\documentclass[12pt,letterpaper]{amsart}

\usepackage[top=2.6cm,
bottom=2.6cm,
left=2.3cm,
right=2.3cm]{geometry}
\usepackage{amssymb}
\usepackage{amsfonts}
\usepackage[hypertexnames=false]{hyperref}

\usepackage[english]{babel}
\usepackage{bm}
\usepackage{setspace}
\usepackage{physics}
\usepackage[style=numeric, maxbibnames=99, maxalphanames=99]{biblatex}
\usepackage{stmaryrd}
\usepackage{upgreek}
\usepackage{setspace}

\usepackage{color}

\usepackage[utf8]{inputenc}
\usepackage{hyperref}
\usepackage{relsize}
\usepackage{bm}
\usepackage{amssymb,amsmath}
\usepackage{dsfont}
\usepackage[all]{xy}
\usepackage{mathtools}
\usepackage{mathrsfs}
\usepackage{comment}
\usepackage{braket}
\usepackage{enumitem}
\usepackage{adjustbox}
\usepackage{combelow}

\renewcommand{\theta}{\uptheta}
\renewcommand{\iota}{\upiota}
\renewcommand{\alpha}{\upalpha}
\renewcommand{\beta}{\upbeta}
\renewcommand{\gamma}{\upgamma}
\renewcommand{\delta}{\updelta}
\renewcommand{\zeta}{\upzeta}
\renewcommand{\pi}{\uppi\hspace{0.05em}}
\renewcommand{\xi}{\upxi}
\renewcommand{\chi}{\upchi}
\renewcommand{\sigma}{\upsigma}
\renewcommand{\Lambda}{\Uplambda}
\renewcommand{\Gamma}{\Upgamma}
\renewcommand{\phi}{\upphi}
\renewcommand{\nu}{\upnu}
\renewcommand{\tau}{\uptau}
\renewcommand{\mu}{\upmu}
\renewcommand{\eta}{\upeta}

\newtheorem{theorem}{Theorem}[section]
\newtheorem{thmx}{Theorem}

\newtheorem{conjx}{Conjecture}

\newtheorem{proposition}[theorem]{Proposition}
\newtheorem{lemma}[theorem]{Lemma}
\newtheorem{conjecture}[theorem]{Conjecture}
\newtheorem{corollary}[theorem]{Corollary}

\theoremstyle{definition}
\newtheorem{definition}[]{Definition}

\theoremstyle{remark}
\newtheorem{example}[theorem]{Example}
\newtheorem{remark}[theorem]{Remark}

\renewcommand{\AA}{\mathbb{A}}

\DeclareMathOperator{\B}{B\!}

\DeclareMathOperator{\Yang}{\mathbf{Y}}

\DeclareMathOperator{\DD}{\mathbb{D}}

\renewcommand{\dd}{\mathbf{d}}
\DeclareMathOperator{\vv}{\mathbf{v}}
\DeclareMathOperator{\ww}{\mathbf{w}}

\DeclareMathOperator{\reldim}{\mathrm{reldim}}

\newcommand{\JH}{\mathtt{JH}}

\DeclareMathOperator{\Grad}{Grad}
\DeclareMathOperator{\Filt}{Filt}

\let \ol=\overline
\let \ul=\underline

\newcommand{\FX}{\mathfrak{X}}

\newcommand{\ff}{\mathbf{f}}

\newcommand{\lazy}{\bm{e}}

\DeclareMathOperator{\opp}{op}
\DeclareMathOperator{\MO}{\mathtt{MO}}

\newcommand{\TT}{\mathbb{T}}

\newcommand{\LL}{\mathbb{L}}
\newcommand{\NN}{\BoN}

\DeclareMathOperator{\GL}{GL}

\DeclareMathOperator{\de}{def}
\DeclareMathOperator{\ob}{ob}
\DeclareMathOperator{\ext}{ext}
\DeclareMathOperator{\Ann}{Ann}

\DeclareMathOperator{\SPL}{Sp}

\newcommand{\Msp}{\mathcal{M}}

\DeclareMathOperator{\sst}{-ss}

\DeclareMathOperator{\nil}{nil}
\DeclareMathOperator{\unit}{unit}
\DeclareMathOperator{\counit}{counit}

\DeclareMathOperator{\vdim}{vdim}

\DeclareMathOperator{\Hilb}{Hilb}

\DeclareMathOperator{\coker}{coker}

\DeclareMathOperator{\crit}{crit}
\DeclareMathOperator{\Dcrit}{\underline{\crit}}

\DeclareMathOperator{\cyc}{cyc}

\DeclareMathOperator{\Hom}{Hom}

\DeclareMathOperator{\End}{End}

\DeclareMathOperator{\Ext}{Ext}

\DeclareMathOperator{\G}{GL}

\DeclareMathOperator{\Map}{Map}
\DeclareMathOperator{\DQM}{\underline{\mathsf{QM}}}
\DeclareMathOperator{\DMap}{\underline{Map}}
\DeclareMathOperator{\QM}{\mathsf{QM}}

\DeclareMathOperator{\Gr}{\mathbf{Gr}}
\DeclareMathOperator{\thGr}{\widehat{\mathbf{Gr}}}

\DeclareMathOperator{\Vect}{Vect}

\DeclareMathOperator{\rep}{Rep}

\DeclareMathOperator{\Ker}{Ker}

\DeclareMathOperator{\BM}{BM}

\newcommand{\fg}{\mathfrak{g}}

\DeclareMathOperator{\Nak}{\mathbf{N}}

\DeclareMathOperator{\bbP}{\mathbb{P}}

\DeclareMathOperator{\Rep}{Rep}

\DeclareMathOperator{\Coh}{Coh}
\DeclareMathOperator{\Bun}{Bun}
\DeclareMathOperator{\FCoh}{\mathfrak{Coh}}
\DeclareMathOperator{\codim}{codim}

\DeclareMathOperator{\Sym}{\mathbf{Sym}}

\DeclareMathOperator{\Spec}{Spec}
\DeclareMathOperator{\Gl}{GL}

\DeclareMathOperator{\gl}{\mathfrak{gl}}
\DeclareMathOperator{\SL}{SL}

\DeclareMathOperator{\id}{id}
\DeclareMathOperator{\Jac}{Jac}

\DeclareMathOperator{\Crit}{Crit}

\DeclareMathOperator{\pt}{pt}
\DeclareMathOperator{\Tot}{Tot}

\DeclareMathOperator{\Char}{char}

\DeclareMathOperator{\rk}{rk}

\DeclareMathOperator{\vir}{vir}
\DeclareMathOperator{\Ob}{Ob}

\DeclareMathOperator{\HO}{\mathbf{H}}

\DeclareMathOperator{\pr}{pr}

\DeclareMathOperator{\sw}{\mathbf{sw}}

\DeclareMathOperator{\Coha}{\mathcal{A}}
\DeclareMathOperator{\HCoha}{\HO\!\mathcal{A}}

\newcommand{\vmult}{\bm{m}}

\newcommand{\vact}{\bm{a}}
\newcommand{\vtau}{\bm{\tau}}

\DeclareMathOperator{\Lie}{Lie}
\renewcommand{\Res}{\mathrm{Res}}
\renewcommand{\ev}{\mathsf{ev}}

\DeclareMathOperator{\Q}{\mathsf{Q}}
\DeclareMathOperator{\Cartan}{\mathsf{C}}

\newcommand{\BPSh}{\mathcal{BPS}}

\newcommand{\Dub}{\mathcal{D}}

\DeclareFontFamily{U}{mathc}{}
\DeclareFontShape{U}{mathc}{m}{it}%
{<->s*[1.03] mathc10}{}
\DeclareMathAlphabet{\mathcal}{U}{mathc}{m}{it}

\newcommand{\BoN}{\mathbf{N}}
\newcommand{\BoD}{\mathbf{D}}
\newcommand{\BoQ}{\mathbf{Q}}
\newcommand{\BoC}{\mathbf{C}}

\newcommand{\BoZ}{\mathbf{Z}}

\newcommand{\CF}{\mathcal{F}}
\newcommand{\CH}{\mathcal{H}}

\newcommand{\calC}{\mathcal{C}}
\newcommand{\calD}{\mathcal{D}}

\newcommand{\CM}{\mathcal{M}}

\newcommand{\CY}{\mathcal{Y}}

\newcommand{\FM}{\mathfrak{M}}
\newcommand{\FP}{\mathfrak{P}}
\newcommand{\FS}{\mathfrak{S}}

\newcommand{\FY}{\mathfrak{Y}}

\newcommand{\Fg}{\mathfrak{g}}

\newcommand{\lfps}{(\mathcal{K})}
\newcommand{\fps}{(\mathcal{O})}
\newcommand{\fpslfps}{(\mathcal{K}, \mathcal{O})}
\newcommand{\lfpslfps}{(\mathcal{K}, \mathcal{K})}

\newcommand{\Rt}{\mathcal{R}}
\newcommand{\Mt}{\mathcal{M}}
\newcommand{\Nt}{\mathcal{N}}
\newcommand{\Tt}{\mathcal{T}}
\newcommand{\Ut}{\mathcal{U}}

\newcommand{\Vt}{\mathcal{V}}

\newcommand{\Coul}{\Sh{C}}
\newcommand{\edge}{\mathtt{edge}}

\newcommand{\fp}{\mathfrak{p}}
\newcommand{\etap}{{}^{\phi}\eta}

\newcommand{\phip}[1]{{}^{\mathfrak{p}}\!\phi_{#1}}

\DeclareMathOperator{\CoHA}{\mathcal{A}}

\DeclareMathOperator{\Stab}{Stab}

\DeclareMathOperator{\fr}{\mathtt{fr}}

\DeclareMathOperator{\Frac}{Frac}

\newcommand{\Sh}[1]{\mathcal{#1}}

\DeclareMathOperator{\loc}{\mathtt{loc}}
\DeclareMathOperator{\gr}{gr}

\DeclareMathOperator{\Hecke}{Hecke}

\DeclareMathOperator{\onil}{\omega-nil}

\DeclareMathOperator{\cofib}{cofib}
\DeclareMathOperator{\DPsi}{\underline{\Psi}}

\usepackage{tikz-cd}
\tikzset{
    labl/.style={anchor=south, rotate=90, inner sep=.5mm}
}
\tikzset{
    labl2/.style={anchor=south, rotate=-90, inner sep=.5mm}
}

\newcommand{\tom}[1]{{\color{blue}Tom: [#1]}}

\DeclareMathOperator{\Quot}{Quot}

\title{Quasimap critical cohomology, Coulomb branches, and quantum groups}
\date{29th August 2026}

\author{Tommaso Maria Botta, Spencer Tamagni}

\address{T. M. Botta: Department of Mathematics, Columbia University, New York City, USA}
\email{tommaso.botta@columbia.edu}
\address{S. Tamagni: Department of Physics, University of California, Berkeley, USA}
\email{stamagni@berkeley.edu}

\newcommand{\paperpart}[1]{%
  \refstepcounter{part}%
  \setcounter{section}{0}%
  \addcontentsline{toc}{part}{Part \thepart.\ \MakeUppercase{#1}}%
  \par\bigskip
  \begin{center}
    {\Large\bfseries PART \thepart.\ #1\par}
  \end{center}
  \medskip
}

\providecommand{\paperpart}[1]{}

\renewcommand{\paperpart}[1]{%
  \refstepcounter{part}%
  \setcounter{section}{0}%
  \addcontentsline{toc}{part}{Part \thepart.\ #1}%
  \par\bigskip
  \begin{center}
    {\normalsize\bfseries
      PART \thepart.\ \MakeUppercase{#1}\par}
  \end{center}
  \medskip
}

\begin{document}
\onehalfspacing

\begin{abstract}
    This article presents a systematic study of the critical cohomology of $\QM^{\xi}(\bbP^1, X)$, the moduli space of based quasimaps from $\bbP^1$ to a Nakajima quiver variety $X$. We develop two derived equivalent presentations of this moduli space as a global critical locus. Accordingly, we obtain two equivalent explicit realizations of the canonical DT sheaf of the moduli space.
    
    Each critical description displays a distinguished aspect of the geometric representation theory of $\QM^{\xi}(\bbP^1, X)$.
    \begin{itemize}
        \item The first model, which builds upon \cite{quasimapcrit}, presents $\QM^{\xi}(\bbP^1, X)$ as the derived critical locus
        of a quiver with potential. We use this description and the theory of cohomological Hall algebras to construct a canonical action of a shifted Yangian on the vanishing cycle cohomology of $\QM^{\xi}(\bbP^1, X)$.
        \item The second model lifts the uniformization theorem for $\Bun_G(\bbP^1)$ to quasimaps $\QM^{\xi}(\bbP^1, X)$, and exhibits the latter as the critical locus of a function on a infinite rank (Tate) bundle over $\Bun_G(\bbP^1)$. We exploit this construction to realize a canonical action of the quantized Coulomb branch defined by Braverman, Finkelberg and Nakajima on the vanishing cycle cohomology of $\QM^{\xi}(\bbP^1, X)$.
        To achieve this, we develop a alternative construction of the BFN Coulomb branch via vanishing cycles.
    \end{itemize}
    We also study the compatibilities between these actions. Specifically, we produce a canonical morphism from the appropriate Hall algebra to the Coulomb branch by geometric means, and show that it extend to a surjective morphism from the shifted Yangian. We conjecture that the shifted Yangian and Coulomb actions are compatible via this morphism, and prove the conjecture under a geometric assumption on the base point $\xi\in X$. We also reduce the conjecture to a sheaf theoretic statement that may be approached using Joyce conjecture.    
\end{abstract}

\maketitle

\setcounter{tocdepth}{2} 
\tableofcontents


\section{Introduction}
\subsection{Overview}
In this paper we study the moduli space of quasimaps $\QM(\mathbb{P}^1, X)$ from $\mathbb{P}^1$ to a Nakajima quiver variety $X$, from the point of view of categorical Donaldson-Thomas (DT) theory and geometric representation theory. From the abstract point of view, the study of $\QM(\mathbb{P}^1, X)$ we undertake in this paper can be put into the framework of DT sheaves on oriented $(-1)$-shifted symplectic stacks (we survey all these notions at more length below). Our main results are of the form that the DT cohomologies of the quasimap spaces are modules over well-studied algebras, namely shifted Yangians and Coulomb branches in the sense of \cite{BFNII}, \cite{BFN_quiver}. 

At the heart of our analysis is the derived algebraic geometry of the space of quasimaps, and in particular the ability to exhibit it as the (derived) critical locus of an explicit function (in fact, in more than one way). A differential-geometric construction of such a critical locus presentation in fact appeared in Nakajima's original proposal \cite{nakajimaCoulombVC} for the definition of the Coulomb branch.\footnote{And is also well-known to physicists studying three-dimensional gauge theories with $\Sh{N} = 4$ supersymmetry.} Let us recall the salient features of this construction, skipping over some finer details. First recall the classical construction of Atiyah and Bott \cite{AtBo83} exhibiting the moduli stack of $G$-bundles on a Riemann surface $C$ in gauge-theoretic language:
\[
\text{Bun}_G(C) = \{ \text{$(0, 1)$ parts of connections} \}/\mathcal{G}
\]
where $\mathcal{G}$ is the (very infinite-dimensional) group of all $C^\infty$ automorphisms of a bundle of given topological type. We will view a $(0, 1)$ part of a connection $A$ via its associated Cauchy-Riemann operator which we call $\overline{\partial}_A$. 

Recall that a Nakajima variety $X$ is open in a stack of the form $T^*(N/G)$, where $G$ is a product of general linear groups and $N$ is a $G$-module arising from a quiver. Therefore, a ($K_C^{1/2}$-twisted) quasimap to $X$ is given by a holomorphic $G$-bundle on $C$, which we parameterize by $\overline{\partial}_A$, together with global holomorphic sections $Q \in H^0(C, \Sh{N} \otimes K_C^{1/2})$, $\widetilde{Q} \in H^0(C, \Sh{N}^* \otimes K_C^{1/2})$ of vector bundles associated to $N$ and $N^*$, which must satisfy the moment map condition $\mu(Q, \widetilde{Q}) = 0$ at all points of $C$, where $\mu$ is the moment map associated to $G$-action on $T^*N$. Moreover $(A, Q, \widetilde{Q})$ must satisfy a stability condition that we suppress for the sake of exposition. Above, $K_C^{1/2}$ denotes a choice of square root of the canonical bundle $K_C$ of $C$. 

In algebraic geometry, a section being holomorphic is often taken for granted, but in the smooth context this condition is in fact enforced by a partial differential equation:
\[
\overline{\partial}_A Q = \overline{\partial}_A \widetilde{Q} = 0.
\]
It is now elementary to observe that the above three conditions (holomorphicity of each section and the moment map) organize into the critical point conditions of the potential function 
\[
\Sh{W} = \int_C \langle \widetilde{Q}, \overline{\partial}_A Q \rangle
\]
where $Q, \widetilde{Q}$ appearing in $\Sh{W}$ are now understood as $C^\infty$-sections of the complex vector bundle $T^* \Sh{N} \otimes K_C^{1/2}$. Unfortunately, the ambient space of triples $(A, Q, \widetilde{Q})$ on which $\Sh{W}$ is defined is too infinite-dimensional to be practically useful to study e.g. the vanishing cycles of $\Sh{W}$. 

In this paper, we study $C = \mathbb{P}^1$ and circumvent this problem by finding global critical locus presentations of the (based) quasimap space which are accessible in finite-dimensional algebraic geometry. One way to approach our critical locus presentations is via the following heuristic. The way to interpret $\Sh{W}$ is that it is associated to the Dolbeault complex resolving the condition of the sections $Q, \widetilde{Q}$ being holomorphic. In a sufficiently homotopy-coherent context, one can imagine replacing the Dolbeault resolution by any other quasi-isomorphic resolution, which might be smaller and more amenable to computation. The quiver model we use in Part \ref{part: quiver model} can be viewed as associated to certain canonical resolutions arising from resolution of the diagonal $\Delta \subset \mathbb{P}^1 \times \mathbb{P}^1$. The Čech model we use in Part \ref{part: Čech model} can be interpreted similarly as arising from Čech resolution using the distinguished coordinate charts on $\mathbb{P}^1$.

In this way, the potentials appearing in our analysis below can be considered as ``the same'' as the above $\Sh{W}$, and correspond to finite-dimensional models of it in algebraic geometry. It would be interesting to make this heuristic more precise. 

We now discuss at more length the precise results and structure of this paper.
\subsection{Quasimaps and their $(-1)$ shifted symplectic structure}
\label{subsec: intro 1}
The main object of study of this article is the moduli space of quasimaps from $\bbP^1$ to a Nakajima quiver variety $X$. Therefore, our presentation must begin with the latter. Fix an arbitrary quiver\footnote{No assumption on the quiver will be ever made throughout this article.} $Q$. Quiver varieties \cite{Nak98, Nak94, Nakajima_tensorI} are, by definition, moduli spaces of framed representations of the preprojective algebra $\Pi_Q$ of the quiver. They are defined as GIT symplectic reductions and—for a generic choice of stability condition—are naturally embedded in an ambient stack quotient $\FX$ as an open substack. The moduli space of quasimaps \cite{CIOCANFONTANINE201417, okounkov2017lectures} to $X$, denoted by $\QM(\bbP^1, X)$, is then defined as the open substack\footnote{In fact, every connected component is a finite type scheme.} of the mapping stack $\Map(\bbP^1,\FX)$ consisting of those maps $f: \bbP^1\to \FX$ that land in $X\subset \FX$ for all but finitely many points in $\bbP^1$. This article concerns the action of certain distinguished symmetry algebras on a \emph{suitable cohomology theory} associated to $\QM(\bbP^1, X)$. Our choice of cohomology theory, or equivalently the choice of its coefficients, is determined by the derived geometry of the moduli space $\QM(\bbP^1, X)$. 

To begin with, we will actually consider a variant of $\QM(\bbP^1, X)$, namely the moduli stack of \emph{based} quasimaps $\QM^{\xi}(\bbP^1, X)$. This is defined as the substack of $\QM(\bbP^1, X)$ consisting of maps $f: \bbP^1\to \FX$ whose evaluation at $\infty\in \bbP^1$ lands in a given point $\xi\in X$. This condition has a certain representation theoretic interpretation, but for us it has an even more fundamental geometric significance. In fact, the canonical derived enhancement of the moduli space $\QM^{\xi}(\bbP^1, X)$ admits a canonical $(-1)$-shifted symplectic structure. This key geometric imput suffices to determine the natural coefficients for a suitable cohomology theory of $\QM^{\xi}(\bbP^1, X)$. Indeed, it is proved in \cite{Br12} that to any $(-1)$-shifted symplectic (oriented) derived algebraic stack $\FY$ it is associated a canonical constructible complex $\varphi_{\FY, o}$ known as the DT sheaf of $\FY$. The second subscript refers to the chosen orientation.

Donaldson-Thomas sheaves of $(-1)$-shifted symplectic stacks are a central objects in both algebraic geometry and geometric representation theory, and are intimately related with categorified enumerative invariants and cohomological Hall algebras \cite{KS2, KS3, da13, DaMe4, kinjo2026cohomologicalhallalgebras3calabiyau}—to which we will come later.

For us, the DT sheaf emerges as the natural coefficient for a cohomology theory of $\QM^{\xi}(\bbP^1, X)$ that is amenable for representation theoretic purposes. In fact their properties are ubiquitous in all our constructions. In essence, this article aims at initiating the study of the geometric representation theory of the DT cohomology
\begin{equation}
    \label{eq: DT cohomology intro}
    \HO_{T}(\QM^{\xi}(\bbP^1, X), \varphi_{\QM^{\xi}(\bbP^1, X)}),
\end{equation}
i.e. of the (suitably $T$-equivariant) derived global sections of the canonical\footnote{The moduli space $\QM^{\xi}(\bbP^1, X)$ admits a canonical orientation (in the CY limit), which we hereby fix once and for all.} DT sheaf of $\QM^{\xi}(\bbP^1, X)$. 

The existence of a canonical $(-1)$-shifted symplectic structure, and hence of a canonical DT sheaf on $\QM^{\xi}(\bbP^1, X)$, can be deduced by very general principles. Firstly, observe that $X$ is a symplectic variety, and indeed it is an open subvariety of the $0$-shifted symplectic stack\footnote{The stack $\FX$ is of the form $T^*( N/G)$, where $N$ is a $G$-representation.} $\FX$. Since $\xi$ is a stable point, we have $\QM^\xi(\bbP^1, X)=\ev^{-1}_{\infty}(\xi)$, where $\ev_{\infty}$ is the evaluation map $\ev_\infty: \Map(\bbP^1, \FX)\to \FX\supset X$. Therefore, $\QM^\xi(\bbP^1, X)$ fits in a pullback diagram
\begin{equation}
    \begin{tikzcd}
        \QM^\xi(\bbP^1, X)\arrow[r]\arrow[d] &  \Map(\bbP^1, \FX)\arrow[d]
        \\
        T_{\xi}\FX\arrow[r] & T\FX=\Map(\Spec \BoC[\epsilon]/\epsilon^2, \FX),
    \end{tikzcd}
\end{equation}
where all the moduli spaces are implicitly upgraded to their natural derived enhancements. Since $\FX$ is $0$-shifted symplectic, we have $T^*\FX \simeq T\FX$ and the map $T_{\xi}\FX\to T\FX$ is naturally a Lagrangian morphism. Moreover, since the pair $(\bbP^1, 2\infty=\Spec \BoC[\epsilon]/\epsilon^2)$ is log Calabi-Yau, it follows that the morphism $\Map(\bbP^1, \FX)\to T\FX$ is also Lagrangian (more specifically, by \cite[Theorem 2.9]{calaquelagrangianmap}). Since the derived intersection of two $n$-shifted Lagrangians is ($n-1$)-shifted symplectic \cite[Theorem 2.9]{pantev2013shiftedsymplecticstructures}, the quasimap space $\QM^{\xi}(\bbP^1, X)$ is naturally $(-1)$-shifted symplectic.

The intrinsic definition of the DT sheaf in terms of the oriented $(-1)$-shifted symplectic structure makes it extremely general and amenable for formal constructions but, as a drawback, it does \emph{not} admit an explicit description. This is a unavoidable compromise of the theory, but effectively a significant limitation for geometric representation theory purposes and a considerable obstacle for practical computations. Therefore, the quest for explicit and computationally effective models for the DT sheaf is a central problem in geometric representation theory, and one of the main results of this article is the construction of \emph{two} independent models for the DT cohomology \eqref{eq: DT cohomology intro}. Each model will highlight an interesting facet of its representation theory. 

The very construction of the DT sheaf \cite{Br12}\footnote{See also \cite{khan2025perversepullbacks} for the relative version, and hence in particular the version over $BT$ that we need to make our statements $T$-equivariant.} hints at a way of producing explicit models. Indeed, a fundamental theorem of Brav-Bussi-Joyce \cite{BBD19}, known as Darboux theorem for $(-1)$-shifted symplectic stacks, states any $(-1)$-shifted symplectic stack $\FY$ is locally modeled as a \emph{derived critical locus}, i.e. as the natural derived enhancements of the critical locus of a function $f: M\to \BoC$ on a smooth ambient stack $M$. Each global critical locus admits a canonical orientation and, on each of them, the DT sheaf is simply defined as the \emph{vanishing cycle functor} $\phip{f}$ of $f$ applied to the (suitably shifted) constant sheaf. In other words, we have\footnote{The sheaf $\BoQ^{\vir}$ is simply the shifted  constant sheaf $\BoQ_M[\dim(M)]$}
\[
\varphi_{\FY, o}=\phip{f}\BoQ^{\vir}_M|_{\Crit(f)}\otimes (\text{local system}).
\]
The twist is determined by the orientation $o$ of $\FX$, which serves to glue the various patches together. A fundamental result of \cite{Br12} is then the statement that this construction is intrinsic of the oriented $(-1)$-shifted symplectic structure, i.e. it is independent of the choice of Darboux charts.

In conclusion, it is the gluing process that makes the DT sheaf inexplicit. And hence its solution is—res ipsa loquitur—avoiding any gluing, that is finding a presentation of the moduli space as a global critical locus. This is in general hopeless for an arbitrary $(-1)$-shifted symplectic stack, but turns out to be possible for the moduli space $\QM^{\xi}(\bbP^1, X)$. In fact, this article builds on on two different global critical models, which we now introduce.

\subsection{Quiver model, Hecke modifications, and Hall algebras}
\label{subsec: intro 2}
The first model, which was originally developed in \cite{quasimapcrit}, is of quiver type, i.e. it exhibits the moduli space of quasimaps as a moduli space of representations of a quiver with potential. Specifically, the moduli space of quasimaps $\QM^{\xi}(\bbP^1, X)$ admits a decomposition in connected components
\[
\QM^{\xi}(\bbP^1, X)=\bigsqcup_{\dd}\QM^{\xi}_{\dd}(\bbP^1, X),
\]
where $\dd$ is the topological type recording the effective curve class, and the main theorem of \cite{quasimapcrit} exhibits a natural isomorphism between $\QM^{\xi}_{\dd}(\bbP^1, X)$ and the space of $\dd$-dimensional representation of a distinguished quiver $Q^{\QM}$ with potential $W$. In other words, it exhibits $\QM^{\xi}_{\dd}(\bbP^1, X)$ as the classical critical locus of the function \footnote{Technically, $\FM_{\dd}(Q^{\QM})$ is the quotient of $\Rep_{Q^{\QM}}(\dd)$ by the reductive group $G=\prod_{i}\Gl_{\dd_i}$, while $\QM^{\xi}_{\dd}(\bbP^1, X)$ is the quotient of the same space by a unipotent extension of $\Gl_{\dd_i}$, cf. \S\ref{subsec: quasimaps to quiver varieties}. However, it will turn out that this does not cause any trouble.}
\[
\tr(W): \FM^{\sst}_{\dd}(Q^{\QM})\to \BoC,
\]
where $ \FM^{\sst}_{\dd}(Q)$ is the (appropriate semistable locus of the) moduli space of $\dd$-dimensional representations of the quiver $Q$. We stress that the pair $(Q^{\QM}, W)$ solely depends on the pair $(X, \xi)$. Since the latter are kept fixed in all our construction, the quiver model provides a uniform description of $\QM^{\xi}(\bbP^1, X)$. 

In a nutshell, the quiver model is obtained as follows. Recall that a quasimap to a quiver variety is a map from $\bbP^1$ to a stack quotient of the form $T^*(N/G)$ where $G$ is a product of general linear groups, and hence is a moduli space parametrizing a collection of vector bundles over $\bbP^1$ with sections (depending on $N$) satisfying a moment map condition. Informally speaking, the quiver model exploits the Beilinson resolution to trade vector bundles with vector spaces and sections with linear maps between them. This gives rise to a moduli space of quiver representations $\FM_{\dd}(Q^{\QM})$. The fundamental observation is then that the relations among the linear maps enforced by the relations among the sections \emph{combined} with the constraints that the linear map themselves come from sections naturally assemble into the critical locus of a potential function $W$. 

The first main result of this paper refines this result by proving that this isomorphism lifts to the category of oriented $(-1)$-shifted symplectic stacks. Specifically, we have
\begin{thmx}
\label{main thm: derived QM quiver}
    The intrinsic $(-1)$-shifted symplectic structure on $\QM^\xi( \bbP^1, X)$ agrees, in the quiver model, with the canonical $(-1)$-shifted symplectic structure on the derived critical locus of the quasimap potential $\tr(W): \FM_{\dd}(Q^{\QM})\to \BoC$. Moreover, the canonical orientations also match, so we have a canonical isomorphism of sheaves $\phip{W}\BoQ^{\vir}_{\FM_{\dd}(Q^{\QM})}\cong \varphi_{\QM^\xi(\bbP^1, X)}$.
\end{thmx}

Thanks to this theorem, the DT cohomology \eqref{eq: DT cohomology intro} can be accessed by equivalently considering the more manageable \emph{critical cohomology} of the quiver with potential $(Q^{\QM}, W)$, that is the cohomology 
\begin{equation}
\label{eq: crit cohomology quiver model intro}
    \HO(\FM_{\dd}(Q^{\QM}), \phip{W}\BoQ^{\vir}).
\end{equation}
This embeds the study of \eqref{eq: DT cohomology intro} into the vast and ever growing literature of critical Hall algebras  \cite{da13}, convolution algebras \cite{VVcriticalconvolution}, and stable envelopes \cite{COZZ1, COZZ2, BDZZ}, and gives a hands-on approach to the representation theory naturally attached to these moduli spaces. However, before getting there, it is useful to deduce the geometric nature of the algebra of correspondences acting on \eqref{eq: DT cohomology intro} or, equivalently \eqref{eq: crit cohomology quiver model intro}. As previously outlined, the moduli space of based quasimaps can be informally described as 
\[
\QM^{\xi}(\bbP^1, X)=\lbrace (\Sh{V},s) \; + \text{moment map condition on $s$}\rbrace /\sim
\]
where $\Sh{V}$ is a tuple of bundles, $s$ is a section of an associated bundle, and the moment map condition is inherited from $T^*(N/G)=\mu^{-1}(0)/G$. It is then natural to anticipate a geometric action of a symmetry algebra induced by the moduli stack of Hecke correspondences
\begin{equation}
    \label{eq: hecke stack intro}
    \Hecke^{\xi}(X)=\lbrace (\Sh{V},\Sh{V}',s, s') \; | \; \Sh{V}' \subset \Sh{V} + \text{ compatibility of sections} \rbrace /\sim
\end{equation}
This stack admits two canonical maps to $\QM^{\xi}(\bbP^1, X)$, corresponding to the two pairs $(\Sh{V},s)$ and $(\Sh{V}',s')$, and hence naturally induces a map between a suitable choice of cohomology theory—which as we anticipated, is DT cohomology. Amusingly, the algebra of symmetry naturally emerges by closely inspecting \eqref{eq: hecke stack intro}. Indeed, $ \Hecke^{\xi}(X)$ admits a third canonical map, this time to $\FCoh_{\Pi_Q, 0}(\AA^1)$ the moduli stack of representations of the preprojective algebra $\Pi_Q$ \emph{valued} in coherent sheaves on $\AA^1$. Concretely, this map simply emerges by the observation that, taking the cokernel to the inclusion $\Sh{V}' \subset \Sh{V}$ we obtain a (rank zero) coherent sheaf on\footnote{Since we consider based quasimaps, the value of the map at $\infty\in \bbP^1$ is fixed and hence the Hecke modification must live on its complement $\AA^1$.} $\AA^1$.
The compatible sections then descend to morphisms on the quotients, and the moment map condition enforces the relation of $\Pi_Q$. Remarkably, the abelian category of $\Pi_Q$-modules valued in rank zero coherent sheaves on $\AA^1$ is equivalent to the category of representations of the Jacobi algebra $\Jac(\tilde W, \tilde Q)$ of the tripled quiver with its canonical cubic potential, cf. \S\ref{subsec: Doubling, looping, tripling, and framing}. This statement
implies that the moduli stack $\FCoh_{\Pi_Q, 0}(\AA^1)$ admits a presentation as a quiver with potential—and indeed one of the most famous and studied ones \cite{preproj, DHSM23, BD2023, SV23O}.

The critical cohomology of this quiver with potential, denoted $\HCoha^T_{\tilde Q, \tilde W}$, possesses a distinguished algebra structure making it into a so-called critical Hall algebra, whose construction was outlined in \cite{KS2} and then extensively studied by Davison \cite{da13}. Not surprisingly, it turns out that this is the the natural symmetry algebra acting on quasimap cohomology via \eqref{eq: hecke stack intro}. Possibly more surprising is that the ambient geometries describing $\HCoha^T_{\tilde Q, \tilde W}$ and \eqref{eq: crit cohomology quiver model intro} are perfectly compatible as the pair $(Q^{\QM},W)$ can be realized as a \emph{framing} of the quiver with potential $(\tilde Q,\tilde W)$ determining $\HCoha^T_{\tilde Q, \tilde W}$. This observation makes it possible to describe the Hall action on \eqref{eq: crit cohomology quiver model intro} via vanishing cycles and critical cohomology.

At this point it is worth stressing that the existence of a Hall action can be also derived abstractly by arguing that the correspondence 
\[
\FCoh_{\Pi_Q, 0}(\AA^1) \times \QM^{\xi}(\bbP^1, X) \mapsfrom \Hecke^{\xi}(X) \mapsto \times \QM^{\xi}(\bbP^1, X)
\]
naturally lifts to a $(-1)$ shifted Lagrangian morphism; the action of the DT cohomology of $\FCoh_{\Pi_Q, 0}(\AA^1)$ on \eqref{eq: DT cohomology intro} is then a formal consequence of  \cite{kinjo2026cohomologicalhallalgebras3calabiyau}. However, the quiver model is the key for extending this Hall action to an action of an even more interesting algebra: a shifted Yangian.

Shifted Yangians, and more generally shifted quantum groups, are generalizations of the well-known Yangians and quantum groups that play a fundamental role in the study of quantum integrable systems and related structures. Concretely, they may be realized as \emph{Drinfeld doubles} of cohomological Hall Algebras, that is they can be expressed as
\begin{equation}
\label{eq: Yangian intro}
    \Yang_{\mu}^{T}(\tilde Q, \tilde W):=\HCoha_{\tilde Q, \tilde W}^{T,\nil}\ast (\text{Cartan})\ast \HCoha_{\tilde Q, \tilde W}^{T,\nil, \opp}/(\text{explicit set of relations}).
\end{equation}
Here, $\HCoha_{\tilde Q, \tilde W}^{T,\nil}$ is the so-called \emph{nilpotent} subalgebra of $\HCoha_{\tilde Q, \tilde W}^T$, ``Cartan'' refers to a certain commutative algebra, and $\ast$ denotes the free product of algebras. We are now ready to state the second main result of this paper.
\begin{thmx}
\label{thmx: Cohas and shifted Yangians action}
    There exists a canonical action of the cohomological Hall algebra $\HCoha_{\tilde Q, \tilde W}^{T,\nil}$ on 
    \[
    \HO_T(\FM_{\dd}(Q^{\QM}), \phip{W}\BoQ^{\vir})\cong \HO_{T}(\QM^{\xi}(\bbP^1, X), \varphi_{\QM^{\xi}(\bbP^1, X)})
    \]
    Moreover, this action naturally extends to an action of the shifted Yangian $\Yang_{\mu}^{T}(\tilde Q, \tilde W)$, where the shift $\mu$ is determined by the dimension and framing vector defining $X$.
\end{thmx}
We reiterate that the quiver presentation of $\QM^{\xi}(\bbP^1, X)$ plays a key role in proving the second statement as the relations in \eqref{eq: Yangian intro} are explicitly checked using an intersection-theoretic argument involving the ambient space $\FM(Q^{\QM})$. Theorem \ref{thmx: Cohas and shifted Yangians action} confirms (and generalizes) a conjecture of Costello \cite{costello2017holographykoszuldualityexample}. 

As a byproduct of this theorem, we obtain the action of various variants of the Hall algebra $\HCoha_{\tilde Q, \tilde W}^{T,\nil}$. Of particular geometric significance is the action of the so-called $\omega$-nilpotent Hall algebra $\HCoha_{\tilde Q, \tilde W}^{T,\onil}$, originally introduced in \cite{preproj}. This algebra fits in a chain of strict inclusions\footnote{For sufficiently generic $T$}
\[
\HCoha_{\tilde Q, \tilde W}^{T,\nil}\subset \HCoha_{\tilde Q, \tilde W}^{T,\onil}\subset \HCoha_{T, \tilde Q, \tilde W}
\]
and, from the geometric perspective, corresponds to imposing that the Hecke modification acting on quasimaps via \eqref{eq: hecke stack intro} is supported at $0\in \bbP^1$. We remark that this is a natural condition as we will always work equivariantly with respect to the natural $\BoC^*$-action on $\bbP^1$ so the only fixed point other than the base point $\infty$ is $0$. This locality condition naturally leads us to the next chapter of this paper.

\subsection{Cech Model and Coulomb branch action}
\label{subsec: intro 3}

The natural geometric object encoding local Hecke modifications of vector bundles and, more generally, principal $G$ bundles, is the affine Grassmannian $\Gr_G$. The affine Grassmannian plays a central role in geometric representation theory through the geometric Satake correspondence, which connects its geometry with representations of the Langlands dual group. A distinctive feature of the affine Grassmannian is that its category of coherent sheaves, and hence, via decategorification, its K-theory and Borel-Moore homology, admits a remarkable convolution product structure, which naturally emerges by concatenating Hecke modifications.

Although abstractly defined as the moduli space parametrizing Hecke modifications of a trivial $G$-bundle on the formal disk $\DD$, the closed points in the affine Grassmannian can be presented as the quotient space\footnote{Most of the literature opts for the right quotient description $G\lfps/G\fps$. Our choice of departing from it is determined by wanting the Coulomb branch acting on \eqref{eq: DT cohomology intro} from the left, see \S\ref{subsec: Coulomb branch action on quasimaps}.} $G\fps\backslash G\lfps$, where $\Sh{O}$ is the coordinate ring on a formal disk $\DD$—morally centered at the support of the Hecke modification— and $\Sh{K}$ is its function field. The correspondence
\begin{equation}
    \label{eq: conv diag gr intro}
    G\fps\backslash G\lfps\times G\fps\backslash G\lfps\leftarrow  G\fps\backslash G\lfps  \times G\lfps \rightarrow G\fps\backslash G\lfps  \times_{G\fps} G\lfps  \to G\fps\backslash G\lfps
\end{equation}
describes the monoidal structure on the Satake category, which is dual to the natural monoidal structure on $\text{Rep}(G^L)$ \cite{Lusztig1983}, \cite{Ginzburg1995}, \cite{MirkovicVilonen2007}. The same convolution diagram induces a convolution algebra structure on the $G\fps$-equivariant Borel-Moore homology 
\[
\Coul_{G}=\HO^{\BM}_{G\fps}(\Gr_G).
\]
The moduli space of principal $G$ bundles on $\mathbb{D}$ can be presented as the mapping stack $\Map(\DD, \B G)$, so it is natural to seek for generalizations of the convolution algebra $\Coul_{G}$ by considering mapping stacks of the form $\Map(\DD, N/G)$, where $N$ is a $G$-representation---that is to say, searching for moduli spaces parameterizing Hecke modifications of bundles \textit{with sections}. This problem has originally been addressed by Braverman, Finkelberg, and Nakajima \cite{BFNII}, \cite{braverman2018coulombbranches3dmathcal}  in order to provide a mathematical definition of the \emph{Coulomb branch} of a 3d $\Sh{N}=4$ supersymemtric gauge theory.\footnote{With matter representation of cotangent type, i.e. excluding the possibility of half-hypermultiplets.} In \emph{loc.cit.} the authors consider a convolution algebra structure on the moduli space $\Rt$ parametrizing triples $(\Sh{P},n,g)$ where $\Sh{P}$ is a $G$-bundle on the formal disk $\DD$, $g$ is a trivialization over the punctured disk $\DD^\times$, and $n$ is a section of the associated bundle $\Sh{P}\times_G N$ that is sent to a regular section of a trivial bundle under $g$. For $N=0$ one gets back $\Rt=\Gr_G$. One of the main results of their work is the construction of a natural convolution algebra structure on $\HO^{\BM}_{G\fps}(\Rt)$ via a generalization of \eqref{eq: conv diag gr intro}. In absence of torus action, this algebra is commutative, and its spectrum
\begin{equation}
    \label{eq: symplectic dual intro}
    X^!:=\Spec(\HO^{\BM}_{G\fps}(\Rt))
\end{equation}
is, by definition, the classical \emph{Coulomb branch} of the theory. The commutativity fails in equivariant cohomology, and the equivariant parameter $\hbar$ for the $\BoC^*$-action on $\DD$ is identified as the quantization parameter. These varieties provide the mathematical definition of the symplectic dual variety to the quiver variety $X$ introduced in \S\ref{subsec: intro 1}.

The following longstanding conjecture, formulated by Bullimore-Dimofte-Gaiotto-Hilburn-Kim \cite{bullimore2018vorticesvermas} and, independently, Okounkov, states the following:
\begin{conjx}
\label{conjx: intro}
    The quantized Coulomb branch algebra $\BoC[X^!]_\hbar$ acts on the equivariant cohomology of the moduli space $\QM^{\xi}(\bbP^1, X)$.
\end{conjx}
In this article, we provide a refinement of this conjecture alongside with its proof. From a mathematical standpoint, it is natural to anticipate $\HO^{\BM}_{G\fps}(\Rt)$ to act on quasimap cohomology. In fact, we have seen in \S\ref{subsec: intro 2} that Hecke correspondences naturally act on quasimaps, and the moduli space $\Rt$ is by definition, a ``Steinberg variety'' of Hecke correspondences. However, two technical points stand in the way of turning this heuristic into a precise statement:
\begin{enumerate}
    \item Firstly, the moduli space $\Rt$ parametrizes Hecke modifications for $\Map(\DD, N/G)$, while $\QM(\bbP^1, X)$ is an open substack of 
    \[
    \Map(\bbP^1, T^* (N/G))= \Map(\bbP^1, \mu^{-1}(0)/G))
    \]
    and \emph{not} of $\Map(\bbP^1, N/G)$.
    \item As advocated in \S\ref{subsec: intro 2}, the natural cohomology theory for $\QM^{\xi}(\bbP^1, X)$ is DT cohomology, and not Borel-Moore homology, which is the underlying theory for the Coulomb branch \eqref{eq: symplectic dual intro}. 
\end{enumerate}
In this work, we solve both problems by lifting the definition of Coulomb branch to critical cohomology and providing a global critical locus model of  $\QM^{\xi}(\bbP^1, X)$ compatible with it. Our construction is similar in spirit to Nakajima's original proposal for the Coulomb branch \cite{nakajimaCoulombVC}, by the heursitic reviewed above. Here, we reinterpret this idea by constructing a space $\tilde \Rt$ together with a potential $\Sh{W}:\tilde \Rt\to \BoC$ satisfying\footnote{$\mathbb{B}$ denotes the ``bubble'' or ``raviolo'' $\mathbb{B} = \DD \cup_{\DD^\times} \DD$.}
\begin{equation}
    \label{eq: crit pot coul}
    \Crit(\Sh{W})/G\fps=\Map(\mathbb{B},\mu^{-1}(0)/G).
\end{equation}
and use it to provide an alternative definition of the Coulomb branch. This leads us to our third main result. 
\begin{thmx}
\label{thmx: critical raviolo}
    There exists a natural convolution algebra structure on the (suitably defined) critical cohomology
    \begin{equation}
        \label{eq: crit Coulomb}
        \HO_{G\fps}(\tilde \Rt,\phip{\Sh{W}}\DD\BoQ)
    \end{equation}
    This convolution algebra structure admits natural quantization and flavor deformation. Additionally, there exists a natural dimensional reduction isomorphism 
    \[
    \HO_{G\fps}(\tilde \Rt,\phip{\Sh{W}}\DD\BoQ)\cong \HO^{\BM}_{G\fps}(\Rt)
    \]
    that is compatible with the algebra structures, quantizations, and flavor deformations of two sides.
\end{thmx}
In plain English, this theorem provides an equivalent description of the Coulomb branch in terms of critical cohomology. Moreover, thanks to \eqref{eq: crit pot coul}, it solves the first of the two points above. We also remark that this construction works for arbitrary theories of cotangent type\footnote{It is not immediate to generalize the construction given in this paper to non-cotangent type, but seems likely such generalization should exist.}, not just quiver gauge theories.

As this stage we would like to remark that, since $\tilde \Rt$ is an ind-scheme of infinite type, and hence in particular it is \emph{not} quasi-compact, the definition of \eqref{eq: crit Coulomb} is understood in terms of finite approximations of $\tilde \Rt$, $G\fps$, and of the potential $\Sh{W}$. The final construction does not depend on such choice of approximation. In this regard, our approach parallels the one in \cite{BFNII}. Additionally, we remark that because of the finite dimensional approximations, the vanishing cycle $\phip{\Sh{W}}$ is not applied to a smooth variety\footnote{The approximations of $\Gr_G$ are given by Schubert varieties, which are not smooth.}, and hence at these stage it cannot be interpreted \eqref{eq: crit Coulomb} as the DT sheaf of $(-1)$-shifted symplectic stack. It would be very interesting to make precise sense of this statement.

We now discuss point (2) above, which is also addressed in a similar fashion. Namely, we exhibit $\QM^{\xi}(\bbP^1, X)$ as the critical locus of a pair $(\Mt, \Sh{W}^{\xi})$ that is \emph{manifestly compatible} with $(\Rt, \Sh{W})$. The fundamental idea behind this construction is based on the Beauville-Laszlo uniformization theorem, which exhibits the moduli space $\Bun_{G, \infty}(\bbP^1)$ of $G$-bundles on $\bbP^1$ with trivialization at infinity as the double quotient 
\[
G\fps \backslash G\lfps/ G_1[\AA^1_{\infty}]=\Gr_G/ G_1[t^{-1}],
\]
where $\AA^1_{\infty}=\bbP^1\setminus \{0\}$ and $G_1(\AA^1_{\infty})$ is the kernel of the natural evaluation $G_1(\AA^1_{\infty})\to G$. In this article, we provide a uniformization theorem for the moduli space $\QM^{\xi}(\bbP^1, X)$ by exhibiting the latter as the critical locus of a function $\Sh{W}^{\xi}:\Mt\to \BoC$. Informally speaking, the space $\Mt$ lives over $\Bun_{G, \infty}(\bbP^1)$ and parametrizes pairs of quasimap data to $(T^*N)/G$ defined on the disk $\DD$ and on the affine line $\AA^1_{\infty}$. The potential enforces the factorization through $T^*(N/G)=\mu^{-1}(0)/G$ as well as the condition that the quasimap data on $\DD$ and $\AA^1_{\infty}$ glue to produce a quasimap on $\bbP^1=\DD\cup_{\DD^\times}\AA^1_{\infty}$. We refer to such critical model as the Čech model. 

Similarly to the case of the Coulomb branch, we produce finite type approximations of the pair $(\Mt, \Sh{W}^{\xi})$, and hence provide a well-defined notion of the critical cohomology 
\[
\HO_T(\Mt, \phip{\Sh{W}}\BoQ^{\vir}_{\Mt}).
\]
We remark that, unlike for the Coulomb branch, our finite type approximations are smooth so $\phip{\Sh{W}}\BoQ^{\vir}$ is a natural candidate for the canonical DT sheaf of $\QM^{\xi}(\bbP^1, X)$. 
\begin{thmx}
    \label{thmx: derived critical Čech}
    The intrinsic $(-1)$-shifted symplectic structure on $\QM^\xi(X)$  agrees, in the Čech model, with the canonical $(-1)$-shifted symplectic structure on the derived critical locus of the Čech potential $\Sh{W}^{\xi}: \Mt\to \BoC$. Moreover, the canonical orientations also match, so we have a canonical isomorphism $\phip{\Sh{W}}\BoQ^{\vir}_{\Mt}\cong \varphi_{\QM^\xi(\bbP^1, X)}$.
\end{thmx}
A remarkable feature of the Čech model is that it is manifestly compatible with the critical description of the Coulomb branch from Theorem \ref{thmx: critical raviolo}. Therefore, building on this result and Theorem \ref{thmx: derived critical Čech} we arrive at the following result:
\begin{thmx}
\label{thmx: Coulomb aciton intro}
    The critical cohomology
    \[
    \HO_T(\Mt, \phip{\Sh{W}}\BoQ^{\vir}_{\Mt})\cong \HO_{T}(\QM^{\xi}(\bbP^1, X), \varphi_{\QM^{\xi}(\bbP^1, X)})
    \]
    carries a natural action of the quantized, flavor deformed BFN Coulomb branch algebra $\Coul^T_{N,G}$ or, equivalently, of the quantized flavor deformed coordinate ring of the symplectic dual variety $X^!$.
\end{thmx}
We regard this result as the adequate affirmative answer to Conjecture \ref{conjx: intro}.

\subsection{Coulomb branches as truncated shifted Yangians and their actions}
Throughout \S\ref{subsec: intro 2} and \S\ref{subsec: intro 3} we have discussed two distinguished algebra actions on the same $\HO_T:=\HO_T(\pt)$-module $\HO_{T}(\QM^{\xi}(\bbP^1, X), \varphi_{\QM^{\xi}(\bbP^1, X)})$. In this section, we discuss their compatibilities. At the algebra level, our results can be seen as a generalization—and geometrization—of the known statement that Coulomb branches of $ADE$ type can be realized as quotients of shifted Yangians \cite{BFN_quiver, KWWY,KTWWY}.

Our results hold for an arbitrary quiver $Q$ and emerge from the following simple geometric observation. For expository purposes, assume momentarily that $G=\GL(n)$. In type $A$, the affine Grassmannian admits a distinguished subvariety $\Gr^+_G\subseteq \Gr_G$, sometimes known as the positive affine Grassmannian, whose cohomology $\HO_{G\fps}(\Gr^+)$ is naturally a subalgebra of $\HO_{G\fps}(\Gr)$. As a moduli space $\Gr^+$ parametrizes pairs $\Sh{V}\subset \Sh{V'}$ of sub-bundles of the same rank on the disk $\DD$. Therefore, by taking quotients, we get a canonical morphism $\Gr^+_G\to \FCoh_{0}(\DD)$, the moduli stack of rank zero coherent sheaves on the disk. The latter has a Hall algebra structure and this construction provides, via pullback, the sought after compatibility the latter and $\HO_{G\fps}(\Gr^+)$. In this article, we show that this picture generalizes to arbitrary quivers. 
\begin{thmx}
\label{thmx: map hall Coulomb}
    There exists a canonical morphism of $\HO_T$-algebras $\theta: \HCoha^{T, \onil}_{\tilde Q, \tilde W} \to \Coul^T_{N,G}$ extending to a morphism of algebras $\Yang_{\mu}^{T}(\tilde Q, \tilde W)\mapsto \Coul^T_{N,G}$ that becomes surjective after applying $\Frac(\HO_T)\otimes(-)$.
\end{thmx}
From a technical standpoint, this result is obtained by first dimensionally reducing $\HCoha^{T, \onil}_{\tilde Q, \tilde W}$ to a two dimensional Hall algebra for $\Rep_{Q}(\Coh_0(\DD))$, the category of representations of $Q$ valued in $\Coh_0(\DD)$, and then proving that the canonical morphism from the positive space of triples $\Rt^+$ to the moduli stack of $\Rep_{Q}(\Coh_0(\DD))$ is quasi-smooth, so there is a natural notion of virtual pullback supplementing the map $\theta$. Currently, we do not see a straightforward way to compare these algebras using operations in critical cohomology.

We also remark that, while drafting the article, the work by Jindal-Negut \cite{jindal2026loopnilpotentcohomologicalhallalgebra} has appeared. In their work, the authors construct a surjective morphism of algebras from a suitable Drinfeld double of $\HCoha^{T, \onil}_{\tilde Q, \tilde W}$ to $\Coul^T_{N,G}$ by combinatorial methods—specifically via shuffle algebras. We expect our construction, which is purely geometric, to agree with theirs. In fact, it is already clear that these two definitions agree on the spherical part of $\HCoha^{T, \onil}_{\tilde Q, \tilde W}$\footnote{This can be deduced by comparing the explicit formulas in \cite{jindal2026loopnilpotentcohomologicalhallalgebra} with those in Lemma \ref{lemma: formulas generators hall in Coulomb}.}. 
 
Combining Theorem \ref{thmx: map hall Coulomb} with Theorem \ref{thmx: Coulomb aciton intro} we obtain an action of the shifted Yangian $\Yang_{\mu}^{T}(\tilde Q, \tilde W)$ on the DT cohomology that is a priori different from the one of Theorem \ref{thmx: Cohas and shifted Yangians action}. However, the same geometric nature of the two actions strongly suggests the following conjecture.

\begin{conjx} $ $

\label{conjx: diagram}
\begin{enumerate}
    \item 
    \label{conj 1 intro}
    The action of the Hall algebra $\HCoha^{T, \onil}_{\tilde Q, \tilde W}$ on $\HO_T(\QM^{\xi}, \varphi_{\QM^{\xi}(\bbP^1, X)})$ induced by the composition of the algebra morphism $\theta: \HCoha^{T, \onil}_{\tilde Q, \tilde W}\mapsto \Coul^T_{G,N}$ with the Coulomb action from Theorem \ref{thmx: map hall Coulomb} coincides with the action from Theorem \ref{thmx: Cohas and shifted Yangians action}. Equivalently, we have a commutative diagram of algebras:
    \begin{equation}
    \label{eqx: diagram conjecture}
        \begin{tikzcd}[column sep=4em]
            \HCoha^{T,\onil}_{\tilde Q, \tilde W} \arrow[r, "\text{Thm. \ref{thmx: Cohas and shifted Yangians action}}"]\arrow[d, swap, "\theta"] & \End\left(\HO_T(\QM^{\xi}(\bbP^1, X), \varphi_{\QM^{\xi}(\bbP^1, X)})\right)
            \\
            \Coul^T_{G,N}\arrow[ur, swap, "\text{Thm. \ref{thmx: Coulomb aciton intro}}"]
        \end{tikzcd}
    \end{equation}
    \item 
    \label{conj 2 intro}
    The action of the Yangian $\Yang^{T}_\mu(\tilde Q, \tilde W)$ induced by the composition of the algebra morphism $\Yang_{\mu}^{T}(\tilde Q, \tilde W)\mapsto \Coul^T_{N,G}$ with the Coulomb action from Theorem \ref{thmx: map hall Coulomb} coincides with the action from Theorem \ref{thmx: Cohas and shifted Yangians action}. Equivalently, the diagram above commutes with  $\HCoha^{T,\onil}_{\tilde Q, \tilde W}$ replaced by $\Yang^{T}_\mu(\tilde Q, \tilde W)$.
\end{enumerate}
\end{conjx}
Restated in sheaf language, the conjecture states that the Hall and Coulomb actions, which are both given by correspondences supported on $\Hecke^{\xi}(X)$, but are realized in terms of \emph{different} global critical models—and hence a priori depend on the ambient spaces—are actually the same. We refer to Remark \ref{Rem: intrinsic conj} for a detailed discussion. This phrasing is very close
to being (a special case of) the Joyce conjecture \cite[Conj.1.1]{JS_2020}. The latter is the object of a work in preparation by Khan-Kinjo-Park-Safronov \cite{KKPSJoyce}. We hope their results will shed light on Conjecture \eqref{conjx: diagram}.

We remark that conjecture \eqref{conj 1 intro}, if appropriately lifted to the sheaf level, implies conjecture \eqref{conj 2 intro} since at the sheaf level the lowering operators are recovered by the raising operators via Verdier duality.

We also have a partial positive answer to Conjecture \eqref{conjx: diagram}, to which we come now. Informally speaking, the subtlety of this conjecture lies in the lack of a uniform language to treat all the morphisms involved in diagram \eqref{eq: diagram conjecture}. In fact, the downward pointing map is constructed in Borel-Moore homology, while the right points ones are constructed in critical cohomology, and in fact with two completely different models—the quiver model and the Čech model. At this point it is important to stress that, to the best of our knowledge, it is  generally \emph{not} possible to dimensionally reduce the critical cohomology  $\HO_T(\QM^{\xi}(\bbP^1, X), \varphi_{\QM^{\xi}(\bbP^1, X)})$ to Borel-Moore homology. There is, however, a particular case when dimensional reduction becomes possible, namely whenever the evaluation point $\xi\in X$ is polarized, that is it lies in a Lagrangian subvariety of the form $L=X\cap N/G$. In this case we have the following:
\begin{thmx}
    \label{propx}
    Assume that $\xi\in X$ is polarized. Then we have a canonical isomorphism 
    \[
    \HO_T(\QM^{\xi}(\bbP^1, X), \varphi_{\QM^{\xi}(\bbP^1, X)})\cong \HO_T(\QM^{\xi}(\bbP^1, L), \DD\BoQ^{\vir}_{\QM^{\xi}(\bbP^1, L)}).
    \]
    Moreover, conjecture \ref{conjx: diagram} holds.
\end{thmx}
The proof is obtained by performing dimensional reduction to all the arrows in \eqref{eqx: diagram conjecture} and hence reducing the statement to Borel-Moore homology. At this stage all maps are constructed geometrically using the geometry of the underlying moduli spaces—without any further reference to the critical models, and the correspondences defining the action can be directly compared. In current work in progress, we plan to prove the conjecture in full generality by constructing an intrinsic notion of stable envelope for $(-1)$-shifted symplectic schemes.

\subsection{Further directions and related works}

\subsubsection{Stable envelopes} 
In this article, we describe Coulomb branch and quantum group actions by explicitly analyzing the generating correspondences. A complementary and very powerful approach is provided by the geometric $R$-matrix formalism of \cite{MO19}, pioneered in critical cohomology in the recent articles \cite{COZZ1}, \cite{COZZ2}. The quiver model of this paper can be thought of as placing the modules 
\[
\HO_T(\QM^\xi(\mathbb{P}^1, X), \varphi_{\QM^\xi(\mathbb{P}^1, X)})
\]
on essentially equal footing with the cohomologies of Nakajima varieties themselves $\HO_T(X)$ which were the focus of \cite{MO19}. In particular, a distinguished role is expected to be played by the geometric $R$-matrices braiding the cohomology of a Nakajima variety itself with the critical cohomology of quasimaps. For quivers of type $A$, these define quantizations of so-called monopole scattering matrices \cite{bullimore2015coulombbranch3dmathcaln4}, otherwise known as matrix descriptions of (generalized) affine Grassmannian slices \cite{BFN_quiver}, \cite{krylov2021dominantgeneralizedslicesconvolution}. 

Assuming the validity of Conjecture \eqref{conjx: diagram}, the $R$-matrix technology of \cite{COZZ2} will define\footnote{Strictly speaking, one also needs some assumption about faithfulness of the Coulomb branch action in these modules to say the Coulomb branch coproduct is uniquely fixed by this consideration, however Conjecture \eqref{conjx: diagram} is enough to relate the Coulomb branch actions in these classes of modules.} a notion of Coulomb branch comultiplication for any quiver gauge theory attached to an antidominantly shifted Yangian, extending certain results of \cite{finkelberg2017comultiplicationshiftedyangiansquantum}, \cite{Frassek_2022}, \cite{kamnitzer2026categorymathcalotruncatedshifted} to more general quivers.

It is interesting to experiment with the dominantly shifted case, as the scattering matrices appear to be well-defined and match certain features visible from the geometry of the quasimap spaces. We expect that the geometric $R$-matrix for $\mathsf{Y}_{+1}(\widehat{\mathfrak{gl}}_1)$ braiding a Fock representation with the module we study in this paper for $X = \text{Hilb}_n(\mathbb{C}^2)$ will describe quantization of the operator studied in Section 5 of \cite{tamagni2026scatteringtransformnoncommutativeinstantons}; more generally, for affine type $A$ quivers resulting operators should be described as quantizations of instanton scattering matrices introduced in that paper.

\subsubsection{K-theory and line operator categories} 
Many of the arguments we make in this paper are geometric in nature, hence work very generally and in particular should extend to critical $K$-theory (we understand critical $K$-theory following \cite{COZZ1}), for which there is also a notion of critical stable envelopes. 

An interesting new degree of freedom in the $K$-theoretic situation is the so-called slope parameter of stable envelopes. It is tempting\footnote{Based on our interpretation of certain proposals of D. Gaiotto and Y. Zenkevich, \cite{zenkevich2026wallcrossingstringnetworks}.} to speculate that corresponding slope $R$-matrices are related to wall-crossing structures arising in the context of cluster charts on $K$-theoretic Coulomb branches. 

An important technical point is that the logic we use in this paper to compare the Coulomb branch and quantum group actions depends on the fact that critical cohomology depends only on the oriented $(-1)$-shifted symplectic structure on corresponding derived critical locus. To our knowledge, this statement is only partially understood in critical $K$-theory \cite{hennion2025gluinginvariantsdonaldsonthomastype}. 

One can speculate even further about the categorification of the $K$-theoretic picture, to produce a categorical action on a matrix factorization category of $\QM^\xi(\mathbb{P}^1, X)$ defined using either the quiver or Čech model of this paper. The object which should be acting is some version of coherent sheaves on the varieties of \cite{BFNII} (see \cite{cautis2018clustertheorycoherentsatake}) which should correspond to a line operator category in a 4d $\Sh{N} = 2$ supersymmetric gauge theory \cite{niu2022localoperators4dmathcaln2}. This would be a far-reaching generalization of the 2-representations of $\mathfrak{sl}_2$ studied by Rouquier \cite{rouquierquasimap}. 

\subsubsection{Vertex functions and Hikita conjecture}
Recently, a reformulation of the so-called quantum Hikita conjecture via the enumerative theory of quasimaps has been proposed in \cite{dinkins2026quantumhikitaconjecturequasimaps}. It is an expectation that graded traces in the Coulomb branch modules constructed in the present article can be used to study twisted traces on Coulomb branch algebras, and should furnish a basis in the space of twisted traces when the torus $T$ acts on $X$ with isolated fixed points and the corresponding equivariant parameters are generic. 

A basic observation that may be useful in this direction is the following. A direct consequence of the analysis of the derived structure on the quasimap spaces summarized in Theorems \ref{main thm: derived QM quiver}, \ref{thmx: derived critical Čech}, together with well-known results relating the graded dimension of critical cohomology with integrals against corresponding virtual class \cite{Behrend}, is that the characters of the Coulomb branch modules we construct here can be identified with the Calabi-Yau specialized vertex functions employed in \cite{dinkins2026quantumhikitaconjecturequasimaps}.

It would be interesting to investigate the extent to which our results here could prove or even categorify the conjectures of \cite{dinkins2026quantumhikitaconjecturequasimaps}. A related point is that the original conjecture of \cite{bullimore2018vorticesvermas} included the statement that the quasimap cohomology should define a Coulomb branch Verma module. In the present work, we do not discuss the structure of resulting module in any generality; it would be of interest to revisit this question in the future.

\subsubsection{Relative version over $X$} 
Consider the stack $\QM^{\text{ns }\infty}(\bbP^1, X)$ on nonsingular quasimaps at $\infty$, defined as the preimage of $X\subset \FX$ under the evaluation map $\ev_{\infty}: \QM(\bbP^1, X)\to \FX$. 
This is, by construction, a $(-1)$-shifted symplectic fibration. By recent work of Kan-Kinjo-Park-Safrnov, \cite{khan2025perversepullbacks} there exists a relative notion of DT sheaf providing the natural coefficient for a cohomology theory of the shifted symplectic fibration $\QM^{\text{ns }\infty}(\bbP^1, X)\to X$. Compared to the current work, the new feature of this cohomology is that, besides the left action by Hecke modifications, it admits a right action by modifications of the quiver variety. These statements will be clarified by a forthcoming work by the authors.

\subsubsection{Related works}
The DT cohomology of $\QM^\xi(\mathbb{P}^1, X)$ can be viewed as a modern generalization of very classical constructions in geometric representation theory. In the very special case $X = T^*(GL_n/P)$ when $X$ is a cotangent bundle to a partial flag variety in type $A$, Theorems \ref{thmx: Cohas and shifted Yangians action}, \ref{thmx: Coulomb aciton intro}, and Theorem \ref{propx} (together with known results on Coulomb branches) recover the results in \cite{Braverman_2011}, \cite{feigin2011yangianscohomologyringslaumon}, \cite{nakajima2011handsawquivervarietiesfinite}, with a completely different proof. The article \cite{Hilburn_2023} constructs Coulomb branch modules which we expect to agree with ours when the evaluation point $\xi$ is polarized, so that the DT cohomology can be dimensionally reduced to Borel-Moore homology. The PhD thesis of H. Liu \cite{Liuphd} (based, in turn, on unpublished work of Nakajima and Okounkov) describes quantum group modules which ``approximate'', in a certain sense, the DT cohomologies studied in this paper.

We have also been informed by E. Vasserot of a work in progress with W. Zhao \cite{VasserotCoulomb} aimed at constructing Coulomb branches on global curves for general affine symplectic targets by means of a vanishing cycle formalism. The authors also consider dimensional reduction in this setting and relate their construction to quasimaps. It would be very interesting to compare our work with theirs. 

While drafting this paper, two independent works by Jindal-Negut\cite{jindal2026loopnilpotentcohomologicalhallalgebra} and Muthiah-Weekes \cite{MWzastavas} relating shifted Yangians to Coulomb branches have appeared. The former uses shuffle algebra methods to present a surjective map from the appropriate double of the $\omega$-nilpotent Hall algebra $\HCoha^{T, \onil}_{\tilde Q, \tilde W}$ to the Coulomb branch. We believe this to be the combinatorial counterpart of our Theorem \ref{thmx: map hall Coulomb}. On the other hand, the work by Muthiah-Weekes gives an injective morphism from the Coulomb branch to certain asymptotic limit of a shifted Yangian. It would be interesting to relate it to our construction, and understand its geometric nature.

\subsection{Structure of the paper}

The paper is divided in four parts.

\begin{enumerate}
    \item In Part\ref{part: quiver model}, we consider the quiver model for quasimaps critical cohomology and construct the relevant Hall algebra and shifted Yangian action. In particular, we state and prove Theorem \ref{thmx: Cohas and shifted Yangians action}.

    \item In Part \ref{part: Čech model}, we develop an alternative description of the BFN Coulomb branch via vanishing cycles, introduce the Čech model for quasimap critical cohomology, and combine these results to construct the action on the Coulomb branch on \eqref{eq: DT cohomology intro}. Theorem \ref{thmx: critical raviolo} and Theorem \ref{thmx: Coulomb aciton intro} are stated and proved.
    
    \item In Part \ref{part: compatibilities}, we state Theorem \ref{thmx: map hall Coulomb} and discuss Conjecture \ref{conjx: diagram}.

    \item In part \ref{part: quasimap shifted symplectic}, we compare the shifted symplectic structures of the two critical models introduced in Part \ref{part: quiver model} and Part \ref{part: Čech model} and prove their equivalence. In particular, Theorem \ref{thmx: derived critical Čech} and \ref{main thm: derived QM quiver} are stated and proved.
\end{enumerate}

We also have three appendices. In Appendix \ref{app: Sheaf operations and critical cohomology} we review known constructions in sheaf theory and critical cohomology. In Appendix \ref{subsec: Proof of long Theorem} we prove Theorem \ref{thmx: map hall Coulomb}. In Appendix \ref{app: derived quasimaps} we present some foundational material on the derived structure of the moduli space $\QM^\xi(\bbP^1, X)$.

\subsection{Conventions}
\label{subsec: Conventions}
\begin{itemize}
\item 
All algebraic stacks are defined over the complex numbers $\BoC$. Accordingly, $\dim(X)$, $\codim_Y(X)$ etc. are the complex dimensions.
For a fixed torus $T$ let $\FM$ be a stack presented as a global quotient $Y/(G\times T)$ for some algebraic group $G$, and some variety $Y$.  Then we define the shifted perverse t-structure on $\Dub_c(X/T)$ and $\Dub_c(\FM)$ by setting ${}^{\fp'}\!\vtau^{\leq i}\coloneqq {}^{\fp}\vtau^{\leq i-\dim(T)}$ and ${}^{\fp'}\!\vtau^{\geq i}\coloneqq {}^{\fp}\vtau^{\geq i-\dim(T)}$.  Objects in the heart of this t-structure are precisely those complexes $\CF$ such that $(Y/G\rightarrow Y/(G\times T))^*\CF$ is perverse.  We then take the shift of the Verdier duality functor $\DD$ preserving this heart.

\item 
For vector bundle $E$ with Chern roots $x_i$ we write $c_u(E) := \prod_i (u - x_i)$ for its Chern polynomial. For a $K$-theory class $E-F$, $c_u(E - F) = c_u(E)/c_u(F)$. 

\item The cohomology of the clssifying stack $\B G$ of an algebraic group $G$ is denoted by $\HO_G$. We always treat the $G$-equivariant cohomology of any space as a $\HO_G$-module. In particular, tensor product are taken in this category.

\item
In this article we will consider various versions of the stack of Hecke modifications of $G$-bundles and related constructions. We use the following convention. We reserve the notation $\Hecke$ for the stack of ``positive'' Hecke modifications of $\GL(n)$ bundles, that is for Hecke modifications that extend to an injective morphism of sheaves. We will use the notation $\Hecke^{\pm}$ for the stack of ``all'' Hecke modifications. i.e. with no condition about injectivity.

\item We choose the right coset presentation for the affine Grassmannian, i.e. we exhibit it as $\Gr=G\fps\backslash G\lfps$. This is different from the convention in \cite{BFNII} as well as from vast part of the lliterature, but has the advantage of making the associated Coulomb branch acting on quasimap cohomology from the left, in agreement with the Yangian action.

\item Our notation $\Coul_{G,N}$ for the Coulomb branch is also nonstandard; the notation introduced in \cite{BFNII} and commonly used in the literature is $\mathcal{A}_{G,N}$. We adopt this convention to distinguish the Coulomb branch from the Hall algebra, which we denote by $\HCoha_{\tilde Q,\tilde W}$, in accordance with the existing literature. We will also make use of the notation $\Coul^{\vir}_{G,N}$ to denote a version of $\Coul_{G,N}$ with cohomological degree shifted by some finite amount. The superscript only refers to this homological shift so, as $\pi_1$-graded algebras, we have a natural isomorphism $\Coul_{G,N}\cong \Coul^{\vir}_{G,N}$

\item In the more computationally intensive parts of this paper, such as Part \ref{part: quasimap shifted symplectic} and Appendix \ref{app: derived quasimaps}, we will suppress all indices related to nodes and edges of the quiver $Q$ in an attempt to reduce clutter. The concerned reader is invited to reinstate the indices.
\end{itemize}

\subsection{Acknowledgments}

We would like to warmly thank Ben Davison and Andrei Okounkov for years of inspiration and encouragement, as well as for their fundamental role in shaping this century's geometric representation theory. We also thank Andres Ibanez Nunez for sharing his knowledge in derived algebraic geometry, which helped us deal with some of the technical parts of this work. Special thanks are also owed to Peng Zhou for giving us the courage and inspiration to use more than one critical chart. We also thank Mina Aganagic, Viktor Alekseev, Alexander Braverman, Yalong Cao, Sabin Cautis, Igor Chaban, Hunter Dinkins, Daniel Halpern-Leistner, Andres Fernandez Herrero, Davide Gaiotto, Shivang Jindal, Joel Kamnitzer, Alexis Leroux-Lapierre, Hiraku Nakajima, Nikita Nekrasov, Tasuki Kinjo, Vasily Krylov, Tony Pantev, Hyeonjun Park, Pavel Safronov, Andrey Smirnov, Eric Vasserot, Alex Weekes, Yaping Yang, and Yehao Zhou for many useful conversations.

\newpage

\paperpart{Quiver model and Hall action}
\label{part: quiver model}

\section{Setup}

\subsection{Quivers and their representations in a category}
\label{subsec: Quivers and their representations in a category}
A quiver $Q=(Q_1,Q_0,s,t)$ consists of a set of arrows $Q_1$, a set of vertices $Q_0$, and two maps $s,t: Q_1\rightarrow Q_0$ taking an arrow to its source and target, respectively. By abuse of notation, we will often denote a quiver by $Q=(Q_1,Q_0)$, ignoring the maps $s$ and $t$. Its adjacency matrix  $\Q\in\End(\BoZ^{Q_0})$ is the matrix with entries
\[
\Q_{ij}=\#\{\text{edges $e\in Q_1$ $|$ $s(e)=i,\, t(e)=j$}\}.
\]
The Cartan matrix $\Cartan$ of the quiver is then defined as
\[
\Cartan=2\id-\Q-\Q^T,
\]
where $\Q^T$ is the transpose of $\Q$.

We denote by $\BoC Q$ the free path algebra of $Q$ over $\BoQ$.  This algebra has a $\BoQ$-basis given by paths in $Q$, including length zero paths $\lazy_i$ for each $i\in Q_0$, which start and end at $i$.
Here, $pq$ denotes the concatenation of $p$ and $q$, and we read paths from right to left (the same way we compose morphisms in a category). Composition is defined by concatenation. It is clear that $\lazy_i$ are pairwise orthogonal idempotents and that $\sum_{i\in Q_0} \lazy_i=1$. More generally, we will fix a two sided ideal $I\subset \BoQ R$ such that $\lazy_i\neq I$ for all $i\in Q_0$ and consider the quotient $\BoC Q/I$. The data $(Q, I, \{\lazy_i\}_{i\in Q_0})$ determine a category $\calC(\BoC Q/I)$, which is $\BoQ$-linear. Its objects are the the idempotents $\lazy_i$ and the morphisms are $\Hom(\lazy_i, \lazy_j)= \lazy_j (\BoC Q/I) \lazy_i$. It is clear that $\calC(\BoQ Q/I)$ is equivalent to the quotient of $\calC(\BoC Q)$ by the congruence relations $\lazy_i, \lazy_j$ for all $i,j\in \BoQ_0$.

Given another linear category $\calD$ over a field extension $K$ of $\BoQ$, a representation of $\BoC Q/I$ in $\calD$ is a linear functor $\calC(\BoC Q)\to \calD$. We will be only interested in the case where $\calD$ is the category of finite dimensional vector spaces over $\BoC$, or, more generally, coherent sheaves on a curve. Explicitly, a representation  of $\BoC Q$ in $\calD$ consists of the data of a collections of objects $\Sh{F}_i\in \Ob(\calD)$ indexed over $Q_0$ and a collection of morphisms $\alpha_a:\Sh{F}_{s(a)}\to \Sh{F}_{t(a)}$ for all $a\in Q_1$ satisfying the relations of $I$.  These functors form a $K$-linear category $\Rep_{\BoC Q/I}(\calD)$, the category of representations of $Q$ with values in $\calD$. Its morphisms are natural transformation between two given functors. More explicitly, a morphism in $\Rep_{\BoC Q/I}(\calD)$
\[
(\{\Sh{F}_i\}_{i\in Q_0}, \{\alpha_a\}_{a\in Q_1})\to (\{\Sh{G}_i\}_{i\in Q_0}, \{g_a\}_{a\in Q_1})
\]
is given collection of morphisms $h_i: \Sh{F}_i\to \Sh{G}_i$ in $\calD$ such that $h_{t(a)}f_{a}=g_{t(a)}h_{s(a)}$ for all $a\in Q_1$. 
To remove clutter, we will write $\Rep_{\BoC Q/I}$ instead of $\Rep_{\BoC Q/I}(\Vect_{\BoC})$, and $\Rep_Q(\calD)$ as a further shortcut for $\Rep_{\BoC Q}(\calD)$.

We remark that in some cases it is possible to identify $\Rep_{\BoQ}(\calD)$ with $\Rep_{\BoC Q_{\calD}/I_{\calD}}$, where $ Q_{\calD}$ and $I_{\calD}$ depend on $Q$ and $\calD$. We will discuss some basic examples in Proposition \ref{prop: equivalence of categories of coh and reps}.


\subsection{Doubling, looping, tripling, framing, and potentials}
\label{subsec: Doubling, looping, tripling, and framing}

Fix a quiver $Q$. The set of dimension vectors is the monoid $\BoN^{Q_0}\subset \BoZ^{Q_0}$. Given a dimension vector $\dd = (\dd_i)_{i \in Q_0}  \in\BoN^{Q_0}$ the space of $\dd$-dimensional representations is defined as
\[
\Rep_{Q}(\dd)\coloneqq \prod_{a\in Q_1}\Hom(\BoC^{\dd_{s(a)}},\BoC^{\dd_{t(a)}}).
\]
A potential for $Q$ is an element $W\in (\BoC Q)_{\cyc}\coloneqq \BoC Q/[\BoC Q,\BoC Q]_{\mathrm{vect}}$ where the subscript means that we take the quotient by the vector space spanned by commutators. For any $\dd\in \BoN^{Q_0}$, we get a $\GL_{\dd}$-invariant function 
\[
\Tr(W): \Rep_Q(\dd)\to \BoC.
\]

Given $Q$, the doubled quiver $\overline{Q}$ is defined by setting $\overline{Q}_0=Q_0$ and setting $\overline{Q}_1=Q_1\coprod Q^{\opp}_1$, where $Q^{\opp}_1=\{a^*\}$ is a copy of $Q_1=\{a\}$ as a set but its arrows $a^*$ have opposite orientation, i.e. $s(a)=t(a^*)$ and $t(a)=s(a^*)$. The trace pairing then gives an isomorphism $\Rep_{\overline Q}(\dd)=T^*\Rep_Q(\dd)$. With respect to the inherited symplectic structure, the  $\GL_{\dd}$ action is Hamiltonian, so it has a moment map 
\begin{equation}
    \label{eq: moment map}
    \mu_{\dd}: \Rep_{\overline Q}(\dd)\to \gl_{\dd}=\prod_{i\in Q_0}\End(\BoC^{\dd_i}) \qquad (a, a^*)\mapsto \sum_{a\in Q_1}[a, a^*].
\end{equation}
We define the preprojective algebra
\begin{equation}
\label{PA_def}
\Pi_Q\coloneqq\BoC\overline{Q}/\langle \sum_{a\in Q_1}[a,a^*]\rangle.
\end{equation}
The space of $\dd$-dimensional representations of $\Pi_Q$ form a subvariety of $\Rep_Q(\dd)$, which is naturally identified with $\mu^{-1}_{\dd}(0)\subset\Rep_Q(\dd)$.

Given a quiver $Q$, we also introduce its \emph{looping} $\hat Q$ and its \emph{tripling} $\tilde Q$. The former is the quiver defined by $\hat Q_0=Q_0$ and $\hat Q_1=Q_1\sqcup \{\omega_i \;|\; i \in Q_0\}$; the latter is instead the quiver
\[
\tilde{Q}_0=Q_0;\quad\quad \tilde{Q}_1=\overline{Q}_1\coprod \{\omega_i\;\lvert \;i\in Q_0\}
\]
with $s(\omega_i)=t(\omega_i)=i$. In particular, the tripled quiver $\tilde Q$ is the looping of the doubled quiver $\bar Q$. The tripled quiver $\tilde Q$ carries the canonical cubic potential $W\in (\BoC \tilde Q)_{\cyc}$ given by
\begin{equation}
\label{CCP}
\tilde{W}\coloneqq \left(\sum_{a\in Q_1} [a,a^*]\right)\left(\sum_{i\in Q_0} \omega_i\right).
\end{equation}
It plays a key role in the realization of quantum groups as Hall algebras.

To both $\hat Q$ and $\tilde Q$, we associate distinguished quotients of their path algebra. We define $\Xi_Q$ be the quotient of the path algebra $\BoC \hat Q$ by the two-sided ideal generated by 
\begin{equation}
    \label{eq: hat quver relation}
    \sum_{a\in Q_1, i\in Q_0} a\omega_i -\omega_i a
\end{equation}
Similarly, we define the \emph{Jacobi algebra} $\Jac(\tilde Q, \tilde W)$. as the the quotient of the path algebra $\BoC \tilde Q$ by the two-sided ideal generated by 
\[
\sum_{a\in Q_1, i\in Q_0} a\omega_i -\omega_i a \qquad \sum_{a\in Q_1, i\in Q_0} a^*\omega_i -\omega_i a^*\qquad \sum_{a\in Q_1}[a, a^*].
\]
The notation refers to the fact that $\Jac(\tilde Q, \tilde W)$ is obtained as the quotient of $\BoC Q$ by the non-commutative derivatives of the cubic potential \eqref{CCP}.

As mentioned in S\ref{subsec: Quivers and their representations in a category}, we will also be interested in representations of $Q$ and $\Pi_Q$ in categories other than $\Vect_{\BoC}$, and particularly in the case when is $\Coh_0(\AA^1)$, namely be the category of rank zero coherent sheaves on $\AA^1$. The following basic proposition plays a key role in the geometric understanding of the categories of representations of $\Xi_Q$ and $\Jac(\tilde Q, \tilde W)$.
\begin{proposition}
\label{prop: equivalence of categories of coh and reps}
    We have natural equivalences of categories
    \[
    \Rep_{Q}(\Coh_0(\AA^1))\simeq \Rep_{\Xi_Q}\qquad \Rep_{\Pi_Q}(\Coh_0(\AA^1))\simeq \Jac(\tilde Q, \tilde W).
    \]
\end{proposition}
We will discuss this equivalence explicitly in \S\ref{subsec: Q decorated coherent sheaves on A1} and \S\ref{subsec: Hecke correspondence on quasimap space}

We end this section with a discussion about framings. Fix a quiver $Q$ and a framing vector $\ww\in \BoQ^{Q_0}$. The framed quiver $Q_{\ww}=(Q_{\ww,0},Q_{\ww,1})$ is defined by setting 
\[
Q_{\ww,0}=Q_0\sqcup \{\infty\}\qquad Q_{\ww,1}=Q_1\sqcup \{ r_{i,m} \;|\; i\in Q_0\;,1\leq m\leq \ww_i\} 
\]
where $seri(r_{i,m})=i$ and $t(r_{i,m})=\infty$ for all $1\leq m\leq \ww_i$. Notice that $\Rep_{Q_{\ww}}(\vv,1)=\Rep_Q(\vv)\oplus \Hom(\BoC^{\ww},\BoC^{\vv})$, where
\[
\Hom(\BoC^{\vv},\BoC^{\ww})=\bigoplus_{i\in Q_0}\Hom(\BoC^{\vv_i}, \BoC^{\ww_i}).
\]
We write $\bar Q_{\ww}$, $\tilde{Q}_{\ww}$ and $\hat Q_{\ww}$ to denote the quivers obtained by framing and \textit{then} doubling, tripling,  or looping, as opposed to the more logical, but also more cumbersome $\ol{Q_{\ff}}$, $\widetilde{Q_{\ff}}$, and $\widehat {Q_{\ww}}$. Additionally, to save on space, we will sometimes write $\AA_{\dd}(Q)$ instead of the slightly more cumbersome $\Rep_{Q}(\dd)$.

\subsection{Moduli spaces of quiver representations}
\label{subsec: moduli of quiver reps}
Fix an arbitrary quiver $Q=(Q_0,Q_1)$ and a dimension vector $\dd\in \BoQ^{Q_0}$. We denote by $\FM_{\dd}(Q)$ the quotient of the stack of $\dd$-dimensional $\BoC Q$-modules. The stack of finite dimensional $\BoC Q$-modules is then given by 
\[
\bigsqcup_{\dd\in \BoQ^{Q_0}}\FM_{\dd}(Q)
\]
We have an equivalence
\[
\FM_{\dd}(Q)\simeq \Rep_Q(\dd)/\GL_{\dd}.
\]
More generally, we will consider the moduli stack of representations of the quotient $A=\BoC Q/I$ associated to the choice of a two sided ideal $I\subset \BoC Q$, cf. \S\ref{subsec: Quivers and their representations in a category}. We denote by $\FM(A)$ the corresponding moduli stack. It is a closed substack of $\FM(Q)$ and each connected component admits a presentation $\FM_{\dd}(A)\simeq \Rep_{A}(\dd)/\GL_{\dd}$. If besides $\GL_{\dd}$ we have an additional torus $T$ acting on $\Rep_{A}(\dd)/\GL_{\dd}$ then we consider the quotient $\FM^T_{\dd}(A)$ of $\FM_{\dd}(A)$ by $T$.

Consider now the coarse moduli space $\CM_{\dd}(Q)=\mathrm{Spec}(\Gamma(\mathcal{O}_{\Rep{A}(\dd)})^{\Gl_{\dd}})$.  Closed points of $\CM_{\dd}(A)$ are in bijection with isomorphism classes of $\dd$-dimensional semisimple $\BoC Q$-modules.  We denote by $\JH_{\dd}\colon \FM_{\dd}(A)\rightarrow \CM_{\dd}(A)$ the affinization morphism.  From the Hilbert-Mumford criterion it follows that a morphism sends a module to its semisimplification, i.e. the direct sum of the subquotients appearing in some (equivalently, any) Jordan--H\"older filtration.  The torus $T$ acts on $\CM_{\dd}(A)$, and we denote by $\CM^T_{\dd}(A)=\Msp_{\dd}(A)/T$ the stack-theoretic quotient. We continue to denote by the same symbols $\JH_{\dd}$ the induced morphism $\JH_{\dd}\colon \FM^T_{\dd}(A)\rightarrow \CM_{\dd}^T(A)$.

We will also consider versions of the above stacks incorporating stability conditions.  For a given (slope) stability condition $\zeta\in\BoZ^{Q_0}$, we denote by $\Rep_{A}(\dd)^{\zeta\sst}\subset \Rep_{A}(\dd)$ the open subvariety whose points correspond to $\zeta$-semistable $\dd$-dimensional $\BoC Q$-modules.  We set $\FM^{T,\zeta\sst}_{\dd}(A)\coloneqq\Rep_{A}(\dd)^{\zeta\sst}/(\Gl_{\dd}\times T)$.  The natural morphism $\FM^{T,\zeta\sst}_{\dd}(A)\rightarrow \FM^{T}_{\dd}(A)$ is an open embedding.  We construct the coarse moduli space $\CM^{\zeta\sst}_{\dd}(A)$ of $\zeta$-semistable $\dd$-dimensional $A$-modules via GIT as in \cite{King}. Namely, given an integral stability condition, we define 
\[
\CM^{\zeta \sst}_{\dd}(Q)\coloneqq\text{Proj}\left(\bigoplus_{n\geq 0} \Gamma\left(\Sh{O}_{\Rep_{\dd}(A)}\otimes \theta_{\zeta}^n\right)^{\GL_\dd}\right)
\]
where $\theta_{\zeta}\in \Char{\GL_{\dd}}$ is the cocharacter given by 
\[
(g_i)_{i\in I}\mapsto \det(g_i)^{\zeta_i}.
\]
Closed points of this quasi-projective variety are in bijection with isomorphism classes of $\zeta$-polystable $\dd$-dimensional $\BoC Q$-modules, and the natural morphism $\pi\colon\CM^{\zeta\sst}_{\dd}(A)\rightarrow \CM_{\dd}(A)$ is a GIT quotient map, and hence projective (see \cite{King}).  By construction, the torus $T$ acts on $\CM^{\zeta\sst}_{\dd}(Q)$, and the morphism $\pi$ is $T$-equivariant.  We define the stack $\CM^{T,\zeta\sst}_{\dd}(A)\coloneqq \CM^{\zeta\sst}_{\dd}(A)/T$.  We denote by 
\[
\JH^{\zeta}_{\dd}\colon \FM^{T,\zeta\sst}_{\dd}(Q)\rightarrow \CM^{T,\zeta\sst}_{\dd}(Q)
\]
the natural morphism. 

\subsection{Nakajima quiver varieties}
\label{subsec: Nakajima quiver varieties}
For any given stability condition $\zeta\in \BoQ^{Q_{\ww,0}}$ for the quiver $Q_{\ww}$, the Nakajima quiver variety $\Nak^{\theta}_Q(\ww,\vv)$ is defined as the GIT quotient 
\begin{align*}
    \Nak^{\zeta}_Q(\ww,\vv)
    &=\CM^{(\zeta) \sst}_{(\vv,1)}(\Pi_{Q_{\ww}}).
\end{align*}
The following is well known:
\begin{lemma}
\label{lemma characterization anticanonical stability}
    The slope $\zeta$ is generic. Therefore, stability is equal to semi-stability and all points have trivial stabilizer. Furthermore, a $\Pi_{Q_{\ww}}$-module iff there is no proper nontrivial submodule supported on $Q_0\subset Q_0\cup\{\infty\}$. Equivalently, for any submodule of dimension $(\vv,0)$, we have $\vv=0$.
\end{lemma}
As a corollary, it follows that that for a generic stability $\zeta$, the quiver variety $\Nak^{\zeta}_Q(\ww,\vv)$ is isomorphic to the stack quotient $\mu_{\ww,\vv}^{-1}(0)^{\zeta\sst} /\GL_{\vv}$.
In particular, we deduce that the canonical map 
\[
\FM^{\zeta\sst}_{(\vv,1)}(\Pi_{Q_{\ww}})\rightarrow \CM^{\zeta\sst}_{(\vv,1)}(\Pi_{Q_{\ww}})=\Nak^{\zeta}_Q(\ww,\vv)
\]
is a trivial $\BoC^*$-gerbe, where the $\BoC^*$ action comes from rescaling the framing vertex $\infty$. In other words, we have 
\begin{equation}
    \label{eq: quiver variety as stack}
    \Nak^{\zeta}_Q(\ww,\vv)/{\BoC^*}=\FM^{\zeta\sst}_{(\vv,1)}(\Pi_{Q_{\ww}}).
\end{equation}
Throughout this paper we set the stability condition for quiver varieties to be $\zeta^-=(0,-1)$. It is a generic stability condition. We will usually denote a tuple of in the prequotient $\mu^{-1}_{\ww}(0)$ by $(\xi_2, \xi_3, \iota, \jmath)$, where $(\xi_2,\xi_3)$ corresponds to a point in $T^*\Rep_Q(\vv)=\Rep_Q(\vv)\oplus (\Rep_Q(\vv))^\vee$ and $(\iota, \jmath)$ to the framing arrows in $\Hom(\BoC^{\ww}, \BoC^{\vv})\oplus (\Hom(\BoC^{\ww}, \BoC^{\vv}))^\vee$. The stability $\zeta^-$ corresponds to the condition that the vector spaces $\BoC^{\vv_i}$ are \emph{co-generated} by $\jmath$ and the action of $\xi_2, \xi_3$.

In the rest of this paper, we will fix a Nakajima variety $X$, understood with stability $\zeta^-$ and some fixed choice of $\vv, \ww$. 

\subsection{Torus action on quiver moduli} 
\label{subsec: torus action quiver varieties}
In this section we gather notations and conventions for the action of tori on various moduli spaces of quiver representations that will be used throughout this paper. Start by considering a Nakajima quiver variety $X = \Nak^\zeta_Q(\ww, \vv)$.  By construction, $X$ arises as a GIT symplectic reduction of $T^* \text{Rep}_Q(\vv, \ww)$. On this cotangent bundle acts a torus $\mathbb{C}^\times_\hbar$, by dilation of the cotangent fibers. $\mathbb{C}^\times_\hbar$ continues to act on $X$ and scales its canonical symplectic form $\omega_X$ with weight $\hbar$. 

On $X$ also acts the group 
\[
\GL_{\ww} \times \SPL(\text{self loops}) \times \GL(\text{edges}) \coloneqq \prod_{i \in Q_0} (\GL_{\ww_i} \times \SPL(\Q_{ii})) \times \prod_{i \neq j \in Q_0} \GL(\Q_{ij}).
\]
where $\Q_{ij}$ is the adjacency matrix, cf. \ref{subsec: Quivers and their representations in a category}.

We set $T \subset \GL_{\ww} \times \SPL(\text{self loops}) \times \GL(\text{edges}) \times \mathbb{C}^\times_\hbar$ to be a maximal torus. Later in this paper, when discussing Hall algebras and quasimaps, we will consider various moduli spaces of quiver representations arising from some framing of the tripled quiver $\tilde{Q}$. On any such moduli space, there is a distinguished torus $\mathbb{C}^\times_{t_1}$ which acts by scaling the self-loops $\omega_i$ (referred to as $B_{1, i}$ in the context of the quasimap quiver, see \S\ref{sec: quiver critical description quasimaps}). On such moduli spaces we always understand the $T$-action via the embedding
\[
\begin{tikzcd}
    T \arrow[r, "\sim"] & \ker(t_1 \cdot \hbar) \arrow[r, hook] & \mathbb{C}^\times_{t_1} \times T. 
\end{tikzcd}
\]
Recall that for each edge $a \in Q_1$ we denote by $a^*$ the opposite edge in the doubled quiver. For ease of writing formulas below, we make the following identifications in dealing with $T$-equivariant cohomology and $K$-theory:
\[
\begin{split}
    K_T(\text{pt}) & \simeq  \mathbb{C}[t_1^{\pm}, t_a^{\pm}, t_{a^*}^\pm, a_{i_k}^\pm,]/(t_1 t_a t_{a^*} - 1)_{a \in Q_1} \\
    \HO_T(\text{pt}) & \simeq \mathbb{C}[\epsilon_1, \epsilon_a, \epsilon_{a^*}, a_{i_k}]/(\epsilon_1 + \epsilon_a + \epsilon_{a^*})_{a \in Q_1}
\end{split}
\]
where indices $i_k$ have $i \in Q_0$, $k = 1, \dots, w_i$. We will sometimes refer to $T$ as the \textit{Calabi-Yau torus} and the condition enforced by the quotients above as the \textit{Calabi-Yau specialization}, borrowing terminology from the equivariant DT theory of threefolds.

\subsection{Vanishing cycles and DT perverse sheaves}
\label{subsec: Vanishing cycles and DT perverse sheaves}

For an oriented $(-1)$-shifted symplectic stack $\FX$, Ben-Bassat-Brav-Bussi-Joyce \cite{Br12} constructed a perverse sheaf $\varphi_{\FX,o}$ called the Donaldson–Thomas sheaf. The subscript $o$ refers to the choice of orientation. Whenever that $\FX$ is a global critical locus of a function $f:\FY\to \BoC$ on a smooth algebraic stack $\FY$, the DT sheaf can be defined as the vanishing cycle functor $\phip{f}$ applied to the shifted constant sheaf $\BoQ^{\vir}_{\FY}=\BoQ_{\FY}[\dim{\FY}]$. In other words, we have $\phip{f}\BoQ^{\vir}_{\FY}=\varphi_{\FX,o_{\text{can}}}$, where $o_{\text{can}}$ is the canonical orientation for the critical locus $\FX=\Crit(f)$. The theory of vanishing cycles is reviewed in Appendix \ref{app: Appendix VC}. 

More generally, there exists a relative version of the DT sheaf for $(-1)$-shifted symplectic fibration $\pi: \FX\to B$ \cite{khan2025perversepullbacks}. In this paper, we will only be concerned with the simplest case of this construction, corresponding to a fibration of the form $\pi: \FX/T\to \B T$, where $\FX$ is $(-1)$ shifted and the $T$-action preserves all the structure. Wist slight abuse of notation, will denote the relative DT sheaf\footnote{Following \cite{khan2025perversepullbacks} this is perverse pullback functor $\pi_{\varphi}$ applied to the constant sheaf $\BoQ_{\B T}$. This sheaf is perverse onnce pulled back along the map $\FX\to \FX/T$, in agreement with our conventions \S\ref{subsec: Conventions}.} by $\varphi_{\FX/T,o}$. If $\FX=\Crit(f)\subset \FY$ and the $T$-action extends on $\FY$ and preserves $f:\FY\to \BoC$, then $\varphi_{\FX/T,o_{\text{can}}}=\phip{f}\BoQ^{\vir}_{\FY/T}$. This is the only case we will consider in this paper. Vanishing cycles are reviewed in Appendix \ref{app: Appendix VC}.

The subscript referring to the orientation will be dropped whenever the canonical orientation associated to a critical locus is understood. 

Despite being only concerned with global critical loci, in this article we will use two derived equivalent but different critical loci model to describe the same sheaf $\varphi_{\FX/T}$ for $\FX=\QM^{\xi}(\bbP^1, X)/T$. To avoid confusion, we will use $\varphi_{\FX/T,o_{\text{can}}}$ to denote the intrinsic DT sheaf, and expressions like $\phip{f}\BoQ^{\vir}_{\FY/T}$ to denote the vanishing cycle construction that gives rise to the DT sheaf in a given critical model $(\FY, f)$.
With slight abuse of notation, we will write $\HO_T(\FX, \varphi_{\FX,o_{\text{can}}})$ to denote the cohomology
\[
\HO(\FX/T, \varphi_{\FX/T,o_{\text{can}}})\cong \HO(\FY/T, \phip{f}\BoQ^{\vir}_{\FY/T}).
\]
We will also informally refer to all these $\HO_T$-modules as critical cohomology or vanishing cycle cohomology of $\FX$.

\section{Cohomological Hall algebras}

\subsection{The three dimensional CoHA}
\label{subsec: The three dimensional CoHA}
Hall algebras can be constructed for arbitrary quiver but, since we are only interested in the tripled quiver $\tilde Q$ with canonical cubic potential $\tilde W$, we restrict to this case. Recall the notation from \S\ref{subsec: moduli of quiver reps}. We will work torus-equivariantly in the Calabi-Yau specialization. In particular, we will assume that the edges $a$ and $a^*$ are scaled by dual weights and the loop $\omega_i\in \tilde Q_1$ is scaled by a weight such that that the potential $\tilde W$ is preserved by $T$. For more details, see \S\ref{subsec: torus action quiver varieties}. 

Consider the constructible complex on $\CM^T_{\dd}(\tilde{Q})$ given by
\begin{equation}
    \label{eq: 3d bps coha graded component}
    \Coha^T_{\tilde{Q}, \tilde{W},\dd}=(\JH_{\dd})_*\:\phip{\Tr(W)}\BoQ_{\FM^{T}_{\dd}(\tilde Q)}^{\vir}
\end{equation}
As a sheaf over $\CM^T(\tilde{Q})$, the 3d-relative Hall algebra is defined by
\[
\Coha^T_{\tilde{Q}, \tilde{W}}\coloneqq\bigoplus_{\dd\in \BoN^{Q_0}}\Coha^T_{\tilde{Q}, \tilde{W},\dd}
\]
The Hall product equips the latter with the structure of an associative algebra with respect to the monoidal structure on $D_c(\CM^T(\tilde Q))$ defined on connected components by 
\begin{equation}
    \label{eq: monoidal structure good moduli space}
    \Sh{F_{\dd'}}\otimes_{\oplus}\Sh{F_{\dd''}}=\oplus_*\left(\Sh{F_{\dd'}}\boxtimes\Sh{F_{\dd''}}\right)
\end{equation}
where $\oplus$ is the map
\[
\oplus: \CM^T_{\dd'}(\tilde Q)\times_{\B T} \CM^T_{\dd''}(\tilde Q)\to \CM^T_{\dd'+\dd''}(\tilde Q).
\]
The multiplication 
\begin{equation}
\label{eq: coha multiplication}
    \Coha^T_{\tilde{Q}, \tilde{W},\dd'}\otimes_{\oplus} \Coha^T_{\tilde{Q}, \tilde{W},\dd''}(\tilde Q)\to \Coha^T_{\tilde{Q}, \tilde{W},\dd'+\dd''}
\end{equation}
is then given by a standard pull-push along the top row of the commutative diagram 
\[
\begin{tikzcd}
    \FM^T_{\dd'}(\tilde Q)\times_{\B T} \FM^T_{\dd''}(\tilde Q)\arrow[d, "\JH_{\dd'}\times_{\B T}\JH_{\dd''}"]& \FM^T_{\dd',\dd''}(\tilde Q)\arrow[l, swap,  "q"]\arrow[r, "p"] & \FM^T_{\dd'+\dd''}(\tilde Q)\arrow[d, "\JH_{\dd}"]
    \\
    \CM^T_{\dd'}(\tilde Q)\times \CM^T_{\dd''}(\tilde Q)\arrow[rr, "\oplus"] & &  \CM^T_{\dd'+\dd''}(\tilde Q)
\end{tikzcd}
\]
Here, $\FM^T_{\dd'+\dd''}(\tilde Q)$ is the moduli stack parametrizing short exact sequences of $\BoC \tilde Q$ modules $0\to M'\to M\to M''\to 0$ and $p$ and $q$ are the canonical maps forgetting various terms of the sequence.
We refer to \cite{da13} for details on the construction. Taking derived global sections we obtain the absolute CoHA for the pair $(\tilde{Q},\tilde{W})$\footnote{This construction works more  generally for an arbitrary quiver, but we will not consider such variants in this article.}. Namely, we define the $\BoN^{Q_0}$-graded, cohomologically graded $\HO_T$-module $\HCoha_{\tilde{Q}, \tilde{W}}^T$ by setting 
\begin{align}
\label{eq: absolute 3d coha}
\HCoha^T_{\tilde{Q}, \tilde{W}}\coloneqq\bigoplus_{\dd\in \BoN^{Q_0}}\HCoha^T_{\tilde{Q}, \tilde{W},\dd}\qquad 
\HCoha_{\tilde{Q}, \tilde{W},\dd}^{T}=\HO(\FM^T_{\dd}(\tilde{Q}),\phip{\Tr(W)}\BoQ_{\FM^{T}_{\dd}(Q)}^{\vir})
\end{align}
and equip it with a algebra structure taking derived global sections of \eqref{eq: coha multiplication}. 
\subsection{BPS sheaves and integrality}
\label{subsec: BPS sheaves and integrality}
By \cite{preproj}, the sheaf $(\JH_{\dd})_*\phip{\Tr(\tilde W)}\BoQ^{\vir}_{\FM^T_{\dd}(\tilde Q)}$ is pure, so it decomposes as 
\begin{equation}
    \label{eq: purity coha}
    (\JH_{\dd})_*\phip{\Tr(\tilde W)}\BoQ^{\vir}_{\FM^T_{\dd}(\tilde Q)}=\bigoplus_{i\in  \BoZ} {}^{\fp'}\!\CH^i\!\left((\JH_{\dd})_*\phip{\Tr(\tilde{W})}\BoQ^{\vir}_{\FM^T_{\dd}(\tilde{Q})} \right)[-i].
\end{equation}
Moreover, the perverse cohomologies vanish for $i\leq 0$. Therefore, we define the \textit{BPS sheaf} by
\[
\BPSh^{T}_{\tilde Q,\tilde W,\dd}\coloneqq {}^{\fp'}\!\CH^1\!\left(\JH_*\:\phip{\Tr( \tilde W)}\BoQ_{\FM^{T}_{\dd}(\tilde  Q)}^{\vir}\right).
\]
By definition is a perverse sheaf on $\CM_{\dd}(\tilde{Q})$.
Define 
\[
\BPSh^{T}_{\tilde Q,\tilde W}=\bigoplus_{\dd\in \BoN^{Q_0}} \BPSh^{T}_{\tilde Q,\tilde W,\dd},
\]
which is a perverse sheaf on $\bigsqcup_{\dd\in \BoN^{Q_0}}\CM^T_{\dd}(\tilde Q)$. The stack $\FM^{T}_{\dd}(\tilde  Q)$ carries a canonical tautological line bundle for every vertex $i\in Q_0$, inducing a canonical map $\eta\times \id\colon \FM^{T}_{\dd}(\tilde Q)\rightarrow \B \BoC^*\times \FM^{T}_{\dd}(\tilde Q)$ where $\eta$ is the classifying morphism of the determinant line bundle. This morphism induces a $\HO(\B\BoC^*, \BoQ^{\vir})$-action on $ \CoHA^{T}_{\tilde Q,\tilde W}$ and hence, via the split inclusion $\BPSh^{T}_{\tilde Q,\tilde W,\dd}\hookrightarrow \CoHA^{T}_{\tilde Q,\tilde W}$ coming from \eqref{eq: purity coha}, a canonical morphism
\begin{equation}
    \label{eq: BPS to COHA}
    \HO(\B\BoC^*, \BoQ^{\vir})\otimes\BPSh^{T}_{\tilde Q,\tilde W}\to \CoHA^{T}_{\tilde Q,\tilde W}.
\end{equation}
We also define
\begin{align*}
\Fg^{T}_{Q,W,\dd}\coloneqq
&\HO^*(\CM^{T}_{\dd}(\tilde Q),\BPSh^{T,\zeta}_{\tilde Q,\tilde W,\dd}[-1]).
\end{align*}
and $\Fg^{T}_{\tilde Q,\tilde W}=\bigoplus_{\dd\in \BoN^{Q_0}} \Fg^{T}_{\tilde Q,\tilde W,\dd}$.
The split inclusion $\BPSh^{T}_{\tilde Q,\tilde W,\dd}\hookrightarrow \CoHA^{T}_{\tilde Q,\tilde W}$ induces an injective  morphism $\iota\colon \Fg^{T}_{\tilde Q,\tilde W}\rightarrow \HO\!\CoHA^{T}_{\tilde Q,\tilde W}$. By \cite[Thm.C]{QEAs} , the image is closed under the commutator and hence defines a Lie subalgebra of $\HO\!\CoHA^{T}_{\tilde Q,\tilde W}$. 

We now come to the integrality morphism for the Hall algebra $\HCoha^T_{\tilde Q, \tilde W}$. By \eqref{eq: BPS to COHA}, we obtain a canonical morphism
\[
\Omega\colon \Sym_{\HO_T}\!\left(\HO(\B\BoC^*)\otimes \BPSh^{T}_{\tilde Q,\tilde W} \right)\rightarrow \CoHA^{T}_{\tilde Q,\tilde W}
\]
and hence, taking derived global sections, we obtain a canonical morphism 
\[
\HO\!\Omega\colon \Sym_{\HO_T}\!\left(\HO_{\BoC^*}\otimes \Fg^{T}_{\tilde Q,\tilde W} \right)\rightarrow \HO\!\CoHA^{T,\zeta}_{\tilde Q,\tilde W}
\]
where $\Sym_{\HO_T}\!\left(\HO_{\BoC^*}\otimes \Fg^{T}_{\tilde Q,\tilde W}\right)$ is the direct sum for $n\geq 0$ of the $\FS_n$-invariant subspaces of $(\HO\!\CoHA^{T,\zeta}_{\tilde Q,\tilde W})^{\otimes _{\HO_T}n}$, where the symmetric group $\FS_n$ acts via the symmetric tensor structure on cohomologically graded, $\BoN^{Q_0}$-graded vector spaces discussed e.g. in \cite[\S2.5]{BD2023}.
\begin{theorem} \cite[Thm.C]{QEAs}
\label{PBW_thm}
The PBW morphism $\Omega$ is an isomorphism. Hence, the same holds for $\HO\!\Omega$.
\end{theorem}

\subsection{2d Preprojective CoHA and Yangians}
Let $p: \FM_{\dd}^T(\tilde{Q})\rightarrow \FM_{\dd}^T(\overline{Q})$ and $\iota:  \FM_{\dd}^T(\Pi_{Q})\rightarrow \FM_{\dd}^T(\overline{Q})$ be the natural morphisms. By Davison's dimensional reduction—Theorem \ref{thm: dimred}—the natural morphism 
\begin{equation}
    \label{eq: dim red coha isom relative}
    \iota_*\BoD\BoQ_{\FM^{T}_{\dd}(\Pi_Q)}^{\vir}\rightarrow p_*\phip{\Tr(\tilde{W})}\BoQ^{\vir}_{\FM_{\dd}(\tilde{Q})}
\end{equation}
is an isomorphism. We set $\BoD\BoQ_{\FM^{T}_{\dd}(\Pi_Q)}^{\vir}\coloneqq (\BoD\BoQ_{\FM^{T}_{\dd}(\Pi_Q)})[(\dd,\dd)_Q]$ and, following \cite{preproj}, define the $\BoN^{Q_0}$-graded, cohomologically graded $\HO_T$-module $\HO\!\CoHA_{\Pi_Q}^T$ by setting 
\begin{align*}
\HCoha_{\Pi_Q,\dd}^{T}
\coloneqq&\HO^*(\FM^{T}_{\dd}(\Pi_Q),\BoD\BoQ_{\FM^{T}_{\dd}(\Pi_Q)}^{\vir})\qquad \HCoha_{\Pi_Q}^{T}
\coloneqq\bigoplus_{\dd\in \BoN^{Q_0}} \HCoha_{\Pi_Q,\dd}^{T}
\end{align*}
We use the isomorphism $\HCoha_{\Pi_Q}^{T}\cong \HCoha^{T}_{\tilde Q, \tilde W}$ induced by \eqref{eq: dim red coha isom relative} to equip the left hand side with an algebra structure. This is known as the \emph{Preprojective CoHA}. 
\begin{remark}
    The algebra structure on $\HCoha_{\Pi_Q,\dd}^{T}$ can be directly defined using the geometry of the stack $\FM(\Pi_Q)$. This approach was originally formulated by \cite{ScVa13} for the Jordan quiver and then extended by \cite{YZ18,ScVa20} for arbitrary quivers. The equivalence\footnote{Up to a sign.} of the two algebra structures is proved in Davison's appendix to \cite{RS17}.
\end{remark}

\begin{theorem}[{\cite{BD2023}}]
\label{thm: BD}
    There is an isomorphism of $\BoZ^{Q_0}$-graded, cohomologically graded Lie algebras $\fg^{T}_{\Pi_Q}\cong \fg^{\MO,T}_Q$. This isomorphism extends to an isomorphism of $\BoN^{Q_0}$-graded, cohomologically graded algebras $\Yang^{\MO, +}_Q \cong \HO\!\CoHA_{\Pi_Q}^T$, intertwining the respective actions on cohomology of Nakajima quiver varieties.
\end{theorem}
The isomorphism $\Yang^{\MO, +}_Q \cong \HO\!\CoHA_{\Pi_Q}^T$ was proved independently in \cite{SV23O}.

The identification between the preprojective Hall algebra $\HCoha^T_{\Pi_Q}$ and the critical Hall algebra $\HCoha^T_{\tilde Q, \tilde W}$ has been of fundamental importance for relating cohomological DT theory to geometric representation theory, and the original proof of Theorem \ref{thm: BD} heavily depends on it too\footnote{An alternative proof, which solely relies on critical cohomology, will appear in \cite{BDZZ}}. However, in this article, we will make extensive use of a \emph{different} dimensional reduction of the critical Hall algebra $\HCoha^T_{\tilde Q, \tilde W}$, see \S\ref{subsec: Dimensional reduction and comparison with the 3d CoHA}. Amusingly, this alternative dimensional reduction will prove crucial to relate Hall algebras and Coulomb branches geometrically, see\S\ref{subsec: From hall to Coulomb}.

\subsection{Nilpotent and $\omega$-nilpotent 3d CoHA}
\label{subsec: Nilpotent and omega-nilpotent 3d CoHA}
In this section, we introduce two variants of the 3d Hall algebra from \S\ref{subsec: The three dimensional CoHA}. From the perspective of the paper, the first variant, which we call the $\omega$-nilpotent 3d CoHA, corresponds to performing Hecke modifications supported at $0\in \bbP^1$ as opposed to an arbitrary number of points in $\bbP^1\setminus \{\infty\}$. The second variant, the nilipotent CoHA, will be related to shifted quantum group actions on quasimaps' vanishing cycle cohomology.
Considering the natural morphisms
\[
\{0\}\xrightarrow{\iota_{\dd,0}}\times \CM_{\dd}^{T,\onil}(\tilde Q)\xrightarrow{\iota_{\dd}} \CM^T_{\dd}(\tilde Q)
\]
Here, $0$ is understood as the zero representation, and $\CM_{\dd}^{T,\onil}(\tilde Q)$ is the zero locus of the functions 
\[
\Tr(\omega^k_i)\qquad k=1\dots, \dd_i\; i\in Q_i.
\]
Taking pullback diagrams, we get closed substacks 
\begin{equation}
    \label{eq: nilp and omega nilp tripled qivers}
    \begin{tikzcd}
     \FM^{T,\nil}_{\dd}(\tilde Q)\arrow[r, "\eta_{\dd, 0}"]\arrow[d] & \FM^{T,\onil}_{\dd}(\tilde Q)  \arrow[r, "\eta_{\dd}"]\arrow[d, "\JH"]& \FM^T_{\dd}(\tilde Q) \arrow[d, "\JH"]\\
     \{0\} \arrow[r, "\iota_{\dd, 0}"]& \CM_{\dd}^{T,\onil}(\tilde Q)\arrow[r, "\iota_{\dd}"] & \CM^T_{\dd}(\tilde Q)
\end{tikzcd}
\end{equation}
Recall that $\JH$ is the semi-simplification morphism. Hence, $\FM^{T, \onil}_{\dd}(\tilde Q)$ is the stack of $\dd$-dimensional $\BoC\tilde Q$-representations with nilpotent loops $\omega_i$ and $\FM^{T,\nil}_{\dd}(\tilde Q)$ is the stack of nilpotent $\BoC\tilde Q$-representations. We then define the 3d relative $\omega$-nilpotent CoHA by setting
\[
\Coha^{T,\onil}_{\tilde{Q}, \tilde{W}}\coloneqq\iota^! \bigoplus_{\dd\in \BoN^{Q_0}}\Coha^T_{\tilde{Q}, \tilde{W},\dd}=\bigoplus_{\dd\in \BoN^{Q_0}}\iota_{\dd}^!\Coha^T_{\tilde{Q}, \tilde{W},\dd}
\]
and 
\[
\HCoha^{T,\onil}_{\tilde{Q}, \tilde{W}}=(\bigoplus_{\dd\in \BoN^{Q_0}} \HO\iota_{\dd}^!\Coha^T_{\tilde{Q}, \tilde{W},\dd})
\]
Similarly, the $3d$ nilpotent CoHA is defined by further pulling back to along $\iota_0$:
\[
\Coha^{T,\nil}_{\tilde{Q}, \tilde{W}}\coloneqq\iota_0^!\iota^! \bigoplus_{\dd\in \BoN^{Q_0}}\Coha^T_{\tilde{Q}, \tilde{W},\dd}=\bigoplus_{\dd\in \BoN^{Q_0}}\iota_{\dd,0}^!\iota_{\dd}^!\Coha^T_{\tilde{Q}, \tilde{W},\dd}
\]
\[
\HCoha^{T,\nil}_{\tilde{Q}, \tilde{W}}=\bigoplus_{\dd\in \BoN^{Q_0}} \HO(\iota_{\dd,0}^0\iota_{\dd}^!\Coha^T_{\tilde{Q}, \tilde{W},\dd}).
\]
Since the map 
\[
\iota: \bigsqcup_{\dd\in \BoN^{Q_0}} \{0\}\times \CM^T_{\dd}(\bar Q)\to \bigsqcup_{\dd\in \BoN^{Q_0}} \CM^T_{\dd}(\tilde Q)
\]
is closed, we can apply the canonical adjunctions $\iota_*\iota^!\to id$ and $\iota_{0}*\iota^!_{0}\to id$  to the sheaf $\Coha^{T,\onil}_{\tilde{Q}, \tilde{W}}$ to get morphisms
\begin{equation}
    \label{eq: map nilp sheafy coha}
\iota_{0,*}\iota_*\Coha^{T,\nil}_{\tilde{Q}, \tilde{W}}\to \Coha^T_{\tilde{Q}, \tilde{W}}\qquad \iota_*\Coha^{T,\onil}_{\tilde{Q}, \tilde{W}}\to \Coha^T_{\tilde{Q}, \tilde{W}}.
\end{equation}
and hence, passing to derived global section, maps
\begin{equation}
    \label{eq: map nilp coha}
    \HCoha^{T,\nil}_{\tilde{Q}, \tilde{W}}\to \HCoha^{T,\onil}_{\tilde{Q}, \tilde{W}}\to \HCoha^T_{\tilde{Q}, \tilde{W}}
\end{equation}
\begin{proposition}
    Both maps \eqref{eq: map nilp sheafy coha} are morphisms of algebras in $D(\CM^T(\tilde Q))$ with respect to the monoidal structure given by \eqref{eq: monoidal structure good moduli space}. Therefore, both maps in \eqref{eq: map nilp coha} are morphisms of $\HO_T$-algebras.
\end{proposition}
\begin{proof}[Proof sketch]
    The logic of the proof is the same for both maps, so we only discuss the case of $\iota_*\Coha^{T,\onil}_{\tilde{Q}, \tilde{W}}\to \Coha^T_{\tilde{Q}, \tilde{W}}$. The space $\FM^{\onil}(\tilde Q)$ is the moduli stack parametrizing representations of $\tilde Q$ whose ``tripling'' loops $\omega_i$ acts nilpotently. The latter form a Serre subcategory of the category of representations of $\tilde Q$. In other words, given an exact sequence 
    \[
    0\to M_1\to M_2\to M_3\to 0
    \]
    of $\BoC Q$-modules, then the loops $\omega_i$ in $M_2$ acts nilpotently iff those in $M_1$ and $M_3$ do so. The result then follows from a standard base change argument.
\end{proof}

Consider now the canonical morphism $m_{\dd}: \AA^1/T\times_{\B T} \CM^T_{\dd}(\Pi_Q)\to \CM^T_{\dd}(\tilde Q)$ obtained by extending a $\Pi_Q$-module to a $\BoC Q$ module by setting the additional loops $\omega_i$ to be equal to the diagonal matrix with constant eigenvalue equal to $x\in \AA^1$. We now recall Davison's support lemma.
\begin{lemma}[{\cite[Lemma 4.1 and \S4.2]{preproj}}]
\label{lma: support lemma}
    The 3d BPS sheaf $\BPSh^{T}_{\tilde Q,\tilde W,\dd}$ is isomorphic to 
    \[
    (m_{\dd})_* \left(\BoQ^{\vir}_{\AA^1/T}\boxtimes_{\B T} \BPSh^{T}_{\Pi_Q, \dd}\right).
    \]
\end{lemma}

\begin{lemma}
\label{lemma: nilpotent BPS}
    The $\omega$-nilpotent BPS sheaf $\iota_{\dd}^! \BPSh^{T}_{\tilde Q,\tilde W,\dd}$ is isomorphic to $\BoQ_0[-1]\boxtimes \BPSh^{T}_{\Pi_Q, \dd}$.
    Moreover, under the isomorphism and Lemma \ref{lma: support lemma}, the natural map
    \begin{equation}
        \label{eq: nilp BPS and BPS}
        \HO(\BPSh^{T}_{\Pi_Q, \dd})=\HO(\iota_{\dd}^!\BPSh^{T}_{\tilde Q,\tilde W,\dd})\to \HO(\BPSh^{T}_{\tilde Q,\tilde W,\dd})=\HO(\BPSh^{T}_{\Pi_Q, \dd}[2])
    \end{equation}
    induced by the adjunction $(\iota_{\dd})_*\iota_{\dd}^! \to id $ is given by the action of the first Chern class of the character $\epsilon_1$ of $T$. In particular, this map is injective.
    \begin{proof}
        Consider now the following Cartesian diagram 
        \[
         \begin{tikzcd}
          \CM^T_{\dd}(\Pi_Q)\arrow[r, "j_{\dd}"]\arrow[d, "\iota_{\dd}^0"] & \CM_{\dd}^{T,\onil}(\tilde  Q)\arrow[d, "\iota_{\dd}"]\\
          \AA^1\times_{\B T} \CM^T_{\dd}(\Pi_Q)\arrow[r, "m_{\dd}"] & \CM^T_{\dd}(\tilde Q)
          \end{tikzcd}
        \]
        By base change, we get
        \begin{align*}
            \iota_{\dd}^! \BPSh^{T}_{\tilde Q,\tilde W,\dd} 
            & = \iota_{\dd}^! (m_{\dd})_* \left(\BoQ^{\vir}_{\AA^1/T}\boxtimes_{\B T}\BPSh^{T}_{\Pi_Q, \dd}\right)
            \\
            & = (j_{\dd})_* (\iota^0_{\dd})^! \BoQ^{\vir}_{\AA^1/T}\boxtimes_{\B T} \BPSh^{T}_{\Pi_Q, \dd}
            \\
            &=\BoQ_{0/T}[-1]\boxtimes_{\B T} \BPSh^{T}_{\Pi_Q, \dd},
        \end{align*}
    where the last equation follows because
    \[
    (\{0\}\hookrightarrow \AA^1/T)^! \BoQ^{\vir}_{\AA^1/T}\boxtimes_{\B T}=(\{0\}\hookrightarrow \AA^1/T)^! \DD\BoQ_{\AA^1}[-1]=\DD\BoQ_{0/T}[-1]=\BoQ_{0/T}[-1].
    \]
    For the second statement, notice that the previous argument implies that the map \eqref{eq: nilp BPS and BPS} is induced by pushforward along the inclusion $\{ 0\}\to \BoC$, which is manifestly given by multiplication by the torus character of $\AA^1$, which is $\epsilon_1$.
    \end{proof}
    
    \begin{corollary}
        The map $\HCoha^{T,\onil}_{\tilde{Q}, \tilde{W}}\to \HCoha^T_{\tilde{Q}, \tilde{W}}$ is injective in $T$-equivariant cohomology. The same is true for the map $\HCoha^{T,\nil}_{\tilde{Q}, \tilde{W}}\to \HCoha^{T,\onil}_{\tilde{Q}, \tilde{W}}$.
    \end{corollary}
    \begin{proof}
        The proof of the first statement follows from the integrality Theorem \ref{PBW_thm} and Lemma \ref{lemma: nilpotent BPS}. For the second statement,we can factor the morphism as $\HCoha^{T,\nil}_{\tilde{Q}, \tilde{W}}\to \HCoha^{T,\onil}_{\tilde{Q}, \tilde{W}}$ as 
        \[
        \HCoha^{T,\nil}_{\tilde{Q}, \tilde{W}}\to \HCoha^{T, \bar Q-\nil}_{\tilde Q, \tilde W} \to \HCoha^{T}_{\tilde Q, \tilde W}.
        \]
        where $\HCoha^{T, \bar Q-\nil}_{\tilde Q, \tilde W}$ is the Hall algebra of $\tilde Q$-representations such that the representation of the sub-quiver $\bar Q\subset \tilde Q$ is nilpotent. The rightmost map is injective by \cite[Rem. 10.4.]{preproj} and the leftmost map by another application of integrality \ref{PBW_thm} and Lemma \ref{lemma: nilpotent BPS}. 
    \end{proof}
\end{lemma}
\begin{remark}
\label{rem: inj omega nilp lie algebra}
    The injection $\HCoha^{T,\onil}_{\tilde{Q}, \tilde{W}}\to \HCoha^T_{\tilde{Q}, \tilde{W}}$ induces a map $ \Fg^{T,\onil}_{\tilde Q, \tilde W}\to \Fg^{T}_{\tilde Q, \tilde W}$.
    This map is injective morphism of Lie algebras, but it is not an isomorphism. This simply follows from the observation that $\Fg^{T}_{\tilde Q, \tilde W}$ is a flat $\HO_T$-module \cite[Thm. 9.6]{preproj} and the image of $\Fg^{T,\onil}_{\tilde Q, \tilde W}\to \Fg^{T}_{\tilde Q, \tilde W}$ lands in degree one for the graded induced by the polynomial ring $\HO_T$. However, the $\BoN^Q_0$-graded Lie algebras $\Fg^{T,\onil}_{\tilde Q, \tilde W}$ and $\Fg^{T}_{\tilde Q, \tilde W}$ become isomorphic upon tensoring with $\Frac(\HO_{T})$. On the other hand, for $T=1$ the Lemma \ref{lemma: nilpotent BPS} and flatness of $\Fg^{T}_{\tilde Q, \tilde W}$ as a $\HO_T$-module imply that the morphism $\Fg^{\onil}_{\tilde Q, \tilde W}\to \Fg^{}_{\tilde Q, \tilde W}$ is zero\footnote{Notice there is no contradiction here because $\mathfrak{g}^T_{\tilde Q, \tilde W}/\Fg^{T,\onil}_{\tilde Q, \tilde W}$ is not flat over $\HO_T$.} but the graded dimensions of the two Lie algebras are the same. It would be interesting to determine the Lie algebra structure of the Lie $\BoQ$-algebra $\Fg^{\onil}_{\tilde Q, \tilde W}$.
\end{remark}

In this article, we will also consider the spherical subalgebra of $\HCoha^{T,\nil}_{\tilde Q, \tilde W}$, defined as the subalgebra $\Sh{SH}^{T,\nil}$ generated by $\HCoha^{T,\nil}_{\tilde Q, \tilde W, \delta_i}$ for all $i\in Q_0$.
Overall, this gives four algebras fitting in a chain of inclusions 
\[
\Sh{SH}^{T,\nil}\subset \HCoha^{T,\nil}_{\tilde{Q}, \tilde{W}}\subset \HCoha^{T, \onil}_{\tilde Q, \tilde W} \subset \HCoha^{T}_{\tilde Q, \tilde W}.
\]
\begin{remark}
\label{rem: signed product hall algebras}
For each of these algebras, it will be useful to consider a variant where the multiplication map is twisted by a sign. It suffices to define the sign twist for $\HCoha^{T}_{\tilde Q, \tilde W}$. We define $\HCoha^{T, \circ}_{\tilde Q, \tilde W}$ as the algebra whose underlying $\HO_T$-module is the same as $\HCoha^{T}_{\tilde Q, \tilde W}$ and such that the summand $\vmult_{\dd',\dd''}$ is multiplied by $(-1)^{(\dd'')^T\Q\dd'}$. It is easy to check that this convention gives a well-defined algebra structure. We denote the various subalgebras accordingly: 
\[
\Sh{SH}^{T,\nil, \circ}\subset \HCoha^{T,\nil, , \circ}_{\tilde{Q}, \tilde{W}}\subset \HCoha^{T, \onil, \circ}_{\tilde Q, \tilde W} \subset \HCoha^{T,\circ}_{\tilde Q, \tilde W}.
\]
This sign convention will be suitable for dimensional reduction statements. We will use the same notation for the relative versions of these algebras, namely $\CoHA^{T,\nil, \circ}_{\tilde{Q}, \tilde{W}}$, $\CoHA^{T, \onil, \circ}_{\tilde Q, \tilde W}$, and  $\CoHA^{T, \circ}_{\tilde Q, \tilde W}$. The same sign convention for the multiplication is understood. 
\end{remark}


\subsection{The Hall algebra of $\Coh_0(\AA^1)$-valued quiver representations}
\label{subsec: Q decorated coherent sheaves on A1}
Let $\Coh_0(\AA^1)$ be the category of zero dimensional coherent sheaves on $\AA^1$.
For a fixed quiver $Q$, set 
\[
\Coh_{Q,0}(\AA^1)\coloneqq\Rep_{Q}(\Coh_0(\AA^1)),
\]
be the abelian category of $\BoC Q$-representations in $\Coh_0(\AA^1)$. Explicitly, the objects are pairs 
\[
(\{\Sh{F}_i\}_{i\in Q_0}, \{\alpha_a\}_{a\in Q_1})
\]
where each $\Sh{F}_i$ is a zero dimensional coherent sheaf on $\AA^1$ and $\alpha_a: \Sh{F}_{s(a)}\to \Sh{F}_{t(a)}$ are morphisms of sheaves. A morphism 
\[
(\{\Sh{F}_i\}_{i\in Q_0}, \{\alpha_a\}_{a\in Q_1})\to (\{\Sh{G}_i\}_{i\in Q_0}, \{\beta_a\}_{a\in Q_1})
\]
is a collection of morphisms of sheaves $h_i: \Sh{F}_i\to \Sh{G}_i$ such that $h_{t(a)}\alpha_{a}=\beta_{t(a)}h_{s(a)}$ for all $a\in Q_1$. It is clear that $\Coh_{Q,0}(\AA^1)$ is an abelian category recovering $\Coh_0(\AA^1)$ for $Q=A_1$. Let $\FCoh_{Q,0}(\AA^1)$ be the associated moduli stack. We have a morphism 
\[
\FCoh_{Q,0}(\AA^1)\to \BoN^{Q_0} \quad (\{\Sh{F}_i\}_{i\in Q_0}, \{\alpha_a\}_{a\in Q_1})\to \text{length}(\Sh{F_i})_{i\in Q_0}
\]
and an associated decomposition
\[
\FCoh_{Q,0}(\AA^1)=\bigsqcup_{\dd\in Q_0}\FCoh_{Q,0,\dd}(\AA^1),
\]
The moduli stacks $\FCoh_{Q,0,\dd}(\AA^1)$ can be described explicitly via Proposition \ref{prop: equivalence of categories of coh and reps}. Firstly, any zero dimensional coherent sheaf on $\AA^1=\Spec(\BoC[x])$ of length $d$ is equivalent to a $\dd$-dimensional vector space $V=\HO(\AA^1, \Sh{F})$ together with an endomorphism $\omega$ prescribing the action of $x\in \BoC[x]$. 
Hence, the category $\Coh=\Coh_{A_1}$ is equivalent to the category of finite dimensional modules of the one loop quiver $\hat A_1$. More generally, by Proposition \ref{prop: equivalence of categories of coh and reps}, we have an equivalence $\Coh_{Q,0}(\AA^1)\cong \Rep_{\Xi_Q}$, cf \S\ref{subsec: Doubling, looping, tripling, and framing}. Notice that the loop variables $\omega_i$ prescribe the action of $x\in \BoC[x]$ on $\Sh{F}_i$ and the evaluation of $a\in Q_1$ prescribe the morphisms $\alpha_a: \Sh{F}_{s(a)}\to \Sh{F}_{t(a)}$. The relations $a\omega_i -\omega_i a$ ensure that the maps $\alpha_a$ are morphisms of sheaves, that is, they commute with the action of $\BoC[x]$. The equivalence of categories induces an equivalence of stacks
\begin{equation}
    \label{eq: stack CohQ and reps}
    \FCoh_{Q,0,\dd}(\AA^1)\simeq \FM_{\dd}(\Xi_Q),
\end{equation}
where $\FM(\Xi_Q)$ is the moduli stack of $\dd$-dimensional representations of the algebra, see \S\ref{subsec: Doubling, looping, tripling, and framing}. Explicitly, $\FM(\Xi_Q) \simeq \AA_{\dd}(\Xi_Q)/\GL_{\dd}$, where $\AA_{\dd}(\Xi_Q)$ is the zero locus of the function
\begin{equation}
    \label{eq: condition Xi_Q prequotient}
    \nu_{\dd}: \AA_{\dd}(\hat Q)\to \bigoplus_{a\in Q_1}\Hom(\BoC^{\dd_{s(a)}}, \BoC^{\dd_{t(a)}})
\end{equation}
given by evaluating \eqref{eq: hat quver relation} on a representation. By assumption, the action of the torus $T$ on $\FM(\tilde Q)$ introduced in \S\ref{subsec: The three dimensional CoHA} preserves $\FM(\bar Q)$ and further restricts\footnote{Since the potential $\tilde W$ is fixed by $T$, our weighting condition implies that and the relation $a\omega a_i-\omega_i a$ is only scaled by a weight.} to $\FM(\Xi _Q)$. We denote by $\FM^T(\Xi_Q)$ and $\FCoh^T_{Q,0,\dd}(\AA^1)$ the corresponding stack quotients.

The category $\Coh_{Q,0,\dd}(\AA^1)$ or, equivalently, $\Rep_{\Xi_Q}$, is naturally a Hall category in the sense of \cite{bu2025intrinsicdonaldsonthomastheoryi}, and the stacks \eqref{eq: stack CohQ and reps} admit quasi-smooth derived enhancements. Consequently, their Borel-Moore cohomology 
\[
\HCoha^T_{\Xi_Q}\coloneqq\bigoplus_{\dd\in \BoN^{Q_0}}\HO(\FM^T_{\dd}(\Xi_Q), \DD\BoQ^T_{\FM_{\dd}(\Xi_Q)})
\]
admits a Hall algebra structure, see. \cite[\S8.1.8]{bu2025cohomologysymmetricstacks}. Notice that the stack $\FM^T_{\dd}(\Xi_Q)$ has virtual dimension zero.
In fact, the Hall algebra can be defined at the level of sheaves on the coarse moduli spaces $\CM^T(\Xi_Q)$. The construction parallels the one for the 3d CoHA of a quiver with potential, cf. \ref{subsec: The three dimensional CoHA}. Consider the monoidal structure on $D_c(\CM^T(\Xi_Q))$ defined on connected components by $\Sh{F_{\dd'}}\otimes_{\oplus}\Sh{F_{\dd''}}=\oplus_*\left(\Sh{F_{\dd'}}\boxtimes\Sh{F_{\dd''}}\right)$, where $\oplus$ is the direct sum map
\[
\oplus: \CM^T_{\dd'}(\Xi_Q)\times_{\B T} \CM^T_{\dd''}(\Xi_Q)\to \CM^T_{\dd'+\dd''}(\Xi_Q).
\]
Set $\Coha^T_{\Xi_{\Q, \dd}}=(\JH_{\dd})_*\DD\BoQ_{\FM^T_{\dd}(\Xi_Q)}^{\vir}$ and $\Coha^T_{\Xi_Q}=\oplus_{\dd} \Coha^T_{\Xi_{Q, \dd}}$, where $\JH_{\dd}: \FM^T_{\dd}(\Xi_Q)\to \CM^T_{\dd}(\Xi_Q)$ is the semi-simplification map. We can explicitly describe multiplication map
\[
m^{\Xi_Q}: \Coha^T_{\Xi_Q}\otimes_{\oplus} \Coha^T_{\Xi_Q}\to \Coha^T_{\Xi_Q}
\]
defining the algebra structure on $\Coha^T_{\Xi_Q}$ using the commutative diagram
\begin{equation}
\label{eq: 2d coha coh multiplication diagram}
\begin{tikzcd}
    \FM^T_{\dd'}(\Xi_Q)\times \FM^T_{\dd''}(\Xi_Q)\arrow[d, "\JH_{\dd'}\times\JH_{\dd''}"]& \FM^T_{\dd',\dd''}(\Xi_Q)\arrow[l, swap,  "\bar q"]\arrow[r, "\bar p"] & \FM^T_{\dd'+\dd''}(\Xi_Q)\arrow[d, "\JH_{\dd}"]
    \\
    \CM^T_{\dd'}(\Xi_Q)\times \CM^T_{\dd''}(\Xi_Q)\arrow[rr, "\oplus"] & &  \CM^T_{\dd'+\dd''}(\Xi_Q)
\end{tikzcd}
\end{equation}
and the embedding $\FM^T(\Xi_Q)\hookrightarrow \FM^T(\hat Q)$ induced by the quotient $\Xi_Q=\BoC \hat Q/(\sum_{a, 
\omega} a\omega -\omega a)$. Explicitly, for any given dimension vector $\dd\in \BoN^{Q_0}$, this embedding is presented as
\[
\Rep_{\Xi_Q}(\dd)/(\GL_{\dd}\times T)\hookrightarrow \Rep_{\hat Q}(\dd)/(\GL_{\dd}\times T).
\]
Set 
\[
Z_{\dd',\dd''}\coloneqq\left(\Rep_{\hat Q}(\dd')\times \Rep_{\hat Q}(\dd'')\times \bigoplus_{a\in Q_1}\Hom(\BoC^{\dd''}_{s(a)}, \BoC^{\dd'}_{t(a)}) \right)/(\GL_{\dd'}\times\GL_{\dd''}\times T)
\]
Consider the diagram 
\begin{equation}
    \label{eq: diagram coha Xi_Q}
\begin{tikzcd}
    Z_{\dd',\dd''} & \FM^T_{\hat Q, \dd', \dd''} \arrow[l, swap, "q"]\arrow[r, "p"] & \FM^T_{\hat Q, \dd'+\dd''}
    \\
    \FM^T_{\Xi_Q, \dd'}\times \FM^T_{\Xi_Q, \dd''} \arrow[u, "k"]& \arrow[u]\FM^T_{\Xi_Q, \dd', \dd''} \arrow[l, swap, "\bar q"]\arrow[r, "\bar p"] & \FM^T_{\Xi_Q, \dd'+\dd''}\arrow[u]
\end{tikzcd}
\end{equation}
where $\FM^T_{\Xi_Q, \dd''}$ and $\FM^T_{\hat Q, \dd', \dd''}$ are the usual stacks of short exact sequences. The map $k$ is specified on the prequotients by
\[
\Rep_{\hat Q}(\dd', \dd'')\xrightarrow{a\times b\times c} \Rep_{\hat Q}(\dd')\times \Rep_{\hat Q}(\dd'')\times \bigoplus_{a\in Q_1}\Hom(\BoC^{\dd''}_{s(a)}, \BoC^{\dd'}_{t(a)}) 
\]
where $a$ and $b$ are the canonical projections and $c$ is given as follows. On the vertices $s(a)$ and $t(a)$ a representation of $\hat Q$ expressed by three matrices 
\[
B_{s(a)}=\begin{pmatrix}
    B'_{s(a)} & E_{s(a)} \\
    0 & B''_{s(a)}
\end{pmatrix}\qquad B_{s(a)}=\begin{pmatrix}
    B'_{t(a)} & E_{t(a)} \\
    0 & B''_{t(a)}
\end{pmatrix}
\qquad 
A=\begin{pmatrix}
    A' & X\\
    0 & A''
\end{pmatrix}
\]
corresponding to the arrows $\omega_{s(a)}$, $\omega_{t(a)}$, and $a$ i the path algebra $\BoC \hat Q$. Then the projection  of $c$ to the factor $\Hom(\BoC^{\dd''}_{s(a)}, \BoC^{\dd'}_{t(a)})$ is given by 
\begin{equation}
    \label{eq: section Coh_Q^0}
    A'E_V-E_WA''-B_W'X-XB_V''.
\end{equation}
The following lemma follows directly from the definitions.
\begin{lemma}
\label{lemma: quasi-smoothness pull 2d coha}
    The left square in diagram \eqref{eq: diagram coha Xi_Q} is Cartesian. 
\end{lemma}
Since $q$ is smooth, we have a canonical morphism $\DD\BoQ_{Z_{\dd,\dd''}}\to q_*\DD\BoQ_{\FM^T_{\hat Q, \dd',\dd''}}$ and hence, applying $k^!$ and performing base change, we obtain a canonical Gysin morphism 
\[
q^!: \DD\BoQ_{\FM^T_{\Xi_Q,\dd'}\times_{\B T}\FM^T_{\Xi_Q,\dd''}}\to \bar q_*\DD\BoQ_{\FM^T_{\Xi_Q, \dd',\dd''}}
\]
Therefore, using properness of $p$ and commutativity of diagram \eqref{eq: 2d coha coh multiplication diagram}, we obtain canonical morphisms as the composition
\begin{align}
\label{eq: sheafy xi mult}
\begin{split}
    \Coha^T_{\Xi_{Q}, \dd'}\otimes_{\oplus}\Coha^T_{\Xi_{Q}, \dd'}
    & \cong \oplus_{*}(\JH_{\dd'}\times \JH_{\dd'})_* \left(\DD\BoQ_{\FM^T_{\dd'}(\Xi_Q)\times\FM^T_{\dd'}(\Xi_Q)} \right)
    \\
    &\to \oplus_{*}(\JH_{\dd'}\times \JH_{\dd'})_*\bar q_*\DD\BoQ_{\FM^T_{\Xi_Q, \dd',\dd''}}
    \\
    & = \JH_{\dd'+\dd''}\bar p_*\bar p^!\DD\BoQ_{\FM^T_{\Xi_Q, \dd'+\dd''}}
    \\
    &\to \JH_{\dd'+\dd''} \DD\BoQ_{\FM^T_{\Xi_Q, \dd'+\dd''}}
\end{split}
\end{align}

We have obtained a multiplication map $\vmult^{\Xi_Q}_{\dd',\dd''}$ $\Coha^T_{\Xi_{Q}, \dd'}\otimes_{\oplus}\Coha^T_{\Xi_{Q}, \dd'}\to \Coha^T_{\Xi_{Q}, \dd'}$. With slight abuse of notation, we use the same notation for denote for the map
\[
\vmult^{\Xi_Q}_{\dd',\dd''}: \HCoha^T_{\Xi_{Q}, \dd'}\otimes \HCoha^T_{\Xi_{Q}, \dd''}\to \HCoha^T_{\Xi_{Q}, \dd'+\dd''}.
\]
obtained by taking derived global section. 
The proof of the following theorem is standard and we omit it. 
\begin{theorem}
\label{thm: 2d CoHA mult}
        The maps $\vmult^{\Xi_Q}_{\dd',\dd''}$ equip $\Coha^T_{\Xi_{Q}}$ with the structure of an associative unital $\BoN^{Q_0}$-graded algebra object in $D_c(\CM^T(\Xi))$. Therefore, $(\HCoha^T_{\Xi_{Q}}, \vmult^{\Xi_Q})$ is an $\BoZ\times \BoN^{Q_0}$ graded unital algebra.
\end{theorem}

\begin{remark}
\label{rem: alternative pullback Coh_Q^0}
We end this section with a slightly different description of the virtual pullback $\bar q^!$, which will turn out useful later in the paper. Let $V^{\de}$ and $V^{\ob}$ be two copies of the vector space $\bigoplus_{a\in Q_1}\Hom(\BoC^{\dd''}_{s(a)}, \BoC^{\dd'}_{t(a)})$ and let $V^{\ext}=\bigoplus_{i\in Q_0}\Hom(\BoC^{\dd''_i}, \BoC^{\dd'_i})$\footnote{The superscript ``ext'' refers to the fact that this vector space encodes the deformations of the torsion sheaves, which are unobstructed. Similarly,  ``def'', and ``ob'', refer, respectively, to the space of deformation and obstruction of the morphisms between the torsion sheaves.}. Set 
\begin{align*}
    E_{\Xi_Q}&=\left(\Rep_{\hat Q}(\dd')\times \Rep_{\hat Q}(\dd'') \times V^{\ext} \times V^{\de}\times V^{\ob}\right)/(P \times T)
    \\
    E_{\Xi_Q}&=\left(\Rep_{\Xi_Q}(\dd')\times \Rep_{\Xi_Q}(\dd'') \times V^{\ext}\times V^{\de}\times V^{\ob}\right)/(P \times T)
    \\
    B_{\Xi_Q}&=\left(\Rep_{\Xi_Q}(\dd')\times \Rep_{\Xi_Q}(\dd'') \times V^{\ext}\times V^{\de}\right)/(P \times T)
    \\
    Z'_{\dd',\dd''}&=\left(\Rep_{\Xi_Q}(\dd')\times \Rep_{\Xi_Q}(\dd'') \times V^{\ext}\times V^{\ob}\right)/(\Gl_{\dd'}\times \GL_{\dd''}\times T)
\end{align*}
and notice that 
\begin{align*}
    Z_{\dd',\dd''}&=\left(\Rep_{\hat Q}(\dd')\times \Rep_{\hat Q}(\dd'') \times  V^{\ob}\right)/(\Gl_{\dd'}\times \GL_{\dd''}\times T)
    \\
    \FM^T_{\hat Q, \dd',\dd''}&=\left(\Rep_{\Xi_Q}(\dd')\times \Rep_{\Xi_Q}(\dd'') \times V^{\ext}\times V^{\de}\right)/(P \times T)
\end{align*}
Consider the diagram 
\[
\begin{tikzcd}
    Z_{\dd',\dd''} & E_{\hat Q}\arrow[l]&  \FM^T_{\hat Q, \dd',\dd''} \arrow[ll, swap, bend right, "q"]\arrow[l]
    \\
    Z'_{\dd',\dd''} \arrow[u]& E_{\Xi_Q}\arrow[l]\arrow[u]&  B_{\Xi_Q} \arrow[l, swap, "s_{\Xi_Q}"]\arrow[u]
    \\
    \FM^T_{\Xi_Q, \dd'}\times \FM^T_{\Xi_Q, \dd''} \arrow[u]\arrow[uu, bend left, shift left=2, "k"]& B_{\Xi_Q}\arrow[l, swap, "\pi"]\arrow[u, "o"]& \arrow[u, "i"]\FM^T_{\Xi_Q, \dd', \dd''} \arrow[ll, swap, bend left, "\bar q"]\arrow[l, "i"]
\end{tikzcd}
\]
All the unlabelled map are the natural inclusions or projections. The map $\pi$ is the canonical projection. The map $o$ is the zero section of the vector bundle $E_{\Xi_Q}\to B_{\Xi_Q}$. The map $s_{\Xi_Q}$ is also a section of the same bundle, and it is given given by \eqref{eq: section Coh_Q^0}. The diagram commutes. As a consequence, a standard base change argument shows that the Gysin morphism $q^!:  \DD\BoQ_{\FM^T_{\Xi_Q,\dd'}\times_{\B T}\FM^T_{\Xi_Q,\dd''}}\to \bar q_*\DD\BoQ_{\FM^T_{\Xi_Q, \dd',\dd''}}$ equals the composition
\begin{multline*}
\DD\BoQ_{\FM^T_{\Xi_Q,\dd'}\times_{\B T}\FM^T_{\Xi_Q,\dd''}}\to \pi_*\pi^* \DD\BoQ_{\FM^T_{\Xi_Q,\dd'}\times_{\B T}\FM^T_{\Xi_Q,\dd''}}
\\
= \pi_*\DD\BoQ_{B_{\Xi_Q}}\xrightarrow{} \pi_*i_* \DD\BoQ_{\FM^T_{\Xi_Q, \dd',\dd''}}=\bar q_*\DD\BoQ_{\FM^T_{\Xi_Q, \dd',\dd''}},
\end{multline*}
where the second map is the virtual pullback associated to the data $(E_{\Xi_Q}, s_{\Xi_Q})$, cf. \ref{app: Refined Gysin pullback}.
\end{remark}

\subsection{Fixed support condition and nilpotency}
\label{subsec: Fixed support condition and nilpotency}

In \S\ref{subsec: Q decorated coherent sheaves on A1} we associated a (equivariant) Hall algebra $\HCoha^T_{\Xi_Q}$ to the category $\Coh_{Q,0}(\AA^1)$ of $Q$-valued zero dimensional sheaves on $\AA^1$. Consider the subcategory $\Coh^0_{Q,0}(\AA^1)$ of torsion sheaves supported at the origin $0\in \AA^1$. It is isomorphic to the category $\Coh^T_{0}(\DD)$ of zero dimensional coherent sheaves on the formal disk $\DD$. Let now $\FCoh^T_{Q,0}(\DD)$ be the associated moduli substack. It is a closed substack of $\FCoh^T_{Q,0}(\AA^1)$. Under the equivalence $\FCoh^T_{Q,0}(\AA^1)\cong \FM^T(\Xi_Q)$, which identifies the support of a sheaf with the eigenvalues of the loops $\omega_i\in \Xi_Q$, we have an identification
\[
\begin{tikzcd}
    \FCoh^T_{Q,0}(\DD)\arrow[d, equal]\arrow[r, hookrightarrow] & \FCoh^T_{Q,0}(\AA^1) \arrow[d, equal]
    \\
    \FM_{\dd}^{T,\onil}(\Xi_Q) \arrow[r, hookrightarrow] & \FM^T_{\dd}(\Xi_Q)
\end{tikzcd}
\]
where $\FM_{\dd}^{T,\onil}(\Xi_Q)=Y_{\dd}^{\onil}/(\GL_{\dd}\times T)$ and $Y_{\dd}^{\onil}$ is the subvariety of $\Rep_{\Xi_Q}(\dd)$ consisting of those representations on which the loops $\omega_i$ act nilpotently for all $i\in Q_0$. Moreover, we have a Cartesian square 
\begin{equation}
    \label{eq: niplotent 2d hall good moduli}
    \begin{tikzcd}
    \FM_{\dd}^{T,\onil}(\Xi_Q)\arrow[d, "\JH_{\dd}"] \arrow[r, hookrightarrow] & \FM^T_{\dd}(\Xi_Q)\arrow[d, "\JH_{\dd}"]
    \\
    \CM_{\dd}^{T,\onil}(\Xi_Q) \arrow[r, hookrightarrow, "\nu_{\dd}"] & \CM^T_{\dd}(\Xi_Q)
\end{tikzcd}
\end{equation}
where $\CM^{T, \onil}_{\dd}(\Xi_Q)$ is the zero locus in $\CM^T_{\dd}(\Xi_Q)$ of the functions $\Tr(\omega^k_i)$ for $k=1\dots, \dd_i\; i\in Q_i$. Set
\[
\Coha^{T, \onil}_{\Xi_Q}\coloneqq\nu_*\nu^!\Coha^T_{\Xi_Q}=\oplus_{\dd}(\nu_{\dd})^!\Coha^T_{\Xi_Q, \dd}
\]
and $\vmult^{\Xi_Q, \onil}\coloneqq\nu^!\vmult_{\Xi_Q}$. Then the pair $(\Coha^{T, \onil}_{\Xi_Q}, \vmult^{\Xi_Q})$
is an algebra object in $D_c(\CM^{T, \onil}_{\Xi_Q})$, where the latter is equipped with the monoidal structure defined similarly to the one on $D_c(\CM^T_{\Xi_Q})$. We also remark that we are crucially using the fact that $\Coh^{0}_{Q,0}(\AA^1)$ is a Serre subcategory of $\Coh^T_{Q,0}(\AA^1)$, i.e. is stable under extensions, subobjects, and quotients), so in particular the direct sum map $\oplus: \CM^{T, \onil}_{\Xi_Q}\times \CM^{T, \onil}_{\Xi_Q}\to \CM^{T, \onil}_{\Xi_Q}$ is well-defined. It then follows at once that the canonical morphism 
\[
\nu_*\Coha^{T, \onil}_{\Xi_Q, \dd}=\nu_*\nu^!\Coha^T_{\Xi_Q}\to \Coha^T_{\Xi_Q}
\]
is a morphism of algebra objects in $D_c(\CM^T_{\Xi_Q})$. Furthermore, since the diagram \eqref{eq: niplotent 2d hall good moduli} is Cartesian, we have a natural isomorphism
\[
{\CoHA}^{T, \onil}_{\dd}=\nu_{\dd}^! {\JH_{\dd}}\DD\BoQ_{\FM^T_{\dd}(\Xi_Q)}\cong {\JH_{\dd}}\DD\BoQ_{\FM^{T, \onil}_{\dd}(\Xi_Q)},
\]
so we obtain an algebra structure on the Borel-Moore cohomology 
\[
\HCoha_{\Xi_Q}^{T, \onil}=\bigoplus_{\dd\in \BoN^{Q_0}}\HO(\FM_{\dd}^{T, \onil}(\Xi_Q), \DD\BoQ_{\FM_{\dd}^{T, \onil}(\Xi_Q)})
\]
by taking derived global sections. We summarize the discussion in this section with the following proposition. The last statement concerning injectivity will be proved in the next section.

\begin{proposition} $ $
\label{prop: map coha to omega nil coha}
\begin{enumerate}
    \item The maps $\vmult^{\Xi_Q, \onil}_{\dd',\dd''}$ equip $\Coha^{T, \onil}_{\Xi_{Q}}$ with the structure of an associative unital $\BoN^{Q_0}$-graded algebra object in $D_c(\CM^{T, \onil}(\Xi))$. Therefore, $(\HCoha^{T, \onil}_{\Xi_{Q}}, \vmult^{\Xi_Q})$ is an $\BoZ\times \BoN^{Q_0}$ graded unital algebra.
    \item The canonical morphism $\nu_*\Coha^{T, \onil}_{\Xi_{Q}}\to \Coha^T_{\Xi_{Q}}$ is a morphism of algebras in $D_c(\CM^T(\Xi_Q))$.
    \item The induced map $\HCoha_{\Xi_Q}^{T, \onil}\to \HCoha^T_{\Xi_Q}$ is a morphism of algebras. Moreover, this morphism is injective in $T$-equivariant cohomology.
\end{enumerate}
\end{proposition}

\subsection{Dimensional reduction and comparison with the 3d CoHA}
\label{subsec: Dimensional reduction and comparison with the 3d CoHA}

The Hall algebra $\HCoha^T_{\Xi_{Q}}$ can be identified, via dimensional reduction, with the Hall algebra $\HCoha^T_{\tilde Q, \tilde W}$ for the tripled quiver \eqref{eq: absolute 3d coha}. Explicitly, the isomorphism is constructed as follows. Set $Y_{\dd}= \Rep_{\Xi_Q}(\dd)$. It is the closed subvariety of $\Rep_{\hat Q}(\dd)$ cut out by the relation \eqref{eq: hat quver relation}. It is clear that $\FM^T_{\dd}(\Xi_Q)=Y_{\dd}/(\GL_{\dd}\times T)$. Let $\ol Y_{\dd}$ be the preimage of $\Rep_{\Xi_Q}(\dd)\subseteq \Rep_{\hat Q}(\dd)$ along the affine bundle $\Rep_{\tilde Q}(\dd)\to \Rep_{\hat Q}(\dd)$. Denote by $k: \ol Y_{\dd}/(\Gl_{\dd}\times T)\hookrightarrow \FM^T_{\dd}(\tilde Q)$, $k': \FM^T(\Xi_Q)\hookrightarrow \FM^T_{\dd}(\hat Q)$, and $p: \FM^T_{\dd}(\tilde Q)\to \FM^T_{\dd}(\hat Q)$ and the associated inclusions and affine fibration, respectively. Consider the dimensional reduction morphism
\begin{equation}
    \label{eq: dim red 3d coha to coh_Q}
    (k')_*\DD\BoQ_{\FM^T_{\dd}(\Xi_Q)}[\dd^T\Q\dd]= p_* k_*\DD\BoQ_{\ol Y_{\dd}/(\GL_{\dd}\times T)}=p_* k_*k^! \DD\BoQ_{\FM^T_{\dd}(\tilde Q)}\to p_*\varphi_{\Tr(W)}\DD\BoQ_{\FM^T_{\dd}(\tilde Q)},
\end{equation}
c.f. Theorem \ref{thm: dimred}. 
\begin{proposition}
    The morphism \eqref{eq: dim red 3d coha to coh_Q} is an isomorphism of complexes.
\end{proposition}

\begin{proof}
    The stack $\ol Y_{\dd}/(\Gl_{\dd}\times T)$ is the substack of $\FM^T(\tilde Q)_{\dd}$ cut out by the relations \eqref{eq: hat quver relation}. These relations are obtained by taking non-commutative derivations of the potential $\tilde W=\sum_{\omega, a}\omega[a, a^*]$ along the arrows $a^*$ for all $a\in Q_1$. Moreover, the potential $\tilde W$ is linear in $a^*$, so the hypotheses of the dimensional reduction theorem \ref{thm: dimred} are satisfied and the statement follows. 
\end{proof}

Pushing forward \eqref{eq: dim red 3d coha to coh_Q} along $\FM^T(\bar Q)\xrightarrow{\JH} \CM(\bar Q)$ and using the diagram 
\begin{equation*}
    \begin{tikzcd}
        \FM^T(\Xi_Q)\arrow[r, "k'"]\arrow[d, "\JH"] & \FM^T(\bar Q) \arrow[d, "\JH"]& \FM^T(\tilde Q)\arrow[l, swap, "p"]\arrow[d, "\JH"]
        \\
         \CM^T(\Xi_Q)\arrow[r, "b"] & \CM^T(\bar Q) & \FM^T(\tilde Q),\arrow[l, swap, "a"]
    \end{tikzcd}
\end{equation*}
where all the morphisms labelled by $\JH$ are the semi-simplification maps, we obtain an isomorphism of complexes
\begin{equation}
    \label{eq: dim red 3d coha to coh_Q 2}
    b_*\Coha^T_{\Xi_Q}\to a_*\Coha^T_{\tilde Q, \tilde W}
\end{equation}
on $\CM(\tilde Q)$. Taking derived global sections, we get a morphism of $\BoZ\times \BoZ^{Q_0}$-graded vector spaces
\[
    dr: \HCoha^T_{\Xi_Q}=\bigoplus_{\dd\in \BoN^{Q_0}}\HCoha^T_{\Xi_Q, \dd}\xrightarrow{\cong} \bigoplus_{\dd\in \BoN^{Q_0}}\HCoha^T_{\tilde{Q}, \tilde{W},\dd}=\HCoha^T_{\tilde Q, \tilde W}.
\]
Recall that $\Coha^T_{\Xi_Q}$ (resp. $\Coha^T_{\tilde Q, \tilde W}$) is an algebra object in $(D_c(\CM^T(\Xi_Q)), \otimes_{\oplus})$ (resp. $(D_c(\CM^T(\tilde Q)), \otimes_{\oplus})$), therefore, both the domain and the target are naturally algebra objects in $(D_c(\CM^T(\bar Q)), \otimes_{\oplus})$. Recall also the sign convention from Remark \ref{rem: signed product hall algebras}.

\begin{theorem}
\label{thm: dr coh 3d coha}
The dimensional reduction isomorphism \eqref{eq: dim red 3d coha to coh_Q 2} intertwines the algebra structures on $b_*\Coha^{T}_{\Xi_Q}$ and $a_*\Coha^{T, \circ}_{\tilde Q, \tilde W}$. 
As a consequence, the same is true for $dr$.
\end{theorem}
The proof of the theorem parallels the argument in the appendix of \cite{RS17}, where the analog statement relating the multiplications of $\HCoha^T_{\tilde Q, \tilde W}$ and $\HCoha^T_{\Pi_Q}$ is proved. Therefore, we omit it. We stress the algebra structures are compatible on the nose because of the sign conventions from Remark \ref{rem: signed product hall algebras}.
In \S\ref{subsec: Dimensionally reduced action hall} we will sketch a proof a simple variant of this theorem. 

Analogous results also hold for the $\omega$-nilpotent versions of $\CoHA^T_{\Xi_Q}$ and $\CoHA^T_{\tilde Q, \tilde W}$. 
Consider the commutative diagrams
\[
\begin{tikzcd}
    \FM^T(\tilde Q)\arrow[r, "p"]\arrow[d, "\JH"] & \FM^T(\bar Q) \arrow[d, "\JH"] 
\\
\CM^T(\tilde Q)\arrow[r, "a"] & \CM^T(\bar Q)
\\
\CM^{T, \onil}(\tilde Q)\arrow[r, "a'"]\arrow[u, "\iota"] & \CM^{T, \onil}(\bar Q) \arrow[u, "\iota'"] 
\end{tikzcd}
\qquad 
\begin{tikzcd}
    \FM^T(\Xi_Q)\arrow[r, "k'"]\arrow[d, "\JH"] & \FM^T(\bar Q) \arrow[d, "\JH"] 
\\
\CM^T(\Xi_Q)\arrow[r, "b"] & \CM^T(\bar Q)
\\
\CM^{T, \onil}(\Xi_Q)\arrow[r, "b'"]\arrow[u, "\nu"] & \CM^{T, \onil}(\bar Q) \arrow[u, "\iota'"] 
\end{tikzcd}
\]

where we have removed the subscripts $\dd$ denoting the topological type to remove clutter. In these two diagrams, all the spaces of the form $X^{\onil}$ are the zero locus in $X$ of the functions $\Tr(\omega^k_i)$ for all $k=1\dots, \dd_i\; i\in Q_i$, so their closed points are representations (or the relevant algebra $\Xi_Q$, $\BoC \bar Q$, or $\BoC \tilde Q$) such that the loops $\omega_i$ act nilpotently for all $i\in Q_0$. All the maps are the canonical ones. Additionally, all the bottom right squares are Cartesian. 

We now apply the functor $(\iota')^!(\FM^T(\bar Q)\xrightarrow{\JH} \CM^T(\bar Q))_*$ to the isomorphism \eqref{eq: dim red 3d coha to coh_Q}. On the target side, by base change, we get 
\begin{align}
  (\iota')^!\JH_* p_*\varphi_{\Tr(
\tilde W)}\DD\BoQ^{\vir}_{\FM^T_{\dd}(\tilde Q)} \cong (a)'_* \iota^!\JH_*\varphi_{\Tr(\tilde W)}\DD\BoQ^{\vir}_{\FM^T_{\dd}(\tilde Q)}.
\end{align}
On the domain side, also via base change, we get
\[
(\iota')^!\JH_*(k')_*\DD\BoQ_{\FM^T_{\dd}(\Xi_Q)}\cong  (b)'_*\nu^!\JH_* \DD\BoQ_{\FM^T_{\dd}(\Xi_Q)}.
\]
Therefore, unwrapping the definitions, we get an isomorphism of complexes $(b)'_*\Coha^{T, \onil}_{\Xi_Q}\to (a')_*\Coha^{T, \onil}_{\tilde Q, \tilde W},$, and hence, passing to derived global sections, an isomorphism
\begin{equation}
    \label{eq: dr nilp version}
    dr: \HCoha^{T, \onil}_{\Xi_Q}\to \HCoha^{T, \onil}_{\tilde Q, \tilde W}.
\end{equation}
\begin{proposition}
\label{prop: dr nilp version}
    The map \eqref{eq: dr nilp version} induces an isomorphism of $\BoZ\times \BoN^{Q_0}$-graded algebras $\HCoha^{T, \onil}_{\Xi_Q}\cong \HCoha^{T, \onil, \circ}_{\tilde Q, \tilde W}$. The sign-twisted Hall algebra $\HCoha^{T, \onil, \circ}_{\tilde Q, \tilde W}$ is defined in Remark \ref{rem: signed product hall algebras}.
\end{proposition}
\begin{proof}
    The multiplication on $\Coha^{T, \onil}_{\Xi_Q}=\nu^!\Coha^T_{\Xi_Q}$ (resp. $\Coha^{T, \onil}_{\tilde Q, \tilde W}=\nu^!\Coha^T_{\tilde Q, \tilde W}$) is defined by applying the functor $\nu^!$ (resp. $\iota^!$) to the multiplication map of $\Coha^T_{\Xi_Q}$ (resp. $\Coha^T_{\tilde Q, \tilde W}$). Therefore, the statement follows from Theorem \ref{thm: dr coh 3d coha}.
\end{proof}


\subsection{Anti-involution}
\label{subsec: anti-involution}
In this section, we define a natural anti-involution on the algebra $\HCoha^T_{\Xi_Q}$. Later we will identify it with the the Chevalley anti-involution on the shifted Yangian introduced in \S\ref{subsec: Shifted Yangian action}.
Consider the category $\Coh_{0}(\AA^1)$ of zero dimensional sheaves on $\AA^1$. We have an auto-equivalence 
\[
\Theta: \Coh_{0}(\AA^1)\to \Coh_{0}(\AA^1)^{\opp}\qquad  Q\mapsto \Ext^1(Q, \Sh{O}_{\AA^1}).
\]
More generally, applying $\Ext^1(Q,-)$ induces an auto-equivalence 
\[
\Theta: \Coh_{Q,0}(\AA^1)\to \Coh_{Q^{\opp},0}(\AA^1)^{\opp}
\]
where $Q^{\opp}$ denotes the opposite quiver. Therefore, we obtain an equivalence of stacks 
\[
\Theta:\FCoh_{Q,0}(\AA^1)\to \FCoh_{Q^{\opp},0}(\AA^1)^{\opp}
\]
It is interesting to check that, under the identification $\FCoh_{Q,0}(\AA^1)\simeq \FM(\Xi_Q)$, the equivalence $\Theta$ is identified with the equivalence $\FM(\Xi_Q)\simeq \FM(\Xi_{Q^{\opp}})$ induced by taking transpose of matrices. To work equivariantly, we need to introduce the automorphism $u:T\to T$ given by $\epsilon_1\mapsto \epsilon_1$, $\epsilon_{a^*}\mapsto \epsilon_{a}$, and $\epsilon_a\mapsto \epsilon_{a^*}$ for all $a\in Q_1$. We denote by 
$\FM^{T^{\opp}}(\Xi_{Q^{\opp}})$ the quotient of $\FM^T(\Xi_{Q^{\opp}})$ by the composition of $u$ with the natural action of $T$ on $\FM^T(\Xi_{Q^{\opp}})$. We then have an equivalence $\Theta: \FM^T(\Xi_Q)\simeq \FM^{T^{\opp}}(\Xi_{Q^{\opp}})$.

\begin{proposition}
\label{prop: involution coha}
    The map $\Theta_{\dd, *}$ induces an isomorphism of $\HO_{T}$-algebras 
    \[
    \Theta_*: \HCoha^T_{\Xi_{Q^{\opp}}}\to \HCoha^{T^{\opp},\opp}_{\Xi_{Q^{\opp}}}.
    \]
    where the superscript $\opp$ in the target stands for opposite algebra structure.
    Moreover, the same construction induces an isomorphism of $\HO_{T}$-algebras 
    \[
    \Theta_*^{\onil}: \HCoha^{T,\onil}_{\Xi_{Q^{\opp}}}\to \HCoha^{T^{\opp}, \onil,\opp}_{\Xi_{Q^{\opp}}}
    \]
    fitting in the diagram 
    \begin{equation*}
        \begin{tikzcd}
            \HCoha^{T,\onil}_{\Xi_{Q^{\opp}}}\arrow{d}\arrow[r, "\Theta_*"] &  \HCoha^{T,\onil^{\opp}, \opp}_{\Xi_{Q^{\opp}}}\arrow{d}
            \\
         \HCoha^T_{\Xi_{Q^{\opp}}}\arrow[r, "\Theta_*^{\onil}"] & \HCoha^{T^{\opp}, \opp}_{\Xi_{Q^{\opp}}}
        \end{tikzcd}
    \end{equation*}
    where the vertical maps are as in Proposition \ref{prop: map coha to omega nil coha}.
\end{proposition}
\begin{remark}
\label{remark identification equivariant parameters involution}
    The $\HO_T$-action on the target of either $\Theta_*$ is modified by precomposition with the automorphism $u_*$ of $\HO_T=\BoQ[\epsilon_1, \epsilon_a, \epsilon_{a^*}]_{a\in Q_1}/(\epsilon_1 + \epsilon_a + \epsilon_{a^*})_{a \in Q_1}$ given by $\epsilon_1\mapsto \epsilon_1$, $\epsilon_a\mapsto \epsilon_{a^*}$ and $\epsilon_{a^*}\mapsto \epsilon_{a}$.
\end{remark}
\begin{proof}
Given a short exact sequence of coherent sheaves $0\to Q_1\to Q_2\to Q_3\to 0$, the long exact sequence for $Ext^*(-, \Sh{O}_{\AA^1})$ implies that we have a short exact sequence $0\to \Theta(Q_3)\to \Theta(Q_2)\to \Theta(Q_1)\to 0$. Therefore, $\Theta$ induces equivalences of stacks fitting in a commutative diagram 
\begin{equation*}
    \begin{tikzcd}
        \FM^T_{\Xi_Q, \dd'}\times \FM^T_{\Xi_Q, \dd''} \arrow[d, swap, "\sw\circ \Theta\times \Theta"]& \arrow[d, "\Theta"]\FM^T_{\Xi_Q, \dd', \dd''} \arrow[l, swap, "\bar q"]\arrow[r, "\bar p"] & \FM^T_{\Xi_Q, \dd'+\dd''}\arrow[d, "\Theta"]
    \\
    \FM^T_{\Xi_Q, \dd''}\times \FM^T_{\Xi_Q, \dd'} &\FM^T_{\Xi_Q, \dd'', \dd'}  \arrow[l, swap, "\bar q"]\arrow[r, "\bar p"] & \FM^T_{\Xi_Q, \dd'+\dd''}
    \end{tikzcd}
\end{equation*}
Here $\sw$ is simply the morphism swapping the order of the tensor factors. Moreover, since $\Theta$ is induced by matrix transposition, it is easily checked that the virtual structure as in \eqref{eq: diagram coha Xi_Q}, which induce the virtual pullback along $\bar q$, are compatible with $\Theta$. Therefore, the result follows from compatibility of virtual pullback and pushforward in Cartesian squares.
\end{proof}

\begin{remark}
    One can also construct an anti-involution $\Theta_*^{3d}:\CoHA^{T,\onil}_{\tilde Q, \tilde W}\to \CoHA^{T_{\opp}, \onil, \opp}_{\tilde Q, \tilde W}$ by interpreting $\CoHA^{T,\onil}_{\tilde Q, \tilde W}$ as the Hall algebra of $\bar Q$-decorated coherent sheaves on $\AA^1$ satisfying the moment map condition. The anti-involution is also induced by the duality functor $Ext^*(-, \Sh{O}_{\AA^1})$, and dimensionally reduced to the morphism $\Theta_{*}$ from \ref{prop: involution coha}. 
\end{remark}

\section{Quiver model and Hall action}

\subsection{Quasimaps to Nakajima quiver varieties}
\label{subsec: quasimaps to quiver varieties}

Fix a quiver $Q$ and dimension and framing vectors $\vv,\ww \in \mathbb{Z}_{\geq 0}^{Q_0}$. Let $X\coloneqq \Nak^{\zeta}_Q(\ww,\vv)$ be the corresponding Nakajima variety. It is an open substack of the stack quotient $\FX= [\mu_{\ww,\vv}^{-1}(0)/\G_{\vv}]$. Consider the mapping stack $\Map(\bbP^1, \FX)$. By definition, the moduli space of quasimaps from $\bbP^1$ to $X$ is the open substack
\begin{equation}
    \QM(\bbP^1, X) \subset \Map(\mathbb{P}^1, \mu^{-1}(0)/ \GL_{\vv} )
\end{equation}
consisting of those maps $f\in \Map(\bbP^1, \FX)$ such that for $\eta\in \bbP^1$ the generic point, we have $f(\eta)\in X$.

Notice that a map $f: \mathbb{P}^1 \dashrightarrow \FX$ from $\mathbb{P}^1$ consists, by definition, of the following data. For each $i \in Q_0$, we have a vector bundle $\Sh{V}_i$ on $\mathbb{P}^1$, together with morphisms of sheaves 
\begin{equation}
\begin{split}
    \Sh{J}_i: & \Sh{V}_i \to W_i \otimes \Sh{O}_{\mathbb{P}^1}  \\
    \Sh{I}_i: &  W_i \otimes \Sh{O}_{\mathbb{P}^1} \to \Sh{V}_i.
\end{split}
\end{equation}
The vector spaces $W_i$ are the framing spaces entering the definition of $X$, so $\dim W_i= \ww_i$ and $\text{rank} \, \Sh{V} = (\text{rank} \, \Sh{V}_i)_{i \in Q_0} = \bf{v}$. For each $e \in Q_1$ we have morphisms 
\begin{equation}
\begin{split}
    \Sh{B}_{2, e} : & \Sh{V}_{s(e)} \to \Sh{V}_{t(e)} \\
    \Sh{B}_{3, e} : & \Sh{V}_{t(e)} \to \Sh{V}_{s(e)} 
\end{split}
\end{equation}
which are required to satisfy the usual moment map condition: for $i \in Q_0$
\begin{equation}
    \sum_{e \in t^{-1}(i)} \Sh{B}_{2, e} \Sh{B}_{3, e} - \sum_{e \in s^{-1}(i)} \Sh{B}_{3, e} \Sh{B}_{2, e} + \Sh{I}_i \Sh{J}_i = 0.
\end{equation}
A map $f: \bbP^1\to \FX$ is a quasimap $f: \bbP^1\to X$, i.e. it belongs to $\QM(\bbP^1, X)$ if in the generic fiber of $\oplus_i \Sh{V}_i$, $\ker \Sh{J}$ has no nonzero $\Sh{B}$-invariant subspace. This is a reflection of the stability condition $\zeta^-$ we have chosen for $X$ throughout this paper, cf. \S\ref{subsec: Nakajima quiver varieties}. 

The degree of a the $f$ is the vector 
\begin{equation}
    \deg(f) \coloneqq (\deg \Sh{V}_i)_{i \in Q_0} \in \mathbb{Z}^{Q_0}
\end{equation}
and for our choice of stability condition, the degree lies in the ample cone $\deg \Sh{V}_i \leq 0$ for each $i$. Therefore, we have $\deg(f)=-\dd$ for some $\dd\in \mathbb{Z}_{\geq 0}^{Q_0}$. Later in the paper, we will often suppress $\mathbb{P}^1$ from the notation and write $\QM(X)$. This shall not generate confusion as $\bbP^1$ is the only source curve considered in this paper.

We choose a distinguished point $\infty \in \mathbb{P}^1$. There is an evaluation map $\text{ev}_\infty: \QM(X) \to \FX$
and we will consider the moduli space 
\begin{equation*}
    \QM^\xi(\bbP^1, X) \coloneqq \text{ev}^{-1}_\infty(\xi)
\end{equation*}
for our chosen point $\xi = (\xi_2, \xi_3, \iota, \jmath) \in X \subset \FX$. By construction $\QM^\xi(\bbP^1, X)$ is a closed substack of $\QM(\bbP^1, X)$. It is convenient to choose $\xi \in X^{T}$ for some torus $T$ acting on $X$, and we will often make such a choice; in particular, if the fixed locus of $T$ is isolated we have finitely many preferred choices of $\xi$ if we wish for $T$ to act on $\QM^\xi$. It will also be convenient to work with quasimaps of fixed degree $-\dd$ in addition to evaluating to $p$ at $\infty$; denote the space of such quasimaps by 
\begin{equation}
    \QM^\xi_{\dd}(X) \subset \QM^\xi(X). 
\end{equation}
\subsection{Quasimap spaces as critical loci} 
\label{sec: quiver critical description quasimaps}
Fix a quiver $Q=(Q_0, Q_1)$, two dimension vectors $\vv,\ww\in \BoQ^{Q_0}$, and let $X=\Nak_{Q}(\ww,\vv)$ be the associated quiver variety. We introduce the quiver $Q^{\QM}_{\ww,\vv}=(Q^{\QM}_{\ww,\vv,0}, Q^{\QM}_{\ww,\vv,1})$
\begin{align*}
    Q_{\ww, \vv, 0}
    =&Q_0\sqcup \{\infty\}
    \\
    Q_{\ww,\vv,1}
    =&\ol{Q_1}\sqcup \{\omega_i\;|\;i\in Q_0\}\sqcup \{ \iota_{i,m} \;|\; i\in Q_0\;, 1\leq m\leq \vv_i\}
    \\
    & 
    \sqcup \{ j_{a,m}, j_{a^*,m} \;|\; a\in Q_1\;,1\leq m\leq \vv_i\} 
    \sqcup\{ \alpha_{i,m} \;|\; i\in Q_0\;,1\leq m\leq \ww_i\}
\end{align*}
The orientation of the edges in $\ol{Q_1}$ is inherited from the quiver $\bar Q$ while the one of the new arrows is as prescribed by the following assignments:
\begin{align*}
    &s(\omega_{i})=t(\omega_i)=i\\
    &s(\iota_{i,m})=\infty\quad t(\iota_{i,m})=i\\
    &s(j_{a,m})=t(a)\quad t(j_{i,m})=\infty\\
    &s(j_{a^*,m})=s(a)\quad t(j'_{i,m})=\infty\\
    &s(\alpha_{i,m})=i\quad t(\alpha_{i,m})=\infty
\end{align*}

The connected components of the space of representations for such quiver are labeled by dimension vectors $(\dd,n)\in \BoN^{Q_0}\times \BoN$ ($\BoN = \mathbb{Z}_{\geq 0}$). In particular, the space of linear representations of dimension $(\dd,1)\in \BoN^{Q_{\ww, \vv, 0}}\cong \BoN^{Q_0}\times \BoN$ is given by 
\begin{equation*}
\begin{split}
    \rep_{Q^{\QM}_{\ww,\vv}}(\dd,1) \coloneqq & \bigoplus_{i \in Q_0} \text{Hom}(U_i, U_i) \oplus \text{Hom}(V_i, U_i) \oplus \text{Hom}(U_i, W_i) \\
    & \bigoplus_{e \in Q_1} \text{Hom}(U_{s(e)}, U_{t(e)}) \oplus \text{Hom}(U_{t(e)}, U_{s(e)}) \oplus \text{Hom}(U_{t(e)}, V_{s(e)}) \oplus \text{Hom}(U_{s(e)}, V_{t(e)}). 
\end{split}
\end{equation*}
where $V_i$ and $W_i$ are the fibers of tautological bundles over $X$ at a point $\xi \in X$, and $U_i$ has dimension $d_i$ fixed by the degree, $i \in Q_0$. Hence, a representation $m\in \rep_{Q^{\QM}_{\ww,\vv}}(\dd,1)$ is equivalent to a collection of $\BoC$-linear maps 
\begin{equation}
\label{eq: explicit quiver data quasimaps}
\begin{split}
    B_{1, i} : & U_i \to U_i \\
    I_{i} : & V_i \to U_i \\
    B_{2, e}: & U_{s(e)} \to U_{t(e)} \\
    B_{3, e}: & U_{t(e)} \to U_{s(e)} \\
    J_{2, e}: & U_{t(e)} \to V_{s(e)} \\
    J_{3, e}: & U_{s(e)} \to V_{t(e)} \\
    \alpha_i: & U_i \to W_i.
\end{split}
\end{equation}

\begin{remark}
The vertex at infinity $\infty\in Q^{\QM}_{\ww,\vv,0}$ can be thought as a framing. In this article, we choose choose King's approach to describing framed quiver representations. Namely, instead of adding many framing vertices, one for each vertex $i\in \bar Q_0$ of the original—unframed—quiver, and a single arrow between $i$ and the corresponding framing node, we add a single framing vertex $\infty$ and many arrows to each of the vertices $i\in \bar Q_0$.  The resulting linear spaces of representations are isomorphic.
\end{remark}

In \cite{quasimapcrit}, the following presentation of $\QM^\xi_{\dd}(X)$ is obtained. The value of the quasimap at infinity is given by specifying vector spaces $V = (V_i)_{i \in Q_0}$ with fixed isomorphisms $V_i \simeq \Sh{V}_i \eval_{\infty}$ for each $i \in Q_0$, and fixed quiver data $(\xi_2, \xi_3, \iota, \jmath)$ defining the point $\xi \in X$. 

We have the group $\GL(U) \ltimes \text{Hom}(U, V) \coloneqq \prod_{i \in Q_0} \GL(U_i) \ltimes \Hom(U_i, V_i)$, an element $(g, x)$ of which acts on a representation $m=(B_1,B_2, B_3, I, J_1, J_2, \alpha)\in \rep_{Q^{\QM}_{\ww,\vv}}(\dd,1)$ by
\begin{equation*}
    (g^{-1} B_1g + g^{-1} I x, g^{-1} B_2g, g^{-1} B_3 g, g^{-1} I, J_2g - \xi_3 x + xg^{-1}B_3 g, J_3g + \xi_2 x - x g^{-1} B_2 g, \alpha g + \jmath x).
\end{equation*}
We then define, for each fixed $\xi \in X$ the potential function $W^\xi: \rep_{Q^{\QM}_{\ww,\vv}}(\dd,1)\to \BoC$ as
\begin{multline}
\label{eq: quasimap potential}
    W^\xi \coloneqq \sum_{e \in Q_1} \left( \tr B_{1, t(e)} B_{2, e} B_{3, e} - \tr B_{1, s(e)} B_{3, e} B_{2, e}\right)\\
    + \left(\tr J_{2, e}(B_{2, e} I_{s(e)} - I_{t(e)} \xi_{2, e} ) + \tr J_{3, e} (B_{3, e} I_{t(e)} - I_{s(e)} \xi_{3, e}) \right) \\
    + \sum_{j \in Q_0} \tr I_j \iota_j \alpha_j 
\end{multline}
We call it quasimap potential. It is interesting to check that $W^\xi$ is $GL(U) \ltimes \text{Hom}(U, V)$-invariant. 

Consider now the slope stability $\zeta^+=(0,1)\in \BoQ^{Q_0}\times \BoQ$ and let $\rep_{Q^{\QM}_{\ww,\vv}}^{\zeta^+\sst}(\dd,1)$ be the subspace of $(\dd,1)$-dimensional (slope)-semistable representations of $Q^{\QM}_{\ww,\vv}$. The following is well-known.
\begin{lemma}
\label{lemma: stabiloity quasimaps space}
    A representation $m=(B_1,B_2, B_3, I, J_1, J_2, \alpha)\in \rep_{Q^{\QM}_{\ww,\vv}}(\dd,1)$ is $\zeta^+$-semistable iff $\mathbb{C} \langle B_1, B_2, B_3 \rangle I(V) = U$.
\end{lemma}
We define 
\[
\FM_{\dd}(Q^{\QM}_{\ww,\vv}) \coloneqq \rep_{Q^{\QM}_{\ww,\vv}}(\dd,1)/\GL(U) \ltimes \Hom(U, V). 
\]
and, similarly, the substack of $\zeta^+$-semistable points
\begin{equation}
    \FM_{\dd}^{\zeta^+\sst}(Q^{\QM}_{\ww,\vv}) \coloneqq \rep_{Q^{\QM}_{\ww,\vv}}^{\zeta^+\sst}(\dd,1)/\GL(U) \ltimes \Hom(U, V). 
\end{equation}
Either of the stacks above is a smooth Artin stack of finite type. In fact, the semistable locus $ \FM_{\dd}^{\zeta^+\sst}(Q^{\QM}_{\ww,\vv}$ is representable by a scheme.\footnote{This follows from Proposition 3.8 of \cite{quasimapcrit} and GIT.} 
As the quiver data $(\xi_2, \xi_3, \iota, \jmath)$ are fixed at $\xi \in X$, we may regard $W^\xi$ as a function
\[
W^\xi: \FM_{\dd}(Q^{\QM}_{\ww,\vv}) \to \BoC
\] 
We use the same notation to denote its restriction to the semistable locus: $W^\xi: \FM_{\dd}^{\zeta^+\sst}(Q^{\QM}_{\ww,\vv})  \to \BoC$.  
\begin{theorem}[\cite{quasimapcrit}]
\label{thm: critical description quasimaps}
We have an equivalence of algebraic stacks $\QM^p_{\dd}(X) \simeq \crit(W^\xi)$.
\end{theorem}
Let now $\rep^\circ_{Q^{\QM}}(\dd)$ be the open subset of $\rep_{Q^{\QM}_{\ww,\vv}}^{\zeta\sst}(\dd)$ such that $\BoC\langle B_1 \rangle I(V) = U$ and denote by $\FM^{\circ}_{\dd}(Q^{\QM}_{\ww,\vv})$ its stack quotient by $\GL(U) \ltimes \Hom(U, V)$. By Lemma \ref{lemma: stabiloity quasimaps space}, it follows that $\FM^{\circ}_{\dd}(Q^{\QM}_{\ww,\vv})$ is an open substack of $\FM_{\dd}^{\zeta^+\sst}(Q^{\QM}_{\ww,\vv})$. Although this inclusion is proper for general $\dd\in \BoN^{Q_0}$, we have:
\begin{lemma}
    The inclusion $\FM^{\circ}_{\dd}(Q^{\QM}_{\ww,\vv})\cap \crit(W^\xi)\subseteq \crit(W^\xi)$ is an isomorphism.
\end{lemma}
\begin{proof}
    Let $m\in \rep_{Q^{\QM}_{\ww,\vv}}(\dd)$ be a representation of the Jacobi algebra $J(Q^{\QM}_{\ww,\vv},W^{\xi})$. Then, among the relations satisfied by $m$, we have
    \begin{align*}
        B_{2,e}I_{s(e)}=I_{t(e)}\xi_{2,e}\\
        B_{3,e}I_{s(e)}=I_{t(e)}\xi_{3,e}\\
        B_{2,e}B_{1, t(e)}-B_{1,s(e)}B_{2,e}+I_{s(e)}J_{2,e} = 0\\
        B_{3,e}B_{1, t(e)}-B_{1,s(e)}B_{3,e}+I_{s(e)}J_{2,e} = 0
    \end{align*}
    which imply at once that $\mathbb{C} \langle B_1, B_2, B_3 \rangle I(V) \subseteq \mathbb{C} \langle B_1\rangle I(V) $, and hence the statement.
\end{proof}

\subsection{Main example: Hilbert scheme}
In the very important special case, $X = \text{Hilb}_n(\mathbb{C}^2)$, this construction reads as follows. We fix the evaluation at infinity to be a partition $\lambda \in \text{Hilb}_n(\mathbb{C}^2)$, so that the torus $T \simeq \mathbb{G}_m^2$ acts on $\QM^\lambda(\bbP^1, X)$. In terms of quiver representations, the data $\iota = 0$ identically, and we have $(\xi_2, \xi_3, \jmath)$ where $\xi_{2, 3} : V \to V$ are two commuting linear maps $\comm{\xi_2}{\xi_3} = 0$, and $\jmath: V \to \mathbb{C}$ is a linear function, where $\dim V = n = |\lambda|$.  We can encode all the linear algebraic data in a quiver 
\begin{equation} \label{bigadhmquiver}
\begin{tikzcd}[row sep = huge]
    \underline{V} \arrow[out=45, in=90, loop, swap, "\xi_2"] 
    \arrow[out=90, in=135, loop, swap, "\xi_3"] 
    \arrow[r, "j"]
    \arrow[d, swap, "I"] & \underline{\mathbb{C}} \\
    U
    \arrow[out=157.5, in=202.5, loop, swap, "B_1"]
    \arrow[out=247.5, in=292.5, loop, swap, "B_2"]
    \arrow[out=337.5, in=22.5, loop, swap, "B_3"] 
    \arrow[u, bend left, "J_2"]
    \arrow[u, bend left, shift left=3ex, "J_3"]
    \arrow[ur, "\alpha"]
\end{tikzcd}
\end{equation}
where $\dim V = d$ for quasimaps of degree $d$. The underlined vector spaces are regarded as fixed framing. The potential function simply reads 
\begin{equation}
    W = \tr B_1 \comm{B_2}{B_3} + \tr J_2(B_2 I - I \xi_2) + \tr J_3(B_3 I - I \xi_3). 
\end{equation}

\subsection{Torus action}
\label{subsec: Torus action crit qm}
On the moduli space of quasimaps $\QM^\xi(X)$, there is an evident action of $\mathbb{C}^\times_{t_1} \times T$, where $T$ is a torus acting on $X$ as in \ref{subsec: torus action quiver varieties} and $\mathbb{C}^\times_{t_1}$ acts by rotation of the source curve $\mathbb{P}^1$ fixing $\infty$ and $0 \in \mathbb{P}^1$. We follow our convention from $\ref{subsec: torus action quiver varieties}$ and consider the same abstract torus $T$ to act on $\QM^\xi(X)$ as a Calabi-Yau torus. We have the following elementary result:
\begin{lemma}
    The $\mathbb{C}^\times_{t_1} \times T$ action extends to the ambient representation spaces $\FM_{\dd}(Q^{\QM}_{\ww, \vv})$, and the quasimap potential $W$ is homogeneous of degree $t_1 \hbar$ with respect to $\mathbb{C}^\times_{t_1} \times T$.
\end{lemma}
In particular, $W$ is $T$-invariant for the Calabi-Yau torus introduced in \S\ref{subsec: torus action quiver varieties}. Notice that in checking the above assertion it is important to use the $T$-module structure of the framing spaces $W_i$, $V_i$ in the quasimap quiver induced from the choice of $T$-fixed point $\xi \in X$. 

Following the conventions of \S\ref{subsec: moduli of quiver reps}, we set $\FM^T_{\dd}(Q^{\QM}_{\ww, \vv}):=\FM_{\dd}(Q^{\QM}_{\ww, \vv})/T$ and denote the induced potential by $W^{\xi}:\FM^T_{\dd}(Q^{\QM}_{\ww, \vv}) \to \BoC$. As a byproduct of the previous lemma we deduce that the equivalence of stacks from Theorem \ref{thm: critical description quasimaps} holds equivariantly, i.e. we have $\QM^{\xi}_{\dd}(X)/T \simeq \crit(W^\xi)/T$. Equivalently, we may see $\QM^{\xi}_{\dd}(X)/T$ as the critical locus of the function $W^{\xi}:\FM^{T, \zeta^+\sst}_{\dd}(Q^{\QM}_{\ww, \vv}) \to \BoC$. 

\subsection{Evaluation point and dimensional reduction}
\label{subsec: ev point and dim red}
  Fix a point
 \[
 p=(\xi_2, \xi_3, i, j)\in X \subset \mu^{-1}(0)/G.
 \]
 We say that $\xi$ is \emph{polarized} if there exists a sub-quiver $Q'\subset \bar Q$ such that $\ol{Q'}=\bar Q$ and $\xi$ belongs to the Lagrangian subvariety $L= \Rep_{Q'_{\ff}}(\dd,1)/G\cap X$. Clearly, up to changing the orientations of some arrows of $Q$, we may always assume that $Q'=Q$.  In coordinates on the prequotient $\Rep_{\bar{Q'_{\ff}}(\dd,1)}=T^*\Rep_{Q'_{\ff}(\dd,1)}$, a polarized fixed point can always be written as $(\xi_2, 0,0,\jmath)$. Notice that the vanishing $\iota=0$ as opposed to $\jmath =0$ is forced by our choice of stability condition on the quiver variety. 

 \begin{remark}
     The fiber at $0\in X_0$
     of the resolution $\pi: X\to X_0$ is a (generally singular and reducible) Lagrangian subvariety of $X$ and any $L= \Rep_{Q'_{\ff}}(\dd,1)/G\cap X$ as above is an irreducible component in it. However, not all irreducible components are of this form. For example, consider $X=\Hilb_3(\BoC^2)$. Then $\pi^{-1}(0)$ admits three connected components, each containing one fixed point of $X$. It is easy to check two out of these three fixed points are polarized, but the one corresponding to the partition $3=2+1$ is not. More generally, for $\Hilb_n(\BoC^2)$, the only polarized fixed points are those corresponding to the partitions whose Young diagram is either just a row or just a column.
 \end{remark}
 
 Given a polarized $\xi \in X $, set $N=\Rep_{Q'_{\ff}}(\dd,1)/G$. Let $\QM(\bbP^1,L)\subset \Map(\bbP^1, N/G)$ be the moduli space of quasimpas to $L\subset N/G$. If $\xi \in X $ is polarized, then the critical cohomology of $\QM(\bbP^1, X)$ can be dimensionally reduced. To prove this statement, it is convenient to introduce a quiver model for $\QM(\bbP^1,L)$.  Let $\hat Q_{\ww,\vv}=(\hat Q_{\ww,\vv,0}, \hat Q_{\ww,\vv,1})$ be the quiver
 \begin{align*}
 	\hat Q_{\ww, \vv, 0}
 	=&Q_0\sqcup \{\infty\}
 	\\
 	\hat Q_{\ww,\vv,1}
 	=&Q_1\sqcup \{\omega_i\;|\;i\in Q_0\}\sqcup \{ \iota_{i,m} \;|\; i\in Q_0\;, 1\leq m\leq \vv_i\}
 	\\
 	& 
 	\sqcup \{ j_{a,m} \;|\; a\in Q_1\;,1\leq m\leq \vv_i\} 
 	\sqcup\{ \alpha_{i,m} \;|\; i\in Q_0\;,1\leq m\leq \ww_i\}
 \end{align*}
 Let $\Rep_{\Xi_{Q_{\ww,\vv}}}(\dd,1)$ be the space of $(\dd,1)$-dimensional representations of the quiver algebra $\Xi_{Q_{\ww,\vv}}$, cf. \S\ref{subsec: Doubling, looping, tripling, and framing}. Explicitly, this is the subvariety of $\Rep_{\hat Q_{\ww,\vv}}(\dd,1)$ given by the zero locus of the equations
 \begin{align}
 	\label{eq: quasimap to lagrangian equations}
 \begin{split}
 	&B_{1, t(e)} B_{2, e} -  B_{2, e}B_{1, s(e)} +  I_{t(e)} J_{3, e}=0
 	\\
 	&B_{2, e} I_{s(e)} - I_{t(e)} \xi_{2, e}.
 \end{split}
 \end{align}
 for all $e\in Q_1$. We further define
 \[
 \FM_{\dd}(\Xi_{Q_{\ww,\vv}}) \coloneqq\Rep_{\Xi_{Q_{\ww,\vv}}}(\dd,1)/\GL(U) \ltimes \Hom(U, V). 
 \]
 and, similarly,
 \begin{equation}
 	\FM^{\zeta^+\sst}_{\dd}(\Xi_{Q_{\ww,\vv}}) \coloneqq\Rep^{\zeta^+\sst}_{\Xi_{Q_{\ww,\vv}}}(\dd,1)/\GL(U) \ltimes \Hom(U, V). 
 \end{equation}
 Either of the stacks above is a quasi-smooth Artin stack of finite type. Moreover, the semistable locus $\FM_{\dd}(\Xi_{Q_{\ww,\vv}})_{\zeta^+\sst}(\mathbf{d})$ is representable by a scheme.  The following result follows from a straightforward adaptation of the argument in \cite{quasimapcrit}.
 \begin{proposition}
 	\label{prop: dimensionally reduced quasimap quiver model}
 	Assume that $\xi$ is polarized. Then we have an equivalence of stacks 
 	\[
 	\FM^{\zeta^+\sst}_{\dd}(\Xi_{Q_{\ww,\vv}})\simeq \QM^\xi_{\dd}(\bbP^1,L).
 	\]
 \end{proposition}
The action of the torus $T$ on $\FM(Q^{\QM_{\ww,\vv}})$ and $\QM(\bbP^1, X)$ discussed in \S\ref{subsec: Torus action crit qm} restricts to an action on $\FM(\Xi_{Q_{\ww,\vv}})$ and $\QM(\bbP^1, L)$. We denote by $\FM^T(\Xi_{Q_{\ww,\vv}})$ and $\QM^{\xi}(\bbP^1, L)/T$ the associated quotient stacks. 

We can now relate the Borel-Moore homology of $\QM^{\xi}(\bbP^1,L)$ with the critical cohomology of $\QM^{\xi}_{\dd}(\bbP^1, X)$. By Proposition \ref{prop: dimensionally reduced quasimap quiver model} we deduce that the virtual dimension of $\QM^{\xi}(\bbP^1,L)$ is $\dd^T(\Q^T -\Q)\vv+\dd^T\ww$. Accordingly, set 
 \[
 \BoQ^{\vir}_{\QM^{\xi}(\bbP^1, L)/T}\coloneqq\BoQ_{\QM(\bbP^1, L)/T}[\dd^T(\Q^T -\Q)\vv+\dd^T\ww)].
 \]
\begin{proposition}
\label{prop: dimred Beilison model}
    Assume that $\xi$ is polarized. Then we have a canonical isomorphism
    \[
    \HO_T(\QM^{\xi}_{\dd}(\bbP^1, X), \varphi_{\QM^{\xi}_{\dd}(\bbP^1, L)})\cong \HO^{\BM}_T(\QM^{\xi}_{\dd}(\bbP^1, L), \BoQ^{\vir}).
    \]
\end{proposition}
\begin{proof}
    We drop the torus $T$ from the notation as the argument is insensitive of it. Up to changing the orientation of $Q$, we may assume that $\xi =(\xi_2, 0,0, \jmath)$. Therefore, the quasimap potential \eqref{eq: quasimap potential} takes the form 
    \begin{multline*}
    \sum_{e \in Q_1} \left( \tr B_{1, t(e)} B_{2, e} B_{3, e} - \tr B_{1, s(e)} B_{3, e} B_{2, e}\right)
    + \left(\tr J_{2, e}(B_{2, e} I_{s(e)} - I_{t(e)} \xi_{2, e} ) + \tr J_{3, e} (B_{3, e} I_{t(e)}) \right)
    \\
    =\sum_{e \in Q_1} \tr  \left(B_{3, e} \left(B_{1, t(e)} B_{2, e} -  B_{2, e}B_{1, s(e)} +  I_{t(e)} J_{3, e}\right)\right)
    + \tr\left( J_{2, e}\left(B_{2, e} I_{s(e)} - I_{t(e)} \xi_{2, e} \right) \right)
\end{multline*}
In particular, it is homogeneous in $B_{3, e}$ and $J_{2,e}$. Moreover, these fields are exactly the Lagrange multipliers of \eqref{eq: quasimap to lagrangian equations}. Set 
\[
\FM^{\zeta^+\sst}_{\dd}(\hat Q_{\ww,\vv}) \coloneqq\Rep^{\zeta^+\sst}_{\hat Q_{\ww,\vv}}(\dd,1)/\GL(U) \ltimes \Hom(U, V). 
\]
and let $E\to \FM^{\zeta^+\sst}_{\dd}(\Xi_{Q_{\ww,\vv}})$ be the restriction of the vector bundle 
\[
\FM^{\zeta^+\sst}_{\dd}(Q^{\QM}_{\ww,\vv})\to \FM^{\zeta^+\sst}_{\dd}(\hat Q_{\ww,\vv})
\]
to $\FM^{\zeta^+\sst}_{\dd}(\Xi_{Q_{\ww,\vv}})\subset \FM^{\zeta^+\sst}_{\dd}(\hat Q_{\ww,\vv})$.
By dimensional reduction \ref{thm: dimred}, there is a canonical isomorphism 
 \begin{equation}
 	\label{eq: in the proof dimred quasimap}
    \HO_T(\QM^\xi_{\dd}(\bbP^1, X), \varphi_{\QM^\xi_{\dd}(\bbP^1, X)})\cong \HO^{\BM}(E, \BoQ[\dim(\FM_{{Q}^{\QM}_{\ww,\vv}})]),
\end{equation}
Therefore, combining  \eqref{eq: in the proof dimred quasimap} and Proposition \ref{prop: dimensionally reduced quasimap quiver model}, we get
\[
\HO_T(\QM^\xi_{\dd}(\bbP^1, X), \varphi_{\QM^\xi_{\dd}(\bbP^1, X)})\cong \HO^{\BM}(\QM^\xi_{\dd}(\bbP^1, L), \BoQ[\dim(\FM_{{Q}^{\QM}_{\ww,\vv}})-2\rk(E)])
\]
The statement now follows from a direct computation of the cohomological shift.
\end{proof}

\section{Critical Hall action on quasimaps via Hecke modifications}

\subsection{Hecke correspondence on quasimap spaces}
\label{subsec: Hecke correspondence on quasimap space}
Let $X$ be a Nakajima quiver variety, and for a given fixed point $\xi\in X$, let $\QM^{\xi}_{\dd}(X)$ be the moduli space of degree ${\dd}$ quasimaps with fixed evaluation ${\ev}_{\infty}=\xi$ with our notations and terminology as usual.
\begin{lemma}
\label{lemma negative vector bundle}
    Let $(\Sh{V},s)\in \QM^{\xi}(X)$ be a quasimap datum. Then every sub-bundle of $\Sh{V}$ has non-positive degree.
\end{lemma}
\begin{proof}
    Fix $i\in Q_0$ and consider a splitting $\Sh{V}_i=\Sh{O}(d_1)\oplus \dots \Sh{O}(d_{\vv_i})$. By quasimap stability for the chosen slope $\zeta$ and Lemma \ref{lemma characterization anticanonical stability}, there exists a non-zero map from $\Sh{O}(d_l)$ to the framing $\Sh{O}^{\oplus \ww}$. This forces $d_l\leq 0$.
\end{proof}
\begin{corollary}
    For any quasimap datum, we have $\HO^0(\bbP^1, \Sh{V}(-1))=0$.
\end{corollary}
Consider now the moduli spaces
$$
\Hecke^{\xi}(X)=\lbrace (\Sh{V},\Sh{V}',s, s') \; | \; \Sh{V}' \subset \Sh{V} + \text{compatibility of sections} \rbrace /\sim
$$
The notation is as follows. $\Sh{V} = (\Sh{V}_i)_{i \in Q_0}$ is shorthand for a tuple of bundles on $\mathbb{P}^1$, $s =  (\mathcal{B}_2, \mathcal{B}_3, \mathcal{I}, \mathcal{J})$ is shorthand for the sections entering the quasimap data; we insist that $(\Sh{V}, s)$ and $(\Sh{V}', s')$ define stable quasimaps to $X$ such that $s(\infty) = s'(\infty) = \xi$. Compatibility of sections we mean that for any edge $a\in Q_1$ and $i\in Q_0$, we have commutative diagrams
\[
\begin{tikzcd}
    \Sh{V}'_{s(a)} \arrow[r, hook]\arrow[d, "\mathcal{B}'_{2, a}"] & \Sh{V}_{s(a)}\arrow[d, "\mathcal{B}_{2, a}"]\\
    \Sh{V}'_{t(a)} \arrow[r, hook] & \Sh{V}_{t(a)} 
\end{tikzcd}
\qquad \qquad 
\begin{tikzcd}
     \Sh{V}'_{i} \arrow[r, hook]\arrow[d, "\mathcal{J}'_i"] & \Sh{V}_i\arrow[d, "\mathcal{J}_i"]\\
    W_i \otimes \Sh{O}_{\mathbb{P}^1} \arrow[r, equal] & W_i \otimes \Sh{O}_{\mathbb{P}^1}
\end{tikzcd}
\]
and similar ones for $\mathcal{B}_{3, a}$ and $\mathcal{I}_i$. We denote by $\Hecke^{\xi}_{\dd, \dd'}(X)$ the component of $\Hecke^{\xi}(X)$ such that $\deg{\Sh{V}}= -\dd$ and $\deg{\Sh{V}'}= -\dd-\dd'$.
It comes with natural forgetful morphisms 
\[
p: \Hecke^{\xi}_{\dd, \dd'}(X) \to \QM^{\xi}_{\dd}(X)\qquad p': \Hecke^{\xi}_{\dd,\dd'}(X) \to \QM^{\xi}_{\dd + \dd'}(X)
\]
Any point $(\Sh{V}, \Sh{V'}, s,s')\in \Hecke^{\xi}_{\dd,\dd'}(X)$ induces the following exact sequence
\begin{equation}
    \label{eq: diagram ses sub quotient general discussion}
        \begin{tikzcd} 0\arrow[r]&\Sh{V}'_{s(a)}\arrow[r]\arrow[d, "s'_a"]& \Sh{V}_{s(a)}\arrow[r]\arrow[d, "s_a"]&\Sh{V}_{s(a)}/\Sh{V'}_{s(a)}\arrow[d, "\beta_a"] \arrow[r]&0
    \\
    0\arrow[r]& \Sh{V}'_{t(a)}\arrow[r]& \Sh{V}_{t(a)}\arrow[r]&    \Sh{V}_{t(a)}/\Sh{V'}_{t(a)}\arrow[r] & 0
    \end{tikzcd}
\end{equation}
for all $a\in Q_1$. Set $Q_i=: \Sh{V}_i/\Sh{V}'_i$. These are torsion sheaves on $\bbP^1$ supported at $\AA^1=\bbP^1\setminus \{\infty\}$. Applying the functor $\HO^0(\bbP^1,-)$ to the pair $(Q, \beta)$ gives a representation of the tripled quiver $\tilde Q$, where the action a loop $\omega_i\in \tilde Q_1$ is induced by the action of a fixed generator of $\BoC[\AA^1]$ on $Q_i$ and the action of $a$ and $a^*$ is induced by the maps $\beta_a$ and $\beta_{a^*}$. We denote such representation by $M$. By additivity of degree in exact sequences, $M$ is a $\dd'$-dimensional representation.
\begin{lemma}
    The representation $M$ satisfies the relation of the Jacobi algebra $\Jac(\tilde Q, \tilde W)$ for the tripled quiver $\tilde Q$ with its canonical cubic potential $\tilde W$. 
\end{lemma}
\begin{proof}. The same argument of \S\ref{subsec: Q decorated coherent sheaves on A1} implies that the relations 
\[
    \sum_{a\in Q_1, i\in Q_0} a\omega_i -\omega_i a\qquad \sum_{a\in Q_1, i\in Q_0} a^*\omega_i -\omega_i a^*
\]
are satisfied. Therefore, it suffices to show that $\gamma$ satisfies the moment map relations $\sum_{a\in Q_1} [a, a^*]$ are satisfied too. By Lemma \ref{lemma negative vector bundle}, we have a commutative diagram
\begin{equation*}
    \begin{tikzcd}
        0\arrow[r]& \HO^0(\bbP^1, \Sh{V}_{s(a)}/\Sh{V'}_{s(a)})\arrow[r]\arrow[d, "\beta_a"]& \HO^1(\bbP^1, \Sh{V}'_{s(a)}(-1)) \arrow[r]\arrow[d, "\HO^1(s'_a)"]& \HO^1(\bbP^1, \Sh{V}_{s(a)}(-1))\arrow[r]\arrow[d, "\HO^1(s_a)"]& 0
        \\
        0\arrow[r]& \HO^0(\bbP^1, \Sh{V}_{t(a)}/\Sh{V'}_{t(a)})\arrow[r]& \HO^1(\bbP^1, \Sh{V}'_{t(a)}(-1)) \arrow[r]& \HO^1(\bbP^1, \Sh{V}_{t(a)}(-1))\arrow[r]& 0
    \end{tikzcd}
\end{equation*}
for any $a\in \bar Q_1$. Since the sections $s'_a$ and $s_e$ satisfy the moment map condition fiberwise, the pair $(\HO^1(s'_a), \HO^1(s'_{a^*}))_{a\in Q_1}$ satisfies the moment map equation. Therefore, the same must hold for $(\beta_a, \beta_{a^*})_{a\in Q_1}$. The proof follows.
\end{proof}
The lemma implies that we have a morphism 
\[
q: \Hecke^{\xi}_{\dd, \dd'}(X)\to \FM_{\dd'}(\Jac(\tilde Q,\tilde W))=\crit_{\dd'}(\tilde W).
\]
and hence a correspondence 
\begin{equation}
\label{eq: Hecke correspondence quasiaps}
\begin{tikzcd}
    \crit_{\dd'}(\tilde W)\times\QM^{\xi}_{\dd}(X)  & \Hecke^{\xi}_{\dd, \dd'}(X)\arrow[l, swap, "q\times p"]\arrow[r, "p'"] & \QM^{\xi}_{\dd+\dd'}(X)
\end{tikzcd}
\end{equation}
It is clear that all the arguments of this section hold equivariantly for the action of the torus $T$ on $\QM^{\xi}(X)$ defined in \ref{subsec: Torus action crit qm}. It will now be explained how this diagram induces a Hall action in the setting of critical cohomology.


\subsection{Critical description of Hecke correspondences}
\label{subsec: Critical description of Hecke correspondences}
In this section, we relate the Hecke correspondence \eqref{eq: Hecke correspondence quasiaps} to cohomological Hall induction for the moduli space $\FM_{\dd}^{\zeta^+\sst}(Q^{\QM}_{\ww,\vv})$. It is easy to check that, for every $\dd\in \BoN^{Q_0}$ the connected components of the stack of graded objects
\[
\Grad(\FM_{\dd}(Q^{\QM}_{\ww,\vv}))
\]
are of the form
\[
\prod_{i=1}^{n-1}\FM_{\dd_i}(\tilde Q) \times \FM_{\dd_n}(Q^{\QM}_{\ww,\vv}),
\]
where $\dd=\dd_1+\dots+\dd_n$ is an arbitrary decomposition of $\dd$. In the case $n=2$, we denote by $\FM_{\dd, \dd'}(Q^{\QM}_{\ww,\vv})$ the preimage of the component $\FM_{\dd'}(\tilde Q) \times \FM_{\dd}(Q^{\QM}_{\ww,\vv})$ of $\Grad(\FM_{\dd + \dd'}(Q^{\QM}_{\ww,\vv}))$ under the canonical map 
\[
\Filt(\FM_{\dd + \dd'}(Q^{\QM}_{\ww,\vv}))\to \Grad(\FM_{\dd+\dd'}(Q^{\QM}_{\ww,\vv})).
\]
Explicitly, $\FM_{\dd, \dd'}(Q^{\QM}_{\ww,\vv})$ is the stack parameterizing short exact sequences of $Q^{\QM}_{\ww,\vv}$-modules
\[
\FM_{\dd , \dd'}(Q^{\QM}_{\ww,\vv})=\{0\to M'\to M'' \to M\to 0\}
\]
where $M'$ is of dimension $(\dd',0)$ and $M''$ is of dimension $(\dd + \dd',1)$. Note that a $Q^{\QM}_{\ww, \vv}$-module, forgetting the quiver arrows, is a short exact sequence for each node $i \in Q_0$
\[
\begin{tikzcd}
    0 \arrow[r] & V_i \arrow[r] & F_i \arrow[r] & U_i \arrow[r] & 0. 
\end{tikzcd}
\]
Then $M''$, $M$ in the short exact sequence of $Q^{\QM}_{\ww, \vv}$-modules above consists of a pair $(F_i'', U_i'')$, $(F_i, U_i)$ fitting into a diagram 
\[ 
\begin{tikzcd}
    0 \arrow[r] & V_i \arrow[d, equal] \arrow[r] & F_i'' \arrow[d, two heads] \arrow[r] & U_i'' \arrow[d, two heads] \arrow[r] & 0 \\ 
    0 \arrow[r] & V_i \arrow[r] & F_i \arrow[r] & U_i \arrow[r] & 0. 
\end{tikzcd}
\]
The canonical projection $F_i'' \to U_i''$ induces an isomorphism $\ker(F_i'' \to F_i) \simeq \ker(U_i'' \to U_i) =: E_i$, so the short exact sequence takes the schematic form 
\[
\begin{tikzcd}
    0 \arrow[r] & E \arrow[r] & (F'', U'') \arrow[r] & (F, U) \arrow[r] & 0
\end{tikzcd}
\] 
where $E$ is a module over the unframed tripled quiver $\tilde{Q}$ and the other terms are $Q^{\QM}_{\ww, \vv}$-modules. The quiver arrows must commute with the inclusions and projections in the obvious way. Fix a cocharacter $\sigma: \mathbb{C}^\times \to \GL_{\dd + \dd'}$ corresponding to the splitting 
\[
U_{\dd + \dd'} = U_{\dd} \oplus U_{\dd'}
\]
scaling the first summand with weight $0$ and the second with weight $+1$. $\sigma$ defines a $\mathbb{C}^\times$-action on $\GL_{\dd + \dd'}$ by conjugation; let $P_{\dd, \dd'} \hookrightarrow \GL_{\dd + \dd'}$ be the parabolic subgroup defined by taking attractors. Let $\rep^{\geq 0}_{Q^{\QM}_{\ww,\vv}}(\dd + \dd', 1) \hookrightarrow \rep_{Q^{\QM}_{\ww, \vv}}(\dd + \dd', 1)$ denote the subspace of nonnegative weight under $\sigma$. We have the elementary 
\begin{lemma} \label{lemma: quotient description of induction space}
    We have an isomorphism of stacks 
    \[
    \FM_{\dd ,\dd'}(Q^{\QM}_{\ww,\vv}) \simeq \rep^{\geq 0}_{Q^{\QM}_{\ww,\vv}}(\dd + \dd', 1)/P_{\dd, \dd'} \ltimes \Hom(U_{\dd}, V).
    \]
    Here we regard $\Hom(U_{\dd}, V)$ as an additive group which acquires a $P_{\dd, \dd'}$-action via pullback under the canonical maps $P_{\dd, \dd'} \to \GL_{\dd} \times \GL_{\dd'} \to \GL_{\dd}$. 
\end{lemma}
\begin{proof}
    Translate the above discussion to matrices. 
\end{proof}

Forgetting the various terms of the sequence, we get maps
\begin{equation}
    \label{eq: hall action correspondence ambient spaces}
    \begin{tikzcd}
    \FM_{\dd'}(\tilde Q) \times \FM_{\dd}(Q^{\QM}_{\ww,\vv}) & \FM_{\dd ,\dd'}(Q^{\QM}_{\ww,\vv})\arrow[l, swap,  "q"]\arrow[r, "p"] & \FM_{\dd + \dd')}(Q^{\QM}_{\ww,\vv}).
\end{tikzcd}
\end{equation}
The maps $p$ and $q$ are, respectively, a proper map and a vector bundle stack. 

\begin{remark}
    In the standard setting of quiver representations, it is a routine argument to verify that $p$ is proper, by factoring it as a closed inclusion followed by a $G/P$ fibration. However, since the group $GL(U) \ltimes \text{Hom}(U, V)$ we quotient by is not reductive, this argument does not a priori apply. In the end, one can see that $p$ is nonetheless proper e.g. by using Proposition 3.8 of \cite{quasimapcrit} to reduce to a situation where the standard argument applies, or using a description of $p$ as in Lemma \ref{lemma: projective bundle description of Nakajima operators} below. 
\end{remark}

Additionally, because of our choice of stability condition $\zeta^+$, there is a well-defined commutative diagram 
\begin{equation}
    \label{eq: hall action diamgram ambient spaces}
    \begin{tikzcd}
    \FM_{\dd'}(\tilde Q) \times \FM_{\dd}(Q^{\QM}_{\ww,\vv}) & \FM_{\dd,\dd'}(Q^{\QM}_{\ww,\vv})\arrow[l, swap,  "q"]\arrow[r, "p"] & \FM_{\dd + \dd'}(Q^{\QM}_{\ww,\vv})
    \\
    \FM_{\dd'}(\tilde Q) \times \FM^{\zeta^+ \sst}_{\dd}(Q^{\QM}_{\ww,\vv}) \arrow[u, hookrightarrow]& \FM^{\zeta^+\sst}_{\dd,\dd'}(Q^{\QM}_{\ww,\vv})\arrow[l, swap, hookrightarrow,  "q^{\zeta^+}"]\arrow[r, "p^{\zeta^+}"]\arrow[u, hookrightarrow] & \FM^{\zeta^+\sst}_{\dd + \dd'}(Q^{\QM}_{\ww,\vv})\arrow[u, hookrightarrow]
\end{tikzcd}
\end{equation}
where the right square is Cartesian. It follows that $p^{\zeta^+}$ is proper and $q^{\zeta^+}$ is smooth. Denote the function $W^\xi_{\dd + \dd'}$ induced by $W$ on the target $\FM^{\zeta^+\sst}_{\dd + \dd'}(Q^{\QM}_{\ww,\vv})$ of the map $p^{\zeta^+}$. 

\begin{theorem} 
\label{thm: hecke as crit}
    We have an equivalence of stacks 
    \[
    (p^{\zeta^+})^{-1}(\Crit(W^\xi_{\dd + \dd'})) \simeq \Hecke^{\xi}_{\dd, \dd'}(X).
    \]
\end{theorem}

This will follow from Propositions \ref{prop: hecke commutative diagram} and \ref{thm: hecke as crit} below. Let $(\Sh{V}, s)$ and $(\Sh{V}', s')$ be a pair of quasimaps related by $\Hecke^{\xi}_{\dd, \dd'}(X)$. By definition, this means we have short exact sequences of sheaves on $\mathbb{P}^1$ 
\[
\begin{tikzcd}
    0 \arrow[r] & \Sh{V}'_i \arrow[r, "g_i"] & \Sh{V}_i \arrow[r] & \text{torsion} \arrow[r] & 0
\end{tikzcd}
\]
where the support of the torsion is disjoint from $\infty \in \mathbb{P}^1$. Applying Lemma 2.2 from \cite{quasimapcrit} to the morphism $g_i$ we get an induced commutative diagram 
\[
\begin{tikzcd}
    0 \arrow[r] & \Sh{V}'_i \arrow[d, "g_i"] \arrow[r] & \HO^1(\mathbb{P}^1, \Sh{V}'_i(-2)) \otimes \Sh{O}_{\mathbb{P}^1} \arrow[d, two heads] \arrow[r] & \HO^1(\mathbb{P}^1, \Sh{V}'_i(-1)) \otimes \Sh{O}_{\mathbb{P}^1}(1) \arrow[d, two heads] \arrow[r] & 0 \\ 
    0 \arrow[r] & \Sh{V}_i \arrow[r] & \HO^1(\mathbb{P}^1, \Sh{V}_i(-2)) \otimes \Sh{O}_{\mathbb{P}^1} \arrow[r] & \HO^1(\mathbb{P}^1, \Sh{V}_i(-1)) \otimes \Sh{O}_{\mathbb{P}^1}(1) \arrow[r] & 0
\end{tikzcd}
\]
exact on the top and bottom rows. The surjectivity of the second and third downward arrows follows from the fact that the support of the cokernel of $g_i$ is zero-dimensional. Taking the fiber over $\infty \in \mathbb{P}^1$ and using the definition of the isomorphism $\Psi$ from Theorem \ref{thm: critical description quasimaps}, we get the data $(F'', U'') \twoheadrightarrow (F, U)$ from the sought-after short exact sequence of $Q^{\QM}_{\ww, \vv}$-modules above. It remains to check the compatibility of the arrows. From the definition of $\Hecke^{\xi}_{\dd, \dd'}(X)$, we have commutative diagrams for each edge $e \in Q_1$
\[ 
\begin{tikzcd}
    \Sh{V}'_{s(e)} \arrow[r, "\Sh{B}'_{2, e}"] \arrow[d, "g_i"] & \Sh{V}'_{t(e)} \arrow[d, "g_i"] \\ 
    \Sh{V}_{s(e)} \arrow[r, "\Sh{B}_{2, e}"] & \Sh{V}_{t(e)}.
\end{tikzcd}
\]
Applying the functor $\HO^1(\mathbb{P}^1, - \otimes \Sh{O}_{\mathbb{P}^1}(-2))$ to this square and the defintion of $\Psi$ gives the compatibility of the quiver arrows $(B_{2, e}, J_{3, e})$ with $(F'', U'') \twoheadrightarrow (F, U)$. A similar argument works for $(B_{3, e}, J_{2, e})$. 

Likewise, from the definition of $\Hecke^{\xi}_{\dd, \dd'}(X)$ we have commutative diagrams
\[ 
\begin{tikzcd}
    \Sh{V}'_i \arrow[r, "\Sh{J}'_i"] \arrow[d, "g_i"] & W_i \otimes \Sh{O}_{\mathbb{P}^1} \arrow[d, equal] \\ 
    \Sh{V}_i \arrow[r, "\Sh{J}_i"] & W_i \otimes \Sh{O}_{\mathbb{P}^1}.
\end{tikzcd}
\]
Applying $\HO^1(\mathbb{P}^1, - \otimes \Sh{O}_{\mathbb{P}^1}(-2))$ gives the compatibility of the quiver data $\alpha_i$. Finally let $z$ denote the operation of multiplication by the affine coordinate on $\mathbb{P}^1$, regarded as a global section of $\Sh{O}_{\mathbb{P}^1}(1)$. There is an obvious commutative square 
\[ 
\begin{tikzcd}
    \Sh{V}'_i(-2) \arrow[r, "z"] \arrow[d, "g_i"] & \Sh{V}'_i(-1) \arrow[d, "g_i"] \\
    \Sh{V}_i(-2) \arrow[r, "z"] & \Sh{V}_i(-1). 
\end{tikzcd}
\]
Applying the functor $\HO^1(\mathbb{P}^1, -)$ to this square gives the compatibility of the quiver data $(B_{1, i}, I_i)$. 

Note this construction works in the universal family as in \cite{quasimapcrit}. This discussion may be summarized as follows. Let $\Psi: \QM^\xi_{\dd}(X) \to \Crit(W^\xi_{\dd})$ denote the isomorphism of Theorem \ref{thm: critical description quasimaps}.

\begin{proposition}
\label{prop: hecke commutative diagram}
We have a commutative diagram 
\begin{equation} \label{eq: hecke as crit}
    \begin{tikzcd}
        \Hecke^{\xi}_{\dd, \dd'}(X) \arrow[r, "p \times p' "] \arrow[d] & \QM^\xi_{\dd}(X) \times \QM^\xi_{\dd  + \dd'}(X) \arrow[d, "\Psi \times \Psi"] \\ (p^{\zeta^+})^{-1}(\Crit(W^\xi_{\dd + \dd'})) \arrow[r] & \Crit(W^\xi_{\dd}) \times \Crit(W^\xi_{\dd + \dd'}).
    \end{tikzcd}
\end{equation}
Moreover, this diagram is Cartesian.
\end{proposition}

In other words, $\Psi$ pulls back to the claimed isomorphism in Theorem \ref{thm: hecke as crit}. 

\begin{proof}
    The first statement follows from the discussion above, so it suffices to show that the diagram is Cartesian. We need to check that if a pair of quiver representations $(m, m'')$ in $\Crit(W^\xi_{\dd}) \times \Crit(W^\xi_{\dd + \dd'})$ are in the image of $(p^{\zeta^+})^{-1}(\Crit(W^\xi_{\dd + \dd'}))$, then the corresponding quasimaps $(\Sh{V}, s; \Sh{V}', s')$ are in the image of $\Hecke^{\xi}_{\dd, \dd'}(X)$. 

    Given $(m, m'')$ in the image of $(p^{\zeta^+})^{-1}(\Crit(W^\xi_{\dd + \dd'}))$, we have a surjection $(F'', U'') \twoheadrightarrow (F, U)$ of $Q^{\QM}_{\ww, \vv}$-modules. By definition of $\Psi$, the vector bundles $\Sh{V}' = (\Sh{V}'_i)_{i \in Q_0}$ fit into short exact sequences 
    \[ 
    \begin{tikzcd}
    0 \arrow[r] & \Sh{V}'_i \arrow[r] & F_i'' \otimes \Sh{O}_{\mathbb{P}^1} \arrow[r, "(B_{1, i}'' - z \, \, \, I''_i)"] & U_i'' \otimes \Sh{O}_{\mathbb{P}^1}(1) \arrow[r] & 0.
    \end{tikzcd}
    \]
    The compatibility of  quiver data $(B_{1, i}'', I_i''; B_{1, i}, I_i)$ implies that the morphisms $g_i: \Sh{V}'_i \to \Sh{V}_i$ defined by
    \[ 
    \begin{tikzcd}
        \Sh{V}'_i \arrow[d, "g_i"] \arrow[r, hook] & F_i'' \otimes \Sh{O}_{\mathbb{P}^1} \arrow[d, two heads] \\ 
        \Sh{V}_i \arrow[r, hook] & F_i \otimes \Sh{O}_{\mathbb{P}^1}
    \end{tikzcd}
    \]
    where the downward arrow is part of the induced surjection $(F'', U'') \twoheadrightarrow (F, U)$ of $Q^{\QM}_{\ww, \vv}$-modules, are well-defined. The compatibility of the rest of the quiver data with $(F'', U'') \twoheadrightarrow (F, U)$ implies the sections $s, s'$ entering the quasimap data fit into the required commutative diagrams with $g$, by simply using the definition of $\Psi$ and Lemma 2.2 of \cite{quasimapcrit}. 

    Then what must be shown is that $g_i$ defined as above is an injective sheaf homomorphism with torsion cokernel, which is equivalent to showing that $g_i$ induces a fiberwise isomorphism on the complement of a finite number of points in $\mathbb{P}^1$. Restricting to the fiber over $z \in \mathbb{P}^1$, this becomes a question of linear algebra. Writing $F_i'' \simeq U_i'' \oplus V_i$, $F_i \simeq U_i \oplus V_i$ as vector spaces, $U''_i \simeq U_i \oplus E_i$ (where $E_i \coloneqq \ker(U_i'' \to U_i) \simeq \ker(F_i'' \to F_i)$), by Lemma \ref{lemma: quotient description of induction space} we may assume 
    \begin{equation}
    \begin{split}
        B_{1, i}'' & = \begin{pmatrix} B_{1, i} && 0 \\ \rho_{1, i} && \gamma_{1, i} \end{pmatrix} \\
        I_i'' & = \begin{pmatrix} I_i \\ I'_i  \end{pmatrix}
    \end{split}
    \end{equation}
    with respect to this decomposition. The equation $(B_{1, i}'' - z) \psi''_i + I_i'' \chi_i = 0$ cutting out the fiber $\Sh{V}'_i |_z \hookrightarrow U_i'' \oplus V_i \ni (\psi''_i, \chi_i)$ becomes equivalent to 
    \begin{equation}
    \begin{pmatrix} (B_{1, i} - z) \psi_i + I_i \chi_i \\ \rho_{1, i} \psi_i + (\gamma_{1, i} - z) s_i + I'_i \chi_i \end{pmatrix} = 0
    \end{equation}
    where $\psi''_i = (\psi_i, s_i) \in U_i \oplus E_i$. $g_i |_z : \Sh{V}'_i |_z \to \Sh{V}_i |_z$ sends $(\psi_i, s_i,  \chi_i)$ solving these equations to $(\psi_i, \chi_i)$ solving the first equation. It is clear that if $z$ is not an eigenvalue of $\gamma_{1, i}$, a solution to the first equation uniquely determines a solution to the second equation, showing that $g_i |_z$ is an isomorphism for $z \in \mathbb{P}^1 \setminus \Spec \gamma_{1, i}$. 
\end{proof}

\begin{remark}
\label{remark: support hecke modifications}
Note a corollary of the proof is that the eigenvalues of the matrices $\gamma_{1, i}$ defining (part of) the $\widetilde{Q}$-module structure on $E = (E_i)_{i \in Q_0}$ determine the support of the Hecke modifications. In particular, the Hecke modification is supported at $0\in \bbP^1$ if and only if the matrices $\gamma_{1, i}$ are nilpotent for all $i\in Q_0$. To compare with Coulomb branches as defined in \cite{BFNII}, it is necessary to restrict to modifications supported at a single point (and this point must be $0 \in \mathbb{P}^1$ to allow for equivariance with respect to the standard $\mathbb{C}^\times \hookrightarrow \text{Aut}(\mathbb{P}^1)$). This explains, from the geometric perspective, the relevance of the $\omega$-nilpotent CoHA for comparison with the Coulomb branch. 
\end{remark}

\subsection{Hall action on critical quasimaps}
\label{subsec: Hall action on critical quasimaps}

To remove clutter in what follows, let $\FX_{\dd}$ denote the critical locus of the function $\tr(W): \FM^T_{\dd}(\tilde Q) \to \BoC$. Likewise, we denote by $\FX_{\dd}^{\fr}$ the critical locus of $W: \FM^T_{\dd}(Q^{\QM}_{\ww,\vv}) \to \BoC$. We also work equivariantly with respect to the action of the torus $T$ on $\QM^{\xi}(X)$ defined in \ref{subsec: Torus action crit qm}, but we do not keep track of it in the notation. The superscript refer to the fact that $Q^{\QM}_{\ww,\vv}$ is a framed version of $\tilde Q$. In particular, by Theorem \ref{thm: critical description quasimaps}, we have an identification 
\[
\QM^{\xi}_{\dd}(X)=\FX_{\dd}^{\fr, \zeta^+\sst},
\]
where the right hand side it the semistable locus of $\FX_{\dd}^{\fr}$ with respect to the chosen stability condition $\zeta^+$.
Fix a decomposition $\dd=\dd'+\dd''$. Restricting \eqref{eq: hall action diamgram ambient spaces} to the relevant critical loci and applying Theorem \ref{thm: hecke as crit}, we get the diagram
\begin{equation}
\label{eq: diagram hall action quasimaps critical}
\begin{tikzcd}
    \FX_{\dd'} \times \FX_{\dd''}^{\fr} & \FX_{\dd',\dd''}^{\fr}\arrow[l, swap,  "q"]\arrow[r, "p"] & \FX_{\dd}^{\fr}
    \\
    \FX_{\dd'} \times \QM^{\xi}_{\dd''}(X) \arrow[u, hookrightarrow, "j"]& \Hecke^{\xi}_{\dd',\dd''}(X)\arrow[l, swap,  "\bar q"]\arrow[r, "\bar p"]\arrow[u, hookrightarrow, "k"] & 
    \QM^{\xi}_{\dd}(X)\arrow[u, hookrightarrow, "l"]
\end{tikzcd}
\end{equation}
Notice that the right diagram is Cartesian and the top row is simply the restriction of correspondence 
\[
\Grad(\FX_{\dd}^{\fr})\leftarrow \Filt(\FX_{\dd}^{\fr})\to \FX_{\dd}^{\fr}.
\]
to the component of topological type $(\dd',\dd'')$. Therefore, by \cite[Thm. B]{kinjo2026cohomologicalhallalgebras3calabiyau}, we have a canonical morphism
\begin{equation}
    \label{eq: Joice conj mor}
    q^* \varphi_{\FX_{\dd'} \times \FX_{\dd''}^{\fr}}\to p^! \varphi_{\FX_{\dd}}[I_{\dd'}].
\end{equation}
where the index $I_{\dd}$ can be explicitly computed as
\begin{multline*}
    I_{\dd'}=
    \left(\dim(\FM_{\dd}(Q^{\QM}_{\ww,\vv}))-\dim(\FM_{\dd',\dd''}(Q^{\QM}_{\ww,\vv}))\right)
    \\ -\left(\dim(\FM_{\dd',\dd''}(Q^{\QM}_{\ww,\vv}))-\dim(\FM_{\dd'}(\tilde Q)\times \FM_{\dd''}(Q^{\QM}_{\ww,\vv}))\right)
\end{multline*}
Notice that, although $I_{\dd'}$ may in principle depend on $\dd''$, an explicit computation using the formulas from \S\ref{subsec: quasimaps to quiver varieties} gives 
\begin{equation}
    \label{eq: index shifted action}
    I_{\dd'}=(\dd')^T(\ww-\Cartan \vv)
\end{equation}
where $\Cartan=2-\Q-\Q^T $ is the Cartan matrix of the quiver $Q$, see \S\ref{subsec: Quivers and their representations in a category}. 

\begin{remark}
The vector $\mu=\ww-\Cartan \vv$ has a clear representation theoretic interpretation. In fact, the cohomology of the disjoint union $\bigsqcup_{\vv\in \BoN^{Q_0}} \Nak(\vv,\ww)$ is a module for the Maulik-Okounkov Lie algebra $\Fg^{\MO}_Q$ (equivalently, BPS Lie algebra) and the cohomology of the summand $\Nak(\vv,\ww)$ is identified with a weight space for the Cartan subalgebra $\mathfrak{h}\subset \Fg^{\MO}_Q$ with weight $\mu=\ww-\Cartan \vv\in \BoN^{Q_0}$, (in the simple root basis). In \S\ref{subsec: Shifted Yangian action} we will also identify $\mu$ with the shift of the shifted quantum group extending the Hall algebra action on the vanishing cycle cohomology of $\bigsqcup_{\dd\in \BoN^{Q_0}}\QM^{\xi}_{\dd}(X)$ discussed in this section.
\end{remark}

Notice that $I_{\dd'}$ is non-zero in general because the stack $\FM_{\dd}(Q^{\QM}_{\ww,\vv})$ is not symmetric.
Since $k$ and $l$ are open inclusions, and hence smooth, we have a canonical morhism $k^*p^!\cong \bar p^! l^*$. Therefore, pulling back the map above along $k$ we obtain a canonical morphism 
\[
 \bar q^* j^* \varphi_{\FX_{\dd'} \times\FX_{\dd''}^{\fr}}=k^*q^* \varphi_{\FX_{\dd'} \times\FX_{\dd''}^{\fr}}\to k^*p^! \varphi_{\FX_{\dd}} [I_{\dd'}]\cong  \bar p^! l^*  \varphi_{\FX_{\dd}}[I_{\dd'}].
\]
Moreover, since the DT sheaf is compatible with restriction on open sets, we obtain a canonical morphism 
\begin{equation}
    \label{eq: joyce conj map for quasimaps}
    \bar q^* \varphi_{\FX_{\dd'} \times \QM^{\xi}_{\dd''}(X)}[I_{\dd'}] \to \bar p^! \varphi_{\QM^{\xi}_{\dd}(X)}.
\end{equation}
Passing to derived global sections we get a canonical morphism via the composition
\begin{align}
\begin{split}
\label{eq: construction action hall quasimaps}
    \HO(\FX_{\dd'} \times \QM^{\xi}_{\dd''}, \varphi_{\FX_{\dd'} \times \QM^{\xi}_{\dd''}})
    &\to \HO(\FX_{\dd'} \times \QM^{\xi}_{\dd''}, \bar q_* \bar q^* \varphi_{\FX_{\dd'} \times \QM^{\xi}_{\dd''}})
    \\
    &=\HO(\Hecke^{\xi}_{\dd',\dd''}(X), \bar q^* \varphi_{\FX_{\dd'} \times \QM^{\xi}_{\dd''}})
    \\
    & \to \HO(\Hecke^{\xi}_{\dd',\dd''}(X), \bar p^! \varphi_{\QM^{\xi}_{\dd}}[I_{\dd'}])
    \\
    &=\HO(\Hecke^{\xi}_{\dd',\dd''}(X), \bar p_!\bar p^! \varphi_{\QM^{\xi}_{\dd}}[I_{\dd'}])\\
    &\to 
    \HO(\QM^{\xi}_{\dd}, \varphi_{\QM^{\xi}_{\dd}}[I_{\dd'}])
\end{split}
\end{align}
Here, the first map is the co-unit $\id\to \bar q_*\bar q^*$, the second is induced by \eqref{eq: joyce conj map for quasimaps} last one is the unit $\bar p_! \bar p^!\to \id$ and along the way we are using properness of $\bar p$. 
In conclusion, applying Thom-Sebastiani and noticing that $\HCoha^T_{\tilde Q, \tilde W, \dd'}=\HO(\FX_{\dd'}, \varphi_{\FX_{\dd'}})$, we get a morphism 
\begin{equation}
    \label{eq: hall action quasimaps}
    \vact_{\dd',\dd''}:\HCoha^T_{\tilde Q, \tilde W, \dd'}\otimes \HO_T(\QM^{\xi}_{\dd''}, \varphi_{\QM^{\xi}_{\dd''}})\to \HO_T(\QM^{\xi}_{\dd}, \varphi_{\QM^{\xi}_{\dd}}[I_{\dd'}]).
\end{equation}
In the last equations we have recast the action of the torus $T$ in the notation. 
Set 
\[
    \HO_T(\QM^{\xi}, \varphi_{\QM^{\xi}})\coloneqq \bigoplus_{\dd\in \BoN^{Q_0}}  \HO_T(\QM^{\xi}_{\dd}, \varphi_{\QM^{\xi}_{\dd}}[I_{\dd'}]).
\]
\begin{theorem}
\label{thm: Hall acts on critical quasimaps}
    The morphisms \eqref{eq: hall action quasimaps} define an action of the Hall algebra $\HCoha^T_{\tilde Q, \tilde W}$ on the cohomology 
    \[
    \HO_T(\QM^{\xi}, \varphi_{\QM^{\xi}})=\bigoplus_{\dd\in \BoN^{Q_0}} \HO_T(\QM^{\xi}_{\dd}, \varphi_{\QM^{\xi}}).
    \]
\end{theorem}
\begin{proof}
    It suffices to check associativity, and the argument is standard. In the language of this section, it can be directly deduced by the third statement of \cite[Prop. 9.11]{kinjo2026cohomologicalhallalgebras3calabiyau}.
\end{proof}

\begin{remark}
\label{rem: alternative description action}
    The action map $\vact_{\dd', \dd''}$ can be alternatively defined using the ambient geometries of the moduli spaces of quiver representations. Formally, consider the correspondence 
    \[
    \begin{tikzcd}
        \FM^T_{\dd'}(\tilde Q) \times_{\B T} \FM^{T,\zeta^+ \sst}_{\dd''}(Q^{\QM}_{\ww,\vv}) & \FM^{T,\zeta^+\sst}_{\dd',\dd''}(Q^{\QM}_{\ww,\vv})\arrow[l, swap, hookrightarrow,  "q^{\zeta^+}"]\arrow[r, "p^{\zeta^+}"]& \FM^{T,\zeta^+\sst}_{\dd' + \dd''}(Q^{\QM}_{\ww,\vv})
    \end{tikzcd}
    \]
    from diagram \eqref{eq: hall action diamgram ambient spaces}. The morphism $q$ is smooth and $p$ is proper, and we have an equality 
    \[
    \tilde{W} \boxplus W^{\xi} \circ q^{\zeta^+} = W^{\xi} \circ p^{\zeta^+}.
    \]
    Therefore—ignoring cohomological shifts—we have maps in critical cohomology given by 
    \begin{multline*}
        \HO(\FM^T_{\dd'}(\tilde Q), \phip{\tilde W}\DD\BoQ) \otimes \HO(\FM^{T,\zeta^+ \sst}_{\dd''}(Q^{\QM}_{\ww,\vv}), \phip{W^{\xi}}\DD\BoQ)
        \\
        \xrightarrow{(q^{\zeta^+})^*} \HO(\FM^{T,\zeta^+\sst}_{\dd,\dd'}(Q^{\QM}_{\ww,\vv}), \phip{\tilde W\boxplus W^{\xi}\circ q^{\zeta^+}}\DD\BoQ)
         \\
         \xrightarrow{(p^{\zeta^+})_*} \HO(\FM^{T,\zeta^+\sst}_{\dd' + \dd''}(Q^{\QM}_{\ww,\vv}), \phip{W^{\xi}}\DD\BoQ).
    \end{multline*}
By Corollary \ref{cor: DT are the same}, we can identify the domain and target of the morphism above with the ones of \eqref{eq: hall action quasimaps}. Moreover, it is shown in \cite[\S7.7]{kinjo2026cohomologicalhallalgebras3calabiyau} that this construction is equivalent to the one inducing \eqref{eq: hall action quasimaps}. We will exploit this alternative characterization in \S\ref{subsec: Shifted Yangian action} and \S\ref{subsec: Dimensionally reduced action hall} to extend the Hal action to the action of a suitably defined shifted Yangian.
\end{remark}

We conclude this section discussing a variant of the previous construction. Recall from \S\ref{subsec: Nilpotent and omega-nilpotent 3d CoHA} that we have a canonical injective morphism of algebras $\HCoha^{T, \onil}_{\tilde Q, \tilde W}\to \HCoha^T_{\tilde Q, \tilde W}$. Combined with Theorem \ref{thm: Hall acts on critical quasimaps}, we obtain an action of the $\omega$-nilpotent Hall algebra $\HCoha^{T, \onil}_{\tilde Q, \tilde W}$ on $\HO_T(\QM^{\xi}, \varphi_{\QM^{\xi}})$. In the rest of the section, we translate this nilpotency condition into the language of Hecke correspondences. As before, we drop the torus $T$ from the notation.

Let $\Hecke^{\xi,0}_{\dd', \dd''}(X)$ be the closed substack of $ \Hecke^{\xi}_{\dd', \dd''}(X)$ parametrizing those Hecke modifications that are supported at the point $0\in \bbP^1$.
Let $\FX^{\onil}_{\dd}$ be the stack fitting in the leftmost of the following Cartesian squares 
\begin{equation}
\label{eq: diagram omega nilp and ambeint spaces}
    \begin{tikzcd}
        \FX^{\onil}_{\dd}\arrow[r, ""]\arrow[d, "\kappa_{\dd}"]& \FM ^{\onil}_{\dd}(\tilde Q)\arrow[d, "\eta_{\dd}"]\arrow[r, "\JH_{\dd}"] &\CM ^{\onil}_{\dd}(\tilde Q)\arrow[d, "\iota_{\dd}"] \\
        \FX_{\dd}\arrow[r] & \FM_{\dd}(\tilde Q)\arrow[r, "\JH_{\dd}"]&  \CM_{\dd}(\tilde Q)
    \end{tikzcd}
\end{equation}
and consider the following enlargement of diagram \eqref{eq: diagram hall action quasimaps critical}:
\begin{equation}
    \label{eq: diagram omega nilp and hecke}
\begin{tikzcd}
    &\FX_{\dd'} \times \QM^{\xi}_{\dd''} & \Hecke^{\xi}_{\dd',\dd''}(X)\arrow[l, swap, hookrightarrow, near start,  "\bar q"]\arrow[r, "\bar p"]& 
    \QM^{\xi}_{\dd}
    \\
    \\
    \FX^{\onil}_{\dd'}\times\QM^{\xi}_{\dd''}\arrow[uur, swap, near start, "\kappa_{\dd}\times \id"] &\Hecke^{\xi,0}_{\dd', \dd''}(X) \arrow[l, swap,  "\bar r"]\arrow[uur, near start, "t"]\arrow[uurr, "\bar s"]
\end{tikzcd}
\end{equation}
All the squares are Cartesian except from those with arrows $(q, \bar q, j,k)$, $(r, \bar r, \tilde j,\tilde k)$, and $(s, \bar s, \tilde k, l)$. In particular, the fact that the bottom square is Cartesian, i.e. that a Hecke correspondence is supported at zero iff the loops $\omega_i$ are nilpotent for all $i\in Q_0$, follows from Remark \ref{remark: support hecke modifications}. Therefore, by base change, we have a natural isomorphism of functors $\bar r^*(\kappa_{\dd}\times \id)^!\to t^! \bar q^*$ which, combined with \eqref{eq: joyce conj map for quasimaps}, gives a canonical morphism
\begin{equation}
\label{eq: joice conj for quasimpas nilp version}
    \bar r^* (\kappa_{\dd}\times \id)^!\varphi_{\FX_{\dd'} \times \QM^{\xi}_{\dd''}(X)} \to 
    t^! \bar q^* \varphi_{\FX_{\dd'} \times \QM^{\xi}_{\dd''}(X)} \to \bar s^! \varphi_{ \QM^{\xi}_{\dd}(X)}[I_{\dd'}]
\end{equation}
of complexes on $\Hecke^{\xi,0}_{\dd', \dd''}$. Notice that a standard base change argument using the diagram \eqref{eq: diagram omega nilp and ambeint spaces} gives $\HO(\FX^{\onil}_{\dd'}, \kappa_{\dd}^!\varphi_{\FX_{\dd'}})\cong \HCoha^{T, \onil}_{\tilde Q, \tilde W}$. Therefore, arguing as in \eqref{eq: construction action hall quasimaps} and recasting the torus $T$ in our notation, we get a canonical morphism 
\begin{equation}
\label{eq: hall action quasimaps nil version}
\vact^{T, \onil}_{\dd',\dd''}: \HCoha^{T, \onil}_{\tilde Q, \tilde W}\otimes \HO_T(\QM^{\xi}_{\dd''}, \varphi_{\QM^{\xi}_{\dd''}})\to \HO_T(\QM^{\xi}_{\dd}, \varphi_{\QM^{\xi}_{\dd}})[I_{\dd'}].
\end{equation}
\begin{proposition}
    The maps \eqref{eq: hall action quasimaps nil version} define an action of the nilpotent Hall algebra $\HCoha^{T,\onil}_{\tilde Q, \tilde W}$ on on the cohomology $\HO_T(\QM^{\xi}, \varphi_{\QM^{\xi}})$. Moreover, this action is compatible with the action of $\HCoha_{T,\tilde Q, \tilde W}$ on $\HO_T(\QM^{\xi}, \varphi_{\QM^{\xi}})$ via the canonical morphism $\HCoha^{T, \onil}_{\tilde Q, \tilde W}\to \HCoha^T_{\tilde Q, \tilde W}$. Equivalently, the following diagram commutes for all $\dd',\dd''\in \BoN^{Q_0}$:
    \begin{equation*}
        \begin{tikzcd}
            \HCoha^{T, \onil}_{\tilde Q, \tilde W, \dd'}\otimes \HO_T(\QM^{\xi}_{\dd''}, \varphi_{\QM^{\xi}_{\dd''}})\arrow[r, "\vact^{\onil}_{\dd',\dd''}"] \arrow[d] & \HO_T(\QM^{\xi}_{\dd}, \varphi_{\QM^{\xi}_{\dd}})[I_{\dd'}]
            \\
            \HCoha_{T, \tilde Q, \tilde W, \dd'}\otimes\HO_T(\QM^{\xi}_{\dd''}, \varphi_{\QM^{\xi}_{\dd''}})\arrow[ur, swap,  "\vact_{\dd',\dd''}"] 
        \end{tikzcd}
    \end{equation*}
\end{proposition}
The proof of the proposition follows from the definitions and standard six-functors operations using the map of the diagram \eqref{eq: diagram omega nilp and hecke}.

\begin{remark}
    In \cite{DHSM23, BD2023} it is shown that the cohomology
    \[
    \bigoplus_{\vv\in \BoN^{Q_0}} \HO_T(\Nak(\ww, \vv), \BoQ^{\vir})
    \]
    is isomorphic to a direct sum of irreducible lowest weight modules for the Maulik-Okounkov Lie algebra $\Fg^{\MO}_{Q}$ or, equivalently, for the generalized Kac-Moody Lie algebra whose positive part is given by the BPS Lie algebra $\Fg_{\tilde Q \tilde W}$. It would be very interesting to describe the critical cohomology of the direct sum 
    \[
     \bigoplus_{\dd\in \BoN^{Q_0}}\HO_T(\QM^{\xi}_{\dd}, \varphi_{\QM^{\xi}_{\dd}})
    \]
    in similar terms. Note in particular that in this case there are two geometrically meaningful actions to consider. The one of the BPS Lie algebra $\Fg^T_{\tilde Q, \tilde W}$, and the one of its $\omega$-nilpotent subalgebra $\Fg^{T,\onil}_{\tilde Q, \tilde W}$, see Remark \ref{rem: inj omega nilp lie algebra}. We also stress that, unlike for the case of quiver varieties, the answer to this question cannot be flat over $\HO_T$ as the embedding $\Fg^{T,\onil}_{\tilde Q, \tilde W}\hookrightarrow \Fg^{T}_{\tilde Q, \tilde W}$ only holds equivariantly.

\end{remark}

\subsection{Nakajima-type operators}
The action of certain generators of the cohomological Hall algebra can be described without reference to the Hall algebra itself. We now spell out the construction, which is a variant of the so-called Nakajima operators. Fix a vertex $i\in Q_0$, fix a dimension vector $\dd\in \BoN^{Q_0}$, and set $\dd'=\dd+\delta_i$. There is a tautological line bundle $\Sh{L}_i$ on $\FM_{\dd, \delta_i}(Q^{\QM}_{\ww,\vv})$ whose fiber at a closed point is given by $\ker(U_i'' \twoheadrightarrow U_i)$. We follow the now-standard convention of not keeping track of the action of the torus $T$ in our notation.

Restricting to the locus $\FM^{\nil}_{\dd, \delta_i}(Q^{\QM}_{\ww, \vv}) \subset \FM_{\dd, \delta_i}(Q^{\QM}_{\ww, \vv})$ where the unframed representation is nilpotent, and further to the stable loci in $\FM_{\dd}(Q^{\QM}_{\ww, \vv})$ and $\FM_{\dd + \delta_i}(Q^{\QM}_{\ww, \vv})$, we obtain a correspondence 
\begin{equation}
\begin{tikzcd}
    \overline{\FP}(\dd, \dd + \delta_i) \arrow[r, "\bar{q}^{\zeta^+} \times p^{\zeta^+}"] & \FM^{\zeta^+\sst}_{\dd}(Q^{\QM}_{\ww, \vv}) \times \FM^{\zeta^+\sst}_{\dd + \delta_i}(Q^{\QM}_{\ww, \vv}).
\end{tikzcd}
\end{equation}
where $\bar q^{\zeta^+}$ is the composition $\overline{\FP}(\dd, \dd + \delta_i) \to \FM^{\nil}_{\delta_i}(\tilde{Q}) \times \FM^{\zeta^+\sst}_{\dd}(Q^{\QM}_{\ww, \vv}) \to \FM^{\zeta^+\sst}_{\dd}(Q^{\QM}_{\ww, \vv})$. Note the above map can actually be identified as a closed inclusion, because any morphism of $\zeta^+$-stable framed quiver representations which is identity on the framing is unique if it exists and is necessarily surjective. This is an immediate consequence of the cyclic stability condition. 

We can then define the operators 
\[
e_{i, r}: \HO_T(\QM_{\dd}^\xi(X), \varphi_{\QM^\xi_{\dd}(X)})\to \HO_T(\QM_{\dd + \delta_i}^\xi(X), \varphi_{\QM^\xi_{\dd + \delta_i}(X)}) 
\]
by setting 
\[
e_{i, r}(\gamma)= (-1)^{t_i}(p^{\zeta^+})_*\left(c_1(\Sh{L_i})^r \cup  (\bar{q}^{\zeta^+})^*\gamma\right).
\]
The sign twist is as in \ref{theorem: shifted Yangian action on critical cohomology} below. Notice that $\FM^{\nil}_{\widetilde{Q}}(\delta_i) \simeq \B\mathbb{C}^\times$ and $\Sh{L}_i$ is obtained from the universal line bundle on $\B \mathbb{C}^\times$ by pullback. We have a well-defined correspondence action on critical cohomology because $(\bar{q}^{\zeta^+})^*(W^\xi_{\dd}) = (p^{\zeta^+})^*(W^\xi_{\dd + \delta_i})$, where pullbacks here denote pullback of functions. 

We may also define the operators 
\[ 
f_{i, r}: \HO_T(\QM_{\dd + \delta_i}^\xi(X), \varphi_{\QM^\xi_{\dd + \delta_i}(X)}) \to \HO_T(\QM^\xi_{\dd}(X), \varphi_{\QM^\xi_{\dd}(X)})
\]
by setting
\[ 
f_{i, r}(\gamma) = (\bar{q}^{\zeta^+})_*(c_1(\Sh{L}_i)^r \cup (p^{\zeta^+})^*\gamma).  
\]
That the pullback and pushforward are well-defined can be seen using the description of $p^{\zeta^+}$ and $\overline{q}^{\zeta^+}$ in Lemma \ref{lemma: projective bundle description of Nakajima operators} below.

Indeed, there is the following explicit description of $p^{\zeta^+}: \overline{\FP}(\dd - \delta_i, \dd) \to \FM^{\zeta^+\sst}_{\dd}(Q^{\QM}_{\ww, \vv})$ and $\bar{q}^{\zeta^+}: \overline{\FP}(\dd, \dd + \delta_i) \to \FM^{\zeta^+\sst}_{\dd}(Q^{\QM}_{\ww, \vv})$, useful for computations. Define the \textit{outgoing tautological complex} of vector bundles over $\FM^{\zeta^+\sst}_{\dd}(Q^{\QM}_{\ww, \vv})$ to be 
\[ 
\Sh{C}_i^+ \coloneqq 
\begin{tikzcd}
   F_i \arrow[r] & U_i \oplus W_i \oplus \bigoplus_{a: i \to j} F_j
\end{tikzcd}
\]
where $a$ runs over all edges in the doubled quiver from $i$ to $j$, and the map is the tautological one built from the quiver arrows. $F_i$ is placed in cohomological degree $0$, and cohomological degree increases from left to right. Recall the vector bundle $F_j$ is the universal extension 
\[ 
\begin{tikzcd}
    0 \arrow[r] & V_j \arrow[r] & F_j \arrow[r] & U_j \arrow[r] & 0. 
\end{tikzcd}
\]
The tautological map $F_i \to U_i$ is $\begin{pmatrix} B_{1, i} && I_i\end{pmatrix}$ and the tautological map $F_i \to W_i$ is $\begin{pmatrix} \alpha_i && \jmath_i \end{pmatrix}$. 

Define the \textit{incoming tautological complex} of vector bundles over $\FM^{\zeta^+\sst}_{\dd}(Q^{\QM}_{\ww, \vv})$ to be 
\[
\Sh{C}^-_i \coloneqq 
\begin{tikzcd}
    F_i \oplus \bigoplus_{b: j \to i} U_j \arrow[rr, "(B_{1, i} \, \, I_i \, \, B_b)"] && U_i
\end{tikzcd}
\]
with notations as above. $U_i$ is placed in cohomological degree $1$. 

Projectivizations of 2-term complexes of vector bundles supported in cohomological degrees $[-1, 0, 1]$ are well-defined as derived Artin stacks. The Hall correspondences are derived stacks in a canonical fashion, whence so are the $\overline{\FP}$ correspondences by restriction. We have the following 
\begin{lemma} \label{lemma: projective bundle description of Nakajima operators}
There are commutative diagrams 
\[
\begin{tikzcd}
   \overline{\FP}(\dd - \delta_i, \dd) \arrow[r, "\sim"] \arrow[d, "p^{\zeta^+}"] & \mathbb{P}(\Sh{C}^+_i) \arrow[dl] \\ \FM^{\zeta^+\sst}_{\dd}(Q^{\QM}_{\ww, \vv}) 
\end{tikzcd} \, \, \, \textnormal{and} \, \, \, 
\begin{tikzcd}
    \overline{\FP}(\dd, \dd + \delta_i) \arrow[r, "\sim"] \arrow[d, "\bar{q}^{\zeta^+}"] & \mathbb{P}((\Sh{C}^-_i)^\vee) \arrow[dl] \\ \FM^{\zeta^+\sst}_{\dd}(Q^{\QM}_{\ww, \vv}).
\end{tikzcd}
\]
The left isomorphism identifies $\Sh{L}_i$ with $\Sh{O}(-1)$, and the right identifies $\Sh{L}_i$ with $\Sh{O}(+1)$. The diagonal maps are the natural projections of the projective bundles.
\end{lemma}

\begin{proof}
    Observe that the fibers of $p^{\zeta^+}$ are identified canonically with 
    \[
    \mathbb{P}\Bigg(\ker \Big(F_i \longrightarrow 
    U_i \oplus W_i \oplus \bigoplus_{a: i \to j} F_j \Big) \Bigg) 
    \]
    and the fibers of $\bar{q}^{\zeta^+}$ are identified canonically with
    \[
    \mathbb{P}\Bigg( \coker \Big(U_i^* \longrightarrow F_i^* \oplus \bigoplus_{b: j \to i} U_j^* \Big) \Bigg)
    \]
    both of which follow e.g. immediately from Lemma \ref{lemma: quotient description of induction space} and the restriction to nilpotent unframed quiver representations. This identifies classical truncations; to see that the derived structure agrees with that arising from the Hall correspondence is a matter of checking definitions. 
\end{proof}

\begin{remark} \label{remark: geometric description of pole cancellation}
    Note that since we assume $(\xi_2, \xi_3, \iota, \jmath)$ is a stable quiver representation, the kernel of the restriction of the differential in $\Sh{C}_i^+$ to $V_i \hookrightarrow F_i$ is trivial. This means that the classical truncation of $\mathbb{P}(\Sh{C}_i^+)$ is fiberwise disjoint from the subvariety $\mathbb{P}(V_i) \hookrightarrow \mathbb{P}(F_i)$. Concretely, this means the poles at the torus weights of $V_i$ will not contribute in residue formulas like \eqref{eq: raising operator residue formula} below. 
\end{remark}

\subsection{Shifted Yangian action}
\label{subsec: Shifted Yangian action}
As discussed above, we have the nilpotent CoHA with a canonical morphism \eqref{eq: map nilp coha} to $\HCoha^T_{\tilde Q, \tilde W}$. We pass further to a subalgebra, the spherical nilpotent CoHA:
\begin{equation}
    \label{embedding spherical nilpotent in coha}
    \Sh{SH}^{T,\nil} \hookrightarrow \HCoha^{T, \nil}_{\tilde Q, \tilde W}.
\end{equation}
It is generated by the operators 
\[
e_{i, r} = (-1)^{t_i}c_1(\Sh{L}_i)^r \cap [\FM_{\widetilde{Q}}(\delta_i)^{\nil}],
\]
where $\Sh{L}_i$ be the tautological line bundle on $\FM_{\tilde Q}^{\nil}(\delta_i)$. The integer $t_i$ is defined in Theorem \ref{theorem: shifted Yangian action on critical cohomology} and plays no role until then.
We also introduce $\Sh{H}^0_Q \coloneqq \mathbb{C}[h_{i, r}]_{i \in Q_0, r \in \mathbb{Z}}$, which is a polynomial algebra in infinitely many variables. It is useful to organize these into generating series: given $i \in Q_0$ write
\[
e_i(u)=\sum_{r \geq 0} e_{i,r}u^{-r-1}\qquad h_i(u)=\sum_{r \in \mathbb{Z}} h_{i,r}u^{-r-1} \qquad f_i(u)=\sum_{r \geq 0} f_{i,r}u^{-r-1}.
\]
The infinitely shifted Yangian\footnote{In \cite{COZZ2}, authors refer to this as $\Sh{DSH}_\infty(\tilde Q, \tilde W)$ to distinguish it from Faddeev-Reshetikhin-Takhtajan Yangians arising from $R$-matrices. In this paper we will only have use for the cohomological Hall algebra version, so we use this notation.} $\Yang^T_\infty(\tilde Q, \tilde W)$ \cite{COZZ2} is the algebra obtained from the free product $\Sh{SH}^{T,\nil} * \Sh{H}^0_Q * \Sh{SH}^{T, \nil, \opp}$ by imposing the relations
\begin{equation} \label{eq: shifted Yangian relations}
\begin{split}
    e_{i, r} f_{j,s} - f_{j, s} e_{i, r} & = \delta_{ij} \gamma_i h_{i, r + s} \\ 
    h_i(u)e_j(v) & = \zeta_{ij}(u) e_j(v) h_i(u) \\ 
    h_i(u) f_j(v) & = \zeta_{ij}(u)^{-1} f_j(v) h_i(u). 
\end{split}
\end{equation}
The notation $\gamma_i \coloneqq \epsilon_1 \times \prod_{a: i \to i} \epsilon_a$, i.e. it is the product of equivariant variables for each self-loop at node $i$. Likewise\footnote{Recall that $a: i \to j$ and $b: j \to i$ refers to all arrows in the \textit{doubled} quiver.},
\[
\zeta_{ij}(u) \coloneqq \Big( \frac{u - \sigma_j - \epsilon_1}{u - \sigma_j + \epsilon_1}\Big)^{\delta_{ij}} \times \frac{\prod_{a: i \to j} ( u - \sigma_j - \epsilon_a)}{\prod_{b: j \to i}(u - \sigma_j + \epsilon_b)}
\]
where $\sigma_j$ acts by $e_{j, r} \mapsto e_{j, r + 1}$, $f_{j, r} \mapsto f_{j, r + 1}$. 

For $\mu \in \mathbb{Z}^{Q_0}$, the \textit{$\mu$-shifted Yangian} $\Yang^T_{\mu}(\tilde Q, \tilde W)$ is obtained from $\Yang^T_\infty(\tilde Q, \tilde W)$ by further imposing the relations 
\[
h_{i, r} = 
\begin{cases}
0 & \text{if $r < - \mu_i - 1$} \\ 
1 & \text{if $ r = - \mu_i - 1$}.
\end{cases}
\]
\begin{remark}
    Assume that $Q$ is an orientation of a Dynkin diagram. The associated Yangian—or, more generally, shifted Yangian—is classically defined as the quotient of a free algebra in $e_{i,r}, f_{i,r}, h_{i,r}$ by the ideal generated by a number of equations, including the Fourier modes of \eqref{eq: shifted Yangian relations}. The relations exclusively involving Fourier modes in $e_{i}(v)$ or in $f_{i}(v)$ are already encoded in $\Sh{SH}^{T,\nil}$ and $\Sh{SH}^{T, \nil, \opp}$, and hence do not appear in our relations. Moreover, it is worth stressing that the second and third relation in \eqref{eq: shifted Yangian relations} are sometimes modified by adding a term that is independent of $v$ (see, e.g., \cite[\S B5]{BFNII}). Since both $e_i(v)$ and $f_i(v)$ have no non-negative modes, this extra term does not modify our relations but, a priori, introduces new ones. However, one can check that these are not new relations but are already determined by those in \eqref{eq: shifted Yangian relations}.
\end{remark}
\begin{remark}
    It is easy to check that the anti involution on $\Sh{SH}^{T,\nil} * \Sh{H}^0_Q * \Sh{SH}^{T, \nil,\opp}$ given by
    \[
    e_i(u)\mapsto f_i(u)\qquad f_i(u)\mapsto e_i(u)\qquad h_i(u)\mapsto h_i(u)
    \]
    preserves relations \eqref{eq: shifted Yangian relations} and hence induces an anti-involution $\mathbf{\iota}: \Yang_{\mu}^{T}(\tilde Q, \tilde W) \to \Yang_{\mu}^{T}(\tilde Q, \tilde W)$. It is easy to check that $\iota$ is induced by the anti-involution $\Theta_{*}$ introduced in \S\ref{subsec: anti-involution}.
\end{remark}

Recall from \S\ref{subsec: Quivers and their representations in a category} the definition of $\Cartan$ Cartan matrix of the quiver $Q$. We have the following.

\begin{theorem} \label{theorem: shifted Yangian action on critical cohomology}
    Recall the Nakajima variety $X = \NN_Q^{\zeta^-}(\vv, \ww)$. Let $\mu = \ww - \Cartan \vv$. Notations for equivariant parameters as in \S\ref{subsec: torus action quiver varieties}. Then the assignments
    \[
    \begin{split}
    e_i(u) & \mapsto (-1)^{t_i} \frac{1}{u - c_1(\Sh{L}_i)} [\overline{\FP}(\dd, \dd+ \delta_i)] \\ 
    f_i(u) & \mapsto \frac{1}{u - c_1(\Sh{L}_i)}[\overline{\FP}(\dd - \delta_i, \dd)]^T \\ 
    h_i(u) & \mapsto c_u\Big( t_1 U_i - t_1^{-1} U_i + \sum_{a: i \to j} t_a U_j - \sum_{b: j \to i} t_b^{-1} U_j +W_i + \sum_{a: i \to j} t_a V_j - (1 + t_1^{-1})V_i\Big) \cup -
    \end{split}
    \]
    define an action of $\Yang_{\mu}^{T}(\tilde Q, \tilde W)$ on $\HO_T(\QM^\xi(X), \varphi_{\QM^\xi(X)})$, extending the action of the spherical nilpotent cohomological Hall algebra induced by Theorem \ref{thm: Hall acts on critical quasimaps} and the canonical morphism \eqref{eq: map nilp coha}. The sign twist is chosen as 
    \[
    t_i = v_i + w_i + \sum_{a: i \to j}(v_j + d_j).
    \]
\end{theorem}
This theorem essentially follows from the computation done in the proof of Theorem 3.3 in \cite{COZZ2}, but for the convenience of the reader we will provide a self-contained argument. The proof is just a computation, but our proof is lengthy because we have explained all the steps. 

\begin{proof}
The most nontrivial commutation relation to check is the first line of \eqref{eq: shifted Yangian relations}, so we begin with this one. Note that the operator 
\[ 
e_{i, r} f_{j, s} - f_{j, s} e_{i, r}
\]
is supported on a cycle in $\FM^{\zeta^+\sst}_{\dd + \delta_i - \delta_j}(Q^{\QM}_{\ww, \vv}) \times \FM^{\zeta^+\sst}_{\dd}(Q^{\QM}_{\ww, \vv}) \ni (U', U)$, where $U'$ and $U$ are a pair of $Q^{\QM}_{\ww, \vv}$-modules in the corresponding moduli spaces. From definitions, $e_{i, r} f_{j, s}$ is supported on the locus where there exists a pair of surjections 
\[
\begin{tikzcd}
    & U \arrow[d, two heads] \\ U' \arrow[r, two heads] & U''
\end{tikzcd}
\]
to a framed stable quiver representation $U'' \in \FM^{\zeta^+\sst}_{\dd - \delta_j}(Q^{\QM}_{\ww, \vv})$. Note that, because we consider framed stable quiver representations, if a morphism to $U''$ exists then it is unique and surjective. On this locus, $e_{i, r} f_{j, s}$ coincides with the operator of multiplication by $c_1(\Sh{L}_i)^r c_1(\Sh{L}_j)^s$, where $\Sh{L}_i = \ker(U' \twoheadrightarrow U'')$ and $\Sh{L}_j = \ker(U \twoheadrightarrow U'')$ are the tautological line bundles. 

Similarly, $f_{j, s} e_{i, r}$ is supported on the locus where there exists a pair of surjections 
\[
\begin{tikzcd}
    U''' \arrow[r, two heads] \arrow[d, two heads] & U \\ U' 
\end{tikzcd}
\]
from a framed stable quiver representation $U''' \in \FM^{\zeta^+\sst}_{\dd + \delta_i}(Q^{\QM}_{\ww, \vv})$. On this locus it coincides with the operator of multiplication by $c_1(\Sh{L}_i)^r c_1(\Sh{L}_j)^s$, where $\Sh{L}_i = \ker(U''' \twoheadrightarrow U)$ and $\Sh{L}_j = \ker(U''' \to U')$ are the tautological line bundles. 

Now observe that if $i \neq j$ or if $i = j$ but $U$ and $U'$ are not isomorphic, then $U'''$ is uniquely determined by $U''$ and vice versa by taking pullbacks and pushouts (it is again important to invoke here that morphisms between framed stable quiver representations are unique and surjective if they exist). Moreover, this construction naturally identifies the tautological line bundles, therefore the operator $\comm{e_{i, r}}{f_{j, s}}$ vanishes if $i \neq j$ and is supported on the diagonal $\Delta \subset \FM^{\zeta^+\sst}_{\dd}(Q^{\QM}_{\ww, \vv}) \times \FM^{\zeta^+\sst}_{\dd}(Q^{\QM}_{\ww, \vv})$ if $i = j$: 
\[
e_{i, r} f_{j, s} - f_{j, s} e_{i, r} = \delta_{ij} \Delta_*(\beta_{i, r, s}) \in \HO_T(\FM^{\zeta^+\sst}_{\dd}(Q^{\QM}_{\ww, \vv}) \times \FM^{\zeta^+\sst}_{\dd + \delta_i - \delta_j}(Q^{\QM}_{\ww, \vv}), \phip{W_2 - W_1} \BoQ^{\text{vir}}_{\FM \times \FM}). 
\]
 The essential observation is that the potential $W_2 - W_1$ vanishes when restricted to $\Delta$, so that $\beta_{i, r, s}$ may be computed in ordinary (i.e. non-critical) equivariant cohomology; it is the class $\comm{e_{i, r}}{f_{i, s}} \cdot 1_{\dd}$, where now we understand the $e$ and $f$ correspondences in ordinary equivariant cohomology of the ambient spaces $\FM^{\zeta^+\sst}_{\dd}(Q^{\QM}_{\ww, \vv})$ of framed stable quiver representations. 

By Lemma \ref{lemma: projective bundle description of Nakajima operators} and Bott residue formula, we have\footnote{For readers uncomfortable with the derived world, all we are using here is that the map $p^{\zeta^+}$ may be factored as the inclusion of the zero locus of a section of a vector bundle, followed by a projective space bundle. We are just getting a numerator from the Euler class of the vector bundle, followed by the usual Bott residue formula for the fiber integration. A similar consideration holds for $\bar{q}^{\zeta^+}$ below.}
\begin{equation} \label{eq: raising operator residue formula}
\begin{split}
e_{i, r} \cdot 1_\dd  & = (-1)^{t_i} \oint_\gamma \frac{dx}{2\pi i} x^r \prod_{m = 1}^{d_i + 1} \frac{u^{(i)}_m - x + \epsilon_1}{u^{(i)}_m - x} \times  \prod_{a: i \to j} \prod_{n = 1}^{d_j + \delta_{ij}} (u^{(j)}_n - x + \epsilon_a) \times \\ 
& \frac{\prod_{a: i \to j} \prod_{l = 1}^{v_j}(v^{(j)}_l - x + \epsilon_a)}{\prod_{l = 1}^{v_i}(v^{(i)}_l - x)} \prod_{\alpha = 1}^{w_i}(a^{(i)}_\alpha - x) \\
& = \Bigg\{ c_u\Big( t_1 U_i - U_i + \sum_{a: i \to j} (t_a U_j + t_a V_j) + W_i - V_i  \Big) \Bigg\}_{u^{-r - 1}} \in \HO_T(\FM^{\zeta^+\sst}_{\dd + \delta_i}(Q^{\QM}_{\ww, \vv}), \BoQ^{\vir}_{\FM}).
\end{split}
\end{equation}
In the first line, the symbols $u^{(i)}_m$ are the Chern roots of the tautological bundle $U_i$ on $\FM^{\zeta^+\sst}_{\dd + \delta_i}(Q^{\QM}_{\ww, \vv})$, and likewise $v^{(i)}_l$ are the $T$-weights of $V_i$, and $a^{(i)}_{\alpha}$ are the $T$-weights of $W_i$. Integral means formal extraction of residue at $u^{(i)}_m$ and $v^{(i)}_l$. The subscript in the second line means take the coefficient of $u^{- r - 1}$ in the $u \to \infty$ expansion. The sign disappears by our choice of $t_i$, and the fact that the number of self-loops at any node of the doubled quiver is always even.

By similar reasoning, there is a formula 
\begin{equation} \label{eq: lowering operator residue formula}
\begin{split}
    f_{i, s} \cdot 1_\dd & = \oint_{\gamma} \frac{dx}{2\pi i} x^s \prod_{m = 1}^{d_i - 1} \frac{x - u^{(i)}_m}{x - u^{(i)}_m + \epsilon_1} \times \frac{1}{\prod_{l = 1}^{v_i} (x - v^{(i)}_l  + \epsilon_1 )} \times \\ & \prod_{b: j \to i} \prod_{n = 1}^{d_j - \delta_{ij}}\frac{1}{x - u^{(j)}_n + \epsilon_b} \\
    & = \Bigg\{c_v\Big(U_i - t_1^{-1} U_i - \sum_{b: j \to i} t_b^{-1} U_j - t_1^{-1} V_i \Big) \Bigg\}_{v^{-s - 1}} \in \HO_T(\FM^{\zeta^+\sst}_{Q^{\QM}_{\ww, \vv}}(\dd - \delta_i), \BoQ^{\vir}_{\FM}).  
\end{split}
\end{equation}
We may iterate this procedure to obtain formulas for the compositions $e_{i, r} f_{i, s} \cdot 1_\dd$ and $f_{i, s} e_{i, r} \cdot 1_\dd$. For the first, bearing in mind the exact sequence 
\[
\begin{tikzcd}
    0 \arrow[r] & \Sh{L}_i \arrow[r] & U_\dd \arrow[r] & U_{\dd - \delta_i} \arrow[r] & 0
\end{tikzcd}
\]
of vector bundles on the support of $e_{i, r}$, we have 
\begin{equation}
\begin{split}
    e_{i, r} f_{i, s} \cdot 1_\dd & = (-1)^{t_i} \oint_\gamma \frac{dx}{2\pi i} x^r (\text{as in RHS of \eqref{eq: raising operator residue formula}}) \times \\ 
    & \Bigg\{ \frac{v - x + \epsilon_1}{v - x} \prod_{b: i \to i} (v - x + \epsilon_b) \times c_v\Big( U_i - t_1^{-1} U_i - \sum_{b: j \to i} t_b^{-1} U_j - t_1^{-1} V_i \Big) \Bigg\}_{v^{-s - 1}} \\
    & = \Bigg\{ \Bigg\{ \frac{v - u + \epsilon_1}{v - u} \prod_{b:i \to i} (v - u + \epsilon_b) \times \\
    & c_u\Big( t_1 U_i - U_i + \sum_{a: i \to j} (t_a U_j + t_a V_j) + W_i - V_i \Big) \times  \\
    & c_v\Big( U_i - t_1^{-1} U_i - \sum_{b: j \to i} t_b^{-1} U_j - t_1^{-1} V_i \Big) \Bigg\}_{v^{-s - 1}} \Bigg\}_{u^{- r - 1}} \in \HO_{T}(\FM^{\zeta^+\sst}_{\dd}(Q^{\QM}_{\ww, \vv}), \BoQ^{\vir}_{\FM}). 
\end{split}
\end{equation}
In a similar fashion, we calculate 
\begin{equation}
\begin{split}
    f_{i, s} e_{i, r} \cdot 1_\dd & = \Bigg\{ \Bigg\{ \frac{v - u + \epsilon_1}{v - u} \prod_{a: i \to i}(v - u + \epsilon_a) \times \\
    & c_u\Big( t_1 U_i - U_i + \sum_{a: i \to j} (t_a U_j + t_a V_j) + W_i - V_i \Big) \times  \\
    & c_v\Big( U_i - t_1^{-1} U_i - \sum_{b: j \to i} t_b^{-1} U_j - t_1^{-1} V_i \Big) \Bigg\}_{u^{-r - 1}} \Bigg\}_{v^{- s- 1}} \in \HO_{T}(\FM^{\zeta^+\sst}_{\dd}(Q^{\QM}_{\ww, \vv}), \BoQ^{\vir}_{\FM}). 
\end{split}
\end{equation}
Now let $g_1(u), g_2(v)$ be rational functions, and let $P(u, v)$ be a polynomial. We have the integral identity 
\begin{equation}
\begin{split}
& \int_{|v| > |u| > \text{poles}} \frac{du}{2\pi i} \frac{dv}{2\pi i} \frac{P(u, v)}{v - u} g_1(u) g_2(v) - \int_{|u| > |v| > \text{poles}} \frac{du}{2\pi i} \frac{dv}{2\pi i} \frac{P(u, v)}{v - u} g_1(u) g_2(v) \\
& = \int_{|u| > \text{poles}} \frac{du}{2\pi i} P(u, u) g_1(u) g_2(u) 
\end{split}
\end{equation}
established by a standard contour deformation argument and the residue theorem. Applied in our situation this gives 
\begin{equation}
\begin{split}
    \beta_{i, r, s} & = \comm{e_{i, r}}{f_{i, s}} \cdot 1_{\dd}
     = \epsilon_1 \prod_{a: i \to i} \epsilon_a \times \\
    &\Bigg\{ c_u \Big( t_1 U_i - t_1^{-1} U_i + \sum_{a: i \to j} t_a U_j - \sum_{b: j \to i} t_b^{-1} U_j + \sum_{a: i \to j} t_a V_j + W_i - (1  + t_1^{-1})V_i\Big) \Bigg\}_{u^{-r - s - 1}}.
\end{split}
\end{equation}
This completes the verification of the first line of \eqref{eq: shifted Yangian relations}. To determine the shift $\mu$, we note it is uniquely determined by the asymptotics $h_i(u) = u^{\mu_i}(1 + O(u^{-1}))$ for $u \to \infty$. Noting that the quiver is doubled, so there is a bijection between $a: i \to j$ and $b: j \to i$, the first four terms in parentheses above do not contribute to $\mu$ and we read off 
\[
\mu_i = \ww_i + \sum_{a: i \to j} \vv_j - 2\vv_i = \ww_i - \sum_{j \in Q_0} \Cartan_{ij} \vv_j
\]
where $\Cartan_{ij}$ are the matrix elements of the Cartan matrix $\Cartan$. 

The second and third line of \eqref{eq: shifted Yangian relations} follow immediately from the fact that Chern polynomials are multiplicative over exact sequences. 
\end{proof}

\begin{remark}
\label{Cartan Coulomb and Hall commute}
We can anticipate the compatibility of the Hall and Coulomb branch action discussed in \S\ref{subsec: intro 3} by comparing the generators expressed here in the quiver language with the ``Coulomb branch'' formulas for the generators appearing in Theorem \ref{eq: generators hall to Coulomb} below. From the defining short exact sequence 
\[
\begin{tikzcd}
    0 \arrow[r] & \Sh{V}_i \arrow[r] & F_i \otimes \Sh{O}_{\mathbb{P}^1 \times \QM} \arrow[r] & U_i \otimes \Sh{O}_{\mathbb{P}^1 \times \QM}(1) \arrow[r] & 0 
\end{tikzcd}
\]
of the quiver description from \S\ref{sec: quiver critical description quasimaps} we deduce the relation 
\[
\Sh{V}_i |_0 = V_i + (1 - t_1)U_i \in K_T(\QM^\xi(X))
\]
where the left hand side is the $K$-class of the fiber of $\Sh{V}_i$ at zero. In terms of $\Sh{V}_i |_0$ the virtual bundle entering the formula for $h_i(u)$ takes a simpler form:
\[
\begin{split}
h_i(u) & =  c_u\Big( t_1 U_i - t_1^{-1} U_i + \sum_{a: i \to j} t_a U_j - \sum_{b: j \to i} t_b^{-1} U_j + \sum_{a: i \to j} t_a V_j + W_i - (1 + t_1^{-1})V_i  \Big) \\
& = \frac{c_u(W_i) \prod_{a: i \to j} c_{u - \epsilon_a}( \Sh{V}_j \eval_0)}{c_u(\Sh{V}_i |_0) c_{u + \epsilon_1}( \Sh{V}_i |_0)}
\end{split}
\]
where we have used that 
\[
\sum_{a: i \to j} t_a U_j - \sum_{b: j \to i} t_b^{-1} U_j = \sum_{a: i \to j} t_a(1 - t_1) U_j 
\]
since summations run over arrows in the doubled quiver and the torus variables are at the Calabi-Yau specialization.
\end{remark}

\begin{remark}
    The action of the shifted Yangian $\Yang^T_{\mu}(\tilde Q, \tilde W)$ is obtained by first restricting the action of the Hall algebra $\HCoha^T_{\tilde Q, \tilde W}$ to the nilpotent spherical subalgebra $\Sh{SH}^{T, \nil}$ and then realizing its Drinfeld double by transposing the Nakajima operators realizing the action of the spherical generators of $\Sh{SH}^{T, \nil}$. It would be very interesting extend this ``doubling'' process to the action of the ``whole'' 3d Hall algebra $\HCoha^T_{\tilde Q, \tilde W}$. This process would realize a shifted analog of the (un-shifted) Maulik-Okounkov Yangian \cite{MO19}.
\end{remark}

\subsection{Dimensionally reduced action}
\label{subsec: Dimensionally reduced action hall}
We conclude the chapter about the Yangian action on quasimap critical cohomology by analyzing the case when the evaluation point $\xi \in X$ is polarized in the sense of \S\ref{subsec: ev point and dim red}. To this end, it will be convenient to use the sign-twisted version of Hall algebras introduced in Remark \ref{rem: signed product hall algebras}. To make the sign-twisted Hall algebra $\HCoha^{T}_{\tilde Q, \tilde W}$ act on quasimaps, we further twist the action map \eqref{eq: hall action quasimaps} by $(-1)^{(\dd''+\vv)^T \Q \dd'}$. With this sign convention, one can check that $\HCoha^{T,\circ}_{\tilde Q, \tilde W}$ acts on $\HO_T(\QM^\xi_{\dd}(\bbP^1, X), \varphi_{\QM^\xi_{\dd}(\bbP^1, X)})$. We will always work with this sign convention whenever considering the action of $\HCoha^{T,\circ}_{\tilde Q, \tilde W}$ instead of the $\HCoha^{T}_{\tilde Q, \tilde W}$-action.

Since $\xi$ is polarized, Proposition \ref{prop: dimred Beilison model},gives a dimensional reduction isomorphism
\[
    dr: \HO_T(\QM^\xi_{\dd}(\bbP^1, X), \varphi_{\QM^\xi_{\dd}(\bbP^1, X)})\cong \HO_T^{\BM}(\QM^\xi_{\dd}(\bbP^1, L), \BoQ^{\vir}),
\]
On the other hand, by Theorem \ref{thm: dr coh 3d coha}, we have an isomorphism of algebras $dr: \HCoha^T_{\Xi_Q}\to \HCoha^{T, \circ}_{\tilde Q, \tilde W}$
Therefore, Theorem \ref{thm: Hall acts on critical quasimaps} gives an action\footnote{We incorporate the sign $(-1)^{(\dd''+\vv)^T \Q \dd'}$ in this action.}
\begin{equation}
	\label{eq: 2d hall action dimensionally reduced quasimaps}
 {\vact}^{\Xi_Q}_{\dd',\dd''}: \HCoha^T_{\Xi_Q,\dd'}\otimes \HO_T^{\BM}(\QM^\xi_{\dd'}(\bbP^1, L), \BoQ^{\vir})\to \HO_T^{\BM}(\QM^\xi_{\dd"+\dd''}(\bbP^1, L), \BoQ^{\vir})[I_{\dd'}].
\end{equation}
In this section, we discuss how to interpret this action by solely relying on the geometry of the moduli spaces $\QM^{\xi}(\bbP^1, L)$ and $\FM(\Xi_Q)$. Let 
\[
\Hecke^{\xi}(L)=\lbrace (\Sh{V},\Sh{V}',s, s') \; | \; \Sh{V}' \subset \Sh{V} + \text{compatibility of sections} \rbrace /\sim
\]
be the moduli space of Hecke modifications for $\QM^\xi(\bbP^1, L)$. Since the stability condition on the variety $L= \Rep_{Q'_{\ff}}(\dd,1)/G\cap X$ is inherited from the quiver variety $X$, it follows that $\HO^0(\Sh{V(}-1))=\HO^0(\Sh{V}'(-1))=0$ and therefore the same argument of \S\ref{subsec: Hecke correspondence on quasimap space} induces a canonical correspondence 
\begin{equation}
	\label{eq: dimensionally redduce hecke correspondence 1}
\begin{tikzcd}
    \FCoh_{Q, 0, \dd}(\AA^1)\times\QM^{\xi}_{\dd'}(L)  & \Hecke^{\xi}_{\dd, \dd'}(L)\arrow[l, swap, "\bar q\times \bar p"]\arrow[r, "\bar p'"] & \QM^{\xi}_{\dd+\dd'}(L)
\end{tikzcd}
\end{equation}
where $\Coh_{Q, 0}(\AA^1)$ is the moduli stack $Q$-valued torsion sheaves on $\AA^1$ introduced in \S\ref{subsec: Q decorated coherent sheaves on A1}. As discussed in \emph{loc. cit.}, we have a natural equivalence of stacks $\Coh_{Q, 0, \dd}(\AA^1)\simeq \FM_{\dd}(\Xi_Q)$. Moreover, by Proposition \ref{prop: dimensionally reduced quasimap quiver model}, the moduli space $\QM^{\xi}_{\dd'}(L)$ is equivalent to $\FM^{\zeta^+\sst}_{\dd'}(\Xi_{Q_{\ww,\vv}})$. Let now $\FM^{\zeta^+\sst}_{\dd',\dd''}(\Xi_{Q_{\ww,\vv}})\subseteq \FM^{\zeta^+\sst}_{\dd',\dd''}(\hat Q_{\ww,\vv})$ be the substack fitting in the Cartesian diagram 
\[
\begin{tikzcd}
      \FM^{\zeta^+\sst}_{\dd',\dd''}(\Xi_{Q_{\ww,\vv}})\arrow[r, hookrightarrow]\arrow[d,"\bar p'"] &\FM^{\zeta^+\sst}_{\dd,\dd'}(\hat Q_{\ww,\vv})\arrow[d, "\bar p"]
      \\
      \FM^{\zeta^+\sst}_{\dd+\dd'}(\Xi_{Q_{\ww,\vv}})\arrow[r, hookrightarrow] &\FM^{\zeta^+\sst}_{\dd+\dd'}(\hat Q_{\ww,\vv})
\end{tikzcd}
\]
Arguing along the lines of Theorem \ref{thm: hecke as crit} and Proposition \ref{prop: map coha to omega nil coha}, one can check that there is a natural equivalence of stacks
\[
 \Hecke^{\xi}_{\dd, \dd'}(L)\simeq \FM^{\zeta^+\sst}_{\dd',\dd''}(\Xi_{Q_{\ww,\vv}})
\]
identifying \eqref{eq: dimensionally redduce hecke correspondence 1} with 
\begin{equation}
\label{eq: correspondence 2d coha action dimensionally reduced quasimap}
\begin{tikzcd}
	\FM_{\dd}(\Xi_Q)\times \FM^{\zeta^+\sst}_{\dd'}(\Xi_{Q_{\ww,\vv}}) & \FM^{\zeta^+\sst}_{\dd',\dd''}(\hat Q_{\ww,\vv}) \arrow[l,  swap, "\bar q\times \bar p"]\arrow[r,  "\bar p'"] & \FM^{\zeta^+\sst}_{\dd'}(\Xi_{Q_{\ww,\vv}}).
\end{tikzcd}
\end{equation}
Informally, this correspondence is just a framed version of the top row of \eqref{eq: 2d coha coh multiplication diagram}. We remark that our choice of stability guarantees that the morphism $\bar p'$ is proper. Quasi-smoothness of $\bar q\times \bar p$ follows from the construction above \eqref{lemma: quasi-smoothness pull 2d coha} and the lemma itself. Therefore, by the same argument outlined in \ref{subsec: Q decorated coherent sheaves on A1},  the correspondence above induces a (shifted) action of the Hall algebra $\HCoha^T_{\Xi_Q}$ on $\HO_T^{\BM}(\QM^\xi_{\dd"+\dd''}(\bbP^1, L), \BoQ^{\vir})$ via pull-push in Borel-Moore homology.
\begin{proposition}
\label{prop: dimensional reduction hall action}
	The $\HCoha_{\Xi_Q}$-action on $\HO^{\BM}(\QM^{\xi}_{\dd"+\dd''}(\bbP^1, L), \BoQ^{\vir})$ induced by \eqref{eq: correspondence 2d coha action dimensionally reduced quasimap} coincides with the action \eqref{eq: 2d hall action dimensionally reduced quasimaps}. Equivalently, we have a commutative diagram
    \begin{equation*}
        \begin{tikzcd}
            \HCoha^T_{\Xi_Q,\dd'}\otimes \HO_T^{\BM}(\QM^\xi_{\dd'}(\bbP^1, L), \BoQ^{\vir})\arrow[r, "\bar p'_*(\bar q\times \bar p)^!"] \arrow[d, "dr\otimes dr"]&\HO_T^{\BM}(\QM^{\xi}_{\dd}(\bbP^1, L), \BoQ^{\vir})[I_{\dd'}]\arrow[d, "dr"]
            \\
           \HCoha_{\tilde Q, \tilde W, \dd'}\otimes \HO_T(\QM^{\xi}_{\dd''}, \varphi_{\QM^{\xi}_{\dd''}})\arrow[r, "\vact_{\dd',\dd''}"] & \HO_T(\QM^{\xi}_{\dd}, \varphi_{\QM^{\xi}_{\dd}})[I_{\dd'}].
        \end{tikzcd}
    \end{equation*}
\end{proposition}
\begin{proof}[Sketch of the proof]
	The proof parallels the proof of Theorem \ref{thm: dr coh 3d coha}, which is in turm base on \cite[Appendix 1]{RS17}. Following \S\ref{subsec: Doubling, looping, tripling, and framing}, we remove clutter by denoting quiver representations by $\AA_{\dd}(Q)$ instead of $\Rep_Q(\dd)$. We will only outline the main steps and let refer the interested reader to \emph{loc.cit.} for further details. Let $\ol \AA_{\dd}(\Xi_{Q_{\ww, \vv}})$ be  the restriction  of the (trivial) fibration  $\AA_{\dd}( Q^{\QM}_{\ww,\vv})\to \AA_{\dd}( \hat Q_{\ww,\vv})$ on $\AA_{\dd}(\Xi_{Q_{\ww, \vv}})\subseteq  \AA_{\dd}( \hat Q_{\ww,\vv})$.  Similarly, let $\ol \AA_{\dd}(\Xi_{Q})$ be  the restriction  of  $\AA_{\dd}( \tilde Q)\to \AA_{\dd}( \hat Q)$ on $\AA_{\dd}(\Xi_{Q})\subseteq  \AA_{\dd}( \hat Q )$.  Additionally, set 
	\[
	Z=\bigoplus_{a\in Q_1}\Hom(\BoC^{\dd_{s(a)}''}\oplus \BoC^{\vv_{s(a)}}, \BoC^{\dd_{t(a)}'})\qquad Z^\vee=\bigoplus_{a\in Q_1}\Hom(\BoC^{\dd_{s(a^*)}'}, \BoC^{\dd_{t(a^*)}''}\oplus \BoC^{\vv_{t(a^*)}}).
	\]
	Notice that since $s(a)=t(a^*)$ and $t(a)=s(a^*)$, the vector spaces $Z$ and $Z^\vee$ are dual to each other via trace pairing. Consider now the following diagram
	\begin{equation}
		\label{eq: big diagram dimred prequotients 1}
		\begin{tikzcd}[column sep=small]
			\AA_{\dd'}(\tilde Q)\times \AA^{\zeta^+\sst}_{\dd''}(Q^{\QM}_{\ww,\vv})  &\arrow[l, swap, "a"] Z\times \AA_{\dd'}(\tilde Q)\times \AA^{\zeta^+\sst}_{\dd''}(Q^{\QM}_{\ww,\vv}) \arrow[d, "d"] & \AA^{\zeta^+\sst}_{\dd',\dd''}(Q^{\QM}_{\ww,\vv}) \arrow[d, "e"]\arrow[l, swap, "g"]
			\\
			Z^{\vee} \times \AA_{\dd'}(\tilde Q)\times \ol \AA^{\zeta^+\sst}_{\dd''}(Q^{\QM}_{\ww,\vv}) \arrow[r, "f"] \arrow[u, "c"]& Z^{\vee}\times Z\times \AA_{\dd'}(\tilde Q)\times \AA^{\zeta^+\sst}_{\dd''}(Q^{\QM}_{\ww,\vv}) & Z^{\vee}\times \AA^{\zeta^+\sst}_{\dd',\dd''}(Q^{\QM}_{\ww,\vv}) \arrow[d, "h"]\arrow[l, swap, "b"]
			\\
			& &\AA^{\zeta^+\sst}_{\dd'+\dd''}(Q^{\QM}_{\ww,\vv})
		\end{tikzcd}
	\end{equation}
	The maps in the diagram are defined as follows. 
		\begin{itemize}
		\item The maps $a$ and $c$ are just the canonical projections.
		\item The map $f$ is induced by the inclusion $\{0 \}\to Z$.
		\item The maps $d$ and $e$ are induced by the inclusion $\{0\}\to Z^{\vee}$.
		\item To define $h$, first consider the $\BoC^*$-action on the vector space $\AA_{\dd'+\dd''}(Q^{\QM}_{\ww,\vv})$ whose attracting locus is $\AA_{\dd',\dd''}(Q^{\QM}_{\ww,\vv})$.  Then the subspace of negative weights contains $Z^{\vee}$ as a summand. Therefore, we have a canonical inclusion 
		\[
		Z^{\vee}\times \AA_{\dd',\dd''}(Q^{\QM}_{\ww,\vv})\hookrightarrow \AA_{\dd'+\dd''}(Q^{\QM}_{\ww,\vv}),
		\]
		which we denote by $h$.
		\item The map $b$ is $\id_Z^{\vee}\times g$.
        \item  The map $g$ was defined as follows. The projection of $g$ on the factor $\AA_{\dd'}(\tilde Q)\times \AA^{\zeta^+\sst}_{\dd''}(Q^{\QM}_{\ww,\vv})$ is just the canonical projection to the fixed locus of the $\BoC^*$ action considered above. On the other hand, given an element\footnote{Notice what we are calling $A_a, A_{a^*}$ here are referred to $B_{2, e}, B_{3, e}$ elsewhere in the paper.}
        \[
            x=(B_{1,i}, A_{a}, A_{a^*}, J_a, J_a^*, I_i, \alpha_i)\in \AA_{\dd',\dd''}(Q^{\QM}_{\vv,\ww})\subseteq \AA_{\dd'+\dd''}(Q^{\QM}_{\vv,\ww})
        \]
       we declare the projection of $g(x)$ on $Z$ to be given by taking the matrix 
       \[
      \begin{pmatrix}
        A_{a} B_{1, s(a)} - B_{1, t(a)} A_{a} + I_{t(a)}J_{a} & A_{a} I_{s(a)} - I_{t(a)}\xi_{a}
      \end{pmatrix}
       \in \Hom(\BoC^{\dd_{s(a)}}\oplus \BoC^{\vv_{s(a)}}, \BoC^{\dd_{t(a)}})
    \]
and projecting it to the factor 
\begin{equation*}
    \label{eq: fiber vbundle virtual pull half hecke quiver version}
    Z=\Hom(\BoC^{\dd_{s(a)}''}\oplus \BoC^{\vv_{s(a)}}, \BoC^{\dd_{t(a)}'})\subseteq \Hom(\BoC^{\dd_{s(a)}}\oplus \BoC^{\vv_{s(a)}}, \BoC^{\dd_{t(a)}}).
\end{equation*}

	\end{itemize}
	It is clear that the diagram commutes and that the rightmost square is Cartesian.  Moreover, one can check that these maps are all equivariant with respect to the following group actions. Set $U_{\dd}=\Hom(\BoC^{\dd},\BoC^{\vv})$. The group $\G_{\dd'}\times \G_{\dd''}\ltimes U_{\dd''}$ acts on all the four spaces in the first two rows of the diagram. For the four spaces at the corners of the left square, this action is induced by the projection The group $G_{\dd'+\dd''}\ltimes U_{\dd'+\dd''}$ acts on $\AA_{\dd'+\dd''}(Q^{\QM}_{\ww,\vv})$. The action of $P_{\dd',\dd''}\ltimes U_{\dd'+\dd''}$ on $\AA_{\dd',\dd''}(Q^{\QM}_{\ww,\vv})$ and $Z^{\vee}\times \AA_{\dd',\dd''}(Q^{\QM}_{\ww,\vv})$ is induced by the action of $G_{\dd'+\dd''}\ltimes U_{\dd'+\dd''}\supset P_{\dd',\dd''}\ltimes U_{\dd''}$ on $\AA_{\dd',\dd''}(Q^{\QM}_{\ww,\vv})$ and the embeddings $h\circ e$ and $e$, whose images are easily seen to be preserved. 
		Each of the spaces in the diagram carries a potential, assigned as follows:
	\begin{itemize}
		\item The potential on $\AA_{\dd'}(\tilde Q)\times \AA_{\dd''}(Q^{\QM}_{\ww,\vv})$ is $\tilde W_{\dd'}+ W_{\dd''}$, where $\tilde W_{\dd'}$ is the canonical cubic potential \eqref{CCP} on  $\AA_{\dd'}(\tilde Q)$ and $W_{\dd''}$ is the quasimap potential \eqref{eq: quasimap potential} on $\AA_{\dd''}(Q^{\QM}_{\ww,\vv})$.
		\item The potentials on $Z^{\vee}\times \AA_{\dd'}(\tilde Q)\times \AA_{\dd''}(Q^{\QM}_{\ww,\vv})$, $Z\times \AA_{\dd'}(\tilde Q)\times \AA_{\dd''}(Q^{\QM}_{\ww,\vv})$, and $\AA_{\dd',\dd''}(Q^{\QM}_{\ww,\vv})$  are pulled back from $\AA_{\dd'}(\tilde Q)\times \AA_{\dd''}(Q^{\QM}_{\ww,\vv})$ via $c$, $a$, and $a\circ g$, respectively.
		\item The potential on $Z^{\vee}\times Z  \times \AA_{\dd'}(\tilde Q)\times \AA_{\dd''}(Q^{\QM}_{\ww,\vv})$ is $\tilde W_{\dd'}+ W_{\dd''}+ \tr(lz)$, where $l$ and $z$ are coordinates on $Z^{\vee}$ and $Z$, respectively.
		\item The potential on $\AA_{\dd'+\dd''}(Q^{\QM}_{\ww,\vv})$ is the quasimap potential $W_{\dd'+\dd''}$.
		\item The potential on $Z^{\vee}\times \AA_{\dd',\dd''}(Q^{\QM}_{\ww,\vv})$ is the pullback of $W_{\dd'+\dd''}$ along $h$.
	\end{itemize}
	 The potentials commute with all the maps in the diagram and are invariant under the relevant group action.  Notice in particular that the equality
	\begin{equation*}
		(\tilde W_{\dd'}+ W_{\dd''}+ \tr(lz))\circ b=W_{\dd'+\dd''}\circ h
	\end{equation*}
	follows from the definition of $b$. 
    Taking quotients of \eqref{eq: big diagram dimred prequotients 1}, we get a diagram 
    \begin{equation*}
		\label{eq: big diagram dimred prequotients}
		\begin{tikzcd}[column sep=small]
			\FM_{\dd'}(\tilde Q)\times \FM^{\zeta^+\sst}_{\dd''}(Q^{\QM}_{\ww,\vv})  &\arrow[l, swap, "a"] \left(\FM_{\dd'}(\tilde Q)\times \FM^{\zeta^+\sst}_{\dd''}(Q^{\QM}_{\ww,\vv})\right)^{Z} \arrow[d, "d"] & \FM^{\zeta^+\sst}_{\dd',\dd''}(Q^{\QM}_{\ww,\vv}) \arrow[d, "e"]\arrow[l, swap, "g"]
			\\
		    \left(\FM_{\dd'}(\tilde Q)\times \FM^{\zeta^+\sst}_{\dd''}(Q^{\QM}_{\ww,\vv})\right)^{Z^{\vee}} \arrow[r, "f"] \arrow[u, "c"]& \left(\FM_{\dd'}(\tilde Q)\times \FM^{\zeta^+\sst}_{\dd''}(Q^{\QM}_{\ww,\vv})\right)^{Z^{\vee}\times Z} &\left(\FM^{\zeta^+\sst}_{\dd',\dd''}(Q^{\QM}_{\ww,\vv}) \right)^{Z^{\vee}}\arrow[d, "h"]\arrow[l, swap, "b"]
			\\
			& &\FM^{\zeta^+\sst}_{\dd'+\dd''}(Q^{\QM}_{\ww,\vv})
		\end{tikzcd}
	\end{equation*}
    where $\left(\FM_{\dd'}(\tilde Q)\times \FM^{\zeta^+\sst}_{\dd''}(Q^{\QM}_{\ww,\vv})\right)^{Z^{\vee}}$ stands for the tautological bundle on $\FM_{\dd'}(\tilde Q)\times \FM^{\zeta^+\sst}_{\dd''}(Q^{\QM}_{\ww,\vv})$ with fiber $Z^{\vee}$ and the same logic is used to defined all the spaces with similar superscripts. From now on the proof follows the exact same argument of Davison's dimensional reduction in \cite[App. 1]{RS17}, with diagram $(7)$ in \emph{loc. cit.} replaced by our diagram \eqref{eq: big diagram dimred prequotients}. 

\end{proof}


\newpage

\paperpart{Cech Model and Coulomb action}
\label{part: Čech model}

\section{The BFN Coulomb branch}

\subsection{Mapping stacks on the formal disk and related geometries}

\label{subsec: Mapping stacks on the formal disk and related geometries}
Let $\DD$ be the formal disk, and $N$ a complex representation of a reductive group $G$. The prototypical example is the mapping stack $\Map(\DD, N/G)=N\fps/G\fps$. 
Set $N^k\fps= N(\Sh{O})/t^k N(\Sh{O})$ and $G^k\fps= G(\Sh{O}/t^k\Sh{O})$. The mapping stack $\Map(\DD, N/G)$ is the moduli space parametrizing pairs $(\Sh{P},s)$ where $\Sh{P}$ is a principal $G$ bundle on the formal disk $\DD$ and $s$ is a section of the associated bundle $N\times_G \Sh{P}$. 
We have canonical morphisms $N^{k+l}\fps/G^{k+l}\fps\to N^k\fps/G^k\fps$ for all $k\geq 0$ and $\Map(\DD, N/G)$ is the inverse limit 
\[
\Map(\DD, N/G)=\varprojlim \Map(\DD_{k}, N/G)=\varprojlim N^k\fps/G^k\fps.
\]

Let now $\Hecke^{\pm}(\DD, N/G)$ be the stack of Hecke modifications of these data. Explicitly, $\Hecke^{\pm}(\DD, N/G)$ is the moduli of parametrizing tuples $(\Sh{P}_1, \Sh{P}_2, n_1, n_2, g)$ where $\Sh{P}_1, \Sh{P}_2$ are $G$-bundles on $\DD$, $g$ is an isomorphism $g: \Sh{P}_1|_{\DD^\times}\xrightarrow{\cong } \Sh{P'}_2|_{\DD^\times}$,  and $n_1$ and $n_2$ are sections of the associated bundles $N\times_G\Sh{P}_1$ and $N\times_G\Sh{P}_2$, respectively, such that $g n_1=n_2$. The moduli space $\Hecke^{\pm}(\DD, N/G)$ is a ind-stack of ind-infinite type. It admits a global quotient presentation
\[
\Hecke^{\pm}(\DD, N/G)=\Rt_{G,N}/G\fps ,
\]
where $\Rt_{G,N}$, sometimes referred to as the \emph{space of triples}, is the moduli space parametrizing the same data $(\Sh{P}_1, \Sh{P}_2,  n_1, g, n_2)$ as in $\Hecke^{\pm}(\DD, N/G)$ but with $\Sh{P}_2=G\fps$, i.e. with $\Sh{P}_2$ being the trivial $G$-torsor. Set $\Tt_{G,N}= N\fps \times_{G\fps}G\lfps$. In other words, $\Tt_{G,N}$ also parametrizes tuples $(\Sh{P}_1, \Sh{P}_2,  n_1, g, n_2)$ as above but now without the requirement that $g n_1=n_2$. Throughout this article, we will often drop $\Sh{P}_1$ and $\Sh{P}_2$ from the notation and denote a closed point in $\Rt_{G,N}$ by $(n_1, g, n_2)$. 

We have a closed embedding $\Rt_{G,N}\hookrightarrow \Tt_{G,N}$, as well as a map $\pi: \Tt_{G,N}\to \Gr_G$. For any $k> 0$, we introduce approximations 
\[
\pi_{k}: \Tt^d_{G,N}= (N\fps/t^kN\fps) \times_{G\fps} G\lfps \to \Gr_G
\]
and further set $\Tt^k_{N,G, \leq \lambda}=(\pi_d)^{-1}(\Gr_{\leq \lambda})$. These are schemes of finite type. We also let $\Tt_{N,G, \leq \lambda}=(\pi)^{-1}(\Gr_{\leq \lambda})$. If $\lambda\leq \mu$ with respect to the dominance order in the coweight lattice, we have inclusions 
\[
\Tt_{N,G, \leq \lambda}\hookrightarrow \Tt_{N,G, \leq \mu}.
\]
The definition of $\Rt^k_{N,G, \leq \lambda}$ and $\Rt_{N,G, \leq \lambda}$ is completely analogous, cf. \cite[\S 2(i)]{BFNII}. 

\subsection{The case of quiver representations}
\label{subsec: The case of quiver representations}
In this article we will be interested in the case where $G=\GL_{\vv}$ is a product of general linear groups and $N=\Rep_Q(\vv, \ww)$ is the moduli of framed representations of a quiver $Q$. In this case $\Map(\DD, N/G)$ is, more explicitly, the moduli stack of tuples $(\Sh{V}, n)$ consisting of rank $\vv_i$ vector bundles $\Sh{V_i}$ and a collection of morphisms coherent sheaves 
\begin{equation}
    \label{eq: ambient alpha}
    n=(n_a, n_i)_{a\in Q_1, i \in Q_0}\in \bigoplus_{a\in Q_1} \Hom(\Sh{V}_{s(a)}, \Sh{V}_{t(a)}) \oplus \bigoplus_{i \in Q_0} \text{Hom}(\Sh{V}_i,  W_i \otimes \Sh{O}_{\DD}).
\end{equation}
In this case $\Map(\DD, N/G)$ can also be seen a substack of the moduli stack $\FCoh_{Q^{\text{fr}}}(\DD)$ of objects in the abelian category $\Coh_{Q^{\text{fr}}}(\DD)$ of framed $Q$-valued sheaves on $\DD$. That is, the category whose objects are tuples of coherent sheaves $\{\Sh{V}_i\}_{i\in Q_0}$ on the disk and morphisms $\{n_a, n_i\}_{a\in Q_1, i \in Q_0}$ such that $n_a: \Sh{V}_{s(a)}\to \Sh{V}_{t(a)}$, $n_i: \Sh{V}_i \to W_i \otimes \Sh{O}_{\mathbb{D}} $. More precisely, $\Map(\DD, N/G)$ is the open substack of the connected component of $\FCoh_Q(\DD)$ parameterzing rank $\vv$ coherent sheaves with $\ww$-dimensional framing such that the coherent sheaves are locally free. 

The stack $\Hecke(\DD, N/G)$ is then the moduli stack of the data $(\Sh{V}_1, \Sh{V}_2, n_1, g, n_2)$, were $\Sh{V}$ and $\Sh{V}'$ are collections of vector bundles of the same rank, $g_i: \Sh{V}_{1,i}|_{\DD^\times }\xrightarrow{\cong }\Sh{V}'_{2,i}|_{\DD^\times }$, and $n_1$ and $n_2$ are sections as in \eqref{eq: ambient alpha} (with $\Sh{V}_2$ replacing $\Sh{V}_1$ in the case of $n_2$) satisfying $g_{t(a)}n_{1,a}=n_{2,a}g_{s(a)}$, $n_{1,i} = n_{2,i} g_i $ on $\DD^\times$. With slight abuse of notations, we will often denote a closed point in $\Hecke(\DD, N/G)$ as a triple $(n_1, g, n_2)$ , thus dropping the bundles $\Sh{V}_1$ and $\Sh{V}_2$ from the notation.

\subsection{The positive space of triples}
\label{subsec: The positive space of triples}

We retain the notation and the assumptions of the previous section. Additionally, set $G\lfps^+=G\lfps\cap \Fg\fps$ and $\Gr^+=G\fps\backslash G\lfps^+$. We have a canonical inclusion $\Gr^+\subset \Gr$.
The positive affine Grassmannian $\Gr^+$ parametrizes tuples $\Sh{V}_i$ of locally free sheaves rank $\vv_i$ on the disk $\DD$ together with injective morphism of sheaves $\Sh{V}_i \hookrightarrow \Sh{O}^{\vv_i}_{\DD}$. It decomposes in connected components $\Gr^+=\bigsqcup_{\dd\in \BoN^{Q_0}} \Gr^+_{\dd}$ according to the lengths $\dd_i$ of the quotients $\Sh{O}^{\vv_i}_{\DD}/\Sh{V_i}$.
We further set 
\begin{align*}
    \Tt^+&=N\fps\times_{G\fps} G\lfps ^+ 
    \\
    \Rt^+&=\Tt^+\cap (\Gr^+\times N\fps)
\end{align*}
Both $\Tt$ and $\Rt$ map canonically to $\Gr^+$ and hence inherit decompositions 
\[
\Tt^+=\bigsqcup_{\dd\in \BoN^{Q_0}} \Tt^+_{\dd} \qquad \Rt^+=\bigsqcup_{\dd\in \BoN^{Q_0}} \Rt^+_{\dd}.
\]
Overall, we can assemble these spaces into the diagram
\[
\begin{tikzcd}
    \Rt^+\arrow[r, hookrightarrow]\arrow[d, hookrightarrow, swap, "\rho"] & \Tt^+ \arrow[r, hookrightarrow] \arrow[d, hookrightarrow]& \Gr^+\times N\lfps\arrow[d, hookrightarrow]
    \\
    \Rt\arrow[r, hookrightarrow] & \Tt \arrow[r, hookrightarrow] & \Gr\times N\lfps
\end{tikzcd}
\]
All the squares are Cartesian, all vertical morphisms are closed, and all the maps are equivariant with respect to the left action of $G\fps$. Moreover, the quotient $G\fps\backslash \Rt^+$ is canonically isomorphic to the stack $\Hecke(\DD, N/G)$, which parametrizes tuples $(g, n_1, n_2)$, where 
\[
g_i: \Sh{V}_{1,i}\hookrightarrow\Sh{V}_{2,i}
\]
is an injective morphism of locally free sheaves of rank $\vv_i$
while $n_{1,a}: \Sh{V}_{1, s(a)}\to \Sh{V}_{1, t(a)}$ and  $n_{2,a}: \Sh{V}_{2,s(a)}\to \Sh{V}_{2, t(a)}$ (and similarly $n_{1,i}, n_{2,i}$) are morphisms such that the following squares commute: 
\begin{equation}
    \label{eq: compatibility positive Coulomb} 
\begin{tikzcd}
    0\arrow[r] &\Sh{V}_{1, s(a)} \arrow[r, "g_{s(a)}"]\arrow[d, "n_{1,a}"] & \Sh{V}_{2,s(a)}\arrow[d, "n_{2,a}"]
    \\
        0\arrow[r] & \Sh{V}_{1, t(a)} \arrow[r, "g_{t(a)}"] & \Sh{V}_{2,t(a)}
\end{tikzcd}
\qquad 
\begin{tikzcd}
    0\arrow[r] &\Sh{V}_{1, i} \arrow[r, "g_{i}"]\arrow[d, "n_{1,i}"] & \Sh{V}_{2,i}\arrow[d, "n_{2,i}"]
    \\
        0\arrow[r] & \Sh{W}_{i} \arrow[r, equal] & \Sh{W}_{i}
\end{tikzcd}
\end{equation}
It is clear that $\Hecke(\DD, N/G)$ decomposes in connected components
\[
\Hecke(\DD, N/G)=\bigsqcup_{\dd\in \BoN^{Q_0}} \Hecke_{\dd}(\DD, N/G).
\]
according to the lengths $\dim(\Sh{V}'_i/\Sh{V}_i)=\dd_i$ for all $i\in Q_0$, and that $\Hecke_{\dd}(\DD, N/G)=\Rt^+ /G\fps$. Moreover, if $\dd=0$ then $\Gr_{\dd}=\Spec(\BoC)$ and $\Hecke_{0}(\DD, N/G)=\Map(\DD, N/G)$.
Finally, notice that we have canonical morphisms 
\begin{equation}
    \label{eq: diagram various maps from hecke on disk}
    \begin{tikzcd}
    & \Hecke_{\dd}(\DD, N/G)\arrow[dl, swap, "p_1"]\arrow[d, "p_2"] \arrow[dr, "p_3"]
    \\
    \Map(\DD, N/G)  & \Map(\DD, N/G)  & \FCoh_{Q,0,\dd}(\DD)
\end{tikzcd}
\end{equation}
sending a tuple $(n_1, g, n_2)$ to $n_1$, $n_2$, and $\beta$, respectively. Here, $\beta$ is the unique tuple of morphism fitting in the diagram 
\begin{equation*}
    \begin{tikzcd}
    0\arrow[r]& \Sh{V}_{1,s(a)} \arrow[r, "g_{s(a)}"]\arrow[d, "n_{1,a}"] & \Sh{V}_{2, s(a)} \arrow[r]\arrow[d, "n_{2,a}"]& \Sh{V}_{2, s(a)}/\Sh{V}_{1, s(a)}\arrow[d, "\beta_a"] \arrow[r] & 0
    \\
    0\arrow[r]&\Sh{V}_{1, t(a)} \arrow[r, "g_{t(a)}"] & \Sh{V}_{2, t(a)}\arrow[r] & \Sh{V}_{2,t(a)}/\Sh{V}_{1,t(a)} \arrow[r] &0
\end{tikzcd}
\end{equation*}
.

\begin{remark}
\label{rem: retraction positive Coulomb 1}
    It is worth noticing that the moduli space $\Rt^+_{G,N}$ is a vector bundle over, in particular contractible to, $\Rt^+_{G,N_0}$, where $N_0\subset N$ is the $G$-subrepresentation that is independent on the framing. Equivalently, $N_0$ is determined by the decomposition $N=N_0\oplus\bigoplus_{i\in Q_0}\Hom(V_i, W_i)$. To see this, observe that for any $i\in Q_0$, the right square in \eqref{eq: compatibility positive Coulomb} implies that $n_{1,i}$ is uniquely determined by $g_i$ and $n_{2,i}$ for all $i$. Moreover, since $g$ is defined on the whole disk $\DD$, the equation $n_{1,i}=n_{2,i}g_i$ imposes no constraint on $n_{2,i}$. Therefore, the canonical morphism $\Rt^+_{G,N}\to \Rt^+_{G,N_0}$ forgetting $n_{2,i}$ and $n_{1,i}$ for all $i\in Q_0$ is a vector bundle projection. We remark that this retraction does not extend over $\Rt_{G, N}$ because for a general $(n_{1,i},g, n_{2,i})\in \Rt_{G,N}$ the Hecke modification $g$ is only defined on the punctured disk $\DD^\times$, so the the equations $n_{1,i}=n_{2,i}g_i$ imposes constraints on $n_{2,i}$: the compositions $n_{2,i}g_i$ must be pole-free for all $i\in Q_0$.
\end{remark}

\subsection{Cohomology of pro and ind-stacks}
\label{subsec: cohomology of pro and ind stacks}
In this article, we will need to consider cohomology of certain pro-stacks and ind-stacks. We review the definition here. The easiest example is the mapping stack $\Map(\DD, N/G)$. The morphisms $N^{k+l}\fps/G^{k+l}\fps\to N^k\fps/G^k\fps$
are locally trivial fibrations of relative dimension $l(\dim(N)-\dim(G))$ so they induce a canonical pullback morphism in Borel-Moore homology 
\begin{equation}
\label{eq: cohomology mappisng stack f im approx comparison}
    \HO^{\BM}(N^k\fps/G^k\fps, \BoQ)\xrightarrow{\cong} \HO^{\BM}(N^{k+l}\fps/G^{k+l}\fps, \BoQ)[-2l(\dim(N)-\dim(G))]
\end{equation}
which is actually an isomorphism. Hence, up to a degree shift, the homologies 
\[
\HO^{\BM}(N^k\fps/G^k\fps, \BoQ)
\]
are all canonically isomorphic. Therefore, one may define the BM homology of $\Map(\DD, N/G)$
as the BM homology of any of the finite dimensional approximations $N^k\fps/G^k\fps$. More precisely, in order to have a well defined notion of cohomological degree, we set 
\begin{equation}
    \label{eq: def cohomology of Mapping stck disk}
    \HO^{\BM}(\Map(\DD, N/G), \BoQ)=\HO^{\BM}(N^k\fps/G^k\fps, \BoQ)[-2k(\dim(N)-\dim(G))].
\end{equation}
for some $k\geq 0$. Equation \eqref{eq: cohomology mappisng stack f im approx comparison} implies that the definition, including the notion of cohomological degree, does not depend on $k$. Notice that although in this case any $k\geq 0$ works, later well have to impose some lower bounds on $k$.

The definition of the BM homology of $\HO^{\BM}_{G\fps}( \Tt_{G,N}, \BoQ)$ is similar. Firstly we define the BM homology of $\Tt_{N,G, \leq \lambda}$ by noticing that we have locally trivial fibrations $G^{k+l}\fps\backslash \Tt^{l+k}_{N,G, \leq \lambda}\to G^{k}\fps\backslash \Tt^{l}_{N,G, \leq \lambda}$ inducing isomorphisms 
\begin{equation}
    \label{eq: shift dimension for different approximations BFN}
    \HO^{\BM}_{G^k\fps}(\Tt^{k}_{N,G, \leq \lambda}, \BoQ)\to \HO^{\BM}_{G^{k+l}\fps}(\Tt^{l+k}_{N,G, \leq \lambda}, \BoQ)[-2l(\dim(N)-\dim(G))]
\end{equation}
which can be bootstrapped to define 
\begin{equation}
    \label{eq: def cohomology of variety Tt}
\HO^{\BM}_{G\fps}( \Tt_{N,G, \leq \lambda}, \BoQ)\coloneqq\HO^{\BM}_{G^k\fps}(\Tt^{k}_{N,G, \leq \lambda}, \BoQ)[-2k(\dim(N)-\dim(G))]
\end{equation}
for some large enough $k$. The independence of $k$ follows from \eqref{eq: shift dimension for different approximations BFN}.
Finally, we define the cohomology 
\[
\HO^{\BM}_{G\fps}( \Tt_{G,N}, \BoQ)=\varinjlim \HO^{\BM}_{G\fps}(\Tt_{N,G, \leq \lambda}, \BoQ).
\]
via the direct limit (for the direct system associated to the dominance order on the co-weight lattice)
\[
\HO^{\BM}_{G\fps}( \Tt_{N,G, \leq \lambda}, \BoQ)\to \HO^{\BM}_{G\fps}( \Tt_{N,G, \leq \mu}, \BoQ)
\]
Here, the morphisms are proper pushforwards associated to the morphisms $ \Tt_{N,G, \leq \lambda}\to  \Tt_{N,G, \leq \mu}$.

The definition of $\HO^{\BM}_{ G\fps}( \Rt_{G,N}, \BoQ)$
is completely similar, and we refer to \cite[\S 2(ii)]{BFNII} for full details. 

Throughout the article, we will also consider some slight variants of these BM homologies, e.g. finite rank vector bundles on $\Tt_{G,N}$ or, more generally, on fiber products of a Artin stack and  $\Tt_{G,N}$. In all such cases, the cohomology of these stacks is defined via a straightforward adaptation of the above constructions. The $G\fps \times G\fps$-equivariant Borel-Moore homology of $G\lfps$ is also considered, and we refer to \cite[\S2(ii)]{BFNII} for a detailed discussion. of its definition. Although we will not spell always spell out all the details of the finite dimensional approximations, we will make sure to highlight the key ideas and potential subtleties.


\subsection{Quantization and Flavor Symmetry} 
\label{subsec: Quantization and Flavor Symmetry}
The space of triples $\Rt_{G,N}$ there admits a natural action of the semidirect product $G\fps\rtimes \BoC^*_{\hbar}$ induced by the natural action of $\BoC^*_{\hbar}$ on the disk $\DD$. The cohomology $\HO^{\BM}_{G\fps\rtimes \BoC^*_{\hbar}}( \Rt_{G,N}, \BoQ)$ is known as the quantized Coulomb branch\footnote{Its name comes  from the fact that, unlike $\HO^{\BM}_{G\fps}( \Rt_{G,N}, \BoQ)$, its convolution algebra structure is non-commutative.}. Additionally, one can introduce further equivariance by considering the action of a ``flavor'' torus $A$ acting on $N$ and hence on $\Rt_{G,N}$. In this article, we will always fully equivariantly with respect to a torus $T=A\times \BoC^*_{\hbar}$ and hence consider the cohomology 
\begin{equation}
    \label{eq: coulomb branch with flavor symmetry}
    \HO^{\BM}_{G\fps\rtimes T}( \Rt_{G,N}, \BoQ)
\end{equation}
If $G$ and $N$ are of quiver type, i.e. we have $N=\Rep_{Q}(\ww, \vv)$ and $G=\prod_{i\in Q_0}\GL(\vv_i)$, we take $T$ to be the torus defined in \S\ref{subsec: torus action quiver varieties}. Notice that although its action was defined on $T^*N$, it preserves the subspace $N$.

To remove clutter, we set $G\fps_T\coloneqq G\fps\rtimes T$ and hence write $\HO^{\BM}_{G\fps_T}( \Rt_{G,N}, \BoQ)$ as a shortcut for \eqref{eq: coulomb branch with flavor symmetry}. Throughout this article, we will only consider the quantized flavored-deformed Coulomb branch $\HO^{\BM}_{G\fps_T}( \Rt_{G,N}, \BoQ)$. Although all our constructions are $T$-equivariant, we will only factor out $T$ whenever convenient for the exposition.

\subsection{The BFN Coulomb branch}
\label{subsec: Coulomb branch basics}

In this section we review the construction of the Coulomb branch developed by Braverman, Finkelberg, and Nakajima in \cite{BFNII}. To remove clutter, we set $\Rt=\Rt_{G,N}$ and $\Tt=\Tt_{G,N}$. The quantized BFN Coulomb branch is, by definition, a convolution algebra whose underlying vector space is $\HO^{\BM}_{G\fps_T}(\Rt, \BoQ)$. The algebra structure is given by a convolution product whose construction we now come to review.

The stack $\Tt$ admits a canonical embedding $\Tt=  N\fps\times_{G\fps}G\lfps\to \Gr\times N\lfps$ induced by $[g,n]\mapsto ([g], gn)$. In this ambient space, we have $\Rt=\Tt\cap (\Gr\times N\fps)$. 

Following \cite{BFNII}, we consider the diagram 
\begin{equation}
    \label{eq: Coulomb branch mult diagram}
\begin{tikzcd}
    \Rt\times \Rt \arrow[d, hookrightarrow]& p^{-1}(\Rt\times \Rt) \arrow[d, hookrightarrow]\arrow[l, swap, "\tilde p"] \arrow[r, "\tilde q"] & q(p^{-1}(\Rt\times \Rt))\arrow[d, hookrightarrow]\arrow[r, "\tilde m"] & \Rt\arrow[d, "i", hookrightarrow]
    \\
    \Rt\times \Tt & \Rt\times G\lfps \arrow[l, swap, " p"] \arrow[r, "q"] & (\Rt\times_{G\fps} G\lfps)\arrow[r, " m"] & \Tt
\end{tikzcd}
\end{equation}
The maps in the second row are defined by 
\begin{equation}
    \label{eq: maps correspondence Coulomb product}
    ([n,g_1], [g_1n,g_2])\mapsfrom ([n,g_1], g_2)\mapsto [[n,g_1],g_2]\mapsto [n,g_1g_2],
\end{equation}
and those on the top row are defined by restriction as all vertical maps are injective.  
\begin{remark}
   Compared with the BFN convolution diagram \cite[(3.2)]{BFNII}, the order of all the factors in the fiber product in \eqref{eq: Coulomb branch mult diagram} appears to be reversed. This is due to our choice to present the affine Grassmannian $\Gr$ as a left quotient, $\Gr=G\fps\backslash G\lfps$, rather than as a right quotient, as in BFN's convention. Our choice is determined by the requirement that the Coulomb branch act on critical quasimaps from the left; with BFN's convention, we would instead obtain a right action.
\end{remark}

The following lemma follows easily by analyzing the topology of $p$. To remove clutter, we introduce the notation $\infty_{G,N}=\dim N\fps-\dim G\fps$. We will also suppress the action of $T$ from the notation.

\begin{lemma}[{\cite[Lemma 3.5]{BFNII}}]
\label{lemma: Lemma pullback BFN}
    There is a canonical isomorphism of $G\fps\times G\fps$-equivariant sheaves
    \begin{align*}
        p^*\DD\BoQ_{\Rt \times \Tt}\cong p^*\DD\BoQ_{\Rt\times G\lfps}[2\infty_{G,N}]
    \end{align*}
\end{lemma}
Here the degree shift $2\infty_{G,N}$ is understood by taking finite dimensional approximation as in \eqref{subsec: cohomology of pro and ind stacks}.

Applying the functor $(\text{id} \times i)^!$ to the unit $\id \to p_*p^* \DD\BoQ_{\Rt \times \Tt}$ and performing base change we get a canonical Gysin morphism 
\begin{equation}
    \label{eq: gysin pullback Coulomb}
    \DD\BoQ_{\Rt\times \Rt}\to \tilde p_{*}\DD\BoQ_{p^{-1}(\Rt\times\Rt)}[2\infty_{G,N}]
\end{equation}
and hence a map in cohomology 
\[
\tilde p^!: \HO_{G\fps}(\Rt, \DD\BoQ_\Rt)\to  \HO_{G\fps}(p^{-1}(\Rt\times \Rt), \DD\BoQ_{p^{-1}(\Rt\times \Rt)}).
\]
The convolution algebra structure on the Coulomb branch $\HO_{G\fps}(\Rt, \DD\BoQ_\Rt)$ is then defined as the composition 
\begin{align}
\label{eq: Coulomb branch mult}
\begin{split}
    \HO_{G\fps}(\Rt, \DD\BoQ_\Rt)\otimes \HO_{G\fps}(\Rt, \DD\BoQ_\Rt) 
    & \xrightarrow{\cong} \HO_{G\fps\times G\fps}(\Rt\times \Rt, \DD\BoQ_{\Rt\times \Rt}) 
    \\
    &\xrightarrow{\tilde p^!}  \HO_{G\fps\times G\fps}(p^{-1}(\Rt\times \Rt), \DD\BoQ_{p^{-1}(\Rt\times \Rt)})
    \\
    &\xrightarrow[\tilde q^*]{\cong}\HO_{ G\fps}(q(p^{-1}(\Rt\times \Rt)), \DD\BoQ_{q(p^{-1}(\Rt\times \Rt))})
    \\
    &\xrightarrow{m_*} \HO_{G\fps}(\Rt, \DD\BoQ_\Rt)
\end{split}
\end{align}
where the last map is well-defined because it is ind-proper.
Following \cite{BFNII}, the (quantized, flavor deformed) BFN Coulomb branch is the associative algebra 
\[
\Coul^T_{G,N}=\HO^{\BM}_{ G\fps_T}(\Rt_{G,N}, \BoQ)
\]
with convolution product defined as in \eqref{eq: Coulomb branch mult}. The unit is given by the fundamental class of the fiber $\Rt^+\to \Gr^+$ at the base point $[1]\in \Gr^+$, corresponding to the trivial Hecke modification.

Recall that the affine Grassmannian $\Gr_G$ has connected components indexed by $\pi_1(G)$. Accordingly, the space $\Rt$ decomposes as $\Rt=\sqcup_{\gamma\in \pi_1(G)}\Rt^{\gamma}$. In terms of the cell decomposition $\Gr_G=\bigsqcup_{\lambda\in X_*(S)_+} \Gr_G^{\lambda}$, the canonical map $\Gr_G\to \pi_1(G)$ sends a point $x\in \Gr_G^{\lambda}$ to the coset $[\lambda]\in\pi_1(G)=X_*(S)/\{\text{coroots}\}$. Accordingly, the cohomology $\HO^{\BM}_{G\fps}(\Rt, \BoQ)$ inherits a $\pi_1(G)$ grading
\[
\HO^{\BM}_{G\fps_T}(\Rt, \BoQ)=\bigoplus_{\gamma\in \pi_1(G)} \HO^{\BM}_{G\fps_T}(\Rt_{\gamma}, \BoQ)
\] 
and it is easy to check that the multiplication $\eta$ is compatible with this grading. Therefore, the pair $(\HO^{\BM}_{G\fps_T}(\Rt,\BoQ), \eta)$ is a cohomologically graded, $\BoZ^{Q_0}$-graded algebra.

\begin{remark}
In the case of a quiver gauge theory, we have a natural embedding
    \[
    \pi_1(G)=\pi_1\left(\prod_{i\in Q_0}\GL(\vv_i)\right)\subseteq \BoZ^{Q_0},
    \]
which is an equality unless $\vv_i=0$ for some $i\in Q_0$. In the case of a strict embedding, we identify $\pi_1(G)$ with the subgroup of $\BoZ^{Q_0}$ consisting of those vectors $\dd\in \BoZ^{Q_0}$ such that $\dd_i=0$. 

In \S\ref{subsec: From hall to Coulomb}, we will match this $\BoZ^{Q_0}$ grading with the $\BoN^{Q_0}$ grading on the Hall algebra $\HCoha^{T}_{\tilde Q, \tilde W}$ and its variants.
\end{remark}

\begin{remark}
\label{rem: virtual Coulomb branch}
    Assume now that the Coulomb branch is of quiver type, i.e. $N=\Rep_Q(\ww,\vv)$ and $G=\prod_{i\in Q_0}\GL(\dd_i)$. To relate such Coulomb branches to shifted Yangians, it is convenient to shift the cohomological degree by setting 
    \[
    \Coul^{T,\vir}_{G,N}\coloneqq\bigoplus_{\dd\in \BoN^{Q_0}}\HO^{\BM}_{G\fps_T}(\Rt_{N,G,\dd}, \BoQ[2s_{\dd}]),
    \]
    where $s_{\dd}=\vv^T\dd-\vv^T\Q^T\dd$. Notice that, since the function $\dd\mapsto s_{\dd}$ is linear, the shift is compatible with the cohomologically graded associative algebra structure. In other words, the multiplication 
    \[
    \HO^{\BM}_{G\fps_T}(\Rt_{N,G, \dd'}, \BoQ[s_{\dd'}])\otimes_{\HO_{T}} \HO^{\BM}_{G\fps_T}(\Rt_{N,G, \dd''}, \BoQ[s_{\dd''}])\to \HO^{\BM}_{G\fps_T}( \Rt_{N,G, \dd+\dd''} \BoQ[s_{\dd'+\dd''}])
    \]
    preserves the cohomological degree.
\end{remark}

\subsection{Coulomb Branch generators}
\label{subsec: coulomb branch generators}
We conclude the review of the BFN Coulomb Branch by discussing its generators in the quiver case. Therefore, we assume $N=\Rep_Q(\vv,\ww)$ and $G=\prod_{i\in Q_0}\GL(\vv_i)$. For a given $i\in Q_0$ and $1\leq n\leq \vv_i$. let $\omega_{i,n}$ be the $i$-th fundamental weight of the $\GL(\vv_i)$-factor in the gauge group $G$ and let $\omega_{i,n}^*=-w_0\cdot \omega_{i,n}$. Here $w_0$ is the longest element in the Weyl group of $\GL(V_i)$. In the standard coordinate basis of $\GL(\vv_i)$, we have 
\[
\omega_{i,n}=(\underbrace{1,\dots, 1}_{\text{n times}}, 0,\dots, 0)\qquad \omega^*_{i,n}=-(0,\dots, 0, \underbrace{1,\dots, 1}_{\text{n times}}).
\]
These are the minuscule weights of $\GL(\vv_i)$ and the minuscule weights of $G=\prod_{i\in Q_0}{\GL(\vv_i)}$ are given by a choice of minuscule weight for each factor in the product. 
Let $\lambda$ be a minuscule weight of $G$. Then the corresponding $G\fps$-orbit $\Gr_{G,\lambda}\subset \Gr_{G}$ is closed. Overall, $\Rt_{\lambda}$ is closed and $G\fps$-stable in $\Rt_{G}$, so it induces a fundamental class $[\Rt_{\lambda}]$ in the Borel-Moore homology $\HO_{G\fps_T}(\Rt, \DD\BoQ_\Rt)$. Similarly,  one has \cite[\S 6]{BFNII}
\[
\HO_{G\fps_T}(\Rt_{\lambda}, \DD\BoQ_\Rt)= \HO_{\Stab_G(\lambda)\times T}\cap [\Rt_{\lambda}].
\]
Given $f\in \HO_{\Stab_G(\lambda)}$ we denote the corresponding class by
\[
M_{f, \lambda}\coloneqq f\cap [\Rt_{\lambda}]\in \HO_{G\fps_T}(\Rt, \DD\BoQ_\Rt).
\]
These classes are known as \emph{dressed monopoles operators}, where dressing refers to the factor $f$. It is known that the Coulomb branch is generated as a $\HO_{T}$-algebra by the dressed monopoles operators \cite[Remark 6.7]{BFNII}\cite[Thm. 4.32]{FT_ShitedQuantum}\cite[Prop. 3.1]{Weekes}. In the paper we will need the following variant of this result by Weekes.
\begin{theorem}[{\cite[Thm. 3.7 and Rem. 3.9]{Weekes}}]
\label{eq: generators Coulomb}
    The localized Coulomb branch $\Coul^T_{N,G}\otimes_{\HO_{\BoC^*_{\hbar}}}\Frac(\HO_{\BoC^*_{\hbar}})$ is generated as a $\HO_{G\fps_T}\otimes_{\HO_{\BoC^*_{\hbar}}}\Frac(\HO_{\BoC^*_{\hbar}})$-algebra by the monopoles operators $M_{1, \omega_{i,1}}$ and $M_{1, \omega^*_{i,1}}$.
\end{theorem}

\subsection{The positive part of the Coulomb branch}
\label{subsec: Positive part of the Coulomb branch}

Recall the definition of the positive space of triples $\Rt^+$ from \S\ref{subsec: The positive space of triples}. In this section we argue that
\begin{equation}
    \label{eq: positive Coulomb}
    \HO^{\BM}_{G\fps_T}(\Rt^+, \BoQ)
\end{equation}
admits a natural convolution algebra structure, and that the canonical pushforward 
\[
\rho_*: \HO^{\BM}_{G\fps_T}(\Rt^+, \BoQ)\to \HO^{\BM}_{G\fps_T}(\Rt, \BoQ)
\]
is an injective morphism of algebras, where the right hand side is the Coulomb branch algebra by \cite{BFNII}. This construction was originally outlined in \cite[\S3(ii)]{BFN_quiver} but, since we did not find a full proof of these statements, we outline it here. The construction itself will be relevant for the proof of Theorem \ref{thmx: map hall Coulomb}. We begin our analysis by listing some properties of \eqref{eq: positive Coulomb}. We drop the action of the torus $T$ from the notation.

\begin{proposition}
\label{prop: first properties positive Coulomb}
    Let $S\subset G$ be a maximal torus. 
    \begin{enumerate}
        \item The cohomologically graded vector space \eqref{eq: positive Coulomb} is concentrated in even degree and flat as a $\HO_G$-module.
        \item The base change morphism 
        \[
        \HO_G\otimes_{\HO_S}\HO_{S\fps}(\Rt^+, \BoD\BoQ_\Rt^+)\to \HO_{G\fps}(\Rt^+, \BoD\BoQ_\Rt^+)
        \]
        is an isomorphism. In particular $\left(\HO_{S\fps}(\Rt^+, \BoD\BoQ_\Rt^+)\right)^W\cong \HO_{G\fps}(\Rt^+, \BoD\BoQ_\Rt^+)$.
    \end{enumerate}
\end{proposition}
\begin{proof}
    The proof is an adaptation of \cite[Lemma 5.3]{BFNII}, which in turn is based on \cite{BFM}. We sketch how to adapt their proof here. It is well known that $\Gr$ has a decomposition into $G\fps$-orbits $\Gr=\sqcup_{\lambda\in X_*(S)^+}\Gr^{\lambda}$ where $\Gr^{\lambda}$ is the $G(\Sh{O})$-orbit of $t^\lambda$ where $\lambda$ is a dominant coweight. Each $\Gr^\lambda$ is an affine fibration over a partial flag variety for $G$, thus the equivariant cohomology of $\Gr^\lambda$ is even and free over $\HO_G(\text{pt})$. Similarly, the positive affine Grassmannian $\Gr^+$ admits a decomposition $\Gr^+=\sqcup_{\lambda\in \Gamma} \Gr^{\lambda}$ where $\Gamma\subset X_*(S)^+$ is a certain subset of dominant coweights. 

    Let now $\Rt^{\lambda}=\pi^{-1}(\Gr^{\lambda})$, where $\pi: \Rt\to \Gr$ is the canonical map. Then we get $\Rt^+=\sqcup_{\lambda\in \Gamma} \Rt^{\lambda}$.
    Moreover, by \cite[Lemma 2.2]{BFNII}, the map $\Rt^{\lambda}\to \Gr^\lambda$ is an affine bundle (of infinite rank). Conclude by the above that $\HO_{G\fps}(\Rt^{\lambda}, \DD\BoQ_{\lambda})$ is also even and free over $\HO_G(\text{pt})$.

    Consider the partial order on the strata $\{\Gr_{\lambda}\}_{\lambda \in \Gamma}$ obtained by taking the transitive closure of the binary relation $\ol{\Gr^{\lambda}}\cap \Gr^{\mu}\neq \emptyset \implies \Gr^{\lambda}\geq \Gr^{\mu}$. Here the closure is taken in $\Gr^+$. It is refined by the dominance order of coweights. Therefore, $\Gr^{+, \leq\lambda}=\sqcup_{\mu\leq \lambda } \Gr_{\lambda}$ and contains $\Gr^{\lambda}$ as an open subset. Let $\Rt^{+, \leq \lambda}=\pi^{-1}(\Gr^{+, \leq\lambda})$ and let $\Gr^{+, <\lambda}$ denote the complement of $\Rt^{\lambda}\hookrightarrow \Rt^{+, \leq\lambda}$. Since $\Rt^{\lambda}$ has no odd cohomology, a standard inductive argument with the long exact sequence for the open-closed embeddings $\Rt^{\lambda} \hookrightarrow \Rt^{+, \leq\lambda}\hookleftarrow \Rt^{+, <\lambda}$ shows that $\Rt^{+, \leq \lambda}$ has no odd cohomology too. Since the cohomology $\HO_{G\fps}(\Rt^+, \DD\BoQ_{\Rt^+})$ is defined as the direct limit of the cohomologies $\Rt^{+, \leq \lambda}$, it follows that the latter has no odd cohomology too. Moreover, since $\HO(\Rt^{+, \leq\lambda}, \DD\BoQ_{\Rt^{+, \leq\lambda}})$ is even, the previous long exact sequence of $\HO_G$-modules splits:
    \[
    0\to \HO_{G\fps}(\Rt^{+, <\lambda}, \DD\BoQ_{\Rt^{+, <\lambda}})\to \HO_{G\fps}(\Rt^{+, \leq\lambda}, \DD\BoQ_{\Rt^{+, \leq\lambda}})\to \HO_{G\fps}(\Rt^{\lambda}, \DD\BoQ_{\Rt^{\lambda}})\to 0
    \]
    Since flatness is preserved by taking extensions,  another inductive argument with respect to the dominance order shows that $\HO_{G\fps}(\Rt^{+, \leq\lambda}, \DD\BoQ_{\Rt^{+, \leq\lambda}})$ is $\HO_G$-flat, and hence the same is true for the limit $\HO_{G\fps}(\Rt^{+, \leq\lambda}, \DD\BoQ_{\Rt^{+, \leq\lambda}})$. This concludes the proof of the first point of the lemma. 
    The proof of the second part follows the same logic, and is proved by induction using the short exact sequence above. The details are left to the reader.
    \end{proof}

We now equip \eqref{eq: positive Coulomb} with an algebra structure. Notice that $\rho$ is ind-proper so the pushforward map $\rho_*$ is well-defined. We now upgrade diagram \eqref{eq: Coulomb branch mult diagram} as follows:
\begin{equation}
    \label{eq: huge diagram Coulomb}
\begin{tikzcd}
    & \Rt\times \Rt \arrow[dd, hookrightarrow, near start, swap, "i\times \id"]& p^{-1}(\Rt\times \Rt) \arrow[dd, hookrightarrow]\arrow[l, swap, "\tilde p"] \arrow[r, "\tilde q"] & q(p^{-1}(\Rt\times \Rt))\arrow[dd, hookrightarrow]\arrow[r, "\tilde m"] & \Rt\arrow[dd, hookrightarrow, "i"]
    \\
    \Rt^+\times \Rt^+\arrow[ur, hookrightarrow, "\rho\times \rho"] \arrow[dd, hookrightarrow, "i^+\times \id"]& (p^+)^{-1}(\Rt^+\times \Rt^+) \arrow[dd, shift left=5pt, hookrightarrow]\arrow[l, swap, "\tilde p^+"] \arrow[r, "\tilde q^+"] \arrow[ur, hookrightarrow]& q^+((p^+)^{-1}(\Rt^+\times \Rt^+))\arrow[dd, shift left=5pt, hookrightarrow]\arrow[r, "\tilde m"] \arrow[ur, hookrightarrow]& \Rt^+\arrow[dd, shift left=5pt, hookrightarrow, near start, "i^+"]\arrow[ur, hookrightarrow, "\rho"]
    \\
    & \Rt\times \Tt & \Rt\times G\lfps \arrow[l, swap, " p"] \arrow[r, "q"] &(\Rt\times_{G\fps}G\lfps)\arrow[r, " m"] & \Tt
    \\
    \Rt^+\times \Tt^+ \arrow[ur, hookrightarrow]& \Rt^+ \times G\lfps^+\arrow[l, swap, " p^+"] \arrow[r, "q^+"]\arrow[ur, hookrightarrow]& (\Rt^+\times_{G\fps}G\lfps^+)\arrow[r, " m^+"] \arrow[ur, hookrightarrow]& \Tt^+\arrow[ur, hookrightarrow]
\end{tikzcd}
\end{equation}
The spaces in the bottom row are all embedded in their analogs without superscript $(-)^+$, so the maps that relate them are formally defined as in \eqref{eq: maps correspondence Coulomb product}. It is easy to check that the diagram is commutative and all the faces in the two leftmost cubes of the diagram are Cartesian. Informally, this is achieved because all the spaces with the superscript $(-)^+$ are the subset of their analogs without the superscript obtained by enforcing that the Hecke modifications encoded in their data are actually embeddings.

To adapt the construction of \S\ref{subsec: Coulomb branch basics} to obtain a compatible algebra structure on the positive half $\HO_{G\fps}(\Rt^+, \DD\BoQ_{\Rt^+})$ is now easy. Firstly, by applying the functor $(\rho\times \rho)^!$ to ]\eqref{eq: gysin pullback Coulomb} and base changing with respect to the square with maps $(\rho\times \rho, \tilde p^+, \tilde p),$ we obtain a canonical morphism 
\begin{equation}
\label{eq: gysin pullback Coulomb +}
    (\tilde p^+)^!: \DD\BoQ_{\Rt^+\times\Rt^+}\to (\tilde p^+)_{*}\DD\BoQ_{(p^+)^{-1}(\Rt^+\times\Rt^+)}[2\infty_{G,N}]
\end{equation}
Therefore, we can define a multiplication map on $\HO_{G\fps}(\Rt^+, \DD\BoQ_{\Rt^+})$ by setting 
\begin{align}
\begin{split}
    \label{eq: gysin pullback Coulomb + 2}
    \HO_{G\fps}(\Rt^+, \DD\BoQ_{\Rt^+})\otimes \HO_{G\fps}(\Rt^+, \DD\BoQ_{\Rt^+})
    &\cong\HO_{G\fps\times G\fps}(\Rt^+\times \Rt^+, \DD\BoQ_{\Rt^+\times \Rt^+}) 
    \\
    &\xrightarrow{(\tilde p^+)^!}  \HO_{G\fps\times G\fps}(\tilde p^{-1}(\Rt^+\times \Rt^+), \DD\BoQ_{p^{-1}(\Rt^+\times \Rt^+)})
    \\
    &\xrightarrow{(\tilde q^{+,*})^{-1}}\HO_{ G\fps}(q(p^{-1}(\Rt^+\times \Rt^+)), \DD\BoQ_{q(p^{-1}(\Rt^+\times \Rt^+))})
    \\
    &\xrightarrow{(m^+)_*} \HO_{G\fps}(\Rt^+, \DD\BoQ_{\Rt^+}).
\end{split} 
\end{align}

Switching on the $T$-equivariance again, we get:
\begin{theorem} $ $
\label{thm: positive Coulomb in Coulomb}
    \begin{enumerate}
        \item The map $\eta^+$ equips $\HO_{G\fps_T}(\Rt^+, \DD\BoQ_{\Rt^+})$ with an associative unital algebra structure, where the unit is given by the fundamental class of the fiber $\Rt^+\to \Gr^+$ at the base point $[1]\in \Gr^+$, corresponding to the trivial Hecke modification.
        \item The morphism
        \begin{equation}
            \label{eq: positive to full Coulomb map}
            \rho_*: \HO_{G\fps_T}(\Rt^+, \BoD\BoQ_{\Rt^+})\to \HO_{G\fps_T}(\Rt, \BoD\BoQ_\Rt)
        \end{equation}
        is an injective morphism of unital algebras. Equivalently, $\rho_*$ preserves the unit elements and the following diagram is commutative
        \[
        \begin{tikzcd}
            \HO_{G\fps_T}(\Rt^+, \BoD\BoQ_{\Rt^+})\otimes_{\HO_{T}} \HO_{G\fps_T}(\Rt^+, \BoD\BoQ_{\Rt^+})\arrow[r, "\rho_*\otimes \rho_*"]
            \arrow[d, "\eta^+"]&\HO_{G\fps_T}(\Rt, \BoD\BoQ_\Rt)\otimes_{\HO_{T}} \HO_{G\fps_T}(\Rt, \BoD\BoQ_{\Rt})\arrow[d, "\eta"]
            \\
            \HO_{G\fps_T}(\Rt^+, \BoD\BoQ_{\Rt^+})\arrow[r, "\rho_*"] & \HO_{G\fps_T}(\Rt, \BoD\BoQ_{\Rt})
        \end{tikzcd}
        \]
        
    \end{enumerate}
\end{theorem}

Before proving the theorem, consider the special case when $G$ is a torus, that is when all the entries of the dimension vector $\vv\in \BoN^{Q_0}$ are either $0$ or $1$. Up to replacing the quiver with the sub-quiver where $\vv$ has full support), we may assume that $\vv_i=1$ for all $i\in Q_0$. Then we have (on $\mathbb{C}$-points) $\Gr_T=X_*(S)=\BoZ^{Q_0}$ and $\Gr_T^+=X_*(S)=\BoN^{Q_0}$. Moreover, we have 
\[
\Rt=\bigsqcup_{\lambda\in \BoZ^{Q_0}} \{\lambda\}\times \left(t^{\lambda} N\fps \cap N\fps \right)\qquad \Rt^+=\bigsqcup_{\lambda\in \BoN^{Q_0}} \{\lambda\}\times \left(t^{\lambda} N\fps \cap N\fps \right)
\]
Therefore, in the abelian case the map $\rho: \Rt^+\to \Rt$ is just the inclusion of some connected components in $\Rt$. Then all the statements of the theorem are clear. In particular, we have:

\begin{lemma}
\label{lemma: inclusion provite Coulomb in injective in abelian case}
    Assume that $G=S$. Then the morphism \eqref{eq: positive to full Coulomb map} is injective.
\end{lemma}

\begin{proof}[Proof of Theorem \ref{thm: positive Coulomb in Coulomb}]
    The proof of the first point is lengthy, but parallels step by step the proof of associativity for $\eta$ given in \cite{BFNII}, so we omit it. The proof of the fact that $\rho$ is compatible with the units follows from functoriality of fundamental classes in Borel-Moore homology together with the observation that $[1]\in \Gr$ is contained in $\Gr^+$, and so its fiber with respect to $\Rt\to \Gr$ factors through the inclusion $\rho: \Rt^+\hookrightarrow \Rt$. The statement that $\rho_*$ is a morphism of algebras follows from standard base-change arguments.
    
    It remains to prove that $\rho_*$ is injective. Let $S\subset G$ be a maximal torus. It suffices to prove that the morphism in $S\fps$-equivariant cohomology
    \[
    \rho_*: \HO_{S\fps}(\Rt^+, \BoD\BoQ_{\Rt^+})\to \HO_{S\fps}(\Rt, \BoD\BoQ_\Rt)
    \]
   as the statement in $G\fps$-equivariant cohomology would then follow from Proposition \eqref{prop: first properties positive Coulomb} and \cite[Lemma 5.3]{BFNII} by taking $W$-invariants. 
   Let $N_S$ be the restriction of the $G$-module $N$ to $T$. In the quiver language, this corresponds to replacing $(Q, \vv)$ with $(Q^{\text{ab}}, \ul{1})$, where $Q^{\text{ab}}$ is the abelianized quiver, obtained by replacing a vertex $i$ of $Q$ with $\vv_i$ vertices. Set $\Rt_{T, N_S}$ be the space of triples for $(T, N_S)$. We have a canonical closed embedding $\iota: \Rt_{S,N_S}\to \Rt$. Since $(\Rt_{S,N_S})^S=\Rt^S$, the localization theorem implies that the pushforward map 
    \begin{equation}
        \label{eq: abelianization Coulomb branch}
        \iota_*: \HO_{S\fps}(\Rt_{S,N_S}, \BoD\BoQ_{\Rt_{S,N_S}}) \to \HO_{S\fps}(\Rt, \BoD\BoQ_{\Rt})
    \end{equation}
    is an isomorphism after tensoring with $\Frac(\HO_S)$.
    
    Consider now the positive halves $\Rt^+$ and $\Rt_{S,N_S}^+$. We have a similar closed embedding $\iota^+:  \Rt_{S,N_S}^+\hookrightarrow \Rt^+$ satisfying $(\Rt^+)^S=(\Rt^+_{S, N_S})^S$, so the associated pushforward morphism
    \[
    \iota^+_*: \HO_{S\fps}(\Rt^+_{S,N_S}, \BoD\BoQ_{\Rt^+_{S,N_S}}) \to \HO_{G\fps}(\Rt^+, \BoD\BoQ_{\Rt^+})
    \]
    becomes an isomorphism after localization.
    Moreover, we have a commutative diagram
    \[
    \begin{tikzcd}
        \HO_{S\fps}(\Rt^+_{S,N_S}, \BoD\BoQ_{\Rt^+_{S,N_S}}) \arrow[r, "\iota^+_*"] \arrow[d, "\rho^T_*"]&  \HO_{S\fps}(\Rt^+, \BoD\BoQ_{\Rt^+})\arrow[d, "\rho_*"]
        \\
        \HO_{S\fps}(\Rt_{S,N_S}, \BoD\BoQ_{\Rt_{S,N_S}}) \arrow[r, "\iota_*"]&  \HO_{S\fps}(\Rt, \BoD\BoQ_{\Rt})
    \end{tikzcd}
    \]
    where $\rho^T_*$ is the pushforward associated to the inclusion $\rho_T: \Rt^+_{S,N_S}\to \Rt_{S,N_S}$. By Lemma \ref{lemma: inclusion provite Coulomb in injective in abelian case}, $\rho^T_*$ is injective. Moreover, both horizontal morphisms are isomorphisms after tensoring with $\Frac(\HO_S)$, so the localization of $\rho_*$
    \[
    \HO_{S\fps}(\Rt^+, \BoD\BoQ_{\Rt^+})\otimes \Frac(\HO_S)\to \HO_{S\fps}(\Rt, \BoD\BoQ_{\Rt})\otimes \Frac(\HO_S)
    \]
    is injective. Since $\HO_S$ is an integral domain, injectivity of $\rho_*$ now follows from flatness of the $\HO_G$-module $\HO_{G\fps}(\Rt^+, \BoD\BoQ_{\Rt^+})$, see Proposition \ref{prop: first properties positive Coulomb}.
\end{proof}

The algebra $(\HO^{\BM}_{G\fps}(\Rt^+, \BoQ), \eta^+)$ has a similar grading, induced by the inclusions $\Gr^+\hookrightarrow \Gr$. Set $\Gr^+_{\dd}= \Gr_{\dd}\cap \Gr^+$. It is easy to check that $\Gr_{\dd}\cap \Gr^+=\emptyset$ unless $\dd\in \BoN^{Q_0}$. In the sheaf language, the component $\Gr^+_{\dd}$ parametrizes those tuples of inclusions $\Sh{V}_i\hookrightarrow \Sh{O}^{\vv_i}$ such that the cokernels $\Sh{O}^{\vv_i}/\Sh{V}_i$ are of length $\vv_i$. Therefore, we obtain a decomposition $\Rt_+=\bigsqcup_{\dd\in \BoN^{Q_0}} \Rt^+_{\dd}$ and hence a $\BoN^{Q_0}$-grading 
\[
\HO^{\BM}_{G\fps}(\Rt^+, \BoQ)=\bigoplus_{\dd\in \BoN^{Q_0}} \HO^{\BM}_{G\fps}(\Rt^+_{\dd}, \BoQ)
\]

\begin{definition}
    The positive part of the (quantized, flavor deformed) Coulomb branch is the $\BoN^{Q_0}$-graded, cohomologically graded, unital algebra $\Coul^{T, +}_{G,N}=(\HO^{\BM}_{G\fps_T}(\Rt^+, \BoQ), \eta^+)$. 
\end{definition}

\begin{remark}
\label{rem: virtual positive Coulomb branch}
Like in the case of the whole Coulomb branch, this grading is compatible with the multiplication $\eta^+$. Additionally, following Remark \ref{rem: virtual Coulomb branch}, we set 
\[
\Coul^{T, \vir,+}_{G,N}=\bigoplus_{\dd\in \BoN^{Q_0}}\HO^{\BM}_{G\fps_T}(\Rt^+_{\dd}, \BoQ[2s_{\dd}]).
\]
This is also a $\BoZ^{Q_0}$-graded, cohomologically graded algebra.
\end{remark}

\begin{remark}
\label{rem: retraction positive Coulomb 2}
    The retraction  $ \Rt^+_{G,N} \to  \Rt^+_{G, N_0}$ outlined in Remark \ref{rem: retraction positive Coulomb 1} induces a canonical isomorphism 
    \[
    \HO^{\BM}_{G\fps_T}(\Rt^+_{G,N})\cong \HO^{\BM}_{G\fps_T}(\Rt^+_{G,N_0}, \BoQ)
    \]
    of $\HO_{G\fps\times\BoC^*}$-modules.
    Moreover, this isomorphism is clearly compatible with the convolution product on either sides of the equation above. Therefore, it induces an isomorphism of the convolution algebras $\Coul^+_{G,N}$ and $\Coul^+_{G,N_0}$. As noted in Remark \ref{rem: retraction positive Coulomb 1}, the  vector bundle $\Rt^{T,+}_{G,N}\to \Rt^{T,+}_{G, N_0}$ does not extend to $\Rt^{T,+}_{G,N}$ and hence there is in general no isomorphism between $\Coul_{G,N}$ and $\Coul_{G,N_0}$.
\end{remark}

We conclude this section with an analysis of the generators of $\Coul^{T,+}_{G,N}$. Recall the general shape of minuscule weights of $G$ from \S\ref{subsec: Coulomb branch basics}. We say that a minuscule weight $\lambda$ is positive if it is of the form $\lambda=(\omega_{i,k_i})_{i\in Q_0}$ for some $1\leq k_i\leq \vv_i$. Notice that for any positive $\lambda$, the subspace $\Rt_{\lambda}\subset \Rt$ is contained in $\Rt^+$. Therefore,  for any $f\in \HO_{\Stab_{G}(\lambda)}$ we have a well defined class 
\[
M_{\lambda, f}=f\cap [\Rt_{\lambda}]\in \HO_{G\fps_T}(\Rt^+, \DD\BoQ_{\Rt^+}).
\]
We call such operators positive monopole operators.
The following result is a straightforward adaptation of \cite[Remark 6.7]{braverman2018coulombbranches3dmathcal}, \cite[Thm. 4.32]{FT_ShitedQuantum}, and \cite[Thm. 3.7]{Weekes}.

\begin{proposition}
\label{prop: positive coulomb branch generators}
    The positive Coulomb branch $\HO_{G\fps_T}(\Rt^+, \DD\BoQ_{\Rt^+})$ is generated as a $\HO_{G\fps_T}$-algebra by the positive monopoles operators $M_{f, \lambda}$. Moreover, the localized positive Coulomb branch $\HO_{G\fps_T}(\Rt^+, \DD\BoQ_{\Rt^+})\otimes_{\HO_{\BoC^*_{\hbar}}}\Frac(\HO_{\BoC^*_{\hbar}})$ is generated as a $\HO_{G\fps_T}\otimes_{\HO_{\BoC^*_{\hbar}}}\Frac(\HO_{\BoC^*_{\hbar}})$-algebra by the monopoles operators $M_{f, \omega_{i,1}}$.
\end{proposition}

\section{The BFN Coulomb branch via vanishing cycles}

\subsection{Critical description of the Coulomb Branch}
\label{subsec: critical Coulomb branch}

In this section we provide a description of the Coulomb branch in terms of critical cohomology. As it will dimensionally reduce to the BFN Coulomb branch, it is not intended as an alternative definition, but rather as an alternative presentation that is better suited to construct the Coulomb branch action on quasimap critical cohomology. We begin with the heuristics of the construction.

We retain the notation of \S\ref{subsec: The case of quiver representations}. In particular, we fix a representation $N$ of a reductive group $G$, and let $\mu: T^*N\to \Fg$ be the associated moment map. We begin with the heuristics of the construction. Consider the mapping stack $\Map(\DD, T^*N/G)$. Its closed points are tuples $(\Sh{P}, n, n^*)$, where $\Sh{P}$ is a principal $G$ bundle and $n$ (resp.) $n^*$ are sections of the associated bundles $\Sh{P}\times_G N$ (resp.  $\Sh{P}\times_G N^*$).
The moment map lifts to a morphism 
\begin{equation}
    \label{eq: moment map doubled Coulomb}
\mu: \Map(\DD, T^*N/G)\to \Map(\DD, \Fg/G)
\end{equation}
fitting in the Cartesian square 
\[
\begin{tikzcd}
    \Map(\DD, \mu^{-1}(0)/G)\arrow[r, hookrightarrow] \arrow[d] & \Map(\DD, T^*N/G)\arrow[d, "\mu"]\\
     \Map(\DD, 0/ G)\arrow[r, hookrightarrow] &  \Map(\DD, \Fg/ G)
\end{tikzcd}
\]
Hence, a closed point in $\Map(\DD, \mu^{-1}(0)/G)$ is a closed point $(\Sh{V}, n, n^*)$ in $\Map(\DD, T^*N/G)$ satisfying the moment map condition $\mu(n, n^*)=0$ in $\prod_{i\in Q_0}\End(\Sh{V}_i)\cong \Fg\fps$.

Let now $\Hecke^{\pm}(\DD,\mu^{-1}(0)/G)$ be the stack of all Hecke modifications of $\Map(\DD, \mu^{-1}(0)/G)$. Formally, we define 
\[
\Hecke^{\pm}(\DD,\mu^{-1}(0)/G)=\Hecke^{\pm}(\DD, T^*N/G)\times_{\Map(\DD, T^*N/G)} \Map(\DD, \mu^{-1}(0)/G).
\]
Notice that although there are two canonical maps $\Hecke^{\pm}(\DD,\mu^{-1}(0)/G)\to \Map(\DD, T^*N/G)$ the two associated fiber products are canonically isomorphic. More explicitly, $\Hecke^{\pm}(\DD, N/G)$ is the moduli of parametrizing tuples $(\Sh{P}_1, \Sh{P}_2, g, n_1, n_1^*, n_2, n_2^*)$ where $g: \Sh{P}|_{\DD^\times}\to \Sh{P}_2|_{\DD^\times}$ is an isomorphism and the sections satisfy $g n_1=n_2$, $gn^*_1 =n_2^*$, as well as the moment map equation $\mu(n_1, n_1^*)=0$. Notice that the latter is then equivalent to $\mu(n_2, n_2^*)=0$ because of the equations $g n_1=n_2$, $gn^*=n_2^*$. Additionally, $n_1$ and $n_1^*$ are uniquely determined by $n_2, n_2^*$ and $g$.

As in \S\ref{subsec: Mapping stacks on the formal disk and related geometries}, $\Hecke^{\pm}(\DD,\mu^{-1}(0)/G)$ admits the global quotient description 
\[
\Hecke^{\pm}(\DD,\mu^{-1}(0)/G)= \Rt_{\mu^{-1}(0), G}/G\fps,
\]
where $ \Rt_{\mu^{-1}(0), G}\subset \Rt_{T^* N, G}$ is the zero locus of the moment map equation $\mu: \Rt_{T^* N, G}\to \Fg\fps$, which descends to $\eqref{eq: moment map doubled Coulomb}$ after quotienting by $G\fps$. We now exhibit $\Rt_{\mu^{-1}(0), G}$ as a critical locus. Firstly, let $\tilde\Rt_{G,N}$ be the moduli space parametrizing  $(\Sh{P}, \Sh{P}_2, g, n_1, n_2, n_2^*)$ where $g: \Sh{P}|_{\DD^\times}\to \Sh{P}_2|_{\DD^\times}$ is an isomorphism but now no compatibility with $n_i$ and $n_i^*$ is assumed, nor the moment map equation $\mu(n_1, n^*_1)=0$.
Explicitly, we have 
\begin{equation}
    \label{critical raviolo ambient space}
    \tilde \Rt_{G,N}\coloneqq N\fps\times_{G\fps} G\lfps \times_{G\fps} (N^*\lfps\times N\fps)
\end{equation}
Notice that $\Rt_{\mu^{-1}(0), G}$ is a closed substack of $\tilde\Rt_{G,N}$. Consider now the potential $\Sh{W}: \tilde\Rt_{G,N}\to \BoC$ given by 
\begin{equation}
    \label{eq: Coulomb branch potential}
    \Sh{W}=\Res\Tr\left( n_2^*\left(gn_1-n_2\right)\right).
\end{equation}
The following result follows from a straightforward variation of the argument of Theorem \ref{thm: Čech critical quasimap}.
\begin{proposition}
\label{prop: critical raviolo}
There is an equivalence of ind-stacks $\Crit(\Sh{W})\cong \Rt_{\mu^{-1}(0), G}$.
\end{proposition}

We conclude this section by discussing torus actions on $\tilde \Rt_{G,N}$. For a general pair $(G,N)$, consider the torus $T=A\times \BoC^*_{\hbar}$, introduced in \S\ref{subsec: Quantization and Flavor Symmetry}. It naturally acts on $T^*N$ by letting $\hbar$ scale the cotangent direction with weight $\hbar^{-1}$. We consider the induced action on $\tilde \Rt_{G,N}$. We remark that now the action of the subtorus $\BoC^*_{\hbar}$ is induced by both its action on the disk $\DD$ and on $N^*\subset T^*N$. Accordingly, incorporating the natural action of $G\fps$ the relevant torus action on $\tilde \Rt_{G,N}$ becomes the one $G\fps\rtimes T$. We will always arrange this action in such a way that the uniformizer $t$ of $\DD=\Spec\BoC[\![t]\!]$ and the cotangent fiber are scaled with opposite weights. This implies that the potential \eqref{eq: Coulomb branch potential} is fixed by the $T$ action\footnote{And $\Rt_{\mu^{-1}(0),G}$ is preserved by it.}.

If $(G,N)$ are of quiver type, then we set $T$ to be the torus introduced in \S\ref{subsec: torus action quiver varieties}.

\subsection{Finite type approximations}
\label{subsec: Finite type approximations critical Coulomb}
We now exploit this map to give a presentation of the Coulomb branch algebra in terms of critical cohomology. This will be achieved in three steps. In the first one, we define formally define the critical cohomology of the pair $(\tilde \Rt_{G,N},\Sh{W})$ and show it dimensionally reduces to $ \HO^{\BM}_{G\fps_T}(\Rt^+, \BoQ)$. Secondly, we equip this critical cohomology with a convolution algebra structure. Finally, we show that the resulting convolution algebra is identified, under dimensional reduction, to the BFN convolution algebra on $\HO^{\BM}_{G\fps_T}(\Rt^+, \BoQ)$.

We now proceed with the first step. As usual, we temporarily drop the torus $T$ from the notation. Consider the infinite rank vector bundle $\tilde\Rt_{G,N}\to \Gr_G$ and let $\tilde\Rt_{N,G, \leq\lambda}$ denote its restriction to a Schubert variety $\Gr_{G, \leq \lambda}$. Explicitly, we have 
\begin{align*}
    \tilde\Rt_{N,G, \leq\lambda}
    &=N\fps\times_{G\fps} G\lfps_{\leq \lambda}\times N^*\lfps\times N\fps
    \\
    &=\Tot_{\Gr_{G, \leq \lambda}}(N\fps)\times N^*\lfps\times N\fps
\end{align*}
where $G\lfps_{\leq \lambda}$ is the preimage of $\Gr_{G, \leq \lambda}$ under $G\lfps\to \Gr_G$. Fix non-negative integers $m, r, d\geq 0$. Let $N^m\fps=N\fps/t^mN\fps$, let $N^{*}_d\lfps\subset N\lfps$ be the $G\fps$ submodule with poles of order at most $d$, and let $N^*\lfps^{r}_d\subset N\lfps$ be its quotient by $ t^r N\fps \subset N^{*}_d\lfps$. Set 
\begin{align*}
    \tilde\Rt^{ m,d,r}_{N,G,\leq\lambda}
    \coloneqq\Tot_{\Gr_{G, \leq \lambda}}(N^m\fps)\times N^*\lfps^{r}_d\times N^d\fps
\end{align*}
It is a finite rank vector bundle over $\Gr_{G, \leq \lambda}$, and hence it is a variety. Notice that the representative $g(t)\in G\lfps_{\leq \lambda}$ of any point in $\Gr_{\leq \lambda}$ has poles bounded by a constant that depends on $\lambda$. As a consequence, for $m\ggg d$ (relative to $\lambda$), the $G\fps$ module $t^m N\fps$ satisfies $g(t) t^m N\fps\subseteq t^dN\fps$. Therefore, the restriction of the potential \eqref{eq: Coulomb branch potential} to 
\[
\Tot_{\Gr_{G, \leq \lambda}}(N\fps)\times N^*_{d}\lfps\times N\fps
\]
factors through the quotient bundle $\tilde\Rt^{ m,d,r}_{N,G,\leq\lambda}$.
\[
\Sh{W}^{m,d,r}_{\lambda}: \tilde\Rt^{ m,d,r}_{N,G,\leq\lambda}\to \BoC
\]
the resulting map. It is $G\fps$ invariant for the residual $G\fps$ action on $\tilde\Rt^{ m,d,r}_{N,G,\leq\lambda}$. Additionally, for $k\ggg 0$\footnote{in particular, $k>d+r$ to make sure that that the $G\fps$ action on $N_d^{*,r}\lfps$ factors through $G^k\fps$.} the $G\fps$ action on $\tilde\Rt^{m,d,r}_{N,G,\leq\lambda}$ factors through $G^k\fps= G(\Sh{O}/t^k \Sh{O})$, so we have a well-defined critical cohomology
\[
\HO_{G^k\fps}(\tilde\Rt^{m,d,r}_{N,G, \leq\lambda }, \phip{\Sh{W}^{m,d}}\DD\BoQ).
\]
\begin{remark}
\label{rem: order of zero dual field Coulomb branch}
    From now on, we further assume that $r\ggg 0$ (relatively to $\lambda$). This will ensure that, for any representative $g\in G\lfps$ of an element in $\Gr_{\leq \lambda}$, the order of the pole of the term $g\alpha$ appearing in \eqref{eq: Coulomb branch potential} is smaller or equal than $r$. The relevance of this condition will become clear in Proposition \ref{eq: direct sum dec crit Coulomb}.
\end{remark}

\begin{lemma}
    For any $m',d',r'\geq 0$ and $l$ large enough, we have canonical isomorphisms
    \begin{multline*}
        \HO_{G_{k+l}}(\tilde \Rt^{m+m',d+d',r+r'}_{N,G, \leq \lambda}, \phip{\Sh{W}^{m+m', d+d',r+r'}_{\lambda}}\DD\BoQ)\cong
        \\\HO_{G_{k}}(\tilde \Rt^{m,d,r+r'}_{N,G, \leq \lambda}, \phip{\Sh{W}^{m, d,r}_{\lambda}}\DD\BoQ)[-2((d'+m'+r')\dim(N)-l\dim(G))].
    \end{multline*}
\end{lemma}
\begin{proof}
    Consider the $G\fps$ equivariant maps 
    \[
    \begin{tikzcd}
        \Rt^{m,d,r}_{N,G, \leq \lambda} & \Tot_{\Gr_{G, \leq \lambda}}(N\fps^{m+m'})\times N\fps^{d+d'}\times N^*\lfps^{r+r'}_{d}\arrow[r, hookrightarrow, "i"] \arrow[l, swap, "j"]& \Rt^{m+m',d+d'}_{N,G, \leq \lambda}.
    \end{tikzcd}
    \]
    The map $j$ is induced by taking quotients of the fibers, and hence is an affine fibration or rank $(d'+m'+r')\dim(N)$. The map $i$ is induced by the inclusion $N_d\lfps\subseteq N_{d+d'}\lfps$, and hence is a closed immersion. To remove clutter, set 
    \[
    Y=\Tot_{\Gr_{G, \leq \lambda}}(N\fps^{m+m'})\times N\fps^{d+d'}\times N^*\lfps^{r+r'}_{d}.
    \]
    Passing to cohomology, we obtain canonical morphisms 
    \begin{align*}
        \HO_{G_{k+l}}(\tilde \Rt^{m,d,r}_{N,G, \leq \lambda}, \phip{\Sh{W}^{m, d,r}_{\lambda}}\DD\BoQ)
        &\xrightarrow{j^*} \HO_{G_{k+l}}( Y, \phip{\Sh{W}^{m, d,r}_{\lambda}\circ j}\DD\BoQ)[-2(d'+m'+r')\dim(N)]
        \\
        &= \HO_{G_{k+l}}(Y, \phip{\Sh{W}^{m+m', d+d'}_{\lambda}\circ i}\DD\BoQ)[-2(d'+m'+r')\dim(N)]\\
        &\xrightarrow{i_*}\HO_{G_{k+l}}(\tilde \Rt^{m+m',d+d',r+r'}_{N,G, \leq \lambda}, \phip{\Sh{W}^{m+m', d+d',r+r'}_{\lambda}}\DD\BoQ)[-2(d'+m'+r')\dim(N)].
        \end{align*}
    Since $j$ is an affine fibration, the pullback $j^*$ is an isomorphism. The fact that $i_*$ is also an isomorphism is a consequence of Theorem \eqref{thm: dimred}. The proof then follows from the fact that 
    \[
    \HO_{G_{k+l}}(\tilde \Rt^{m,d,r}_{N,G, \leq \lambda}, \phip{\Sh{W}^{m, d,r}_{\lambda}}\DD\BoQ)=\HO_{G_{k}}(\tilde \Rt^{m,d,r}_{N,G, \leq \lambda}, \phip{\Sh{W}^{m, d,r}_{\lambda}}\DD\BoQ)[2l\dim(G)].
    \]
\end{proof}
The previous lemma implies that we can set, without ambiguity, 
\begin{equation}
    \label{eq: def crit Coulomb branch over orbit closure}
    \HO( \tilde\Rt_{N,G, \leq\lambda},\phip{\Sh{W}_{\lambda}}\DD\BoQ)\coloneqq\HO_{G_{k}}(\tilde \Rt^{m,d,r}_{N,G, \leq \lambda}, \phip{\Sh{W}^{m, d}_{\lambda}}\DD\BoQ)[-2(d+m+r)\dim(N)-k\dim(G))].
\end{equation}
for $m,d,r,k\ggg 0$. Moreover, the closed embeddings $\tilde \Rt_{N,G\leq \lambda}\hookrightarrow \Rt_{N,G\leq \mu}$ gives canonical pushforward morphisms
\[
\HO( \tilde\Rt_{N,G, \leq\lambda},\phip{\Sh{W}_{\lambda}}\DD\BoQ)\to \HO( \tilde\Rt_{N,G, \leq\mu},\phip{\Sh{W}_{\mu}}\DD\BoQ)
\]
which are compatible with the finite approximations on the domain and target\footnote{Where the approximations are given by the same choice of $m\ggg d>0$ and $k\ggg0$ on the domain and target.}. This finally leads us to the definition of the critical Coulomb branch:
\begin{equation}
    \label{eq: def crit Coulomb}
    \HO_{G\fps}(\tilde \Rt_{G,N}, \phip{\Sh{W}}\DD\BoQ)\coloneqq\lim_{\to} \HO_{G\fps}(\tilde \Rt_{N,G, \leq \lambda}, \phip{\Sh{W}_{\lambda}}\DD\BoQ).
\end{equation}

By construction, a class in \eqref{eq: def crit Coulomb} is of degree $n$ if it is represented by a class of degree $n-2(d+m+r)\dim(N)-k\dim(G))$ in $\HO_{G_{k}}(\tilde \Rt^{m,d,r}_{N,G, \leq \lambda}, \phip{\Sh{W}^{m,d,r}_{\lambda}}\DD\BoQ)$ for some dominant $\lambda$ and large enough $m,d,r$, and $k$.

\begin{remark}
\label{rem: critical Coulomb branch degree}
    In the rest of the paper, we will denote the shift appearing in \eqref{eq: def crit Coulomb branch over orbit closure} by
    \[
    -2\infty_{M,G}\coloneqq-2\dim(N^*\lfps)-2\dim N\fps+2\dim G\fps,
    \]
    which we understand as follows: the term $-2\dim(N^*\lfps)$ provides, in a finite approximation, a factor of $-2(d+r)\dim(N)$; the term $-2\dim N\fps$ provides a factor of $-2m\dim(N)$, and the term $2\dim G\fps$ provides a factor $-2k\dim(G)$. The reason why a larger factor appears in front of the shift associated to  $-2\dim(N\lfps)$ is due to the fact that $N^*\lfps$ ``counts twice'' since it's infinite in both the positive and negative direction.
\end{remark}

\begin{remark}
    We stress here that, unlike for the case of quasimaps discussed in \S\ref{subsec: DT sheaf Čech model}, we do \emph{not} interpret the critical locus of the function $\Sh{W}_{\lambda}^{m,d,r}$ on $\tilde \Rt^{m,d,r}_{N,G, \leq \lambda}$ as a derived critical locus, yet alone as a $(-1)$ shifted symplectic scheme. In fact, this statement not even make sense as the orbit $\Gr_{G, \leq \lambda}$ and hence $ \Rt^{m,d,r}_{N,G, \leq \lambda}$ is not smooth in general. On the other hand, in \S\ref{subsec: Coulomb branch action on quasimaps} we will show that the vanishing cycle formalism on $\Rt^{m,d,r}_{N,G, \leq \lambda}$ combined with Proposition \ref{prop: dr Coulomb vector sp} below provides an alternative—but equivalent—construction of the BFN Coulomb branch.
\end{remark}
Notice also that the composition $\tilde \Rt_{G,N}\to \Gr_{G}\to \pi_1(G)$
induces a decomposition 
\[
\tilde \Rt_{G,N}=\bigsqcup_{\gamma \in \pi_1(G)}  \tilde \Rt_{N,G, \gamma}
\]
and hence a decomposition
\begin{equation}
    \label{eq: direct sum dec crit Coulomb}
    \HO_{G\fps}(\tilde \Rt_{G,N}, \phip{\Sh{W}}\DD\BoQ)=\bigoplus_{\gamma\in \pi_1(G)} \HO_{G\fps}(\tilde \Rt_{N,G, \gamma}, \phip{\Sh{W}}\DD\BoQ).
\end{equation}
Reinstating the torus equivariance, we get:
\begin{proposition}
\label{prop: dr Coulomb vector sp}
    We have a canonical isomorphism
    \[
    dr: \HO_{G\fps_T}(\tilde \Rt_{G,N}, \phip{\Sh{W}}\DD\BoQ)\xrightarrow{\cong} \HO_{G\fps_T}(\Rt_{G,N},\DD\BoQ),
    \]
   of $\BoZ\times \pi_1(G)$-graded vector spaces with grading determined by the conventions in \S\ref{subsec: cohomology of pro and ind stacks} and Remark \ref{rem: critical Coulomb branch degree}, respectively.
\end{proposition}
\begin{proof}
    By definition of both sides, the proposition is really about the existence of a dimensional reduction morphism
    \[
    \HO_{G^k\fps}(\tilde \Rt^{m,d,r}_{N,G, \leq \lambda}, \phip{\Sh{W}^{m,d,r}_{\lambda}}\DD\BoQ)\cong \HO_{G^k\fps}(\Rt^m_{N,G, \leq \lambda},\DD\BoQ[2d\dim(N^{*})]). 
    \]
    for the finite approximations.
    Let $\bar \Rt^{d}_{N, G, \leq \lambda}$ be the subspace of $\tilde \Rt^{m,d,r}_{N,G, \leq \lambda}$ with $n_2^*=0$. Then $\tilde \Rt^d_{N, G, \leq \lambda}$ is naturally vector bundle over $\bar \Rt^{m,d,r}_{N, G, \leq \lambda}$ of rank $(d+r)\dim(N^{*})$. Via residue pairing, the assignment
    \[
    \bar \Rt^d_{N, G, \leq \lambda} \ni ([n_1, g],n_2)\mapsto gn_1-n_2\in N\lfps/t^dN\fps
    \]
    gives a section of the dual bundle $(\tilde \Rt^{m,d,r}_{N, G, \leq \lambda})^\vee$ over $\bar \Rt^{d}_{N, G, \leq \lambda}$. By Remark \ref{rem: order of zero dual field Coulomb branch} the zero locus of this section is 
    the locus of $\bar \Rt^{d}_{N, G, \leq \lambda}$ consisting of those triples $([n_1, g],n_2)$ satisfying $n_2-gn_1=0\mod t^dN\fps$ in $ N\lfps$. Since $n_2\in N\fps/t^dN\fps$ is holomorphic, this condition enforces $gn_1$ to be holomorphic and uniquely determines $n_2$ in terms of $n_1$ and $g$. Therefore, the zero locus of $n_2-gn_1=0\mod t^dN\fps$ is canonically isomorphic to a rank $2(d+r)\dim(N^{*})$ vector bundle over $\Rt^{m}_{N,G, \leq \lambda}$. Hence, the statement follows from dimensional reduction \ref{thm: dimred}.
\end{proof}

\subsection{Critical convolution algebra structure}

\label{subsec: critical Coulomb branch convolution algebra}
Next, we define a convolution algebra structure on $ \HO_{G\fps_T}(\tilde \Rt_{G,N}, \phip{W}\DD\BoQ)$ by recasting BFN's construction in critical cohomology. 
For $i=1,\dots, 4$ let $\tilde\Rt^{(2),i}_{G,N}$ be the moduli spaces
\begin{align}
\label{eq: three spaces critical Coulomb mult}
\begin{split}
    \tilde\Rt^{(2),0}_{G,N}&=N\fps\times_{G\fps} G\lfps \times N^*\lfps\times N\fps\times  G\lfps \times N^*\lfps \times N\fps
    \\
    \tilde\Rt^{(2),1}_{G,N}&=N\fps\times_{G\fps} G\lfps \times_{G\fps} \left(N^*\lfps\times N\fps\times G\lfps\right) \times N^*\lfps \times N\fps
    \\
    \tilde\Rt^{(2),2}_{G,N}&=\left(N^*\lfps\times N\lfps\times N\fps\right) \times_{G\fps}G\lfps \times_{G\fps} G\lfps \times (N^*\lfps \times N\fps)
    \\
    \tilde\Rt^{(2),3}_{G,N}&=\left(N^*\lfps\times N\fps\right) \times_{G\fps}G\lfps \times_{G\fps} G\lfps \times (N^*\lfps \times N\fps)
    \\
    \tilde\Rt^{(2),4}_{G,N}&= N\fps \times_{G\fps}G\lfps \times_{G\fps} G\lfps \times (N^*\lfps \times N\fps)
\end{split}
\end{align}
Closed points in $\tilde\Rt^{(2),0}_{G,N}, \tilde\Rt^{(2),1}_{G,N}, \tilde\Rt^{(2),2}_{G,N}, \tilde\Rt^{(2),3}_{G,N}$, and $\tilde\Rt^{(2),4}_{G,N}$ are denoted by 
\[
\begin{split}
    & \left( [n_1, g], n_2^*, n_2, g', n_3^*, n_3\right)
    \\
    &[n_1, g, n_2^*, n_2, g', n_3^*, n_3]
    \\
    & [ w^*, w,n_1, g, g', n_3^*, n_3]
    \\
    & [w^*,n_1, g, g', n_3^*, n_3]
    \\
    & [n_1, g, g', n_3^*, n_3]
\end{split}
\]
respectively. Notice that we ordered the entries of the cosets in a way that naturally matches the terms in the fiber product. We assemble this spaces into the following diagram 
\begin{equation}
    \label{eq: critical Coulomb branch correspondence}
    \begin{tikzcd}
        \tilde\Rt_{G,N}\times \tilde\Rt_{G,N} & \tilde\Rt^{(2),0}_{G,N}\arrow[l, swap, "\tilde p"]\arrow [r, "\tilde q_1"] & \tilde\Rt^{(2),1}_{G,N}\arrow[r, "\tilde j"]&  \tilde\Rt^{(2),2}_{G,N}
        \\
        & & & \tilde\Rt^{(2),3}_{G,N}\arrow[u, swap, "\tilde q_2"]\arrow[d, "\tilde q_3"]
        \\
        & & & \tilde\Rt^{(2),4}_{G,N} \arrow[r, "\tilde m"] & \tilde\Rt^{}_{G,N} 
\end{tikzcd}
\end{equation}
The maps $\tilde p$ and $\tilde q_1$ in the diagram are defined by
\begin{multline}
    \label{eq: first two maps critical ravolo correspondence}
     ([n_1, g, n_2^*, n_2], [n_2, g', n_3^*, n_3])\mapsfrom \left( [n_1, g], n_2^*, n_2, g', n_3^*, n_3\right)
     \mapsto [n_1, g, n_2^*, n_2, g', n_3^*, n_3]
\end{multline}
The map $\tilde j$ is defined by
\[
[n_1, g, n_2^*, n_2, g', n_3^*, n_3]\mapsto [g'gn_3^*-gn_2^*, gn_2-n_1,n_1, g, g', n_3^*, n_3].
\]
The maps $\tilde q_2$ and $\tilde q_3$ are defined by 
\[
[w^*, 0, n_1, g, g', n_3^*, n_3] \mapsfrom [w^*,n_1, g, g', n_3^*, n_3]\mapsto [n_1, g, g', n_3^*, n_3]
\]
Finally, the map $\tilde m$ is defined by 
\[
[n_1, g, g', n_3^*, n_3]\mapsto [n_1,  g'g, n_3^*, n_3]
\]
Notice that $\tilde p$ is the composition of a fiber bundle map with a closed embedding (induced by the diagonal map $n_2\mapsto (n_2,n_2)$); the map $\tilde q_1$ is a quotient by $G\fps$; the map $\tilde j$ is a closed immersion; the map $\tilde q_2$ is the zero section of a vector bundle; the map $q_3$ is a vector bundle; the map $\tilde m$ is ind-proper.

\begin{remark}
To help the reader navigate through this section, we clarify in advance that the correspondence \eqref{eq: critical Coulomb branch correspondence} will play the role of diagram \eqref{eq: Coulomb branch mult diagram}—or rather of its last row. Once understood appropriately, all morphisms $\tilde q_i$ will all induce isomorphisms in critical cohomology so the Coulomb branch multiplication will be effectively induced by a pullback along $\tilde p$, a pushforward along $\tilde j$ and a further pushforward along $\tilde m$. We encourage those readers that are still daunted by diagram \eqref{eq: critical Coulomb branch correspondence} to carefully read Remark \ref{rem: ideal coulomb branch mult}. 

Finally, we stress that, since the the expressions $g'gn_3^*-gn_2^*$ and $gn_2-n_1$ may exhibit poles, they naturally live in $N^*\lfps$ and $N\lfps$; they would not be defined as elements in $N^*\fps$ and $N\fps$. 
\end{remark}

\begin{example}
\label{ex: quiver picture critical raviolo correspondence}
In the quiver case, the various spaces in \eqref{eq: critical Coulomb branch correspondence} can be conveniently visualized by means of certain diagrams. We spell out the details in the case of an $A_2$ quiver. The various spaces $\tilde \Rt^{(2), i}_{G,N}$ all live over certain prequotients—depending on $i$— of $G\fps\backslash G\lfps\times_{G\fps}G\lfps$, the moduli space of pairs of Hecke modifications with a bundle trivialized\footnote{in the case of $\tilde \Rt^{(2), 0}_{G,N}$ two bundles out of three are trivialized}. The diagrams will not be sensible of such trivializations. Therefore, we visualize the data coming from the base $G\fps\backslash G\lfps\times_{G\fps}G\lfps$ as
\[
\begin{tikzcd}
    \Sh{V} \arrow[r, dashed, "g"] & \Sh{V}' \arrow[r, dashed, "g'"] & \Sh{V''}
\end{tikzcd}
\]
Here and below, a dotted map refers that it is only defined on the punctures disk $\DD^{\times}$. The remaining data determining the moduli spaces in $\tilde \Rt^{(2), i}_{G,N}$ are then given by tuples of sections between these bundles. Explicitly, we have 
\[
n_1\in \bigoplus_{a\in Q_1} \Hom(\Sh{V}_{s(a)}, \Sh{V}_{t(a)}) \qquad  n_2\in \bigoplus_{a\in Q_1} \Hom(\Sh{V}'_{s(a)}, \Sh{V}'_{t(a)})\qquad n_3 \in \bigoplus_{a\in Q_1} \Hom(\Sh{V}''_{s(a)}, \Sh{V}''_{t(a)})
\]
\[
n_2^*\in \bigoplus_{a\in Q_1} \Hom(\Sh{V}'_{t(a)}, \Sh{V}'_{s(a)}) \qquad n_3^*\in \bigoplus_{a\in Q_1} \Hom(\Sh{V}_{t(a)}, \Sh{V}'_{s(a)})
\]
With these premises, we can visualize the map $\tilde p$ as
\begin{equation}
\label{eq: diagrammatic presentation Coulomb branch correspondence}
\begin{tikzcd}
    \Sh{V}\arrow[d, "n_1"]\arrow[r, dashed, "g"] & \Sh{V}'\arrow[d, "n_2"]
    \\
    \Sh{V}\arrow[r, dashed, "g"] & \Sh{V}'\arrow[u, shift left, dashed, "n_2^*"]
\end{tikzcd}
\times
\begin{tikzcd}
     \Sh{V}'\arrow[d, "n_2"]\arrow[r, dashed, "g'"] & \Sh{V}''\arrow[d, "n_3"]
    \\
    \Sh{V}'\arrow[r, dashed, "g'"]& \Sh{V}''\arrow[u, shift left, dashed, "n_3^*"]
\end{tikzcd}
\mapsfrom
\begin{tikzcd}
    \Sh{V}\arrow[d, "n_1"]\arrow[r, dashed, "g"] & \Sh{V}'\arrow[d, "n_2"]\arrow[r, dashed, "g'"] & \Sh{V}''\arrow[d, "n_3"]
    \\
    \Sh{V}\arrow[r, dashed, "g"] & \Sh{V}'\arrow[u, shift left, dashed, "n_2^*"]\arrow[r, dashed, "g'"]& \Sh{V}''\arrow[u, shift left, dashed, "n_3^*"]
\end{tikzcd}
\end{equation}
The map $\tilde j$ is given by
\[
\begin{tikzcd}
    \Sh{V}\arrow[d, "n_1"]\arrow[r, dashed, "g"] & \Sh{V}'\arrow[d, "n_2"]\arrow[r, dashed, "g'"] & \Sh{V}''\arrow[d, "n_3"]
    \\
    \Sh{V}\arrow[r, dashed, "g"] & \Sh{V}'\arrow[u, shift left, dashed, "n_2^*"]\arrow[r, dashed, "g'"]& \Sh{V}''\arrow[u, shift left, dashed, "n_3^*"]
\end{tikzcd}
\mapsto 
\begin{tikzcd}
    \Sh{V}\arrow[d, dashed, "gn_2-n_1"]
    \\
    \Sh{V}\arrow[u, shift left, dashed,  "g'gn_3^*-gn_2^*"]
\end{tikzcd}
\times 
\begin{tikzcd}
    \Sh{V}\arrow[d, "n_1"]\arrow[r, dashed, "g"] & \Sh{V}'\arrow[r, dashed, "g'"] & \Sh{V}''\arrow[d, "n_3"]
    \\
    \Sh{V}\arrow[r, dashed, "g"] & \Sh{V}'\arrow[r, dashed, "g'"]& \Sh{V}''\arrow[u, shift left, dashed, "n_3^*"]
\end{tikzcd}
\]
Here, for visualization purposes, we are introducing another copy of $\Sh{V}$ although the two copies are to be thought as identified. Equivalently, the fiber product should be thought over $\B G\fps$.

The maps $\tilde q_1, \tilde q_2$ and $\tilde q_3$ are elementary, so we do not give a pictorial representation. Finally, the map $\tilde m$ is given by 
\[
\begin{tikzcd}
    \Sh{V}\arrow[d, "n_1"]\arrow[r, dashed, "g"] & \Sh{V}'\arrow[r, dashed, "g'"] & \Sh{V}''\arrow[d, "n_3"]
    \\
    \Sh{V}\arrow[r, dashed, "g"] & \Sh{V}'\arrow[r, dashed, "g'"]& \Sh{V}''\arrow[u, shift left, dashed, "n_3^*"]
\end{tikzcd}
\mapsto 
\begin{tikzcd}
    \Sh{V}\arrow[d, "n_1"]\arrow[r, dashed, "g'g"] & \Sh{V}''\arrow[d, "n_3"]
    \\
    \Sh{V}\arrow[r, dashed, "g'g"] & \Sh{V}''\arrow[u, shift left, dashed,  "n_3^*"]
\end{tikzcd}
\]
We remark that none of the squares in these diagrams are required to commute. The commutativity of these  squares will be eventually enforced by the vanishing cycle functors.
\end{example}

We now introduce potentials on all the spaces in in diagram \eqref{eq: critical Coulomb branch correspondence}. Let 
\[
\Tr:N\lfps\times N^*\lfps\to \BoC\lfps, (n, n^*)\mapsto \Tr\Res (n^*n).
\]
be the natural pairing. The assignment is as follows:
\begin{itemize}
    \item The potential on $\tilde\Rt_{G,N}\times \tilde\Rt_{G,N}$ is the direct sum potential $\Sh{W}\boxplus \Sh{W}$, where each $\Sh{W}$ is a copy of \eqref{eq: Coulomb branch potential}.
    \item The potential on $\tilde\Rt^{(2),2}_{G,N}$ is 
    \begin{equation}
    \label{eq: potential W_0}
        \Sh{W}_2\coloneqq(\Sh{W}\boxplus \Sh{W})\circ \tilde p=\Res\Tr\left(n_2^*(gn_1-n_2)\right) +\Res\Tr\left(n_3^*(g'n_2-n_3)\right)
    \end{equation}
    the pullback of $\Sh{W}\boxplus \Sh{W}$ along $\tilde p$. Since it is $G\fps$ equivariant, it descends to a potential $\Sh{W}_1$ on $\tilde\Rt^{(2),1}_{G,N}$.
    \item The potential $\Sh{W}_3$ on $\tilde\Rt^{(2),3}_{G,N}$ is given by 
    \begin{equation}
     \label{eq: potential W_1}
    \Sh{W}_3 \coloneqq \Res\Tr(w^*w)+\Res\Tr\left( n_3^*\left(g'g n_1-n_3\right)\right).
    \end{equation}
    \item The potential on $\tilde\Rt^{(2),4}_{G,N}$, denoted by $\Sh{W}_4$, is the composition 
    \[
    \Sh{W}_4\coloneqq\Sh{W}_2\circ \tilde q_2=\Res\Tr\left( n_3^*\left(g'gn_1-n_3\right)\right)
    \]
    \item The potential on $\tilde\Rt^{(2),4}_{G,N}$ is 
    \[
    \Sh{W}_4\coloneqq\Sh{W}\circ \tilde m=\Res\Tr\left( n_3^*\left(g'gn_1-n_3\right)\right).
    \] 
\end{itemize}
We use these potentials to associate critical cohomology to each of these spaces. As for the definition of $\HO_{G\fps}(\tilde \Rt_{G,N}, \phip{\Sh{W}}\DD\BoQ)$, the accurate definition of these critical cohomologies involves considering finite dimensional approximations of the spaces and of the associated potentials. The construction of such finite approximations parallels the one discussed in detail in \S\ref{subsec: Coulomb branch action on quasimaps}, where an action of the critical Coulomb branch $\HO_{G\fps}(\tilde \Rt_{G,N}, \phip{\Sh{W}}\DD\BoQ)$ on quasimap cohomology\footnote{The argument here would actually be simpler, because of the lack of global structure in the Coulomb branch, compared to the moduli of quasimaps $\QM^{\xi}(\bbP^1, X)$}. Since the finite dimensional approximations are discussed there, we do not replicate them here.

From the definition of the potentials on the various spaces in \eqref{eq: critical Coulomb branch correspondence}, it is clear that all the potentials are compatible with the morphisms, except for those on the domain and codomain of $\tilde q_1$. This compatibility is settled by the next lemma. To remove clutter, from now on, we remove the subscripts $G,N$ in all the spaces of this section.
\begin{lemma}
\label{lemma: critical raviolo 1}
    The equality $\Sh{W}_2\circ \tilde j=\Sh{W}_1$ holds.
\end{lemma}
\begin{proof}
     Combining \eqref{eq: first two maps critical ravolo correspondence} and \eqref{eq: potential W_1}, the evaluation of the left and side of the statement at a point $[n_1, g, n_2^*, n_2, g', n_3^*, n_3]$ gives
     \begin{multline*}
         \Res\Tr((g'gn_3^*-gn_2^*)(gn_2-n_1)) + \Res\Tr(n_3^*(g'gn_1-n_3))
         \\
         =\Res\Tr\left((n_3^*(g'n_2-n_3)\right)+\Res\Tr\left(n_2^*(gn_1-n_2)\right),
     \end{multline*}
The equality follows from the invariance of the pairing $\Res\Tr(-)$ under conjugation by $G\lfps$. But the right hand side is exactly \eqref{eq: potential W_0}.
\end{proof}
\begin{lemma}
\label{lemma: critical raviolo 2}
    We have a dimensional reduction isomorphism
    \[
    \tilde q_{2,*}\circ \tilde q_3^*: \HO_{G\fps}(\tilde \Rt^{(2),4}_{G,N}, \phip{\Sh{W}_3}\DD\BoQ)\xrightarrow{\cong} \HO_{G\fps}(\tilde \Rt^{(2),2}_{G,N}, \phip{\Sh{W}_1}\DD\BoQ)
    \]
\end{lemma}
\begin{proof}
    We will show that both maps are isomorphisms. Consider the closed immersion $\tilde q_2: \tilde \Rt^{(2),3}_{G,N}\hookrightarrow \tilde \Rt^{(2),2}_{G,N}$. Applying the vanishing cycle $\phip{\Sh{W}_2}$ to the co-unit $\tilde q_{2,*}\tilde q_2^!\to \id$ gives a natural morphism 
    \[
    \tilde q_{2,*}\phip{\Sh{W}_3}\DD\BoQ_{ \tilde \Rt^{(2),3}_{G,N}}\to \phip{\Sh{W}_2}\DD\BoQ_{ \tilde \Rt^{(2),2}_{G,N}}
    \]
    which, by Theorem \ref{prop: def dim red smooth case}, becomes an isomorphism after passing to derived global sections. This proves that $\tilde q_{2,*}$ is an isomorphism.
    
    On the other hand, the fact the fact that $\tilde q_2^*$ is an isomorphism follows from the observation that $\tilde q_3: \tilde \Rt^{(2),3}\to \Rt^{(2),4}$ is a infinite rank vector bundle and the potential $\Sh{W}_3$ is pulled back from $\Sh{W}_4$ via $\tilde q_3$. Therefore, we get the desired isomorphism isomorphism 
    \[
    \tilde q_2^* : \HO(\Rt^{(2),3}, \phip{\Sh{W}_3}\DD\BoQ)\to \HO(\Rt^{(2),2}, \phip{\Sh{W}_2}\DD\BoQ)  
    \]
    We remark that, despite the unit $\id\to \tilde q_{3,*}\tilde q_3^*$ applied to $\DD\BoQ_{\Rt^{(2),1}}$ induces a \emph{shifted} isomorphism 
    \[
    \DD\BoQ_{\Rt^{(2),1}}\to \tilde q_{2,*}\DD\BoQ_{\Rt^{(2),2}}[-2\dim N^*\lfps]
    \]
    with an infinite shift in our formal description—but an actual finite shift in any finite approximation. The finite approximations of this shift are built in in the definition of the cohomology $\HO(\Rt^{(2),2}, \phip{\Sh{W}_2}\DD\BoQ)$, and hence do not appear in the statement of the lemma.
\end{proof}
We can now complete the definition of the critical convolution product. Firstly, notice that since both $\tilde \Rt\times \tilde \Rt$ and $\tilde \Rt^{(2),0}$ are smooth, we have a canonical isomorphism 
\[
\tilde p^*\DD\BoQ_{\tilde \Rt\times \tilde \Rt}=\DD\BoQ_{\tilde \Rt^{(2),0}}[2\infty_{G,N}],
\]
where as usual we have set $\infty_{G,N}=\dim N\fps-\dim G\fps$.
Therefore, we obtain a canonical morphism 
\begin{equation*}
    \begin{tikzcd}
        \phip{\Sh{W}\boxplus \Sh{W}}\DD\BoQ_{\tilde \Rt\times \tilde \Rt}\arrow[r] & \phip{\Sh{W}\boxplus \Sh{W}} \tilde p_*\tilde p^* \DD\BoQ_{\tilde \Rt\times \tilde \Rt}\arrow[dl, equal] 
        \\
    \phip{\Sh{W}\boxplus \Sh{W}} \tilde p_*\DD\BoQ_{\tilde \Rt^{(2),0}_{G,N}}[2\infty_{G,N}]\arrow[r] & \tilde p_*\phip{\Sh{W}_1}\DD\BoQ_{\tilde \Rt^{(2),0}_{G,N}}[2\infty_{G,N}],
    \end{tikzcd}
\end{equation*}
where the first map is induced by the unit $\id\to \tilde p_*\tilde p^*$ and the last one by base-changing the vanishing cycle. Passing to derived global section gives a pullback morphism\footnote{Similarly to Lemma \ref{lemma: critical raviolo 2}, the shift by $\infty_{G,N}$ is absorbed, after passing to finite approximations, in the definition of the target of the map \eqref{eq: critical convolutionl Coulomb 1}. Similar conventions are implicitly made below.}
\begin{equation}
    \label{eq: critical convolutionl Coulomb 1}
    \tilde p^*: \HO_{G\fps\times G\fps}(\tilde \Rt\times \tilde \Rt, \phip{\Sh{W}\boxplus \Sh{W}}\DD\BoQ)\to \HO_{{G\fps}\times {G\fps}}(\tilde \Rt^{(2),0}, \phip{\Sh{W}_0}\DD\BoQ).
\end{equation}
Moreover, the $G\fps$-invariance of the potential $\Sh{W}\boxplus \Sh{W}\circ p$ gives a pullback isomorphism 
\begin{equation}
    \label{eq: critical convolutionl Coulomb 2}
     \tilde q^*_1: \HO_{G\fps}(\tilde \Rt^{(2),1}, \phip{\Sh{W}_1}\DD\BoQ) \xrightarrow{\cong} \HO_{G\fps\times G\fps}(\tilde \Rt^{(2),0}, \phip{\Sh{W}_0}\DD\BoQ)
\end{equation}
Finally, since $\tilde j$ $\tilde m$ are proper representable we have pushforward morphisms in critical cohomology
\begin{align}
\label{eq: critical convolutionl Coulomb 3}
\begin{split}
    &\tilde j_*:\HO_{{G\fps}}(\tilde \Rt^{(2),1}, \phip{\Sh{W}_1}\DD\BoQ)\to \HO_{{G\fps}}(\tilde \Rt^{(2),2}, \phip{\Sh{W}_2}\DD\BoQ)
    \\
    &\tilde m_*: \HO_{{G\fps}}(\tilde \Rt^{(2),4}, \phip{\Sh{W}_3}\DD\BoQ)\to \HO_{G\fps}(\tilde \Rt, \phip{\Sh{W}}\DD\BoQ)
\end{split}
\end{align}
Recall the isomorphism $ \tilde q_{2,*}\circ \tilde q_3^*$ from Lemma \ref{lemma: critical raviolo 2}. We now define
\begin{equation}
    \label{eq: formula critical Coulomb branch multiplication}
    \etap\coloneqq\tilde m_* \circ (\tilde q_{2,*}\circ \tilde q_3^*)^{-1}\circ \tilde j_*\circ (\tilde q^*_1)^{-1} \circ \tilde p^*
\end{equation}
This composition defines a natural map 
\begin{equation}
    \label{eq: critical Coulomb branch multiplication}
    \etap: \HO(\tilde \Rt, \phip{\Sh{W}}\DD\BoQ)\otimes \HO(\tilde \Rt, \phip{\Sh{W}}\DD\BoQ)\to \HO(\tilde \Rt, \phip{\Sh{W}}\DD\BoQ).
\end{equation}
In the next section, we identify this morphism will with the BFN convolution convolution product.

\begin{proposition}
\label{prop: crit convolution structure}
    Let $e$ be the identity element in $\pi_1(G)$. The morphism $\etap$ defines a $\pi_1(G)$-graded cohomologically graded associative algebra structure on 
    \[
    {}^{\phi}\Coul^T_{N,G}\coloneqq\HO_{G\fps_T}(\tilde \Rt_{G,N}, \phip{\Sh{W}^{(2)}}\DD\BoQ).
    \]
    The unit in this algebra is the fundamental class
    \[
    [\Rt_{N,G, e}\times N^*\lfps]\in \HO_{G\fps_T}(\Rt_{N,G, e}\times N^*\lfps,\DD\BoQ)\cong  \HO_{G\fps}(\tilde \Rt_{N,G, e}, \phip{W_{\lambda}}\DD\BoQ)
    \]
    of the subvariety $N^*\lfps \times \Rt_{N,G, e}\subseteq N^*\lfps \times \tilde \Rt_{N,G, e}$.
\end{proposition}

The proof of this theorem is analog to the proof of \cite[Theorem 3.10]{BFNII}. In fact, the proof is fundamentally easier as it solely relies on standard functoriality and compatibility properties of pullback and pushforward in critical cohomology. In particular, unlike the BFN construction, it does not rely on any ambient geometry. Since the proof is elementary, we do not give it here. An alternative proof follows from the identification of $\etap$ with the BFN multiplication via dimensional reduction, to which we come to now. In particular, the latter argument ensures that the critical convolution structure introduced above is independent on the choice of approximations of the spaces involved in the correspondence. In \S\ref{subsec: Coulomb branch action on quasimaps} we will use the presentation of the Coulomb branch provided by Proposition \ref{prop: crit convolution structure} to construct an action of the Coulomb branch on $\HO_T(\QM^{\xi}(X), \varphi_{\QM^{\xi}(X)})$. The proof that the action map is compatible with the Coulomb branch multiplication can be easily rearranged to prove \eqref{prop: crit convolution structure} by solely relying on operations in critical cohomology.

\begin{remark}
\label{rem: virtual crit Coulomb}
    Paralleling Remark \ref{rem: virtual Coulomb branch}, if $G,N$ come from a quiver gauge theory, we define 
    \[
    {}^{\phi}\Coul^{T, \vir}_{N,G}\coloneqq\bigoplus_{ \dd\in \BoN^{Q_0}}\HO_{G\fps_T}(\tilde \Rt_{N,G, \dd}, \phip{\Sh{W}}\DD\BoQ[-2s_{\dd}]),
    \]
    where $s_{\dd}=\vv^T\dd-\vv^T\Q^T\dd$. For the reasons explained in Remark \ref{rem: virtual Coulomb branch}, this is still a $\pi_1(G)$-graded cohomologically graded associative algebra. By Theorem \ref{thm: dr Coulomb agebra}, it will be clear that ${}^{\phi}\Coul^{T, \vir}_{N,G}$ dimensionally reduces to $\Coul^{T, \vir}_{N,G}$.

\end{remark}

\begin{remark}
\label{rem: ideal coulomb branch mult}
    Diagram \eqref{eq: first two maps critical ravolo correspondence} can be rewritten more compactly as follows. Set $M\lfps=N^*\lfps\times N\lfps$ and denote its closed points by $(w^*, w)$. We also set
    \[
    \tilde\Rt^{M\lfps}_{G,N}= \left(M\lfps \times N\fps\right) \times_{G\fps}G\lfps \times_{G\fps} G\lfps \times (N^*\lfps \times N\fps).
    \]
    This is naturally a vector bundle over $\tilde\Rt_{G,N}$ with fiber $M\lfps$. We assign it the potential $\Sh{W}^{\xi}+\Sh{Q}$, where $\Sh{Q}=\Res\Tr(w^*w)$. This is just the pullback of the potential $\Sh{W}^{\xi}$ on the base $\tilde\Rt_{G,N}$ plus the natural quadratic potential on the fiber $M\lfps$.
    
    Let now $\tilde m' : \tilde \Rt^{(2),1}_{G,N}\to\tilde\Rt^{M\lfps}_{G,N}$ be the morphism defined by
    \[
    [n_1, g, n_2^*, n_2, g', n_3^*, n_3]\mapsto [g'gn_3^*-gn_2^*, gn_2-n_1, n_1, g'g, n_3^*, n_3].
    \]
    By the same argument of Lemma \ref{lemma: critical raviolo 1}, it follows that $(\Sh{W}^{\xi}+\Sh{Q})\circ \tilde m'=\Sh{W}_1$. So we have a natural correspondence 
    \begin{equation}
        \label{eq: ideal critical raviolo}
        \begin{tikzcd}
        \tilde\Rt_{G,N}\times \tilde\Rt_{G,N} & \tilde\Rt^{(2),0}_{G,N}\arrow[l, swap, "\tilde p"]\arrow [r, "\tilde q_1"] & \tilde\Rt^{(2),1}_{G,N} \arrow[r, "\tilde m'"]& \tilde\Rt^{M\lfps}_{G,N}.
    \end{tikzcd}
    \end{equation}
    that is compatible with the potentials on these spaces.     Ideally, this would be the critical analog of the BFN correspondence. In fact, one could then just define the convolution algebra structure on $\HO_{G\fps}(\tilde\Rt_{G,N}, \phip{\Sh{W}}\DD\BoQ)$ as the composition
    \begin{equation}
        \label{eq: ideal critical moltiplication}
        c \circ  \tilde m'_* \circ (\tilde q^*)^{-1} \circ \tilde p^*,
    \end{equation}
    where $c$ is the dimensional reduction isomorphism 
    \[
    \HO_{G\fps_T}(\tilde\Rt^{M\lfps}_{G,N}, \phip{\Sh{W}+Q}\DD\BoQ)\xrightarrow[\cong]{c} \HO_{G\fps_T}(\tilde\Rt_{G,N}, \phip{\Sh{W}}\DD\BoQ).
    \]
    The reason why we we opt for the more cumbersome diagram \eqref{eq: critical Coulomb branch correspondence} is due to some technical details with the finite dimensional approximations of these maps that can be circumvented by ordering the maps differently (see in particular the proof of Proposition \ref{prop: Coulomb action def} and the choice of finite type approximations there). However, assuming the well-definiteness of all the maps in \eqref{eq: ideal critical moltiplication}, then the latter would be equal as \eqref{eq: critical Coulomb branch multiplication}. In fact, the statement would easily follow from  the commutativity of the diagram
    \[
    \begin{tikzcd}
        \tilde\Rt_{G,N}\times \tilde\Rt_{G,N} & \tilde\Rt^{(2),0}_{G,N}\arrow[l, swap, "\tilde p"]\arrow [r, "\tilde q_1"] & \tilde\Rt^{(2),1}_{G,N}\arrow[r, "\tilde j"]\arrow[rr, bend left, "\tilde m'"]&  \tilde\Rt^{(2),2}_{G,N} \arrow[r]& \tilde\Rt^{M\fps}_{G,N}
        \\
        & & & \tilde\Rt^{(2),3}_{G,N}\arrow[u, swap, "\tilde q_2"]\arrow[d, "\tilde q_3"]\arrow[r]& \tilde\Rt^{N^*\fps}_{G,N}\arrow[u]\arrow[d]
        \\
        & & & \tilde\Rt^{(2),4}_{G,N} \arrow[r, "\tilde m"] & \tilde\Rt^{}_{G,N}
\end{tikzcd}
    \]
    whose squares are Cartesian, together with the fact that the dimensional reduction isomorphism $c$ is induced by the rightmost vertical maps (by the same argument of Lemma \ref{lemma: critical raviolo 2}).

    In conclusion, the reason why we work with \eqref{eq: critical Coulomb branch correspondence} and not with \eqref{eq: ideal critical raviolo} is rather technical, and we strongly believe that, with a more powerful sheaf theory at hand, equation \eqref{eq: ideal critical moltiplication} could be taken as the actual definition of the multiplication on ${}^{\phi}\Coul_{N,G}$. The same simplification could also be applied to the definition of the Coulomb branch action on $\HO_T(\QM^{\xi}(X), \varphi_{\QM^{\xi}(X)})$ defined in \S\ref{subsec: Coulomb branch action on quasimaps}.
\end{remark}

\subsection{Identification with the BFN Coulomb branch}

In this section, we show that the critical version of the Coulomb branch introduced in \S\ref{subsec: critical Coulomb branch convolution algebra} reduces, via dimensional reduction, to the BFN Coulomb branch.

\begin{theorem}
    \label{thm: dr Coulomb agebra}
    The dimensional reduction isomorphism
    \[
     dr: \HO_{G\fps_T}(\tilde \Rt_{G,N}, \phip{\Sh{W}}\DD\BoQ)\xrightarrow{\cong} \HO_{G\fps_T}(\Rt_{G,N},\DD\BoQ)
    \]
    from Proposition \ref{prop: dr Coulomb vector sp} induces an isomorphism of $\pi_1(G)$-graded, cohomologically graded, unital algebras 
    \[
    \left(\HO_{G\fps-T}(\tilde \Rt_{G,N}, \phip{\Sh{W}}\DD\BoQ), \etap\right)\cong \left(\HO_{G\fps_T}(\Rt_{G,N},\DD\BoQ), \eta\right)
    \]
    and hence a canonical identification of the left hand side with he BFN Coulomb branch.
\end{theorem}
The proof of this theorem requires some notation. Firstly, we introduce the critical analog of the space $\Tt_{G,N}$ used in the BFN construction, cf. \S\ref{subsec: Coulomb branch basics}. To this end, let 
\[
\tilde \Tt_{G,N}=N\fps \times_{G\fps} G\lfps   \times(N\lfps\times N^*\lfps)
\]
We denote a point in the prequotient of $\tilde \Tt_{G,N}$ as a tuple $(n, g, n_2^*,n_2)$, where $n_2^*$ and $n_2$ are meromorphic while $n_1$ is holomorphic. In the diagrammatic language introduced in the previous section, we represent $\tilde \Tt_{G,N}$ by
\[
\tilde \Tt_{G,N}=
\begin{tikzcd}
    \Sh{V}\arrow[d, "n_1"]\arrow[r, dashed, "g"] & \Sh{V}'\arrow[d, dashed, "n_2"]
    \\
    \Sh{V}\arrow[r, dashed, "g"] & \Sh{V}'\arrow[u, shift left, dashed, "n_2^*"]
\end{tikzcd}
\]
where the dashed arrows serve to denote meromorphic fields.
On this space, we introduce the quadratic potential 
\[
Q: \tilde \Tt_{G,N}\to \BoC\qquad ([n_1,g], n_2^*, n_2)\mapsto \Res\Tr(n_2^* n_2).
\]
and consider the closed immersion 
\[
i: \tilde \Rt_{G,N}\hookrightarrow \tilde \Tt_{G,N} \qquad ([n_1,g], n_2^*, n_2)\mapsto ([n_1,g], n_2^*, n_2) ([n_1, g], gn_1-n_2, n_2^*)
\]
Diagrammatically, we represent it as 
\[
\begin{tikzcd}
    \Sh{V}\arrow[d, "n_1"]\arrow[r, dashed, "g"] & \Sh{V}'\arrow[d, "n_2"]
    \\
    \Sh{V}\arrow[r, dashed, "g"] & \Sh{V}'\arrow[u, shift left, dashed, "n_2^*"]
\end{tikzcd}
\mapsto 
\begin{tikzcd}
    \Sh{V}\arrow[d, "n_1"]\arrow[r, dashed, "g"] & \Sh{V}'\arrow[d,dashed, "n_2-gn_1"]
    \\
    \Sh{V}\arrow[r, dashed, "g"] & \Sh{V}'\arrow[u, shift left, dashed, "n_2^*"]
\end{tikzcd}
\]
Notice that we have $i^*Q=\Sh{W}$, where the right hand side is the Coulomb branch potential \eqref{eq: Coulomb branch potential}. Therefore, the potentials on $\tilde \Rt_{G,N}$ and $\tilde \Tt_{G,N}$ are compatible along $j$. 
\begin{remark}
The same argument of Lemma \ref{lemma: critical raviolo 2} gives a dimensional reduction isomorphism 
\[
\HO_{G\fps}(\Tt_{G,N}, \DD\BoQ)\xrightarrow{\cong} \HO_{G\fps}(\tilde \Tt_{G,N}, \phip{Q}\DD\BoQ) .
\]
fitting\footnote{Both the horizontal and vertical maps are pushforward morphisms in critical cohomology, so commutativity of the diagram reduces to associativity of pushforward in critical cohomology.} in a commutative square 
\[
\begin{tikzcd}
    \HO_{G\fps}(\Rt_{G,N}, \DD\BoQ)\arrow[r] \arrow[d, "\cong"] & \HO_{G\fps}(\Tt_{G,N}, \DD\BoQ) \arrow[d, "\cong"]
    \\
\HO_{G\fps}(\tilde \Rt_{G,N}, \phip{\Sh{W}}\DD\BoQ)\arrow[r]  & \HO_{G\fps}(\tilde \Tt_{G,N}, \phip{Q}\DD\BoQ),
\end{tikzcd}
\]
where the horizontal maps are pushforward morphisms. Hence, the pair $(\tilde \Tt_{G,N}, Q)$ can be viewed as the critical version of $\Tt_{G,N}$. However, as observed in \cite[Rem. 3.9.]{BFNII}, it does not seem natural to equip $\HO_{G\fps}(\Tt_{G,N}, \DD\BoQ)$ with a convolution algebra structure. The critical description does not seem to help in this regard.
\end{remark}

We further define 
\begin{align*}
\tilde \Tt^{(2),0}_{G,N}&= N\fps\times_{G\fps} G\lfps\times  N^*\lfps\times N\fps \times G\lfps \times  (N^*\lfps \times N\lfps)
\\
\tilde \Tt^{(2),1}_{G,N}&= N\fps\times_{G\fps} G\lfps\times_{G\fps} \left(N^*\lfps \times N\fps \times G\lfps\right) \times  (N^*\lfps \times N\lfps).
\end{align*}
Diagrammatically, we have 
\[
\tilde \Tt^{(2,0)}_{G,N}=
\begin{tikzcd}
    \Sh{V}\arrow[d, "n_1"]\arrow[r, dashed, "g"] & \Sh{V}'\arrow[d, "n_2"]\arrow[r, dashed, "g'"] & \Sh{V}''\arrow[d, dashed, "n_3"]
    \\
    \Sh{V}\arrow[r, dashed, "g"] & \Sh{V}'\arrow[r, dashed, "g'"]\arrow[u, shift left, dashed, "n_2^*"]& \Sh{V}''\arrow[u, shift left, dashed, "n_3^*"]
\end{tikzcd}
\]
Consider now the maps
\[
\begin{tikzcd}
       \tilde\Rt_{G,N}\times \tilde \Tt_{G,N}&\tilde\Tt^{(2),0}_{G,N}\arrow[l, swap, "\tilde p"]\arrow[r, "\tilde q_1"] &\tilde\Tt^{(2),1}_{G,N}
\end{tikzcd}
\]
defined by 
\[
\left(([n_1, g], n_2^*, n_2), ([n_2, g'], n_2^*,  n_2) \right)\mapsfrom \left([n_1, g], n_2^*, n_2, g', n_3^*, n_3\right)\mapsto [n_1, g, n_2^*, n_2, g', n_3^*, n_3].
\]
Informally speaking, this map simply separates the two squares in the diagram above. 
The potential on $\tilde\Tt^{(2,0)}_{G,N}$ is set to be the pullback of the one on $\tilde\Rt_{G,N}\times \tilde \Tt_{G,N}$. Explicitly, it is given by 
\[
\Res\Tr(n_2^*(gn_1-n_2))+\Res\Tr(n_3^*n_3).
\]

To remove clutter, we drop the subscripts $G,N$ and we also set $\bar\Rt=\Rt\times N^*\lfps$ and $\bar \Tt=\Tt\times N^*\lfps$. Points in $\bar\Rt$ and $\bar\Tt$ are denoted by $([n_1,g], n_2^*)$ We have an embedding $j: \bar \Rt \hookrightarrow \tilde \Rt$ sending $[n_1, g], n_2^*)\in \bar \Rt_{G,N}$ to $([n_1, g],n_2^*,gn_1)\in \tilde \Rt_{G,N}$. Similarly, let $k: \bar \Tt\hookrightarrow \tilde \Tt$ be the closed embedding whose sending $([n_1,g], n_2^*)$ to $([n_1,g], n_2^*, 0)$. These maps fit in a Cartesian square 
\begin{equation}
    \label{eq: various T and R spaces}
\begin{tikzcd}
   \bar \Rt\arrow[d, hookrightarrow, "i"] \arrow[r, hookrightarrow, "j"]& \tilde \Rt\arrow[d,hookrightarrow, "i"]
   \\
   \bar \Tt \arrow[r, hookrightarrow, "k"]& \tilde \Tt
\end{tikzcd}
\end{equation}
Notice also that the potentials on $\tilde \Tt$ and $\tilde \Rt$ pullback to zero on the subspaces $\bar \Tt$ and $\bar \Rt$.
 Consider now the following diagram

\begin{equation}
\label{eq: big diagram dimred Coulomb}
\begin{tikzcd}
    & \tilde\Rt\times \tilde\Rt \arrow[dd, hookrightarrow, near start, swap, "i\times \id"]& \tilde\Rt^{(2),0} \arrow[dd, hookrightarrow]\arrow[l, swap, "\tilde p"] \arrow[r, "\tilde q_!"] & \tilde\Rt^{(2)}\arrow[dd, hookrightarrow]
    \\
    \bar\Rt\times \bar\Rt\arrow[ur, hookrightarrow, "j\times j"] \arrow[dd, hookrightarrow, "i\times \id"]& \tilde p^{-1}( \bar\Rt\times \bar\Rt) \arrow[dd, shift left=5pt, hookrightarrow]\arrow[l, swap, "\tilde p"] \arrow[r, "\tilde q_1"] \arrow[ur, hookrightarrow]& \tilde q_1((\tilde p)^{-1}(\bar\Rt\times \bar\Rt)\arrow[dd, shift left=5pt, hookrightarrow]\arrow[ur]
    \\
    & \tilde\Rt \times \tilde\Tt & \tilde\Tt^{(2),0} \arrow[l, swap, "\tilde p"] \arrow[r, "\tilde q"] &\tilde\Tt^{(2)}
    \\
    \bar\Rt \times \bar\Tt\arrow[ur, "k\times j"]&  \bar\Rt \times (G\lfps\times N^*\lfps) \arrow[l, swap, "p"] \arrow[r, "q"]\arrow[ur, hookrightarrow, "l"]&  \bar\Rt \times_{G\fps}(G\lfps\times N^*\lfps)\arrow[ur]
\end{tikzcd}
\end{equation}
The maps on the leftmost face are induced by \eqref{eq: various T and R spaces} via pullback. The horizontal maps on the back face are defined above and in \S\ref{subsec: critical Coulomb branch convolution algebra}. The vertical ones are instead obtained by pullback and taking $G\fps$ quotients. All the remaining maps are also defined by taking cartesian squares (those in the leftmost cube) and $G\fps$ quotients (those in the rightmost cube).
All the squares are Cartesian too. Additionally, all the potentials defined on the back face become zero after restriction to the front face.

By construction, all spaces on the front page are vector bundles with fibers $N^*\lfps^2$ (recall that all spaces with an overline are bundles with fiber $N^*\lfps$). Passing to the zero sections of these bundles, we recover the diagram 
\begin{equation*}
\begin{tikzcd}
    \Rt\times \Rt \arrow[d, hookrightarrow, "i\times \id"]& p^{-1}(\Rt\times \Rt) \arrow[d, hookrightarrow]\arrow[l, swap, "\tilde p"] \arrow[r, "\tilde q"] & q(p^{-1}(\Rt\times \Rt))\arrow[d, hookrightarrow]
    \\
    \Rt\times \Tt & \Rt\times  G\lfps \arrow[l, swap, " p"] \arrow[r, "q"] & \Rt \times_{G\fps}  G\lfps
\end{tikzcd}
\end{equation*}
which is precisely the two leftmost squares of diagram \eqref{eq: Coulomb branch mult diagram} used in the BFN construction of the convolution product. We can now complete the proof of the main result of this section.

\begin{proof}[Proof of Theorem \ref{thm: dr Coulomb agebra}] The units match by definition, so it suffices to prove compatibility of the multiplications. The proof is mostly a standard compatibility argument for sheaf operations in Cartesian squares. The main subtle point is the compatibility involving the BFN pullback \eqref{eq: gysin pullback Coulomb}, so we start from there. From the left square of the bottom face of diagram \eqref{eq: big diagram dimred Coulomb} we get a commutative diagram
\[
\begin{tikzcd}
    (k\times j)_*(k\times j)^! \arrow[d]\arrow[r] &  \tilde p_* l_*l^! \tilde p^*\arrow[d]
    \\
   \id \arrow[r] &   \tilde p_*  \tilde p^*
\end{tikzcd}
\]
Applying $\DD\BoQ_{\tilde \Rt\times \tilde \Tt}$ to the right and the functor $(i\times \id)^!$ to the left we obtain 
\[
\begin{tikzcd}
    (j\times j)_*\DD\BoQ_{\bar\Rt\times \bar\Rt} \arrow[d]\arrow[r] &   (i\times id)^!\tilde p_* l_*l^! \tilde p^*\DD\BoQ_{\tilde \Rt\times \tilde \Tt}\arrow[d]
    \\
   \DD\BoQ_{\bar \Rt\times \bar \Rt} \arrow[r] &   (i\times id)^!\tilde p_* \tilde p^*\DD\BoQ_{\tilde \Rt\times \tilde \Tt}
\end{tikzcd}
\]
By construction, the  leftmost map induces the dimensional reduction isomorphism of Proposition \ref{prop: dr Coulomb vector sp}. On the other hand, using smoothness of $\tilde \Rt\times \tilde \Tt $ we can identify the bottom map with the unit map
\[
\DD\BoQ_{\tilde \Rt\times \tilde \Rt} \to   \tilde p_*  \tilde p^*\DD\BoQ_{\tilde \Rt\times \tilde \Rt}=\tilde p_*  \DD\BoQ_{\tilde \Rt^{(2),0}}[2\infty_{G,N}]
\]
Additionally smoothness of (all) the spaces on the bottom face of \eqref{eq: big diagram dimred Coulomb} implies that the top arrow agrees with the application of the functor $(i\times \id)^!$ to
\[
\DD\BoQ_{\tilde\Rt\times \tilde\Tt}\to \tilde p_* \tilde p^* \DD\BoQ_{\tilde\Rt\times \tilde\Tt}
\]
which induces the BFN pullback \eqref{eq: gysin pullback Coulomb} (more precisely its pullback along the trivial fibration $\bar\Rt\times\bar\Rt\to \Rt\times\Rt$). Therefore, applying the vanishing cycle functor $\phip{\Sh{W}\boxplus Q}$ we arrive at the commutative diagram
\[
\begin{tikzcd}
    (j\times j)_*\DD\BoQ_{\bar\Rt\times\bar\Rt} \arrow[d]\arrow[r] &   (j\times j)_* \tilde p_{*}\DD\BoQ_{p^{-1}(\Rt\times\Rt)}[2\infty_{G,N}]\arrow[d]\arrow[dr]
    \\
   \phip{\Sh{W}\boxplus Q}\DD\BoQ_{\tilde \Rt\times \tilde \Rt} \arrow[r] &   \phip{\Sh{W}\boxplus Q} \tilde p_*\tilde p^*\DD\BoQ_{\tilde \Rt\times \tilde \Rt} \arrow[r]  & \tilde p_* \phip{\Sh{W}\boxplus Q\circ \tilde p}  \DD\BoQ_{\tilde \Rt^{(2),0}}[2\infty_{G,N}]
\end{tikzcd}
\]
Taking derived global sections of the outer frame, contracting the trivial $N^*\lfps$-directions, and using the conventions involving the shifts, we arrive to the commutative diagram 
\begin{equation}
    \label{eq: proof dimred Coulomb 1}
\begin{tikzcd}[column sep=small]
    \HO_{G\fps}(\tilde \Rt,\phip{\Sh{W}}\DD\BoQ)\otimes \HO_{G\fps}(\tilde \Rt,\phip{\Sh{W}}\DD\BoQ) \arrow[r, "\tilde p^*"]   &  \HO_{G\fps\times G\fps}(\tilde \Rt^{(2),0},\phip{\Sh{W}\boxplus \Sh{W}\circ \tilde p}\DD\BoQ)
    \\
    \HO_{G\fps}(\Rt,\DD\BoQ)\otimes \HO_{G\fps}(\Rt,\DD\BoQ) \arrow[r, "p^!"]  \arrow[u, "\cong"]&  \HO_{G\fps\times G\fps}(p^{-1}(\Rt\times \Rt),\DD\BoQ)\arrow[u, "\cong"]
\end{tikzcd}
\end{equation}
On the other hand, we clearly have a commutative diagram
\begin{equation}
    \label{eq: proof dimred Coulomb 2}
\begin{tikzcd}[column sep=small]
      \HO_{G\fps\times G\fps}(\tilde \Rt^{(2),0},\phip{\Sh{W}_0}\DD\BoQ) \arrow[r, "(\tilde q_1^*)^{-1}"] &  \HO_{G\fps}(\tilde \Rt^{(2),1},\phip{\Sh{W}_1}\DD\BoQ)
      \\
      \HO_{G\fps\times G\fps}(p^{-1}(\Rt\times \Rt),\DD\BoQ)  \arrow[r, "(q^*)^{-1}"]\arrow[u, "\cong"] &  \HO_{G\fps}(q(p^{-1}(\Rt\times \Rt)),\DD\BoQ) \arrow[u, "\cong", "a"']
\end{tikzcd}
\end{equation}
The new label $a$ will be useful to keep track of the morphisms. Now consider the diagram
\[
\begin{tikzcd}
    \tilde \Rt^{(2),1}\arrow[r, "\tilde j"] & \Rt^{(2),2}
    \\
    \tilde q_1((\tilde p)^{-1}(\bar\Rt\times \bar\Rt)\arrow[r, hookrightarrow, "\bar j"]\arrow[u] &\Rt^{(2),3}\arrow[u, hookrightarrow, "\tilde q_2"]
\end{tikzcd}
\]
where the map $\bar j$ is defined via $[n_1, g, n_2^*, g', n_3^*]\mapsto [g'gn_3^*-gn_2^*, n_1, g, g', n_3^*, g'g n_1]$.
It is easy to check that this diagram commutes using the definition of morphisms $\tilde q_1$ and $\tilde j$ from \S\ref{subsec: critical Coulomb branch convolution algebra}. Taking pushforwards in critical cohomology we get a commutative diagram 
\begin{equation}
    \label{eq: proof dimred Coulomb 3}
\begin{tikzcd}[column sep=small]
      &\HO_{G\fps}(\tilde \Rt^{(2),1},\phip{\Sh{W}_1}\DD\BoQ)\arrow[r, "\tilde j_*"] & \HO_{G\fps}(\tilde \Rt^{(2),2},\phip{\Sh{W}_1}\DD\BoQ)
      \\
      \HO_{G\fps}(q(p^{-1}(\Rt\times \Rt)),\DD\BoQ)\arrow[r, equal]\arrow[rr, bend right, "b"]\arrow[ur, "\cong"', "a"]&\HO_{G\fps}(\tilde q_1(\tilde p^{-1}(\Rt\times \Rt)),\DD\BoQ) \arrow[u, "\cong"]\arrow[r, "\bar j_*"] & \HO_{G\fps}(\tilde \Rt^{(2),3},\phip{\Sh{W}_3}\DD\BoQ)\arrow[u, "\tilde q_2"]
\end{tikzcd}
\end{equation}
where the equality comes from contracting the $N^*\lfps$ fibers. We introduce the label $b$ for future refference. Finally, consider the diagrams
\begin{equation}
\label{eq: two diagrams at once proof dr Coulomb}
    \begin{tikzcd}
        Y \arrow[r, hookrightarrow, "j'"]\arrow[d, "\tilde m"]& \tilde \Rt^{(2),3}\arrow[d, "\tilde m"]
        \\
        \bar \Rt\arrow[r, hookrightarrow, "j"] &  \tilde \Rt
    \end{tikzcd}
    \qquad
    \begin{tikzcd}
        & \tilde \Rt^{(2),3}\arrow[dd, bend left, "\tilde q_3"]
        \\
        \tilde q_1((\tilde p)^{-1}(\bar\Rt\times \bar\Rt)\arrow[ur, hookrightarrow, "\bar j"]\arrow[r, hookrightarrow, "\bar j'"]\arrow[d, "t'"] &X\arrow[u, "\cong", "s"']\arrow[d, "t"]
        \\
        Y \arrow[r, hookrightarrow, "j'"]& \tilde \Rt^{(2),4}
    \end{tikzcd}
\end{equation}
Here $Y$ is the sub-bundle of $\tilde q_1((\tilde p)^{-1}(\bar\Rt\times \bar\Rt)$ consisting of tuples of the form $[n_1, g, 0, g', n_3^*]$ (as opposed to $[n_1, g, n_2^*, g', n_3^*]$ with arbitrary $n_2^*$ for an element in $\tilde q_1((\tilde p)^{-1}(\bar\Rt\times \bar\Rt)$); we denote its points by $[n_1, g, g', n_3^*]$.  The maps $t'$ and $j'$ are given by 
\[
[n_1, g, n_2^*, g', n_3^*]\mapsto [n_1, g, g', n_3^*]\mapsto [n_1, g, g', n_3^*, g'g n_1].
\]
It is easy to check that the leftmost diagram commutes. As for the right diagram, set
\[
X= N\fps\times_{G\fps}G\lfps \times_{G\fps} N^*\lfps\times_{G\fps} G\lfps \times (N^*\lfps\times N\fps)
\]
and define the maps $s$ and $t$ by 
\[
[g'g n_3^*-gn_2^*, n_1, g, g', n_3^*, n_3]\mapsfrom [n_1, g, n_2^*,g, n_3^*, n_3]\mapsto [n_1, g,g, n_3^*, n_3].
\]
It is clear that $s$ is an isomorphism. We set $\bar j'=s^{-1}\circ \bar j$. It is now easy to check that the bottom square is Cartesian. We now claim that the diagram
\begin{equation}
\label{eq: proof dimred Coulomb 4}
        \begin{tikzcd}
        & \HO_{G\fps}(\tilde q_1((\tilde p)^{-1}(\bar\Rt\times \bar\Rt),\DD\BoQ)\arrow[d, "(t'^*)^{-1}"]\arrow[r, "\bar j_*"]\arrow[dr, "A"]& \HO_{G\fps}(\tilde \Rt^{(2),3}, \phip{\Sh{W}_3}\DD\BoQ)\arrow[d, "(\tilde q_3^*)^{-1}", "\cong"']
        \\
        \HO_{G\fps}(qp^{-1}(\Rt\times\Rt), \DD\BoQ)\arrow[equal]\arrow[d, "m"]\arrow[dr, phantom, "B"]\arrow[ur, "b"]\arrow[r, equal]&\HO_{G\fps}(Y, \DD\BoQ)\arrow[r, "j_*"]\arrow[d, "\tilde m"] \arrow[dr, phantom, "C"]& \HO_{G\fps}(\tilde \Rt^{(2),4}, \phip{\Sh{W}_4}\DD\BoQ)\arrow[d, "\tilde m"]
        \\
        \HO_{G\fps}(\Rt, \DD\BoQ)\arrow[r, equal]& \HO_{G\fps}(\bar \Rt\DD\BoQ)\arrow[r, "j"] &\HO_{G\fps}(\tilde \Rt, \phip{\Sh{W}}\DD\BoQ)
    \end{tikzcd}
\end{equation}
commutes. Commutativity of square C follows from commutativity of the left diagram in \eqref{eq: two diagrams at once proof dr Coulomb}. In B, the horizontal equalities are given by contraction of the fiber $N\lfps$ and hence commutativity reduces to compatibility of pullback and pushforward in Cartesian squares. For  A, we argue as follows. Consider the second diagram in \eqref{eq: two diagrams at once proof dr Coulomb}. Since $s$ is an isomorphism, from the upper triangle we deduce that $\bar j_*=s_* \circ \bar j'_*= (s^*)^{-1}\bar j'_*$. Since its bottom square is Cartesian, we also have $t^* \circ j'_*= \bar j'_*\circ t'^*$. Finally, from the remaining triangle we get $s^*\circ \tilde q_3^*=t^*$. The sought after commutativity follows by combining thee three equations.

The proof now follows by combining the definitions of $\eta$ and $\etap$ (see \S\ref{subsec: Coulomb branch basics} and \S\ref{subsec: critical Coulomb branch}) with diagrams \eqref{eq: proof dimred Coulomb 1}, \eqref{eq: proof dimred Coulomb 2}, \eqref{eq: proof dimred Coulomb 3}, and \eqref{eq: proof dimred Coulomb 4}.
\end{proof}

\section{Cech model for quasimaps}

\subsection{The infinite dimensional Čech model}
\label{subsec: infinite Čech model}
As we have seen in Theorem \ref{theorem: shifted Yangian action on critical cohomology}, the critical locus description of the quasimap moduli space described in Section \ref{sec: quiver critical description quasimaps} is useful for constructing the action of the shifted Yangian associated to the underlying quiver $Q$ on the quasimap critical cohomology. There is an alternative critical locus description in which the Coulomb branch action will become manifest. It will be verified in Appendix \ref{app: derived quasimaps} and Part \ref{part: quasimap shifted symplectic} that either critical locus description can be used as a model for the moduli space of quasimaps as a $(-1)$-shifted symplectic stack (and that the orientation data agree, see Appendix \ref{ap: orientations}), so that we can use either as a model for the critical cohomology, see Corollary \ref{cor: DT are the same}

\subsubsection{Reminder on vector bundles on $\mathbb{P}^1$} \label{subsubsect: Čech review}
We consider $\mathbb{P}^1$ with uniformizing coordinate $t$, $U_\infty = \text{Spec} \, \mathbb{C}[t^{-1}] \subset \mathbb{P}^1$ the $\mathbb{A}^1$ chart containing $\infty$, and $\mathbb{D}_0 = \text{Spec} \, \mathbb{C}[[t]]$ the formal neighborhood of $0 \in \mathbb{P}^1$; we have $U_\infty \cap \mathbb{D}_0 = \mathbb{D}^\times_0 = \text{Spec} \, \mathbb{C}((t))$, the formal punctured disc at zero. For brevity, set $\Sh{O} = \mathbb{C}[[t]]$, and $\Sh{K} = \mathbb{C}((t))$.

Describing bundles on $\mathbb{P}^1$ by transition functions, we have the well-known isomorphism 
\begin{equation} \label{eq: BunG uniformization}
\text{Bun}_G(\mathbb{P}^1) \simeq G(\Sh{O}) \backslash G(\Sh{K}) / G[t^{-1}],
\end{equation}
and we write a transition function as $g(t)$. $G(\Sh{O}) \times G[t^{-1}] \ni (g_0(t), g_\infty(t))$ acts on transition functions by $g(t) \mapsto g^{-1}_0(t) g(t) g_\infty(t)$. For a vector bundle $\mathscr{E}$ on $\mathbb{P}^1$, its coherent sheaf cohomology $H^\bullet(\mathbb{P}^1, \mathscr{E})$ is computed as the cohomology of the (formal) Čech complex\footnote{Observe that, by the splitting theorem for bundles on $\mathbb{P}^1$, it suffices to check this when $\mathscr{E} = \Sh{O}_{\mathbb{P}^1}(k)$.}
\[
\Gamma(U_\infty, \mathscr{E}) \oplus \Gamma(\mathbb{D}_0, \mathscr{E}) \to \Gamma(\mathbb{D}^\times_0, \mathscr{E}). 
\]
Trivializing $\mathscr{E}$ over each chart, the Čech differential is given in terms of the transition function by $\delta(s_\infty(t), s_0(t)) = s_\infty(t) - g^{-1}(t) \cdot s_0(t)$. 

\subsubsection{Main construction}
Now we return to the setting of $\QM^\xi(X)$. Denote $N \coloneqq \text{Rep}_Q(\vv, \ww)$, regarded as a representation of $G \coloneqq \GL_{\vv} = \prod_{i \in Q_0} \GL_{v_i}$. There is a natural map 
\[
\QM^\xi(X) \to \text{Bun}_G(\mathbb{P}^1)
\]
and we view the target via \eqref{eq: BunG uniformization} as collections $g(t) = (g_i(t))_{i \in Q_0}$ of transition functions on rank $v_i$ vector bundles $\Sh{V}_i$. Consider the space 
\[
\Sh{M}_{G, N} \coloneqq  N(\Sh{O}) \times_{G(\Sh{O})} G(\Sh{K}) \times_{G[t^{-1}]} (N^*(\Sh{K}) \times N[t^{-1}]). 
\]
There is a natural map $\Sh{M}_{G, N} \to \text{Bun}_G(\mathbb{P}^1)$; its fiber over $\Sh{V} = (\Sh{V}_i)_{i \in Q_0}$ can be interpreted intrinsically as 
\[
\Gamma(\mathbb{D}_0, \Sh{N}) \oplus (\Gamma(\mathbb{D}_0^\times, \Sh{N}^*) \oplus \Gamma(U_\infty, \Sh{N}))/\text{Aut}(\Sh{V})
\]
where $\Sh{N}$ is the vector bundle associated to $\Sh{V}$ via the representation $N$, and $\text{Aut}(\Sh{V}) \subseteq G
(\Sh{O}) \times G[t^{-1}]$ can be viewed as the stabilizer of the transition function. Explicitly, a point of $\Sh{M}_{G, N}$ may be described in terms of a tuple 
\[
( \Sh{B}_{3, 0}(t), \Sh{J}(t), g(t), \Sh{B}_2(t), \Sh{I}(t), \Sh{B}_{3, \infty}(t), \Sh{J}_\infty(t)).  
\]
acted on by $G(\Sh{O}) \times G[t^{-1}] \ni (g_0(t), g_\infty(t))$ via
\begin{align*}
& (g^{-1}_0(t) \Sh{B}_{3, 0}(t) g_0(t), \Sh{J}_0(t) g_0(t), g_0^{-1}(t) g(t) g_\infty(t),   \\
& g^{-1}_\infty(t) \Sh{B}_2(t) g_\infty(t), g^{-1}_\infty(t) \Sh{I}(t), g^{-1}_\infty(t) \Sh{B}_{3, \infty}(t) g_\infty(t), \Sh{J}_\infty(t) g_\infty(t)). 
\end{align*}
Let $\text{Bun}_{G, \infty}(\mathbb{P}^1)$ be the moduli stack of $G$-bundles on $\mathbb{P}^1$ with a trivialization of the fiber over $\infty$; the canonical map $\text{Bun}_{G, \infty}(\mathbb{P}^1) \to \text{Bun}_G(\mathbb{P}^1)$ is a principal $G$-bundle. Define $\widetilde{\Sh{M}}_{G, N}$ by the pullback square
\[
\begin{tikzcd}
    \widetilde{\Sh{M}}_{G, N} \arrow[r] \arrow[d] & \Sh{M}_{G, N} \arrow[d] \\ 
    \text{Bun}_{G, \infty}(\mathbb{P}^1) \arrow[r] & \text{Bun}_G(\mathbb{P}^1).
\end{tikzcd}
\]
Consider the closed substack 
\begin{equation} \label{eq: define MxiGN}
\Sh{M}^\xi_{G, N} \hookrightarrow \widetilde{\Sh{M}}_{G, N}
\end{equation}
of data satisfying the condition\footnote{Notice that we are not entitled to fix the evaluation of $\Sh{B}_2, \Sh{I}$ at $\infty$ yet, because they a priori are only defined over the formal punctured disk at zero. Instead, the rest of the evaluation data will enter through the potential $\Sh{W}$ below.} $\Sh{B}_{3, \infty}(\infty) = \xi_3$, $\Sh{J}_\infty(\infty) = \jmath$. Now consider the formal expression 
\begin{equation} \label{eq: formula for Čech potential}
\begin{split}
\Sh{W} & = \Res \tr (\xi_2 \Sh{B}_{3, \infty}(t)) + \Res \tr(\iota \Sh{J}_\infty(t)) \\
& - \Res \tr \Sh{B}_2(t)( \Sh{B}_{3, \infty}(t) - g^{-1}(t) \Sh{B}_{3, 0}(t) g(t)) \\ 
& - \Res \tr \Sh{I}(t)(\Sh{J}_{3, \infty}(t) - \Sh{J}_0(t) g(t))
\end{split}
\end{equation}
where $\text{Res}$ is the linear functional extracting the coefficient of $t^{-1}$. Notice that the torus $T$ introduced at the end of \S\ref{subsec: critical Coulomb branch} naturally acts on $\Mt^{\xi}_{N,G}$ and preserves the potential \eqref{eq: formula for Čech potential}. We will always work equivariantly with respect to this torus action.

\begin{lemma} \label{lemma: Čech potential is well-defined}
$\Sh{W}$ descends to a well-defined function on $\Sh{M}^\xi_{G, N}$ if and only if the quiver data $(\xi_2, \xi_3, \iota, \jmath)$ satisfy the moment map equation $\comm{\xi_2}{\xi_3} + \iota \jmath = 0$. 
\end{lemma}

\begin{proof}
    First, notice that by construction, the degree in $t$ of each factor in the formula for $\Sh{W}$ is bounded from below, so that the residue is well-defined. It suffices to then check that $\Sh{W}$ is $G(\Sh{O}) \times G_1[t^{-1}]$-invariant, where $G_1[t^{-1}] \subset G[t^{-1}]$ is the kernel of evaluation at $t = \infty$. For the last two lines of \eqref{eq: formula for Čech potential}, this is obvious by the definition of the group action. For the first line, we reason as follows. Expand 
    \[
    g_\infty(t) = 1 + \frac{g^{(1)}_\infty}{t} + O(t^{-2}).
    \]
Then upon acting by $G(\Sh{O}) \times G_1[t^{-1}]$, the first line of \eqref{eq: formula for Čech potential} is (upon dropping terms that do not contribute to the residue, and using an obviously similar notation for the expansion coefficients of $\Sh{B}_{3, \infty}, \Sh{J}_\infty$):
\begin{align*}
& \text{Res} \tr \xi_2(1 - g^{(1)}_\infty t^{-1})(\xi_3 + \Sh{B}^{(1)}_{3, \infty} t^{-1} )( 1 + g^{(1)}_\infty t^{-1}) \\
& + \text{Res} \tr  \iota (\jmath + \Sh{J}^{(1)}_\infty t^{-1})(1 + g^{(1)}_\infty t^{-1}) \\
& =  \tr \xi_2 \Sh{B}^{(1)}_{3, \infty} + \tr \iota \Sh{J}^{(1)}_\infty + \tr(( \comm{\xi_2}{\xi_3} + \iota \jmath) g^{(1)}_\infty).
\end{align*}
$\Sh{W}$ is invariant if and only if the third term in the last line vanishes, implying the result.
\end{proof}

We always assume that $(\xi_2, \xi_3, \iota, \jmath)$ satisfy the moment map condition and the stability condition, so that $\xi \in X$.
The following is the main result of this section.

\begin{theorem} 
\label{thm: Čech critical quasimap}
There is an isomorphism of algebraic stacks $\QM^\xi(X) \simeq \textnormal{Crit}(\Sh{W})$.
\end{theorem}
In Appendix \ref{app: derived quasimaps} it will be shown that this lifts to an equivalence of derived stacks. The proof consists essentially of unravelling the definition of $\text{Crit}(\Sh{W})$. 

\begin{proof}
The sought-after isomorphism will be constructed by filling in the top arrow in the diagram 
\[
\begin{tikzcd}
    \QM^\xi(X) \arrow[d] \arrow[r, "\sim", dashed] & \text{Crit}(\Sh{W}) \arrow[d] \\ 
    \text{Bun}_G(\mathbb{P}^1) \arrow[r, "\sim"] & G(\Sh{O}) \backslash G(\Sh{K}) / G[t^{-1}]
\end{tikzcd}
\]
where the bottom arrow is the uniformization isomorphism \eqref{eq: BunG uniformization} sending a bundle to its transition function. By definition, the fiber of $\QM^\xi(X)$ over a given bundle $\Sh{V}$ is the space of global sections
\[
H^0(\mathbb{P}^1, \Sh{N} \oplus \Sh{N}^*)
\]
evaluating to $\xi$ at infinity\footnote{Notice that the quasimap stability condition is automatically satisfied if we assume infinity is sent to a stable point of $\FX$.} and satisfying the moment map equation 
\[
\comm{\Sh{B}_2}{\Sh{B}_3} + \Sh{I} \Sh{J} = 0 \in H^0(\mathbb{P}^1, \mathscr{E}nd(\Sh{V}))
\]
taken up to automorphisms of $\Sh{V}$. Let us compare this with the critical point condition $\delta \Sh{W} = 0$. From \eqref{eq: formula for Čech potential}, varying $\Sh{W}$ with respect to $(\Sh{B}_2(t), \Sh{I}(t)) \in N^*(\Sh{K})$ gives the $N(\Sh{K})$-valued equations 
\begin{align*}
\Sh{B}_{3, \infty}(t) - g^{-1}(t) \Sh{B}_{3, 0}(t) g(t) & = 0 \\
\Sh{J}_\infty(t) - \Sh{J}_0(t)g(t) & = 0
\end{align*}
which are precisely the conditions for the Čech cochain $(\Sh{B}_{3, \infty}, \Sh{J}_\infty; \Sh{B}_{3, 0}, \Sh{J}_0)$ to define an element of $H^0(\mathbb{P}^1, \Sh{N})$. Likewise, the part of $\delta \Sh{W}$ involving $(\delta \Sh{B}_3, \delta \Sh{J})$ is 
\begin{align*}
\delta \Sh{W} & = \text{Res} \tr (\Sh{B}_2(t) - \xi_2) \delta \Sh{B}_{3, \infty}(t) + \text{Res} \tr (\Sh{I}(t) - \iota) \delta \Sh{J}_\infty(t) \\
& - \text{Res} \tr ( g^{-1}(t) \Sh{B}_2(t) g(t)) \delta \Sh{B}_{3, 0}(t) - \text{Res} \tr g^{-1}(t) \Sh{I}(t) \delta \Sh{J}_0(t) + \dots 
\end{align*}
Notice that, because we take $\Sh{W}$ to be defined on $\Sh{M}^\xi_{G, N}$, $(\delta \Sh{B}_{3, \infty}(t), \delta \Sh{J}_\infty(t)) \in t^{-1}N[t^{-1}]$. Vanishing of $\delta \Sh{W}$ then enforces the conditions: 
\begin{align*}
(\Sh{B}_2(t) - \xi_2, \Sh{I}(t) - \iota) & \in t^{-1}N^*[t^{-1}] \\ 
(g(t) \Sh{B}_2(t) g^{-1}(t), g(t) \Sh{I}(t)) & \in N^*(\Sh{O}). 
\end{align*}
These once again assert that $(\Sh{B}_2(t), \Sh{I}(t))$ extends to define an element of $H^0(\mathbb{P}^1, \Sh{N}^*)$; notice also that $(\Sh{B}_2, \Sh{B}_3, \Sh{I}, \Sh{J})$ evaluates to $(\xi_2, \xi_3, \iota, \jmath)$ at $\infty \in \mathbb{P}^1$ on the critical locus of $\Sh{W}$. Finally, the part of $\delta \Sh{W}$ involving variation of the transition function is 
\[
\delta \Sh{W} = \text{Res} \tr( \comm{\Sh{B}_2(t)}{g^{-1}(t) \Sh{B}_{3, 0}(t) g(t)} + \Sh{I}(t) \Sh{J}_0(t)g(t) ) g^{-1}(t) \delta g(t) + \dots 
\]
Using that $g^{-1}(t) \delta g(t) \in \Fg(\Sh{K})$ and the other critical point conditions, we see that we must have 
\[
\comm{\Sh{B}_2(t)}{\Sh{B}_{3, \infty}(t)} + \Sh{I}(t) \Sh{J}_\infty(t) = 0
\]
implying that the sections satisfy the moment map condition. The isomorphism $\QM^\xi(X) \to \text{Crit}(\Sh{W})$ then consists of sending $(\Sh{B}_3, \Sh{J})$ to the corresponding Čech cochain, restricting $(\Sh{B}_2, \Sh{I})$ to $U_\infty$ and then completing at $\mathbb{D}^\times_0$.
\end{proof}

Notice that the action of the torus $T$ on $\QM^\xi(X)$ extends to $\Mt^{\xi}_{G,N}$ and preserves the potential. Moreover, it is clear that the argument of the theorem holds $T$-equivariantly. Therefore, the statement of Theorem \ref{thm: Čech critical quasimap} holds $T$-equivariantly too.


\subsection{Finite type approximations of the Čech model}
\label{subsec: Finite type approximations of the Čech model}

A practical limitation of the Čech model discussed in \S\ref{subsec: infinite Čech model} is that it presents the moduli space of quasimaps, which is a scheme locally of finite type\footnote{And of finite type in each degree.} as the critical locus of a function on a space that is \emph{not} locally of finite type. Therefore, the usual sheaf theory package is not available—and indeed the analytic topology is not even naively defined on $\Sh{M}^\xi_{G, N}$. In fact, $\Sh{M}^\xi_{G, N}$ is naturally an infinite rank affine bundle over the moduli space $\Bun_{G, \infty}(\bbP^1)= G(\Sh{O})\backslash G\lfps/G_1[t^{-1}]$ parameterizing $G$-bundles on $\bbP^1$ together with a trivialization at $\infty\in \bbP^1$. More precisely, $\Sh{M}^\xi_{G, N}\to \Bun_{G, \infty}(\bbP^1)$ is a Tate bundle \cite{drinfeld2004infinitedimensionalvectorbundlesalgebraic} over a smooth scheme locally of finite type, which informally means that can be realized by taking pro and ind limit of finite rank affine bundles over the base. The goal of this section is to provide an explicit presentation of these finite type approximations together with a compatible approximation of the potential $\Sh{W}: \Sh{M}^\xi_{G, N}\to \BoC$. This will be crucial to obtain a presentation of the DT sheaf on $\QM^{\xi}(X)$ that is compatible with the Coulomb branch action, cf. \S\ref{subsec: DT sheaf Čech model}.

Firstly, notice that 
\begin{align}
\label{eq: presentation Čech ambeint space as tot}
\begin{split}
    \Sh{M}^{\xi}_{G,N} 
    \coloneqq&N\fps\times_{G\fps} G\lfps \times_{G_1[t^{-1}]}\left(N^*\lfps \times N^{\xi}[t^{-1}])\right)
    \\
    &=\Tot_{\Bun_{G,\infty}(\bbP^1)}\left(N\fps\times N_\infty^*\lfps \times N^\xi[t^{-1}]\right)
\end{split}
\end{align}
where, with slight abuse of notation, we are using the same notation for the vector bundles and their fibers, and we have added a subscript to $N^*(\Sh{K})$ to indicate that we regard it as a $G_1[t^{-1}]$-space. We also rewrite the potential \eqref{eq: formula for Čech potential} more compactly as 
\begin{equation}
    \label{eq: potential Čech compact version}
    \Sh{W}=\Res\left(n^*_{\infty}(t)\left(g^{-1}(t)n_0(t)-n_{\infty}(t)\right)\right)+\Res\Tr((\xi_2, \iota) n_{\infty}(t))
\end{equation}
where $(g(t), n_0(t), n^*_{\infty}(t),n_{\infty}(t))$ denotes a tuple in $ G\lfps \times N\fps \times N^*\lfps \times N^{\xi}[t^{-1}]$. 
\begin{remark}
    It may be convenient to visualize these data by using a diagram of the form 
    \[
    \begin{tikzcd}
    \Sh{V}|_{\DD_0}\arrow[d, "n_0"]\arrow[r, dashed, "g"] & \Sh{V}|_{U_{\infty}}\arrow[d, "n_{\infty}"]
    \\
    \Sh{V}|_{\DD_0}\arrow[r, dashed, "g"] & \Sh{V}|_{U_{\infty}}\arrow[u, shift left, dashed, "n^*"]
    \end{tikzcd}
    \]
    where a dashed arrow indicates that the map labelled is defined over the punctured disk $\DD^\times_0$.
\end{remark}

Fix a pair of non-negative integers $l ,m \in \mathbb{Z}_{\geq 0}$. We further assume $l \neq 0$. Recall that $N^{\xi}[t^{-1}]$ is the the $G_1[t^{-1}]$-variety $N^{\xi}[t^{-1}]=(\xi_3, \iota )+t^{-1}N[t^{-1}]$. Set 
\[
N^{\xi}_l[t^{-1}]=(\xi_3,\iota )+t^{-1}N[t^{-1}]/t^{-l}N[t^{-1}].
\]
This is also a representation of $G_1[t^{-1}]$, and hence induces a \emph{finite rank} vector bundle on $\Bun_{G, \infty}(\bbP^1)$ of the same rank. Similarly, we define 
\[
N^m\fps=N\fps/t^m\fps
\]
which is a $G\fps$-representation and hence induces another \emph{finite rank} vector bundle on $\Bun_{G, \infty}(\bbP^1)$.

It remains to truncate $N^*_\infty\lfps$. To begin with, notice that any element $g(t)\in G\lfps$ induces a map
\[
g^{-1}(t): N\fps\hookrightarrow N_\infty\lfps\qquad n_0(t)\mapsto g^{-1}(t) n_0(t),
\]
which descends to an injective morphism of vector bundles over $\Bun_{G, \infty}(\bbP^1)$, which we denote in the same way. Set 
\[
N^{*}_m\lfps\coloneqq \Ann\left(g^{-1}(t)t^m N\fps\right) \subset N^*\lfps 
\]
where the annihilator is induced by the residue pairing 
\[
\Res: N^{*}\lfps\times N\lfps\qquad (n(t),n^*(t))\mapsto \Res( n(t),n^*(t))
\]
Notice that $N^{*}_{m}\lfps$ is still a bundle of infinite rank. Now consider the sheaf on $\Bun_{G, \infty}(\bbP^1)$ fiberwise given by 
\[
N^{*}_{m, l}\lfps\coloneqq \Ann\left(g^{-1}(t)t^m N\fps\right) \cap \Ann(t^{-l} N[t^{-1}])\subset N^*\lfps.
\]
Let $\Bun_{G,\infty}^{\leq 0}(\bbP^1)$ be the open substack of $\Bun_{G,\infty}(\bbP^1)$ parametrizing stably negative bundles, that is bundles whose Birkhoff factorization only contains non-positive line bundles. Openness follows from upper semi-continuity of cohomology. We further denote by  $\Bun_{G,\infty, \dd}^{\leq 0}(\bbP^1)$ the connected component parametrizing bundles of degree $\dd\in \BoZ^{Q_0}$.
\begin{lemma}
\label{lemma: approxiamtion Čech model}
    For large enough $m, l$, the restriction of $N^{*}_{m, l}\lfps$ to the substack $\Bun_{G,\infty, \dd}^{\leq 0}(\bbP^1)\subseteq \Bun_{G,\infty}(\bbP^1)$ is a vector bundle of finite rank.
\end{lemma}

\begin{proof}
    Consider the Harder–Narasimhan stratification of $\Bun_{G}(\bbP^1)$. It induces, via pullback along the principal $G$-bundle $\Bun_{G,\infty}(\bbP^1)\to \Bun_{G}(\bbP^1)$ a stratification on $\Bun_{G,\infty}(\bbP^1)$, and so on the open substack $\Bun_{G,\infty, \dd}^{\leq 0}(\bbP^1)\subseteq \Bun_{G,\infty}(\bbP^1)$. In terms of the uniformization theorem, these resulting Harder–Narasimhan strata correspond to the $G\fps \times G[t^{-1}]$-orbits $G\fps t^{-\lambda}G[t^{-1}]$ for dominant cocharacters $\lambda$ with no negative entries and such that $|\lambda|=\dd$. In particular, there are finitely many strata. Moreover, by $G[t^{-1}]$-invariance of the residue pairing, it follows that $N^*\lfps_{m, l}$ has constant rank on each of these strata. Therefore, it suffices to compute the rank for $g(t)=t^{-\lambda}$ and show that the result does not jump as $\lambda$ varies across the strata. 
    
    Let $N=\bigoplus_{\mu} N_{\mu} $ be the weight decomposition of $N$ with respect to the maximal torus $S\subset G$. We have
    \[
    t^{\lambda}t^m N\fps= \bigoplus_{\mu} t^{\lambda+m}N_{\mu}\fps =  \bigoplus_{\mu} t^{\langle\lambda, \mu\rangle+m} N_{\mu}\fps
    \]
    Assume now that $m> \max_{\lambda, \mu}(-\langle\lambda, \mu\rangle)$, where $\mu$ ranges among the (finite set of) weights of $N$ and $\lambda$ among the (finite set of) strata of $\Bun_{G,\infty, \dd}^{\leq 0}(\bbP^1)\subseteq \Bun_{G,\infty}(\bbP^1)$, and that $l$ is sufficiently large. Computing the annihilator, we get 
    \[
     \Ann\left(g^{-1}(t)t^m N\fps\right)=\widehat\bigoplus_{n\geq -\langle\lambda, \mu\rangle-m} (N_{\mu})^* t^n \subset N^*\lfps.
    \]
    Similarly, we have 
    \[
    \Ann(t^{-l} N[t^{-1}])=\bigoplus_{\mu}\widehat\bigoplus_{n\leq l-2} (N_{\mu})^* t^n\subset N^*\lfps,
    \]
    and thus, choosing $m$ and $l$ large enough we have
    \[
    N^{*}_{m,l}=\bigoplus_{\mu}\bigoplus_{\langle\lambda, \mu\rangle+m \leq n \leq l-2} (N_{\mu})^* t^n.
    \]
    Therefore, we deduce that the rank of $N^*_{m, l}(\Sh{K})$ on the stratum $G\fps t^{-\lambda}G[t^{-1}]$ is given by
    \begin{align}
    \label{eq: computation dimension fibers Čech model}
    \begin{split}
        \sum_{\mu}\sum_{-\langle\lambda, \mu\rangle-m \leq n< l-2} \dim(N_{\mu})
        &=\sum_{\mu}(l -1+m)\dim(N_{\mu})+\sum_{\mu} \dim(N_{\mu})\langle\lambda, \mu\rangle
        \\
        &=(l-1+m)\dim(N)+\langle\lambda, \chi_{\det(N)}\rangle
    \end{split}
    \end{align}
    where $\chi_{\det(N)}$ is the unique $S$-weight of the determinant representation. To conclude that the right hand side only depends on $m, l$ and the degree $\dd$, it suffices to show that the map 
    \[
    X_*(S)\to \BoZ \qquad \lambda\mapsto \langle\lambda, \chi_{\det(N)}\rangle
    \]
    factors through the quotient
    \[
    X_*(S)/\{\text{coroots}\}\cong \pi_1(G).
    \]
    Indeed, any coroot $\lambda: \BoC^*\to S\subset G$ factors through a copy of $\SL(2)$ in $G$ so the claim reduces to the fact that the character of the determinant representation associated to an irreducible $\SL(2)$ representation is trivial. In conclusion, we have shown that for $m$ and $l$ large enough the sheaf $N^*\lfps_{m,l}$ has constant rank fibers across the Harder–Narasimhan strata, and hence is a \emph{finite rank} vector bundle on $\Bun_{G,\infty, \dd}^{\leq 0}(\bbP^1)$.
    \end{proof}
    
We can now complete the construction of the finite dimensional approximation of \eqref{eq: presentation Čech ambeint space as tot}. Define
\begin{equation}
    \label{eq: approximated Čech model}
    \Sh{M}^{\xi, m, l}_{G, N, \dd} \coloneqq\Tot_{\Bun^{\leq 0}_{G, \infty, \dd}(\bbP^1)}\left(N^m\fps\times N^*\lfps_{m,l} \times N^\xi_l[t^{-1}]\right)
\end{equation}
By Lemma \ref{lemma: approxiamtion Čech model}, it follows that for $l$ and $m$ large enough, $\Sh{M}^{\xi, m, l}_{G, N, \dd}$ is a finite rank vector bundle over $\Bun^{\leq 0}_{G, \infty, \dd}(\bbP^1)$, and hence is a smooth scheme of finite type.  Moreover, it follows from the construction that the restriction of the potential \eqref{eq: potential Čech compact version} to the subspace
\[
\Tot_{\Bun^{\leq 0}_{G, \infty, \dd}(\bbP^1)}\left(N\fps \times N^*\lfps_{m, l} \times N^\xi[t^{-1}]\right)\subset \Sh{M}^{\xi}_{G, N, \dd}
\]
factors through the quotient $\Sh{M}^{\xi, m, l}_{G, N, \dd}$, giving a commutative diagram 
\[
\begin{tikzcd}
    \Sh{M}^{\xi}_{G, N, \dd}\arrow[r, "\Sh{W}"] & \AA^1
    \\
    \Tot_{\Bun^{\leq 0}_{G, \infty, \dd}(\bbP^1)}\left(N^*\lfps_{m, l} \times N^\xi[t^{-1}]\times N\fps\right)\arrow[u, hookrightarrow] \arrow[r, twoheadrightarrow]& \Sh{M}^{\xi, m, l}_{G, N, \dd}\arrow[u, swap, "\Sh{W}^{m, l}"]
\end{tikzcd}
\]
To show that the pair $(\Sh{M}^{\xi, m, l}_{G, N, \dd}, \Sh{W}^{m,l})$ provides a finite dimensional approximation $(\Sh{M}^{\xi}_{G, N, \dd}, \Sh{W})$, it remains to prove the compute the critical locus of $\Sh{W}^{m, l}$.
\begin{proposition}
\label{prop: approx crit Čech}
    For large enough $m, l$, we have a natural isomorphism
    \[
    \crit(\Sh{W}^{m, l}) \simeq \crit(\Sh{W}).
    \]
\end{proposition}

Before beginning the proof, notice that our choice of stability condition on the quiver variety $X$ implies that $\crit(\Sh{W})$ lies entirely over $\text{Bun}^{\leq 0}_{G, \infty, \dd}$, see Lemma \ref{lemma negative vector bundle}.

\begin{proof}
    We will give the consideration for the critical point condition arising from variation of $n^*(t)$, as the others are similar. In the finite approximation, this critical locus condition reads
    \begin{equation}
        \label{eq: approximated critical condition}
        n_{\infty}(t)+t^{-l}(\text{Polynomial in }t^{-1})=g^{-1}(t)\left(n_0(t) +t^m(\text{Power series in $t$}\right).
    \end{equation}
    Denote by $\bar{n}_\infty(t)$ and $\bar{n}_0(t)$ the unique lifts of minimal degree to $N^\xi[t^{-1}]$, $N(\Sh{O})$ respectively. To prove the claim it suffices to check that, for $m$ and $l$ large enough, this equation implies 
    \begin{equation}
    \label{eq: full critical condition}
    \bar{n}_{\infty}(t)= g^{-1}(t)\bar{n}_0(t)
    \end{equation}
    Since this equation is $G\fps \times G[t^{-1}]$-equivariant, we may assume that $g=t^{-\lambda}$, for some $\lambda$ such that $G(\Sh{O})t^{-\lambda} G[t^{-1}]$ appears among the Harder–Narasimhan strata of $\Bun^{\leq 0}_{G, \dd}(\bbP^1)$. Let $n^{\mu}_{0,\infty}(t)$ be the component of $n_{0,\infty}(t)$ in the weight space $N_{\mu}\subset N$. Then the previous equation restricts to 
    \[
     \bar{n}^{\mu}_{\infty}(t)+t^{-l}(\text{Polynomial in $t^{-1}$})=t^{\langle\lambda, \mu\rangle}\bar{n}_0^{\mu}(t) +t^{\langle\lambda, \mu\rangle+m}(\text{Power series in }t)
    \]
    Since the left hand side has only negative modes in $t$, it follows that for $m\geq -\langle\lambda, \mu\rangle$ the power series on the right is zero. Similarly, for $l\gg 0$ the polynomial on the right hand side has to be zero. Running this argument for all weights $\mu$ and for all\footnote{These are finitely many in a given degree component.} strata indexed by $\lambda$, we deduce that for $m,l\gg 0$ equation \eqref{eq: approximated critical condition} implies \eqref{eq: full critical condition}.
\end{proof}

\section{Coulomb branch action on critical quasimaps}

\subsection{The quasimap DT sheaf in the Čech model}
\label{subsec: DT sheaf Čech model}
Fix a quiver $Q$ and dimension vectors $\vv,\ww\in \BoN^{Q_0}$. Let $X=\Nak(\ww,\vv)$ be the associated quiver variety and let $\xi\in X$ be a torus fixed point. In \S\ref{subsec: Shifted Yangian action}, we exploited the presentation of the moduli space of based quasimaps $\QM^{\xi}(X)$ as a quiver with potential to construct an action of the shifted Yangian $\Yang_{\mu}^{T}(\tilde Q, \tilde W)$ on the vanishing cycle cohomology $\HO_T(\QM^{\xi}(X), \varphi_{\QM^{\xi}(X)})$. In this section, we will use the Čech model developed in \S\ref{subsec: infinite Čech model} to obtain a different presentation of the DT sheaf $\varphi_{\QM^{\xi}(X)}$. This presentation will be crucial in the next section, where we will exploit it to construct an action of the Coulomb branch on $\HO_T(\QM^{\xi}(X), \varphi_{\QM^{\xi}(X)})$. Since the Čech model presents $\QM^{\xi}(X)$ as the critical locus of a function $\Sh{W}: \QM^{\xi}(X)\to \BoC$, in order to effectively apply sheaf theory we need to pass to the finite approximations of the Čech model discussed in \S\ref{subsec: Finite type approximations of the Čech model}. 

Let $\QM_{\dd}^{\xi}(X)\subset \QM^{\xi}(X)$ the connected component of degree $\dd$ quasimaps. It admits a canonical morphism to $\Bun_{G, \infty, \dd}(\bbP^1)$, the smooth scheme parametrizing degree $\dd$ principal $G=\prod_{i\in Q_0}\GL_{\dd_i}$-bundles on $\bbP^1$ equipped with a trivializtion at $\infty\in \bbP^1$. As discussed in  \S\ref{subsec: Finite type approximations of the Čech model}, for sufficiently large $l,m\in \BoN$, there exists a finite rank affine bundle 
\[
\Sh{M}^{\xi, m, l}_{G, N, \dd}=\Sh{M}^{\xi,m, l}_{G, N, \dd}=\Tot_{\Bun_{G, \infty}(\bbP^1)}\left(N^*\lfps_{m,l} \times N[t^{-1}]^{l}\times N\fps^{m}\right)
\]
together with a regular function 
\[
\Sh{W}^{m, l}: \Sh{M}^{\xi, m, l}_{G, N, \dd}\to \BoC
\]
such that $\crit(\Sh{W}^{m, l})\cong QM_{\dd}^{\xi}(X)$. Consider the sheaf 
\begin{equation}
\label{eq: approximated DT sheaf Čech model}
    \phip{\Sh{W}^{m, l}}\BoQ^{\vir}_{\Sh{M}^{\xi, m, l}_{G, N, \dd}}
\end{equation}
where, as usual, $\BoQ^{\vir}_{\Sh{M}^{\xi, m, l}_{G, N, \dd}}=\BoQ_{\Sh{M}^{\xi, m, l}_{G, N, \dd}}[\dim(\Sh{M}^{\xi, m, l}_{G, N, \dd})]$. 
Bu Proposition \ref{prop: approx crit Čech}, the critical locus of $\Sh{W}^{m, l}$ is supported on $\QM_{\dd}^{\xi}(X)$, so we may see \eqref{eq: approximated DT sheaf Čech model} as a complex in $D_c(\QM_{\dd}^{\xi}(X))$.
\begin{proposition}
    For all $l'\geq l\gg 0$ and $m'\geq m\gg 0$, we have a canonical isomorphism of sheaves over $\QM_{\dd}^{\xi}(X)$
    \[
    \left(\phip{\Sh{W}^{m, l}}\BoQ^{\vir}_{\Sh{M}^{\xi, m, l}_{G, N, \dd}}\right)|_{\QM_{\dd}^{\xi}(X)}\cong  \left(\phip{\Sh{W}^{m', l'}}\BoQ^{\vir}_{\Sh{M}^{\xi, m', l'}_{G, N, \dd}}\right)|_{\QM_{\dd}^{\xi}(X)}.
    \]
\end{proposition}
\begin{proof}
Set 
\[
\Sh{M}^{\xi, m',m,l', l}_{G, N, \dd} \coloneqq\Tot_{\Bun^{\leq 0}_{G, \infty}(\bbP^1)}\left(N^*\lfps_{m,l} \times N[t^{-1}]^{l'}\times N\fps^{m'}\right)
\]
The quotient maps 
\begin{align*}
    N[t^{-1}]^{l'}&=N[t^{-1}]/t^{-l'}N[t^{-1}]\twoheadrightarrow N[t^{-1}]/t^{-l'}N[t^{-1}] =N[t^{-1}]^{l}
    \\
    N\fps^{m'}&=N\fps/t^{m'}N\fps \twoheadrightarrow N\fps/t^{m}N\fps=N\fps^{m}
\end{align*}
are, respectively, $G_1[t^{-1}]$ and $G\fps$ equivariant and hence induce a fiberwise-surjective morphism of bundles $j: \Sh{M}^{\xi, m',m,l', l}_{G, N, \dd}\twoheadrightarrow \Sh{M}^{\xi, m, l}_{G, N, \dd}$. Similarly, the inclusions
\begin{align*}
    g(t)t^{m'} N\fps \hookrightarrow g(t)t^{m} N \qquad t^{-l'}G[t^{-1}]\hookrightarrow t^{-l}G[t^{-1}]
\end{align*}
are, respectively, $G_1[t^{-1}]$ and $G\fps$ equivariant, and induce equivariant inclusions 
\begin{align*}
    \Ann(g(t)t^{m} N\fps) &\hookrightarrow \Ann(g(t)t^{m'} N\fps)
    \\
    \Ann(t^{-l} N[t^{-1}]) &\hookrightarrow \Ann(t^{-l'} N[t^{-1}]).
\end{align*}
Therefore, we obtain equivariant inclusions 
\begin{multline*}
    N^*\lfps_{m,l}=\Ann(g(t)t^{m} N\fps)\cap \Ann(t^{-l} N[t^{-1}])
    \\
    \hookrightarrow \Ann(g(t)t^{m'} N\fps)\cap \Ann(t^{-l'} N[t^{-1}])=N^*_{m', l'}\lfps
\end{multline*}
inducing a fiberwise injective morphism of affine bundles $i: \Sh{M}^{\xi, m',m,l', l}_{G, N, \dd}\hookrightarrow \Sh{M}^{\xi, m',l'}_{G, N, \dd}$. Moreover, it is easy to check that the critical locus $\crit(\Sh{W}_{m', l'})\subset \Sh{M}^{\xi,m',l'}_{G, N, \dd}$ factors through $i$ and that 
\begin{equation}
    \label{eq: compat potential comparison different approximations}
    \Sh{W}_{m,l}\circ j=\Sh{W}_{m', l'}\circ i.
\end{equation}
Overall, we have constructed canonical maps fitting in a commutative diagram 
\begin{equation}
    \label{eq: diagram DT Čech indep of approx}
\begin{tikzcd}
    \Sh{M}^{\xi,m, l}_{G, N, \dd}&\Sh{M}^{\xi, m',m,l', l}_{G, N, \dd}\arrow[l, twoheadrightarrow, swap, "j"]\arrow[r, hookrightarrow, "i"]& \Sh{M}^{\xi, m',l'}_{G, N, \dd}
    \\
     \QM^{\xi}(X)\arrow[u, hookrightarrow]& 
     j^{-1}(\QM^{\xi}(X))\arrow[l, twoheadrightarrow]\arrow[u, hookrightarrow, "s"]
     \\
     &\QM^{\xi}(X)\arrow[r, equal]\arrow[ul, equal]\arrow[u, hookrightarrow, "t"]& \QM^{\xi}(X)\arrow[uu, hookrightarrow, "u"]
\end{tikzcd}
\end{equation}
where both squares are Cartesian.
Here, $t$ is the zero section of the affine fibration $j^{-1}(\QM^{\xi}(X))\to \QM^{\xi}(X)$ and the leftmost and rightmost vertical maps are the natural inclusion of the critical loci for the functions $\Sh{W}_{m,l}$ and $\Sh{W}_{m', l'}$, respectively.
By non-degeneracy of the residue pairing, we have 
\begin{equation}
    \label{eq: invariance VC Čech appprox 0}
    \reldim(j)=\codim_{ \Sh{M}^{\xi, m',l'}_{G, N, \dd}}(\Sh{M}^{\xi, m',m,l', l}_{G, N, \dd}).
\end{equation}
Set $k\coloneqq\reldim(j)$. Notice that since $j$ is an affine fibration, we have a natural isomorphism
\begin{equation}
\label{eq: invariance VC Čech appprox 1}
    \phip{\Sh{W}_{m,l}\circ j}\BoQ_{\Sh{M}^{\xi, m',m,l', l}_{G, N, \dd}}\cong j^*\phip{\Sh{W}_{m,l}} \BoQ_{\Sh{M}^{\xi,m, l}_{G, N, \dd}}.
\end{equation}
Moreover, applying the vanishing cycle to the co-unit $i_*i^!\to \id$ we obtain a natural morphism of complexes
\begin{align*}
\begin{split}
     i_*\phip{\Sh{W}_{m,l}\circ j }\BoQ_{\Sh{M}^{\xi,m',m,l',l}_{G, N, \dd}}
    &=
    i_*\phip{\Sh{W}_{m', l'}\circ i }\BoQ_{\Sh{M}^{\xi,m',m,l',l}_{G, N, \dd}}
    \\
    &\cong i_*\phip{\Sh{W}_{m', l'}}i^!\BoQ_{\Sh{M}^{\xi,m',l'}_{G, N, \dd}}[-2k]
    \\
    &\cong \phip{\Sh{W}_{m', l'}}i_*i^!\BoQ_{\Sh{M}^{\xi,m',l'}_{G, N, \dd}}[-2k]
    \\
    &\to i_*\phip{\Sh{W}_{m', l'}}\BoQ_{\Sh{M}^{\xi,m',l'}_{G, N, \dd}}[-2k]
\end{split}
\end{align*}
where in first step follows from \eqref{eq: compat potential comparison different approximations}, the second one from smoothness of $\Sh{M}^{\xi,m',m,l',l}_{G, N, \dd}$ and $\Sh{M}^{\xi,m',l'}_{G, N, \dd}$ and the third one by properness of $i$. Applying $u^*$ we obtain, via proper base change, a canonical morphism 
\[
t^*s^*\phip{\Sh{W}_{m,l}\circ j }\BoQ_{\Sh{M}^{\xi,m',m,l',l}_{G, N, \dd}}\to u^*\phip{\Sh{W}_{m', l'}}\BoQ_{\Sh{M}^{\xi,m',l'}_{G, N, \dd}}[-2k]
\]
Combining this morphism with \eqref{eq: invariance VC Čech appprox 1}, taking into account the virtual shifts, and using \eqref{eq: invariance VC Čech appprox 0}, we obtain a canonical morphism 
\begin{equation}
\label{eq: invariance VC Čech appprox 2}
\left(\phip{\Sh{W}^{m, l}}\BoQ^{\vir}_{\Sh{M}^{\xi, m, l}_{G, N, \dd}}\right)|_{\QM_{\dd}^{\xi}(X)}\to  \left(\phip{\Sh{W}^{m', l'}}\BoQ^{\vir}_{\Sh{M}^{\xi, m', l'}_{G, N, \dd}}\right)|_{\QM_{\dd}^{\xi}(X)}.
\end{equation}
It remains to show that this map is an isomorphism. It suffices to work locally over $\Bun^{\leq 0}_{G, \infty}(\bbP^1)$ (in the  lisse-e\'tale topology). All the spaces in the top row of the diagram \eqref{eq: diagram DT Čech indep of approx} are total spaces of vector bundles over $\Bun^{\leq 0}_{G, \infty}(\bbP^1)$. On any small enough smooth cover $q: U\to\Bun^{\leq 0}_{G, \infty}(\bbP^1)$ of a point, these bundles trivialize and the top row of the diagram can be written as
\[
\begin{tikzcd}
    \Sh{M}^{\xi,m, l}_{G, N, \dd}|_U&\Sh{M}^{\xi,m, l}_{G, N, \dd}|_U \oplus \AA^k_U\arrow[l, twoheadrightarrow, swap, "j"]\arrow[r, hookrightarrow, "i"]& \Sh{M}^{\xi,m, l}_{G, N, \dd}|_U \oplus \AA^{2k}_U
\end{tikzcd}
\]
where $\AA^{2k}_U$ is a trivial vector bundle of rank $2k$ on $U$ and $\AA^{k}$ is a rank $k$ sub-bundle. Moreover, on this restriction the potential $\Sh{W}_{m', l'}$ takes the form $\Sh{W}_{m', l'}|_U=\Sh{W}_{m,l}|_U+Q$, where $\Sh{W}_{m,l}|_U$ is pulled back from the first factor and $Q$ is a nondegenerate quadratic form on $\AA^{2k}_U$ that restricts to zero on $\AA^{k}_U$. The non-degeneracy of $Q$ implies that the critical locus of $Q$ is $0\in \AA^{2k}_U$, the zero bundle. Hence, it follows from Proposition \ref{prop: def dim red smooth case} that the restriction of \eqref{eq: invariance VC Čech appprox 2} is an isomorphism. The theorem follows.

\end{proof}

Since we have shown that $ \phip{\Sh{W}^{m, l}}\BoQ^{\vir}_{\Sh{M}^{\xi, m, l}_{G, N, \dd}}$, seen as a sheaf on $\Crit(\Sh{W}^{m, l})=\QM^{\xi}_{\dd}(X)$ does \emph{not} depend on the choices of $m, l\ggg 0$, we will denote it by $\phip{\Sh{W}}\BoQ^{\vir}_{\Sh{M}^{\xi}_{G, N, \dd}}$. However, we stress that the construction of this implicitly involves a choice of approximation $\Sh{M}^{\xi, m, l}_{G, N, \dd}$. We do \emph{not} define $\phip{\Sh{W}}\BoQ^{\vir}_{\Sh{M}^{\xi}_{G, N, \dd}}$ by applying the formal definition of vanishing cycle functor to the constant sheaf\footnote{whose perverse exact shift is not even defined by infinite dimensionality.} on the \emph{infinite} dimensional space $\Sh{M}^{\xi}_{G, N}$ as this would create pathological behaviors. Likewise, the derived global sections of $\phip{\Sh{W}}\BoQ^{\vir}_{\Sh{M}^{\xi}_{G, N, \dd}}$ will be denoted by 
\[
\HO_T(\Sh{M}^{\xi}_{G, N, \dd}, \phip{\Sh{W}}\BoQ^{\vir}_{\Sh{M}^{\xi}_{G, N, \dd}})
\]
although we really understand this expression by taking $l,m\ggg 0$ and setting
\[
\HO_T(\Sh{M}^{\xi}_{G, N, \dd}, \phip{\Sh{W}}\BoQ^{\vir}_{\Sh{M}^{\xi}_{G, N, \dd}})\coloneqq\HO_T(\Sh{M}^{\xi, m, l}_{G, N, \dd}, \phip{\Sh{W}^{m,l}}\BoQ^{\vir}_{\Sh{M}^{\xi, m, l}_{G, N, \dd}})
\]
Since it is not illuminating to constantly work with finite dimensional approximations, we will usually drop them from the notation. However, all the statements and constructions of this work involve them. 

We also remark that, by Corollary \ref{cor: DT are the same}, we have a canonical isomorphism 
\[
\HO_T(\Sh{M}^{\xi}_{G, N, \dd}, \phip{\Sh{W}}\BoQ^{\vir}_{\Sh{M}^{\xi}_{G, N, \dd}})\cong \HO_T(\QM^{\xi}_{\dd}(X), \varphi_{\QM^{\xi}_{\dd}(X)}).
\]
Therefore,the pair $(\Sh{M}^{\xi}_{G, N},\Sh{W}^{\xi})$ provides a global critical locus model for the DT cohomology of $\QM^{\xi}_{\dd}(X)$.

\subsection{Dimensional reduction in the Čech model}
\label{subsec: dim red Čech model}
In \S\ref{subsec: ev point and dim red} it is proved that if the base-point $p=(\xi_2,\xi_3, \imath,\jmath)$ is polarized, i.e. it can we written as $p=(\xi_2,0, 0,\jmath)$, then we have a canonical isomorphism 
\begin{equation}
    \label{eq: dimred quasimap Čech model formula}
    \HO_T(\QM^{\xi}_{\dd}(X), \varphi_{\QM^{\xi}_{\dd}(X)})\cong \HO_T^{\BM}(\QM^{\xi}(L),\BoQ^{\vir}).
\end{equation}
It is instructive to recover this isomorphism using the Čech model above. Since $\xi_3=0$ and $\imath=0$, the potential \eqref{eq: potential Čech compact version} on $\Mt^{\xi}_{G,N}$ reads 
\[
\Sh{W}=\Res\left(n^*_{\infty}\left(gn_0-n_{\infty}\right)\right)
\]
In particular, it is homogeneous in $n^*_{\infty}$ and hence can be dimensionally reduced. Formally, set 
\begin{align*}
    \Sh{U}^{\xi}_{G, N} 
    \coloneqq&N\fps\times_{G\fps} G\lfps \times_{ G[t^{-1}] }N^{\xi}[t^{-1}])
    \\
    =&\Tot_{\Bun_{G, \infty}(\bbP^1)}\left(N\fps\times N[t^{-1}]\right)
\end{align*}
A point in the prequotient of $\Sh{U}^{\xi}_{G, N} $ is denoted by $(n_0(t), g(t), n_{\infty}(t))$. Consider also the zero zero locus 
\[
\Nt^{\xi}_{G,N} \coloneq Z(gn_0-n_{\infty})\subset \Ut^{\xi}_{G,N}
\]
Recalling the notation from\S\ref{subsec: infinite Čech model} we further introduce finite approximations of these space: we set 
\[
\Ut^{\xi,l,m}_{G,N}=\Tot_{\Bun_{G, \infty}(\bbP^1)}\left(N^m\fps\times N^l[t^{-1}]\right)
\]
and let $\Nt^{\xi,l,m}_{G,N}\subset \Ut^{\xi,l,m}_{G,N}$ be the zero locus of $gn_0-n_{\infty}$.

The proof of the follows Proposition parallels the one of Theorem \ref{thm: Čech critical quasimap} and Proposition \ref{prop: approx crit Čech}.
\begin{proposition}
\label{prop: dimensionaly reduced quasimap Čech model}
We have a natural equivalence $\Nt^{\xi}_{G,N}\simeq \QM^{\xi}(L)$.
In particular, the left hand side is a scheme of finite type. Moreover, we also have natural equivalences $\Nt^{\xi,l,m}_{G,N}\simeq \QM^{\xi}(L)$ for any $l, m\ggg 0$.
\end{proposition}
We are now in the position to recover \eqref{eq: dimred quasimap Čech model formula} using the Čech model. To argue formally, fix $l,m\ggg 0$. Let $E^{m,l}\to \QM^{\xi}(L)$ be the restriction of the (finite rank) vector bundle 
\begin{align*}
    \Mt^{\xi, m,l}_{G,N}
    &=N^m\fps\times_{G\fps} G\lfps \times_{ G[t^{-1}] }\left( N^*\lfps \times N^{\xi}_{l}[t^{-1}])\right)
    \\
    &\to N^m\fps\times_{G\fps} G\lfps \times_{G[t^{-1}]}N^{\xi}_{l}[t^{-1}])=\Ut^{\xi, m,l}_{G,N}
\end{align*}
to $\QM^{\xi}(L)\simeq \Nt^{\xi,m,l}_{G,N}\subseteq \Ut^{\xi,m,l}_{G,N}$. By the dimensional reduction theorem \ref{thm: dimred}, we obtain a canonical isomorphism
\[
\HO_T(\Sh{M}^{\xi, m, l}_{G, N, \dd}, \phip{\Sh{W}^{m,l}}\BoQ^{\vir}_{\Sh{M}^{\xi, m, l}_{G, N, \dd}})=\HO_T^{\BM}(\QM^{\xi}(L), \BoQ[\dim(\Sh{M}^{\xi, m, l}_{G, N, \dd}, \phip{\Sh{W}^{m,l}})-2\rk(E_{m,l})]),
\]
where now $E_{m,l}$ is a vector bundle over $\QM^{\xi}(L)$ with fiber $N^*\lfps_{m,l}$, see \S\ref{subsec: Finite type approximations of the Čech model} for the notation. Therefore it remains to check that the resulting cohomological shift equals the one from Proposition \ref{prop: dimred Beilison model}, and in particular that it is independent of $m,l$. We sketch the argument here. Using \eqref{eq: approximated Čech model}, the fact that $\dim(\Bun_{G, \infty}(\bbP^1))=0$, and equation \eqref{eq: computation dimension fibers Čech model}, one sees that 
\[
\dim(\Sh{M}^{\xi, m, l}_{G, N, \dd}, \phip{\Sh{W}^{m,l}})-2\rk(E_{m,l})=-\langle \dd, \chi_{\det(N)}\rangle,
\]
where we are identifying $\dd$ with an element in $X_*(S)/\{\text{coroots}\}\cong \pi_1(G)$. But we have 
\[
\langle \dd, \chi_{\Hom(V_i, W_i)}\rangle=\sum_{i}-\dd_i\ww_i \qquad \langle \dd, \chi_{\Hom(V_i, V_j)}\rangle=\sum_{i\to j}\dd_j\vv_i-\dd_i\vv_j,
\]
whence 
\begin{equation}
    \label{eq: computation shift dimred quasimap Čech}
    -\langle \dd, \chi_{\det(N)}\rangle=\dd^T(\Q^T -\Q)\vv+\dd^T\ww,
\end{equation}
in agreement with Proposition \ref{prop: dimred Beilison model}.



\subsection{Coulomb branch action on quasimaps}
\label{subsec: Coulomb branch action on quasimaps}

Fix a quiver $Q=(Q_0,Q_1)$ and set $N=\Rep_Q(\ww,\vv)$ and $G=\GL_{\dd}$. In this section, we construct the Coulomb branch action on
\begin{equation}
    \label{eq: DT cohomology quasimaps and critical rpesentation}
    \HO_T(\QM^{\xi}_{\dd}(X), \varphi_{\QM^{\xi}_{\dd}(X)})=\HO_T(\Sh{M}^{\xi}_{G, N, \dd}, \phip{\Sh{W}}\BoQ^{\vir}_{\Sh{M}^{\xi}_{G, N, \dd}})
\end{equation}
To achieve this, we combine the Čech model for the the quasimap DT sheaf \S\ref{subsec: DT sheaf Čech model} with the critical description of the Coulomb branch developed in \S\ref{subsec: critical Coulomb branch} and \S\ref{subsec: critical Coulomb branch convolution algebra}. Recall that the Čech presentation of $\QM^{\xi}(X)$ is given in terms of the space 
\begin{align}
\label{eq: ambient Čech action section}
    \Sh{M}^{\xi}_{G, N} 
    \coloneqq&N\fps\times_{G\fps} \times G\lfps \times_{G[t^{-1}]} \left(N^*\lfps \times N^{\xi}[t^{-1}]\right)
\end{align}
with potential
\begin{equation}
\label{eq: Čech potential action section}
    \Sh{W}^{\xi}=\Res\left(n^*_{\infty}\left(gn_2-n_{\infty}\right)\right)+\Res\Tr((\xi_2, \iota) n_{\infty})
\end{equation}
Here, $(n_2, g, n^*_{\infty},n_{\infty})$ denotes a point in $ G\lfps \times N\fps \times N^*\lfps \times N^{\xi}[t^{-1}]$. The reason why we label the field $n_2$ with such subscript will be apparent in the course of this section. 
Recall also that the critical presentation of the Coulomb branch is given by the pair space with potential $(\tilde \Rt_{G,N}, \Sh{W})$, where
\begin{equation}
\label{eq: ambient Coulomb action section}
    \tilde \Rt_{G,N}= N\fps\times_{G\fps} G\lfps\times (N^*\lfps\times N\fps)
\end{equation}
and
\[
\Sh{W}: \tilde\Rt_{G,N}\to \BoC \qquad \Res\Tr\left( n_2^*\left(gn_1-n_2\right)\right).
\]
where $(n_1, g, n_2^*, n_2)$ is a tuple in $ N\fps\times G\lfps\times (N^*\lfps\times N\fps)$.

\begin{remark}
    All the sheaf theory operations in this section will be understood in terms of the finite dimensional approximations of these spaces. The approximation of the pair $(\Sh{M}^{\xi}_{G, N} , \Sh{W}^{\xi})$ is presented in \S\ref{subsec: Finite type approximations of the Čech model}, while the one of the pair $(\Rt_{G,N}, \Sh{W})$ is given in \S\ref{subsec: Finite type approximations critical Coulomb}. We remark here that there is  a fundamental difference between these two approximations. The approximation $\Sh{M}^{\xi, m, l}_{G, N}$ of $\Sh{M}^{\xi}_{G, N}$ is obtained by truncating the fibers of the vector bundle 
    \[
    \Sh{M}^{\xi}_{G, N}\to G\fps\backslash G\lfps^{\leq 0} /G_1[t^{-1}]=\Bun^{\leq 0}_{G, \infty}(\bbP^1)
    \]
    which is a infinite rank vector bundle over the quasi-compact space $\Bun^{\leq 0}_{G, \infty}(\bbP^1)$, which is locally of finite type. In particular, the base $\Bun_{G, \infty}(\bbP^1)$ is \emph{not} approximated.

    On the other hand, the approximation $\tilde \Rt^{d}_{G,N}$ of $\tilde \Rt_{G,N}$ is obtained by first restricting the bundle $\tilde \Rt_{G,N}\to \Gr_G$ to a Schubert variety $\Gr_{G,\leq \lambda}$, then truncating the fibers, and finally noticing the $G\fps$ action on this space factors through a quotient $G_d=G(\Sh{O}/t^d \Sh{O})$. This is unavoidable as $\Gr_G$ is not a finite type object. Despite defining the correspondence inducing the Coulomb branch action formally using infinite dimensional geometries, the proof of Proposition \ref{prop: Coulomb action def}, where the actual map is defined, uses finite type approximations.
\end{remark}

We now introduce the moduli spaces 
\begin{align}
\label{eq: three spaces critical Coulomb action}
\begin{split}
    \Mt^{(2),0, \xi}_{G,N}&=N\fps\times_{G\fps} G\lfps \times \left(N^*\lfps\times N\fps\times  G\lfps^{\leq 0}\right)\times_{G_1[t^{-1}]} \left(N^*\lfps \times N^{\xi}[t^{-1}]\right)
    \\
    \Mt^{(2),1, \xi}_{G,N}&=N\fps\times_{G\fps} G\lfps \times_{G\fps} \left(N^*\lfps\times N\fps\times  G\lfps^{\leq 0}\right) \times_{G_1[t^{-1}]} \left(N^*\lfps \times N^{\xi}[t^{-1}]\right)
    \\
    \Mt^{(2),2,\xi}_{G,N}&=\left(N^*\lfps\times N\lfps\times N\fps\right) \times_{G\fps}G\lfps \times_{G\fps} G\lfps^{\leq 0} \times_{G_1[t^{-1}]} \left(N^*\lfps \times N^{\xi}[t^{-1}]\right)
    \\
    \Mt^{(2),3,\xi}_{G,N}&=\left(N^*\lfps\times N\fps\right) \times_{G\fps}G\lfps \times_{G\fps} G\lfps^{\leq 0} \times_{G_1[t^{-1}]} \left(N^*\lfps \times N^{\xi}[t^{-1}]\right)
    \\
    \Mt^{(2),4,\xi}_{G,N}&= N\fps \times_{G\fps}G\lfps \times_{G\fps} G\lfps^{\leq 0} \times_{G_1[t^{-1}]} \left(N^*\lfps \times N^{\xi}[t^{-1}]\right).
\end{split}
\end{align}
Closed points in each of thse spaces are denoted by 
\[
\begin{split}
    & \left( [n_1, g], n_2^*, n_2, [g', n_\infty^*, n_\infty]\right)
    \\
    &  [n_1, g, n_2^*, n_2, g', n_\infty^*, n_\infty]
    \\
    & [w^*, w, n_1, g, g', n_\infty^*, n_\infty]
    \\
    & [ w^*,n_1, g, g', n_\infty^*, n_\infty]
    \\
    & [n_1, g, g', n_\infty^*, n_\infty].
\end{split}
\]
respectively. Notice that we ordered the entries of the cosets in a way that naturally matches the terms in the fiber product. We assemble this spaces into the following diagram 
\begin{equation}
    \label{eq: correspondence Coulomb branch action}
    \begin{tikzcd}
        \tilde\Rt_{G,N}\times \Mt^{\xi}_{G,N} & \Mt^{(2),0, \xi}_{G,N}\arrow[l, swap, "\tilde p^{\xi}"]\arrow [r, "\tilde q_1^{\xi}"] & \Mt^{(2),1, \xi}_{G,N}\arrow[r, "\tilde j^{\xi}"]& \Mt^{(2),2, \xi}_{G,N}
        \\
        & & & \Mt^{(2),3, \xi}_{G,N}\arrow[u, swap, "\tilde q_2^{\xi}"]\arrow[d, "\tilde q_3^{\xi}"]
        \\
        & & &\Mt^{(2),4, \xi}_{G,N} \arrow[r, "\tilde m^{\xi}"] & \Mt^{\xi}_{G,N} 
\end{tikzcd}
\end{equation}
Formally, all the maps are defined as in \S\ref{subsec: critical Coulomb branch convolution algebra}. Explicitly, The maps $\tilde p^{\xi}$ and $\tilde q^{\xi}_1$ are defined by
\begin{multline}
    \label{eq: first two maps critical ravolo action}
     (([n_1, g], n_2^*, n_2), [n_2, g', n_\infty^*, n_\infty])\mapsfrom \left( [n_1, g], n_2^*, n_2, [g', n_\infty^*, n_\infty]\right)
     [n_1, g, n_2^*, n_2, g', n_\infty^*, n_\infty]
\end{multline}
The map $\tilde j$ is defined by
\begin{equation}
\label{eq: j map Coulomb action}
    [n_1, g, n_2^*, n_2, g', n_\infty^*, n_\infty]\mapsto [g'gn_\infty^*-gn_2^*,gn_2-n_1,n_1, g, g', n_\infty^*, n_\infty].
\end{equation}
The maps $\tilde q^{\xi}_2$ and $\tilde q^{\xi}_3$ are defined by 
\begin{equation}
    \label{eq: second two pair of maps critical ravolo action}
    [w^*, 0,n_1, g, g', n_\infty^*, n_\infty] \mapsfrom [w^*, n_1, g, g', n_\infty^*, n_\infty]\mapsto [n_1, g, g', n_\infty^*, n_\infty]
\end{equation}
Finally, the map $\tilde m$ is defined by 
\begin{equation}
    \label{eq: last map critical ravolo action}
    [n_1, g, g', n_\infty^*, n_\infty]\mapsto [n_1,  g'g, n_\infty^*, n_\infty]
\end{equation}

\begin{remark}
    By definition, all the spaces in \eqref{eq: three spaces critical Coulomb action} are of the form $X\times_{G_1[t^{-1}]}(N^*\lfps\times N^{\xi}[t^{-1}])$, where $X$ is a space carrying a $G_1[t^{-1}]$-action. Replacing the factor $(-)\times_{G_1[t^{-1}]}(N^*\lfps\times N^{\xi}[t^{-1}])$ with $(-)\times (N^*\lfps\times N\fps)$ while keeping $X$ fixed\footnote{Technically, one also replaces the factor $G\lfps$ with $G\lfps^{\leq 0}$.} one gets the four spaces \eqref{eq: three spaces critical Coulomb mult} appearing in the critical Coulomb branch correspondence \eqref{eq: critical Coulomb branch correspondence}.  In fact, the correspondence itself can be thought as obtained from \eqref{eq: critical Coulomb branch correspondence} by replacing the term $(-)\times_{G_1[t^{-1}]}(N^*\lfps\times N^{\xi}[t^{-1}])$ with $(-)\times (N^*\lfps\times N\fps)$ in all the spaces involved while keeping the formal definition of the morphisms the same.
\end{remark}

We now introduce potentials on these spaces as follows:
\begin{itemize}
    \item The potential on $\tilde\Rt_{G,N}\times \Mt^{\xi}_{G,N}$ is the direct sum potential $\Sh{W}\boxplus \Sh{W}^{\xi}$.
    \item The potential on $\Mt^{(2),0, \xi}_{G,N}$ is 
    \begin{multline}
    \label{eq: potential W_0 xi}
        \Sh{W}_0^{\xi}\coloneqq(\Sh{W}\boxplus \Sh{W}^{\xi})\circ \tilde p^{\xi}
        \\
        =\Res\Tr\left(n_2^*(gn_1-n_2)\right) +\Res\Tr\left(n_\infty^*(g'n_2-n_\infty)\right)+\Res\Tr((\xi_2, \iota) n_{\infty}),
    \end{multline}
    the pullback of $\Sh{W}\boxplus \Sh{W}^{\xi}$ along $\tilde p^{\xi}$. Since it is $G\fps$ equivariant, it descends to a potential $\Sh{W}_1$ on $\Mt^{(2),1}_{G,N}$.
    \item The potential $\Sh{W}_2^{\xi}$ on $\Mt^{(2),2, \xi}_{G,N}$ is given by 
    \begin{equation}
     \label{eq: potential W_1 xi}
    \Sh{W}_2^{\xi}\coloneqq\Res\Tr(w^*w)+\Res\Tr\left( n_\infty^*\left(g'g n_1-n_\infty\right)\right)+\Res\Tr((\xi_2, \iota) n_{\infty}).
    \end{equation}
    \item The potential on $\Mt^{(2),3, \xi}_{G,N}$, denoted by $\Sh{W}_3$, is the composition 
    \[
    \Sh{W}_3^{\xi}\coloneqq\Sh{W}_2\circ \tilde q_2^{\xi}=\Res\Tr\left( n_\infty^*\left(g'gn_1-n_\infty\right)\right)+\Res\Tr((\xi_2,\iota) n_{\infty}).
    \]
    \item The potential on $\Mt^{(2),4, \xi}_{G,N}$ is 
    \[
    \Sh{W}_4^{\xi}\coloneqq\Sh{W}^{\xi}\circ \tilde m^{\xi}=\Res\Tr\left( n_\infty^*\left(g'gn_1-n_\infty\right)\right)+\Res\Tr((\xi_2,\iota) n_{\infty})
    \] 
\end{itemize}
The same argument of Lemma \ref{lemma: critical raviolo 2} shows that $\Sh{W}_2^{\xi}\circ \tilde j^{\xi}=\Sh{W}_1^{\xi}$, hence all the morphisms in \eqref{eq: three spaces critical Coulomb action} are compatible with the potentials assigned to the vertices of the diagram. We now proceed with constructing the correspondence defining the Coulomb action on \eqref{eq: DT cohomology quasimaps and critical rpesentation}. Recall the definition of the index $I_{\dd}$ from \eqref{eq: index shifted action}. 

\begin{proposition}
\label{prop: Coulomb action def}
    The correspondence induces a canonical shifted action of the Coulomb branch on $\HO_T(\Mt_{G,N}, \phip{\Sh{W}^{\xi}}\BoQ^{\vir})$. Specifically, the restriction of this action to a summand ${}^{\phi}\Coul^{T,\vir}_{G,N, \dd}$ is a morphism 
    \[
    {}^{\phi}\eta^{\xi}:{}^{\phi}\Coul^{T,\vir}_{G,N, \dd} \otimes \HO_T(\Mt_{G,N}, \phip{\Sh{W}^{\xi}}\BoQ^{\vir})\to \HO_T(\Mt_{G,N}, \phip{\Sh{W}^{\xi}}\BoQ^{\vir})[I_{\dd}].
    \]
\end{proposition}
\begin{proof}
Informally, the map ${}^{\phi}\eta^{\xi}$ is given by the composition 
\begin{equation}
\label{eq: legit action}
    \tilde m^{\xi}_* \circ (\tilde q^{\xi}_{2,*}\circ \tilde q^{\xi, *}_3)^{-1}\circ \tilde j^{\xi}_*\circ (\tilde q^{\xi,*}_1)^{-1} \circ \tilde p^{\xi,*}
\end{equation}
where all the morphisms are formally defined as those in \eqref{eq: critical Coulomb branch multiplication}. However, here we choose to discuss the details of all the finite dimensional approximations. To remove clutter, we drop the subscripts $G,N$ from all the spaces below. Fix a dominant cocharacter $\lambda$ and let $\dd$ be its conjugacy class in $\pi_1(G)=\BoZ^{Q_0}$. Recall the finite type approximations $\tilde\Rt^{m,d,r}_{\leq\lambda}$ and $\Mt^{\xi, d,l}_{\dd'}$ from \S\ref{subsec: Finite type approximations critical Coulomb} and \S\ref{subsec: Finite type approximations of the Čech model}, respectively. Following the assumptions from those sections, we assume that $m,d,r,l\ggg 0$ and $m\ggg d$ (relative to $\lambda$).

Consider the approximation of $\Mt^{(2),0,\xi}$ given by\footnote{Technically, the definition below hides a slight abuse of notation since $N^*\lfps_{d,l}$ was defined in \S\ref{subsec: Finite type approximations of the Čech model} as a vector bundle over $\Bun^{\leq 0}_{G, \infty}$. Here, we are identifying it with its fiber. Since it acts passively in most of the construction, this shall not cause confusion.}
\begin{multline*}
    \Mt^{(2),0,\xi,m,d,r,l}_{\leq \lambda, \dd'}=N^m\fps\times_{G\fps} G\lfps_{\leq \lambda} 
    \\
    \times \left(N^{*,r}_d\lfps\times N^d\fps \times G\lfps^{\leq 0}_{\dd'}\right) \times_{G_1[t^{-1}]} \left(N^*\lfps_{d,l} \times N^{\xi}_{l}[t^{-1}]\right)
\end{multline*}
Then assignment \eqref{eq: first two maps critical ravolo action} gives a well defined map $\tilde p^{\xi}: \Mt^{(2),0,m,d,r,l}_{\leq \lambda, \dd'}\to \tilde\Rt^{m,d,r}_{G,N, \leq\lambda}\times \Mt^{\xi, d,l}_{G,N,\dd'}$. Similarly, we define an approximation by $\Mt^{(2),1,\xi}$  by
\begin{multline*}
    \Mt^{(2),1,\xi,m,d,r,l}_{\leq \lambda, \dd'}=N^m\fps\times_{G\fps} G\lfps_{\leq \lambda} 
    \\
    \times_{G\fps} \left(N^{*,r}_d\lfps\times N^d\fps \times G\lfps^{\leq 0}_{\dd'}\right) \times_{G_1[t^{-1}]} \left(N^*\lfps_{d,l} \times N^{\xi}_{l}[t^{-1}]\right).
\end{multline*}
We clearly have a well defined map $\tilde q_1: \Mt^{(2),0,\xi,m,d,r,l}_{\leq \lambda, \dd'}\to \Mt^{(2),1,\xi,m,d,r,l}_{\leq \lambda, \dd'}$.
To define a finite approximation of $\Mt^{(2),2,\xi}$, we proceed as follows. Consider the morphism of vector bundles
\begin{multline*}
    f: N^m\fps\times_{G\fps} G\lfps_{\leq \lambda} \times_{G\fps}\left(N\lfps\times N^*\lfps\times_{G\fps}  G\lfps^{\leq 0}_{\dd'} \right) \times_{G_1[t^{-1}]} \left(N^*_{d, l}\lfps \times N^{\xi}_{l}[t^{-1}]\right)
    \\
    \to
    \left(N^*\lfps\times N\lfps\times N^m\fps)\right) \times_{G\fps}G\lfps_{\leq \lambda} \times_{G\fps} G\lfps^{\leq 0}_{\dd'} \times_{G_1[t^{-1}]} \left(N^*_{d,l}\lfps \times N^{\xi}[t^{-1}]\right)
\end{multline*}
 over the base
 \begin{equation}
 \label{eq: base morphism some long proof}
     N^m\fps\times_{G\fps} G\lfps_{\leq \lambda} \times_{G\fps} G\lfps^{\leq 0}_{\dd'} \times_{G_1[t^{-1}]} \left(N^*_{d, l}\lfps \times N^{\xi}_{l}[t^{-1}]\right)
 \end{equation}
 given by
\[
[n_1, g, n_2^*, n_2, g', n_\infty^*, n_\infty]\mapsto [gn_2^*, gn_2,n_1, g, g', n_\infty^*, n_\infty].
\]
The domain of this map contains a finite chain of subbundles $E^{d,r}\subset F^{d,r}$ defined as
\begin{align*}
     E^{d,r}&\coloneqq N^m\fps\times_{G\fps} G\lfps_{\leq \lambda} \times \left(t^rN\fps\times t^dN^*\fps\times  G\lfps^{\leq 0}_{\dd'} \right) \times_{G_1[t^{-1}]} \left(N^*_{d, l}\lfps \times N^{\xi}_{l}[t^{-1}]\right)
    \\
    F^{d,r}&\coloneqq N^m\fps\times_{G\fps} G\lfps_{\leq \lambda} \times \left(N_d\lfps\times N^*_r\lfps\times  G\lfps^{\leq 0}_{\dd'} \right) \times_{G_1[t^{-1}]} \left(N^*_{d, l}\lfps \times N^{\xi}_{l}[t^{-1}]\right).
\end{align*}
We set $\Mt^{(2),1,\xi,m,d,r,l}_{\leq \lambda, \dd'}=f(F^{d,r})/f(E^{d,r})$. Since $\rk(F^{d,r}/E^{d,r})$ has rank $2(r+k)$ and \eqref{eq: base morphism some long proof} is of finite type, it follows that $\Mt^{(2),1,\xi,m,d,r,l}_{\leq \lambda, \dd'}$ is also of finite type and of the same rank. It is now easy to check that the assignment \eqref{eq: j map Coulomb action} induces a morphism $\tilde j^{\xi}: \Mt^{(2),1,\xi,m,d,r,l}_{\leq \lambda, \dd'}\to \Mt^{(2),2,\xi,m,d,r,l}_{\leq \lambda, \dd'}$. The construction of the remaining approximations is similar. Specifically, $\Mt^{(2),3,\xi,m,d,r,l}_{\leq \lambda, \dd'}$ is the subbundle of $\Mt^{(2),2,\xi,m,d,r,l}_{\leq \lambda, \dd'}$ realized as the quotient $f(\bar F^{d,r})/f(\bar E^{d,r})$, where 
\begin{align*}
     \bar E^{d,r}&\coloneqq N^m\fps\times_{G\fps} G\lfps_{\leq \lambda} \times \left( t^rN^*\fps\times_{G\fps}  G\lfps^{\leq 0}_{\dd'} \right) \times_{G_1[t^{-1}]} \left(N^*_{d, l}\lfps \times N^{\xi}_{l}[t^{-1}]\right)
    \\
    \bar F^{d,r}&\coloneqq N^m\fps\times_{G\fps} G\lfps_{\leq \lambda} \times \left( N^*_d\lfps\times_{G\fps}  G\lfps^{\leq 0}_{\dd'} \right) \times_{G_1[t^{-1}]} \left(N^*_{d, l}\lfps \times N^{\xi}_{l}[t^{-1}]\right).
\end{align*}
The map $\tilde q^{\xi}_2: \Mt^{(2),3,\xi,m,d,r,l}_{\leq \lambda, \dd'}\hookrightarrow \Mt^{(2),2,\xi,m,d,r,l}_{\leq \lambda, \dd'}$ is the inclusion of the domain as a subbundle of the target. The approximation $\Mt^{(2),4,\xi,m,d,r,l}_{\leq \lambda, \dd'}$ is simply given by \eqref{eq: base morphism some long proof}, and $\tilde q^{\xi}_3: \Mt^{(2),3,\xi,m,d,r,l}_{\leq \lambda, \dd'}\to \Mt^{(2),4,\xi,m,l}_{\leq \lambda, \dd'}$ is a finite rank vector bundle. Therefore, we have a commutative diagram 
\begin{equation}
\label{eq: approximated correspondence}
    \begin{tikzcd}
        \tilde\Rt^{m,d,r}_{\leq\lambda}\times \Mt^{\xi, d,l}_{\dd'} & \Mt^{(2),0,m,d,r,l}\arrow[l, swap, "\tilde p^{\xi}"]\arrow [r, "\tilde q_1^{\xi}"] & \Mt^{(2),1,\xi,m,d,r,l}_{\leq \lambda, \dd'}
        \\
        & &  \Mt^{(2),2,\xi,m,d,r,l}_{\leq \lambda, \dd'}\arrow[u, swap, "\tilde q_2^{\xi}"]\arrow[d, "\tilde q_3^{\xi}"]
        \\
        & & \Mt^{(2),3,m,l,\xi}_{\leq \lambda, \dd'} \arrow[r, "\tilde m^{\xi}"] & \Mt^{\xi, m,l}_{\dd+\dd'} 
\end{tikzcd}
\end{equation}
All the spaces of finite type except from the source of $\tilde p^{\xi}$ and $\tilde q_1^{\xi}$. The latter are both principal $G\fps$-bundles. To further turn them into finite type morphism one may quotient such morphisms by the groups $K_i=\Ker(G\fps\mapsto G^i\fps)$ for $i$ large enough. The precise procedure is explained in \cite[\S2(ii)]{BFNII}, so we do not discuss it here.
Note also that the map $\tilde m^{\xi}$ in the target is given by the composition of the proper map\footnote{Notice that we dropped the superscript $\leq 0$ on the term $G\lfps_{\dd'+\dd''}$ appearing in the target of this map. This is because, after applying a Hecke modification to a stably negative bundle, the resulting bundle may not be stably negative. The support condition on $G\lfps^{\leq0}$ will be recovered after applying the vanishing cycle functor, which is supported on the critical locus and hence on stably negative bundles.}
\begin{multline*}
    N^m\fps\times_{G\fps} G\lfps_{\leq \lambda} \times G\lfps^{\leq 0}_{\dd'} \times_{G_1[t^{-1}]} \left(N^*_{d, l}\lfps \times_{G\fps} N^{\xi}_{l}[t^{-1}]\right)
    \\\to N^m\fps\times_{G\fps} G\lfps_{\dd'+\dd''} \times_{G_1[t^{-1}]} \left(N^*_{d, l}\lfps \times N^{\xi}_{l}[t^{-1}]\right)
\end{multline*}
with the closed immersion
\begin{multline*}
    N^m\fps\times_{G\fps} G\lfps_{\dd'+\dd''} \times_{G_1[t^{-1}]} \left(N^*_{d, l}\lfps \times N^{\xi}_{l}[t^{-1}]\right)
    \\\hookrightarrow N^m\fps\times_{G\fps} G\lfps_{\dd'+\dd''} \times_{G_1[t^{-1}]} \left(N^*_{m, l}\lfps \times N^{\xi}_{l}[t^{-1}]\right),
\end{multline*}
which is well defined because of the assumption $m \ggg d$. 

We leave to the reader to check that all the potentials descend to these finite dimensional approximations and are compatible with all th maps in \eqref{eq: approximated correspondence}. In particular, this implies checking the analog of Lemma \ref{lemma: critical raviolo 1} in the finite approximations. Once this is established, formula \eqref{eq: legit action} gives a well defined morphism 
\[
    {}^{\phi}\eta^{\xi}:{}^{\phi}\Coul^{T,\vir}_{G,N,\dd} \otimes \HO_T(\Mt_{\dd'}, \phip{\Sh{W}^{\xi}}\BoQ^{\vir})\to \HO_T(\Mt_{\dd+\dd'}, \phip{\Sh{W}^{\xi}}\BoQ^{\vir})[s].
\]
for some shift $s$ that we shall determine. It is easy to check that the ``base directions'', namely those corresponding to $\Gr_G$ and $\Bun_{G,\infty}$, do not contribute to any shift\footnote{This is basically a consequence of $\Bun_{G,\infty}$ having dimension zero, which implies that its dimension does not contribute to the virtual shift in $\phip{\Sh{W}^{\xi}}\BoQ^{\vir}$.} so we only compute the shift from the ``fiber directions'', that is the terms depending on $N$ and $N^*$. Let $n=\dim(N)$ and $n^*=\dim(N^*)$\footnote{Of course, $n=n^*$ but we prefer to use a distinct notation to make it easier to understand where each of the shifts in the formulas below come from.}. The domain of the map ${}^{\phi}\eta^{\xi}$ is a vanishing cycle applied to the dualizing complex of the product $\tilde\Rt^{m,d,r}_{\leq\lambda}\times \Mt^{\xi, d,l}_{\dd'}$ , shifted (in the fiber directions, which are the only ones we compute) by 
\[
a\coloneqq\underbrace{-2mn-2dn-2rn}_{\text{from }\tilde\Rt^{m,d,r}_{\leq\lambda}}\underbrace{-2\vv^T\dd+2\vv^T\Q^T\dd}_{\text{vir. shift of }\Coul^{\vir}_{G, N,\dd}}\underbrace{-dn-ln-dn^*-ln^*-\langle \dd', \det(N)\rangle}_{\text{from} \Mt^{\xi, d,l}_{\dd'}}.
\]
The virtual shift in the Coulomb branch is discussed in Remark \ref{rem: virtual crit Coulomb}. The pullback along $\tilde p^{\xi}$ introduces a shift by $2dn$. The maps $\tilde j^{\xi}$ and $(\tilde q^{\xi}_{2,*})^{-1}$ do not introduce shifts. The map $(\tilde q^{\xi, *}_{3})^{-1}$ introduces a shift by $2(r+d)n^*$. Finally, the target is the dualizing complex of the product $\Mt^{\xi, d,l}_{\dd+\dd'}$, shifted by 
\[
b\coloneqq-mn-ln-mn^*-ln^*-\langle \dd'+\dd, \chi_{\det(N)}\rangle.
\]
Therefore, the overall shift $s$ is given by 
\begin{align*}
    s
    &=a+2dn+2(r+d)n^*-b 
    \\
    &=\langle \dd', \chi_{\det(N)}\rangle-2\vv^T\dd+2\vv^T\Q^T\dd
    \\
    &=\dd^T(\Q^T -\Q)\vv+\dd^T\ww-2\vv^T\dd+2\vv^T\Q^T\dd
    \\
    &=\dd^T(\ww-\Cartan\vv)
    \\
    &=I_{\dd}.
\end{align*}
This concludes the proof.
\end{proof}

\begin{theorem}
\label{thm: Coulomb action QM}
    The morphism ${}^{\phi}\eta^{\xi}$ induces an action of the Coulomb branch ${}^{\phi}\Coul^{\vir}_{G,N}$ on the quasimap critical cohomology
    \[
    \HO_T(\Sh{M}^{\xi}_{G, N, \dd}, \phip{\Sh{W}}\BoQ^{\vir}_{\Sh{M}^{\xi}_{G, N, \dd}})\cong \HO_T(\QM^{\xi}_{\dd}(X), \varphi_{\QM^{\xi}_{\dd}(X)}).
    \]
    Equivalently, the unit $[\Rt_{e}\times N^*\lfps]\in \Coul^{T, \vir}_{N,G,0}$ acts as the identity element and, for all pairs $\dd',\dd''\in \pi_1(G)$, the diagram
    \[
    \begin{tikzcd}
        {}^{\phi}\Coul^{T,\vir}_{G,N, \dd'}\otimes {}^{\phi}\Coul^{T,\vir}_{G,N, \dd''}\otimes \HO_{T}(\Mt^{\xi}_{G,N}, \phip{\Sh{W}^{\xi}}\DD\BoQ)\arrow[r,"\id\otimes {}^{\phi}\eta^{\xi}"] \arrow[d, swap, "{}^{\phi}\eta\otimes \id"] & {}^{\phi}\Coul^{T,\vir}_{G,N, \dd'+\dd''}\otimes \HO_T(\Mt^{\xi}_{G,N}, \phip{\Sh{W}^{\xi}}\DD\BoQ)\arrow[d, "{}^{\phi}\eta^{\xi}"]
        \\
        {}^{\phi}\Coul^{T,\vir}_{G,N, \dd'}\otimes \HO_T(\Mt^{\xi}_{G,N}, \phip{\Sh{W}^{\xi}}\DD\BoQ)[I_{\dd''}]\arrow[r, "{}^{\phi}\eta^{\xi}"] &  \HO_T(\Mt^{\xi}_{G,N}, \phip{\Sh{W}^{\xi}}\DD\BoQ)[I_{\dd'+\dd''}].
    \end{tikzcd}
    \]
    is commutative.
\end{theorem}

\begin{proof}
To remove clutter, we write $\tilde\Rt$ and $\Mt^{\xi}$ instead of  $\tilde\Rt_{G,N}$ and $\Mt^{\xi}_{G,N}$. Similarly, we drop the subscripts $\xi$ from all the spaces fitting in the correspondence \eqref{eq: correspondence Coulomb branch action}. We will also not keep track of the cohomological degree. We first prove compatibility with the Coulomb branch multiplication. To this end, we introduce an ``interpolating'' correspondence 
\begin{equation}
    \label{eq: interpolating corespondence}
\begin{tikzcd}
    \tilde\Rt\times \tilde\Rt\times \Mt & \Mt^{(3),0,0} \arrow[l, swap, "\tilde p_{123}"]\arrow[r, "\tilde q_{123, 1}"]& \Mt^{(3),1,1} \arrow[r, "\tilde j_{123}"] & \Mt^{(3),2,2}
    \\
    & & & \Mt^{(3),3,3}\arrow[u, "\tilde q_{3,123}"]\arrow[d, "\tilde q_{123,4}"]
    \\
    & & & \Mt^{(3),4,4}\arrow[r, "\tilde m_{123}"]& \Mt^{\xi}
\end{tikzcd}
\end{equation}
The spaces and the morphisms are defined as follows. Firstly, to remove clutter, set $N\fpslfps=N^*\lfps\times N\fps$ and $N\lfpslfps=N^*\lfps\times N\lfps$. We define
\begingroup
\small
\begin{align}
\label{eq: three spaces critical Coulomb action}
\begin{split}
    \Mt^{(3),0,0}_{G,N}&=N\fps\times_{G\fps} G\lfps\times\left(N\fpslfps \times G\lfps\right) \times \left(N\fpslfps\times  G\lfps^{\leq 0}\right)\times_{G_1[t^{-1}]} \left(N^*\lfps \times N^{\xi}[t^{-1}]\right)
    \\
    \Mt^{(3),1,1}_{G,N}&=N\fps\times_{G\fps} G\lfps\times_{G\fps}\left(N\fpslfps \times G\lfps\right) \times_{G\fps} \left(N\fpslfps\times  G\lfps^{\leq 0}\right)\times_{G_1[t^{-1}]} \left(N^*\lfps \times N^{\xi}[t^{-1}]\right)
    \\
    \Mt^{(3),2,2}_{G,N}&=\left(N\lfpslfps^2 \times N\fps\right) \times_{G\fps}G\lfps \times_{G\fps}G\lfps \times_{G\fps} G\lfps^{\leq 0} \times_{G_1[t^{-1}]} \left(N^*\lfps \times N^{\xi}[t^{-1}]\right)
    \\
    \Mt^{(3),3,3}_{G,N}&=\left(N^*\lfps^2\times N\fps\right) \times_{G\fps}G\lfps \times_{G\fps}G\lfps \times_{G\fps} G\lfps^{\leq 0} \times_{G_1[t^{-1}]} \left(N^*\lfps \times N^{\xi}[t^{-1}]\right)
    \\
    \Mt^{(3),4,4}_{G,N}&= N\fps\times_{G\fps}G\lfps \times_{G\fps}G\lfps \times_{G\fps} G\lfps^{\leq 0} \times_{G_1[t^{-1}]} \left(N^*\lfps \times N^{\xi}[t^{-1}]\right).
\end{split}
\end{align}
\endgroup
As usual, we denote a point in each of these moduli spaces as a conjugacy class of tuples in the prequotient, where each element of the tuple corresponds to a term in the prequotient, read from left to right.  The map $p_{123}$ is then given by 
\[
([n_1, g], n_2^*, n_2)\times ([n_2, g], n_3^*, n_3)\times[n_3, g, n_\infty^*, n_\infty] \mapsfrom ([n_1, g], n_2^*, n_2, g', n_3^*, n_3,[g'', n_\infty^*, n_\infty])
\]
The map $\tilde q_{123,1}$ is just the quotient by $G\fps\times G\fps$. The map $\tilde j_{123}$, which will be crucial for us, is given by
\[
[n_1, g, n_2^*, n_2, g', n_3^*, n_3,g'', n_\infty^*, n_\infty]\mapsto [w^*, w, \bar w^*, \bar w, n_1, g, g',g'', n_\infty^*, n_\infty],
\]
where the assignments
\[
(w^*,w)=(g'gn_3^*-gn_2^*, gn_2-n_1)
\quad 
(\bar w^*, \bar w)=(g''g'gn_\infty^*-g'gn_3^*,g'gn_3-n_1)
\]
are understood.
The map $\tilde q_{123, 2}$ is given by 
\[
[w^*, w, \bar w^*, w, n_1, g, g',g'', n_\infty^*, n_\infty]\mapsto[w^*, \bar w^*, n_1, g, g',g'', n_\infty^*, n_\infty]
\]
and $\tilde q_{123, 3}$ is given by 
\[
[w^*, \bar w^*, n_1, g, g',g'', n_\infty^*, n_\infty]\mapsto [n_1, g, g',g'', n_\infty^*, n_\infty]
\]
Finally, $\tilde m_{123}$ is given by 
\[
[n_1, g, g',g'', n_\infty^*, n_\infty]\mapsto [n_1, g''g'g, n_\infty^*, n_\infty].
\]

Following Example \ref{ex: quiver picture critical raviolo correspondence}, one may visualize the maps $\tilde p_{123}$, $\tilde j_{123}$, and $\tilde m_{123}$ as follows. The map $\tilde p_{123}$ is given by 
\[
\begin{tikzcd}
    \Sh{V}\arrow[d, "n_1"]\arrow[r, dashed, "g"] & \Sh{V}'\arrow[d, "\alpha'"]
    \\
    \Sh{V}\arrow[r, dashed, "g"] & \Sh{V}'\arrow[u, shift left, dashed, "n_2^*"]
\end{tikzcd}
\times
\begin{tikzcd}
     \Sh{V}'\arrow[d, "n_2"]\arrow[r, dashed, "g'"] & \Sh{V}''\arrow[d, "n_3"]
    \\
    \Sh{V}'\arrow[r, dashed, "g'"]& \Sh{V}''\arrow[u, shift left, dashed, "n_3^*"]
\end{tikzcd}
\times
\begin{tikzcd}
     \Sh{V}''\arrow[d, "n_3"]\arrow[r, dashed, "g''"] & \Sh{V}'''\arrow[d, "n_\infty"]
    \\
    \Sh{V}''\arrow[r, dashed, "g''"]& \Sh{V}'''\arrow[u, shift left, dashed, "n_\infty^*"]
\end{tikzcd}
\mapsfrom
\begin{tikzcd}
    \Sh{V}\arrow[d, "n_1"]\arrow[r, dashed, "g"] & \Sh{V}'\arrow[d, "n_2"]\arrow[r, dashed, "g'"] & \Sh{V}''\arrow[d, "n_3"] \arrow[r, dashed, "g''"] & \Sh{V}'''\arrow[d, "n_\infty"]
    \\
    \Sh{V}\arrow[r, dashed, "g"] & \Sh{V}'\arrow[u, shift left, dashed, "n_2^*"]\arrow[r, dashed, "g'"]\arrow[u]& \Sh{V}''\arrow[r, dashed, "g''"]\arrow[u, shift left, dashed, "n_3^*"]& \Sh{V}'''\arrow[u, shift left, dashed, "n_\infty^*"]
\end{tikzcd}
\]
The map $\tilde q_{1,123}$ is given by
\[
\begin{tikzcd}
    \Sh{V}\arrow[d, "n_1"]\arrow[r, dashed, "g"] & \Sh{V}'\arrow[d, "n_2"]\arrow[r, dashed, "g'"] & \Sh{V}''\arrow[d, "n_3"] \arrow[r, dashed, "g''"] & \Sh{V}'''\arrow[d, "n_\infty"]
    \\
    \Sh{V}\arrow[r, dashed, "g"] & \Sh{V}'\arrow[u, shift left, dashed, "n_2^*"]\arrow[r, dashed, "g'"]\arrow[u]& \Sh{V}''\arrow[r, dashed, "g''"]\arrow[u, shift left, dashed, "n_3^*"]& \Sh{V}'''\arrow[u, shift left, dashed, "n_\infty^*"]
\end{tikzcd}
\mapsto        
\begin{tikzcd}
    \Sh{V}\arrow[d, dashed, "w"]
    \\
    \Sh{V}\arrow[u, shift left, dashed,  "w^*"]
\end{tikzcd}
\times 
\begin{tikzcd}
    \Sh{V}\arrow[d, dashed, "\bar w"]
    \\
    \Sh{V}\arrow[u, shift left, dashed,  "\bar w^*"]
\end{tikzcd}
\times
\begin{tikzcd}
    \Sh{V}\arrow[d, "n_1"]\arrow[r, dashed, "g"] & \Sh{V}'\arrow[r, dashed, "g'"] & \Sh{V}''\arrow[r, dashed, "g''"] & \Sh{V}'''\arrow[d, "n_\infty"]
    \\
    \Sh{V}\arrow[r, dashed, "g"] & \Sh{V}'\arrow[r, dashed, "g'"]& \Sh{V}''\arrow[r, dashed, "g''"]& \Sh{V}'''\arrow[u, shift left, dashed, "n_\infty^*"]
\end{tikzcd}
\]
and the map $\tilde m^{\xi}_{123}$ is given by 
\[
\begin{tikzcd}
    \Sh{V}\arrow[d, "n_1"]\arrow[r, dashed, "g"] & \Sh{V}'\arrow[r, dashed, "g'"] & \Sh{V}''\arrow[r, dashed, "g''"] & \Sh{V}'''\arrow[d, "n_\infty"]
    \\
    \Sh{V}\arrow[r, dashed, "g"] & \Sh{V}'\arrow[r, dashed, "g'"]& \Sh{V}''\arrow[r, dashed, "g''"]& \Sh{V}'''\arrow[u, shift left, dashed, "n_\infty^*"]
\end{tikzcd}
\mapsto
\begin{tikzcd}
    \Sh{V}\arrow[d, "n_1"]\arrow[r, dashed, "g''g'g"] & \Sh{V}'''\arrow[d, "n_\infty"]
    \\
    \Sh{V}\arrow[r, dashed, "g''g'g"] & \Sh{V}'''\arrow[u, shift left, dashed, "n_\infty^*"]
\end{tikzcd}
\]
We now assign potentials along the lines of \S\ref{subsec: Coulomb branch action on quasimaps}.
\begin{itemize}
    \item The potential on $\tilde\Rt_{G,N}^2\times \Mt_{G,N}$ is the direct sum potential $\Sh{W}\boxplus\Sh{W}\boxplus\Sh{W}$.
    \item The potential on $\tilde\Rt^{(3),0,0}_{G,N}$ is 
    \begin{multline*}
        \Sh{W}_{0,0}\coloneqq(\Sh{W}\boxplus \Sh{W}\boxplus \Sh{W})\circ \tilde p_{123}
        =\Res\Tr\left(n_2^*(gn_1-n_2)\right) + \Res\Tr\left(n_3^*(g'n_2-n_3)\right) \\+\Res\Tr\left(n_\infty^*(g''n_3-n_\infty)\right)+\Res\Tr((\xi_2,j) n_{\infty}).
    \end{multline*}
    Since it is $G\fps\times G\fps$ equivariant, it descends to a potential $\Sh{W}_{1,1}$ on $\Mt^{(2),1}_{G,N}$.
    \item The potential $\Sh{W}_{2,2}$ on $\Mt^{(2),2,2}_{G,N}$ is given by 
    \begin{equation}
    \label{eq: W_22}
    \Sh{W}_{2,2}\coloneqq\Res\Tr(w^*w)+\Res\Tr(\bar w^*\bar w)+\Res\Tr\left( n_\infty^*\left(g''g'g n_1-n_\infty\right)\right)+\Res\Tr((\xi_2,\iota) n_{\infty}).
    \end{equation}
    \item The potential on $\Mt^{(2),3,3}_{G,N}$, denoted by $\Sh{W}_3$, is the composition 
    \[
    \Sh{W}_{3,3}\coloneqq\Sh{W}_2\circ \tilde q_{123,2}=\Res\Tr\left( n_\infty^*\left(g''g'gn_1-n_\infty\right)\right)+\Res\Tr((\xi_2,\iota) n_{\infty}).
    \]
    \item The potential on $\Mt^{(2),4,4}_{G,N}$ is 
    \[
    \Sh{W}_{4,4}\coloneqq\Sh{W}\circ \tilde m_{123}=\Res\Tr\left( n_3^*\left(g'gn_1-n_3\right)\right)+\Res\Tr((\xi_2,\iota) n_{\infty})
    \] 
\end{itemize}
A simple generalization of the argument of Lemma \ref{lemma: critical raviolo 2} shows that $\Sh{W}_{2,2}\circ \tilde j_{123}=\Sh{W}_{1,1}$. It then follows that all the potentials assigned to the vertices of the diagram \eqref{eq: interpolating corespondence} are compatible with the morphisms of the latter.
Therefore, recasting the argument of Proposition \ref{prop: Coulomb action def} in this setting, we get a well defined morphism
\begin{equation}
    \label{eq: map induced by interpolating correspondence}
    {}^{\phi}\eta_{123}\coloneqq \tilde m_{123,*} \circ(\tilde q_{3,123}^*)^{-1} \circ (\tilde q_{2,123,*})^{-1}  \circ \tilde j_{123,*}\circ (\tilde q_{1,123}^*)^{-1} \circ \tilde p_{123}^*
\end{equation}
It fits on the diagonal of the diagram 
\begin{equation}
    \label{eq: splitted diagram action Coulomb on quasimaps}
    \begin{tikzcd}
        {}^{\phi}\Coul^{T,\vir}_{G,N}\otimes {}^{\phi}\Coul^{T,\vir}_{G,N}\otimes \HO_{T}(\Mt_{G,N}, \phip{\Sh{W}^{\xi}}\DD\BoQ)\arrow[r,"\id\otimes {}^{\phi}\eta^{\xi}"] \arrow[d, swap, "{}^{\phi}\eta\otimes \id"] \arrow[dr, "{}^{\phi}\eta^{\xi}_{123}"]& {}^{\phi}\Coul^{T,\vir}_{G,N}\otimes \HO_T(\Mt_{G,N}, \phip{\Sh{W}^{\xi}}\DD\BoQ)\arrow[d, "{}^{\phi}\eta^{\xi}"]
        \\
        {}^{\phi}\Coul^{T,\vir}_{G,N}\otimes \HO_T(\Mt_{G,N}, \phip{\Sh{W}^{\xi}}\DD\BoQ)\arrow[r, "{}^{\phi}\eta^{\xi}"] &  \HO_T(\Mt_{G,N}, \phip{\Sh{W}^{\xi}}\DD\BoQ)
    \end{tikzcd}
\end{equation}
Therefore, to prove the theorem, it suffices to show that both triangles of \eqref{eq: splitted diagram action Coulomb on quasimaps} commute. We now present the full argument ensuring commutativity of the lower triangle. The one for the upper triangle is essentially the same, so it will only be sketched.  Consider the diagram
\[
\begin{tikzcd}[column sep=1em]
\tilde\Rt^2\times\Mt& \tilde\Rt^{(2),0}\times\Mt \arrow[l, swap, "\tilde p_{12}"] \arrow[r, "\tilde q_{12,1}"] & \tilde\Rt^{(2),1}\times \Mt  \arrow[r, "\tilde j_{12}"] & \tilde\Rt^{(2),2}\times \Mt  & \tilde\Rt^{(2),3}\times \Mt \arrow[l, swap, "\tilde q_{2,12}"] \arrow[r, "\tilde q_{3,12}"]  &\tilde\Rt^{(2),4}\times \Mt   \arrow[r, "\tilde m_{12}"] & \tilde\Rt\times \Mt
\\
&
\Mt^{(3),0,0}\arrow[u]\arrow[ul, "\tilde p_{123}"]\arrow[r]\arrow[dr, swap, "\tilde q_{1, 123}"]& \Mt^{(3),1,0}\arrow[u]\arrow[d]\arrow[r]&\Mt^{(3),2,0} \arrow[u]\arrow[d]&\Mt^{(3),3,0}\arrow[u]\arrow[l] \arrow[d]\arrow[r]&\Mt^{(3),4,0}\arrow[u]\arrow[r]\arrow[d]&\Mt^{(2),0} \arrow[u, swap, "\tilde p"] \arrow[d, "\tilde q_1"]
\\
&& \Mt^{(3),1,1} \arrow[r]\arrow[dr,  swap, "\tilde j_{123}"]&\Mt^{(3),2,1} \arrow[d]&\Mt^{(3),3,1}\arrow[l]\arrow[r]\arrow[d]&\Mt^{(3),4,1}\arrow[r]\arrow[d]&\Mt^{(2),1} \arrow[d, "\tilde j"]
\\
&&&
\Mt^{(3),2,2}
&
\Mt^{(3),3,2}
  \arrow[l]
  \arrow[r]
&
\Mt^{(3),4,2}
  \arrow[r]
&
\Mt^{(2),2}
\\
&&&&
\Mt^{(3),3,3}
  \arrow[ul, "\tilde q_{2, 123}"]
  \arrow[u]
  \arrow[r]
  \arrow[dr, swap, "\tilde q_{3, 123}"]
&
\Mt^{(3),4,3}
  \arrow[u]
  \arrow[r]
  \arrow[d]
&
\Mt^{(2),3}
  \arrow[u, swap, "\tilde q_2"]
  \arrow[d, "\tilde q_3"]
\\
&&&&&
\Mt^{(3),4,4}
  \arrow[r]
  \arrow[dr, swap, "\tilde m_{123}"]
&
\Mt^{(2),4}
  \arrow[d, "\tilde m"]
\\
&&&&&&
\Mt
\end{tikzcd}
\]

The top row is the correspondence \eqref{eq: critical Coulomb branch correspondence}, multiplied by $\Mt$ from the right. The rightmost column is \eqref{eq: correspondence Coulomb branch action}. The main diagonal is \eqref{eq: interpolating corespondence}. A useful way to thinking about the diagram is the following. Any vertical arrow in the first row is a closed embedding followed by a $G\fps$ quotient. Any map from a space of the from $X^{\bullet,i,0}$ (resp. $X^{\bullet,0,i}$) to a space of the form $X^{\bullet,i,1}$ (resp. $X^{\bullet,1,i}$) is a $G\fps$ quotient. Any map from a space of the from $X^{\bullet,i,1}$ (resp. $X^{1,i}$) to a space of the form $X^{\bullet,i,2}$ (resp. $X^{\bullet,2,i}$) is a closed embedding. Any map from $X^{\bullet,i,2}$ (resp. $X^{\bullet,2,i}$) to $X^{\bullet,i,3}$ (resp. $X^{\bullet,3,i}$) is a closed embedding. Any map from $X^{\bullet,i,3}$ (resp. $X^{\bullet, 3,i}$) to $X^{\bullet,i,4}$ (resp. $X^{\bullet,4,i}$) is a vector bundle map. Any map from $X^{\bullet,i,4}$ (resp. $X^{\bullet,i,4}$) to $X^{\bullet,i}$ is a ind-proper map. All squares are Cartesian and the triangles commute. All the vertices and arrows of this diagram are determined by these conditions.

The commutativity of the lower triangle in \eqref{eq: splitted diagram action Coulomb on quasimaps} follows from \eqref{eq: map induced by interpolating correspondence}, the definitions of ${}^{\phi}\eta$ and ${}^{\phi}\eta^{\xi}$—see \eqref{eq: formula critical Coulomb branch multiplication} and \eqref{eq: legit action}—and the standard compatibility between pullback and pushforward in Cartesian squares.

 The argument for the commutativity of the upper triangle is similar, so we only focus on its subtle part. In order to straightforwardly adapt the argument above for the lower triangle, it is convenient to use a variant of \eqref{eq: interpolating corespondence}, where the map $\tilde j_{123}$ is replaced by a map $\tilde k_{123}$ defined by 
 \[
[n_1, g, n_2^*, n_2, g', n_3^*, n_3,g'', n_\infty^*, n_\infty]\mapsto [u^*, u, \bar u^*, \bar u, n_1, g, g',g'', n_\infty^*, n_\infty],
\]
where now
\[
(w^*,w)=(g'gn_3^*-gn_2^*,gn_2- n_1)
\quad 
(\bar w^*, \bar w)=(g''g'gn_\infty^*-g'gn_3^*,g'gn_3-n_1)
\]
are understood.
where now
\[
(u^*,u)=(g''g'gn_\infty^*-gn_3^*, g'gn_2-gn_1) 
\qquad 
(\bar u^*, \bar u)=(g''g'gn_\infty^*-gn_2^*, gn_2-n_1) 
\]
A straightforward adaptation of the argument above implies the map 
\begin{equation}
    \label{eq: map induced by interpolating correspondence v2}
    {}^{\phi}\theta_{123}\coloneqq \tilde m_{123,*} \circ(\tilde q_{3,123}^*)^{-1} \circ (\tilde q_{2,123,*})^{-1}  \circ \tilde k_{123_*}\circ (\tilde q_{1,123}^*)^{-1} \circ \tilde p_{123}^*
\end{equation}
satisfies ${}^{\phi}\theta_{123}={}^{\phi}\eta^{\xi}\circ (\id\otimes {}^{\phi}\eta^{\xi})$. Therefore, to prove the lemma it suffices to show that ${}^{\phi}\theta_{123}={}^{\phi}\eta_{123}$. We will prove the stronger statement 
\begin{equation}
    \label{eq: last eq to prove associativity}
    (\tilde q_{3,123}^*)^{-1} \circ (\tilde q_{2,123,*})^{-1}  \circ \tilde j_{123_*}=(\tilde q_{3,123}^*)^{-1} \circ (\tilde q_{2,123,*})^{-1}  \circ \tilde k_{123,*},
\end{equation}
which implies that ${}^{\phi}\theta_{123}={}^{\phi}\eta_{123}$ by \eqref{eq: map induced by interpolating correspondence}  and \eqref{eq: map induced by interpolating correspondence v2}. Consider the diagram 
\[
\begin{tikzcd}
    \Mt^{(3),1}\arrow[r, "\tilde j_{123}"] \arrow[dr, swap, "\tilde k_{123}"] & \Mt^{(3),1,1} \arrow[d, "c"] & \Mt^{(3),2,2}\arrow[d, "c_0"] \arrow[l, swap, "\tilde q_{2,123}"]\arrow[r, "\tilde q_{3,123}"]& \Mt^{(3),3,3}
    \\
    & \Mt^{(3),1,1} &\arrow[l, "\tilde q_{2,123}"] \Mt^{(3),2,2}\arrow[ur, swap, "\tilde q_{3,123}"]
\end{tikzcd}
\]
where $c$ is given by\footnote{Using the notation from the definition of $\tilde j_{123}$ and $\tilde k_{123}$, the map $c$ is just the change of variables $(w^*, w)=(u^*, u-\bar u)$ and $(\bar w^*, \bar w)=(u^*+\bar u^*, \bar u)$.}
\[
[w^*, w, \bar w^*, \bar w, n_1, g, g',g'', n_\infty^*, n_\infty]\to [\bar w^*, \bar w-w, w^*+\bar w^*, w, n_1, g, g',g'', n_\infty^*, n_\infty]
\]
The map $c_0$ is the induced by pullback map. It is easy to check that $\Sh{W}_{2,2}\circ c=\Sh{W}_{2,2}$, using \eqref{eq: W_22}. Therefore, passing to critical cohomology, we get $\tilde k_{123,*}=\tilde j_{123,*}\circ c_*$, $c_*\circ \tilde q_{2, 123}=\tilde q_{2, 123}\circ c_{0,*}$, and $c_{0,*}\circ \tilde q_{3,123}^*=\tilde q_{3,123}^*$ (the latter because $c_0$ is an isomorphism, so $c_{0,*}=c_0^*$). Equation \eqref{eq: last eq to prove associativity} follows. 

It remains to prove that the unit $[\Rt_{e}\times N^*\lfps]\in {}^{\phi}\Coul^{\vir}_{N,G}$ acts by the neutral element. Let $e\in \pi_1(G)$ be the unit. Recall that $\Gr_e=\pt$ and hence $\Rt_e\cong N\fps$. Recall also that the class $[\Rt_{e}\times N^*\lfps]\in {}^{\phi}\Coul^{\vir}_{N,G}$ is defined as the pushforward of the fundamental class of $[\Rt_{e}\times N^*\lfps]$ under the inclusion $\Rt_e\times N^*\lfps\mapsto \tilde \Rt_e$, whose image consists of tuples $[n_1, g, n_2, n_2^*]$ of the form $[n_1, \id, n_1, n_2^*]$. Consider the following diagram
\begin{equation}
\label{eq: diagram unit is a unit}
    \begin{tikzcd}
        \tilde\Rt_{e}\times \Mt & \Mt^{(2),0}_{e}\arrow[l, swap, "\tilde p"]\arrow [r, "\tilde q_1"] & \Mt^{(2),1}_{e}\arrow[r, "\tilde j"]& \Mt^{(2),2}_{e}
        \\
        (\Rt_e\times N^*\lfps)\times \Mt \arrow[ur, phantom, "A"]\arrow[u, "j\times \id"]\arrow[d, swap, "p\times \id"] & N^*\lfps \times \widehat \Mt \arrow[u, "j'"]\arrow[l, swap, "\bar p"]\arrow[r, "\bar q_1"]\arrow[d]\arrow[ur, phantom, "B"] & N^*\lfps \times_{G\fps} \widehat \Mt\arrow[r, "\bar q_2"]\arrow[u, "j''"]\arrow[d] \arrow[ur, phantom, "C"]& \Mt^{(2),3}_{e}\arrow[u, swap, "\tilde q_2"]\arrow[d, "\tilde q_3"]
        \\
        \Mt &  \widehat \Mt \arrow[l, swap, "p"]\arrow[r, "p"]&\Mt \arrow[r, equal]  &\Mt^{(2),4}_{e} \arrow[r, equal, "\tilde m"] & \Mt.
\end{tikzcd}
\end{equation}
Here, we have defined
\begin{align*}
    \widehat \Mt&\coloneqq  N\fps\times G\lfps^{\leq 0} \times_{G[t^{-1}]} (N^*\lfps\times N^{\xi}),
\end{align*}
so that $\Mt=\widehat \Mt/G\fps$. Let $p:\widehat \Mt\to \Mt$ and $\bar q_1: N^*\lfps\times\widehat \Mt\to N^*\lfps\times_{G\fps} \widehat \Mt$ be the quotient maps. The maps $j'$ and $\bar p$ are defined to make square A Cartesian. Therefore, $j'$ is given by 
\begin{equation*}
    \label{eq: map j' in proof unit is unit}
    (n_2^*, n_1, [g, n_\infty^*, n_\infty])\mapsto (n_1, \id, n_2^*, n_1, [g, n_\infty^*, n_\infty]).
\end{equation*}
The map $j''$ is induced by $j'$ by taking $G\fps$-quotients. The map $\bar q_2$ is defined via the assignment\footnote{Notice that there is no variable $g$ in this argument; in fact we can always set $g=\id$ as a representative of an element in $\Gr_e=G\fps\backslash G\fps=\pt$. However, we keep using the variable $g'$ to denote a variable in the prequotient of $\Mt$ to match formulas \eqref{eq: first two maps critical ravolo action}-\eqref{eq: last map critical ravolo action}.}
\[
[n_2^*, n_1, g', n_\infty^*, n_\infty]\mapsto [g'n_3^*-n_2^*, n_1, \id, g', n_\infty^*, n_\infty].
\]
Commutativity of square C follows by explicit computation:
\[
\begin{tikzcd}[arrows=mapsto]
    \left[n_1, \id, n_2^*, n_1, g', n_\infty^*, n_\infty\right]\arrow[r, "\tilde j"] & \left[g'n_\infty^*-n_2^*,0,n_1, \id, g', n_\infty^*, n_\infty\right]
    \\
    \left[n_2^*, n_1, g', n_\infty^*,n_\infty\right]\arrow[r, "\bar q_2"]\arrow[u, "j''"] & \left[g'n_\infty^*-n_2^*, n_1, \id, g', n_\infty^*, n_\infty\right]\arrow[u, swap, "\tilde q_2"].
\end{tikzcd}
\]
Notice that here we have crucially used the fact, in general, the second entry of the targe of $\tilde j$ is assigned to $gn_2-n_1$ (cf. \eqref{eq: second two pair of maps critical ravolo action}), but here $g=\id$ and $n_1=n_2$ (up to conjugations by $G\fps$ that cancel each other). The same observation also implies that $\Mt=\Mt^{(2),4}_e$ and $\tilde m=\id $, as stated in \eqref{eq: diagram unit is a unit}. Checking that the lower squares commute is left to the reader. All the spaces in the outer frame of \eqref{eq: diagram unit is a unit} admit a potential, and we assign a potential to the remaining space $N^*\lfps \times \widehat \Mt$ and $N^*\lfps \times_{G\fps} \widehat \Mt$ by pullback along any of the maps in the diagram. We also leave to the reader to check that this assignment is unambiguous. 

The proof is now a computation. By definition, the action of the unit $[\Rt_{e}\times N^*\lfps]\in {}^{\phi}\Coul^{\vir}_{N,G}$ is given by the composition 
\begin{equation}
\label{eq: unit action}
    \tilde m_* \circ (\tilde q_3^*)^{-1}\circ (\tilde q_{2,*})^{-1}\circ \tilde j_*\circ (\tilde q_1^*)^{-1}\circ \tilde p^*\circ (j\times \id )_*\circ (p\times \id)^*.
\end{equation}
Since the squares A and B are Cartesian, we obtain 
\[
    (\tilde q_1^*)^{-1}\circ \tilde p^*\circ (j\times \id )_*= j''_*\circ (\bar q_1^*)^{-1}\circ \bar p^*,
\]
while from the commutativity of $C$ we get $(\tilde q_{2,*})^{-1}\circ \tilde j_*\circ j''_*=\bar q_{2,*}$. Combining the last two equations, we deduce that \eqref{eq: unit action} is equal to 
\begin{equation}
\label{eq: unit action 2}
\tilde m_* \circ (\tilde q_3^*)^{-1} \circ \bar q_{2,*}\circ (\bar q_1^*)^{-1}\circ \bar p^* \circ (p\times \id)^*.
\end{equation}
But $\bar q_2$ is an isomorphism, so $\bar q_{2,*}=(\bar q_2)^{-1}$. Therefore, using the commutativity of the lower squares in \eqref{eq: diagram unit is a unit} and the previous observation that $\tilde m=\id$, we deduce that \eqref{eq: unit action 2} is equal to $(p^*)^{-1}\circ p^*=\id$. Therefore, the action of the unit is the identity map. The proof is now complete.

\end{proof}

\subsection{Minuscule monopoles action}
\label{subsec: Minuscule monopoles action}

In this section, we give a simple formula for the action of the monopole operators $M_{f, \lambda}\in \Coul^{T,\vir}_{N,G}$, which, by Proposition \ref{subsec: coulomb branch generators}, provide a set of generators for the Coulomb branch algebra. Fix a minuscule weight $\lambda$ and let $\Rt_{\lambda}\subset \Rt$ be the corresponding variety. 

Let $\Mt^{(2), \xi}_{\lambda}$ be the subvariety of $\Mt^{(2),4, \xi}$ consisting of those points $[n_1, g, g', n_\infty, n_\infty^*]$ such that $gn_1$ is holomorphic and $g\in G\lfps$ is in the $G\fps\times G\fps$ orbit determined by $\lambda$. We have natural morphisms 
\begin{equation}
    \label{eq: simple critical correspondence mu=inuscule}
\begin{tikzcd}
   \Rt_{\lambda}\times \Mt^{\xi} & \Mt^{(2), \xi}_{\lambda}\arrow[r, "q_{\lambda}"]\arrow[l, swap, "p_{\lambda}"]& \Mt^{\xi}
\end{tikzcd}
\end{equation}
explicitly given by 
\[
\begin{tikzcd}[arrows=mapsto]
    (\left[n_1,g\right], \left[gn_1,g', n_\infty, n_\infty^*\right])& \left[n_1, g, g', n_\infty, n_\infty^*\right] \arrow[r]\arrow[l]& \left[n_1, g'g, n_\infty, n_\infty^*\right].
\end{tikzcd}
\]
The potential on $\Rt_{\lambda}\times \Mt^{\xi}$ is pulled back from $ \Mt^{\xi}$. It is easy to check that $p_{\lambda}^* \Sh{W}^{\xi}=q_{\lambda}^* \Sh{W}^{\xi}$.
Consider the commutative diagram 
    \begin{equation}
        \label{eq: diagram action minuscules}
    \begin{tikzcd}
         \tilde \Rt_{\delta_i}\times\Mt^{\xi} & \Mt^{(2),0,\xi}\arrow[r, "\tilde q"]\arrow[l, swap, "\tilde p"] & \Mt^{(2),1,\xi}\arrow[r, "\tilde j"]& \Mt^{(2),2,\xi}
         \\
         (\Rt_{\lambda}\times N^*\lfps) \times\Mt^{\xi}\arrow[u, "j\times \id"]\arrow[d] &  N^*\lfps \times \widehat \Mt^{(2), \xi}_{\lambda}\arrow[u, "j'"]\arrow[l, swap, "\bar p"]\arrow[r, "\bar q_1"]\arrow[d]& N^*\lfps \times_{G\fps} \widehat \Mt^{(2), \xi}_{\lambda}\arrow[u, "j''"]\arrow[d] \arrow[r, "\bar q_2"]& \Mt^{(2),3, \xi}\arrow[u, swap, "\tilde q_2"]\arrow[d, "\tilde q_3"]
         \\
        \Rt_{\lambda}\times\Mt^{\xi} &  \widehat \Mt^{(2), \xi}_{\lambda} \arrow[l, swap, "p"]\arrow[r, "q"]&\Mt^{(2), \xi}_{\lambda}\arrow[r]\arrow[dr, swap, "q_{\lambda}"]  \arrow[ll, bend left, "p_{\lambda}"]&\Mt^{(2), 4,\xi} \arrow[d, "\tilde m"] 
        \\
        &&& \Mt^{\xi}
    \end{tikzcd}
        \end{equation}
    This is the natural generalization of \eqref{eq: diagram unit is a unit}. That diagram served to express the action of the unit of the Coulomb branch, while this one expresses the action of a general monopole $\Rt_{\lambda}$. The map $j'$ is defined by making the top leftmost square Cartesian and the map $j''$ by taking the $G\fps$ quotient of $j'$. The two lower leftmost squares are constructed in the same way. The map $\bar q_2$ is given by
    \[
       [n_1, g, n_2^*, g', n_\infty^*, n_\infty]\mapsto [g'gn_3^*-gn_2^*, n_1, g, g', n_\infty^*, n_\infty].
    \]
    All the spaces in the outer frame of the diagram are equipped with a potential, and we one to the remaining spaces by pullback along any of the maps in the diagram. As in the previous case, this definition is unambiguous.
    The four leftmost square are Cartesian and the rightmost squares commute. The same argument of the proof of Theorem \ref{thm: Coulomb action QM} implies that the action of the class $[\Rt_{\lambda}\times N^*\lfps]\in {}^{\phi}\Coul^{\vir}$ is given by the assignment 
    \begin{equation}
    \label{eq: monopole action}
        \rho\mapsto q_{\lambda, *} p_{\lambda}^*(\rho)
    \end{equation} 
    By Theorem \ref{thm: dr Coulomb agebra}, this is equivalent to the statement that monopole class $M_{1, \lambda}\in [\Rt_{\lambda}]\in {}^{\phi}\Coul^{\vir}$ acts on $\HO_T(\Sh{M}^{\xi}, \phip{\Sh{W}}\BoQ^{\vir}_{\Sh{M}^{\xi}})$ by \eqref{eq: monopole action}. We remark that, for the well definiteness of the various pullback morphisms in critical cohomology we crucially use the fact that, since $\lambda$ is minuscule, all the spaces in \eqref{eq: diagram action minuscules} are smooth. 
    
    More generally, consider a dressed monopole $M_{f,\lambda}$. The dressing $f$ is a linear combination of Chern classes on $\Rt_{\lambda}$. By abuse of notation, we also denote by $f$ its pullback on $\Mt^{(2)}_{\lambda}$ via $q_{\lambda}$. With this notation, we obtain: 
    \begin{proposition}
        \label{prop: monopole action}
        The action of the dressed monopole operator $M_{f,\lambda}$ on  $\rho\in \HO_T(\Sh{M}^{\xi}, \phip{\Sh{W}^{\xi}}\BoQ^{\vir}_{\Sh{M}^{\xi}})$ is given by $\rho\mapsto q_{\lambda, *} \left(f\cap p_{\lambda}^*(\rho)\right)$.
    \end{proposition}

\subsection{Dimensionally reduced action in the Čech model}

Assume now that the base-point $p=(\xi_2,x_3, i,j)$ is polarized, i.e. it can we written as $p=(0,x_3, 0,j)$, and hence we have a canonical isomorphism $\HO_T(\QM^{\xi}_{\dd}(X), \varphi_{\QM^{\xi}_{\dd}(X)})\cong \HO_T^{\BM}(\QM^{\xi}(L),\BoQ^{\vir})$, cf. \S\ref{subsec: ev point and dim red} and \S\ref{subsec: dim red Čech model}. In \S\ref{subsec: Dimensionally reduced action hall} we have proved that, under the polarization hypothesis, the action of the 3d Hall algebra $\HCoha^T_{\tilde Q, \tilde W}$ on $\HO_T(\QM^{\xi}_{\dd}(X), \varphi_{\QM^{\xi}_{\dd}(X)})$ can be dimensionally reduced to an action of $\CoHA^T_{\Xi_Q}$ on $\HO_T^{\BM}(\QM^{\xi}(L),\BoQ^{\vir})$. In this section, we argue that under the same hypothesis, the action of the critical Coulomb branch ${}^{\phi}\Coul^{T,\vir}_{G,N}$ on $\HO_T(\QM^{\xi}_{\dd}(X), \varphi_{\QM^{\xi}_{\dd}(X)})$ from Theorem \ref{thm: Coulomb action QM} dimensionally reduces to an action of the BFN Coulomb branch $\Coul^{T, \vir}_{G,N}$ on $\HO_T^{\BM}(\QM^{\xi}(L),\BoQ^{\vir})$.

Recall the notation from \S\ref{subsec: dim red Čech model}. In the rest of the section we drop the subscripts $G$ and $N$ from the notation. Consider the following variant of the BFN diagram
\begin{equation}
\label{eq: raviolo correspondence quasimaps dimred}
\begin{tikzcd}
    \Rt\times \Nt^{\xi} \arrow[d, hookrightarrow]& p^{-1}(\Rt\times G\fps/G[t^{-1}] ) \arrow[d, hookrightarrow]\arrow[l, swap, "p"] \arrow[r, "q"] & q(p^{-1}(\Rt\times G\fps/G[t^{-1}] )\arrow[d, hookrightarrow]\arrow[r, "m"] & \Nt^{\xi}\arrow[d, hookrightarrow]
    \\
   \Rt\times \Ut^{\xi} & \Rt\times G\fps/G[t^{-1}] \arrow[l, swap, "\tilde p"] \arrow[r, "\tilde q"] & \Rt\times_{G\fps} G\fps/G[t^{-1}] \arrow[r, "\tilde m"] & \Ut^{\xi}
\end{tikzcd}
\end{equation}
It is obtained from \eqref{eq: Coulomb branch mult diagram} by replacing the spaces $\Rt\subseteq \Tt=N\fps\times_{G\fps } G\lfps$
with 
\[
\Nt^{\xi}=\{ n_{\infty}(t)-g(t)n_{0}(t)\}\subseteq \Ut^{\xi}=N\fps\times_{G\fps} G\lfps\times_{G_1[t^{-1}]} N^{\xi}[t^{-1}])
\]
Notice that since $n_{\infty}(t)$ is uniquely determined by $g(t)$ and $n_{0}(t)$, we may denote a point in $\Ut^{\xi}$ as a coset $[n_0(t), g(t)]$.
The maps in the second row are then diagrammatically defined by 
\begin{equation*}
    ([n,g_1], [g_1n,g_2])\mapsfrom ([n,g_1], g_2)\mapsto [[n,g_1],g_2]\mapsto [n,g_1g_2],
\end{equation*}
The exact same argument of \S\ref{subsec: Coulomb branch basics} combined with Proposition \ref{prop: dimensionaly reduced quasimap Čech model} gives a canonical map
\begin{equation}
    \label{eq: dim red Coulomb action}
    \Coul^{T,\vir}_{G,N, \dd}\otimes\HO_T^{\BM}(\QM^{\xi}(L),\BoQ^{\vir})\to \HO_T^{\BM}(\QM^{\xi}(L),\BoQ^{\vir})[I_{\dd}]
\end{equation}
where $I_{\dd}$ was defined in \eqref{eq: index shifted action}.
Moreover, a straightforward modification of the proof of \cite[Thm. 3.10]{BFNII} implies that \eqref{eq: dim red Coulomb action} is compatible with the multiplication on $\Coul^{\vir}_{G,N}$ and hence induces an action of $\Coul^{T,\vir}_{G,N}$ on $\HO_T^{\BM}(\QM^{\xi}(L),\BoQ^{\vir})$. Informally, this is the first statement of the following proposition. Since the proof is an easy adaptation of Theorem \ref{thm: dr Coulomb agebra}, we omit it.

\begin{proposition}
    We have a commutative diagram
    \[
    \begin{tikzcd}
    {}^{\phi}\Coul^{T,\vir}_{G,N, \dd}\otimes\HO_T(\QM^{\xi}_{\dd}(X), \varphi_{\QM^{\xi}_{\dd}(X)})\arrow[r]\arrow[d, swap,  "dr\otimes dr"] & \HO_T(\QM^{\xi}_{\dd}(X), \varphi_{\QM^{\xi}_{\dd}(X)})[I_{\dd}]\arrow[d, "dr"]
        \\
        \Coul^{T,\vir}_{G,N,\dd}\otimes\HO_T^{\BM}(\QM^{\xi}(L),\BoQ^{\vir})\arrow[r] & \HO_T^{\BM}(\QM^{\xi}(L),\BoQ^{\vir}[I_{\dd}])
    \end{tikzcd}
    \]
    where the vertical maps are the dimensional reduction morphisms from \eqref{eq: dimred quasimap Čech model formula} and Theorem \ref{thm: dr Coulomb agebra}.
    Therefore, the action of $\Coul^{T, \vir}_{G,N}$ on $\HO_T^{\BM}(\QM^{\xi}(L),\BoQ^{\vir})$ induced by \eqref{eq: dim red Coulomb action} is identified, via dimensional eduction, to the action of the critical Coulomb branch ${}^{\phi}\Coul^{T,\vir}_{G,N}$ on $\HO_T(\QM^{\xi}_{\dd}(X), \varphi_{\QM^{\xi}_{\dd}(X)})$ from Theorem \ref{thm: Coulomb action QM}.
\end{proposition}


\newpage

\paperpart{Coulomb branches, Yangians, and quasimap actions}
\label{part: compatibilities}

\section{From Hall algebras to Coulomb branches}
		
    \subsection{More on the positive Coulomb branch}
		\label{subsec: more on the positive Coulomb branch}
		Set $G=\GL_{\vv}$. In this section we analyze the structure of the spaces $\Gr^+$ and $\Rt^+$ introduced in \S\ref{subsec: Positive part of the Coulomb branch}. By definition $\Gr^+$ is the closed substack of $\Gr$ parametrizing a tuple of rank $\vv_i$ vector bundles $\Sh{V}_i$ on the disk together with injective morphisms of sheaves $\Sh{V}_i\hookrightarrow \Sh{O}^{\vv_i}$. For every $\dd\in \BoN^{Q_0}$, the component $\Gr^+_{\dd}\subset \Gr^+$ parametrizes those tuples $\Sh{V}_i\hookrightarrow \Sh{O}^{\vv_i}$ such that the cokernel has length $\dd_i$. Recall that $\Gr_{\dd}$ is in general an ind-scheme. On the other hand, the next proposition shows that each component $\Gr^+_{\dd}$ is a finite type scheme. 
        
        Given two non-negative integers $l,n\in \BoN$, let $\Quot(n,l)$ be the quot scheme on $\AA^1$, parametrizing length $l$ rank zero quotients $\Sh{O}^n_{\AA^1}\twoheadrightarrow Q$ of the structure sheaf $\Sh{O}^n_{\AA^1}$. It is a smooth $ln$-dimensional variety. It admits a canonical map $\Quot(n,l)\to \FCoh_0(\AA^1)$ to the stack of rank zero torsion sheaves. We identify $\FCoh_0(\DD)$ with the closed substack $\FCoh_0(\DD)\subset \FCoh_0(\AA^1)$ of rank zero sheaves on $\AA^1$ supported at zero. Now let
		\[
		\Quot_0(n,l)=\Quot(n,l)\times_{\FCoh_0(\AA^1)} \FCoh_0(\DD)
		\]
		the associated punctual quot scheme, parametrizing quotients $\Sh{O}^n_{\AA^1}\twoheadrightarrow Q$ with $Q$ supported at zero.
		\begin{proposition}
			\label{prop: quot and gr}
			The moduli space $\Gr^+_{\dd}$ is represented by a scheme of finite type. In fact, we have 
			\[
			\Gr^+_{\dd}\cong \prod_{i\in Q_0}\Quot_0(\vv_i, \dd_i),
			\]
			where $\Quot_{0}(\vv_i, \ww_i)$ is the punctual quot scheme of points on $\AA^1$.
		\end{proposition}
		\begin{proof}
			The proof follows from that any quotient $\BoC\fps^{ n}\twoheadrightarrow  Q$ factors through $(\BoC[t](t^N))^{\oplus n}$ for any $N\geq \text{length}(Q)$, therefore the moduli functors defining the two sides of the claim naturally equivalent. 
		\end{proof}
		Let $\FCoh_{0, \dd_i}(\DD)$ be the component of $\FCoh_{0}(\DD)$ of length $\dd_i$ quotients. Set $\FCoh_{0, \dd}(\DD)=\prod_{i\in Q_0} \FCoh_{0, \dd_i}(\DD)$. We have a canonical map
		\[
		b: \Gr^+_{\dd}\to \FCoh_{0, \dd}(\DD).
		\]
		The following observation will also be very useful. 
		\begin{lemma}
			\label{lemma: gr to coh is smooth}
			The morphism $b:\Gr^+_{\dd}\to \FCoh_{0,\dd}(\DD)$ is smooth of relative dimension $\vv^T\dd$.
		\end{lemma}
		\begin{proof}
			By Proposition \ref{prop: quot and gr}, it suffices to check that the map $\Quot_0(n, l)\to \Coh_{0,l}(\DD)$ is smooth. By base change, it further suffices to show that $\Quot(n, l)\to \Coh_{0,l}(\AA^1)$ is smooth. This is probably standard but we include the proof for completeness. The moduli space $\Coh_{0,l}(\AA^1)$ admits an universal sheaf $\Sh{Q}$ over $\AA^1 \times \Coh_{l}(\AA^1)$, which is finite flat over $\Coh_{l}(\AA^1)$. Consider the maps 
			\[
			\begin{tikzcd}
				\AA^1 &  \AA^1\times \Coh_{0,l}(\AA^1)\arrow[l, swap, "p"]\arrow[r, "q"] & \Coh_{0,l}(\AA^1)
			\end{tikzcd}
			\]
			Then the sheaf $\Sh{Hom}(p^*\Sh{O}^n, \Sh{Q})$ is flat and finite, and hence $q_*\Sh{Hom}(p^*\Sh{O}^n, \Sh{Q})$ is a vector bundle over $\Coh_{0,l}(\AA^1)$. Let $E\to \Coh_{0,l}(\DD)$ be its total space. Then its fibers over a torsion sheaf $Q$ are the vector spaces $\Hom(\Sh{O}^n, Q)$. In particular, the map $\Quot(n, l)\to \Coh_{0,l}(\AA^1)$ is the composition of the open inclusion $\Quot(n, l)\hookrightarrow E$ with the projection $E\to \Coh_{0,l}(\DD)$. Therefore, it is smooth.
		\end{proof}
		\begin{remark}
			The map $\Quot(n, l)\to \FCoh_{0,l}(\AA^1)$, and hence $\Gr^+_{\dd}\to \FCoh_{0,\dd}(\DD)$ are smooth but not surjective. In fact, $\Quot(n, l)\to \FCoh_{0,l}(\AA^1)$ is just surjective onto the open substack of $\FCoh_{0,\dd}(\DD)$ of sheaves with at most $n$ generators.
		\end{remark}
        \begin{remark}
            The Quot scheme $\Quot(n, l)$ can be explicitly described as the quotient  
            \[
            \{(A, v)\in \End(\BoC^l)\times \Hom(\BoC^n, \BoC^l)\; |\; \BoC[A] v=\BoC^n \}/\GL_l
            \]
            and the map $\Quot(n, l)\to \Coh_{0,l}(\AA^1)=\End(\BoC^l)/\GL(l)$ is identified with the projection $(A,v)\mapsto A$. Hence, the punctual quot scheme $\Quot_0(n, l)$ is the substack
            \[
            \{(A, v)\in \Sh{N}_l\times \Hom(\BoC^n, \BoC^l)\; |\; \BoC[A] v=\BoC^n \}/\GL_l
            \]
            where $\Sh{N}_l$ is the nilpotent cone. This description, which we do not need in what follows, makes smoothness of the morphism $\Quot_0(n, l)\to \FCoh_{0,l}(\DD)=\Sh{N}_l/\GL_l$ particularly manifest.
        \end{remark}
		
		Now set $N=\Rep_Q(\vv)$ and consider the moduli space $\Tt^+=G\lfps^+ \times_{G\fps} N\fps$ and the closed subspace $\Rt^+$, cf. \S\ref{subsec: Positive part of the Coulomb branch}. We also fix connected components $\Tt^+_{\dd}$ and $\Rt^+_{\dd}$. Notice that the canonical map 
		\[
		\Tt^+_{\dd}\to \Gr^+_{\dd}
		\]
		is a pro-vector bundle over the scheme $\Gr^+_{\dd}$. Explicitly, it is the inverse limit of affine bundles $\Tt^{+,k}\coloneqq G\lfps^+ \times_{G\fps} N\fps/t^kN\fps \to \Gr^+$. These spaces fit in a diagram 
		\begin{equation}
			\label{eq: various maps involved in B_{dd}}
			\Rt^+_{\dd}\hookrightarrow \Tt^+_{\dd}\xrightarrow{} \Gr^+_{\dd}\xrightarrow{b} \FCoh_{0,\dd}(\DD).
		\end{equation}
		Recall also the moduli space $\FCoh^T_{Q,0}(\DD)$ from \S\ref{subsec: Fixed support condition and nilpotency}. It parametrizes tuples of rank zero sheaves $Q=\{Q_i\}_{i\in Q_0}$ on $\DD$ together with morphisms $Q_{s(a)}\to Q_{t(a)}$ for all $a\in Q_0$. We denote by $\FCoh_{Q,0, \dd}(\DD)$ the component of length $\dd=(\dd_i)_{i\in Q_0}$ sheaves. We have canonical maps $\FCoh_{Q,0, \dd}(\DD)\to \FCoh_{0, \dd}(\DD)$. Our next goal is to analyze the geometry of the map $\Gr^+_{\dd}\to \FCoh_{Q,0, \dd}(\DD)$. Set
		\begin{equation}
			\label{eq: def B_{dd}}
			B_{\dd} \coloneq N\fps\times \Gr^+_{\dd}\times_{\FCoh_{0,\dd}(\DD)} \FCoh_{Q,0,\dd}(\DD)
		\end{equation}
		Notice that $B_{\dd}$ is the moduli space parametrizing tuples $(Q, g, \alpha, \beta)$ consisting of a tuple of rank zero length $\dd$ sheaves $Q=\{Q_i\}_{i\in Q_0}$ together with tuples of morphisms of $\Sh{O}_{\DD}$-modules such that $q=\{q_i\}_{i\in Q_0}$, $\alpha=\{\alpha_a\}_{a\in Q_1}$ and $\beta=\{\beta_{a}\}_{a\in Q_1}$ are tuples of morphisms 
		\[
		q_i: \Sh{O}^{\vv_i}\twoheadrightarrow Q_i \qquad \alpha_a: \Sh{O}^{\vv_{s(a)}}\to \Sh{O}^{\vv_{s(a)}}\qquad \beta_a: Q_{s(a)}\to Q_{t(a)}. 
		\]
		The moduli space $\Rt^+_{\dd}$ is the closed subspace of $B_{\dd}$ such that $q_{t(a)}\alpha_a-\beta_a q_{s(a)}$ for all $a\in Q_1$, i.e. such that the following diagram commutes:
		\[
		\begin{tikzcd}
			\Sh{O}^{\vv_{s(a)}}\arrow[d, "\alpha_a"]\arrow[r, twoheadrightarrow, "q_{s(a)}"] & Q_{s(a)}\arrow[d, "\beta_a"]
			\\
			\Sh{O}^{\vv_{t(a)}}\arrow[r, twoheadrightarrow, "q_{t(a)}"] & Q_{t(a)}
		\end{tikzcd}
		\]
		More precisely, we have the following result.
		\begin{proposition}
			\label{prop: virtual pullback data hall and Coulomb}
			There exists a rank $\vv^T\Q^T\dd$ bundle $E_{\dd}\to B_{\dd}$ together with a section $s: B_{\dd}\to E_{\dd}$ such that the composition $(s)^{-1}(0)\hookrightarrow B_{\dd} \to \Tt^+_{\dd}$ is an isomorphism onto $\Rt_{\dd}\hookrightarrow\Tt^+_{\dd}$:
			\begin{equation}
				\label{eq: positive Coulomb and inverse system}
				\begin{tikzcd}
					B_{\dd} \arrow[r] & \Tt^+_{\dd}
					\\
					(s)^{-1}(0)\arrow[r, "\cong"] \arrow[u, hookrightarrow] & \Rt^+_{\dd}\arrow[u, hookrightarrow]
				\end{tikzcd}
			\end{equation}
		\end{proposition}
		\begin{proof}
			The proof is constructive. Consider the moduli space $\FCoh_{0,\dd}(\DD)$ and the projections
			\[
			\begin{tikzcd}
				\AA^1 &  \AA^1\times \FCoh_{0,\dd}(\AA^1)\arrow[l, swap, "p"]\arrow[r, "q"] & \FCoh_{0,\dd}(\AA^1)
			\end{tikzcd}
			\]
			Over $\AA^1\times \FCoh_{0,\dd}(\DD)$ we have universal sheaves $\Sh{Q}_i$ for all $i\in Q_0$. Fix an edge $a\in Q_1$. By the same logic of the proof of Lemma \ref{lemma: gr to coh is smooth}, we obtain vector bundles 
			\[
			E'_a\coloneqq \text{Tot}\left( q_*\Sh{Hom}(p^*\Sh{O}^{\vv_{s(a)}}, \Sh{Q}_{t(a)})\right)\to \FCoh_{0,\dd}(\DD).
			\]
			Notice that $E_a$ is a vector bundle of rank $\dd_{t(a)}^T \vv_{s(a)}$. Set 
			\[
			E'_{\dd}\coloneqq\bigoplus_{a\in Q_1} E'_a.
			\]
			By construction, it is a vector bundle of rank $\vv^T \Q^T \dd$ with fibers over a tuple of coherent sheaves $\{Q_i\}_{i\in Q_0}$ given by\footnote{We do not distinguish between $\Hom(\Sh{O}_{\DD}^{\vv_{s(a)}}, Q_{t(a)})$ and $\Hom(\Sh{O}_{\AA^1}^{\vv_{s(a)}}, Q_{t(a)})$ as these are canonically isomorphic.} 
			\[
			\bigoplus_{a\in Q_1} \Hom(\Sh{O}_{\DD}^{\vv_{s(a)}}, Q_{t(a)})
			\]
			Let $E_{\dd}$ be the pullback of $E'_{\dd}$ on $B_{\dd}$ along the composition $B_{\dd}\to \Tt_{\dd}\to \FCoh_{0,\dd}(\DD)$. The section $s$ is then simply the tautological section sending a tuple $(Q, q, \alpha, \beta)\in B_{\dd}$ to the tuple of homomorphisms $q_{t(a)}\alpha_a-\beta_a q_{s(a)}\in \Hom(\Sh{O}_{\DD}^{\vv_{s(a)}}, Q_{t(a)})$.
		\end{proof}
		
		\begin{remark}
			\label{remark: virtual pullback data hall Coulomb approximations}
			The scheme $B_{\dd}$ is, by definition, the inverse limit of a system of morphisms $B^{k+1}_{\dd}\to B^k_{\dd}$ induced by the inverse system $\Tt^{+,k+1}_{\dd}\to \Tt^{+,k}_{\dd}$. 
			Moreover, the vector bundle $E_{\dd}$ on $B_{\dd}$ is also pulled back from the finite approximations $B^k_{\dd}$ for any sufficiently large $k$ (depending on $\dd$). In fact, as soon as $k\geq \max_{i\in Q_0}\{\dd_i\}$, the hom groups $\Hom(\BoC\fps^{ \vv_{s(a)}}, Q_{t(a)})$ are invariant under translation of $\BoC\fps^{ \vv_{s(a)}}$ by $t^k\BoC\fps^{ \vv_{s(a)}}$ and likewise the same is true of $E$, which thus descends to $B^k_{\dd}$. The same argument ensures the well-definiteness of a compatible system of sections $s^k: B_{\dd}^k\to E^k_{\dd}$. In conclusion, for all $k\geq \max_{i\in Q_0}\{\dd_i\}$ we have well defined morphisms 
			\[
			\begin{tikzcd}
				B_{\dd}^k \arrow[r] & \Tt^{+,k}_{\dd}
				\\
				(s^k_{\dd})^{-1}(0)\arrow[r, "\cong"] \arrow[u, hookrightarrow] & \Rt^+_{\dd, k}\arrow[u, hookrightarrow]
			\end{tikzcd}
			\]
			and diagram \eqref{eq: positive Coulomb and inverse system} emerges as the inverse limit  of these diagrams. 
		\end{remark}
		We end this section with a further analysis of the canonical morphism $ b: B_{\dd}\to \FCoh_{Q,0,\dd}(\DD)$. By equations \eqref{eq: def B_{dd}} \eqref{eq: various maps involved in B_{dd}}, we have a canonical factorization
		\[
		B_{\dd}\xrightarrow{b_1} \Gr^+_{\dd}\times_{\FCoh_{0,\dd}(\DD)} \FCoh_{Q, 0,\dd}(\DD)\xrightarrow{b_2} \FCoh_{Q, 0,\dd}(\DD)
		\]
		where the second map is the projection. Informally, the map $b_1$ sends a closed point $(Q, g, \alpha, \beta) \in B_{\dd}$ to $(Q, g, \beta)$, and $b_2$ projects further to $(Q, \beta)$.
		\begin{lemma} $ $
			\label{eq: lemma factorization pullback hall Coulomb 1}
			\begin{enumerate}
				\item The morphism $b_1$ is an inverse limit of finite rank vector bundles 
				\[
				B^k_{\dd}\xrightarrow{b_1^k} \Gr^+_{\dd}\times_{\FCoh_{0,\dd}(\DD)} \FCoh_{Q, 0,\dd}
				\]
				\item The morphism $b_2$ is a smooth morphism of algebraic stacks of relative dimension $\vv^T\dd$.
			\end{enumerate}
		\end{lemma}
		\begin{proof}
			The proof of the first item follows directly from the definition of $B_{\dd}^k$, which keeps track of the morphisms $\alpha_i: \BoC\fps^{\vv_{s(a)}}\to \BoC\fps^{\vv_{t(a)}}$ up to degree $k$. For the second item, it suffices to notice that we have a Cartesian diagram 
			\[
			\begin{tikzcd}
				\Gr^+_{\dd}\times_{\FCoh_{0,\dd}(\DD)}\FCoh_{Q,0,\dd}(\DD)\arrow[d, ""]\arrow[r]& \Gr_{\dd}^+\arrow[d]
				\\ 
				\FCoh_{Q,0,\dd}(\DD) \arrow[r] & \FCoh_{0,\dd}(\DD).
			\end{tikzcd}
			\]
			Hence the statement follows from Proposition \ref{prop: quot and gr} and Lemma \ref{lemma: gr to coh is smooth}.
		\end{proof}

        \subsection{From Hall algebras to Coulomb branches}
		\label{subsec: From hall to Coulomb}

		As discussed in sections \S\ref{subsec: Q decorated coherent sheaves on A1} and \S\ref{subsec: Fixed support condition and nilpotency}, the Borel-Moore homology the stack of $Q$-valued rank zero coherent sheaves $\FCoh^T_{Q,0}(Q)$ forms a Hall algebra, which is canonically isomorphic to
		\[
		\HCoha^{T, \onil}_{\Xi_{Q}}: =\bigoplus_{\dd\in \BoN^{Q_0}}\HO(\FM_{\dd}^{T, \onil}(\Xi_Q), \DD\BoQ_{\FM_{\dd}^{T, \onil}(\Xi_Q)}),
		\]
		namely the Hall algebra of $\omega$-nilpotent representations of the quiver $\hat Q$, cf. \S\ref{subsec: Fixed support condition and nilpotency}. Notice that the virtual dimension of $\FM_{\dd}(\Xi_Q)$ is zero so we have $\DD\BoQ^{\vir}_{\FM_{\dd}^{T, \onil}(\Xi_Q)}=\DD\BoQ_{\FM_{\dd}^{T, \onil}(\Xi_Q)}$. In this section, we construct a morphism of algebras from the Hall algebra to the positive Coulomb branch. In other to make the morphism compatible with the cohomological grading, it is convenient to define the following cohomologically shifted ariant of the positive Coulomb branch algebra  
		\[
		 \HO^{\BM}_{G\fps_T}( \Rt^+, \BoQ)=\bigoplus_{\dd\in \BoN^{Q_0}}  \HO^{\BM}_{G\fps_T}( \Rt^+_{\dd}, \BoQ).
		 \]
        In this section, it will be important to consider the shifted version of the Coulomb branch $\Coul^+_{G,N}$ introduced in Remark \ref{rem: virtual positive Coulomb branch}, which is given by 
        \[
        \Coul^{T,\vir, +}_{G,N}\coloneqq\bigoplus_{\dd\in \BoN^{Q_0}} \HO^{\BM}_{G\fps}( \Rt^+_{\dd}, \BoQ^{\vir}[2s_{\dd}]).
        \]
		for $s_{\dd}=\vv^T\dd-\vv^T\Q^T\dd$. With slight abuse of notation, we will also refer to $\Coul^{T,\vir, +}_{G,N}$ as the positive Coulomb branch algebra. Its sole difference from the ordinary positive Coulomb branch $\Coul^{T,+}_{G,N}$ \ref{eq: positive Coulomb} is a shift in the cohomological degree. The goal this section is to define a morphism of  $\BoN^{Q_0}$-graded, cohomologically graded $\HO_T$-algebras
		\begin{equation}
			\label{eq: map hall to Coulomb}
			\theta: \HCoha^{T, \onil}_{\Xi_{Q}}\to \Coul^{T, \vir, +}_{G,N}.
		\end{equation}
		Fix $\dd\in \BoN^{Q_0}$. We begin by defining a morphism of $\HO_T$-modules 
		\begin{equation}
			\label{eq: map hall to Coulomb for a component}
			\theta_{\dd}: \HO(\FM_{\dd}^{T, \onil}(\Xi_Q), \DD\BoQ^{\vir}_{\FM_{\dd}^{T, \onil}(\Xi_Q)})\to \HO_{G\fps_T}( \Rt^+_{\dd}, \DD\BoQ_{ \Rt^+_{\dd}}[-2s_{\dd}]).
		\end{equation}
		so that \eqref{eq: map hall to Coulomb} will be defined as the direct sum $\theta=\bigoplus_{\dd\in \BoN^{Q_0}} \theta_{\dd}$.
		
		Using the notation of the previous section, we have a chain of morphisms
		\begin{equation}
		    \label{eq: maps from Coulomb to coh}
            \begin{tikzcd}
			b_1^*:\Rt^+_{\dd}\arrow[r, hookrightarrow, "a"]& B_{\dd}\arrow[r, "b_1"] &\Gr^+_{\dd}\times_{\FCoh_{0,\dd}(\DD)} \FCoh_{Q, 0,\dd}(\DD) \arrow[r, "b_2"] & \FCoh_{Q, 0,\dd}(\DD)
		\end{tikzcd}
		\end{equation}
		To remove clutter, set $C_{\dd}\coloneqq\Gr^+_{\dd}\times_{\FCoh_{0,\dd}(\DD)} \FCoh_{Q, 0,\dd}(\DD)$. It is easy to check that all the maps above are equivariant with respect to the action of $ G\fps$.  Its action on $\FCoh_{Q, 0,\dd}(\DD)$ is trivial, while all the other actions all induced by the definition of these spaces in terms of $\Rt^+_{\dd}$ and $\Gr^+_{\dd}$.
		By the second part of Lemma \ref{eq: lemma factorization pullback hall Coulomb 1}, the map $b_1$ is smooth of relative dimension $\vv^T\dd$, so it induces a pullback morphism in Borel-Moore homology
		\begin{equation}
			\label{eq: construction pullback hall Coulomb step 1}
			b_2^*: \HO(\FCoh^T_{Q, 0,\dd}, \DD\BoQ_{\FCoh_{Q, 0,\dd}})\to \HO_{G\fps_T}( C_{\dd}, \DD\BoQ_{ C_{\dd}}[-2\vv^T\dd])
		\end{equation}
		As usual, the pullback as well as the cohomological shifts need to be understood in terms of finite dimensional approximations. The well definiteness of the latter follows from the observation that for any $k\ggg 0$ the groups $K_k=\ker(G\fps\to G^k\fps)$ act trivially on $C_{\dd}$.
		
		By the first part of Lemma \ref{eq: lemma factorization pullback hall Coulomb 1}, the map $b_2$ is approximated by finite rank vector bundles and hence, with the usual conventions on cohomological shifts (see \S\ref{subsec: cohomology of pro and ind stacks}), defines a pullback morphism 
		\begin{equation}
			\label{eq: construction pullback hall Coulomb step 2}
			\HO_{G\fps_T}(C_{\dd}, \DD\BoQ_{C_{\dd}}[-2\vv^T\dd])\to
			\HO_{G\fps_T}(B_{\dd}, \DD\BoQ_{ B_{\dd}}[-2\vv^T\dd]).
		\end{equation}
		
		Finally, by Proposition \ref{prop: virtual pullback data hall and Coulomb} and the construction of Appendix \ref{app: Refined Gysin pullback}, we obtain a virtual pullback
		\begin{equation}
			\label{eq: construction pullback hall Coulomb step 3}
			a^!: \HO_{G\fps_T}(B_{\dd}, \DD\BoQ_{ B_{\dd}}[-2\vv^T\dd])\to
			\HO_{G\fps_T}( \Rt^+_{\dd}, \DD\BoQ_{ \Rt^+_{\dd}}[-2s_{\dd}])
		\end{equation}
        We remark that \eqref{eq: construction pullback hall Coulomb step 3} is really defined in terms of the finite dimensional approximations from Remark \ref{remark: virtual pullback data hall Coulomb approximations} although we suppress this extra piece of notation for the sake of clarity. We can now define the morphism \eqref{eq: map hall to Coulomb for a component} as the composition of \eqref{eq: construction pullback hall Coulomb step 1}, \eqref{eq: construction pullback hall Coulomb step 2}, and \eqref{eq: construction pullback hall Coulomb step 3}. The apparent discrepancy in the cohomological shifts is simply due to the fact that the actual cohomology $\HO_{G\fps_T}( \Rt^+_{\dd}, \DD\BoQ_{ \Rt^+_{\dd}})$ is defined in terms of finite dimensional approximations, each of which is downshifted by $2\dim{G^k\fps}+2\dim{N^k\fps}$. Therefore, in the limit, the shift by $-2\infty_{G,N}$ is built in in the definition of $\HO_{G\fps_T}( \Rt^+_{\dd}, \DD\BoQ_{ \Rt^+_{\dd}})$.
		
		\begin{example}
			As a sanity check for the shift conventions, it is useful to look closely at the example where $\dd=\ul 0$. In this case $\Coh_{Q,0,\ul 0}=\Gr^+_{\ul 0}=\Spec(\BoC)$ and $B_{\dd}=C_{\dd}=\Rt^+_{\dd}=N\fps$, so that $G\fps\backslash B_{\dd}=G\fps\backslash C_{\dd}=G\fps\backslash\Rt^+_{\dd}=\Map(\DD, N/G)$. Then 
			\[
			\HCoha_{\Xi_Q, \ul 0}^{T, \onil}= \HO(\B T, \BoQ_{\B T})=\HO_T
			\]
			and 
			\[
			\HO_{G\fps_T}^{\BM}( R^+_{\ul 0}, \BoQ)= \HO_{G\fps_T}^{\BM}(N\fps, \BoQ)
			\]
			The right hand side is a free $\HO_T$-module of rank one. It is defined up to an overall shift in the cohomological degree. More explicitly, for each $k\geq 0$ we have a finite type approximations $N^k\fps/G^k\fps=\Map(\DD_{k}, N/G)$, cf. \S\ref{subsec: Mapping stacks on the formal disk and related geometries}, and the cohomologies $\HO_T^{\BM}(N^k\fps/G^k\fps, \BoQ)$ are one dimensional $\HO_T$ modules generated by the fundamental class $[N^k\fps/G^k\fps]$, which sits in cohomological degree $-2\dim(N_k)+2\dim(G^k\fps)$. For any $k,l\geq 0$, the cohomologies $\HO_T^{\BM}(N^k\fps/G^k\fps, \BoQ)$ and $\HO_T^{\BM}(N^k\fps/G^k\fps, \BoQ)$ are identified as $\HO_T$-modules by sending the generator $[N^k\fps/G^k\fps]$ to $[N^l\fps/G^l\fps]$. The map $\theta_{\ul 0}$ is to be understood in terms of these approximations, and, for a given $k$, sends the unit $1\in \HO_T=\HCoha_{\Xi_Q, \ul 0}^{T, \onil}$ to the fundamental class $[N^k\fps/G^k\fps]\in \HO_T^{\BM}(N^k\fps/G^k\fps, \BoQ)$ (which is the unit in the Coulomb branch algebra).
		\end{example}
		We can now state the main result of this section. Its proof is lengthy and is deferred to Appendix \ref{sec: Proof of Theoorem Hall Coulomb app}.
		\begin{theorem}
        \label{thm: main thm positive Coulomb hall}
			The morphism \eqref{eq: map hall to Coulomb} is a morphism of $\BoN^{Q_0}$-graded cohomologically graded algebras. In other words, for any $\dd\in \BoN^{Q_0}$ and any decomposition $\dd=\dd'+\dd''$, the diagram 
			\[
			\begin{tikzcd}
				 \HCoha^{T, \onil}_{\Xi_{Q}, \dd'}\otimes \HCoha^{T, \onil}_{\Xi_{Q}, \dd''} \arrow[d, "\theta_{\dd'}\otimes \theta_{\dd''}"] \arrow[rr, "{\vmult^{\omega\nil}}{\dd',\dd''}"] & & \HCoha^{T, \onil}_{\Xi_{Q}, \dd}\arrow[d,  "\theta_{\dd}"]\\
                \Coul^{T,\vir, +}_{G,N, \dd'}\otimes\Coul^{T,\vir, +}_{G,N, \dd''}\arrow[rr, "\eta^+"] &  & \Coul^{T,\vir, +}_{G,N, \dd}
			\end{tikzcd}
			\]
			is commutative. 
		\end{theorem}

        \subsection{Explicit formula for spherical elements}
        Fix $i\in Q_0$ and consider $\HCoha^{T,\nil}_{\Xi_{Q},\delta_i}\subset \HCoha^{T, \onil}_{\Xi_{Q}}$, .cf \S\ref{subsec: Fixed support condition and nilpotency}. Since $\FM^{T, \onil}_{\Xi_Q,\delta_i}=\AA^{g_i}/\BoC^*$, where $g_i$ is the number of loops of the quiver $Q$ at $i\in Q_0$, the cohomology $\HCoha^{T, \onil}_{\Xi_{Q},\delta_i}$ is spanned by the classes 
        \begin{equation}
        \label{eq: spherical generators omega nilpotent coha}
            c_1(\Sh{L}_i)^r \cap [\FM_{\delta_i}(\Xi_Q)^{T, \onil}]
        \end{equation}
        for $r\geq 0$.
        Notice that the orbit $\Rt_{\omega_{i,1}}\subset \Rt^+$ is isomorphic to the Grassmannian $\Gr(1,\vv_i)$ parametrizing one dimensional quotients of $\BoC^{\vv_i}$. Let $\Sh{Q}_i$ denote the associated tautological line bundle on $\Gr(1,\vv_i)$. In the language of the affine Grassmannian, the orbit $\Rt_{\omega_{i,1}}$ parametrizes lattices $L$ such that $t\Sh{O}^{\vv_i}\subset L\subset \Sh{O}^{\vv_i}$ such that $\dim(\Sh{O}^{\vv_i}/L)=1$, and $\Sh{Q_i}$ is the tautological bundle with fiber $\Sh{O}^{\vv_i}/L$.
        \begin{lemma}
        \label{image generators hall in positive Coulomb}
            The image of \eqref{eq: spherical generators omega nilpotent coha} under the morphism $\eta^+$ is $c_1(\Sh{Q}_i)^r \cap [\Rt_{\omega_{i,1}}]$.
        \end{lemma}
        \begin{proof}
            Notice that for $\dd=\delta_i$, we have $\Rt_{\delta_i}=\Rt_{\omega_{i,1}}$. Therefore, pullback under the morphism $\Rt_{\delta_i}\to \FCoh^{Q,0,\delta_i}=\FM^{T, \onil}_{\delta_i}(\Xi_{Q})$ sends the fundamental class $[\FM_{\delta_i}(\Xi_Q)^{T, \onil}]$ to $[\Rt_{\omega_{i,1}}]$. The statement then follows by compatiblity of pullback with action by Chern classes. 
        \end{proof}

        \subsection{Extension to the shifted Yangian}
        \label{subsec: Extension to the shifted Yangian}

         Recall the notation from Remark \ref{rem: signed product hall algebras}. Combining the morphism \eqref{eq: map hall to Coulomb for a component} with Theorem \ref{thm: positive Coulomb in Coulomb} and Proposition  \ref{prop: dr nilp version} we get a canonical algebra morphism $\HCoha^{T, \circ}_{\tilde Q,\tilde W}\to \Coul^{T, \vir}_{G,N}$ and hence, composing with the canonical map $\Sh{SH}^{T,\nil, \circ} \hookrightarrow \HCoha^{T, \circ}_{\tilde Q, \tilde W}$, a morphism of algebras 
        \[
        \phi^+: \Sh{SH}^{T,\nil, \circ} \to \Coul^{T, \vir}_{G,N}.
        \]
        The goal of this section is to extend this morphism from the Yangian. In view of the sign twists in the definition of the multiplication of $\Sh{SH}^{T,\nil, \circ}$, it is convenient to incorporate the same signs in the Yangian. Hence, we define 
        \[
        \Yang_{\mu}^{T, \circ}(\tilde Q, \tilde W)=\Sh{SH}^{T,\nil, \circ} * \Sh{H}^0_Q * \Sh{SH}^{T,\nil, \circ,  \opp}/I
        \]
        where $I$ is \emph{same} two-sided ideal encoding the defining relations of the shifted Yangian $\Yang_{\mu}^{T}(\tilde Q, \tilde W)$ discussed in \S\ref{subsec: Shifted Yangian action}. Informally speaking, the only difference between $ \Yang_{\mu}^{T, \circ}(\tilde Q, \tilde W)$ and the shifted yangian $\Yang_{\mu}^{T}(\tilde Q, \tilde W)$ defined in \S\ref{subsec: Shifted Yangian action} are the sign conventions in the Serre relations. 
        
        The goal of this section is to define a morphism of algebras $ \Yang_{\mu}^{T, \circ}(\tilde Q, \tilde W)\to \Coul^{T,\vir}_{G,N}$.
        To achieve this, we will first construct morphisms 
        \[
         \phi^-: \Sh{SH}^{T,\nil, \circ, \opp}\to \Coul^{T,\vir}_{G,N} \qquad \phi^0: \Sh{H}_Q^0\to \Coul^{T,\vir}_{G,N}.
        \]
        and then show that the induced map $\phi: \Sh{SH}^{T,\nil} * \Sh{H}^0_Q * \Sh{SH}^{T,\nil, \circ, \opp}/I \to \Coul^{T,\vir}_{G,N}$ satisfies the relations of $I$. The first step is achieved in the following proposition. Recall that the orbit $\Rt_{\omega^*_{i,1}}\subset\Rt_{G,N}$ parametrizes lattices $M$ such that $\Sh{O}^{\vv_i}\subset M\subset t^{-1}\Sh{O}^{\vv_i}$ satisfying $\dim(M/\Sh{O}^{\vv_i})=1$, and $\Sh{S_i}$ is the tautological bundle with fiber $M/\Sh{O}^{\vv_i}$.
        \begin{proposition}
        \label{prop: map opposite hall to Coulomb generators}
            There exists a canonical morphism of $\HO_T$ algebras $\HCoha^{T,\onil, \opp}_{\Xi_Q} \to \Coul^{T,\vir}_{G,N}$ sending the spherical elements $c_1(\Sh{L}_i)^r \cap [\FM^{T, \onil}_{\delta_i}(\Xi_Q)]$ to $c_1(\Sh{S}_i)^r \cap [\Rt_{\omega^*_{i,1}}]$.
        \end{proposition}
        \begin{proof}
            By Proposition \ref{prop: involution coha} we have a morphism of algebras $\Theta_*: \HO \Sh{A}^{T,\onil, \opp}_{\Xi_Q}\to \HO \Sh{A}^{T^{\opp}, \onil}_{\Xi_{Q^{\opp}}}$ and therefore, combining Proposition \eqref{eq: map hall to Coulomb for a component} applied to the quiver $Q^{\opp}$ and \ref{thm: positive Coulomb in Coulomb}, we get a composite morphism 
            \begin{equation}
                \label{eq: involution coha}
                \HO \Sh{A}^{T,\onil, \opp}_{\Xi_Q}\to \HO \Sh{A}^{T^{\opp}, \onil, \opp}_{\Xi_{Q^{\opp}}}\to \Coul^{T^{\opp}, \vir}_{G,N^*}.
            \end{equation}
            Notice that the the target is, canonically, the Coulomb branch for the dual representation $N^*$ because $(\Rep_{Q}(\vv,\ww))^*=\Rep_{Q^{\opp}}(\vv,\ww)=N^*$.
            This is by construction a morphism of algebras. Moreover, by \cite[Remark 3.6.]{BFN_quiver}, we have an isomorphism of $\HO_T$-algebras 
            \begin{equation}
                \label{eq: isom Coulombs dual reps}
                \Coul^{T^{\opp}, \vir}_{G,N^*}\to \Coul^{T, \vir}_{G,N}
            \end{equation}
            induced by the isomorphism $\iota: \Rt_{G,N}\to \Rt_{G, N^*}$ defined by $(n, g)\mapsto (n^t, (g^{t})^{-1})$. Here the notation $T^{\opp}$ refers to the fact that the $\HO_T$-action is induced by the composition of the automorphism of $\HO_T$ discussed in Remark \ref{remark identification equivariant parameters involution}. Composing \eqref{eq: involution coha} with \eqref{eq: isom Coulombs dual reps} gives the sought after map, which is $\HO_T$-linear.
            
         It remains to evaluate this map on $c_1(\Sh{L}_i)^r \cap [\FM^{T, \onil}_{\delta_i}(\Xi_Q)]$. Firstly, notice that the orbit $\Rt_{\omega_{i,1}}\subset \Rt_{G,N}$ is mapped by $\iota$ isomorphically onto $\Rt_{\omega^*_{i,1}}$, so we have an induced isomorphism $\iota': \Rt_{\omega^*_{i,1}}\to \Rt_{\omega_{i,1}}$. Therefore, applying Lemma \ref{image generators hall in positive Coulomb} we deduce that the map $\HCoha^{\onil, \opp}_{\Xi_Q} \to \Coul^{T,\vir}_{G,N}$ sends the spherical elements $c_1(\Sh{L}_i)^r \cap [\FM^{T, \onil}_{\delta_i}(\Xi_Q)]$ to $c_1((\iota')^*\Sh{Q}_i)^r \cap [\Rt_{\omega^*_{i,1}}]$. Therefore, it remains to prove that $(\iota')^*\Sh{Q_i}=\Sh{S_i}$. For this, it suffices to assume that $N=0$. But then $\iota: \Gr_G\to \Gr_G$ is the map sending a lattice $t\Sh{O}^{\vv_i}\subset L\subset \Sh{O}^{\vv_i}$ to its dual $\Sh{O}^{\vv_i}\subset L^\vee\subset t^{-1}\Sh{O}^{\vv_i}$, Therefore, the statement follows from the fiberwise definition of $\Sh{Q}_i$ and $\Sh{S}_i$.
        \end{proof}

        We now proceed with the construction the morphism $\phi: \Yang_{\mu}^{T, \circ}(\tilde Q, \tilde W)\to  \Coul^{T,\vir}_{G,N}$. Fix a maximal torus $S$ of $G$. Following \cite[Remark II.5.23]{BFNII}, we have an algebra embedding 
        \begin{equation}
            \label{eq: inclusion Coulomb branch in diff op}
            z^* \circ (i_*)^{-1}: \Coul^{T,\vir}_{G,N}\to \Coul^{T,\vir}_{S,0,\loc},
        \end{equation}
        where $\Coul^{T,\vir}_{S,0}$ is the Coulomb branch for the pure gauge theory for $S$ and the localization $\Coul^{T,\vir}_{S,0,\loc}$ can be explicitly described as
        \[
        \Coul^{T,\vir}_{S,0, \loc}= \langle \epsilon_1^{\pm}, \epsilon_i, w_{i,n}, u^\pm_{i,n}, (\alpha +m\epsilon_1)^{-1}\rangle_{\alpha: \text{root of $G$}, m\in \BoZ}.
        \]
        Here $\{w_{i,n}\}_{\substack{i \in Q_0 \\ 1\leq n\leq \vv_i}}$ are coordinates on $\Lie(S)$, $\epsilon_i$ are the extra torus equivariant parameters, and $u^\pm_{i,n}$ are the shift operators 
        \[
        u^\pm_{i,n}: w_{i,n}\mapsto w_{i,n}\pm\epsilon_1
        \]
        In other words, $\Coul^{T,\vir}_{S,0, \loc}$ is a localization of the ring of difference operators on $\Lie(S)$ valued in the polynomial ring $\BoQ[\epsilon_1^{\pm}, \epsilon_i]$.
        The following lemma is a straightforward adaptation \cite[\S A(ii)]{BFN_quiver}. The formulas there are computed via equivariant localization. The slightly differences in our formulas compared to those in \emph{loc.cit.} are due to our choices of equivariant weightings\footnote{To recover the conventions from \cite{BFN_quiver}, it suffices to set $\epsilon_{a}=-\epsilon_1/2$ for all $a\in Q_1$.}
      \begin{lemma}
      \label{lemma: formulas generators hall in Coulomb}
      Let $f$ be a polynomial function. The following formulas hold.
      \begin{multline*}
          z^* \circ (i_*)^{-1}\left(f(c_1(\Sh{Q}_i)) \cap [\Rt_{\omega_{i,1}}]\right)
          \\
          =\sum_{r=1}^{\vv_i}\frac{ f(w_{i,r})\prod_{a\in Q_1: \; s(a)=i} \prod_{\substack{s=1 \\ t(a)\neq i \text{ or } s\neq r}}^{\vv_{t(a)}} (w_{t(a),s}-w_{i,r}+\epsilon_{a})}{\prod_{i\in Q_0, \; s\neq r}(w_{i,r}-w_{i,s})} u_{i,r}
      \end{multline*}
      \begin{multline*}
          z^* \circ (i_*)^{-1}\left( f(c_1(\Sh{S}_i)) \cap [\Rt_{\omega^*_{i,1}}]\right)
          \\
          =\sum_{r=1}^{\vv_i} f(w_{i,r}-\epsilon_1) \prod_{k: \; i_k=i}(w_{i,r}-a_k - \epsilon_1)\frac{{\prod_{a\in Q_1: \; t(a)=i} \prod_{\substack{s=1 \\ s(a)\neq i \text{ or } s\neq r}}^{\vv_{s(a)}} (w_{i,r}-w_{s(a),s}+\epsilon_{a})}}{\prod_{i\in Q_0, \; s\neq r}(w_{i,s}-w_{i,r})}u^{-1}_{i,r}
      \end{multline*}
      \end{lemma} 
       Following \cite[\S B(ii)]{BFN_quiver}, we set
        \[
        W_i(z)=\prod_{r=1}^{\vv_i}(z-w_{i,r})\qquad W_{i,r}(z)=\prod_{\substack{s=1} \\ s\neq r}^{\vv_i}(z-w_{i,s})\qquad Z_i(z)=\prod_{k: i_k=i}(z-a_k)
        \]
        By the lemma above we get 
        \begin{equation}
            \label{eq: generators e(z) and E(z)}
            z^* \circ (i_*)^{-1}\left(\frac{1}{z-c_1(\Sh{Q}_i)} \cap [\Rt_{\omega_{i,1}}]\right)=(-1)^{\sum_{i\to j} \vv_{j}} E_i(z)
        \end{equation}
        where
        \begin{equation}
          E_i(z) \coloneqq -\sum_{r=1}^{\vv_i}\frac{\prod_{\substack{a\in Q_1
          \\ s(a)=i}} \left(\prod_{t(a)\neq i } W_{t(a)}(w_{i,r}-\epsilon_a)\prod_{t(a)=i} W_{i,r}(w_{i,r}-\epsilon_a)\right)}{(z-w_{i,r})W_{i,r}(w_{i,r})} u_{i,r}
      \end{equation}
      and
      \begin{equation}
          \label{eq: generators f(z) and F(z)}
          z^* \circ (i_*)^{-1}\left(\frac{1}{z-c_1(\Sh{S}_i)} \cap [\Rt_{\omega^*_{i,1}}]\right)=(-1)^{\vv_i - 1} F_{i}(z)
      \end{equation}
      where the sum in the exponent runs over all arrows in the undoubled quiver, and
      \begin{equation}
          F_i(z) \coloneqq \sum_{r=1}^{\vv_i} Z_i(w_{i, r} - \epsilon_1) \times \frac{\prod_{\substack{a\in Q_1
          \\ t(a)=i}} \left(\prod_{s(a)\neq i } W_{s(a)}(w_{i,r}+\epsilon_a)\prod_{s(a)=i} W_{i,r}(w_{i,r}+\epsilon_a)\right)}{(z-w_{i,r}+\epsilon_1)W_{i,r}(w_{i,r})} u_{i,r}^{-1}.
      \end{equation}
      Let $\tilde \phi^+: \Sh{SH}^{T,\nil, \circ} \to \Coul^{T,\vir}_{S,0, \loc}$ be the composition of $\phi^+: \Sh{SH}^{T,\nil,\circ} \to \Coul^{T, \vir}_{G,N}$ with the inclusion \eqref{eq: inclusion Coulomb branch in diff op}. Set also $\gamma^{\edge}_i=\prod_{\substack{a\in \bar{Q}_1\\ s(a)=t(a)}} \epsilon_a$.
      \begin{theorem}
      \label{eq: commutators coulomb}
      \label{eq: generators hall to Coulomb}
          The assignments
          \begin{equation}
              \label{eq: map yangian Coulomb +-}
              e_{i}(z) \mapsto (-1)^{\sum_{i \to j} \vv_j}\gamma^{\edge}_iE_i(z)\qquad f_i(z)\mapsto (-1)^{\vv_i - 1} \gamma^{\edge}_i  F_i(z)
          \end{equation}
          and 
          \begin{equation}
          \label{eq: map yangian Coulomb cartan}
              h_i(z)\mapsto (-1)^{1 + \sum_{i \to j, j \neq i} \vv_j}\gamma^{\edge}_i
              \frac{Z(z)\prod_{\substack{a\in Q_1
          \\ t(a)=i}}  W_{s(a)}(z+\epsilon_a+\epsilon_1) \times\prod_{\substack{a\in Q_1
          \\ s(a)=i}} W_{t(a)}(z-\epsilon_a)}{W_{i}(z+\epsilon_1)W_{i}(z)}
          \end{equation}
          uniquely extend the algebra morphisms \[
          \tilde \phi^+: \Sh{SH}^{T,\nil,\circ} \to \Coul^{T, \vir}_{S,0, \loc} \qquad \tilde \phi^-: \Sh{SH}^{T,\nil, \circ,\opp} \to \Coul^{T, \vir}_{S,0,\loc}
          \]
          to a morphism $\tilde \phi: \Yang_{\mu}^{T,\circ}(\tilde Q, \tilde W)\to \Coul^{T, \vir}_{S,0}$. Moreover, this morphism factors through the inclusion \eqref{eq: inclusion Coulomb branch in diff op}, and hence induces an algebra morphism 
          \[
          \phi: \Yang^T_\mu(\widetilde{Q}, W)\to \Coul^{T, \vir}_{G,N}
          \]
          The restriction of $\phi$ to $\Sh{SH}^{T,\nil,\circ},   \Sh{SH}^{T,\nil,\circ, \opp}\subset\Yang_{\mu}^{T,\circ}(\tilde Q, \tilde W)$ recovers $\phi^+: \Sh{SH}^{T,\nil,\circ}\to \Coul^{T, \vir}_{G,N}$ and $\phi^-: \Sh{SH}^{T,\nil,\circ, \opp}\to \Coul^{T, \vir}_{G,N}$, respectively.
      \end{theorem}
      \begin{proof}
      Fix $i\in Q_0$. The canonical map $\Sh{SH}^{T,\nil,\circ}_{\delta_i}\to \HCoha_{\delta_i}^{T,\onil,\circ}$, which is given by the composition of the inclusion $\Sh{SH}_{\delta_i}^{T,\nil,\circ}\hookrightarrow\HCoha^{T,\onil,\circ}_{\tilde Q, \tilde W, \delta_i}$ and the dimensional reduction isomorphism $\HCoha^{T,\onil,\circ}_{\tilde Q, \tilde W, \delta_i}\cong \HCoha^{T,\onil}_{\Xi_Q, \delta_i}$, sends the Chevalley generator $e_i(z)\in \Sh{SH}^{T,\nil}[\![z^{-1}]\!]$ to 
      \[
      \frac{\gamma^{\edge}_i}{z-c_1(\Sh{L}_i)}\cap [\FM^{T,\onil}_{\delta_i}(\Xi_Q)].
      \]
      where the equivariant weight comes from the nilpotency condition. 
      Therefore, combining Lemma \ref{image generators hall in positive Coulomb}, Lemma \ref{lemma: formulas generators hall in Coulomb} and \eqref{eq: generators e(z) and E(z)}, we deduce that
          \[
          \tilde \phi^+(e_{i}(z))= (-1)^{\sum_{i \to j} \vv_j} \gamma^{\edge}_iE_i(z)
          \]
        Since $\Sh{SH}^{T,\nil}$ is, by definition, spherically generated, it follows that $\tilde \phi^+$ is uniquely determined by this assignment. Similarly, combining Proposition \ref{prop: map opposite hall to Coulomb generators}, Lemma \ref{lemma: formulas generators hall in Coulomb} and \eqref{eq: generators f(z) and F(z)}, we deduce that
          \[
          \tilde \phi^-(f_{i}(z))= (-1)^{\vv_i - 1} \gamma^{\edge}_iF_i(z)
          \]
          and hence, by the same logic, $\tilde \phi^-$ is uniquely determined by this assignment too. This observation concludes the proof of the first statement of the theorem. For the second statement, recall that $\Yang_{\mu}^{T, \circ}(\tilde Q, \tilde W)$ is the quotient of $\Sh{SH}^{T,\nil, \circ} * \Sh{H}^0_Q * \Sh{SH}^{T,\nil, \circ, \opp}$ by the two sided ideal generated by the relations \eqref{eq: shifted Yangian relations}. By the previous considerations, the assignments \eqref{eq: map yangian Coulomb +-} and \eqref{eq: map yangian Coulomb cartan} uniquely determined an algebra morphism
          \[
          \Sh{SH}^{T\nil,\circ} * \Sh{H}^0_Q * \Sh{SH}^{T,\nil, \circ, \opp}\to \Coul^T_{S,0,\loc},
          \]
          so it suffices to check that the relations \eqref{eq: shifted Yangian relations} are in its kernel. We begin with $e_{i, r} f_{j,s} - f_{j, s} e_{i, r} = \delta_{ij} \gamma_i h_{i, r + s}$. Set $\bar h_{i}(z)=\sum_{r>0} h_{i,r} z^{-r}$ the principal part of $h_{i}(z)$. Using $\gamma_i=\epsilon_1\gamma_i^{\edge}$, one can check that this set of relations is equivalent to
      \[
      (z-y)[E_i(z), F_{j}(y)]=\delta_{i,j} \gamma_i(\bar h_i(y)-\bar h_{i}(z)).
      \]
      We now follow the argument from \cite[\S B]{BFN_quiver}.
      It is easy to check that the bracket is zero for $i\neq j$. For $i=j$, a direct computation gives
      \[
      (z-y)[E_i(z), F_{j}(y)]=\sum_{i=1}^r (L_i(y)-R_i(y)) -(L_i(z)-R_i(z)),
      \]
      where 
      \begin{multline*}
          R_i(z)=\frac{Z_i(w_{i,r})\prod_{\substack{a\in Q_1
          \\ t(a)=i}} \left(\prod_{s(a)\neq i } W_{s(a)}(w_{i,r}+\epsilon_a + \epsilon_1)\prod_{s(a)=i} W_{i,r}(w_{i,r}+\epsilon_a + \epsilon_1)\right)
          }{(z-w_{i,r}) W_{i,r}(w_{i,r})W_{i,r}(w_{i,r} + \epsilon_1)}
          \\
          \times \prod_{\substack{a\in Q_1
          \\ s(a)=i}} \left(\prod_{t(a)\neq i } W_{t(a)}(w_{i,r}-\epsilon_a)\prod_{t(a)=i} W_{i,r}(w_{i,r}-\epsilon_a)\right)
      \end{multline*}
      and 
      \begin{multline*}
          L_i(z)=\frac{Z_i(w_{i,r} - \epsilon_1)\prod_{\substack{a\in Q_1
          \\ t(a)=i}} \left(\prod_{s(a)\neq i } W_{s(a)}(w_{i,r}+\epsilon_a)\prod_{s(a)=i} W_{i,r}(w_{i,r}+\epsilon_a)\right)
          }{(z-w_{i,r} + \epsilon_1) W_{i,r}(w_{i,r})W_{i,r}(w_{i,r}-\epsilon_1)}
          \\
          \times \prod_{\substack{a\in Q_1
          \\ s(a)=i}} \left(\prod_{t(a)\neq i } W_{t(a)}(w_{i,r}-\epsilon_a-\epsilon_1)\prod_{t(a)=i} W_{i,r}(w_{i,r}-\epsilon_a-\epsilon_1)\right).
      \end{multline*}
      Therefore, it remains to prove that $L_i(y)-R_i(y)= \gamma_i \bar h_i(y)$. This is a direct generalization\footnote{The only difference are the choices of equivariant weighting and the fact that we allow loops in the quiver, so  a bunch of extra terms appear in the numerators.} of the argument in \cite[Theorem 4.5]{KWWY} and \cite[\S B(iii)]{BFNII}. This proves the first set of relations in \eqref{eq: shifted Yangian relations}. The proof of the second and third line is easier, and it is also spelled out in \emph{loc.cit.} The proof of the second statement is now complete. For the last statement, it suffices to observe that, by construction, $\tilde \phi^{\pm}$ factors through $\Coul^{T, \vir}_{G,N}$. It is also clear that the the assignment of $h_i(z)$ factors through $\Coul^{T, \vir}_{G,N}$ as the coefficients of the right hand side of \eqref{eq: map yangian Coulomb cartan} are polynomial in the variables $w_{i,r}$ and hence elements of $\HO_{G\fps\rtimes T}$. Hence the last statement of the proposition also follows and the proof is complete.
      \end{proof}

      We end this section with the following simple obsrvation. 
      \begin{proposition}
      \label{prop: surj map yangian to Coulomb}
          The localized morphism ${\phi_{\loc}}: \Yang_{\mu, \loc}^{T, \circ}(\tilde Q, \tilde W)\to \Coul^{T, \vir}_{G,N, \loc}$ is surjective.
      \end{proposition}
      \begin{proof}
          By Theorem \ref{eq: generators Coulomb}, we know that the right hand side is generated by the minuscules monopoles $M_{1, \omega_{i,1}}, M_{1, \omega_{i,1}}$ and $\HO_{G{\fps}_T}$. The former are in the image of $\phi_{\loc}$ by Theorem \eqref{eq: generators hall to Coulomb}. The fact that all classes in $\HO_{G\fps}$ are in the image can be seen by inverting \eqref{eq: map yangian Coulomb cartan} to express the Chern polynomial $W_i(z)$ in terms of $h_i(z)$. This is explained in \cite[Thm. B.15]{BFN_quiver}, which is in turn based on \cite[Thm. 4.5]{KWWY}. 
      \end{proof}

\section{Yangians, Coulomb branches and quasimap actions}

\subsection{Statement of the conjecture}

In this section, we discuss the relation between the Hall action and the Coulomb action on $\HO_T(\QM^{\xi}, \varphi_{\QM^{\xi}})$. Observe that, throughout the article, we have obtained the following results:
\begin{enumerate}
\item
\label{item: conj discussion 2}
By Theorem \ref{thm: Coulomb action QM} we have an action
\[
{}^{\phi}\Coul^{T,\vir}_{N,G, \dd'}\otimes\HO_T(\QM^{\xi}_{\dd''}, \varphi_{\QM^{\xi}_{\dd''}}) \to \HO_T(\QM^{\xi}_{\dd'+\dd''}, \varphi_{\QM^{\xi}_{\dd'+\dd''}})
\]
and hence, via Theorem \ref{thm: positive Coulomb in Coulomb}, also an action of the positive subalgebra ${}^{\phi}\Coul^{T,\vir,+}_{G,N}\subset {}^{\phi}\Coul^{T,\vir}_{G,N}$. 
\item
\label{item: conj discussion 3}
In \S\ref{subsec: From hall to Coulomb} we have constructed a morphism $\HCoha^{T, \onil}_{\Xi_{Q}}\to \Coul^{T,c}_{G,N}$ that, via the dimensional reductions $\HCoha^{T, \onil}_{\Xi_{Q}}\cong \HCoha^{T,\onil, \circ}_{\tilde Q, \tilde W}$ and $\Coul^{T,.\vir}_{N,G, \dd'}\cong {}^{\phi}\Coul^{T,\vir}_{G,N}$ from Proposition \ref{prop: dr nilp version} and Theorem \ref{thm: dr Coulomb agebra}, respectively, induce an algebra morphism 
\[
\theta: \HCoha^{T,\onil, \circ}_{\tilde Q, \tilde W}\to {}^\phi \Coul^{T,\vir,+}_{G,N}.
\]
\item
\label{item: conj discussion 1}
On the Hall side, by Theorem \ref{thm: Hall acts on critical quasimaps}, the critical Hall algebra $\HCoha^{T, \circ}_{\tilde Q, \tilde W}$ acts on quasimaps via morphisms\footnote{Because of the sign twist in the definition of the multiplication of $\HCoha^{T, \circ}_{\tilde Q, \tilde W}$, cf. Remark \ref{rem: signed product hall algebras}, to make this action well defined we twist the action map \eqref{eq: hall action quasimaps} by $(-1)^{(\dd''+\vv)^T\Q \dd'}$. This is consistent with the sign conventions in \S\ref{subsec: Dimensionally reduced action hall}.}
\begin{equation}
    \label{eq: hall action JCM section}
    \vact_{\dd',\dd''}: \HCoha^{T, \circ}_{\tilde Q, \tilde W, \dd'}\otimes \HO_T(\QM^{\xi}_{\dd''}, \varphi_{\QM^{\xi}_{\dd''}})\to \HO_T(\QM^{\xi}_{\dd}, \varphi_{\QM^{\xi}_{\dd}})[I_{\dd'}].
\end{equation}
and hence so does the subalgebras 
\[
\Sh{SH}^{T,\nil,\circ}\subseteq \HCoha^{T,\nil, \circ}_{\tilde Q, \tilde W}\subseteq  \HCoha^{T,\onil, \circ}_{\tilde Q, \tilde W}\subseteq \HCoha^{T, \circ}_{\tilde Q, \tilde W},
\]
cf. \S\ref{subsec: Nilpotent and omega-nilpotent 3d CoHA} and Remark \ref{rem: signed product hall algebras}. Moreover, in \S\ref{subsec: Shifted Yangian action}, it is shown that the $\Sh{SH}^{T,\nil, \circ}$-action on $\HO_T(\QM^{\xi}_{\dd''}, \varphi_{\QM^{\xi}_{\dd''}})$ extends to an action of the shifted Yangian $\Yang_{\mu}^{T, \circ}(\tilde Q, \tilde W)$, which contains $\Sh{SH}^{T,\nil, \circ}$ (see \S\ref{subsec: Extension to the shifted Yangian} for the definition of $\Yang_{\mu}^{T,\circ}(\tilde Q, \tilde W)$) as a positive subalgebra.
\end{enumerate}
Informally speaking, the following conjecture states that these three constructions are all compatible.
\begin{conjecture} $ $
\label{conj: the conjecture}

\begin{enumerate}
    \item \label{item 1 conj main body}
    The action of the Hall algebra $\HCoha^{T,\onil, \circ}_{\tilde Q, \tilde W}$ on $\HO_T(\QM^{\xi}, \varphi_{\QM^{\xi}})$ induced by the composition of the algebra morphism $\theta: \HCoha^{T,\onil,\circ}_{\tilde Q, \tilde W}\mapsto \Coul^{T, \vir}_{N,G}$ from \eqref{item: conj discussion 3} with the Coulomb action from Theorem \eqref{item: conj discussion 2} coincides with the Hall action from Theorem \eqref{item: conj discussion 3}. Equivalently, we have a commutative diagram of algebras:
    \begin{equation}
    \label{eq: diagram conjecture}
        \begin{tikzcd}[column sep=4em]
            \HCoha^{T,\onil, \circ}_{\tilde Q, \tilde W} \arrow[r, "\text{\eqref{item: conj discussion 1}}"]\arrow[d, swap, "\theta"] & \End\left(\bigoplus_{\dd\in \BoN^{Q_0}} \HO_T(\QM^{\xi}_{\dd}, \varphi_{\QM^{\xi}_{\dd}})\right)
            \\
            \Coul^{T, \vir}_{G,N}\arrow[ur, swap, "\text{\eqref{item: conj discussion 2}}"]
        \end{tikzcd}
    \end{equation}
    \item \label{item 2 conj main body}
    The restriction of \eqref{eq: diagram conjecture} to the spherical nilpotent subalgebra $\Sh{SH}^{T,\nil,\circ}$ extends to a commutative diagram
    \begin{equation*}
        \begin{tikzcd}
            \Yang_{\mu}^{T,\circ}(\tilde Q, \tilde W) \arrow[r, "\text{\eqref{item: conj discussion 1}}"]\arrow[d] & \End\left(\bigoplus_{\dd\in \BoQ^N} \HO_T(\QM^{\xi}_{\dd}, \varphi_{\QM^{\xi}_{\dd}})\right)
            \\
            {}^\phi\Coul^{T,\vir}_{G,N}\arrow[ur, swap, "\text{\eqref{item: conj discussion 2}}"]
        \end{tikzcd}
    \end{equation*}
    where the labels have the same meaning as above.
\end{enumerate}
\end{conjecture}

We can  prove Conjecture \ref{conj: the conjecture} assuming $\xi\in X$ is polarized, and hence the Hall all the three corners od diagram \eqref{eq: diagram conjecture} can be dimensionally reduced.
\begin{proposition}
\label{prop: proof conj reducible case}
    Assume that $\xi \in X$ is polarized in the sense of \S\ref{subsec: ev point and dim red}. Then Conjecture \eqref{conj: the conjecture} holds. 
\end{proposition}

\begin{remark}

    \label{Rem: intrinsic conj}
Conjecture \ref{eq: diagram conjecture} can be easily lifted at the sheaf level. As discussed in \S\ref{subsec: Shifted Yangian action}, the action of the positive Yangian generators $c_1(\Sh{L}_i)^r\cap [\FM^{T,\nil}_{\delta_i}(\tilde Q)]\in \HO\CoHA^{T,\nil}_{\tilde Q, \tilde W, \delta_i}$
is induced by a fundamental class 
\[
[\overline{\FP}(\dd, \dd + \delta_i)]\in \HO_T(\FM^{\zeta^+\sst}_{\dd}(Q^{\QM}_{\ww,\vv})\times \FM^{\zeta^+\sst}_{\delta_i+\dd}(Q^{\QM}_{\ww,\vv}), \phip{-W_{\dd}}\BoQ^{\vir}\boxtimes\phip{W_{\delta_i+\dd}}\BoQ^{\vir})[I_{\delta_i}]
\]
Similarly, the action of $c_1(\Sh{L}_i)\cap [\FM^{T, \onil}_{\delta_i}(\tilde Q)]\in \HO\CoHA^{T, \onil}_{\tilde Q, \tilde W, \delta_i}$ is induced by a similar class $[\FP(\dd, \dd + \delta_i)]$ living in the same cohomology. Explicitly, this fundamental class can be realized as follows\footnote{The same holds for $[\overline{\FP}(\dd, \dd + \delta_i)]$, with obvious modifications.}.Let $\FM^{\onil, \zeta^+\sst}_{\delta_i,\dd}(Q^{\QM}_{\ww,\vv})$ be the preimage of $\FM^{ \onil}_{\delta_i}(\tilde Q)\subset \FM_{\delta_i}(\tilde Q)$ with respect to the canonical morphism $\FM^{\zeta^+\sst}_{\delta_i,\dd}(Q^{\QM}_{\ww,\vv})\to \FM_{\delta_i}(\tilde Q)$ (see Remark \ref{rem: alternative description action} for the notation). Now consider the canonical morphism
    \[
    \tau:\FM^{\onil, \zeta^+\sst}_{\delta_i,\dd}(Q^{\QM}_{\ww,\vv})\to \FM^{\zeta^+\sst}_{\dd}(Q^{\QM}_{\ww,\vv})\times \FM^{\zeta^+\sst}_{\delta_i+\dd}(Q^{\QM}_{\ww,\vv}) 
    \]
    and denote by $p_{\dd}$ and $p_{\delta_i+\dd}$ the morphisms to the two factors. 
    Applying the vanishing cycle functor $\phip{W_{\delta_i+\dd}-W_{\dd}}$ to the co-unit $\tau_!\tau^!\to \id$, using \ref{eq: base change VC} and \ref{eq: TS}, and taking into account virtual shifts, we obtain a morphism 
    \begin{equation}
    \label{eq: FC generator hall}
        \tau_!\BoQ_{\FM^{\zeta^+\sst}_{\delta_i,\dd}(Q^{\QM}_{\ww,\vv})}[-I_{\delta_i}]\to \phip{-W_{\dd}}\BoQ^{\vir}_{\FM^{\zeta^+\sst}_{\dd}(Q^{\QM}_{\ww,\vv})}\boxtimes \phip{W_{\delta_i+\dd}}\BoQ^{\vir}_{\FM^{\zeta^+\sst}_{\delta_i+\dd}(Q^{\QM}_{\ww,\vv})},
    \end{equation}
    with $I_{\delta_i}$ defined in \eqref{eq: index shifted action}. Since $\FM^{\zeta^+\sst}_{\delta_i,\dd}(Q^{\QM}_{\ww,\vv})$ is smooth, it admits a fundamental class, and its image under the map induced by \eqref{eq: FC generator hall} by taking derived global sections\footnote{This is well defined because $\tau$ is proper representable.}, we get $[\FP(\dd, \dd + \delta_i)]$. 
    
    An alternative but yet equivalent description of the action of $[\FP(\dd, \dd + \delta_i)]$ is given by taking the adjoint of \eqref{eq: FC generator hall} and observing  \cite[\S5.2]{Kin3} that it induced a map
    \begin{align}
    \label{eq: FC generator hall 2}
    p_{\dd}^*\phip{W_{\dd}}\BoQ^{\vir}_{\FM^{\zeta^+\sst}_{\dd}(Q^{\QM}_{\ww,\vv})}\to p_{\delta_i+\dd}^!\phip{W_{\delta_i+\dd}}\BoQ^{\vir}_{\FM^{\zeta^+\sst}_{\dd''}(Q^{\QM}_{\ww,\vv})}[I_{\dd}].
    \end{align}
    induces the action of $[\FM^{T, \onil}_{\tilde Q}(\delta_i)]$ via the usual pull push (as $p_{\delta_i+\dd}$ is proper). More generally, the action of $c_1(\Sh{L}_i)^r\cap[\FM^{T, \onil}_{\tilde Q}]$ is obtained by pre-composing the corresponding morphism \eqref{eq: FC generator hall} with the action of the Chern class $c_1(\Sh{L})^r$. 
    We also remark that \eqref{eq: FC generator hall 2} gives, by dualization, the action of the Yangian lowering operators $f_{i,r}$ too. 
    
    A similar analysis can be run for the Coulomb branch action too. Fix the minuscule cocharacter $\omega_{i,1}$ and recall the definition of $\Mt^{(2),\xi}_{G, N,\omega_{i,1}}$ from \S\ref{subsec: Minuscule monopoles action}. Consider the natural morphism
    \begin{equation}
    \label{eq: sigma map}
        \sigma : \Mt^{(2),\xi}_{G, N,\omega_{i,1}}\to \Mt^{\xi}_{G,N, \dd}\times \Mt^{\xi}_{G,N, \delta_i+\dd}
    \end{equation}
    given by
    \[
    [n_1, g, g', n_\infty, n_\infty^*]\mapsto ([gn_1, g, n_3, n_\infty^*], [n_1, g'g, n_\infty, n_\infty^*])
    \]
    Compared to the notation of \S\ref{subsec: Minuscule monopoles action}, we have added subscripts labelling connected components.
    The projection onto the two factors recover the maps $p_{\omega_{i,1}}$ and  $q_{\omega_{i,1}}$ in \eqref{eq: diagram action minuscules}, which determine the action of the monopole operator $\Rt_{\omega_{i,1}}$, see Proposition \ref{prop: monopole action}. This action can be constructed at the sheaf level as follows. Applying the potential $\phip{-\Sh{W}^{\xi}_{-\dd}+\Sh{W}^{\xi}_{\delta_i+\dd}}$ to the co-unit $\sigma_!\sigma^!\to \id$ and arguing as above, we obtain a morphism 
    \begin{equation}
        \label{eq: FC generator Coulomb}
        \sigma_!\BoQ_{\Mt^{(2),\xi}_{\dd, \dd+\delta_i}}\to \phip{-\Sh{W}^{\xi}_{\dd}} \BoQ^{\vir}_{\Mt^{\xi}_{G,N, \dd}}\boxtimes_{\B T}\phip{\Sh{W}^{\xi}_{\delta_i+\dd}} \BoQ^{\vir}_{\Mt^{\xi}_{G,N, \delta_i+\dd}}[I_{\delta_i}]
    \end{equation}
    and hence a canonical morphism 
    \begin{equation}
        \label{eq: FC generator Coulomb 2}
        q_{\dd}^*\phip{\Sh{W}^{\xi}_{\dd}}\BoQ^{\vir}_{\Mt^{\xi}_{G,N,\dd}}\to q_{\dd}^!\phip{\Sh{W}^{\xi}_{\delta_i+\dd}}\BoQ^{\vir}_{\Mt^{\xi}_{G,N, \delta_i+\dd}}[I_{\delta_i}]
    \end{equation}
    recovering \eqref{eq: monopole action} via pull-push.
    We remark that the left hand side of \eqref{eq: FC generator Coulomb} has no potential because the difference $-\Sh{W}^{\xi}_{\dd}+\Sh{W}^{\xi}_{\dd}$ restricted to $\Mt^{(2),\xi}_{\dd, \dd+\delta_i}$ is zero. This is a simple computation using the definition of $\sigma$. We also stress that the shift $I_{\delta_i}$ can be computed using the same argument as \eqref{prop: Coulomb action def}.  Moreover, the action of the dual generators $\Rt_{\omega^*_{i,1}}$ is obtained by dualizing \eqref{eq: FC generator Coulomb 2}. And the actions of $c_1(\Sh{S})^r\cap \Rt_{\omega^*_{i,1}}$ (resp. $c_1(\Sh{Q})^r\cap \Rt_{\omega_{i,1}}$) is given by pre-composing \eqref{eq: FC generator Coulomb} with the action of the fundamental class $c_1(\Sh{Q})^r$ (resp. $c_1(\Sh{Q})^r$).

    An interesting observation is that all the sheaf operations considered above can be naturally restricted to the support of the DT sheaf, i.e. on $\QM^{\xi}_{\dd}\times \QM^{\xi}_{\delta_i+\dd}$. Specifically, restricting \eqref{eq: FC generator hall 2} and \eqref{eq: FC generator Coulomb 2}—or rather their generalizations encoding the relevant first Chern class action— and applying Corollary \ref{cor: DT are the same}, we obtain canonical morphisms 
    \begin{equation*}
        \rho_!\BoQ_{\Hecke^{\xi}_{\dd, \delta_i+\dd}}[-I_{\delta_i}-2r]\to \varphi_{\QM^{\xi}_{\dd}}\boxtimes \varphi_{\QM^{\xi}_{\delta_i+\dd}}
    \end{equation*}
    where $\rho: \Hecke^{\xi}_{\dd, \delta_i+\dd}\to \QM^{\xi}_{\dd}\times \QM^{\xi}_{\delta_i+\dd}$. To identify the domain with the constant sheaf on the moduli space of Hecke modifications we have used \ref{thm: hecke as crit} in the quiver model and a straightforward computation in the Čech model. Hence, letting $r_{\dd}$ and $r_{\delta_i+\dd}$ be the projections onto the two factors, we get canonical morphism 
    \begin{equation}
        \label{eq: intrinsic action operators}
        r_{\dd}^*\varphi_{\QM^{\xi}_{\dd}}\to r_{\delta_i+\dd}^!\varphi_{\QM^{\xi}_{\delta_i+\dd}}[I_{\delta_i}+2r]
    \end{equation}
    providing an intrinsic definition for the actions of the spherical raising operators $c_1(\Sh{L}_i)\cap [\FM^{T, \onil}_{\tilde Q}(\delta_i)]$ and $c_1(\Sh{Q})^r\cap \Rt_{\omega_{i,1}}$. Applying Verdier duality provides, instead, the lowering operators. 
    \end{remark}
    As a result of our analysis, we conclude that the sheaf theoretic refinement of Conjecture \ref{conj: the conjecture} is the following 
    \begin{conjecture}
    \label{conj 2}
        The intrinsic actions of the operators $c_1(\Sh{L}_i)\cap [\FM^{T, \onil}_{\tilde Q}(\delta_i)]$ and $c_1(\Sh{Q})^r\cap \Rt_{\omega_{i,1}}$, understood as morphisms of the form \eqref{eq: intrinsic action operators}, are equal.
    \end{conjecture}
    Notice in particular that Conjecture \eqref{conj 2} automatically implies the second statement of Conjecture \eqref{conj: the conjecture} as the action of the lowering operators is simply obtained by taking Verdier duals The action of the Cartan generators match by the explicit formulas in Remark \ref{Cartan Coulomb and Hall commute}, Theorem \ref{eq: commutators coulomb} and the proof of Theorem \ref{theorem: shifted Yangian action on critical cohomology}.

\subsection{Proof of Proposition \ref{prop: proof conj reducible case}}

\label{subec: proof conj lagrangian case}

We first illustrate the idea of the proof in the case where we replace $\QM^{\xi}(L)$ with $\Bun^{\leq 0}_{G, \infty}=G\fps\backslash G\lfps^{\leq 0}/G_1[t^{-1}]$. We remark that although the latter is not a moduli space of quasimaps, it serves as an avatar to present the logic of the argument. Recall that $\Gr^+=G\fps\backslash G\lfps^+$. To remove clutter, let $\thGr^{\leq}\coloneqq G\fps^{\leq 0}/G[t^{-1}]$ denote the negative part of the thick affine Grassmannian. Then the two rows of \eqref{eq: raviolo correspondence quasimaps dimred} are equal and get replaced by 
\[
\begin{tikzcd}
&\Gr^+\times \Bun^{\leq 0} _{G, \infty} \arrow[dl, "b"]& \Gr^+\times \thGr^{\leq}\arrow[l, swap, "p"]\arrow[dl] \arrow[r, "q"]  & \Gr^+\times_{G\fps} \thGr^{\leq} \arrow[dlll, bend left, "r"]
\\
\FCoh_{0}(\DD) \times \Bun^{\leq 0} _{G, \infty} & \arrow[l] \FCoh_{0}(\DD)\times \thGr^{\leq} \arrow[l] \end{tikzcd}
\]
where the map $r$ is just the morphism $ \Hecke_{0}\to \FCoh_{0}(\DD)\times \Bun^{\leq 0} _{G, \infty}$ sending a positive Hecke modification of a $G$ bundle to its associated graded. In this case, all the morphisms in the diagram are smooth so the first statement of Conjecture \ref{conj: the conjecture} is equivalent to the statement that $(q^*)^{-1}p^*\circ b^*=r^*$ or, equivalently, $p^*\circ b^*=q^*\circ r^*$. But this simply follows from the commutativity of the diagram and associativity of smooth pullback. 

We now prove the general case, which involves a finer analysis of the various pullback morphisms–as these are not smooth pullbacks. We drop the torus $T$ from the notation. Set $\widehat \Vt= N\fps\times \thGr^{\leq}$ and
\begin{equation}
    \label{eq: def hat B}
    \widehat B_{\dd}=\FCoh_{Q,0}\times_{\FCoh_{0}} \Gr^+_{\dd}\times\widehat \Vt= \FCoh_{Q,0}\times_{\FCoh_{0}} \Gr^+_{\dd}\times N\fps\times \thGr^{\leq}
\end{equation}
Recalling the notation from \S\ref{subsec: dim red Čech model}, consider the commutative diagram 
\begin{equation}
\label{eq: big diagram proof compatibility dimred}
\begin{tikzcd}
    &\Rt^+\times \Nt^{\xi}\arrow[dl, swap, "(b\circ a)\times\id_{\Nt^{\xi}}"] \arrow[dd, hookrightarrow, near start, "\id\times i"]& p^{-1}(\Rt^+\times \thGr^{\leq} ) \arrow[dd, hookrightarrow, "i'"]\arrow[l, swap, "\tilde p"] \arrow[r, "\tilde q"] \arrow[dl, "b'\circ a'"]& q(p^{-1}(\Rt^+\times \thGr^{\leq})\arrow[dd, hookrightarrow]\arrow[llld, bend right, "\tilde r"]
    \\
    \FCoh_{Q,0}\times \Nt^{\xi}\arrow[dd, hookrightarrow, "\id\times i"] &  p_0^{-1}(\FCoh_{Q,0}\times \Nt^{\xi})\arrow[dd, shift left, hookrightarrow, near end, "i''"]\arrow[l, swap, "\tilde p_0"]
    \\
   &\Rt^+\times \Vt\arrow[dl, swap, "(b\circ a)\times\id_{\Vt}"] & \Rt^+\times \thGr^{\leq} \arrow[l, swap, " p"] \arrow[r, "q"] \arrow[dl, "\hat b\circ \hat a"]& \Rt^+\times_{G\fps} \thGr^{\leq} \arrow[llld, bend left, "r"]
   \\
   \FCoh_{Q,0}\times \Vt &  \FCoh_{Q,0}\times \widehat\Vt \arrow[l, swap, "p_0"]
\end{tikzcd}
\end{equation}
The back square of the diagram consists of the two leftmost squares in \eqref{eq: raviolo correspondence quasimaps dimred}; the morphism $b\circ a: \Rt^+\to \FCoh_{Q,0}$ is the composite map \eqref{eq: maps from Coulomb to coh}; the morphisms
    \[
    \begin{tikzcd}
    	    \Rt^+\times \thGr^{\leq} \arrow[r,  "\hat a"] & \widehat B_{\dd}\arrow[r,  "\hat b"] & \FCoh_{Q,0}\times \widehat\Vt 
    \end{tikzcd} 
    \]
    are, respectively, $a\times \id_{\thGr^{\leq}}$ and the natural projection, cf. \eqref{eq: def hat B}. The maps $b'$ and $a'$ are obtained by pullback of $\hat b$, and $\hat a$. The inclusions $i'$ and $i''$ are also obtained by pullback. The morphism $r$ is given by the assignment 
\[
[[n, g], g']\to ((b\circ a)([n,g]), [n,g'])
\]
and the morphism $\tilde r$ is obtained from $r$ by pullback. Notice that we have a canonical isomorphism 
\[
q(p^{-1}(\Rt^+\times \thGr^{\leq})\cong \Hecke^{\xi}(L)
\]
and, under this identification, the map $r$ is the canonical morphism 
\[
q(p^{-1}(\Rt^+\times \thGr^{\leq})\cong \Hecke^{\xi}(L)\to \FCoh_{Q,0}\times \QM^{\xi}(L)\cong \FCoh_{Q,0}\times \Nt^{\xi}
\]
obtained by taking the associated graded of the positive Hecke modification. All vertical squares are Cartesian, as well as 
\begin{equation*}
\label{eq: two diagrams comparison coha and Coulomb aciton dimred}
	\begin{tikzcd}
		 \Rt^+\times \thGr^{\leq} \arrow[r, "q"]\arrow[d, swap, "(b\circ a)\times \id"] & \Rt^+\times_{G\fps} \thGr^{\leq}\arrow[d, "r"]
		 \\
		 \FCoh_{Q,0}\times \widehat\Vt \arrow[r, "p_0"]& \FCoh_{Q,0}\times \Vt 	\end{tikzcd}
	\qquad 
	\begin{tikzcd}
		 p^{-1}(\Rt^+\times \thGr^{\leq}) \arrow[r, "\tilde q"]\arrow[d, swap, "b'\circ a'"] & q(p^{-1}(\Rt^+\times\thGr^{\leq}))\arrow[d, "\tilde r"]
		 \\
		  p_0^{-1}(\FCoh_{Q,0}\times \Nt^{\xi}) \arrow[r, "\tilde p_0"]
 &  \FCoh_{Q,0}\times \Nt^{\xi} 	
 \end{tikzcd}
\end{equation*} 
Recall from \S\ref{subsec: From hall to Coulomb} that the map $a$ fits in a diagram 
\begin{equation}
\label{eq: virtual structure gys pullback theta}
		 \begin{tikzcd} 
		 B_{\dd}\arrow[r, "o"]& E_{\dd}\arrow[d] 
		 \\
			\arrow[u, "a"]\Rt^+_{\dd}\arrow[r, hookrightarrow, "a"]& B_{\dd}		\arrow[u, "s"] (\DD)
		\end{tikzcd}
\end{equation}
where $E_{\dd'}\to B_{\dd'}$ is a finite rank vector bundle, $s$ is a section, and $o$ is the zero section. The morphism $b$ is infinite dimensional smooth and can be approximated by finite type smooth maps, cf. Lemma \ref{eq: lemma factorization pullback hall Coulomb 1}. The structure was used in \S\ref{subsec: From hall to Coulomb} to define $\theta: \HCoha^{T, \onil}_{\Xi_Q}\to \Coul^{\vir, +}_{G,N}$ as the composition of the smooth pullback $b^*$ with the virtual pullback $a^!$. The latter used the data of $E_{\dd'}$ and $s_{\dd'}$. The same structure is inherited, via pullback, by the morphisms  $b'\circ a'$ and $\hat b\circ \hat a$. Therefore, recasting the construction of $\theta$ for these morphisms, we obtain pullback morphisms
\begin{align*}
	(b')^!\circ (a')^*: \HO(p_0^{-1}(\FCoh_{Q,0}\times \Nt^{\xi}), \DD\BoQ) \to \HO(p^{-1}(\Rt^+\times \thGr^{\leq}), \DD\BoQ[-2s_{\dd}])
\end{align*} 
where $s_{\dd}$ is as in \S\ref{rem: virtual Coulomb branch}.
The proof will be the consequence of the following two lemmas.
\begin{lemma}
    The diagram
	\label{lemma 1: dimred compatibility}
	\[
	\begin{tikzcd}
	\HO(\FCoh_{Q,0}, \DD\BoQ)\otimes \HO(\Nt^{\xi},\DD\BoQ)
     \arrow[r, "\theta\otimes id"] \arrow[d, "(\tilde p_0)^*"] & \HO_{G\fps}(\Rt^+, \DD\BoQ[-2s_{\dd}])\otimes \HO(\Nt^{\xi},\DD\BoQ) \arrow[d, "(\tilde p)^!"] 
	\\
	\HO_{G\fps}(p_0^{-1}(\FCoh_{Q,0}\times \Nt^{\xi}), \DD\BoQ)
    \arrow[r, swap, "(b')^!\circ (a')^*"]
    & \HO_{G\fps}(p^{-1}(\Rt^+\times \thGr^{\leq}), \DD\BoQ[-2s_{\dd}]) 
	\end{tikzcd}
	\]

	is commutative.
\end{lemma}

\begin{lemma}
\label{lemma 2: dimred compatibility}
Under the natural identifications $q(p^{-1}(\Rt^+\times \thGr^{\leq})\cong \Hecke^{\xi,+}(L)$, $\FCoh_{Q,0}=\FM(\Xi_Q)$, and $\Nt^{\xi}=\QM^{\xi}(L)=\FM(\Xi_{Q_{\vv,\ww}})$, the morphism $(\bar p\times \bar q)^!$ from \eqref{eq: correspondence 2d coha action dimensionally reduced quasimap} fits in a commutative diagram 
	\begin{equation*}
		\begin{tikzcd}
			\HO_{G\fps}(p_0^{-1}(\FCoh_{Q,0}\times \Nt^{\xi}), \DD\BoQ)
            \arrow[r, "(b')^!\circ (a')^*"]
             & \HO_{G\fps}(p^{-1}(\Rt^+\times \thGr^{\leq}), \DD\BoQ[-2s_{\dd}]) 
                         \\
            \HO(\FCoh_{Q,0},\DD\BoQ)\otimes \HO(\Nt^{\xi}, \DD\BoQ)\arrow[r, "(\bar p\times \bar q)^!"]\arrow[u, "(\tilde p_0)^*"] & \HO(q(p^{-1}(\Rt^+\times \thGr^{\leq}), \DD\BoQ[-2s_{\dd}])\arrow[u, "(\tilde q)^*"]
		\end{tikzcd}
	\end{equation*}
\end{lemma}

Before proving the lemmas, we show how they imply the proposition.
\begin{proof}[Proof of Proposition \ref{prop: proof conj reducible case}]

Combining Lemma \ref{lemma 1: dimred compatibility} and \ref{lemma 2: dimred compatibility}, we get
\begin{align*}
		(\bar p\times \bar q)^!
		&=(\tilde q^*)^{-1} \circ (b')^!\circ (a')^*\circ (\tilde p_0)^*	
		\\
		& = (\tilde q^*)^{-1}\circ (\tilde p)^!\circ (\theta\otimes \id).
\end{align*}
	This equality is equivalent to the commutativity of the top square in the diagram 
	\begin{equation*}
	\begin{tikzcd}
	\HO(\FCoh_{Q,0}, \DD\BoQ)\otimes \HO(\FM(\Xi_{Q_{\vv,\ww}}), \DD\BoQ)
     \arrow[r, "\theta\otimes \id"] \arrow[d, "(\bar p\times \bar q)^!"] & \HO_{G\fps}(\Rt^+, \DD\BoQ[])\otimes \HO(\Nt^{\xi},\DD\BoQ) \arrow[d, "(\tilde q^*)^{-1}\circ (\tilde p)^!"] 
     \\
     \HO(\Hecke^{\xi}(L), \DD\BoQ)\arrow[d, "\bar p'_*"] \arrow[r, equal] & \HO(q(p^{-1}(\Rt^+\times \thGr^{\leq}), \DD\BoQ)\arrow[d, "\tilde m"]
     \\
     \HO(\FM(\Xi_{Q_{\vv,\ww}}), \DD\BoQ) \arrow[r, equal] & \HO(\Nt^{\xi}, \DD\BoQ)
    \end{tikzcd}
	\end{equation*}
	The second square commutes as the maps are literally the same up to the identification of source and target. The leftmost column recovers the action of $\HCoha_{\Xi_Q}$ on $\HO_T(\QM^{\xi}(L))$, while the right column recovers the action of $\Coul^+_{G,N}$. The proposition follows.
\end{proof}


\begin{proof}[Proof of Lemma \ref{lemma 1: dimred compatibility}]
Consider the leftmost bottom square of the diagram \eqref{eq: big diagram proof compatibility dimred}. Expanding it, we get a commutative diagram of the form
\begin{equation}
\label{eq: diagram def p1}
	\begin{tikzcd}
	\FCoh_{Q,0}\times \widehat\Vt\arrow[d, "p_0"] & \widehat B_{\dd'}\arrow[l, swap,  "\hat b"]\arrow[d, "p_1"] & \Rt^+\times \thGr^{\leq} \arrow[d, "p"]\arrow[l, swap,  "\hat a"]
	\\
	 \FCoh_{Q,0}\times \Vt & B_{\dd'}\times \Vt\arrow[l, swap, "b\times\id_{\Vt}"] & \Rt^+\times \Vt\arrow[l, swap, "a\times\id_{\Vt}"]
	\end{tikzcd}
\end{equation}
with right square Cartesian. The first square induces a commutative square
\[
\begin{tikzcd}
(p_0)_*\DD\BoQ_{B_{\dd}\times \widehat \Vt}[2\dim(G\fps)]\arrow[r] & (p_0)_*\hat b_*\DD\BoQ_{\widehat B_{\dd}}[-2\vv^T\dd'-2\infty_{G,N}]
\\
\DD\BoQ_{\FCoh_{Q,0}\times \Vt}\arrow[r]\arrow[u] & (b\times \id_{\Sh{V}})_*\DD\BoQ_{B_{\dd}\times \Vt}[-2\vv^T\dd'-2\dim(N\fps)]\arrow[u]
\end{tikzcd}
\]
Moreover, the second square gives a commutative square 
\[
\begin{tikzcd}
(b\times \id_{\Sh{V}})_*a_* \DD\BoQ_{\widehat B_{\dd'}}\arrow[r] & (b\times \id_{\Sh{V}})_*(p_1)_*\hat a_* \DD\BoQ_{\Rt^+\times \thGr}[2\vv^T\Q\dd']
\\
(b\times \id_{\Sh{V}})_*\DD\BoQ_{B_{\dd'}\times \Vt}\arrow[u] \arrow[r] 	& (b\times \id_{\Sh{V}})_*(p_1)_*\DD\BoQ_{\Rt^+\times \Vt}[2\vv^T\Q\dd']\arrow[u]
\end{tikzcd}
\]
where the horizontal morphisms are constructed as in Appendix \S\ref{app: Refined Gysin pullback} using \eqref{eq: virtual structure gys pullback theta} and its pullback on $\widehat B_{\dd}$. The shift in the last diagram comes from $\rk(E_{\dd'})=2\vv^T\Q\dd'$. The statement of the lemma is now obtained by stacking together the two diagrams, applying the functor $i^!$, base changing the latter with all the pushforwards, and taking derived global sections.
\end{proof}

\begin{proof}[{Proof of Lemma \ref{lemma 2: dimred compatibility}}]
	Consider the morphism $\tilde r$ in the rightmost diagram in \ref{eq: two diagrams comparison coha and Coulomb aciton dimred}. It can be factored as $\tilde r=\tilde b\circ \tilde a$, where $\tilde a$ and $\tilde b$ are uniquely determined by the requirement of fitting in Cartesian squares 
    \[
    \begin{tikzcd}
		 p^{-1}(\Rt^+\times \thGr^{\leq}) \arrow[d, "\tilde q"]\arrow[r, "a'"] & p_1^{-1}(\widehat B_{\dd'})\arrow[r, "b'"]\arrow[d]& p_0^{-1}(\FCoh_{Q,0}\times \Nt^{\xi})\arrow[d, "\tilde p_0"]
		 \\
		  q(p^{-1}(\Rt^+\times\thGr^{\leq})) \arrow[r, "\tilde a"]& \widehat p_1^{-1}(\widehat B_{\dd'})/G\fps \arrow[r, "\tilde b"]&  \FCoh_{Q,0}\times \Nt^{\xi} 	
      \end{tikzcd}
    \]
    Here, $p_1$ is the morphism in \eqref{eq: diagram def p1}.
    From \eqref{eq: def hat B} it follows that $\widehat p_1^{-1}(\widehat B_{\dd})/G\fps$ is the moduli space 
    \[
    (\widehat B_{\dd})/G\fps=\{(\Sh{V}, \Sh{V'},g, s', \alpha) \;|\; g:\Sh{V}\hookrightarrow \Sh{V'}\}/\sim
    \]
    where $\alpha=\{\alpha_a\}_{a\in Q_1}$ is a collection of morphism $\alpha_a: \Sh{V}'_{s(a)}/\Sh{V}_{s(a)}\to \Sh{V}'_{t(a)}/\Sh{V}_{t(a)}$. No compatibility between $n$ and $\alpha$ is required. The map $\tilde b$ sends a tuple $(\Sh{V}, \Sh{V'},g, n', \alpha)$ to $(\alpha, (\Sh{V'},s'))$. The map $\tilde a$ is the closed embedding sending a Hecke modification $(\Sh{V}, \Sh{V'},g, n, n')$ to the tuple $(\Sh{V}, \Sh{V'},g, n', \alpha)$ such that $\alpha$ fits in the exact sequence 
    \[
   \begin{tikzcd}
    0\arrow[r]& \Sh{V}_{s(a)}\arrow[d, "n_a"]\arrow[r] & \Sh{V}'_{s(a)}\arrow[d, "n'_a"]\arrow[r] & \Sh{V}'_{s(a)}/\Sh{V}_{s(a)}\arrow[d, "\alpha_a"]\arrow[r] & 0
   \\
   0\arrow[r]& \Sh{V}_{t(a)}\arrow[r] & \Sh{V}'_{t(a)}\arrow[r] & \Sh{V}'_{t(a)}/\Sh{V}_{t(a)}\arrow[r] & 0
   \end{tikzcd}
    \]
    The same construction of \S\ref{subsec: From hall to Coulomb} exhibits $\tilde a$ as the inclusion of the zero locus $\tilde s^{-1}(0)$ of a section $\tilde s$ of  a vector bundle $\tilde E$ over $ (\widehat B_{\dd})/G\fps$ with fiber $\Hom(\Sh{V}'_{s(a)}, \Sh{V}'_{t(a)}/\Sh{V}_{t(a)})$. Repeating the same argument of the proof of \eqref{lemma 1: dimred compatibility}, one arrives to the commutative diagram
    \begin{equation*}
		\begin{tikzcd}
			\HO_{G\fps}(p_0^{-1}(\FCoh_{Q,0}\times \Nt^{\xi}), \DD\BoQ)
            \arrow[r, "(a')^!\circ (b')^*"]
             & \HO_{G\fps}(p^{-1}(\Rt^+\times \thGr^{\leq}), \DD\BoQ[-2s_{\dd}]) 
                         \\
            \HO(\FCoh_{Q,0},\DD\BoQ)\otimes \HO(\Nt^{\xi}, \DD\BoQ) \arrow[r, "\tilde a^!\circ \tilde b^*"]\arrow[u, "(\tilde p_0)^*"] & \HO(q(p^{-1}(\Rt^+\times \thGr^{\leq}), \DD\BoQ[-2s_{\dd}])\arrow[u, "(\tilde q)^*"]
		\end{tikzcd}
	\end{equation*}
    so it remains to show that $\tilde a^!\circ \tilde b^*=(\bar p\times \bar q)^!$. This can be seen as a consequence that both virtual pullbacks are induced by the same derived structure, but the geometry is simply enough to approach the problem in the following elementary way. To further ease the notation, we set $X=\Hecke^{\xi}(L)$ and $Z=(\widehat B_{\dd'})/G\fps$. The morphism $\tilde a^!$ is induced by the virtual pullback data $(\tilde E, \tilde s)$. The morphism $\bar p\times \bar q$ factors as 
    \[
    X\hookrightarrow B_{\Xi_{Q_{\vv,\ww}}}\xrightarrow{\pi} Z \xrightarrow{\tilde b}\FCoh_{Q,0}\times \Nt^{\xi}.
    \]
    Here
    \[
    B_{\Xi_{Q_{\vv,\ww}}}=\left(\Rep_{\Xi_Q}(\dd')\times \Rep_{\Xi_{Q_{\vv,\ww}}}(\dd'') \times V^{\ext}\times V^{\de}\right)^{\zeta^+\sst}/P 
    \]
    with $V^{\ext}=\bigoplus_{i\in Q_0}\Hom(\BoC^{\dd''_i}\oplus \BoC^{\vv_i}, \BoC^{\dd'})$ and $V^{\de}=\bigoplus_{e\in Q_1}\Hom(\BoC^{\dd''_{s(a)}}, \BoC^{\dd'_{t(a)}})$. Note that $B_{\Xi_{Q_{\vv,\ww}}}$ is a framed semistable version of the space $B_{\Xi_Q}$ defined in in Remark \ref{rem: alternative pullback Coh_Q^0}. The morphism $\pi$ is the vector bundle forgetting the term $V^{\ext}$. Using Remark \ref{rem: alternative pullback Coh_Q^0} and Remark \ref{rem: gysin alternative version}, it follows that $(\bar p\times \bar q)^!=s_{\Xi_{Q_{\ww,\vv}}}^! \pi^* \tilde b^*$, where $s_{\Xi_{Q_{\ww,\vv}}}$ is the virtual pullback associated to the data $(E_{\Xi_{Q_{\ww,\vv}}}, s_{\Xi_{Q_{\ww,\vv}}})$, cf. Remark \ref{rem: alternative pullback Coh_Q^0}. Therefore, it suffices to check that 
    \begin{equation}
        \label{eq: equation to check proof compat dimred}
        s_{\Xi_{Q_{\ww,\vv}}}^! \pi^*=\tilde a^!
    \end{equation}
    Firstly, we claim that we have surjective morphism of vector bundles 
   \begin{equation}
       \label{eq: surj mor bundles proof compatibility actions}
       q: E_{\Xi_{Q_{\ww,\vv}}} \to \pi^* \tilde E\to 0
   \end{equation}
   This easily follows from the observation that that $\tilde E$ is fiberwise given by $\Hom(\Sh{V'_{s(a)}}, \Sh{V'_{t(a)}}/\Sh{V_{t(a)}})$; therefore, the long exact sequence for the functor $\Hom(-,\Sh{V'_{t(a)}}/\Sh{V_{t(a)}})$  applied to the short exact sequence
   \[
       \begin{tikzcd}
           0 \arrow[r] & \Sh{V}_i \arrow[r] & (\BoC^{\dd''_i}\oplus \BoC^{\vv_i}) \otimes \Sh{O}_{\mathbb{P}^1} \arrow[r] & \BoC^{\dd''_i} \otimes \Sh{O}_{\mathbb{P}^1}(1) \arrow[r] & 0
       \end{tikzcd}
    \]
   gives a morphism \eqref{eq: surj mor bundles proof compatibility actions} fiberwise given by 
   \[
   \Hom(\BoC^{\dd''_{s(a)}}\oplus \BoC^{\vv_{s(a)}}, \BoC^{\dd'}_{t(a)})\to \Hom(\Sh{V}'_{s(a)},\Sh{V'_{t(a)}}/\Sh{V_{t(a)}})\to 0.
   \]
   Let now $K$ be the Kernel of \eqref{eq: surj mor bundles proof compatibility actions} and set $Y=(\pi^* \tilde s)^{-1}(0)\subset B_{\Xi_{Q_{\vv,\ww}}}$. We leave to the reader to check that $\pi^* \tilde s =s_{\Xi_{Q_{\ww,\vv}}}\circ q$. It follows that the restriction of $s_{\Xi_{Q_{\ww,\vv}}}$ on $Y$ factors through $K|_Y$, and hence defines a section $s': Y\to K|_Y$. Overall, we have obtained a diagram 
   \[
   \begin{tikzcd}
   (\pi^*\tilde s)^{-1}(0)=Y \arrow[d, "\pi|_Y"] \arrow[r, hookrightarrow]&  B_{\Xi_{Q_{\vv,\ww}}}\arrow[d, "\pi"]
       \\
       \tilde s^{-1}(0)=X\arrow[r, hookrightarrow, "\tilde a"] & Z
   \end{tikzcd}
   \]
   By Lemma \ref{lemma: functoriality Gysin exact sequence bundles} below, we deduce that 
   \[
   s_{\Xi_{Q_{\ww,\vv}}}^! \circ \pi^*=(s')^!\circ (\pi^*\tilde s)^!\circ \pi^*=(s')^!\circ (\pi|_Y)^*\circ \tilde a^!.
   \]
   By construction, the embedding $X=(s')^{-1}(0)\hookrightarrow Y$ is transverse, hence $(s')^!\circ (\pi|_Y)^*=\id$. Therefore, the previous equation recovers \eqref{eq: equation to check proof compat dimred} and hence proves statement \eqref{item 1 conj main body} of Conjecture \ref{conj: the conjecture}. In particular, it proves the first statement of the conjecture. To prove the second statement, it suffices to check that the action of the Yangian generators $e_{i,r}, h_{ir}, f_{i,r}$ factor through $\Coul^{T,\vir}_{G,N}$. The statement for $e_{i,r}$ follows from the first part of the conjecture. The action of $f_{i,r}$ is, up to a sign, the transpose of the action of $e_{i,r}$, and hence their factoring follows from a straightforward  adaptation of the previous argument run one generator at a time. The statement for the Cartan generators $h_{i,r}$ match by the explicit formulas in Remark \ref{Cartan Coulomb and Hall commute}, Theorem \ref{eq: commutators coulomb} and the proof of Theorem \ref{theorem: shifted Yangian action on critical cohomology}.
\end{proof}

\begin{lemma}
\label{lemma: functoriality Gysin exact sequence bundles}
Assume that we are given a short exact sequence $0\to K\to E\xrightarrow{q} \tilde E\to 0$ of vector bundles over $B$ and a section $s: B\to E$. Set $s'=s\circ q$, $X=s^{-1}(0)$ and $Y=s_q^{-1}(0)$. Moreover, the composition $Y\hookrightarrow X\xrightarrow{s} E\xrightarrow{q} E/E'$ is zero so $s$ restricts to a canonical section $s': Y\to K|_Y$. We then have three Cartesian squares 
\[
\begin{tikzcd}
        B\arrow[r, "o"] & E/E'\\
        Y\arrow[u, hookrightarrow] \arrow[r, hookrightarrow] & B\arrow[u, "s_q"]
\end{tikzcd}
\qquad 
\begin{tikzcd}
        B\arrow[r, "o"] & E\\
        X\arrow[u, hookrightarrow] \arrow[r, hookrightarrow] & B\arrow[u, "s"]
\end{tikzcd}
\qquad 
\begin{tikzcd}
        Y\arrow[r, "o"] & E'|_{Y}\\
        X\arrow[u, hookrightarrow] \arrow[r, hookrightarrow] & Y\arrow[u, "s'"]
\end{tikzcd}
\]
and, correspondingly, three refined pullbacks. The following diagram commutes
    \[
    \begin{tikzcd}
        \HO(B, \DD\BoQ_{B})\arrow[r, "s_q^!"]\arrow[rr, bend left, "s^!"]  & \HO(Y, \DD\BoQ_{Y}[2\rk(\tilde E)])\arrow[r, "(s')^!"] & \HO(X, \DD\BoQ_{X}[2\rk(E)]).
    \end{tikzcd}
    \]
\end{lemma}
The proof of this lemma is a standard $6$-functors argument, and hence is omitted. It may be also seen as a special case of Proposition \ref{prop: quasismooth gysin functorial}.

\newpage

\paperpart{Derived symplectic geometry of quasimap moduli} \label{part: quasimap shifted symplectic}

\section{Overview}
In this part of the paper, we will undertake a detailed study of the $(-1)$-shifted symplectic structure on $\QM^\xi(X)$. As discussed in \S\ref{subsec: intro 1} and in more detail in \S\ref{subsec: intrinsic shifted symplectic} below, $\QM^\xi(X)$ has a natural derived enhancement with a canonical $(-1)$-shifted symplectic structure. 

On the other hand, in \S\ref{sec: quiver critical description quasimaps} and \S\ref{subsec: infinite Čech model} we have given \textit{two} different critical locus descriptions of $\QM^\xi(X)$. Each of these has a natural enhancement to a derived critical locus, and derived critical loci carry their \textit{own} $(-1)$-shifted symplectic structures. 

In Appendix \ref{app: derived quasimaps}, we show that the derived structures in the quiver and Čech models, viewed as derived critical loci, agree with the intrinsic derived enhancement of $\QM^\xi(X)$. This is reduced to an analysis of the respective tangent complexes and follows as a corollary of Propositions \ref{prop: quiver tangent complex=quasimap tangent complex}, \ref{prop: Čech tangent complex=quasimap tangent complex}. The purpose of the section is to prove the following two results. 

\begin{theorem} \label{thm: shifted symplectic structure quiver model}
    The intrinsic $(-1)$-shifted symplectic structure on $\QM^\xi(X)$ (described in Corollary \ref{cor: quasimaps is shifted symplectic} below) agrees, in the quiver model, with the canonical $(-1)$-shifted symplectic structure on the derived critical locus of the quasimap potential (described in Lemma \ref{lemma: shifted symplectic form quiver model} below). 
\end{theorem} 

\begin{theorem} \label{thm: shifted symplectic structure Čech model}
   The intrinsic $(-1)$-shifted symplectic structure on $\QM^\xi(X)$ (described in Corollary \ref{cor: quasimaps is shifted symplectic} below) agrees, in the Čech model, with the canonical $(-1)$-shifted symplectic structure on the derived critical locus of the Čech potential (described in Lemma \ref{lemma: shifted symplectic form Čech model}).  
\end{theorem}

The main result of \cite{Br12} assigns a canonical perverse sheaf to any $(-1)$-shifted symplectic stack equipped with an orientation datum. In the case of a derived critical locus with its canonical shifted symplectic structure and orientation data, the construction of \cite{Br12} reproduces the vanishing cycle sheaf of the corresponding function. In Appendix \ref{ap: orientations} we show that the various orientation data on $\QM^\xi(X)$ agree, so we conclude 

\begin{corollary}
\label{cor: DT are the same}
    The vanishing cycle sheaves in the quiver and Čech model are canonically isomorphic, by way of each being canonically isomorphic to the intrinsic DT sheaf of \cite{Br12}:
    \[
    \begin{tikzcd}
        \phip{W^{\xi}} \BoQ^{\vir}_{\FM_{Q^{\QM}}}  & \varphi_{\QM^\xi(X)} \arrow[l, "\sim", swap] \arrow[r, "\sim"] & \phip{\Sh{W}^{\xi}} \BoQ^{\vir}_{\Sh{M}^\xi}.
    \end{tikzcd}
    \]
\end{corollary}

This identification also holds $T$-equivariantly by the naturality of the constructions.

This part of the paper is organized as follows. We introduce our background and conventions/notation for derived algebraic geometry in \S\ref{sec: dag background}. Then in \S\ref{sec: dg quasimap} we explain how to model global functions and global (closed) differential forms in the sense of \cite{pantev2013shiftedsymplecticstructures} via explicit dg-algebras. Finally, \S\ref{sec: compare shifted symplectic structures} contains the proof of Theorems \ref{thm: shifted symplectic structure quiver model}, \ref{thm: shifted symplectic structure Čech model}. The main idea of the proof is simply to unravel the definition of the intrinsic shifted symplectic structure from Corollary \ref{cor: quasimaps is shifted symplectic} in either explicit dg-model from \S\ref{sec: dg quasimap}. 

Section \ref{sec: dag background} has been written for readers who have a background in representation theory and quiver moduli, but assumes no prior knowledge of derived algebraic geometry. Our understanding and presentation of the material owes much to the introductory surveys \cite{calaqueexample}, \cite{pantevvezzosiintro}. Because our description of the quasimap spaces is so explicit, we can reason rather in a rather straightforward way, which is unavailable in the context of more general moduli problems. 

\section{Background} \label{sec: dag background}
\subsection{Sign conventions} We will be dealing extensively with differential graded algebras, for which we use the following sign conventions. A differential graded algebra $A = (A^\bullet, \delta)$ is an associative algebra with a $\mathbb{Z}$-grading $A = \oplus_i A^i$ called cohomological grading, such that the multiplication preserves the grading and is commutative in the following sense: for homogeneous elements $a, b$ of degrees $|a|$, $|b|$,
\[
ab = (-1)^{|a| |b|} ba. 
\]
The differential $\delta$ maps $\delta: A^i \to A^{i + 1}$, satisfies $\delta^2 = 0$, and the graded Leibniz rule with respect to multiplication
\[
\delta(ab) = (\delta a) b + (-1)^{|a|} a (\delta b). 
\]
With this convention, the tensor product of differential graded algebras $(A^\bullet, \delta_A)$, $(B^\bullet, \delta_B)$ may be formed by assigning 
\[
(A \otimes B)^i = \bigoplus_{j +k = i} A^j \otimes B^k.
\]
This is a graded commutative algebra with product defined by $(a \otimes b) (a' \otimes b') = (-1)^{|a'| |b|} aa' \otimes bb'$; the differentials are extended via $\delta_A(a \otimes b) = \delta_A a \otimes b$, $\delta_B(a \otimes b) = (-1)^{|a|} a \otimes \delta _B b$, so that $A \otimes B$ becomes a differential graded algebra with differential $\delta_{A \otimes B} = \delta_A + \delta_B$.

We will have occasion to consider matrix algebras with coefficients in a commutative differential graded algebra; for these, we take into account the sign rule in the definition of the matrix commutator bracket by assigning, for square matrices $a, b$ with entries of homogeneous degrees $|a|, |b|$,
\[
\comm{a}{b} \coloneqq ab + (-1)^{|a| |b| + 1} ba. 
\]
With this assignment, 
\[
\delta \comm{a}{b} = \comm{\delta a}{b} + (-1)^{|a|} \comm{a}{\delta b}. 
\]
To summarize all of the above more succinctly, we apply the Koszul rule of signs \textit{religiously}. Note that, with these conventions, a polynomial algebra on generators of varying cohomological degrees has underlying vector space given by the tensor product of the symmetric algebra on the generators of even degree with the exterior algebra on the generators of odd degree. When we write exterior algebras of graded objects below, a similar convention is understood.

\begin{remark} \label{rmk: complex to space}
    Below, we will recall some foundational definitions in derived algebraic geometry, in which the $\infty$-category of \textit{spaces} $\mathbb{S}$ (in the sense of homotopy theory) will feature. By definition, an object of $\mathbb{S}$ is a specific type of simplicial set called a Kan complex, see \cite{spaceskerodon} and tags referenced therein for details. There is a natural notion of weak homotopy equivalence of Kan complexes, and weakly equivalent Kan complexes are considered equivalent in $\mathbb{S}$. The homotopy theory of $\mathbb{S}$ encodes the homotopy theory of topological spaces; simplicial sets provide a more tractable combinatorial approximation to them.

    Chain complexes $C^\bullet$ determine objects of $\mathbb{S}$ by the following construction\footnote{This is a rational homotopy-theoretic version of one direction of the Dold-Kan correspondence, see \cite{doldkankerodon}, Construction 2.5.6.3.}. Let $\Omega^\bullet(\Delta^n)$ denote the algebraic de Rham complex of the algebraic $n$-simplex $\Delta^n = \text{Spec} \, \mathbb{C}[t_0, \dots, t_n]/(\sum_i t_i - 1)$. Assign the coordinate functions $t_i$ cohomological degree $0$ and their differentials $dt_i$ cohomological degree $+1$. Given a chain complex $C^\bullet$, form the tensor product complex $\Omega^\bullet(\Delta^n) \otimes C^\bullet$ with differential $\delta = \delta_C + d_{\text{dR}}$. Associate to $C^\bullet$ the Kan complex $\mathsf{Map}(\mathbb{C}, C^\bullet)$, which assigns to the $n$-simplex
    \[
    \mathsf{Map}(\mathbb{C}, C^\bullet)(\Delta^n) \coloneqq \text{Hom}_{\text{Cplx}}(\mathbb{C},  \Omega^\bullet(\Delta^n) \otimes C^\bullet)
    \]
    where the right hand side is the usual Hom of complexes and $\mathbb{C}$ is regarded as a chain complex concentrated in degree zero with zero differential. It is instructive to verify that the zero simplices are degree zero cocycles in $C^\bullet$, the one-simplices encode\footnote{Via integration over the simplex and Stokes' theorem.} chain homotopies between zero-simplices, the two-simplices encode homotopies between homotopies, and so on. 

    Most of the objects of $\mathbb{S}$ that we must deal with directly in this paper will arise via this construction, so one can think of the homotopy-theoretic language as a convenient organizing device for the algebraic data encoded in coherent higher chain homotopies.
\end{remark}

\subsection{Basic definitions}
We will quickly define the basic objects of derived algebraic geometry that will be made use of in this paper. This is mainly included to introduce notations, and as a reference point for how our subsequent explicit homological-algebraic constructions fit in with the general theory of \cite{ToenVezzosi}. 

The category of derived affine schemes is declared to be opposite to the category of commutative differential graded algebras (dgas), with zero graded components in degrees\footnote{Grading degree assignments are (of course) purely conventional, but the convention here is chosen to align with the usual degree assignments in the cotangent complex. Specifically, degree zero contains the ``variables'', degree $-1$ represents ``equations/obstructions'', and degrees $\leq -2$ are the higher obstructions. Positive degrees contain information related to automorphisms, thus are absent for derived affine schemes.} $> 0$ (these dgas are called connective). Importantly, these categories are taken to be $\infty$-categories. We will not survey this construction in detail here, but highlight the two most important features. First, two commutative dgas $A$, $B$ are taken to be equivalent if they are quasi-isomorphic as differential graded algebras. Second, Hom-spaces are taken to be \textit{spaces}, i.e. objects of $\mathbb{S}$ as in Remark \ref{rmk: complex to space} above. Denote the $\infty$-category of connective dgas by $\mathbf{cdga}^{\leq 0}$, the $\infty$-category of derived affine schemes by $\mathbf{dAff}$, and the derived affine scheme assigned to a dga $A$ as $\text{Spec} \, A$. 

For readers interested in fine details, the construction \eqref{eq: dgahomspace} below describes dgas as a simplicially enriched category; to obtain from this an $\infty$-category in the sense of \cite{lurie2009higher} one further takes its homotopy coherent nerve, see e.g. \cite{hcnervekerodon}. For a classical reference on the homotopy theory of dgas see \cite{bousfieldgugenheim}. 

Explicitly, the Hom-space between dgas $A$ and $B$ may be described by specifying its $n$-simplices, which is done as follows. Consider $\Omega^\bullet(\Delta^n)$ as in Remark \ref{rmk: complex to space}. The value of the simplicial set $\mathsf{Hom}_{\mathbf{cdga}^{\leq 0}}(A, B)$ on the $n$-simplex is
\begin{equation} \label{eq: dgahomspace}
\mathsf{Hom}_{\mathbf{cdga}^{\leq 0}}(A, B)(\Delta^n) = \text{Hom}_{\mathbf{cdga}^{\leq 0}}(QA, \tau_{\leq 0}(\Omega^\bullet(\Delta^n) \otimes B)).
\end{equation}
Notations are as follows. The Hom on the right hand side is the usual hom of dgas, $\tau_{\leq 0}$ denotes the functor of truncation of a complex to its nonpositive degree part, and $QA \to A$ is a cofibrant replacement\ of commutative dgas. A dga morphism $QA \to A$ is a cofibrant replacement (sometimes also called quasi-free replacement) if $i)$ it is a quasi-isomorphism, and $ii)$ the underlying graded algebra of $QA$ is a polynomial algebra on a finite number of generators of various grading degrees. Informally, $QA$ being cofibrant means ``all the equations have been moved into the differential''. 

There is a functor $t_0: \mathbf{dAff} \to \text{Aff}$ from derived affine schemes to usual affine schemes, called classical truncation, sending 
\[
\text{Spec} \, A \mapsto \text{Spec} \, H^0(A)
\]
where $H^0(A)$ denotes the $0$-th cohomology of $A$, which is an ordinary commutative algebra. $\mathbf{dAff}$ is given a Grothendieck topology\footnote{By definition, a Grothendieck topology on an $\infty$-category is a Grothendieck topology on its homotopy category.} by declaring a collection of maps $\text{Spec} \, A_i \to \text{Spec} \, A$ to be an \'{e}tale covering family if $i)$ each map $A \to A_i$ is derived \'{e}tale (see \cite{pantevvezzosiintro}, Definition 2.23) and $ii)$ $\text{Spec} \, H^0(A_i) \to \text{Spec} \, H^0(A)$ is an \'{e}tale covering family of usual schemes.

A derived stack $X$ is by definition a functor $\mathsf{Map}(-, X): \mathbf{dAff}^{\text{op}} \to \mathbb{S}$ from the opposite category of derived affine schemes to the $\infty$-category of spaces, satisfying the following descent/sheaf condition with respect to the \'{e}tale topology on $\mathbf{dAff}$: given the Čech nerve of an \'{e}tale cover $\text{Spec} \, A_\bullet \to \text{Spec} \, A$, the canonical map
\[
\mathsf{Map}(\text{Spec} \, A, X) \to  \lim \mathsf{Map}(\text{Spec} \, A_\bullet, X)
\]
is an equivalence in $\mathbb{S}$. The limit on the right hand side is taken in the $\infty$-category of spaces, meaning it corresponds to a homotopy limit in the usual sense. Derived stacks themselves form an $\infty$-category called $\mathbf{dStk}$. 

All derived stacks $X$ under consideration in this paper will be Artin, meaning geometric (see \cite{pantevvezzosiintro}, Definition 4.5) and locally of finite presentation. This ensures that $X$ admits a cotangent complex $\mathbb{L}_X$, which is dualizable and whose dual is the tangent complex $\mathbb{T}_X$ (this will be obvious in all examples under consideration).

\subsection{Differential forms on derived Artin stacks}
With the background and context established, we introduce differential forms on derived Artin stacks following \cite{pantev2013shiftedsymplecticstructures}, \cite{pantevvezzosiintro}. 

Recall that given a derived affine scheme $X = \text{Spec} \, A$, its cotangent complex $\mathbb{L}_X$ may be described as follows. Choose a cofibrant replacement $QA \to A$ (the terminology is explained above), then 
\[
\mathbb{L}_X = \Omega^1_{QA} \otimes_{QA} A
\]
where $\Omega^1_{QA}$ denotes the $QA$-module of K\"ahler differentials. Tensoring with $A$ gives an equivalence $\text{QCoh}(X) \simeq QA-\text{mod}$ (where the latter is the $\infty$-category of dg-modules), so that we may identify $\mathbb{L}_X$ with $\Omega^1_{QA}$. K\"ahler differentials $\Omega^1_{QA}$ is given an internal differential $\delta$ defined uniquely by insisting that it is a dg-module over $QA$, and that $d \delta + \delta d = 0$ where $d: QA \to \Omega^1_{QA}$ is the de Rham differential. Then we may define the complex of $n$-shifted $p$-forms on $X$ to be\footnote{The complex does not depend on the choice of $QA$ up to quasi-isomorphism.}
\[
\Sh{A}^p(X; n) \coloneqq (S^p \mathbb{L}_X[-1])[n] = (S^p (\Omega^1_{QA}[-1]))[n] = \Omega^p_{QA}[n - p].
\]
The differential is defined as the unique extension of $\delta$ obeying the graded Leibniz rule. Note that via the construction in Remark \ref{rmk: complex to space}, we may regard this as an object of $\mathbb{S}$. In \cite{pantev2013shiftedsymplecticstructures} the corresponding object is denoted $\Sh{A}^p(X; n)$. We will use the Dold-Kan correspondence to abuse notation and write $\Sh{A}^p(X; n)$ also for the complex; in practice, this creates no trouble because we will always work with the complex.

The complex of $n$-shifted \textit{closed} $p$-forms is described as follows. Note that by construction, 
\[
\bigoplus_{q \geq p} S^q (\mathbb{L}_X[-1])[n]
\]
forms a double complex, with the two differentials being the internal differential $\delta$ and de Rham differential $d$. Pictorially ($\Omega^{p, n}_{QA}$ denotes the cohomological degree $n$ part of $\Omega^p_{QA}$)
\[
\begin{tikzcd}
\Omega^{p, n}_{QA} \arrow[r, "d"] & \Omega^{p + 1, n}_{QA} \arrow[r, "d"] & \dots \\
\Omega^{p, n - 1}_{QA} \arrow[u, "\delta"] \arrow[r, "d"] &  \Omega^{p + 1, n - 1}_{QA} \arrow[u, "\delta"] \arrow[r, "d"] & \dots \\ 
\vdots \arrow[u, "\delta"] & \vdots \arrow[u, "\delta"]
\end{tikzcd}
\]
The cohomological degree in the total complex is defined such that $\Omega^{q, m}_{QA}$ is in degree $q + m - n - p$. Then the complex of $n$-shifted closed $p$-forms is defined as the total complex:
\[
\Sh{A}^{p, cl}(X; n) \coloneqq \text{Tot}^\Pi\Big( \bigoplus_{q \geq p} \Omega^q_{QA}, d + \delta \Big)[n-p]. 
\]
The superscript denotes product totalization, meaning we allow for cocycles to have perhaps infinitely many nonzero components in fixed total cohomological degree. Again, in \cite{pantev2013shiftedsymplecticstructures}, $\Sh{A}^{p, cl}(X; n)$ is regarded as an object in $\mathbb{S}$. An $n$-shifted closed $p$-form $\omega$ is a cocycle in this complex; explicitly, it is an infinite collection $\omega = \{ \omega_i \}_{i \geq 0}$, $\omega_i \in \Omega^{p + i, n - i}_{QA}$ such that 
\begin{equation*}
\begin{split}
    \delta \omega_0 & = 0 \\
    d \omega_i + \delta \omega_{i + 1} & = 0. 
\end{split}
\end{equation*}
The $\omega_i$ for $i \geq 1$ can be viewed as homotopy coherence data witnessing the closure of the underlying $n$-shifted $p$-form $\omega_0$. In this paper, we will often be able to configure things such that most or all of $\omega_i = 0$ for $i \geq 1$. 

This completes the description of forms and closed forms for derived affine schemes. To define them formally for all derived stacks, \cite{pantev2013shiftedsymplecticstructures} show that the functor(s) 
\[
\Sh{A}^{p, (cl)}(-; n): \mathbf{dAff}^{\text{op}} \to \mathbb{S} 
\]
satisfy descent for the \'{e}tale topology, so define derived stacks. Then for any derived stack $X$ we may formally define (where again the rightmost object is a homotopy limit)
\[
\Sh{A}^{p, (cl)}(X; n) \coloneqq \mathsf{Map}_{\mathbf{dStk}}(X, \Sh{A}^{p, (cl)}(-; n)) = \lim_{\text{Spec} \, A \to X} \Sh{A}^{p, (cl)}(\text{Spec} \, A; n). 
\]
The complex of closed zero-shifted zero-forms is called the derived de Rham complex and denoted $\text{DR}^\bullet(X)$. Clearly all the other $n$-shifted $p$-form complexes can be recovered from knowledge of $\text{DR}^\bullet(X)$ by degree shift and truncation.

\begin{remark}
    $\text{DR}^\bullet(X)$ is often regarded as a graded mixed complex (e.g. in \cite{pantev2013shiftedsymplecticstructures}), which differs from a double complex by a regrading. Importantly, the total cohomological grading in $\text{DR}^\bullet(X)$ is the \textit{sum} of the form degree and internal cohomological degree. 
\end{remark}

\subsection{Explicit description for quotient stacks}
We consider the following scenario: $X = [Y/G]$ is a derived global quotient stack, $Y = \text{Spec} \, B$ is an affine derived scheme described by the cofibrant dga $B$, and $G$ is an affine algebraic group fitting into a short exact sequence 
\begin{equation} \label{eq: reductive and unipotent}
\begin{tikzcd}
    1 \arrow[r] & U \arrow[r] & G \arrow[r] & G^{\text{red}} \arrow[r] & 1
\end{tikzcd}
\end{equation}
with $U$ unipotent and $G^{\text{red}}$ reductive, with Lie algebras $\mathfrak{u} = \text{Lie} \, U$, $\mathfrak{g} = \text{Lie} \, G$, $\mathfrak{g}_{\text{red}} = \text{Lie} \, G^{\text{red}}$. Differential forms have the following explicit description in this situation. Since $X$ is Artin, by Proposition 1.14 of \cite{pantev2013shiftedsymplecticstructures}, ignoring the differential we have 
\[
\text{DR}^\bullet(X) \simeq  \Gamma(X, S^\bullet(\mathbb{L}_X[-1]))
\]
where the right hand side is the complex of derived global sections. Note that pullback gives an equivalence $\text{QCoh}(X) \simeq \text{QCoh}_G(Y)$ between quasicoherent sheaves on $X$ and $G$-equivariant quasicoherent sheaves on $Y$, under which we may identify $\mathbb{L}_X$ as
\[
\mathbb{L}_X \simeq  \text{hofib} \Big( \mathbb{L}_Y \to \mathfrak{g}^* \otimes \Sh{O}_Y \Big) 
\]
where the morphism is the dual of the infinitesimal action map. We may also factor the pushforward $X \to \text{pt}$ as $X \to \B G \to \text{pt}$, arriving at (again, for now ignoring differentials)
\begin{equation*}
\begin{split}
   \text{DR}^\bullet([Y/G]) & \simeq \Gamma(\B G, \Gamma(Y, S^\bullet(\mathbb{L}_Y[-1] \oplus \mathfrak{g}^*[-2]))) \\
   & \simeq \Gamma(\B G, \Omega^\bullet_B \otimes S^\bullet(\mathfrak{g}^*[-2])).
\end{split}
\end{equation*}
It remains to describe taking global sections along $\B G$. Essentially by definition, $\Gamma(\B G, -)$ is the derived functor of $G$-invariants, and by the short exact sequence \eqref{eq: reductive and unipotent} we have a natural equivalence of functors 
\[
\Gamma(\B G, -)  \simeq \Gamma(BU, -)^{G^{\text{red}}}
\]
since taking invariants for a reductive group is an exact functor. Now since $U$ is unipotent, for any complex of $U$-modules $(V, \delta_V)$ we may model $\Gamma(BU, V)$ by a Chevalley-Eilenberg complex (e.g. by Hochschild-Mostow spectral sequence, arguing as in \cite{feiginfuchs} Chapter 1 \S3, Corollary to Theorem 6)
\[
\Gamma(BU, V) \simeq \text{CE}(\mathfrak{u}, V) \coloneqq ( V \otimes S^\bullet (\mathfrak{u}^*[-1]), \delta_V + \delta_{\text{CE}})
\]
with the standard Chevalley-Eilenberg differential $\delta_{\text{CE}}$. We arrive at 
\[
\text{DR}^\bullet([Y/G]) \simeq (\Omega^\bullet_B \otimes S^\bullet (\mathfrak{u^*}[-1]) \otimes S^\bullet(\mathfrak{u}^*[-2]) \otimes S^\bullet(\mathfrak{g}^*_{\text{red}}[-2]))^{G^{\text{red}}}. 
\]
It remains to describe the internal and de Rham differentials on $\text{DR}^\bullet([Y/G])$. First observe that we may view 
\[
\Omega^\bullet_B \otimes S^\bullet(\mathfrak{u}^*[-1]) \otimes S^\bullet(\mathfrak{u}^*[-2]) \simeq \Omega^\bullet_{\text{CE}(\mathfrak{u}, B)}
\]
as the algebraic differential forms on the Chevalley-Eilenberg dga 
\[
\text{CE}(\mathfrak{u}, B) = (B \otimes S^\bullet(\mathfrak{u}^*[-1]), \delta_B + \delta_{\text{CE}}). 
\]
In this description, both the de Rham differential $d$ and internal differential $\delta$ on $\Omega^\bullet_{\text{CE}(\mathfrak{u}, B)}$ are uniquely characterized as in the discussion above. Both of these are extended to $\text{DR}^\bullet([Y/G])$ to act trivially on the $S^\bullet(\mathfrak{g}^*_{\text{red}}[-2])$ factor. The total differential on $\text{DR}^\bullet([Y/G])$ receives only one additional contribution, which may be understood as follows. 

Because $G^{\text{red}}$ acts on $Y$ and $U$, there is an induced action of $\mathfrak{g}_{\text{red}}$ on $\text{CE}(\mathfrak{u}, B)$ by derivations. We may regard the action as a canonical element $$V(-) \in \text{Der}(\text{CE}(\mathfrak{u}, B)) \otimes \mathfrak{g}^*_{\text{red}}.$$ 

Identify $S^\bullet(\mathfrak{g}_{\text{red}}^*[-2])$ with polynomial functions on $\mathfrak{g}_{\text{red}}[2] \ni \phi^{\text{red}}$, and view $V(\phi^{\text{red}})$ as a linear function on $\mathfrak{g}_{\text{red}}[2]$ valued in $\text{Der}(\text{CE}(\mathfrak{u}, B))$. Denote by $\iota_{V(\phi^{\text{red}})}$ the operator of contraction with this vector field. $\iota_{V(\phi^{\text{red}})}$ acts on $\Omega^\bullet_{\text{CE}(\mathfrak{u}, B)} \otimes S^\bullet(\mathfrak{g}^*_{\text{red}}[-2])$ and satisfies $\iota_{V(\phi^{\text{red}})} \delta + \delta \iota_{V(\phi^{\text{red}})} = 0$ as a formal consequence of the $G^{\text{red}}$-equivariance of the differential on $\Omega^\bullet_{\text{CE}(\mathfrak{u}, B)}$ and our sign conventions. Then we finally arrive at the following complete ``Cartan-type" model of the derived de Rham complex: 
\begin{proposition} \label{prop: de rham on quotient}
For $Y$ an affine derived scheme modeled by a cofibrant dga $B$ and $G$ an affine algebraic group fitting into an exact sequence of the form \eqref{eq: reductive and unipotent}, the de Rham complex of $[Y/G]$ may be modeled as
\begin{equation}
    \textnormal{DR}^\bullet([Y/G]) \simeq \Big ( (\Omega^\bullet_{\textnormal{CE}(\mathfrak{u}, B)} \otimes S^\bullet(\mathfrak{g}^*_{\textnormal{red}}[-2]))^{G^{\textnormal{red}}}, \, d + \iota_{V(\phi^{\textnormal{red}})} + \delta \Big). 
\end{equation}
The differential $\delta$ is the sum of the differential $\delta_B$ on $B$ and the Chevalley-Eilenberg differential $\delta_{\textnormal{CE}}$, uniquely extended to $\Omega^\bullet_{\textnormal{CE}(\mathfrak{u}, B)}$, and acts trivially on the $S^\bullet(\mathfrak{g}^*_{\textnormal{red}}[-2])$ factor. $d$ is the de Rham differential on $\Omega^\bullet_{\textnormal{CE}(\mathfrak{u}, B)}$, also extended to act trivially on $S^\bullet(\mathfrak{g}^*_{\textnormal{red}}[-2])$. 
\end{proposition}
The total differential indeed squares to zero on $G^{\text{red}}$-invariants by Cartan's magic formula $\acomm{d}{\iota_{V(\phi^{\text{red}})}} = \mathscr{L}_{V(\phi^{\text{red}})}$ (note that both $d$ and $\iota_V$ anticommute with $\delta$, the former by construction and the latter as explained above). The sum $\iota_{V{(\phi^{\text{red}})}} + \delta$ should be viewed as the internal differential on $\text{DR}^\bullet([Y/G])$, while of course $d$ is the de Rham differential. 

\begin{remark}
    A more geometric way to understand Proposition \ref{prop: de rham on quotient} is as follows. By the exact sequence \eqref{eq: reductive and unipotent}, we have an equivalence $[Y/G] = [Z/G^{\text{red}}]$ where $Z = [Y/U]$. Now $[Y/U]$ is an affine stack (see e.g. \cite{toen2006affinestackschampsaffines}, \cite{mathewmondal} Corollary 4.4 and Proposition 4.11). Moreover $\Gamma([Y/U], \Sh{O}_{[Y/U]}) \simeq \text{CE}(\mathfrak{u}, B)$ as we have computed above, and similarly $\text{DR}^\bullet([Y/U]) \simeq \Omega^\bullet_{\text{CE}(\mathfrak{u}, B)}$. To pass further to the quotient by $G^{\text{red}}$, we then just employ the standard Cartan model for reductive quotients. The affine stack point of view will be useful to understand certain constructions below. 
\end{remark}

\subsection{Symplectic and Lagrangian structures}
The final necessary piece of preparatory material is the definition of shifted symplectic and Lagrangian structures. Notice that there is a canonical map 
\[
\Sh{A}^{2, cl}(X; n) \to \Sh{A}^2(X; n)
\]
called the underlying two-form (in terms of $\{ \omega_i \}_{i \geq 0}$, it just projects to $\omega_0$). Any $n$-shifted two-form on an Artin stack $X$ induces a map 
\[
\mathbb{T}_X \to \mathbb{L}_X[n]
\]
and is called \textit{nondegenerate} if this is a quasi-isomorphism (more precisely, equivalence in the $\infty$-category $\text{QCoh}(X)$). An an $n$-shifted closed $2$-form $\omega_X$ whose underlying two-form is nondegenerate is called a shifted symplectic structure on $X$. Notice that the underlying two form does not determine the shifted symplectic structure uniquely, however the underlying two-form together with the homotopy coherence data $\{ \omega_i \}_{i \geq 1}$ does by definition. 

If $X$ has a shifted symplectic structure $\omega$, there is also a notion of Lagrangian structure on a morphism $f: L \to X$ which is the following. An $\textit{isotropic structure}$ on $f$ is a homotopy $\vartheta$ from $f^* \omega$ to zero in the space $\Sh{A}^{2, cl}(L; n)$ of $n$-shifted closed $2$-forms on $L$. Relative to any complex models for $\text{DR}^\bullet(X)$ and $\text{DR}^\bullet(L)$, this is concretely a collection $\{ \vartheta_i \}_{i \geq 0}$, satisfying the conditions
\begin{equation*}
\begin{split}
    f^* \omega_0 & = \delta \vartheta_0 \\ 
    f^* \omega_i & = \delta \vartheta_i + d \vartheta_{i - 1}, \, \, \, i \geq 1.
\end{split}
\end{equation*}
A convenient way to reformulate this is as follows (we follow \cite{calaquelagrangianmap}, Remark 2.1). $f$ induces a canonical map $df: \mathbb{T}_L \to f^*\mathbb{T}_X$, with a dual $df^\vee: f^* \mathbb{L}_X \to \mathbb{L}_L$. We also have the pullback $f^*\mathbb{T}_X \to f^*\mathbb{L}_X[n]$ of the morphism determined by the underlying two-form $\omega_0$ on $X$, and the composite $df^\vee[n] \circ \omega_0 \circ df$ is the pullback $f^*\omega_0$ of $\omega_0$ regarded as a differential $2$-form. Then the leading component $\vartheta_0$ of the isotropic structure provides data witnessing the homotopy commutativity of the diagram
\[
\begin{tikzcd}
    \mathbb{T}_L \arrow[r] \arrow[d, "df"] & 0 \arrow[d] \\ 
    f^* \mathbb{T}_X \arrow[r, "df^\vee \circ \omega_0", swap] & \mathbb{L}_L[n]
\end{tikzcd}
\]
in $\text{QCoh}(L)$. A pair $(f, \vartheta)$ is called $\textit{nondegenerate}$ if this diagram is moreover a homotopy pullback, that is to say the canonical map 
\[
\begin{tikzcd}
\mathbb{T}_L \arrow[r] & \text{hofib} \Big( f^* \mathbb{T}_X \arrow[r, "df^\vee \circ \omega_0"] & \mathbb{L}_L[n] \Big) 
\end{tikzcd}
\]
is an equivalence. A Lagrangian structure on a morphism $f: L \to X$ to an $n$-shifted symplectic derived Artin stack is a nondegenerate isotropic structure. 

Given two morphisms $f_i: L_i \to X$ to an $n$-shifted symplectic stack $X$, together with Lagrangian structures on each $f_i$, form the homotopy pullback/derived fiber product
\[
\begin{tikzcd}
    L_1 \times_X L_2 \arrow[r, "\pi_2"] \arrow[d, "\pi_1"] & L_2 \arrow[d, "f_2"] \\ L_1 \arrow[r, "f_1"] & X. 
\end{tikzcd}
\]
We have the following key result. 
\begin{theorem}{{\cite[Theorem 2.9]{pantev2013shiftedsymplecticstructures}}}
\label{thm: PTVV lag intersection}
Notations as above, $L_1 \times_X L_2$ carries a canonical $(n - 1)$-shifted symplectic structure.
\end{theorem}
The proof is by explicit construction of the shifted symplectic structure given the isotropy data, and we recall the salient points as we will need this construction below. By hypothesis, there is a homotopy $f_2 \circ \pi_2 \sim f_1 \circ \pi_1$ which induces a homotopy $\sigma = \{ \sigma_i \}_{i \geq 0}$ in the space $\Sh{A}^{2, cl}(L_1 \times_X L_2, n)$:
\begin{equation*}
\begin{split}
    \pi_1^* f_1^* \omega_0 - \pi_2^* f_2^* \omega_0 & = \delta \sigma_0 \\
    \pi_1^* f_1^* \omega_i - \pi_2^* f_2^* \omega_i & = \delta \sigma_i + d\sigma_{i - 1}.
\end{split}
\end{equation*}
On the other hand, using the Lagrangian structures we have $f_1^* \omega_i = \delta \vartheta^{(1)}_i + d\vartheta^{(1)}_{i - 1}$, $f_2^* \omega_i = \delta \vartheta^{(2)}_i + d \vartheta^{(2)}_{i - 1}$, so the above equation becomes 
\[
\delta(\pi_1^* \vartheta^{(1)}_i - \pi_2^* \vartheta^{(2)}_i - \sigma_i) + d(\pi_1^* \vartheta^{(1)}_{i - 1} - \pi_2^* \vartheta^{(2)}_{i - 1} - \sigma_{i - 1}) = 0
\]
but this precisely says that the collection 
\[
\Omega_i = \pi_1^* \vartheta^{(1)}_i - \pi_2^*\vartheta^{(2)}_i - \sigma_i
\]
defines a cocycle $\Omega$ in $\Sh{A}^{2, cl}(L_1 \times_X L_2; n - 1)$, in other words an $(n - 1)$-shifted closed two-form. \cite{pantev2013shiftedsymplecticstructures} check that the nondegeneracy of $\vartheta^{(1)}_0$, $\vartheta^{(2)}_0$ implies the nondegeneracy of $\Omega_0$, so that $\Omega$ is an $(n-1)$-shifted symplectic structure. 

\subsection{Example: derived symplectic reduction} \label{subsect: dga for X}
To illustrate the above procedure for constructing shifted symplectic structures, we will describe explicitly the symplectic structure on the stacky quiver variety $\FX = [\mu^{-1}(0)/G]$ by viewing it as a Lagrangian intersection. This example goes back to \cite{calaquelagrangianmap}, \cite{Safronov_2016}. This example will also be technically useful in our analysis below.

As usual we consider $T^*\text{Rep}_Q(\vv, \ww)$ with the canonical symplectic form. We consider $\FX = [\mu^{-1}(0)/G]$ as the quotient of the derived zero locus of $\mu$. $\FX$ satisfies the assumptions of Proposition \ref{prop: de rham on quotient}; before spelling out $\text{DR}^\bullet(\FX)$, let us spell out $\mu^{-1}(0)$ as an affine derived scheme. 

By definition, $\mu^{-1}(0)$ fits into a homotopy pullback diagram of affine derived schemes:
\begin{equation} \label{eq: pullback for moment map}
\begin{tikzcd}
    \mu^{-1}(0) \arrow[r] \arrow[d] & T^*\text{Rep}_Q(\vv, \ww) \arrow[d, "\mu"] \\
    0 \arrow[r, hook] & \mathfrak{g}^*_{\vv}.
\end{tikzcd}
\end{equation}
This means it is equivalent to any dga $$\Sh{O}_{\mu^{-1}(0)} = \Sh{O}_0 \otimes^L_{\Sh{O}_{\mathfrak{g}^*_{\vv}}} \Sh{O}_{T^* \text{Rep}_Q(\vv, \ww)}$$ modeling the derived tensor product of function dgas (here, as all other schemes entering the fiber product diagram are underived, we regard their coordinate rings as dgas with zero differential concentrated entirely in degree zero).

Given dgas and dga morphisms $C \to A $, $C \to B$, the derived tensor product $A \otimes_C^L B$ may be modeled via $\widetilde{A} \otimes_C B$, where $C \to \widetilde{A}$ is any cofibrant replacement of $C \to A$. By definition, this means the following: there is a factorization $C \to \widetilde{A} \to A$ such that $\widetilde{A} \to A$ is a quasi-isomorphism, and such that, as a graded algebra $\widetilde{A} \simeq \text{Sym}_C(V)$ for some nonpositively graded vector space $V$. 

Let us write $\Sh{O}_{\mathfrak{g}^*_{\vv}} \simeq \mathbb{C}[x]$ for some $\mathfrak{g}^*_{\vv}$-valued variable $x$. In this presentation the morphism $\Sh{O}_{\mathfrak{g}^*_{\vv}} \to \Sh{O}_0$ is just the quotient map $\mathbb{C}[x] \to \mathbb{C}$ setting $x = 0$. A cofibrant replacement for this map is given by the first arrow in the factorization
\[
\begin{tikzcd}
    \mathbb{C}[x] \arrow[r, hook] & (\mathbb{C}[x, \chi], \delta) \arrow[r] & \mathbb{C}
\end{tikzcd}
\]
The middle term denotes the symmetric algebra (again, in the graded sense---as a vector space, $\mathbb{C}[x, \chi] \simeq S^\bullet(\mathfrak{g}_{\vv}) \otimes \Lambda^\bullet(\mathfrak{g}_{\vv})$) generated by $x$ in degree zero and another $\mathfrak{g}^*_{\vv}$-valued generator $\chi$ in degree $-1$. The differential $\delta$ is characterized by $\delta \chi = x$, $\delta x = 0$. The projection map setting $x$ and $\chi$ to zero is easily verified to be a quasi-isomorphism. Write $\widetilde{\Sh{O}}_0 \coloneqq (\mathbb{C}[x, \chi], \delta)$, then we arrive at the dga model 
\begin{equation} \label{eq: model for zero locus mu}
\Sh{O}_{\mu^{-1}(0)} = \widetilde{\Sh{O}}_0 \otimes_{\Sh{O}_{\mathfrak{g}^*_{\vv}}} \Sh{O}_{T^*\text{Rep}_Q(\vv, \ww)}. 
\end{equation}
To make this completely explicit, notice that $\Sh{O}_{T^* \text{Rep}_Q(\vv, \ww)} \simeq \mathbb{C}[\xi_2, \xi_3, \iota, \jmath]$, that is to say polynomial functions of the matrix coefficients of the quiver data $(\xi_2, \xi_3, \iota, \jmath)$. Then as a dga, $\Sh{O}_{\mu^{-1}(0)}$ is fully characterized by the following: it has generators given by (all matrix coefficients of) $\xi_2, \xi_3, \iota, \jmath$ in degree zero, a $\mathfrak{g}^*_{\vv}$-valued generator $\chi$ in degree $-1$, and differential characterized by 
\begin{equation*}
\begin{split}
    \delta \chi & = \comm{\xi_2}{\xi_3} + \iota \jmath \\ 
    \delta \xi_2 = \delta \xi_3 = \delta \iota = \delta \jmath & = 0. 
\end{split}
\end{equation*}
We write $\Sh{O}_{\mu^{-1}(0)} \simeq (\Sh{O}_{T^* \text{Rep}_Q(\vv, \ww)} \otimes S^\bullet(\mathfrak{g}_{\vv}[1]), \delta)$ with $\delta$ as above---notice this is just the standard Koszul complex for $\mu = 0$, the only purpose of the above discussion was to explain how it arises from the general formalism. 

Since the group $G_{\vv} = \prod_{i \in Q_0} GL(v_i)$ is reductive, Proposition \ref{prop: de rham on quotient} with $\mathfrak{u} = 0$ applies to give a Cartan model for $\text{DR}^\bullet(\FX)$: 
\[
\text{DR}^\bullet(\FX) \simeq ( (\Omega^\bullet_{\mu^{-1}(0)} \otimes S^\bullet(\mathfrak{g}^*_{\vv}[-2]))^{G_{\vv}}, d + \iota_{V(\phi)} + \delta)
\]
where we have written $G = G^{\text{red}} = G_{\vv}$ and dropped the ``red'' notation because the entire group is reductive. 

We expect that $\FX$ is a zero-shifted symplectic stack, and the most obvious guess for its symplectic form is 
\[
\omega'_{\FX} = \tr(d \xi_2 \wedge d \xi_3) + \tr (d \iota \wedge d \jmath)
\]
(and all $\omega_i = 0$ for $i \geq 1$). This is clearly $d$ and $\delta$-closed, but it is $\iota_{V(\phi)} + \delta$ that is the internal differential on $\text{DR}^\bullet(\FX)$, and $\omega'_{\FX}$ is \textit{not} closed for $\iota_{V(\phi)} + \delta$:
\[
(\iota_{V(\phi)} + \delta)\omega'_{\FX} = \iota_{V(\phi)} \omega'_{\FX} = -d( \tr \phi(\comm{\xi_2}{\xi_3} + \iota \jmath)). 
\]
However, its failure to be closed is $\delta$-exact, $\tr \phi(\comm{\xi_2}{\xi_3} + \iota \jmath) = \delta ( \tr \phi \chi )$, moreover we have $\delta( \tr \phi d \chi) = (\iota_{V(\phi)} + \delta)(\tr \phi d \chi)$ because $\iota_{V(\phi)} d\chi = \comm{\chi}{\phi}$ and  $\comm{\phi}{\phi} = 0$. Then we see that 
\begin{equation*}
\begin{split}
\omega_{\FX} & = \omega'_{\FX} - \tr \phi d \chi \\
& = \tr( d\xi_2 \wedge d \xi_3) + \tr( d \iota \wedge d \jmath) - \tr \phi d\chi
\end{split}
\end{equation*}
defines a cocycle in $\text{DR}^\bullet(\FX)$, and it is $\omega_{\FX}$ that represents the symplectic form in this model.

Let us now see how this ``correction term'' is incoporated in the general formalism. The definining pullback diagram \eqref{eq: pullback for moment map} is $G_{\vv}$-equivariant, thus determines (and is in fact equivalent to) a homotopy pullback diagram of derived quotient stacks
\begin{equation}
\begin{tikzcd}
    \left[ \mu^{-1}(0)/G_{\vv} \right] \arrow[r] \arrow[d] & \left[ T^* \text{Rep}_Q(\vv, \ww)/G_{\vv} \right] \arrow[d, "\mu"] \\ 
    \B G_{\vv} \arrow[r, hook] & \left[ \mathfrak{g}^*_{\vv}/G_{\vv} \right] 
\end{tikzcd}
\end{equation}
where we recognize the classifying stack $\B G_{\vv} \simeq [0/G_{\vv}]$. We would now like to interpret this as a fiber product of Lagrangian morphisms. First notice that $[\mathfrak{g}^*_{\vv}/G_{\vv}]$ is naturally 1-shifted symplectic; for example, it may be realized as $T^*[1]\B G_{\vv}$, and shifted cotangent stacks are always shifted symplectic \cite{Calaque_2019}. More concretely, using our model $\Sh{O}_{\mathfrak{g}^*_{\vv}} = \mathbb{C}[x]$ and the corresponding Cartan model for $\text{DR}^\bullet([\mathfrak{g}^*_{\vv}/G_{\vv}])$ from Proposition \ref{prop: de rham on quotient}, the 1-shifted symplectic form is 
\[
\omega_{T^*[1]\B G_{\vv}} = \tr \phi dx. 
\]
This is indeed closed because $\iota_{V(\phi)} \tr \phi dx = \tr \phi \comm{x}{\phi} = 0$ (the higher closure data $\omega_i$, $i \geq 1$, all may be taken to vanish in this model), and is clearly nondegenerate. Next we verify that the morphisms $f_1 : \B G_{\vv} \to [\Fg^*_{\vv}/G_{\vv}]$ and $f_2 = \mu$ admit canonical Lagrangian structures. We will present them relative to the model \eqref{eq: model for zero locus mu} for $\Sh{O}_{\mu^{-1}(0)}$ above, in other words model $\B G_{\vv}$ as $\text{Spec} \, \widetilde{\Sh{O}}_0/G_{\vv}$. Then we have 
\[
f_1^* \omega_{T^*[1] \B G_{\vv}} = \tr \phi d(\delta \chi)= (\iota_{V(\phi)} + \delta) \vartheta^{(1)}_0
\]
where $\vartheta^{(1)}_0 = - \tr \phi d \chi$, and the higher coherence data vanishes. Similarly, 
\[
\begin{split}
f_2^* \omega_{T^*[1] \B G_{\vv}} & = \tr \phi \mu^*dx \\ 
& = \tr \phi d(\comm{\xi_2}{\xi_3} + \iota \jmath) \\ 
& = \iota_{V(\phi)} \vartheta^{(2)}_0
\end{split}
\]
with $\vartheta^{(2)}_0 = -\omega'_{\FX} = - \tr(d \xi_2 \wedge d \xi_3) - \tr (d\iota \wedge d\jmath)$, and $\vartheta^{(2)}_i = 0$ when $i \geq 1$ (note in this model the internal differential $\delta$ on $\Sh{O}_{T^* \text{Rep}_Q(\vv, \ww)}$ vanishes). It is easily verified that both of these isotropic structures are nondegenerate, so that each of $f_1, f_2$ have $1$-shifted Lagrangian structures. 

In the notation above the previous section, clearly we may take $\sigma_i = 0$ (this is because we are using a strict model for $\Sh{O}_{\mu^{-1}(0)}$, where it is defined by a strictly commutative pushout diagram). Then the canonical 0-shifted symplectic structure on $\FX$ viewed as a 1-shifted Lagrangian intersection is given by 
\[
\omega_{\FX} = \pi_1^* \vartheta^{(1)}_0 - \pi_2^* \vartheta^{(2)}_0 - \sigma_0 = \omega'_{\FX} - \tr \phi d \chi
\]
as before. 

While all the derived geometry collapses to formal consequences of the definition of the moment map for the $G_{\vv}$-action on $T^* \text{Rep}_Q(\vv, \ww)$ in this simple example, the logic used here will continue to work in more interesting situations. 

Finally, notice that inclusion of the stable locus $\iota_{\text{st}}: X \hookrightarrow \FX$ induces a natural pullback map
\[
\iota^*_{\text{st}}: \text{DR}^\bullet(\FX) \to \text{DR}^\bullet(X)
\]
though the target of this map is awkward to describe explicitly as a complex, because $\mu^{-1}(0)^{\text{st}}$ is no longer affine. This means that, if we make a complex model fully explicit, the differential receives an additional contribution from the Čech differential relative to some covering\footnote{To say it more precisely, $\text{DR}^\bullet(X)$ is the homotopy limit of de Rham complexes on each open, intersection, triple intersection, ... of some good cover of $X$, since $\text{DR}^\bullet(-)$ is a stack.} of the stable locus. It is therefore easiest to either restrict attention to forms on $X$ pulled back from $\FX$, or work locally on $X$ to ignore the Čech differential. Of course, in the latter approach one eventually has to glue the local constructions, but in fortunate situations this may be done in a naive fashion. We will find similar phenomena when writing explicit cocycles in $\text{DR}^\bullet(\QM^\xi(X))$ below. 

\section{DG models for quasimap moduli} \label{sec: dg quasimap}
\subsection{Quiver model of $\QM^\xi(X)$ via dgas}
Now we will apply this machinery to explicitly describe the $(-1)$-shifted symplectic structure on $\QM^\xi(X)$ in the quiver model. It will also be convenient for later purposes to consider as an intermediate step $\QM^{\text{ns $\infty$}}(X) \coloneqq \text{ev}_\infty^{-1}(X)$, which fits into an obvious pullback diagram 
\[
\begin{tikzcd}
    \QM^\xi(X) \arrow[d] \arrow[r] & \QM^{\text{ns $\infty$}}(X) \arrow[d, "\text{ev}_\infty"] \\ 
    \xi \arrow[r, hook] & X. 
\end{tikzcd}
\]
Each connected component of $\QM^{\text{ns $\infty$}}(X)$ sits as a stable locus inside a larger derived stack $\Sh{Q}^\infty$ defined in \eqref{eq: define Qinfty} below, which is easier to describe in the quiver language (thought it will still require some preparation). Our presentation of the material in this section was influenced by the papers \cite{Anel_2022}, \cite{Benini_2021}.

Fix a connected component $\QM^{\text{ns $\infty$}}_d(X)$ of degree $d \in \mathbb{Z}_{\geq 0}^{Q_0}$. Consider the vector space 
\[
\text{Rep}_{Q^{\QM}_{\vv, \ww}}(d) \times T^* \text{Rep}_Q(\vv, \ww)
\]
and form the derived quotient stack 
\begin{equation} \label{eq: define M}
M \coloneqq [ (\text{Rep}_{Q^{\QM}_{\vv, \ww}}(d) \times \mu^{-1}(0))/P ]
\end{equation}
where the parabolic subgroup $P \subset GL(U \oplus V) = \prod_{i \in Q_0} GL(U_i \oplus V_i)$ denotes the collection of matrices of the form 
\[
P \coloneqq \Bigg\{ \begin{pmatrix} g_{U_i} && 0 \\ x_i && g_{V_i} \end{pmatrix}_{i \in Q_0} \Bigg\} \subset \prod_{i \in Q_0} GL(U_i \oplus V_i)
\]
with $g_{U_i} \in GL(U_i), g_{V_i} \in GL(V_i)$, $x_i \in \text{Hom}(U_i, V_i)$. Given our earlier notation $G_{\vv} = \prod_{i \in Q_0} GL(V_i)$, it woudl be more logical to call $P$ as $P_{d, \vv}$, but we will omit the subscripts in an attempt to reduce clutter. In particular, $G_{\vv}$ appears as a Levi factor of $P$.

The quiver data in $\text{Rep}_{Q^{\QM}} \times \mu^{-1}(0)$ is just the complete list of $(B_1, I, B_2, J_3, B_3, J_2, \alpha, \xi_2, \xi_3, \iota, \jmath)$ of all arrows appearing in the quasimap quiver of Section \ref{sec: quiver critical description quasimaps}, with the evaluation data satisfying the moment map equation (in the derived sense). $p \in P$ acts on $\text{Rep}_{Q^{\QM}} \times \mu^{-1}(0)$ via 
\begin{equation} \label{eq: P group action}
\begin{split}
    p \cdot \begin{pmatrix} B_{1, i} && I_i \end{pmatrix}  & = g_{U_i}^{-1} \begin{pmatrix} B_{1, i} && I_i \end{pmatrix} \begin{pmatrix} g_{U_i} && 0 \\ x_i && g_{V_i} \end{pmatrix} \\
    p \cdot \begin{pmatrix} B_{2, e} && 0 \\ J_{3, e} && \xi_{2, e} \end{pmatrix} & = \begin{pmatrix} g_{U_{t(e)}} && 0 \\ x_{t(e)} && g_{V_{t(e)}} \end{pmatrix}^{-1}\begin{pmatrix} B_{2, e} && 0 \\ J_{3, e} && \xi_{2, e} \end{pmatrix} \begin{pmatrix} g_{U_{s(e)}} && 0 \\ x_{s(e)} && g_{V_{s(e)}} \end{pmatrix} \\
    p \cdot \begin{pmatrix} B_{3, e} && 0 \\ -J_{2, e} && \xi_{3, e} \end{pmatrix} & = \begin{pmatrix} g_{U_{s(e)}} && 0 \\ x_{s(e)} && g_{V_{s(e)}} \end{pmatrix}^{-1} \begin{pmatrix} B_{3, e} && 0 \\ -J_{2, e} && \xi_{3, e} \end{pmatrix} \begin{pmatrix} g_{U_{t(e)}} && 0 \\ x_{t(e)} && g_{V_{t(e)}} \end{pmatrix} \\
    p \cdot \begin{pmatrix} \alpha_i && \jmath_i \end{pmatrix} & = \begin{pmatrix} \alpha_i && \jmath_i \end{pmatrix}\begin{pmatrix} g_{U_i} && 0 \\ x_i && g_{V_i} \end{pmatrix} \\ 
    p \cdot \iota_i & = g_{V_i}^{-1} \iota_i. 
\end{split}
\end{equation}

We begin by describing the function dga of $M$ from \eqref{eq: define M} in detail. 

\begin{proposition} \label{prop: functions on M}
    Let $\Gamma(M, \Sh{O}_M)$ denote the differential graded algebra of derived global sections of the structure sheaf on $M$. Then
    \[
    \Gamma(M, \Sh{O}_M) \simeq \textnormal{CE}(\mathfrak{u}, \Sh{O}_{\textnormal{Rep}_{Q^{\QM}}} \otimes \Sh{O}_{\mu^{-1}(0)})^{GL(U) \times GL(V)}
    \] 
    with $\Sh{O}_{\mu^{-1}(0)}$ modeled as the Koszul dga \eqref{eq: model for zero locus mu} above. 
\end{proposition}

\begin{proof}
From definitions, 
\[
\Gamma(M, \Sh{O}_M) \simeq \Gamma(BP, \Sh{O}_{\text{Rep}_{Q^{\QM}}} \otimes \Sh{O}_{\mu^{-1}(0)})
\]
 Now $P$ fits into an exact sequence of the form \eqref{eq: reductive and unipotent} with $G^{\text{red}} = \prod_{i \in Q_0} GL(U_i) \times GL(V_i)$ and $U = \prod_{i \in Q_0} \text{Hom}(U_i, V_i)$, the latter regarded as an additive group. Then we may describe the functor $\Gamma(BP, -)$ as in the proof of Proposition \ref{prop: de rham on quotient} to get a model 
\[
\Gamma(M, \Sh{O}_M) \simeq \text{CE}(\mathfrak{u}, \Sh{O}_{\text{Rep}_{Q^{\QM}}} \otimes \Sh{O}_{\mu^{-1}(0)})^{GL(U) \times GL(V)}. 
\]
\end{proof}        
The Chevalley-Eilenberg dga featuring in Proposition \ref{prop: functions on M} admits the following explicit description: as a graded algebra, 
\[
\textnormal{CE}(\mathfrak{u}, \Sh{O}_{\textnormal{Rep}_{Q^{\QM}}} \otimes \Sh{O}_{\mu^{-1}(0)}) \simeq \Sh{O}_{\textnormal{Rep}_{Q^{\QM}}} \otimes \Sh{O}_{\mu^{-1}(0)} \otimes S^\bullet(\mathfrak{u}^*[-1]). 
\]
It is generated in degree zero by all matrix coefficients of all quiver data $$(B_1, I, B_2, J_3, B_3, J_2, \alpha, \xi_2, \xi_3, \iota, \jmath)$$ 
from $\text{Rep}_{Q^{\QM}}$(we suppress quiver indices in what follows to declutter the already rather involved notation). In degree $-1$, there is the $\mathfrak{g}_{\vv}^*$-valued Koszul generator $\chi$ featuring in \eqref{eq: model for zero locus mu} above. Finally, in degree $1$ we have a generator $\zeta$ valued in $\text{Hom}(U, V)$, where the latter is now regarded as an abelian Lie algebra.  The differential is uniquely characterized by the formulas 
\begin{alignat}{4} \label{eq: dga for M}
& \underline{\text{Degree 1}} \quad && \delta \zeta = 0 \quad && \nonumber\\
& \underline{\text{Degree 0}} \quad && \delta B_1 = I \zeta \quad && \delta I = 0 \nonumber \\
& \quad && \delta B_2 = 0 \quad && \delta J_3 = \xi_2 \zeta - \zeta B_2 \quad && \delta \xi_2 = 0\nonumber \\
& \quad && \delta B_3 = 0 \quad && \delta J_2 = \zeta B_3 - \xi_3 \zeta \quad && \delta \xi_3 = 0 \\
& \quad && \delta \alpha = \jmath \zeta \quad && \delta \jmath = 0 \nonumber \\ 
& \quad && \delta \iota = 0 && \quad && \nonumber \\
& \underline{\text{Degree $-1$}} \quad && \delta \chi = \comm{\xi_2}{\xi_3} + \iota \jmath \quad && \quad && \quad \nonumber
\end{alignat}
and the graded Leibniz rule. The differential is read off from restriction of the $P$-action in \eqref{eq: P group action} to $U \subset P$ and further passing to the $\text{Lie} \, U$-action (though in this case the last step is trivial), and summing the result with the differential on $\Sh{O}_{\mu^{-1}(0)}$. 

\begin{remark}
We can interpret $\text{Spec} \, \text{CE}(\mathfrak{u}, \Sh{O}_{\text{Rep}_{Q^{\QM}}} \otimes \Sh{O}_{\mu^{-1}(0)})$ as the affine derived stack $[ \text{Rep}_{Q^{\QM}_{\vv, \ww}}(d) \times \mu^{-1}(0)/U]$. The affine stack terminology and point of view will make the subsequent discussion much more economic: certain diagrams of quotient stacks by $P$ can be read as $GL(U) \times GL(V)$-equivariant diagrams of affine stacks, which are fully characterized by the induced diagram on global functions. 
\end{remark}

We would like to say that our usual potential 
\[
W \coloneqq \tr B_1 \comm{B_2}{B_3} + \tr J_2(B_2 I - I \xi_2) + \tr J_3(B_3 I - I\xi_3) + \tr I \iota \alpha 
\]
defines a morphism $W: M \to \mathbb{A}^1$. But this is \textit{false}, because by the formulas \eqref{eq: dga for M} above, we find that 
\[
\begin{split}
\delta W & = \tr I \zeta \comm{B_2}{B_3} + \tr (\zeta B_3 - \xi_3 \zeta)(B_2 I - I \xi_2) + \tr(\xi_2 \zeta - \zeta B_2)(B_3 I - I \xi_3) + \tr I \iota \jmath \zeta \\
& = \tr I( \comm{\xi_2}{\xi_3} + \iota \jmath) \zeta \neq 0.
\end{split}
\]
Therefore, $W$ does not give rise to a well-defined cocycle in $\Gamma(M, \Sh{O}_M)$, so corresponds to no well-defined global function on the derived stack $M$. Fortunately, its failure to be $\delta$-closed is $\delta$-exact, and we easily obtain the following 
\begin{lemma} \label{lemma: improved W}
The expression
\begin{equation} \label{eq: improved W}
\overline{W} \coloneqq W - \tr I \chi \zeta
\end{equation}
is of total cohomological degree zero and satisfies $\delta \overline{W} = 0$, so it defines a degree zero cocycle in $\Gamma(M, \Sh{O}_M)$. Thus, there is a well-defined map $\overline{W}: M \to \mathbb{A}^1$, and its differential defines a section $$d \overline{W} : M \to T^*M$$ where $T^*M = \textnormal{Spec} \, \textnormal{Sym}_{\Sh{O}_M}( \mathbb{T}_M)$ denotes the cotangent stack to $M$.
\end{lemma}

\begin{remark}
Conceptually, what is going on here is the following. In the underived setting, in checking the $P$-invariance of $W$ one must use the equation $\comm{\xi_2}{\xi_3} + \iota \jmath = 0$. On the derived zero locus $\mu^{-1}(0)$, this equation is only enforced in a homotopy-coherent sense, so $W$ is only $P$-invariant in a homotopy coherent sense. The correction term in $\overline{W}$ can be thought of as a witness to the $P$-invariance of the original $W$.
\end{remark} 

We in addition have a canonical projection $p: M \to \FX$, which induces a natural morphism $\mathbb{T}_M \to p^*\mathbb{T}_{\FX}$ on tangent complexes. Writing 
\[
p^*T^*_{\FX} \coloneqq \text{Spec} \, \text{Sym}_{\Sh{O}_M}(p^* \mathbb{T}_{\FX})
\]
the morphism $\mathbb{T}_M \to p^* \mathbb{T}_{\FX}$ induces a morphism $g : p^*T^*_{\FX} \to T^*M$, which is the inclusion of a closed substack. We may finally define $\Sh{Q}^\infty$ as the derived fiber product 
\begin{equation} \label{eq: define Qinfty}
\begin{tikzcd}
    \Sh{Q}^\infty \arrow[r] \arrow[d] & p^*T^*_{\FX} \arrow[d, "g", hook] \\
    M \arrow[r, "d \overline{W}"] & T^*M. 
\end{tikzcd}
\end{equation}
In English, this diagram says that $\Sh{Q}^\infty$ is the critical locus of $\overline{W}$ on $M$ taken relative to $\FX$, i.e. it is the total space of the family of corresponding fiberwise critical loci. Denote the composition 
\[
\begin{tikzcd}
    \Sh{Q}^\infty \arrow[r] & M \arrow[r, "p"] & \FX 
\end{tikzcd}
\]
by $e$. The following is an immediate consequence of (the proof of) Theorem \ref{thm: spencer's result derived}. 
\begin{lemma} \label{lemma: QM open in quiver}
We have 
 \[   
\QM^{\textnormal{ns $\infty$}}(X) \simeq  (\textnormal{quiver data satisfy $\mathbb{C}[B_1]I(V) = U$}) \cap e^{-1}(X) \hookrightarrow \Sh{Q}^\infty.
\]
where the arrow denotes inclusion of an open substack. Likewise we may define $\Sh{Q}^\xi \coloneqq e^{-1}(\xi)$ as the derived fiber of $e$ over $\xi \in \FX$. If $\xi \in X$, then we similarly have 
\[
\QM^\xi(X) \simeq (\textnormal{quiver data satisfy $\mathbb{C}[B_1]I(V) = U$}) \cap \Sh{Q}^\xi \hookrightarrow \Sh{Q}^\infty.
\]
\end{lemma}

\begin{remark}
    The advantage of working with $\Sh{Q}^\infty$ is that, as we will show below, we may apply Proposition \ref{prop: de rham on quotient} to describe differential forms \textit{globally} on $\Sh{Q}^\infty$, whereas in working with $\QM^{\text{ns $\infty$}}(X)$ directly we would need to glue local descriptions. It will be more convenient to write equations using the model we get for $\Sh{Q}^\infty$, and pull them back to the stable locus by hand.
\end{remark}

Next we aim to give explicit models for $\Gamma(\Sh{Q}^\infty, \Sh{O}_{\Sh{Q}^\infty})$ and $\text{DR}^\bullet(\Sh{Q}^\infty)$ via differential graded algebras. The idea will be to recognize every term in the defining diagram \eqref{eq: define Qinfty} as a global quotient by $P$, and then interpret the diagram as a $GL(U) \times GL(V)$-equivariant diagram of affine stacks. This will again allow us to build the sought-after dga models out of Koszul and Chevalley-Eilenberg dgas. 

First, notice that because $M$ in \eqref{eq: define M} is a global quotient by $P$, we have an equivalence of derived stacks
\[
T^*M \simeq [\nu^{-1}_P(0)/P]
\]
where $\nu_P: T^*(\text{Rep}_{Q^{\QM}_{\vv, \ww}}(d) \times \mu^{-1}(0)) \to (\text{Lie} \, P)^*$ is the moment map for the induced $P$-action on the cotangent bundle to the prequotient in \eqref{eq: define M}. This establishes that the bottom right term in \eqref{eq: define Qinfty} is a global quotient by $P$, but it will be convenient to develop its dga model in greater detail. 

We start with $T^*(\text{Rep}_{Q^{\QM}_{\vv, \ww}}(d) \times \mu^{-1}(0))$ itself.
\begin{proposition} \label{prop: functions on T*Rep}
The differential graded algebra 
\[
\Gamma(T^*(\textnormal{Rep}_{Q^{\QM}} \times \mu^{-1}(0)), \Sh{O}_{T^*(\textnormal{Rep}_{Q^{\QM}} \times \mu^{-1}(0))})
\]
whose spectrum is the affine stack $T^*(\textnormal{Rep}_{Q^{\QM}} \times \mu^{-1}(0))$, admits the following explicit characterization. In degree zero we have the usual generators given by all matrix coefficients of $(B_1, B_2, B_3, J_2, J_3, I, \alpha, \xi_2, \xi_3, \iota, \jmath)$. Also in degree zero, we introduce a dual set of generators consisting of all matrix coefficients of
\[
\begin{split}
(p_{B_i}, p_{J_i}, p_I, p_\alpha) & \in (\textnormal{Rep}_{Q^{\QM}_{\vv, \ww}}(d))^\vee\\
(p_{\xi_2}, p_{\xi_3}, p_{\iota}, p_{\jmath}) & \in (T^*\textnormal{Rep}_Q(\vv, \ww))^\vee
\end{split}
\]
where the vector spaces on the RHS should be regarded canonically as cotangent fibers at the respective points. In degree $-1$ we have the Koszul generator $\chi$ valued in $\mathfrak{g}_{\vv}^*$, and in degree $+1$ we have its dual $p_\chi$ valued in $\mathfrak{g}_{\vv}$. The differential is uniquely determined by the formulas 
\begin{alignat}{5} \label{eq: dga for T*Rep}
& \underline{\textnormal{Degree 1}} \quad && \delta p_\chi = 0 \quad && \quad && \quad && \quad \nonumber \\
& \underline{\textnormal{Degree 0}} \quad && \delta B_i = 0  \quad && \delta p_{B_i} = 0 \quad  && \delta J_i = 0 \quad && \delta p_{J_i} = 0 \nonumber \\ 
& \quad && \delta I = 0 \quad && \delta p_I = 0 && \delta \alpha = 0 \quad && \delta p_{\alpha} = 0 \nonumber  \\  
& \quad && \delta \xi_2 = 0 \quad && \delta p_{\xi_{2}} = \comm{p_{\chi}}{\xi_3} \quad && \delta \xi_3 = 0 \quad && \delta p_{\xi_3} = \comm{\xi_2}{p_{\chi}} \quad \\ 
& \quad && \delta \iota = 0 \quad && \delta p_{\iota} = -\jmath p_\chi \quad && \delta \jmath = 0 && \delta p_{\jmath} = -p_\chi \iota \quad \nonumber \\ 
& \underline{\textnormal{Degree $-1$}} \quad && \delta \chi = \comm{\xi_2}{\xi_3} + \iota \jmath \quad && \quad && \quad && \quad \nonumber
\end{alignat}
and the graded Leibniz rule. 
\end{proposition}

\begin{proof}
By definition we have $T^*(\text{Rep}_{Q^{\QM}_{\vv, \ww}}(d) \times \mu^{-1}(0)) = \text{Spec} \, \text{Sym}_{\Sh{O}_{\text{Rep}} \otimes \Sh{O}_{\mu^{-1}(0)}}( \mathbb{T}_{\text{Rep} \times \mu^{-1}(0)})$, and since $\text{Rep}_{Q^{\QM}} \times \mu^{-1}(0)$ is an affine derived scheme we conclude its cotangent bundle is an affine stack. The tangent complex may be modeled as 
\[
\mathbb{T}_{\text{Rep}_{Q^{\QM}} \times \mu^{-1}(0)} \simeq \begin{tikzcd} T_{\text{Rep}_{Q^{\QM}}} \oplus T_{T^* \text{Rep}_Q} \eval_{\mu^{-1}(0)} \arrow[r, "d\mu"] & \mathfrak{g}^*_{\vv} \otimes \Sh{O}_{\text{Rep}_{Q^{\QM}} \times \mu^{-1}(0)}.
\end{tikzcd}
\]
Resolving this as a dg-module over the Koszul dga \eqref{eq: model for zero locus mu} modeling $\Sh{O}_{\mu^{-1}(0)}$ and writing out $d\mu$ in coordinates, we obtain the expressions above. 
\end{proof}

We may use the dga model of Proposition \ref{prop: functions on T*Rep} to get a model for $\text{DR}^\bullet(T^*(\text{Rep}_{Q^{\QM}} \times \mu^{-1}(0)))$, relative to which the canonical $0$-shifted symplectic form on the cotangent bundle is
\begin{equation} \label{eq: symplectic form T*Rep}
\begin{split}
\omega_{T^*(\text{Rep} \times \mu^{-1}(0))} & = \sum_{i = 1}^3 \tr dp_{B_i} \wedge dB_i + \sum_{i = 2}^3 \tr dp_{J_i} \wedge dJ_i + \tr dp_I \wedge dI + \tr dp_{\alpha} \wedge d\alpha \\
& + \sum_{i = 2}^3 \tr dp_{\xi_{i}} \wedge d\xi_i + \tr dp_{\iota} \wedge d\iota + \tr dp_{\jmath} \wedge d\jmath + \tr dp_\chi \wedge d\chi. 
\end{split}
\end{equation}
The higher closure data $\omega_i = 0$ for $i \geq 1$, and it is instructive as a check on various sign rules that $\delta \omega_{T^*(\text{Rep} \times \mu^{-1}(0))} = 0$.

Given \eqref{eq: symplectic form T*Rep} and \eqref{eq: P group action}, it is a routine exercise to compute the moment map $\nu_P$ for the $P$-action on $T^*(\text{Rep}_{Q^{\QM}_{\vv, \ww}} \times \mu^{-1}(0))$. Notice that, as a vector space:
\[
\text{Lie}(P)^* \simeq \text{End}(U) \oplus \text{Hom}(V, U) \oplus \text{End}(V)
\]
relative to which the moment map decomposes into components as $\nu_P = (\nu_\beta, \nu_\gamma, \nu_\delta)$. One readily computes for the moment map: 
\begin{equation} \label{eq: parabolic moment maps}
\begin{split}
    \nu_\beta & = \sum_{i = 1}^3 \comm{p_{B_i}}{B_i} + \sum_{i = 2}^3 p_{J_i} J_i + p_\alpha \alpha - I p_I \\
    \nu_\gamma & = B_3 p_{J_2} - p_{J_2} \xi_3 + p_{J_3} \xi_2 - B_2 p_{J_3} + p_\alpha \jmath + p_{B_1} I \\ 
    \nu_\delta & = p_I I - \sum_{i = 2}^3 J_i p_{J_i} + \sum_{i = 2}^3 \comm{p_{\xi_i}}{\xi_i} - \iota p_\iota + p_\jmath \jmath - \comm{p_\chi}{\chi}. 
\end{split}
\end{equation}
Recall our sign convention for commutators, so that $\comm{p_\chi}{\chi} = p_\chi \chi + \chi p_\chi$. 

\begin{lemma} \label{lemma: zero locus parabolic moment map}
    The affine substack $\nu_P^{-1}(0) \hookrightarrow T^*(\textnormal{Rep}_{Q^{\QM}} \times \mu^{-1}(0))$ may be viewed as the spectrum of the following differential graded algebra. To \eqref{eq: dga for T*Rep}, we add $\textnormal{Lie}(P)^*$-valued Koszul generators $\chi_P = (\chi_\beta, \chi_\gamma, \chi_\delta)$ in degree $-1$, and uniquely extend the differential from $\eqref{eq: dga for T*Rep}$ by the assignments $\delta \chi_\beta = \nu_\beta$, $\delta \chi_\gamma = \nu_\gamma$, $\delta \chi_\delta = \nu_\delta$, with $\nu_P = (\nu_\beta, \nu_\gamma, \nu_\delta)$ given by \eqref{eq: parabolic moment maps}. 
\end{lemma}

\begin{proof}
    We have a pullback diagram 
    \[
    \begin{tikzcd}
        \nu_{P}^{-1}(0) \arrow[r, hook] \arrow[d] & T^*(\text{Rep}_{Q^{\QM}} \times \mu^{-1}(0)) \arrow[d, "\nu_P"] \\ 
        0 \arrow[r, hook] & (\text{Lie}(P))^*.
    \end{tikzcd}
    \]
    This reduces the question to modeling the derived tensor product $\Sh{O}_0 \otimes^L_{\Sh{O}_{(\text{Lie}(P))^*}} \Sh{O}_{T^*(\text{Rep} \times \mu^{-1}(0)}$. Choosing a cofibrant replacement for $\Sh{O}_{0 \in \text{Lie}(P)^*}$ as in the argument leading to \eqref{eq: model for zero locus mu} gives the result. 
\end{proof}

We are really interested in $T^*M \simeq [\nu_P^{-1}(0)/P]$, and we have the following. 
\begin{proposition} \label{prop: functions on T*M}
The differential graded algebra of global functions on $T^*M$ may be modeled as 
\[
\Gamma(T^*M, \Sh{O}_{T^*M}) \simeq \textnormal{CE}(\mathfrak{u}, \Sh{O}_{\nu^{-1}_P(0)})^{GL(U) \times GL(V)}
\]
where $\Sh{O}_{\nu^{-1}_P(0)}$ is shorthand for the dga described in Lemma \ref{lemma: zero locus parabolic moment map}. The Chevalley-Eilenberg dga $\textnormal{CE}(\mathfrak{u}, \Sh{O}_{\nu^{-1}_P(0)})$ models global functions on the affine stack $[\nu^{-1}_P(0)/U]$, and admits the following explicit description. 

In degree zero, we have the usual generators $(B_i, J_i, I, \alpha, \xi_i, \iota, \jmath)$ and their duals $(p_{B_i}, p_{J_i}, p_I, p_\alpha, p_{\xi_i}, p_\iota, p_\jmath)$ as in Proposition \ref{prop: functions on T*Rep}. In degree $-1$, we have the $\mathfrak{g}^*_{\vv}$-valued Koszul generator $\chi$, and the $\textnormal{Lie}(P)^*$-valued Koszul generators $\chi_P = (\chi_\beta, \chi_\gamma, \chi_\delta)$. In degree $+1$, we have the $\mathfrak{g}_{\vv}$-valued generator $p_\chi$ and the $\textnormal{Hom}(U, V)$-valued Chevalley-Eilenberg generator $\zeta$. The differential is uniquely determined by the following assignments
\begin{alignat}{4} \label{eq: dga for T*M}
& \underline{\textnormal{Degree 1}} \quad && \delta \zeta = 0 \quad && \delta p_\chi = 0 \quad && \nonumber \\
& \underline{\textnormal{Degree 0}} \quad && \delta B_1 = I \zeta \quad && \delta I = 0 \quad && \nonumber\\
& \quad && \delta B_2 = 0 \quad && \delta J_3 = \xi_2 \zeta - \zeta B_2 \quad && \delta \xi_2 = 0 \nonumber \\
& \quad && \delta B_3 = 0 \quad && \delta J_2 = \zeta B_3 - \xi_3 \zeta  \quad && \delta \xi_3 = 0 \nonumber \\
& \quad && \delta \alpha = \jmath \zeta \quad && \delta \jmath = 0 \quad && \delta \iota = 0 \nonumber \\
& \quad && \delta p_{B_1} = 0 \quad && \delta p_I = - \zeta p_{B_1} \quad && \quad \\ 
& \quad && \delta p_{B_2} = p_{J_3} \zeta \quad && \delta p_{J_3} = 0 \quad && \delta p_{\xi_2} = \comm{p_\chi}{\xi_3} - \zeta p_{J_3} \nonumber \\ 
& \quad && \delta p_{B_3} = - p_{J_2} \zeta \quad && \delta p_{J_2} = 0 \quad && \delta p_{\xi_3} = \comm{\xi_2}{p_\chi} + \zeta p_{J_2} \nonumber \\
& \quad && \delta p_\alpha = 0 \quad && \delta p_\jmath = - p_\chi \iota - \zeta p_\alpha \quad && \delta p_{\iota} = -\jmath p_\chi \nonumber \\
& \underline{\textnormal{Degree $-1$}} \quad && \delta \chi = \comm{\xi_2}{\xi_3} + \iota \jmath && \quad && \nonumber \\
& \quad && \delta \chi_\beta = \nu_\beta - \chi_\gamma \zeta \quad && \delta \chi_\gamma = \nu_\gamma \quad && \delta \chi_\delta = \nu_\delta - \zeta \chi_\gamma \nonumber
\end{alignat}
and the graded Leibniz rule. $\nu_\beta, \nu_\gamma, \nu_\delta$ are given by \eqref{eq: parabolic moment maps} above. 
\end{proposition}

\begin{proof}
    From $T^*M \simeq [\nu^{-1}_P(0)/P]$ and definitions, we have 
    \[
    \Gamma(T^*M, \Sh{O}_{T^*M}) \simeq \Gamma(BP, \Gamma(\nu^{-1}_P(0), \Sh{O}_{\nu^{-1}_P(0)}))
    \]
    and the first part of the claim follows by the same reasoning as the proof of Proposition \ref{prop: functions on M}. For the second part of the claim, the $\zeta$-independent part of the differential follows by Proposition \ref{prop: functions on T*Rep} and Lemma \ref{lemma: zero locus parabolic moment map}. The $\zeta$-dependent part of the differential follows from \eqref{eq: P group action}, the induced action on the cotangent bundle, and the coadjoint action on $\text{Lie}(P)^*$, all restricted to the subgroup $U \subset P$ and further linearized along $\mathfrak{u} = \text{Lie}(U)$. 
\end{proof}
From Propositions \ref{prop: functions on T*M} and \ref{prop: de rham on quotient}, we get a Cartan model for $\text{DR}^\bullet(T^*M)$. Relative to the notation in Proposition \ref{prop: de rham on quotient}, write $\phi^{\text{red}} = (\phi_U, \phi_V)$. The Liouville one-form on $T^*M$ is then
\begin{equation} \label{eq: Liouville on T*M}
\begin{split}
    \lambda_{T^*M} & = \sum_{i = 1}^3 \tr (p_{B_i} dB_i) + \sum_{i = 2}^3 \tr(p_{J_i} dJ_i) + \tr(p_I dI) + \tr(p_\alpha d\alpha) \\
    & + \sum_{i = 2}^3 \tr(p_{\xi_i} d\xi_i) + \tr (p_\iota d \iota) + \tr(p_\jmath d\jmath) + \tr(p_\chi d \chi) \\
    & - \tr( \phi_U \chi_\beta) + \tr( d\zeta \chi_\gamma) - \tr(\phi_V \chi_\delta)
\end{split}
\end{equation}
in this model. 

The 0-shifted symplectic form on $T^*M$ in the Cartan model is described as $\omega_{T^* M} = d\lambda_{T^* M}$ (and as usual the higher closure data $\omega_i = 0$, $i \geq 1$). It is straightforward to verify that $(\iota_{V(\phi)} + \delta) \lambda_{T^*M} = 0$, essentially as a formal consequence of the definition of $\nu_P$, the $P$-invariance of $\lambda_{T^*M}$, and sign conventions. 

We are at last in a position to begin analyzing the defining diagram \eqref{eq: define Qinfty} of $\Sh{Q}^\infty$. Relative to the coordinate presentations \eqref{eq: dga for T*M}, \eqref{eq: dga for M}, the section $d\overline{W}: M \to T^*M$ described in Lemma \ref{lemma: improved W} is defined as a dga map as follows. We set 
\[
\lambda_{T^*M} = d\overline{W}
\]
and form the quotient of \eqref{eq: dga for T*M} by the relations obtained by reading off components. It is amusing to verify that this gives rise to a well-defined $GL(U) \times GL(V)$-equivariant morphism of dgas from \eqref{eq: dga for T*M} to \eqref{eq: dga for M}; this is essentially a formal consequence of $\delta \overline{W} = 0$ and the $GL(U) \times GL(V)$-invariance of $\overline{W}$.

Relative to the coordinate presentation \eqref{eq: dga for T*M}, the closed embedding $g: p^* T^*_\FX \to T^*M$ corresponds to setting $p_{B_i} = p_{J_i} = p_I = p_\alpha = \chi_\beta = \chi_\gamma =0$, and is induced by the corresponding quotient map on dgas. Denote this dga quotient as 
\[
\begin{tikzcd}
    \text{CE}(\mathfrak{u}, \Sh{O}_{\nu^{-1}_P(0)}) \arrow[r, "g^*", two heads] & \text{CE}(\mathfrak{u}, \Sh{O}_Y)
\end{tikzcd}
\]
and correspondingly write $p^* T^*_\FX \simeq  [Y/P]$ for certain affine stack $Y$. As a final bit of preparation, define a new dga $\widetilde{\text{CE}}(\mathfrak{u}, \Sh{O}_Y)$ by adding to the dga \eqref{eq: dga for T*M} generators (matrix coefficients of)
\[
\begin{split}
    (\eta_{B_1}, \eta_I)  & \in \text{End}(U) \oplus \text{Hom}(U, V) \\
    (\eta_{B_2}, \eta_{J_3}) & \in \text{Hom}(U, U) \oplus \text{Hom}(V, U) \\ 
    (\eta_{B_3}, \eta_{J_2}) & \in \text{Hom}(U, U) \oplus \text{Hom}(V, U) \\
    \eta_\alpha & \in \text{Hom}(W, U)
\end{split}
\]
all in cohomological degree $-1$, and generators (matrix coefficients of) 
\[
(\beta, \gamma) \in \text{End}(U) \oplus \text{Hom}(V, U)
\]
in cohomological degree $-2$. The differential is uniquely extended by insisting
\begin{alignat}{2} \label{eq: kill the p's}
& \delta \eta_{B_1} = p_{B_1} \qquad && \delta \eta_I = p_I -\zeta \eta_{B_1} \nonumber\\
& \delta \eta_{B_2} = p_{B_2} - \eta_{J_3} \zeta \qquad && \delta \eta_{J_3} = p_{J_3} \\
& \delta \eta_{B_3} = p_{B_3} + \eta_{J_2} \zeta \qquad && \delta \eta_{J_2} = p_{J_2} \nonumber\\
& \delta \eta_{\alpha} = p_\alpha \nonumber
\end{alignat}
on the new degree $-1$ generators, and 
\begin{equation} \label{eq: kill the chi's}
\begin{split}
    \delta \beta & = \sum_{i = 1}^3 \comm{\eta_{B_i}}{B_i} + \sum_{i = 2}^3 \eta_{J_i} J_i + \eta_\alpha \alpha - I \eta_I + \gamma \zeta - \chi_\beta \\
    \delta \gamma & = B_3 \eta_{J_2} - \eta_{J_2} \xi_3 + \eta_{J_3} \xi_2 - B_2 \eta_{J_3} + \eta_\alpha \jmath + \eta_{B_1} I - \chi_\gamma
\end{split}
\end{equation}
on degree $-2$ generators. The following is easily verified. 
\begin{lemma} \label{lemma: resolve g}
    Let $\widetilde{\textnormal{CE}}(\mathfrak{u}, \Sh{O}_Y)$ be defined as above, then
    \begin{enumerate}
        \item The canonical projection $\widetilde{\textnormal{CE}}(\mathfrak{u}, \Sh{O}_Y) \twoheadrightarrow \textnormal{CE}(\mathfrak{u}, \Sh{O}_Y)$ is a quasi-isomorphism. 
        \item The composition $$\textnormal{CE}(\mathfrak{u}, \Sh{O}_{\nu^{-1}_P(0)}) \hookrightarrow \widetilde{\textnormal{CE}}(\mathfrak{u}, \Sh{O}_Y) \twoheadrightarrow \textnormal{CE}(\mathfrak{u}, \Sh{O}_Y)$$ coincides with $g^*$, and defines a cofibrant replacement of it. 
    \end{enumerate}
\end{lemma}

At long last we arrive at 
\begin{proposition} \label{prop: dga for Qinfty}
    $\Sh{Q}^\infty$ defined by the pullback diagram \eqref{eq: define Qinfty} is characterized as $\Sh{Q}^\infty \simeq [\widetilde{\Sh{Q}}^\infty/GL(U) \times GL(V)]$ where $\widetilde{\Sh{Q}}^\infty$ is the affine stack defined as the spectrum of the following dga (all notations are as above). In degree $1$, we have a $\textnormal{Hom}(U, V)$-valued generator $\zeta$. In degree zero, we have generators $(B_i, J_i, I, \alpha, \xi_i, \iota, \jmath)$. In degree $-1$, we have generators $(\eta_{B_i}, \eta_{J_i}, \eta_I, \eta_{\alpha}, \chi)$. In degree $-2$, generators $\beta$, $\gamma$. The differential is uniquely determined by the assignments
    \begin{alignat}{4} \label{eq: dga for Qinfty}
    & \underline{\textnormal{Degree 1}} \quad && \delta \zeta = 0 \quad && \quad && \nonumber \\ 
    & \underline{\textnormal{Degree 0}} \quad && \delta B_1 = I \zeta \quad && \delta I = 0 \quad && \nonumber \\
    & \quad && \delta B_2 = 0 \quad && \delta J_3 = \xi_2 \zeta - \zeta B_2 \quad && \delta \xi_2 = 0 \nonumber \\
    & \quad && \delta B_3 = 0 \quad && \delta J_2 = \zeta B_3 - \xi_3 \zeta \quad && \delta \xi_3 = 0 \nonumber \\ 
    & \quad && \delta \alpha = \jmath \zeta \quad && \delta \jmath = 0 \quad && \delta \iota = 0 \nonumber  \\ 
    & \underline{\textnormal{Degree $-1$}} \quad && \delta \eta_{B_1} = \comm{B_2}{B_3} \quad && \delta \eta_I = \sum_i (J_i B_i - \xi_i J_i) + \iota \alpha - \chi \zeta - \zeta \eta_{B_1} \quad && \delta \chi = \comm{\xi_2}{\xi_3} + \iota \jmath \\ 
    & \quad && \delta \eta_{B_2} = \comm{B_3}{B_1} + IJ_2 - \eta_{J_3} \zeta \quad && \delta \eta_{J_3} = B_3 I - I \xi_3 \quad && \nonumber \\ 
    & \quad && \delta \eta_{B_3} = \comm{B_1}{B_2} + IJ_3 + \eta_{J_2} \zeta \quad && \delta \eta_{J_2} = B_2 I - I \xi_2 \quad && \nonumber \\ 
    & \quad && \delta \eta_\alpha = I \iota \quad && \quad && \nonumber 
    \end{alignat}
    as well as the assignments 
    \[
    \begin{split}
        \delta \beta & = \sum_{i = 1}^3 \comm{\eta_{B_i}}{B_i} + \sum_{i = 2}^3 \eta_{J_i} J_i + \eta_\alpha \alpha - I \eta_I + \gamma \zeta \\
        \delta \gamma & = B_3 \eta_{J_2} - \eta_{J_2} \xi_3 + \eta_{J_3} \xi_2 - B_2 \eta_{J_3} + \eta_\alpha \jmath + \eta_{B_1} I - I \chi
    \end{split}
    \]
    on degree $-2$ generators, and the graded Leibniz rule. 
\end{proposition}

\begin{proof}
By the discussion above, the defining pullback diagram \eqref{eq: define Qinfty} is equivalent to the following $GL(U) \times GL(V)$-equivariant pullback diagram of affine derived stacks: 
\[
\begin{tikzcd}
    \widetilde{\Sh{Q}}^\infty \arrow[r] \arrow[d] & \left[ Y/U \right] \arrow[d, "g"] \\ 
    \left[ (\text{Rep}_{Q^{\QM}} \times \mu^{-1}(0))/U \right] \arrow[r, "d \overline{W}"] & \left[ \nu^{-1}_P(0)/U \right]
\end{tikzcd}
\]
This means that $\widetilde{\Sh{Q}}^\infty$ is equivalent to the spectrum of any dga modeling the derived tensor product
\[
\Sh{O}_{\widetilde{\Sh{Q}}^\infty} \simeq \text{CE}(\mathfrak{u}, \Sh{O}_{\text{Rep} \times \mu^{-1}(0)}) \otimes_{\text{CE}(\mathfrak{u}, \Sh{O}_{\nu^{-1}_P(0)})}^L \text{CE}(\mathfrak{u}, \Sh{O}_Y).
\]
Using the cofibrant replacement for $g^*$ afforded by Lemma \ref{lemma: resolve g}, we may model the derived tensor product as the naive tensor product
\[
\Sh{O}_{\widetilde{\Sh{Q}}^\infty} \simeq \text{CE}(\mathfrak{u}, \Sh{O}_{\text{Rep} \times \mu^{-1}(0)}) \otimes_{\text{CE}(\mathfrak{u}, \Sh{O}_{\nu^{-1}_P(0)})} \widetilde{\text{CE}}(\mathfrak{u}, \Sh{O}_Y). 
\]
By definition of the morphism $d\overline{W}$ from Lemma \ref{lemma: improved W}, this naive tensor product is the quotient of $\widetilde{\text{CE}}(\mathfrak{u}, \Sh{O}_Y)$ by the relations obtained from reading off components of $\lambda_{T^*M} = d \overline{W}$. From \eqref{eq: Liouville on T*M} and the explicit formula for $\overline{W}$, these relations are easily seen to be
\begin{alignat*}{3}
& p_{B_1} = \comm{B_2}{B_3} \quad && p_I = \sum_{i = 2}^3 (J_i B_i - \xi_i J_i) + \iota \alpha - \chi \zeta \quad && \quad \\ 
& p_{B_2} = \comm{B_3}{B_1} + IJ_2 \quad && p_{J_3} = B_3 I - I \xi_3 \quad && p_{\xi_2} = - J_2 I \quad \\ 
& p_{B_3} = \comm{B_1}{B_2} + IJ_3 \quad && p_{J_2} = B_2 I - I \xi_2 \quad && p_{\xi_3} = -J_3 I \\ 
& p_\alpha = I \iota \quad && p_\jmath = 0 \quad && p_\iota = \alpha I \\
& p_\chi = - \zeta I \quad && \quad && \quad \\ 
& \chi_\beta = 0 \quad && \chi_\gamma = I \chi \quad && \chi_\delta = 0.
\end{alignat*}
The claim now follows from \eqref{eq: kill the p's} and \eqref{eq: kill the chi's}.
\end{proof}

We immediately deduce the following 
\begin{corollary} \label{cor: dga for Qxi}
    $\Sh{Q}^\xi = e^{-1}(\xi) \simeq [\widetilde{\Sh{Q}}^\xi/GL(U)]$ for an affine derived stack $\widetilde{\Sh{Q}}^\xi$ defined as the spectrum of the following differential graded algebra. Take the dga from Proposition \ref{prop: dga for Qinfty}, and form the quotient by the relations setting $\xi_2, \xi_3, \iota, \jmath$ to the fixed value determined by $\xi$, and set $\chi = 0$. 
\end{corollary}

\begin{proof}
    From the definition $\Sh{Q}^\xi = e^{-1}(\xi)$, and the fact that the dga in Proposition \ref{prop: dga for Qinfty} is already cofibrant over $\FX$, so the fiber over $\xi$ can be computed naively. 
\end{proof}

Via Proposition \ref{prop: dga for Qinfty}, Corollary \ref{cor: dga for Qxi} and (the proof of) Proposition \ref{prop: de rham on quotient}, we get explicit models for $\text{DR}^\bullet(\Sh{Q}^\infty)$, $\text{DR}^\bullet(\Sh{Q}^\xi)$. Notice that $\Sh{Q}^\xi$ may be identified as the pullback (i.e. a derived critical locus) 
\[
\begin{tikzcd}
    \Sh{Q}^\xi \arrow[r] \arrow[d] & \left[ \text{Rep}_{Q^{\QM}}/GL(U) \ltimes \text{Hom}(U, V) \right] \arrow[d, "z"] \\
    \left[ \text{Rep}_{Q^{\QM}}/GL(U) \ltimes \text{Hom}(U, V) \right] \arrow[r, "dW^\xi"] & T^*\left[ \text{Rep}_{Q^{\QM}}/GL(U) \ltimes \text{Hom}(U, V) \right]
\end{tikzcd}
\]
where $z$ denotes the zero section. Each of the morphisms in the fiber product are evidently $0$-shifted Lagrangian. Thus $\Sh{Q}^\xi$ acquires a canonical $(-1)$-shifted symplectic structure $\Omega_{\Sh{Q}^\xi}$ via Theorem \ref{thm: PTVV lag intersection}.

\begin{lemma} \label{lemma: shifted symplectic form quiver model}
Relative to the model for $\textnormal{DR}^\bullet(\Sh{Q}^\xi)$ obtained from Corollary \ref{cor: dga for Qxi} and Proposition \ref{prop: de rham on quotient}, the canonical $(-1)$-shifted symplectic structure on $\Sh{Q}^\xi$ takes the form 
\[
\begin{split}
    \Omega_{\Sh{Q}^\xi} & = \sum_{i = 1}^3 \tr(d \eta_{B_i} \wedge dB_i) + \sum_{i = 2}^3 \tr(d \eta_{J_i} \wedge dJ_i) + \tr(d \eta_I \wedge dI) + \tr(d \eta_\alpha \wedge d \alpha) \\ 
    & + \tr \phi_U d\beta - \tr (d\zeta \wedge d\gamma),
\end{split}
\]
with higher closure data $\Omega_i = 0$ for $i \geq 1$. 
\end{lemma}
We leave the straightforward proof to the reader. By Lemma \ref{lemma: QM open in quiver} we have an open inclusion $\iota_{\QM}: \QM^\xi(X) \to \Sh{Q}^\xi$, which gives rise to a well-defined pullback map 
\[
\iota_{\QM}^*: \text{DR}^\bullet(\Sh{Q}^\xi) \to \text{DR}^\bullet(\QM^\xi(X)). 
\]
Then $\iota_{\QM}^* \Omega_{\Sh{Q}^\xi}$ defines a $(-1)$-shifted symplectic structure on $\QM^\xi(X)$, though thus far we cannot conclude that it is canonical because it seems to depend on the identification of Theorem \ref{thm: spencer's result derived}. We will see below that in fact it is canonical.

\begin{remark}
    There is a perhaps more instructive way to describe the shifted symplectic structure on $\Sh{Q}^\xi$ via the construction above, for more on this point of view see \cite{Anel_2022}. Notice there is a canonical Lagrangian morphism $q^*T_\FX \to T^*M \times T^*\FX$, and $\Sh{Q}^\infty$ can be viewed as its fiber product with the Lagrangian morphism $d\overline{W}: M \to T^*M$. As the fiber product of two Lagrangian morphisms is always Lagrangian, we find there is a canonical Lagrangian structure on $\Sh{Q}^\infty \to T^* \FX$. Further taking fiber product with the obvious cotangent fiber Lagrangian $T^*_\xi \to T^*\FX$ makes $\Sh{Q}^\xi \to \text{pt}$ shifted Lagrangian, that is to say, describes a shifted symplectic structure on $\Sh{Q}^\xi$. We will see a similar structure emerge in our discussion of the canonical shifted symplectic structure on $\QM^\xi(X)$. 
\end{remark}

\subsection{Cech model via dgas}
Now we will repeat this discussion in the Čech model. Because the details are very similar to the quiver model, we will omit them and simply state the results. 

\begin{remark} \label{remark: finite approx Čech dga}
To make contact with the sheaf theory in the rest of this paper, it is important that we do not work in the superficially infinite type Čech model, but rather with the finite-dimensional approximations from \ref{subsec: Finite type approximations of the Čech model}. These approximations fit into the setting of Proposition \ref{prop: de rham on quotient} and the computations of this appendix, by the following rationale. We will explain it for the pair $(\Sh{M}^\xi_{G, N}, \Sh{W}^\xi)$ for concreteness, though it can be used for other spaces and functions entering at intermediate stages.

The vector bundle $N^*_{\infty, \ell, m}(\Sh{K})$ considered in Section \ref{subsec: Finite type approximations of the Čech model} can be regarded as a finite rank $G(\Sh{O}) \times G[t^{-1}]$-equivariant bundle over $G_d(\Sh{K})$ (subscript $d$ denotes a choice of connected component). Likewise we have finite-rank equivariant vector bundles $N_{\infty, \ell}[t^{-1}]$, $N_{0, m}(\Sh{O})$, on which the $G[t^{-1}]$ (resp. $G(\Sh{O})$) action factors through some quotient $G_j^\infty  = G(\mathbb{C}[t^{-1}]/t^{-j}\mathbb{C}[t^{-1}])$ (resp. $G^0_i = G(\Sh{O}/t^i\Sh{O})$, for some large enough $i, j$. Let $K^\infty_j, K^0_i$ denote the kernels. 

Then $N^*_{\infty, \ell, m}(\Sh{K})$ can be regarded as a $K^0_i \times K^\infty_j$-equivariant bundle on $G_d(\Sh{K})$, and upon pulling back everything to the open subset $G^{\leq 0}_d(\Sh{K}) \to \text{Bun}_G^{\leq 0}(\mathbb{P}^1)_d$ of stably negative bundles of given degree, we have that $K^0_i\backslash N^*_{\infty, \ell, m}(\Sh{K})/K^\infty_j$ is smooth and of finite type, and that $\Sh{M}^{\xi, \ell, m}_{G, N} = G^0_i \backslash E_{i, j, \ell, m}/ G^\infty_j$ where $E$ is again something smooth and of finite type. Now group $G^0_i \times G^\infty_j$  fits into an exact sequence of the form \eqref{eq: reductive and unipotent} with $G^\text{red} = G \times G$ and the reductive quotient identified with evaluation at $(\infty, 0)$, and it will be sufficient for our purposes to notice that Proposition \ref{prop: de rham on quotient} can be used locally on $E$ to model forms on $\text{Crit}(\Sh{W}^{\xi}_{\ell, m})$ for any large enough $\ell, m$, and the quasi-isomorphism class of such models stabilizes for $\ell, m$ large enough to what we write below.

Following the convention of the rest of this paper, rather than explicitly making reference to finite approximations we will write things formally relative to the original presentation of $(\Sh{M}^\xi_{G, N}, \Sh{W}^\xi)$. Readers with a more conservative attitude towards infinite-dimensionality issues may read this as nothing more than a very convenient shorthand capturing the stabilizing quasi-isomorphism class of dga models on any sufficiently large finite approximation\footnote{To say the same but with less jargon: on any given approximation, only finitely many Fourier modes of the formulas we write below will have non-vacuous content.}. We do not believe there is any essential difficulty making sense of things directly in infinite type, but this is not logically necessary.
\end{remark}

Notice first that the spaces $\Sh{M}^\xi_{G, N}$ fit into a family over $\FX$: begin by considering
\[
\Sh{M}_{G, N} = G(\Sh{O}) \backslash(N^*_\infty (\Sh{K}) \times N[t^{-1}] \times N(\Sh{O}) \times G(\Sh{K}) \times N^*_\infty)/G[t^{-1}]
\]
Subscripts signify that the corresponding factors are considered as $G[t^{-1}]$-spaces. In the usual notation, denote the data parameterized by the prequotient by 
$$(\Sh{B}_2(t), \Sh{I}(t), \Sh{B}_{3, \infty}(t), \Sh{J}_\infty(t), \Sh{B}_{3, 0}(t), \Sh{J}_0(t), g(t), \xi_2, \iota).$$
Denote by $\xi_3$, $\jmath$ the constant terms:
\[
\begin{split}
\Sh{B}_{3, \infty}(t) & = \xi_3 + O(t^{-1}) \\ 
\Sh{J}_\infty(t) & = \jmath + O(t^{-1}). 
\end{split}
\]
We have a map $\mu: \widetilde{\Sh{M}}_{G, N} \to \mathfrak{g}^*$ sending a tuple to $\comm{\xi_2}{\xi_3} + \iota \jmath$. Define $\Sh{M}_{G, N}$ as its homotopy fiber 
\[
\begin{tikzcd}
    \Sh{M}_{G, N} \arrow[r, hook] \arrow[d] & \widetilde{\Sh{M}}_{G, N} \arrow[d, "\mu"] \\ 
    0 \arrow[r, hook] & \mathfrak{g}^*.  
\end{tikzcd}
\]
Then $\Sh{M}^\xi_{G, N}$ is recovered as the fiber of $\Sh{M}_{G, N}$ over $\xi \in \FX$. 

The proofs of the following imitate the arguments in the quiver case, thus are omitted. When we use the term ``dga model'' in the statement of the following, we mean that dga models for global functions or differential forms can be easily read off from the following data, as in the quiver model.
\begin{proposition} \label{prop: dga for M(G, N)}
We have the following dga model for $\Sh{M}_{G, N}$. In degree $1$, we have Chevalley-Eilenberg generators $\zeta^\infty(t) \in t^{-1} \mathfrak{g}[t^{-1}]$, $\zeta^0(t) \in t \mathfrak{g}[[t]]$. In degree zero, we have generators given by all matrix coefficients of Fourier modes of the standard tuple $$(\Sh{B}_2(t), \Sh{I}(t), \Sh{B}_{3, \infty}(t), \Sh{J}_\infty(t), \Sh{B}_{3, 0}(t), \Sh{J}_0(t), g(t), \xi_2, \iota).$$ In degree $-1$, a $\mathfrak{g}^*$-valued Koszul generator $\chi$. The differential is uniquely characterized by the assignments 
\begin{alignat*}{3}
& \underline{\textnormal{Degree 1}} \quad && \delta \zeta^\infty(t) = - \frac{1}{2} \comm{\zeta^\infty(t)}{\zeta^\infty(t)} \quad && \delta \zeta^0(t) = - \frac{1}{2} \comm{\zeta^0(t)}{\zeta^0(t)} \nonumber \\
& \underline{\textnormal{Degree 0}} \quad && \delta \Sh{B}_2(t) = \comm{\Sh{B}_2(t)}{\zeta^\infty(t)} \quad && \delta \Sh{I}(t) = - \zeta^\infty(t) \Sh{I}(t) \nonumber \\ 
& \quad && \delta \Sh{B}_{3, \infty}(t) = \comm{\Sh{B}_{3, \infty}(t)}{\zeta^\infty(t)} \quad && \delta \Sh{J}_\infty(t) = \Sh{J}_\infty(t) \zeta^\infty(t) \nonumber \\ 
& \quad && \delta \Sh{B}_{3, 0}(t) = \comm{\Sh{B}_{3, 0}(t)}{\zeta^0(t)} \quad && \delta \Sh{J}_0(t) = \Sh{J}_0(t) \zeta^0(t) \nonumber \\
& \quad && \delta g(t) = - \zeta^0(t) g(t) + g(t) \zeta^\infty (t) \quad &&  \\ 
& \quad && \delta \xi_2 = 0 \quad && \delta \iota = 0 \nonumber \\ 
& \underline{\mathrm{Degree} -1} \quad && \delta \chi = \comm{\xi_2}{\xi_3} + \iota \jmath
\end{alignat*}
and the graded Leibniz rule.
\end{proposition}

\begin{lemma} \label{lemma: improved W Čech}
    Denote by 
    \[
    \begin{split}
    \Sh{W} & = \Res \tr \xi_2 \Sh{B}_{3, \infty}(t) + \Res \tr \iota \Sh{J}_\infty(t) \\
    & -\Res \tr \Sh{B}_2(t)(\Sh{B}_{3, \infty}(t) - g(t) \Sh{B}_{3, 0}(t) g^{-1}(t))  \\
    & -\Res \tr \Sh{I}(t)(\Sh{J}_\infty(t) - \Sh{J}_0(t) g^{-1}(t))
    \end{split}
    \]
    the standard potential in the Čech model. Then 
    \[
    \overline{\Sh{W}} = \Sh{W} - \Res \tr \chi \zeta^\infty(t)
    \]
    defines a degree zero cocycle in $\Gamma(\Sh{M}_{G, N}, \Sh{O}_{\Sh{M}_{G, N}})$, thus a global function $\overline{\Sh{W}}: \Sh{M}_{G, N} \to \mathbb{A}^1$. The differential defines a section $d \overline{\Sh{W}}: \Sh{M}_{G, N} \to T^*\Sh{M}_{G, N}$. 
\end{lemma}

\begin{proof}
    We need only verify $\delta \overline{\Sh{W}} = 0$. Notice 
    \[
    \begin{split}
    \delta \Sh{W} & = \Res \tr \xi_2 \comm{\Sh{B}_{3, \infty}(t)}{\zeta^\infty(t)} + \Res \tr \iota \Sh{J}_\infty(t) \zeta^\infty(t) \\
    & = \Res \tr \xi_2 \comm{\xi_3}{\zeta^\infty(t)} + \Res \tr \iota \jmath \zeta^\infty(t) \\
    & = \Res \tr(\comm{\xi_2}{\xi_3} + \iota \jmath) \zeta^\infty(t)
    \end{split}
    \]
    where we noticed only the constant term can contribute to the residue because $\zeta^\infty(t)$ has vanishing constant term. The claim follows immediately upon noting 
    \[
    \Res \tr \chi \comm{\zeta^\infty(t)}{\zeta^\infty(t)} = 0
    \]
    again because $\zeta^\infty(t)$ has vanishing constant term.
\end{proof}
The canonical map $q: \Sh{M}_{G, N} \to \FX$ induces a closed embedding $h: q^* T^*_\FX \hookrightarrow T^* \Sh{M}_{G, N}$ and we define the relative critical locus $\text{Crit}_\FX(\overline{\Sh{W}})$ by the pullback diagram 
\[
\begin{tikzcd}
    \text{Crit}_\FX(\overline{\Sh{W}}) \arrow[r] \arrow[d] & q^* T^*_\FX \arrow[d, "h", hook] \\ 
    \Sh{M}_{G, N} \arrow[r, "d\overline{\Sh{W}}"] & T^* \Sh{M}_{G, N}. 
\end{tikzcd}
\]
For later use it is convenient to note that we have an evident pullback diagram 
\[
\begin{tikzcd}
    \QM^{\text{ns $\infty$}}(X) \arrow[r, hook] \arrow[d] & \text{Crit}_\FX(\overline{\Sh{W}}) \arrow[d] \\ 
    X \arrow[r, hook] & \FX.
\end{tikzcd}
\]
Notice also that $\Res \tr$ induces isomorphisms 
\[
\begin{split}
(t^{-1}\mathfrak{g}[t^{-1}])^\vee & \simeq \mathfrak{g}(\Sh{K})/t^{-1}\mathfrak{g}[t^{-1}] \\
(\mathfrak{g}(\Sh{O}))^\vee & \simeq \mathfrak{g}(\Sh{K})/\Fg(\Sh{O}).
\end{split}
\]
The following is likewise demonstrated as in the quiver model. 
\begin{proposition} \label{prop: dga for Qinfty in the Čech model}
We have the following dga model for $\textnormal{Crit}_\FX(\overline{\Sh{W}})$. In degree 1, we have Chevalley-Eilenberg generators (all matrix coefficients of variables in) $\zeta^\infty(t) \in t^{-1}\mathfrak{g}[t^{-1}]$, $\zeta^0(t) \in t \Fg[[t]]$. In degree zero, we have all matrix coefficients of Fourier modes of the standard variables $$(\Sh{B}_2(t), \Sh{I}(t), \Sh{B}_{3, \infty}(t), \Sh{J}_\infty(t), \Sh{B}_{3, 0}(t), \Sh{J}_0(t), g(t), \xi_2, \iota).$$ 
In degree $-1$, we have generators (all matrix coefficients of) 
\[
\begin{split}
    (\eta_{2}(t), \eta_{\Sh{I}}(t)) & \in N(\Sh{K}) \\ 
    (\eta_{3, \infty}(t), \eta_{\Sh{J}_\infty}(t)) & \in N^*(\Sh{K)}/t^{-1}N^*[t^{-1}] \\ 
    (\eta_{3, 0}(t), \eta_{\Sh{J}_0}(t)) & \in N^*(\Sh{K})/N^*(\Sh{O}) \\ 
    \eta_g(t) & \in \mathfrak{g}(\Sh{K}) 
\end{split}
\]
as well as the $\mathfrak{g}^*$-valued Koszul generator $\chi$. Finally in degree $-2$ we have generators $\beta_\infty(t) \in \Fg(\Sh{K})/t^{-1}\Fg[t^{-1}]$, $\beta_0(t) \in \Fg(\Sh{K})/\Fg(\Sh{O})$. The differential is uniquely determined by the assignments 
\begin{alignat*}{3}
& \underline{\textnormal{Degree 1}} \quad && \delta \zeta^\infty(t) = -\frac{1}{2}\comm{\zeta^\infty(t)}{\zeta^{\infty}(t)} \quad && \delta \zeta^0(t) = - \frac{1}{2} \comm{\zeta^0(t)}{\zeta^0(t)} \\ 
& \underline{\textnormal{Degree 0}} \quad && \delta \Sh{B}_2(t) = \comm{\Sh{B}_2(t)}{\zeta^\infty(t)} \quad && \delta \Sh{I}(t) = - \zeta^\infty(t) \Sh{I}(t) \\
& \quad && \delta \Sh{B}_{3, \infty}(t) = \comm{\Sh{B}_{3, \infty}(t)}{\zeta^\infty(t)} \quad && \delta \Sh{J}_\infty(t) = \Sh{J}_\infty(t) \zeta^\infty(t) \\ 
& \quad && \delta \Sh{B}_{3, 0}(t) = \comm{\Sh{B}_{3, 0}(t)}{\zeta^0(t)} \quad && \delta \Sh{J}_0(t) = \Sh{J}_0(t) \zeta^0(t) \\ 
& \quad && \delta g(t) = - \zeta^0(t) g(t) + g(t)\zeta^\infty(t) \quad && \quad \\ 
& \quad && \delta \xi_2 = 0 \quad && \delta \iota = 0, 
\end{alignat*}
in degree $-1$ we have ($[ - ]_{\geq 0}$, $[ - ]_{\leq -1}$ denote the canonical projections)
\begin{alignat*}{2}
&  \delta \eta_2(t) = \Sh{B}_{3, \infty}(t) - g^{-1}(t) \Sh{B}_{3, 0}(t) g(t) - \comm{\eta_2(t)}{\zeta^\infty(t)} \quad && \delta \eta_{\Sh{I}}(t) = \Sh{J}_\infty(t) - \Sh{J}_0(t) g(t) - \eta_{\Sh{I}}(t) \zeta^\infty(t) \\
& \delta \eta_{3, \infty}(t) = \Big[ \Sh{B}_2(t) - \xi_2 \Big]_{\geq 0} - \comm{\eta_{3, \infty}(t)}{\zeta^\infty(t)} \quad && \delta \eta_{\Sh{J}_\infty}(t) = \Big[ \Sh{I}(t) - \iota \Big]_{\geq 0} - \zeta^\infty(t) \eta_{\Sh{J}_\infty}(t) \\
& \delta \eta_{3, 0}(t) = [ g(t)\Sh{B}_2(t)g^{-1}(t) ]_{\leq -1} - \comm{\eta_{3, 0}(t)}{\zeta^0(t)} \quad && \delta \eta_{\Sh{J}_0}(t) = [g(t)\Sh{I}(t)]_{\leq -1} - \zeta^0(t) \eta_{\Sh{J}_0}(t) \\
& \delta \eta_g(t) = \comm{\Sh{B}_2(t)}{\Sh{B}_{3, \infty}(t)} + \Sh{I}(t) \Sh{J}_\infty(t) - \comm{\eta_g(t)}{\zeta^\infty(t)} \quad && \quad \\ 
& \delta \chi = \comm{\xi_2}{\xi_3} + \iota \jmath \quad && \quad
\end{alignat*}
and in degree $-2$ we have 
\[
\begin{split}
    \delta \beta_\infty(t) & = \Big[ \comm{\eta_{3, \infty}(t)}{\Sh{B}_{3, \infty}(t)} + \eta_{\Sh{J}_\infty}(t) \Sh{J}_\infty(t) - \eta_g(t) \Big]_{\geq 0} + \comm{\beta_\infty(t)}{\zeta^\infty(t)} + \chi \\
    \delta \beta_0(t) & = \Big[ g(t)(\eta_g(t) +  \comm{\eta_2(t)}{\Sh{B}_2(t)} - \Sh{I}(t) \eta_{\Sh{I}}(t))g^{-1}(t) - \comm{\eta_{3, 0}(t)}{\Sh{B}_{3, 0}(t)} \\
    & - \eta_{\Sh{J}_0}(t) \Sh{J}_0(t)  \Big]_{\leq -1} + \comm{\beta_0(t)}{\zeta^0(t)}. 
\end{split}
\]
\end{proposition}
\begin{corollary} \label{cor: dga for Qxi in the Čech model}
To get a dga model for $\textnormal{Crit}(\Sh{W}^\xi)$, take the dga from Proposition \ref{prop: dga for Qinfty in the Čech model} and form the quotient setting $\xi_2, \xi_3, \iota, \jmath$ to fixed values dictated by the chosen point $\xi$, and set $\chi = 0$. 
\end{corollary}

Assume $\xi \in X$, so that $\QM^\xi(X) \simeq \text{Crit}(\Sh{W}^\xi)$. Then we have
\begin{lemma} \label{lemma: shifted symplectic form Čech model}. 
Relative to the dga model for $\textnormal{Crit}(\Sh{W}^\xi)$ constructed in Proposition \ref{prop: dga for Qinfty in the Čech model} and Corollary \ref{cor: dga for Qxi in the Čech model}, the canonical $(-1)$-shifted symplectic structure on the derived critical locus $\textnormal{Crit}(\Sh{W}^\xi)$ is given by\footnote{Notice that the signs differ in the $\eta_{\Sh{B}_{3, 0}}, \eta_{\Sh{J}_0}$ terms and there is a shift in the $\eta_g(t)$ term because in the dga model, $\delta$ acting on these is not literally the derivative of $\Sh{W}^\xi$ with respect to corresponding dual variable.} 
\[
\begin{split}
    \Omega_{\textnormal{Crit}(\Sh{W}^\xi)} & = -d\Bigg\{ \Res \tr \eta_2(t) d\Sh{B}_2(t) + \Res \tr \eta_{\Sh{I}}(t) d\Sh{I}(t) + \Res \tr \eta_{3, \infty}(t) d \Sh{B}_{3, \infty}(t) \\
    & + \Res \tr \eta_{\Sh{J}_\infty}(t) d\Sh{J}_\infty(t) - \Res \tr \eta_{\Sh{B}_{3, 0}}(t) d\Sh{B}_{3, 0}(t) - \Res \tr \eta_{\Sh{J}_0}(t) d\Sh{J}_0(t)  \\
    & - \Res \tr \Big( \eta_g(t) + \comm{\eta_2(t)}{\Sh{B}_2(t)} - \Sh{I}(t) \eta_{\Sh{I}}(t) \Big) g^{-1}(t) dg(t)  \\
    & + \Res \tr \phi_0 \beta_0(t) - \Res \tr d \zeta^0(t)  \beta_0(t) - \Res \tr d \zeta^\infty(t) \beta_\infty(t) \Bigg\}.
\end{split}
\]
$\phi_0$ denotes the equivariant variable in the Cartan complex associated to the reductive quotient $G(\Sh{O}) \to G$. The higher closure data vanish $\Omega_i = 0$ for $i \geq 1$. 
\end{lemma}

\section{Shifted symplectic structures} \label{sec: compare shifted symplectic structures}
\subsection{Intrinsic shifted symplectic structure on based mapping stacks} \label{subsec: intrinsic shifted symplectic}
Thus far, we have constructed two a priori different $(-1)$-shifted symplectic structures on $\QM^\xi(X)$, $\iota^*_{\QM} \Omega_{\Sh{Q}^\xi}$ and $\Omega_{\text{Crit}(\Sh{W}^\xi)}$. These are defined with respect to different derived critical locus presentations of $\QM^\xi(X)$, and it is difficult to compare them directly, because the corresponding chain-level models for $\text{DR}^\bullet(\QM^\xi(X))$ with respect to which they are constructed are very different. 

Nonetheless, we will argue here that $\iota^*_{\QM} \Omega_{\Sh{Q}^\xi}$ and $\Omega_{\text{Crit}(\Sh{W}^\xi)}$ are nonetheless equivalent, by the following rationale. We will show that there is an \textit{canonical and intrinsically defined} $(-1)$-shifted symplectic structure $\Omega_{\QM^\xi(X)}$ on $\QM^\xi(X)$, and unravel its definition relative to each of the dg-models for $\text{DR}^\bullet(\QM^\xi(X))$ described above. We will find that it specializes, in each presentation, to the canonical $(-1)$-shifted symplectic structure on the corresponding derived critical locus. 

We briefly recall some general notions and setup from \cite{pantev2013shiftedsymplecticstructures}. Let $D$ be a derived Artin stack, which is $\Sh{O}$-compact (see Definition 2.1 of \cite{pantev2013shiftedsymplecticstructures}). $\Sh{O}$-compactness is a technical assumption which guarantees the existence of a canonical morphism 
\[
\kappa_{F, D}: \text{DR}^\bullet(F \times D) \to \text{DR}^\bullet(F) \otimes \Gamma(D, \Sh{O}_D)
\]
for any other derived stack $F$. Assume further that $D$ carries a ``$d$-orientation'', that is a morphism of complexes 
\[
\eta_D : \Gamma(D, \Sh{O}_D) \to \mathbb{C}[-d] 
\]
satisfying the following assumption. For any perfect complex $E$ on $D$, the duality pairing 
\[
\begin{tikzcd}
    \Gamma(D, E) \otimes \Gamma(D, E^\vee) \arrow[r] & \Gamma(D, \Sh{O}_D) \arrow[r, "\eta_D"] & \mathbb{C}[-d]
\end{tikzcd}
\]
induces a morphism $\Gamma(D, E) \to \Gamma(D, E^\vee)^\vee[-d]$, which is required to be a quasi-isomorphism.

For a derived Artin stack $\FX$, let $\text{Map}(D, \FX)$ denote the mapping stack, and 
\[
u: D \times \text{Map}(D, \FX) \to \FX
\]
the universal evaluation. We have the following result. 

\begin{theorem} \label{thm: ptvv symplectic form}
    (\cite{pantev2013shiftedsymplecticstructures}, Theorem 2.5). Let $D$ be an $\Sh{O}$-compact derived Artin stack with a $d$-orientation $\eta_D$. Let $\FX$ be a derived Artin stack equipped with an $n$-shifted symplectic structure $\omega_\FX$. Assume $\textnormal{Map}(D, \FX)$ is itself derived Artin, then the composite morphism 
    \begin{multline}
        \textnormal{DR}^\bullet(\FX) \xrightarrow{u^*}\textnormal{DR}^\bullet(D \times \textnormal{Map}(D, \FX)) \xrightarrow{\kappa} \Gamma(D, \Sh{O}_D) \otimes \textnormal{DR}^\bullet(\textnormal{Map}(D, \FX)) 
        \\
        \xrightarrow{\eta_D}\textnormal{DR}^\bullet(\textnormal{Map}(D, \FX))[-d]
    \end{multline}
    applied to $\omega_\FX$ defines an $(n - d)$-shifted symplectic structure on $\textnormal{Map}(D, \FX)$. 
\end{theorem}
We will refer to shifted symplectic structures on mapping stacks constructed via this procedure as AKSZ symplectic structures. 

In \cite{calaquelagrangianmap}, Calaque studied Lagrangian structures on mapping stacks with AKSZ shifted symplectic structures. Let us further recall some notions required to state his result. Let $D$ be $\Sh{O}$-compact and $d$-oriented as above, and suppose there is a morphism $f: D \to C$ to some other $\Sh{O}$-compact derived stack $C$. A ``boundary structure'' on $f$ is a choice of nullhomotopy of the composition 
\[
\begin{tikzcd}
    \Gamma(C, \Sh{O}_C) \arrow[r, "f^*"] & \Gamma(D, \Sh{O}_D) \arrow[r, "\eta_D"] & \mathbb{C}[-d]. 
\end{tikzcd}
\]
A boundary structure is called nondegenerate if it satisfies the following assumption. For a perfect complex $E$ on $C$, we have maps
\[
\begin{tikzcd}
\Gamma(C, E) \otimes \Gamma(D, f^* E^\vee) \arrow[r, "f^*"] & \Gamma(D, f^*E) \otimes \Gamma(D, f^* E^\vee) \arrow[r] & \Gamma(D, \Sh{O}_D) \arrow[r, "\eta_D"] & \mathbb{C}[-d]. 
\end{tikzcd}
\]
Then there is an induced map $\Gamma(D, f^*E^\vee) \to \Gamma(C, E)^\vee[-d]$, and \cite{calaquelagrangianmap} writes $\Gamma(C, E)^\perp$ for its homotopy fiber. Given a boundary structure on $f$, there is a canonical map $\Gamma(C, E^\vee) \to \Gamma(C, E)^\perp$. The boundary structure is called nondegenerate if it is a quasi-isomorphism. Calaque proves the following result. 

\begin{theorem} \label{thm: calaque lagrangian structure}
    (\cite{calaquelagrangianmap}, Theorem 2.9). Let $\FX$ and $D$ satisfy the assumptions of Theorem \ref{thm: ptvv symplectic form}, let $C$ be an $\Sh{O}$-compact derived stack, and assume we are given a morphism $f: D \to C$ equipped with a nondegenerate boundary structure. Assume that the mapping stack $\textnormal{Map}(C, \FX)$ is also derived Artin. Then there is a canonical Lagrangian structure on the morphism 
    \[
    f^\star : \textnormal{Map}(C, \FX) \to \textnormal{Map}(D, \FX)
    \]
    where the target carries the $(n - d)$-shifted symplectic structure of Theorem \ref{thm: ptvv symplectic form}. 
\end{theorem}
We will be interested in applying Theorem \ref{thm: calaque lagrangian structure} in the situation where $C = \mathbb{P}^1$, $D = \text{Spec} \, \mathbb{C}[t^{-1}]/t^{-2}$, and $f: D \to C$ is the evident inclusion of $\text{Spec} \, \mathbb{C}[t^{-1}]/t^{-2}$ into $\mathbb{P}^1$ as a closed subscheme. 

Notice that $\text{Spec} \, \mathbb{C}[t^{-1}]/t^{-2}$ is naturally zero-oriented, with the orientation given by 
\[
\begin{split}
    \eta_D :  \, \, & \mathbb{C}[t^{-1}]/t^{-2} \to \mathbb{C} \\
                    & a + bt^{-1} \mapsto b = \Res(a + bt^{-1})
\end{split}
\]
i.e. taking residue. The nondegeneracy condition that an orienation must satisfy follows from the observation that $\Res(a + bt^{-1})(c + dt^{-1}) = ad + bc$ is a nondegenerate pairing. 

There is moreover a canonical boundary structure on $f$: if we model 
\begin{equation*}
\Gamma(\mathbb{P}^1, \Sh{O}_{\mathbb{P}^1}) \simeq 
\begin{tikzcd}
 \mathbb{C}[t] \oplus \mathbb{C}[t^{-1}] \arrow[r] & \mathbb{C}[t, t^{-1}]
\end{tikzcd}
\end{equation*}
as a Čech complex, then there is a canonical nullhomotopy of the composite $\Gamma(\mathbb{P}^1, \Sh{O}_{\mathbb{P}^1}) \to \mathbb{C}[t^{-1}]/t^{-2} \to \mathbb{C}$ given by the dashed arrow in the following diagram 
\[
\begin{tikzcd}
    \mathbb{C}[t, t^{-1}] \arrow[dr, dashed] & \\ 
    \mathbb{C}[t] \oplus \mathbb{C}[t^{-1}] \arrow[u] \arrow[r, "\Res_{t=\infty}", swap] & \mathbb{C}
\end{tikzcd}
\]
which is just the trivial observation that taking residue at $t = \infty$ factors through the Čech differential. Nondegeneracy of the boundary structure follows from Serre duality on $\mathbb{P}^1$ (e.g. arguing as in the proof of Claim 3.3 in \cite{calaquelagrangianmap})---the key fact that makes this work is that the ideal sheaf of $D$ is $\Sh{O}_{\mathbb{P}^1}(-2) \simeq K_{\mathbb{P}^1}$. 

Then we may apply Theorem \ref{thm: calaque lagrangian structure} to conclude that $\text{Map}(\mathbb{P}^1, \FX) \to \text{Map}(D, \FX)$ has a natural Lagrangian structure, the target being $n$-shifted symplectic if $\FX$ is. However, we can be even more explicit due to the simple nature of our situation. Notice that 
\[
\text{Map}(D, \FX) = \text{Map}(\text{Spec} \, \mathbb{C}[t^{-1}]/t^{-2}, \FX) \simeq T\FX
\]
is the tangent stack to $\FX$. Moreover it is an amusing exercise to verify that the $n$-shifted AKSZ symplectic structure on $T\FX$ is described as follows. The $n$-shifted symplectic structure $\omega_\FX$ on $\FX$ defines an isomorphism $T\FX \simeq T^*[n]\FX$ from the tangent stack to the shifted cotangent stack, and the symplectic structure on $T\FX$ viewed as a mapping stack is identified with the pullback of the canonical $n$-shifted symplectic structure on $T^*[n]\FX$ \cite{Calaque_2019} under this isomorphism.

In the second description, it is evident that for any choice of point $\xi \in \FX$, the tangent fiber $T_\xi \to T\FX$ is naturally Lagrangian (because shifted conormal stacks are Lagrangian \cite{Calaque_2019}). Combining these observations we arrive at the following. 

\begin{proposition} \label{prop: based maps are symplectic}
    Let $\FX$ be a derived Artin stack with an $n$-shifted symplectic structure, and $\xi \in \FX$ a geometric point. Assume $\textnormal{Map}(\mathbb{P}^1, \FX)$ is a derived Artin stack. Then the stack of based maps $\textnormal{Map}^\xi(\mathbb{P}^1, \FX)$ defined by the pullback square 
    \[
    \begin{tikzcd}
        \textnormal{Map}^\xi(\mathbb{P}^1, \FX) \arrow[r, hook] \arrow[d] & \textnormal{Map}(\mathbb{P}^1, \FX) \arrow[d, "\textnormal{ev}_\infty"] \\ 
        \xi \arrow[r, hook] & \FX
    \end{tikzcd}
    \]
    carries a canonical $(n-1)$-shifted symplectic structure. 
\end{proposition}

\begin{proof}
    Let $\pi: T\FX \to \FX$ be the canonical projection, and consider the diagram
    \[
    \begin{tikzcd}
        \text{Map}^\xi(\mathbb{P}^1, \FX) \arrow[r, hook] \arrow[d] & \text{Map}(\mathbb{P}^1, \FX) \arrow[d, "f^\star", swap] \arrow[dd, bend left, "\text{ev}_\infty"] \\
        T_\xi \arrow[r, hook] \arrow[d] & T\FX \arrow[d, "\pi", swap] \\ 
        \xi \arrow[r, hook] & \FX.
    \end{tikzcd}
    \]
    The bottom square and outer rectangle are evidently homotopy pullbacks, so conclude the top square is a homotopy pullback. On the other hand the top square presents $\text{Map}^\xi(\mathbb{P}^1, \FX)$ as the fiber product of two Lagrangian morphisms to an $n$-shifted symplectic target, so we conclude it carries a canonical $(n - 1)$-shifted symplectic structure by Theorem \ref{thm: PTVV lag intersection}.  
\end{proof}
Now we return to the concrete setting of the paper and study $\FX = T^*[\text{Rep}_Q(\vv, \ww)/G_{\vv}]$. If the chosen point $\xi \in X \subset \FX$ lies in the stable locus $X$, then we have $\text{Map}^\xi(\mathbb{P}^1, \FX) \simeq \QM^\xi(X)$, and deduce as a corollary of the proof of Proposition \ref{prop: based maps are symplectic} that 
\begin{corollary} \label{cor: quasimaps is shifted symplectic}
     $\QM^\xi(X)$ carries a canonical $(-1)$-shifted symplectic structure. Moreover the natural morphism 
    \[
    f^\star: \QM^{\textnormal{ns $\infty$}}(X) \to \textnormal{Map}(\textnormal{Spec} \, \mathbb{C}[t^{-1}]/t^{-2}, X) \simeq TX
    \]
    carries a canonical Lagrangian structure, and the $(-1)$-shifted symplectic structure on $\QM^\xi(X)$ is that induced by the pullback diagram 
    \begin{equation} \label{eq: QM as Lagrangian intersection}
    \begin{tikzcd}
        \QM^\xi(X) \arrow[r, hook] \arrow[d] & \QM^{\textnormal{ns $\infty$}}(X) \arrow[d, "f^\star"] \\ 
        T_\xi \arrow[r, hook] & TX.
    \end{tikzcd}
    \end{equation}
\end{corollary}

\begin{remark}
    The isotropic structure on the morphism $f^\star$ is easy to understand at the underived level: denoting the quasimap section data as $\Sh{B}_2(t), \Sh{B}_3(t), \Sh{I}(t), \Sh{J}(t)$ as usual, we have 
    \[
    (f^\star)^* \omega_{TX} = \Res_{t = \infty} \Big( \tr d \Sh{B}_2(t) \wedge d \Sh{B}_3(t) + \tr d \Sh{I}(t) \wedge d \Sh{J}(t) \Big).
    \]
    On the other hand, if $(\Sh{B}_2(t), \Sh{B}_3(t), \Sh{I}(t), \Sh{J}(t))$ extend to regular sections of the corresponding vector bundles over $\mathbb{P}^1$, the expression in parenthesis has the $t$-dependence of a global section of $\Sh{O}_{\mathbb{P}^1}$, namely constant so the residue vanishes. Theorem \ref{thm: calaque lagrangian structure} in the special case relevant for us can be interpreted as saying that this continues to work in the derived setting, with everything interpreted homotopy-coherently. 
\end{remark}

Denote the corresponding $(-1)$-shifted symplectic structure on $\QM^\xi(X)$ by $\Omega_{\QM^\xi(X)}$. The advantage of the above abstract construction is that it makes the following clear: $\Omega_{\QM^\xi(X)}$ is completely intrinsic, and does not depend on a choice of a presentation or dg model of $\QM^\xi(X)$. 

We will now invert this logic. Given the abstract construction of $\Omega_{\QM^\xi(X)}$ and an explicit dg model for $\text{DR}^\bullet(\QM^\xi(X))$, we can unwind the definition of $\Omega_{\QM^\xi(X)}$ in the corresponding dg model and find an explicit cocycle representing it. This can be used to compare the shifted symplectic structures in the quiver and Čech models for $\QM^\xi(X)$, by way of comparing them each to the intrinsic intermediate object $\Omega_{\QM^\xi(X)}$. We will demonstrate below that in the quiver model, $\Omega_{\QM^\xi(X)}$ specializes to the shifted symplectic structure of Lemma \ref{lemma: shifted symplectic form quiver model}, and in the Čech model $\Omega_{\QM^\xi(X)}$ specializes to that of Lemma \ref{lemma: shifted symplectic form Čech model}. 

Before beginning this analysis, it is useful to state the following. As usual, when we say ``dg model'', we mean that dg models for functions (and, via Proposition \ref{prop: de rham on quotient}, differential forms) can be read off from the following data, working globally on $\FX$ or locally on $X.$

\begin{lemma} \label{lemma: dga for TX}
We have the following dga model for $T\FX = \textnormal{Map}(\textnormal{Spec} \, \mathbb{C}[t^{-1}]/t^{-2}, \FX)$. All $t$-dependent expressions are taken modulo $t^{-2}$ in the following. In degree 1, we have a $\mathfrak{g}_{\vv}$-valued generator $c$. In degree zero, we introduce generators given by matrix coefficients of $(\bm{\xi}_2(t), \bm{\xi}_3(t), \bm{\iota}(t), \bm{\jmath}(t))$. In degree $-1$, we introduce the Koszul generator $\bm{\chi}(t)$. The constant terms are denoted by $(\xi_2, \xi_3, \iota, \jmath, \chi)$ as usual. The differential is uniquely characterized by the assignments 
\begin{alignat*}{3}
& \underline{\textnormal{Degree 1}} \quad && \delta c = 0 \quad && \quad \\ 
& \underline{\textnormal{Degree 0}} \quad && \delta \bm{\xi}_2(t) = t^{-1} \comm{\xi_2}{c} \quad && \delta \bm{\xi}_3(t) = t^{-1} \comm{\xi_3}{c} \\ 
& \quad && \delta \bm{\iota}(t) = -t^{-1}c \iota \quad && \delta \bm{\jmath}(t) = t^{-1} \jmath c \\  
& \underline{\mathrm{Degree } -1} \quad && \delta \bm{\chi}(t) = \comm{\bm{\xi}_2(t)}{\bm{\xi}_3(t)} + \bm{\iota}(t) \bm{\jmath}(t) - t^{-1} \comm{\chi}{c}
\end{alignat*}
and the graded Leibniz rule.
\end{lemma}

\begin{proof}
    Notice that $T\FX$ can be viewed as $[\mu^{-1}(0)/G_{\vv}]$ taken over the base ring $\mathbb{C}[t^{-1}]/t^{-2}$. The derived zero locus $\mu^{-1}(0)$ is then described by the Koszul dga \eqref{eq: model for zero locus mu} taken over $\mathbb{C}[t^{-1}]/t^{-2}$ which gives the $c$-independent part of the dga above. Then we have to take the quotient by $G_{\vv}(\mathbb{C}[t^{-1}]/t^{-2}) \simeq G_{\vv} \ltimes \mathfrak{g}_{\vv}$, where the second factor is regarded as an additive group with a $G_{\vv}$-action. The group we are quotienting by is of the form \eqref{eq: reductive and unipotent} so we may apply the standard logic. 
\end{proof}

\subsection{$\Omega_{\QM^\xi(X)}$ in the quiver model} \label{subsec: prove quiver model symplectic thm}
Let us now unravel the defining diagram \eqref{eq: QM as Lagrangian intersection} in terms of the quiver model. The most challenging point is to describe the morphism $f^\star$ concretely. It will turn out to be enough to analyze the diagram \eqref{eq: QM as Lagrangian intersection} locally, in which case we can make some simplifying assumptions. 

First, recall the defining fiber sequence of the quiver presentation:
\begin{equation} \label{eq: quiver fiber sequence}
    \begin{tikzcd}[ampersand replacement=\&]
        \Sh{V} \arrow[r] \& F \otimes \Sh{O}_{\QM^{\text{ns $\infty$}}(X) \times \mathbb{P}^1} \arrow[r, "{\begin{psmallmatrix} B_1 - t & I \end{psmallmatrix}}"] \& U \otimes \Sh{O}_{\QM^{\text{ns $\infty$}}(X) \times \mathbb{P}^1}(1)
    \end{tikzcd}
\end{equation}
over $\QM^{\text{ns $\infty$}}(X) \times \mathbb{P}^1$. Restricting to $\tilde U \times \text{Spec} \, \mathbb{C}[t^{-1}]/t^{-2}$ where $\tilde U$ is some small enough open locus in $\QM^{\text{ns $\infty$}}(X)$, we may assume 
\begin{enumerate}
    \item $\Sh{V}$ is trivial over $\tilde U \times \text{Spec} \, (\mathbb{C}[t^{-1}]/t^{-2})$. 
    \item \eqref{eq: quiver fiber sequence} splits over $\tilde U \times \text{Spec} \, (\mathbb{C}[t^{-1}]/t^{-2})$.
\end{enumerate}
Moreover $\QM^{\text{ns $\infty$}}(X)$ can be covered with such $\tilde U$. Relative to the trivialization, $\Sh{V} \simeq V[t^{-1}]/t^{-2}$ (where $V$ is the tautological bundle pulled back from $X$), the first arrow in \eqref{eq: quiver fiber sequence} can be described as a tuple 
\begin{equation*}
    \begin{pmatrix} \psi^\infty(t) \\ \chi^\infty(t) \end{pmatrix}: V[t^{-1}]/t^{-2} \to F[t^{-1}]/t^{-2}
\end{equation*}
solving 
\begin{equation*}
    \Big( \frac{B_1}{t} - 1 \Big) \psi^\infty(t) + \frac{1}{t}I \chi^\infty(t) = 0 \mod t^{-2}.
\end{equation*}
Notice that by adjusting the trivialization of $\Sh{V}$, we may as well assume $\chi^\infty(t) = 1 \mod t^{-2}$, at which point $\psi^\infty(t)$ is completely specified:
\[
\psi^\infty(t) =  \frac{1}{t}I + O(t^{-2}). 
\]
The splitting gives us another tuple 
\begin{equation*}
    \begin{pmatrix} \pi^\infty(t) && \rho^\infty(t)  \end{pmatrix}: F[t^{-1}]/t^{-2} \to V[t^{-1}]/t^{-2}
\end{equation*}
which must solve 
\begin{equation*}
    \pi^\infty(t) \psi^\infty(t) + \rho^\infty(t) \chi^\infty(t) = 1 \mod t^{-2}. 
\end{equation*}
Notice that if we write $\rho^\infty(t) = \rho^\infty_0 + \rho^\infty_1 t^{-1} + O(t^{-2})$, $\pi^\infty(t) = \pi^\infty_0 + \pi^\infty_1 t^{-1} + O(t^{-2})$, this equation implies $\rho^\infty_0 = 1$, $\rho^\infty_1 = -\pi_0^\infty I$. 

To understand the form of the map $f^\star$ over $\tilde U$, it is worth recalling the main construction of \cite{quasimapcrit}. Consider for concreteness the quasimap section $\Sh{B}_2$; it fits together with its ADHM data into a commutative diagram 
\[
\begin{tikzcd}[ampersand replacement=\&]
        \Sh{V} \arrow[r] \arrow[d, "\Sh{B}_2"] \& F \otimes \Sh{O}_{\QM^{\text{ns $\infty$}}(X) \times \mathbb{P}^1} \arrow[r] \arrow[d, "{\begin{psmallmatrix}
    B_2 && 0 \\ J_3 && \xi_2
\end{psmallmatrix}}"] \& U \otimes \Sh{O}_{\QM^{\text{ns $\infty$}}(X) \times \mathbb{P}^1}(1) \arrow[d, "B_2"] \\
        \Sh{V} \arrow[r] \& F \otimes \Sh{O}_{\QM^{\text{ns $\infty$}}(X) \times \mathbb{P}^1} \arrow[r] \& U \otimes \Sh{O}_{\QM^{\text{ns $\infty$}}(X) \times \mathbb{P}^1}(1).
    \end{tikzcd}
\]
Now working over $U \times \text{Spec} \, \mathbb{C}[t^{-1}]/t^{-2}$, we can use the splitting to ``solve for $\Sh{B}_2$'':
\[
\Sh{B}_2(t) = \begin{pmatrix} \pi^\infty(t) && \rho^\infty(t) \end{pmatrix} \begin{pmatrix} B_2 && 0 \\ J_3 && \xi_2 \end{pmatrix} \begin{pmatrix} \psi^\infty(t) \\ \chi^\infty(t) \end{pmatrix}.
\]
In a similar way one can describe how isomorphisms of bundles $\Sh{V}$ in \eqref{eq: quiver fiber sequence} are related with isomorphisms of the quiver data. By a routine argument with Lemma 2.2 of \cite{quasimapcrit}, we deduce the following. For any bundle isomorphism $g: \Sh{V}' \to \Sh{V}$, there is a unique isomorphism 
\[
\begin{pmatrix}
    g_U && 0 \\ x && g_V
\end{pmatrix}: F' \to F
\]
of quiver data, making the following diagram commute
\[
\begin{tikzcd}[ampersand replacement=\&]
        \Sh{V}' \arrow[r] \arrow[d, "g"] \& F' \otimes \Sh{O}_{\QM^{\text{ns $\infty$}}(X) \times \mathbb{P}^1} \arrow[r, "{\begin{psmallmatrix} B_1' - t & I' \end{psmallmatrix}}"] \arrow[d, "{\begin{psmallmatrix}
    g_U && 0 \\ x && g_V
\end{psmallmatrix}}"] \& U' \otimes \Sh{O}_{\QM^{\text{ns $\infty$}}(X) \times \mathbb{P}^1}(1) \arrow[d, "g_U"] \\
        \Sh{V} \arrow[r] \& F \otimes \Sh{O}_{\QM^{\text{ns $\infty$}}(X) \times \mathbb{P}^1} \arrow[r, "{\begin{psmallmatrix} B_1 - t & I \end{psmallmatrix}}"] \& U \otimes \Sh{O}_{\QM^{\text{ns $\infty$}}(X) \times \mathbb{P}^1}(1).
    \end{tikzcd}
\]
If we now restrict over $\tilde U \times \text{Spec} \, (\mathbb{C}[t^{-1}]/t^{-2})$ and choose trivializations for $\Sh{V}, \Sh{V}'$ as above, the left square of this diagram reads 
\[
\begin{pmatrix} I /t \\ 1\end{pmatrix} g(t) = \begin{pmatrix} g_U && 0 \\ x && g_V \end{pmatrix} \begin{pmatrix} I'/t \\ 1 \end{pmatrix}
\]
so the bottom row says\footnote{One can interpret what we are saying here more formally as specifying $f^\star$ as a functor between the corresponding action groupoids, at least locally.} $g(t) = g_V + xI' t^{-1} + O(t^{-2})$. 

Then if $\tilde U$ is such a trivializing open patch and $\tilde V$ a compatible open patch in $TX$, $f^\star(\tilde U) \subset \tilde V$, the pullback map $$(f^\star)^*: \Gamma(V, \Sh{O}_{TX}) \to \Gamma(U, \Sh{O}_{\QM^{\text{ns $\infty$}}(X)})$$ 
on functions is completely specified, using the dga models of Proposition \ref{prop: dga for Qinfty} and Lemma \ref{lemma: dga for TX}:
\begin{equation} \label{eq: pullback in quiver model}
\begin{split}
    (f^\star)^*(c) & = \zeta I \\
    (f^\star)^*(\bm{\xi}_2(t)) & = \begin{pmatrix} \pi^\infty(t) && \rho^\infty(t)\end{pmatrix} \begin{pmatrix} B_2 && 0 \\ J_3 && \xi_2 \end{pmatrix} \begin{pmatrix} \psi^\infty(t) \\ \chi^\infty(t) \end{pmatrix} - \frac{1}{t} \zeta\eta_{J_2} \\ 
    & = \xi_2 + \frac{1}{t}(\pi_0^\infty(B_2 I - I \xi_2) + J_3I - \zeta \eta_{J_2}) + O(t^{-2}) \\
    (f^\star)^*(\bm{\xi}_3(t)) & = \begin{pmatrix} \pi^\infty(t) && \rho^\infty(t)\end{pmatrix} \begin{pmatrix} B_3 && 0 \\ -J_2 && \xi_3 \end{pmatrix} \begin{pmatrix} \psi^\infty(t) \\ \chi^\infty(t) \end{pmatrix} - \frac{1}{t} \zeta \eta_{J_3}\\
    & = \xi_3 + \frac{1}{t}(\pi_0^\infty(B_3 I - I\xi_3) - J_2 I - \zeta \eta_{J_3}) + O(t^{-2}) \\
    (f^\star)^* (\bm{\iota}(t)) & = \begin{pmatrix} \pi^\infty(t) && \rho^\infty(t)\end{pmatrix} \begin{pmatrix} 0 \\ \iota  \end{pmatrix} + \frac{1}{t} \zeta \eta_\alpha = \iota - \frac{1}{t} (\pi^\infty_0 I \iota - \zeta \eta_\alpha) + O(t^{-2}) \\
    (f^\star)^*(\bm{\jmath}(t)) & = \begin{pmatrix} \alpha && \jmath \end{pmatrix} \begin{pmatrix} \psi^\infty(t) \\ \chi^\infty(t)  \end{pmatrix} = \jmath + \frac{1}{t} \alpha I + O(t^{-2}) \\
    (f^\star)^*(\bm{\chi}(t)) & = \chi + \frac{1}{t} \Big( \comm{\xi_2}{\pi^\infty_0 \eta_{J_3}} + \comm{\pi^\infty_0 \eta_{J_2}}{\xi_3} - \pi^\infty_0 \eta_\alpha \jmath \\
    & + \eta_I I - J_2 \eta_{J_2} - J_3 \eta_{J_3} - \zeta \gamma \Big) + O(t^{-2}). 
\end{split}
\end{equation}
One understands these formulas as follows: the assignments on degree zero generators are determined from the assignment of quasimap data to quiver data from \cite{quasimapcrit}. In degree 1 the map is fixed by the above discussion of the relation of automorphisms in the two models. In degree $-1$ the choice is the minimal one completing $(f^\star)^*$ to a morphism of dgas. Notice that varying the auxiliary splitting data $\pi^\infty_0$ affects $(f^\star)^*$ only by a chain homotopy.

\eqref{eq: pullback in quiver model} also fixes the pullback map on $\text{DR}^\bullet$ over the corresponding patches. Relative to the dga model of Lemma \ref{lemma: dga for TX} and the Cartan model of Proposition \ref{prop: de rham on quotient}, the AKSZ symplectic form on $TX$ is given by 
    \[
    \omega_{TX} = \Res \Big( \tr d \bm{\xi}_2(t) \wedge d\bm{\xi}_3(t) + \tr d \bm{\iota}(t) \wedge d \bm{\jmath}(t) \Big) - \Res \tr \Big( \phi_V - \frac{dc}{t}\Big)d \bm{\chi}(t).
    \]
As usual, all the higher closure data $\omega_{T X, i} = 0$ for $i \geq 1$. 

\begin{proof}[Proof of Theorem \ref{thm: shifted symplectic structure quiver model}]
Now we prove Theorem \ref{thm: shifted symplectic structure quiver model} by analyzing the defining diagram \eqref{eq: QM as Lagrangian intersection} locally over $\tilde U$ relative to the dg-models explicated above. The rationale is as follows: in our model, diagram \eqref{eq: QM as Lagrangian intersection} commutes in the strict sense so $\sigma_i = 0$ in the notation of the discussion below Theorem \ref{thm: PTVV lag intersection}; likewise since our dg model for $\tilde U$ is already cofibrant relative to $TX$ we may take the intersection with $T_\xi$ naively, thus $\vartheta^{(1)}_i = 0$ as well. Then the full shifted symplectic structure is determined by $\vartheta^{(2)}_i$, which we now compute: we will find that $\vartheta^{(2)}_i = 0$ for $i \geq 1$ and that $\vartheta^{(2)}_0$ is given by an explicit formula. Restricting this formula to $\tilde U$, we will find the restriction of $\Omega_{\QM^\xi(X)}$ to $\tilde U$, which will turn out to be enough to conclude the result. 

Throughout this proof, we will use the relations from Proposition \ref{prop: dga for Qinfty} without comment.

Via \eqref{eq: pullback in quiver model} we find 
\begin{equation*}
\begin{split}
(f^\star)^* \omega_{TX} & = \tr d\xi_2 \wedge d( \pi_0^{\infty} (B_3 I - I\xi_3) -J_2 I - \zeta \eta_{J_3}) + \tr d(\pi^\infty_0(B_2 I - I \xi_2) + J_3 I - \zeta \eta_{J_2}) \wedge d \xi_3  \\
& - \tr d(\pi^\infty_0 I \iota - \zeta \eta_\alpha) \wedge d\jmath + \tr d \iota \wedge d(\alpha I) \\
& - \tr \phi_V d( \comm{\xi_2}{\pi_0^\infty \eta_{J_3}} + \comm{\pi^\infty_0 \eta_{J_2}}{\xi_3} - \pi_0^\infty \eta_\alpha \jmath + + \eta_I I - J_2 \eta_{J_2} - J_3 \eta_{J_3} - \zeta \gamma ) + \tr d(\zeta I)  \wedge d\chi. 
\end{split}
\end{equation*}
We will now compute that $(f^\star)^*\omega_{TX} = (\iota_{V(\phi)} + \delta) \vartheta$ for some $\vartheta$ with $d \vartheta = 0$. It is possible to partially sidestep the following by a less direct argument, but we believe it is instructive to ``discover'' $\vartheta$ by hand.  

First consider the terms 
\[
\begin{split}
    \tr d(\zeta \eta_\alpha) \wedge d \jmath + \tr d \iota \wedge d(\alpha I) & = \tr d(\zeta \eta_\alpha) \wedge d\jmath + \tr d \iota (\alpha dI + d\alpha I) \\
    & = \tr d(\zeta \eta_\alpha) \wedge d\jmath + \tr d(\iota \alpha) \wedge dI + \tr (d(I \iota) \wedge d \alpha) \\
    & = - \delta( \tr  d\eta_\alpha  \wedge d\alpha) - \tr d\eta_\alpha \wedge  d(\jmath \zeta) + \tr d(\zeta \eta_\alpha) \wedge d\jmath + \tr d(\iota \alpha) \wedge dI \\
    & = - \delta(\tr d\eta_\alpha \wedge d\alpha) - \tr d\zeta \wedge d(\eta_\alpha \jmath) + \tr d(\iota \alpha) \wedge dI. 
\end{split}
\]
Similarly observe 
\[
\begin{split}
\tr d(J_2 I) \wedge d \xi_2 - \tr d(\zeta \eta_{J_2}) \wedge d \xi_3 & = \tr dJ_2\wedge (d(I \xi_2 - B_2 I) - dI \xi_2  + \tr d(B_2 I)) \\
& - \tr d\xi_2 \wedge J_2 d I  - \tr d(\zeta \eta_{J_2}) \wedge d \xi_3 \\
& =  \tr d(B_2 I - I \xi_2) \wedge dJ_2 - \tr d(\xi_2 J_2) \wedge dI - \tr d(\zeta \eta_{J_2}) \wedge d \xi_3 \\
& + \tr dJ_2 B_2 \wedge dI  + \tr d(IJ_2) \wedge dB_2 - \tr dI \wedge J_2 dB_2 \\
& = - \delta (\tr d \eta_{J_2} \wedge dJ_2)  - \tr d \eta_{J_2} \wedge d(\zeta B_3 - \xi_3 \zeta) - \tr d(\zeta \eta_{J_2}) \wedge d \xi_3 \\ 
& +  \tr d(J_2 B_2 - \xi_2 J_2) \wedge dI + \tr d(IJ_2) \wedge dB_2 \\ 
& = - \delta( \tr d \eta_{J_2} \wedge d J_2) + \tr d \zeta \wedge d(\eta_{J_2} \xi_3) - \tr d \eta_{J_2} \wedge d(\zeta B_3) \\ 
& + \tr d(J_2 B_2 - \xi_2 J_2) \wedge dI + \tr d(IJ_2) \wedge dB_2. 
\end{split}
\]
By a nearly identical consideration, 
\[
\begin{split}
    \tr d(J_3 I) \wedge d \xi_3 + \tr d(\zeta \eta_{J_3})  \wedge d \xi_2 & = - \delta (\tr d \eta_{J_3} \wedge dJ_3) - \tr d \zeta \wedge d(\eta_{J_3} \xi_2) + \tr d \eta_{J_3} \wedge d(\zeta B_2) \\
    & + \tr d(J_3 B_3 - \xi_3 J_3) \wedge dI + \tr d(IJ_3) \wedge dB_3. 
\end{split}
\]
Next note
\[
\begin{split}
    & \tr d(IJ_2) \wedge dB_2 + \tr d(IJ_3) \wedge dB_3 - \tr d \eta_{J_2} \wedge d ( \zeta B_3) + \tr d \eta_{J_3} \wedge d(\zeta B_2) \\
    & = - \delta \Big( \tr d\eta_{B_2} \wedge d B_2 + \tr d \eta_{B_3} \wedge d B_3 \Big) + \tr d(\eta_{J_3} \zeta) \wedge dB_2 - \tr d(\eta_{J_2} \zeta) \wedge dB_3 \\
    & -\tr d \eta_{J_2} \wedge d(\zeta B_3) + \tr d \eta_{J_3} \wedge d(\zeta B_2) - \tr d(\comm{B_3}{B_1}) \wedge dB_2 - \tr d(\comm{B_1}{B_2} ) \wedge dB _3 \\
    & = - \delta \Big( \tr d \eta_{B_2} \wedge dB_2 + \tr d \eta_{B_3} \wedge dB_3 \Big) + \tr d\zeta \wedge d(B_2 \eta_{J_3}) - \tr d \zeta \wedge d(B_3 \eta_{J_2}) + \tr d(\comm{B_2}{B_3}) \wedge dB_1 \\
    & = - \delta \Big(\sum_{i = 1}^3 \tr d \eta_{B_i} \wedge dB_i  \Big) + \tr d\zeta \wedge d(B_2 \eta_{J_3}) - \tr d \zeta \wedge d(B_3 \eta_{J_2}) - \tr d \eta_{B_1} \wedge d(I \zeta) 
\end{split}
\]
where we have used 
\[
0 = d^2 \tr B_1 \comm{B_2}{B_3} = \tr d(\comm{B_2}{B_3}) \wedge dB_1 + \text{cyclic}.
\]
Thus far it has been established that 
\[
\begin{split}
(f^\star)^* \omega_{TX} & = - \delta \Big( \sum_{i = 1}^3 \tr d \eta_{B_i} \wedge d B_i + \sum_{i = 2}^3 \tr d \eta_{J_i} \wedge dJ_i + \tr d \eta_\alpha \wedge d \alpha \\ 
& + \tr d(\pi^\infty_0 \eta_{J_2}) \wedge d\xi_3 - \tr d(\pi_0^\infty \eta_{J_3}) \wedge d\xi_2 - \tr d(\pi^\infty_0 \eta _\alpha) \wedge d\jmath \Big) \\
& - \tr d\zeta \wedge d(B_3 \eta_{J_2} - \eta_{J_2} \xi_3 + \eta_{J_3} \xi_2 - B_2 \eta_{J_3} + \eta_\alpha \jmath) \\
& + \tr d (J_2 B_2 + J_3 B_3 - \xi_2 J_2 - \xi_3 J_3 + \iota \alpha) \wedge dI - \tr d \eta_{B_1} \wedge d(I \zeta) + \tr d(\zeta I) \wedge d \chi \\
& -\tr \phi_Vd( \comm{\xi_2}{\pi^\infty_0 \eta_{J_3}} + \comm{\pi^\infty_0 \eta_{J_2}}{\xi_3} - \pi^\infty_0 \eta_\alpha \jmath + \eta_I I - J_2 \eta_{J_2} - J_3 \eta_{J_3} - \zeta\gamma). 
\end{split}
\]
Notice now 
\[
\begin{split}
& \tr d(J_2 B_2 +J_3 B_3 - \xi_2 J_2 - \xi_3 J_3 + \iota \alpha) \wedge dI - \tr d \eta_{B_1} \wedge d(I \zeta) + \tr d(\zeta I) \wedge d \chi \\
& = - \delta (\tr d \eta_I \wedge dI) + \tr d(\chi \zeta + \zeta \eta_{B_1}) \wedge dI - \tr d \eta_{B_1} \wedge d(I \zeta) + \tr d(\zeta I) \wedge d \chi \\
& = - \delta( \tr d \eta_I \wedge dI)  + \tr d \zeta \wedge ( dI \chi + I d \chi - d \eta_{B_1} I + \eta_{B_1} dI) 
\end{split}
\]
so altogether we have at this point 
\[
(f^\star)^* \omega_{TX} = - \delta \vartheta' - \tr \phi_V d( \comm{\xi_2}{\pi^\infty_0 \eta_{J_3}} + \comm{\pi^\infty_0 \eta_{J_2}}{\xi_3} - \pi^\infty_0 \eta_\alpha \jmath + \eta_I I - J_2 \eta_{J_2} - J_3 \eta_{J_3} - \zeta\gamma)
\]
with 
\[
\begin{split}
\vartheta' & = \sum_{i = 1}^3 \tr d \eta_{B_i} \wedge dB_i + \sum_{i = 2}^3 \tr d \eta_{J_i} \wedge dJ_i+  \tr d\eta_I \wedge dI + \tr d \eta_\alpha \wedge d\alpha - \tr d \zeta \wedge d \gamma \\
& + \tr d(\pi^\infty_0 \eta_{J_2}) \wedge d\xi_3 - \tr d(\pi^\infty_0 \eta_{J_3}) \wedge d \xi_2 - \tr d( \pi^\infty_0 \eta_\alpha) \wedge d\jmath. 
\end{split}
\]
As the reductive part of the group in the quotient description of Proposition \ref{prop: dga for Qinfty} is $G^{\text{red}} = GL(U) \times GL(V)$ in the notation of \eqref{eq: reductive and unipotent}, introduce $\phi^{\text{red}} = (\phi_U, \phi_V)$ as notation for the equivariant variables in the Cartan model of Proposition \ref{prop: de rham on quotient}. We compute at once 
\[
\begin{split}
\iota_{V(\phi)} \vartheta' & = \sum_{i = 1}^3 (\tr \comm{\eta_{B_i}}{\phi_U} d B_i + \tr d\eta_{B_i} \comm{B_i}{\phi_U}) + \sum_{i = 2}^3 \tr (\eta_{J_i} \phi_V - \phi_U \eta_{J_i}) dJ_i + \tr d \eta_{J_i} (J_i \phi_U - \phi_VJ_i) \\
& + \tr(\eta_I \phi_U - \phi_V \eta_I) dI + \tr d \eta_I (I \phi_V - \phi_U I) - \tr \phi_U \eta_\alpha d \alpha + \tr d \eta_\alpha \alpha \phi_U - \tr( \zeta \phi_U - \phi_V \zeta) d \gamma \\
& - \tr d \zeta (\gamma \phi_V - \phi_U \gamma) + \tr \comm{\pi^\infty_0 \eta_{J_2}}{\phi_V} d\xi_3 + \tr d(\pi_0^\infty \eta_{J_2}) \comm{\xi_3}{\phi_V} \\
& - \tr \comm{\pi^\infty_0 \eta_{J_3}}{\phi_V} d\xi_2 - \tr d(\pi^\infty_0 \eta_{J_3}) \comm{\xi_2}{\phi_V} + \tr \phi_V \pi^\infty_0 \eta_\alpha d \jmath - \tr d(\pi^\infty_0 \eta_\alpha) \jmath \phi_V \\
& =  \tr \phi_U d\Big( \sum_{i = 1}^3 \comm{\eta_{B_i}}{B_i} + \sum_{i = 2}^3 \eta_{J_i} J_i - I \eta_I + \eta_\alpha \alpha + \gamma \zeta \Big) \\ 
& + \tr \phi_V d\Big( \eta_I I - \sum_{i = 2}^3 J_i \eta_{J_i} - \zeta \gamma + \comm{\pi^\infty_0 \eta_{J_2}}{\xi_3} + \comm{\xi_2}{\pi^\infty_0 \eta_{J_3}} - \pi^\infty_0 \eta_\alpha \jmath \Big). 
\end{split}
\]
Finally noticing that $\iota_{V(\phi)} \tr \phi_U d\beta = \tr \phi_U \comm{\beta}{\phi_U} = 0$ we find 
\[
(f^\star)^* \omega_{TX} = -(\iota_{V(\phi)} + \delta) \vartheta
\]
with 
\[
\begin{split}
\vartheta & = \sum_{i = 1}^3 \tr d \eta_{B_i} \wedge dB_i + \sum_{i = 2}^3 \tr d \eta_{J_i} \wedge dJ_i+  \tr d\eta_I \wedge dI + \tr d \eta_\alpha \wedge d\alpha \\
& + \tr \phi_U d \beta - \tr d \zeta \wedge d \gamma \\
& + \tr d(\pi^\infty_0 \eta_{J_2}) \wedge d\xi_3 - \tr d(\pi^\infty_0 \eta_{J_3}) \wedge d \xi_2 - \tr d( \pi^\infty_0 \eta_\alpha) \wedge d\jmath. 
\end{split}
\]
Pulling this back under\footnote{Notice the pullback does not depend on the splitting $\pi^\infty_0$ defined locally on $\tilde U$, and equals the restriction of $\Omega_{\Sh{Q}^\xi}$ to $\tilde U$---in particular, there is no nontrivial coherence data involved in gluing from the local $\tilde U$'s.} $\QM^\xi(X) \to \QM^{\text{ns $\infty$}}(X)$, comparing with Lemma \ref{lemma: shifted symplectic form quiver model} we conclude Theorem \ref{thm: shifted symplectic structure quiver model}. 

\end{proof}

\subsection{$\Omega_{\QM^\xi(X)}$ in the Čech model} \label{subsec: prove Čech model symplectic thm}
Now we will repeat this argument in the Čech model. The qualitative structure and outline of the argument is nearly identical, so we will be more brief and only describe the new features that arise. In this section, we will repeatedly use the notation and action of the differential described in Proposition \ref{prop: dga for Qinfty in the Čech model} without comment. 

The morphism $f^\star: \QM^{\text{ns $\infty$}}(X) \to TX $ induces a pullback on function dgas and de Rham complexes, fully determined by the assignments (notations as in Lemma \ref{lemma: dga for TX} and Proposition \ref{prop: dga for Qinfty in the Čech model}):
\begin{equation*} 
\begin{split}
    (f^\star)^*(c) & = \Res_u \zeta^\infty(u) \\ 
    (f^\star)^*(\bm{\xi}_2(t)) & = \xi_2 + \frac{1}{t} \Big( \Res_u( \Sh{B}_2(u) - \comm{\eta_{3, \infty}(u)}{\zeta^\infty(u)}) \Big) \\ 
    (f^\star)^*(\bm{\xi}_3(t)) & = \xi_3 + \frac{1}{t} \Res_u(\Sh{B}_{3, \infty}(u)) \\ 
    (f^\star)^*(\bm{\iota}(t)) & = \iota + \frac{1}{t} \Big( \Res_u(\Sh{I}(u) - \zeta^\infty(u) \eta_{\Sh{J}_\infty}(u)) \Big) \\ 
    (f^\star)^*(\bm{\jmath}(t)) & = \jmath + \frac{1}{t} \Res_u(\Sh{J}_\infty(u)) \\ 
    (f^\star)^*(\bm{\chi}(t)) & = \chi + \frac{1}{t} \Big( \Res_u(\eta_g(u) - \comm{\eta_{3, \infty}(u)}{\Sh{B}_{3, \infty}(u) - \xi_3} - \eta_{\Sh{J}_\infty}(u)(\Sh{J}_\infty(u) - \jmath) - \comm{\beta_\infty(u)}{\zeta^\infty(u)}) \Big).  
\end{split}
\end{equation*}
It is instructive to check directly that $(f^\star)^*$ intertwines the respective differentials. 

\begin{proof}[Proof of Theorem \ref{thm: shifted symplectic structure Čech model}]
We follow the same strategy explicated at the beginning of the proof of Theorem \ref{thm: shifted symplectic structure quiver model}. In what follows we will use the relations from Proposition \ref{prop: dga for Qinfty in the Čech model} without comment.

The pullback of $\omega_{TX}$ now reads 
\begin{equation} \label{eq: pullback AKSZ in Čech model}
\begin{split}
(f^\star)^* \omega_{TX} & = \Res_u \Big( \tr d \xi_2 \wedge d(\Sh{B}_{3, \infty}(u)) + \tr d(\Sh{B}_2(u) - \comm{\eta_{3, \infty}(u)}{\zeta^\infty(u)}) \wedge d\xi_3 \\
& + \tr d \iota \wedge d(\Sh{J}_\infty(u)) + \tr d(\Sh{I}(u) - \zeta^\infty(u) \eta_{\Sh{J}_\infty}(u)) \wedge d\jmath \\ 
& - \tr \phi_\infty d(\eta_g(u) - \comm{\eta_{3, \infty}(u)}{\Sh{B}_{3, \infty}(u) - \xi_3} - \eta_{\Sh{J}_\infty}(u) (\Sh{J}_\infty(u) - \jmath) - \comm{\beta_\infty(u)}{\zeta^\infty(u)}) \\
& + \tr d(\zeta^\infty(u)) \wedge d \chi \Big) 
\end{split}
\end{equation}
We follow a similar procedure as before. First observe that 
\[
\begin{split}
& \Res \tr d \xi_2 \wedge d(\Sh{B}_{3, \infty}(u)) + \Res \tr d(\Sh{B}_2(u) - \comm{\eta_{3, \infty}(u)}{\zeta^\infty(u)}) \wedge d \xi_3 \\
& = \Res \tr d\Sh{B}_2(u) \wedge d \Sh{B}_{3, \infty}(u) - \Res \tr d(\Sh{B}_2(u) - \xi_2) \wedge d(\Sh{B}_{3, \infty}(u) - \xi_3) - \Res \tr d(\comm{\eta_{3, \infty}(u)}{\zeta^\infty(u)}) \wedge d \xi_3
\end{split}
\]
Now 
\[
- \Res \tr d(\Sh{B}_2(u) - \xi_2) \wedge d(\Sh{B}_{3, \infty}(u) - \xi_3) = -\Res \tr d[ \Sh{B}_2(u) - \xi_2]_{\geq 0} \wedge d (\Sh{B}_{3, \infty}(u) - \xi_3)
\]
because $O(u^{-2})$ terms do not contribute to residue, and we then have 
\[
\begin{split}
& - \Res \tr d[\Sh{B}_2(u) - \xi_2]_{\geq 0} \wedge d (\Sh{B}_{3, \infty}(u) - \xi_3) = \\
& \Res \tr \delta(d \eta_{3, \infty}(u)) \wedge d(\Sh{B}_{3, \infty}(u) - \xi_3) - \Res \tr d(\comm{\eta_{3, \infty}(u)}{\zeta^\infty(u)}) \wedge d(\Sh{B}_{3, \infty}(u) - \xi_3) \\
& = \delta \Big( \Res \tr d \eta_{3, \infty}(u) \wedge d(\Sh{B}_{3, \infty}(u) - \xi_3) \Big) + \Res \tr d \eta_{3, \infty}(u) \wedge d(\comm{\Sh{B}_{3, \infty}(u)}{\zeta^\infty(u)}) \\
& - \Res \tr d(\comm{\eta_{3, \infty}(u)}{\zeta^\infty(u)}) \wedge d(\Sh{B}_{3, \infty}(u) - \xi_3) \\
& = \delta \Big( \Res \tr d \eta_{3, \infty}(u) \wedge d(\Sh{B}_{3, \infty}(u) - \xi_3) \Big) + \Res \tr d\zeta^\infty(u) \wedge d(\comm{\eta_{3, \infty}(u)}{\Sh{B}_{3, \infty}(u)}) \\
& + \Res \tr d(\comm{\eta_{3, \infty}(u)}{\zeta^\infty(u)}) \wedge d \xi_3. 
\end{split}
\]
Altogether this shows 
\[
\begin{split}
& \Res \tr d \xi_2 \wedge d(\Sh{B}_{3, \infty}(u)) + \Res \tr d(\Sh{B}_2(u) - \comm{\eta_{3, \infty}(u)}{\zeta^\infty(u)}) \wedge d \xi_3 \\
& = \delta \Big( \Res \tr d \eta_{3, \infty}(u) \wedge d(\Sh{B}_{3, \infty}(u) - \xi_3) \Big) + \Res \tr d\zeta^\infty(u) \wedge d(\comm{\eta_{3, \infty}(u)}{\Sh{B}_{3, \infty}(u)}) \\
& + \Res \tr d\Sh{B}_2(u) \wedge d \Sh{B}_{3, \infty}(u)
\end{split}
\]
so it remains to deal with the third term. For this we manipulate as follows 
\begin{equation}
\begin{split} \label{eq: formula tr dB2dB3}
    & \Res \tr d \Sh{B}_2(u) \wedge d \Sh{B}_{3, \infty}(u) \\
    & = d \Big(  \Res \tr\Sh{B}_2(u) d(g^{-1}(u) \Sh{B}_{3, 0}(u) g(u))+ \Res \tr \Sh{B}_2(u) d(\delta \eta_{2}(u)) + \Res \tr \Sh{B}_2(u) d(\comm{\eta_2(u)}{\zeta^\infty(u)}) \Big) \\
    & = d\Big( \Res \tr [g(u) \Sh{B}_2(u) g^{-1}(u)]_{\leq -1} d \Sh{B}_{3, 0}(u) - \Res \tr \Sh{B}_2(u) \comm{g^{-1}(u)dg(u)}{g^{-1}(u) \Sh{B}_{3, 0}(u) g(u)} \\
    & + \Res \tr \Sh{B}_2(u) d(\comm{\eta_2(u)}{\zeta^\infty(u)}) \Big) + \Res \tr \delta( d \eta_2(u)) \wedge d\Sh{B}_2(u). 
\end{split}
\end{equation}
It will be convenient to deal with the first two terms in parenthesis one at a time. The first, simpler case is 
\[
\Res \tr [g(u) \Sh{B}_2(u) g^{-1}(u)]_{\leq -1} d\Sh{B}_{3, 0}(u) = \Res \tr \delta \eta_{3, 0}(u) d\Sh{B}_{3, 0}(u) + \Res \tr \comm{\eta_{3, 0}(u)}{\zeta^0(u)} d\Sh{B}_{3, 0}(u)
\]
and therefore 
\[
\begin{split}
& d \Big( \Res \tr [g(u) \Sh{B}_2(u) g^{-1}(u)]_{\leq -1} d\Sh{B}_{3, 0}(u) \Big) = - \delta \Big( \Res \tr d \eta_{3, 0}(u) \wedge d \Sh{B}_{3, 0}(u) \Big) \\
&   - \Res \tr d \eta_{3, 0}(u) \wedge d(\comm{\Sh{B}_{3, 0}(u)}{\zeta^0(u)}) + \Res \tr d( \comm{\eta_{3, 0}(u)}{\zeta^0(u)}) \wedge d \Sh{B}_{3, 0}(u) \\
& = - \delta \Big( \Res \tr d \eta_{3, 0}(u) \wedge d \Sh{B}_{3, 0}(u) \Big) - \Res \tr d \zeta^0(u) \wedge d (\comm{\eta_{3, 0}(u)}{\Sh{B}_{3, 0}(u)}).
\end{split}
\]
Next we have 
\[
\begin{split}
    & -\Res \tr \Sh{B}_2(u) \comm{g^{-1}dg(u)}{g^{-1}(u) \Sh{B}_{3, 0}(u) g(u)} =  \Res \tr \comm{\Sh{B}_2(u)}{g^{-1}(u) \Sh{B}_{3, 0}(u) g(u)} g^{-1}dg(u) \\
    & = -\Res \tr \comm{\Sh{B}_2(u)}{\delta \eta_2(u)} g^{-1}dg(u) - \Res \tr \comm{\Sh{B}_2(u)}{\comm{\eta_2(u)}{\zeta^\infty(u)}} g^{-1}dg(u) \\
    & + \Res \tr \comm{\Sh{B}_2(u)}{\Sh{B}_{3, \infty}(u)} g^{-1}dg(u) . 
\end{split}
\]
Moreover 
\[
\begin{split}
& \delta (g^{-1}dg(u) ) \\
& = -g^{-1}(u)(-\zeta^0(u) g(u) + g(u) \zeta^\infty(u)) g^{-1}dg(u) -g^{-1}d(- \zeta^0(u) g(u) + g(u) \zeta^\infty(u))  \\
& = - d\zeta^\infty(u) + g^{-1}(u) d\zeta^0(u) g(u) - \comm{\zeta^\infty(u)}{g^{-1}dg(u)}
\end{split}
\]
from which we find 
\[
\begin{split}
& \delta \Big( \Res \tr \comm{\Sh{B}_2(u)}{\eta_2(u)} g^{-1}dg(u) \Big) \\
& = \Res \tr \comm{\comm{\Sh{B}_2(u)}{\zeta^\infty(u)}}{\eta_2(u)} g^{-1}dg(u)  + \Res \tr \comm{\Sh{B}_2(u)}{\delta \eta_2(u)} g^{-1}dg(u) \\
& + \Res \tr \comm{\Sh{B}_2(u)}{\eta_2(u)} d\zeta^\infty(u) + \Res \tr \comm{\Sh{B}_2(u)}{\eta_2(u)} \comm{\zeta^\infty(u)}{g^{-1}dg(u)} \\
& - \Res \tr \comm{\Sh{B}_2(u)}{\eta_2(u)} g^{-1}(u) d\zeta^0(u) g(u) \\
& = \Res \tr \comm{\Sh{B}_2(u)}{\delta \eta_2(u)} g^{-1}dg(u) + \Res \tr \comm{\Sh{B}_2(u)}{\comm{\eta_2(u)}{\zeta^\infty(u)}} g^{-1}dg(u) \\
& - \Res \tr \comm{\eta_2(u)}{\Sh{B}_2(u)} d\zeta^\infty(u) + \Res \tr g(u) \comm{\eta_2(u)}{\Sh{B}_2(u)}g^{-1}(u) d\zeta^0(u) 
\end{split}
\]
where we have used identity 
\[
\begin{split}
0 & = \Res \tr \comm{\zeta^\infty(u)}{\comm{\Sh{B}_2(u)}{\eta_2(u)} g^{-1}dg(u)} \\
& = \Res \tr \comm{\zeta^\infty(u)}{\comm{\Sh{B}_2(u)}{\eta_2(u)}}g^{-1}dg(u)  - \Res \tr \comm{\Sh{B}_2(u)}{\eta_2(u)} \comm{\zeta^\infty(u)}{g^{-1}dg(u)} \\
& = -\Res \tr \comm{\comm{\Sh{B}_2(u)}{\zeta^\infty(u)}}{\eta_2(u)} g^{-1}dg(u)  - \Res \tr \comm{\Sh{B}_2(u)}{\eta_2(u)} \comm{\zeta^\infty(u)}{g^{-1}dg(u)} \\
& + \Res \tr \comm{\Sh{B}_2(u)}{\comm{\eta_2(u)}{\zeta^\infty(u)}} g^{-1}dg(u).
\end{split}
\]
It follows that 
\[
\begin{split}
& -\Res \tr \Sh{B}_2(u) \comm{g^{-1}dg(u)}{g^{-1}(u) \Sh{B}_{3, 0}(u) g(u)}  \\
& = \delta \Big( \Res \tr \comm{\eta_2(u)}{\Sh{B}_2(u)} g^{-1}dg(u) \Big) - \Res \tr \comm{\eta_2(u)}{\Sh{B}_2(u)} d\zeta^\infty(u) \\
& + \Res \tr g(u) \comm{\eta_2(u)}{\Sh{B}_2(u)} g^{-1}(u) d\zeta^0(u) + \Res \tr \comm{\Sh{B}_2(u)}{\Sh{B}_{3, \infty}(u)} g^{-1}dg(u).
\end{split}
\]
From the observation
\[
\begin{split}
& \Res \tr \delta( d \eta_2(u)) \wedge d \Sh{B}_2(u) \\
& = \delta \Big( \Res \tr d \eta_2(u) \wedge d \Sh{B}_2(u) \Big) + \Res \tr d \eta_2(u) \wedge d( \comm{\Sh{B}_2(u)}{\zeta^\infty(u)})
\end{split}
\]
together with the final line of \eqref{eq: formula tr dB2dB3} and subsequent computations, we may conclude 
\[
\begin{split}
    & \Res \tr d\Sh{B}_2(u) \wedge d \Sh{B}_{3, \infty}(u) \\
    & = \delta \Big( \Res \tr d \eta_2(u) \wedge d \Sh{B}_2(u) - \Res \tr d \eta_{3, 0}(u) \wedge d \Sh{B}_{3, 0}(u) - d\Res \tr \comm{\eta_2(u)}{\Sh{B}_2(u)} g^{-1}dg(u)\Big) \\
    & - \Res \tr d \zeta^0(u) \wedge d(\comm{\eta_{3, 0}(u)}{\Sh{B}_{3, 0}(u)} - g(u) \comm{\eta_2(u)}{\Sh{B}_2(u)}g^{-1}(u)) + d\Res \tr \comm{\Sh{B}_2(u)}{\Sh{B}_{3, \infty}(u)} g^{-1}dg(u).
\end{split}
\]

Our considerations have established 
\[
\begin{split}
& \Res \tr d \xi_2 \wedge d(\Sh{B}_{3, \infty}(u)) + \Res \tr d(\Sh{B}_2(u) - \comm{\eta_{3, \infty}(u)}{\zeta^\infty(u)}) \wedge d \xi_3 \\
& = \delta \Big( \Res \tr d \eta_2(u) \wedge d \Sh{B}_2(u) +  \Res \tr d \eta_{3, \infty}(u) \wedge d(\Sh{B}_{3, \infty}(u) - \xi_3) - \Res \tr d \eta_{3, 0}(u) \wedge d \Sh{B}_{3, 0}(u) \\
& - d \Res \tr \comm{\eta_2(u)}{\Sh{B}_2(u)} g^{-1}dg(u) \Big) + \Res \tr d\zeta^\infty(u) \wedge d(\comm{\eta_{3, \infty}(u)}{\Sh{B}_{3, \infty}(u)}) \\
& - \Res \tr d \zeta^0(u) \wedge d(\comm{\eta_{3, 0}(u)}{\Sh{B}_{3, 0}(u)} - g(u) \comm{\eta_2(u)}{\Sh{B}_2(u)}g^{-1}(u)) + d(\Res \tr \comm{\Sh{B}_2(u)}{\Sh{B}_{3, \infty}(u)} g^{-1}dg(u)).
\end{split}
\]
Performing the same set of moves on the other terms in \eqref{eq: pullback AKSZ in Čech model} gives 
\[
\begin{split}
& \Res \tr d \iota + d \Sh{J}_\infty(u) + \Res \tr d(\Sh{I}(u) - \zeta^\infty(u) \eta_{\Sh{J}_\infty}(u)) \wedge d \jmath \\
& = \delta \Big( \Res \tr d \eta_{\Sh{I}}(u) \wedge d \Sh{I}(u) + \Res \tr d \eta_{\Sh{J}_\infty}(u) \wedge d (\Sh{J}_\infty(u) - \jmath) - \Res \tr d \eta_{\Sh{J}_0}(u) \wedge d \Sh{J}_0(u) \\
& + d \Res \tr \Sh{I}(u) \eta_{\Sh{I}}(u) g^{-1}dg(u) \Big) + \Res \tr d \zeta^\infty(u) \wedge d(\eta_{\Sh{J}_\infty}(u) \Sh{J}_\infty(u)) \\
& - \Res \tr d\zeta^0(u) \wedge d(\eta_{\Sh{J}_0}(u) \Sh{J}_0(u) + g(u) \Sh{I}(u) \eta_{\Sh{I}}(u) g^{-1}(u)) + d( \Res \tr \Sh{I}(u) \Sh{J}_\infty(u) g^{-1}dg(u)).
\end{split}
\]
Thus far we have expressed \eqref{eq: pullback AKSZ in Čech model} in the form 
\[
\begin{split}
(f^\star)^* \omega_{TX} & = \delta( \dots) + \Res \tr d \zeta^\infty(u) \wedge d( \comm{\eta_{3, \infty}(u)}{\Sh{B}_{3, \infty}(u)} + \eta_{\Sh{J}_\infty}(u) \Sh{J}_\infty(u) + \chi ) \\
& + \Res \tr d\zeta^0(u) \wedge d(g(u)( \comm{\eta_2(u)}{\Sh{B}_2(u)} - \Sh{I}(u) \eta_{\Sh{I}}(u)) g^{-1}(u) - \comm{\eta_{3, 0}(u)}{\Sh{B}_{3, 0}(u)} - \eta_{\Sh{J}_0}(u) \Sh{J}_0(u)) \\
& + d( \Res \tr( (\comm{\Sh{B}_2(u)}{\Sh{B}_{3, \infty}(u)} + \Sh{I}(u) \Sh{J}_\infty(u)) g^{-1}dg(u) )) \\
& - \Res \tr \phi_\infty d(\eta_g(u) - \comm{\eta_{3, \infty}(u)}{\Sh{B}_{3, \infty}(u) - \xi_3} - \eta_{\Sh{J}_\infty}(u) (\Sh{J}_\infty(u) - \jmath) - \comm{\beta_\infty(u)}{\zeta^\infty(u)}).
\end{split}
\]
Finally, use that 
\[
\begin{split}
    & \Res \tr (\comm{\Sh{B}_2(u)}{\Sh{B}_{3, \infty}(u)} + \Sh{I}(u) \Sh{J}_\infty(u)) g^{-1}dg(u) \\
    & = \Res \tr \delta \eta_g(u) g^{-1}dg(u) + \Res \tr \comm{\eta_g(u)}{\zeta^\infty(u)} g^{-1}dg(u) \\
    & = \delta( \Res \tr \eta_g(u) g^{-1}dg(u)) - \Res \tr \eta_g(u) (d \zeta^\infty(u) + \comm{\zeta^\infty(u)}{g^{-1}dg(u)} - g^{-1}(u) d\zeta^0(u) g(u)) \\
    & + \Res \tr \comm{\eta_g(u)}{\zeta^\infty(u)} g^{-1}dg(u)  \\
    & = \delta( \Res \tr \eta_g(u) g^{-1}dg(u)) - \Res \tr d \zeta^\infty(u) \eta_g(u) + \Res \tr d\zeta^0(u) (g(u) \eta_g(u) g^{-1}(u))
\end{split}
\]
so that in particular 
\[
\begin{split}
    (f^\star)^* \omega_{TX} & = \delta( \dots) + \Res \tr d \zeta^\infty(u) \wedge d( \comm{\eta_{3, \infty}(u)}{\Sh{B}_{3, \infty}(u)} + \eta_{\Sh{J}_\infty}(u) \Sh{J}_\infty(u) - \eta_g(u) + \chi) \\
    & + \Res \tr d \zeta^0(u) \wedge d( g(u)( \comm{\eta_2(u)}{\Sh{B}_2(u)} - \Sh{I}(u) \eta_{\Sh{I}}(u) + \eta_g(u))g^{-1}(u) - \comm{\eta_{3, 0}(u)}{\Sh{B}_{3, 0}(u)} \\
    & - \eta_{\Sh{J}_0}(u) \Sh{J}_0(u)) \\
    & - \Res \tr \phi_\infty d(\eta_g(u) - \comm{\eta_{3, \infty}(u)}{\Sh{B}_{3, \infty}(u) - \xi_3} - \eta_{\Sh{J}_\infty}(u) (\Sh{J}_\infty(u) - \jmath) - \comm{\beta_\infty(u)}{\zeta^\infty(u)})
\end{split}
\]
and at last we have shown 
\[
(f^\star)^* \omega_{TX} = \delta \vartheta' - \Res \tr \phi_\infty d(\eta_g(u) - \comm{\eta_{3, \infty}(u)}{\Sh{B}_{3, \infty}(u) - \xi_3} - \eta_{\Sh{J}_\infty}(u) (\Sh{J}_\infty(u) - \jmath) - \comm{\beta_\infty(u)}{\zeta^\infty(u)})
\]
with 
\[
\begin{split}
\vartheta' & = \Res \tr d \eta_2(u) \wedge d \Sh{B}_2(u) + \Res \tr d \eta_{\Sh{I}}(u) \wedge d \Sh{I}(u) + \Res \tr d \eta_{3, \infty}(u) \wedge d (\Sh{B}_{3, \infty}(u) - \xi_3) \\
& + \Res \tr d \eta_{\Sh{J}_\infty}(u) \wedge d (\Sh{J}_\infty(u) - \jmath)  - \Res \tr d \eta_{3, 0}(u) \wedge d \Sh{B}_{3, 0}(u) - \Res \tr d \eta_{\Sh{J}_0}(u) \wedge d \Sh{J}_0(u) \\
& - d \Res \tr \Big( \eta_g(u) + \comm{\eta_2(u)}{\Sh{B}_2(u)} - \Sh{I}(u) \eta_{\Sh{I}}(u)    \Big) g^{-1}dg(u)  \\
& - \Res \tr d \zeta^\infty(u) \wedge d \beta_\infty(u) - \Res \tr d \zeta^0(u) \wedge d \beta_0(u). 
\end{split}
\]
Write $\phi^\text{red} = (\phi_0, \phi_\infty)$ as the equivariant variables in the Cartan model associated to the quotient $G[t^{-1}] \times G(\Sh{O}) \to G \times G$ given by evaluation at $(\infty, 0)$. Compute at once 
\[
\begin{split}
\iota_{V(\phi)} \vartheta' & = \Res \tr \phi_\infty d\Big( \comm{\eta_{3, \infty}(u)}{\Sh{B}_{3, \infty}(u) - \xi_3} + \eta_{\Sh{J}_\infty}(u) (\Sh{J}_\infty(u) - \jmath) - \eta_g(u) + \comm{\beta_\infty(u)}{\zeta^\infty(u)} \Big) \\
& + \Res \tr \phi_0 d \Big( g(u)( \eta_g(u) + \comm{\eta_2(u)}{\Sh{B}_2(u)} - \Sh{I}(u) \eta_{\Sh{I}}(u))g^{-1}(u) -\comm{\eta_{3, 0}(u)}{\Sh{B}_{3, 0}(u)} \\
& - \eta_{\Sh{J}_0}(u) \Sh{J}_0(u) + \comm{\beta_0(u)}{\zeta^0(u)} \Big) 
\end{split}
\]
where we used that, due to the $G \times G$-invariance, 
\[
\begin{split}
&  \iota_{V(\phi)} d \Res \tr \Big( \eta_g(u) + \comm{\eta_2(u)}{\Sh{B}_2(u)} - \Sh{I}(u) \eta_{\Sh{I}}(u) \Big) g^{-1}dg(u) \\
& = - d\iota_{V(\phi)} \Res \tr \Big( \eta_g(u) + \comm{\eta_2(u)}{\Sh{B}_2(u)} - \Sh{I}(u) \eta_{\Sh{I}}(u) \Big) g^{-1}dg(u) \\
& = d \Res \tr \Big( \eta_g(u) + \comm{\eta_2(u)}{\Sh{B}_2(u)} - \Sh{I}(u) \eta_{\Sh{I}}(u) \Big) (\phi_\infty - g^{-1}(u) \phi_0g(u)). 
\end{split}
\]
We conclude that $(f^\star)^* \omega_{TX} = (\delta + \iota_{V(\phi)}) \vartheta$, with 
\[
\begin{split}
\vartheta & = \Res \tr d \eta_2(u) \wedge d \Sh{B}_2(u) + \Res \tr d \eta_{\Sh{I}}(u) \wedge d \Sh{I}(u) + \Res \tr d \eta_{3, \infty}(u) \wedge d (\Sh{B}_{3, \infty}(u) - \xi_3) \\
& + \Res \tr d \eta_{\Sh{J}_\infty}(u) \wedge d (\Sh{J}_\infty(u) - \jmath)  - \Res \tr d \eta_{3, 0}(u) \wedge d \Sh{B}_{3, 0}(u) - \Res \tr d \eta_{\Sh{J}_0}(u) \wedge d \Sh{J}_0(u) \\
& - d \Res \tr \Big( \eta_g(u) + \comm{\eta_2(u)}{\Sh{B}_2(u)} - \Sh{I}(u) \eta_{\Sh{I}}(u)    \Big) g^{-1}dg(u)  \\
& + \Res \tr \phi_0 d \beta_0 - \Res \tr d \zeta^\infty(u) \wedge d \beta_\infty(u) - \Res \tr d \zeta^0(u) \wedge d \beta_0(u). 
\end{split}
\]
Pulling this back to $\QM^\xi(X)$ and comparing with Lemma \ref{lemma: shifted symplectic form Čech model} we conclude the result. 
\end{proof}


\newpage

\addcontentsline{toc}{part}{Appendices}

\begin{appendix} 

\section{Sheaf operations and critical cohomology}
\label{app: Sheaf operations and critical cohomology}

\subsection{Vanishing cycles and critical cohomology}
\label{app: Appendix VC}

Let $X$ be a smooth Artin stack with a function $f:X\to \BoC$. Set $X_0=f^{-1}(0)$ and let $c:
\widetilde{\BoC^*}\to \BoC^*$ be the universal cover. 
Consider the Cartesian diagram 
\[
\begin{tikzcd}
    & \widetilde{X}\arrow[r]\arrow[dd, "\tilde c"] &  \widetilde{\BoC*}\arrow[d, "c"]\\
    & & \BoC^*\arrow[d, hookrightarrow]\\
    X_0 \arrow[r, hookrightarrow, "i"]& X\arrow[r, "f"] & \BoC
\end{tikzcd}
\]
Given a constructible complex $\Sh{F}$ on $X$, is nearby cycle $\psi_f\Sh{F}$ id defined as $\psi_f\Sh{F}\coloneqq i^*\tilde c_*\tilde c^*$.
The standard adjunction gives a canonical mporphims $i^*\Sh{F}\to \psi_f\Sh{F}$, and the vanishing cycle $\varphi_f(F)$ is defined as the cone $i^*\Sh{F}\to \psi_f\Sh{F}\to \varphi_f{F}$. These constructions are functorial, so we get functors
\[
\psi_f: D^+(X)\to D^+(X_0)\qquad \varphi_f  : D^+(X)\to  D^+(X_0)
\]
It is convenient to define $\phip{f}\Sh{F}\coloneqq\varphi_f\Sh{F}[-1]$
as this functor is
exact with respect to the perverse t-structure on the derived category of constructible complexes. A distinctive property of the vanishing cycle $\phip{f}$ is its commutation with Verdier duality : $\phip{f} \DD\simeq \DD \phip{f}$. We also have $\phip{f=0}=\id$.

Given a morphism $h: X\to Y$ and a function $f: Y\to \BoC$, one has the following natural transformations:
\begin{align}
\begin{split}
\label{eq: base change VC}
    &h_! \phip{f\circ h}\to \varphi_{f}  h_!\qquad \phip{f}  h_* \to  h_* \phip{f\circ h}
    \\
    & h^* \phip{f}\to \phip{f\circ h}  h^* \qquad  \phip{f\circ h}  h^!\to h^* \phip{f}
    \end{split}
\end{align}
If $h$ is smooth (resp. proper and representable), then $h^* \phip{f}\to \phip{f\circ h}  h^*$ (resp. $\phip{f}  h_* \to  h_* \phip{f\circ h}$) is an isomorphism. 

Given a  (possibly non-smooth) stack $X$ and a function $f: X\to \BoC$, we define the (Borel-Moore) critical cohomology of the pair $(X,f)$ as
\[
\HO(X, \phip{f}\DD\BoQ).
\]
Given the action of a torus $T$ on a scheme $X$ and a $T$-invariant function on $X$, we will denote by $\HO_T(X, \phip{f}\DD\BoQ)$ the derived global section of the sheaf $\phip{f}\DD\BoQ_{X/T}$.

Critical cohomology is functorial for smooth pullback and representable proper pushforward. Assume that $h:X\to Y$ is proper representable. Then applying the vanishing cycle functor $\phip{f\circ h}$ to the counit $h_!h^!=h_*h^!\to \id$, using $ h_{*}\phip{f\circ h} \cong  \phip{f}  h_*$ and applying $\DD\BoQ_Y$ to the right we get a canonical pushforward morphism 
\[
h_*: \HO(X, \phip{f\circ h}\DD\BoQ_X)\to \HO(Y, \phip{f}\DD\BoQ_Y).
\]
Similarly, if $h$ is smooth, then applying $\phip{f\circ h}$ to the unit $\id\to h_*h^*$ and composing with the morphism the natural transformation $\phip{f}  h_*\to h_*\phip{f\circ h}$ we get a morphism of complexes
\[
\phip{f}\DD\BoQ_Y \to \phip{f}  h_*h^*\DD\BoQ_Y \xrightarrow{} h_*\phip{f\circ h} h^*\DD\BoQ_Y \cong  h_*\phip{f\circ h} \DD\BoQ_X[2\dim(Y)-2\dim(X)]
\]
and thus a pullback morphism 
\[
h^*: \HO(Y, \phip{f}\DD\BoQ_X)\to \HO(X, \phip{f\circ h}\DD\BoQ_X[2\dim(Y)-2\dim(X)]).
\]
We conclude this section by discussing the compatibility of vanishing cycles and tensor product. Given two $T$-equivariant functions $f_i: X_i\to \BoC$ and two bounded complexes $\Sh{F}_i$, there is a natural isomorphism
\begin{equation}
    \label{eq: TS}
    \phip{f_1+f_2}\Sh{F}_1\boxtimes_{\B T} \Sh{F}_2 \cong \phip{f}\Sh{F}_1\boxtimes_{\B T} \phip{f}\Sh{F}_2
\end{equation}
known as the Thom-Sebastiani isomorphism \cite{Ma01}. It follows that there is a canonical morphism 
\[
\HO_T(X_1,\phip{f_1}\DD\BoQ)\otimes\HO_T(X_2,\phip{f_2}\DD\BoQ)\to \HO_T(X_1\times X_2, \phip{f_1+f_2}\DD\BoQ).
\]

\subsection{Dimensional reduction}

Let $X$ be smooth Artin stacks and let $E\to X$ be a vector bundle with a section $s: X\to E$.  Let $\bar s: E^\vee\to \BoC$ be the function associated to $X$. Let $Z(s)\subset X$ be the vanishing locus of $s$ and consider the fiber square
\[
\begin{tikzcd}
    \bar Z(s)\arrow[d, "\pi"]\arrow[r] &  E\arrow[d] 
    \\
    Z(s)\arrow[r, "i"] & X
\end{tikzcd}
\]
Consider the canonical morphism of functors
\begin{equation}
\label{eq: sheafy di red}
    \pi_* (i_* i^! \to \phip{\bar s}\iota^! )\pi^!
\end{equation}
obtained by applying the vanishing cycle $\phip{\bar s}$ to the adjunction and noticing $\phip{\bar s} i_*=i_* \phip{0}=i_*$.

\begin{theorem}[{\cite{da13, Kin22}}]
\label{thm: dimred}
    The morphism $\pi_* i_*\DD\BoQ_{\bar Z(s)}\to \pi_* \phip{\bar s}\DD\BoQ_{E}$ obtained by applying the morphism \eqref{eq: sheafy di red} to the dualizing complex $\DD\BoQ_{X}$ is an isomorphism. In particular, the induced map
    \[
    dr: \HO(Z(s), \DD\BoQ_{Z(s)}[2\rk(E)])\to \HO(E, \phip{\bar s}\DD\BoQ_{E})
    \]
\end{theorem}
We will refer to the morphism $dr$ as the dimensional reduction map. 

The dimensional reduction theorem admits variants informally referred to as partial or deformed dimensional reduction\cite{DaPa20}\cite{COZZ1}. Rather than stating the most general version of the theorem, we present what is sufficient for our purposes. Let $X$ be an algebraic variety with functions $f_1,\dots, f_n,g:X\to \BoC$. Consider the function $f=g+\sum_{i}f_it_i: X\times \AA^n\to \BoC$, where $t_i$ are coordinate functions on $\AA^n$. Let $\pi: X\times\AA^n \to X$ be the projection. Let $Z\hookrightarrow X$ denote the inclusion of the zero locus of $f_1,\dots f_n$. Arguing as above, we obtain a canonical morphism 
\begin{equation}
\label{deformed dim red iso}
    \pi_* (i_* \phip{g} i^! \to \phip{f}\iota^! )\DD\BoQ_{X\times \AA^n}
\end{equation}
\begin{proposition}[{\cite[Thm. C.1.]{COZZ1}}]
\label{prop: def dim red smooth case}
Assume that $Z\subset X$ is smooth of dimension $n$. Then the morphism \eqref{deformed dim red iso} is an isomorphism. 
\end{proposition}
In the main result of \cite{DaPa20}, the smoothness assumption is replaced with a quasi-homogeneity condition on $f$.

\subsection{Refined Gysin pullback}
\label{app: Refined Gysin pullback}
We now review the general construction of refined (or virtual) Gysin pullback. We focus on the level of generality required in this paper. For a more general treatment, see \cite{KhanBM}. Let $B$ be a (possibly non-smooth) algebraic stack and let $\pi: E\to B$ be a vector bundle on it. Let $s: B \to E$ be a section and let $X\coloneqq s^{-1}(0)=B\times_E B\subset B$. Then $s$ fits in the Cartesian square
\begin{equation}
\label{eq: diagram general gysin pullback}
    \begin{tikzcd}
        B\arrow[r, "o"] & E\\
        X\arrow[u, hookrightarrow, "i"] \arrow[r, hookrightarrow, "i"] & B\arrow[u, "s"]
    \end{tikzcd}
\end{equation}
where $o$ is the zero section. 
Since $E$ is a vector bundle over $B$, we have an orientation isomorphism $o^*\DD\BoQ_{E}\cong\DD\BoQ_B[2\rk(E)]$. Therefore, applying the sheaf $\DD\BoQ_{E}$ to the canonical morphism
\begin{equation}
    \label{eq: nat trans inducing Gysin pullback}
    s^!(\id\to o_*o^*)\cong (s^!\to i_!i^!o^*)
\end{equation}
we get a canonical morphism of complexes $\DD\BoQ_{B}\to i_*\DD\BoQ_X[2\rk(E)]$, and thus a canonical morphism of functors
\begin{equation}
    \label{eq: sheafy refined Gysin}
    s^!: \DD\BoQ_{B}\to i_*\DD\BoQ_{X}[2\rk(E)].
\end{equation}
The refined Gysin pulback (or virtual pullback) is the induced map of derived global sections 
\begin{equation}
    \label{eq: gysin pullback}
    s^!: \HO(B, \DD\BoQ_{B})\to \HO(X, \DD\BoQ_{X}[2\rk(E)]).
\end{equation}
and hence \eqref{eq: gysin pullback}. We call the pair $(E,s)$ virtual pullback data for $i: X\to B$.
\begin{remark}
\label{rem: virtual pull: swap o and s}
    We remark that the construction is symmetric in $s$ and $o$, i.e. if their roles are swapped in the construction, the resulting refined Gysin morphism is the same. This easily follows from the fact that the homeomorphism $E\to E$ given by $(b,v)\mapsto (b,s(b)-v)$ swaps $s$ and $o$.
\end{remark}
\begin{remark}
\label{rem: gysin + ordinary at once}
    In the article, we will systematically consider a morphism $h=g\circ i$ obtained as the composition 
    \[
    X\xrightarrow{i} B\xrightarrow{g} Y
    \]
    where $g$ is smooth of relative dimension $r$ and $i$ admits virtual pullback data $(E,s)$. Then, we can define the virtual pullback $h^!$ from $Y$ to $X$ as the composition
    \[
    h^!: \HO(Y, \DD\BoQ_Y)\xrightarrow{g^*} \HO(B, \DD\BoQ_B[-2r])\xrightarrow{s^!} \HO(X, \DD\BoQ_B[-2\vdim(h)])
    \]
    where $\vdim(h)=2r-2\rk(E)$.
\end{remark}
\begin{remark}
\label{rem: gysin alternative version}
    Often, we will be in the situation where there exists a vector bundle $F$ on $Y$ such that $E$ is the pullback of $F$ along $g: B\to Y$. In this case, we have a commutative diagram
    \[
    \begin{tikzcd}
        B\arrow[r, "s"] & E\arrow[r, "t"] & F\\
        X\arrow[u, hookrightarrow, "i"] \arrow[r, hookrightarrow,"i"] & B\arrow[u, "o"]\arrow[r, "g"] & Y\arrow[u, "o"]
    \end{tikzcd}
    \]
    with both squares Cartesian. Let $u=s\circ t$ and $\pi: F\to Y$. It is easy to check that the canonical morphism $o^!(\pi^*\DD\BoQ_{Y}\to u_*u^* \pi^*\DD\BoQ_{Y})$ can be identified via base change with a morphism $\Sh{F}\to i_*h^!\Sh{F}[-2\vdim(h)]$ that, upon taking derived global sections, this morphism recovers \eqref{eq: sheafy refined Gysin}.
\end{remark}
From the commutative diagram 
\begin{equation}
\label{eq: sheafy gysin definition}
    \begin{tikzcd}
        s_*s^!\arrow[d, swap, "\counit"]\arrow[r] & s_* i_*i^! o^*=o_*i_*i^! o^*\arrow[d, "\counit"]
        \\ 
        \id \arrow[r, "\unit"] & o_*o^*
    \end{tikzcd}
\end{equation}
together with the fact that $o^*: \HO(\DD\BoQ_{E}\to o_*o^*\DD\BoQ_{E})$ is the inverse of the pullback morphism $\pi^*: \HO(E, \DD\BoQ_E)\to \HO(B, \DD\BoQ_B[2\rk(E)])$, we deduce that the Gysin pullback fits the following commutative diagram
\[
\begin{tikzcd}
    \HO(B, \DD\BoQ_{B})\arrow[rr, "s^!"]\arrow[d, "s_*"] & & \HO(X, \DD\BoQ_{X}[2\rk(E)])\arrow[d, "i_*"]
    \\
    \HO(E, \DD\BoQ_{E}) \arrow[rr, "(\pi^*)^{-1}=o^*"]& &\arrow[ll, "\pi^*"]  \HO(B, \DD\BoQ_{B}[2\rk(E)])
\end{tikzcd}
\]
where $s_*$ and $i_*$, are ordinary pushforward morphisms.

Virtual pullbacks are compatible with smooth pullback and proper pushforward under the appropriate assumptions. For us, the following statement will suffice. Let $i: X\hookrightarrow B$ and $i': X'\hookrightarrow B'$ be two morphisms with virtual pullback data $(E,s)$ and $(E',s')$. Let also $f: B'\to B$ and $g: X'\to X$ be morphisms such that $f\circ i=i'\circ  g$. We say that the virtual pullback data $(E,s)$ and $(E',s')$ are compatible along $f$ if $(E',s')$ is isomorphic to the pullback of $(E',s')$ along $f$. Explicitly, this means that we have a commutative diagram 
\begin{proposition}
\label{prop: compat gysin and pull push Cartesian squares}
    Assume that $(E,s)$ and $(E',s')$ are compatible along $f$.
    \begin{enumerate}
        \item If $f$ is smooth (and hence $g$ is smooth), then $(s')^!\circ f^*= g^* \circ s^!$.
        \item If $f$ is proper (and hence $g$ is proper), then $s^! \circ f_*= g_* \circ (s')^!$.
    \end{enumerate}
\end{proposition}

    The refined pullback can also be interpreted in therms of bivariant homology theories \cite{KhanBM}. Recall that given any morphism $f: X\to Y$, one defines
    \[
    \HO^*(X/Y)=\HO^*(X, f^!\BoQ_Y)=\Hom_{D(Y)}(\BoQ_Y, f_*f^!\BoQ_Y[*])=\Hom_{D(X)}(\BoQ_X,f^!\BoQ_Y[*])
    \]
    so that we recover usual Borel-Moore homology as $H(X, \DD\BoQ_X)=H(X/\pt)$.
    There is a multiplicative structure 
    \begin{equation}
        \label{eq: mult structure bivariant theory}
        \HO^{d}(X/Y)\otimes  \HO^{e}(Y/Z)\to  \HO^{d+e}(X/Z) (\alpha, \beta)\mapsto \alpha\circ \beta
    \end{equation}
    given by 
    \[
    (\BoQ_X\xrightarrow{\alpha} f^!\BoQ_Y)\otimes (\BoQ_Y\xrightarrow{\beta} g^!\BoQ_Z)\mapsto (\BoQ_X\xrightarrow{\alpha \circ \beta} f^!g^!\BoQ_Z).
    \]
    Given a smooth stack $p: X\to \pt$ of dimension $d$, its fundamental class $[X]=[X/\pt]\in H^{-2d}(X)$ is the class of the canonical purity morphism $\BoQ_{X}\xrightarrow{\cong} p^!\BoQ_X[-2d]$. More generally, given a smooth morphism $f: X\to Y$ of relative dimension $d$, the canonical isomorphism $\BoQ_X\xrightarrow{\cong} f^!\BoQ[-2d]$ induced by purity $f^!\cong f^*[-2d]$ gives a canonical relative fundamental class $[X/Y]\in H^{-2d}(X/Y)$. Moreover, the pullback morphism 
    \[
    f^*: \HO^{e}(Y, \DD\BoQ_Y) \to \HO^{e-2d}(Y, \DD\BoQ_Y)
    \]
    is given by the action of the fundamental class 
    \[
    H^{e}(Y/\pt)\xrightarrow{}\HO^{-2d}(X/Y)\otimes H^{e}(Y/\pt)\to H^{e-2d}(X/\pt) \qquad \alpha\to [X/Y]\circ \alpha.
    \]
    This picture can be generalized to refined Gysin pullbacks. Consider the set up of diagram \eqref{eq: diagram general gysin pullback}. Then the (dual of the) orientation isomorphism $\BoQ_{E}\xrightarrow{\cong} o^!\BoQ_B[-2\rk(E)]$ gives a virtual class $[X/B]\in H(X/B)$ defined by
    \[
    \BoQ_X\xrightarrow{\cong }i^*\BoQ_B\to i^*o^!\BoQ_{E}[2\rk(E)]\xrightarrow{\text{b.ch.}} i^!s^* \BoQ_{E}[2\rk(E)]=i^! \BoQ_{B}[2\rk(E)]
    \] 
    where $\text{b.ch.}$ stands for the canonical base change morphism $i^*o^!\to i^!s^*$. Then the Gysin pullback can be equivalently defined as the map 
    \[
    [X/B]\circ : \HO^{e}(B\to \pt)\to  \HO^{e-2\rk{E}}(X\to \pt) \qquad \alpha\mapsto [X/B]\circ \alpha
    \]
    induced by the fundamental class and the multiplicative structure \eqref{eq: mult structure bivariant theory}.
\begin{remark}
    Virtual fundamental classes can be defined more generally for quasi-smooth morphisms $f: X\to B$ of derived Artin stacks. The general construction involves a deformation to the normal cone construction \cite[\S3.1]{KhanBM}.  The local picture from diagram \eqref{eq: diagram general gysin pullback} can be seen as the local model for a quasi smooth immersion $i: X\to B$ and suffices for our purposes. 
    
    The following result is equivalent relates associativity of the virtual pullback with compatibility of the derived structures. \begin{proposition}[{\cite[Prop 3.12]{KhanBM}}]
\label{prop: quasismooth gysin functorial}
    Assume we are given two quasi-smooth morphisms $f: X\to Y$ and $g: Y\to Z$ of derived algebraic stacks, then we have $[X/Y]\circ [Y/Z]=[X/Z]$. Therefore, we have $f^! \circ g^!=(f\circ g)^!$.
\end{proposition}
\end{remark}

\newpage

\section{Proof of Theorem \ref{thm: main thm positive Coulomb hall}}
\label{subsec: Proof of long Theorem}

\label{sec: Proof of Theoorem Hall Coulomb app}
\subsection{Statement and conventions}
\label{subsec: APP Statement and conventions}
This appendix is devoted to the proof of Theorem \ref{thm: main thm positive Coulomb hall}. It is convenient to we restate it here by expressing all the terms in terms of moduli of sheaves over the formal disk $\DD$. We also drop the torus $T$ from the notation. it will be clear that all arguments hold $T$-equivariantly.
\begin{theorem}
\label{thm: main thm positive Coulomb hall app}
The morphism \eqref{eq: map hall to Coulomb} is a morphism of $\BoN^{Q_0}$-graded cohomologically graded algebras. In other words, for any $\dd\in \BoN^{Q_0}$ and any decomposition $\dd=\dd'+\dd''$, the diagram 
			\[
			\begin{tikzcd}
			 \HO^{\BM}(\FCoh_{Q,0,\dd'}(\DD), \BoQ)\otimes \HO^{\BM}(\FCoh_{Q,0,\dd''}(\DD), \BoQ)\arrow[d, "\theta_{\dd'}\otimes \theta_{\dd''}"]\arrow[r, "\vmult^{\onil}"] & \HO^{\BM}(\FCoh^0_{Q,\dd}(\DD), \BoQ)
        \arrow[d, "\theta_{\dd'+\dd''}"]
        \\
                \HO^{\BM}_{G\fps}( \Rt^+_{\dd'}, \BoQ[s_{\dd'}])\otimes \HO^{\BM}_{G\fps}( \Rt^+_{\dd''}, \BoQ[s_{\dd''}])\arrow[r, "\eta^+"] & \HO^{\BM}_{G\fps}( \Rt^+_{\dd}, \BoQ[s_{\dd}])
			\end{tikzcd}
			\]
		is commutative. 
\end{theorem}
Recall that the map $\theta_{\dd}$ is defined in \S\ref{subsec: From hall to Coulomb} using the following chain of maps,
\begin{equation}
\label{eq: factorization map hall Coulomb app}
		\begin{tikzcd}
			\Rt^+_{\dd}\arrow[r, hookrightarrow, "a"]& B_{\dd}\arrow[r, "b_1"] &C_{\dd} \arrow[r, "b_2"] & \FCoh_{Q, 0,\dd}(\DD)
		\end{tikzcd}
\end{equation}
where $B_{\dd}\coloneqq\Tt_{\dd}\times_{\FCoh_{0,\dd}(\DD)} \FCoh_{Q,0,\dd}(\DD)$ and $C_{\dd}\coloneqq\Gr^+_{\dd}\times_{\FCoh_{0,\dd}(\DD)} \FCoh_{Q, 0,\dd}(\DD)$. All the maps are the obvious ones. It will be convenient to consider an alternative factorization of the map $\Rt^+_{\dd}\to \FCoh_{Q, 0,\dd}(\DD)$. Specifically, we consider 
\begin{equation}
\label{eq: alternative factorization map hall Coulomb}
    \begin{tikzcd}
    \Rt^+_{\dd}\arrow[r, hookrightarrow, "a"]  &  B_{\dd} \arrow[r, "b"] & N\fps \times \FCoh_{Q, 0,\dd}(\DD) \arrow[r] & \FCoh_{Q, 0,\dd}(\DD)
\end{tikzcd}
\end{equation}
where the projection of $b$ to $N\fps$ sends a tuple $(Q, g, n, \beta)\in  B_{\dd}$ to $n$ and the projection to $\FCoh_{Q, 0,\dd}(\DD)$ sends $(Q, g, n, \beta)$ to $(Q, \beta)\in \FCoh_{Q, 0,\dd}(\DD)$. Taking quotients by $G\fps$, we obtain the sequence 
\begin{equation}
\label{eq: alternative factorization map hall Coulomb after quotienting}
    \begin{tikzcd}
    \Hecke_{\dd}(\DD, N/G)\arrow[r, hookrightarrow, "a"]  &  D_{\dd} \arrow[r, "b"] & \Map(\DD, N/G) \times \FCoh_{Q, 0,\dd}(\DD) \arrow[r] & \FCoh_{Q, 0,\dd}(\DD)
\end{tikzcd}
\end{equation}
where we have set 
\begin{equation}
    \label{eq: var triples and filt}
    \Hecke_{\dd}(\DD, N/G)\coloneqq G\fps\backslash \Rt^+_{\dd}\qquad D_{\dd}\coloneqq G\fps\backslash B_{\dd}.
\end{equation}
We also define $c: \Hecke_{\dd}(\DD, N/G)\to \Map(\DD, N/G)\times \FCoh_{Q, 0,\dd}(\DD)$ to be the composition of $c=b\circ a$.
Notice that $\Hecke_{\dd}(\DD, N/G)$ is the moduli space parametrizing triples $(\Sh{V}, \Sh{V'}, g, n, n')$ where $\Sh{V}=\{\Sh{V}_i \}_{i\in Q_0}$ and $\Sh{V}=\{\Sh{V}'_i \}_{i\in Q_0}$ are tuples of rank $\vv_i$ vector bundles, $g=\{g_i\}_{i\in Q_0}$ is a collection of injective morphisms $g: \Sh{V}_i\hookrightarrow\Sh{V}'_i$ of locally free sheaves on the formal disk $\DD$
with torsion cokernels $\Sh{V}'_i/\Sh{V}_i$ such that $\text{length}(\Sh{V}'_i/\Sh{V}_i)=\dd_i$ and $ n=\{n_a\}_{a\in Q_1}$ and $ n'=\{n'_a\}_{a\in Q_1}$ are two collections of morphisms of sheaves
\[
n_a: \Sh{V}_{s(a)}\to \Sh{V}_{t(a)}\qquad n'_a: \Sh{V'}_{s(a)}\to \Sh{V}'_{t(a)}
\]
fitting in the commutative diagram 
\[
\begin{tikzcd}
    0\arrow[r]&  \Sh{V}_{s(a)} \arrow[d, "n_a"]\arrow[r, "g_{s(a)}"] & \Sh{V}'_{s(a)}\arrow[d, "n_a'"]
    \\
     0\arrow[r] & \Sh{V}_{t(a)}\arrow[r, "g_{t(a)}"] & \Sh{V}_{t'(a)}
\end{tikzcd}
\]
The composite map $\Hecke_{\dd}(\DD, N/G)\to \FCoh_{Q, 0,\dd}(\DD)$ sends $(\Sh{V}, \Sh{V'}, g, n', n)$ to $(\Sh{V}'/\Sh{V}, \beta)$, where $\Sh{V}'/\Sh{V}$ stands for the tuple of rank zero quotients $\{ \Sh{V}_i/\Sh{V}_i\}_{i\in Q_0}$ and $\beta=\{\beta_a\}_{a\in Q_1}$ is a collection of morphisms of sheaves $\beta_a: \Sh{V}'_{s(a)}/\Sh{V}_{s(a)}\to \Sh{V}'_{t(a)}/\Sh{V}_{t(a)}$ fitting in a commutative square 
\[
\begin{tikzcd}
    0\arrow[r]  & \Sh{V}_{s(a)}\arrow[d, "n_a"]\arrow[r, "g_{s(a)}"] & \Sh{V}'_{s(a)}\arrow[d, "n'_a"] \arrow[r] &\Sh{V}'_{s(a)}/\Sh{V}_{s(a)}\arrow[d, dashed, "\beta_a"]\arrow[r] & 0
    \\
    0\arrow[r]  &\Sh{V}_{t(a)}\arrow[r, "g_{t(a)}"] & \Sh{V}'_{t(a)}\arrow[r] &\Sh{V}'_{t(a)}/\Sh{V}_{t(a)}\arrow[r] & 0
\end{tikzcd}
\]
Similarly, the moduli space $D_{\dd}$ parametrizes tuples $(\Sh{V}', Q, q, n', \beta)$ where $\Sh{V}'$ and $n'$ is as above, $Q=\{ Q_i\}_{i\in Q_0}$ is a tuple of rank zero length $\dd=(\dd_i)_{i\in Q_0}$ coherent sheaves, $\beta'$ is a tuple $Q_{s(a)}\to Q_{t(a)}$ and $q=\{q_i\}_{i\in Q_0}$ is a tuple of quotient maps $q_i: \Sh{V}'_{i}\twoheadrightarrow Q_i$. The commutativity of the diagram 
\[
\begin{tikzcd}
\Sh{V}'_{s(a)}\arrow[d, "n'_a"] \arrow[r] &Q_{s(a)}\arrow[d, "\beta_a"]\arrow[r] & 0
\\
\Sh{V}'_{t(a)}\arrow[r] &Q_{t(a)}\arrow[r] & 0
\end{tikzcd}
\]
is \emph{not} required at this stage. The substack $\Hecke_{\dd}(\DD, N/G)$ is naturally identified with the substack consisting of those tuples such that the diagram above $\emph{does}$ commute.

Recall now that the morphism 
\[
\theta_{\dd}:  \HO^{\BM}(\FCoh^0_{Q,\dd}(\DD), \BoQ)\to \HO^{\BM}_{G\fps}( \Rt^+_{\dd}, \BoQ^{\vir})
\]
is defined as the composition of pullback morphisms along the smooth morphisms $b_1$ and $b_2$ from \eqref{eq: factorization map hall Coulomb app} and then further via virtual pullback along the quasi-smooth morphism $a$. By functoriality of the pullback morphism, we can equivalently define $\theta_{\dd}$ as the following composition 
\begin{equation}
    \label{eq: diagram virtaul pull hall Coulomb app}
    \begin{tikzcd}
    \HO^{\BM}(\FCoh_{Q,0,\dd}(\DD), \BoQ)\arrow[r, equal]\arrow[d, "\theta_{\dd}"] & \HO^{\BM}(\Map(\DD, N/G)\times \FCoh^0_{Q,\dd}(\DD), \BoQ)\arrow[d, "b^*"]
    \\
     \HO^{\BM}(\Hecke_{\dd}(\DD, N/G), \BoQ^{\vir}) & \arrow[l, swap, "a^!"] \HO^{\BM}(D_{\dd}, \BoQ[-\vv^T\dd])
\end{tikzcd}
\end{equation}
where $b^*$ is smooth pullback along $b: D_{\dd}\to \Map(\DD, N/G)\times \FCoh^0_{Q,\dd}(\DD)$ and $a^!$ is the virtual pullback induced by the virtual data from Proposition \ref{prop: virtual pullback data hall and Coulomb}. We also set $c^!=a^!\circ b^*$. Here, by slight abuse of notation, we are making the identification 
\[
 \HO^{\BM}(\Hecke_{\dd}(\DD, N/G), \BoQ^{\vir})=\HO^{\BM}_{G\fps}( \Rt^+_{\dd}, \BoQ^{\vir})
\]
although the left hand side is really defined as the right hand side. We believe that the benefits of this slight abuse of notation outweigh the drawbacks. To make precise sense of the Borel-Moore cohomology of all the spaces it suffices to pass to a finite approximation $\Spec(\BoC[t]/t^k)$ of the formal disk $\DD$ and hence a finite type approximation of the mapping stack $\Map(\DD, N/G)$ in \eqref{eq: alternative factorization map hall Coulomb after quotienting}. This is, indeed, one of the advantages of the construction of $\theta_{\dd}$ in terms of the factorization \eqref{eq: alternative factorization map hall Coulomb}—or equivalently \eqref{eq: alternative factorization map hall Coulomb after quotienting}. In fact it is easy to check that we have a compatible system of morphisms
\[
\begin{tikzcd}
    \Hecke_{\dd}(\DD, N/G)\arrow[d] \arrow[r, hookrightarrow, "a"]  &  D_{\dd} \arrow[r, "b"] \arrow[d]& \Map(\DD, N/G) \arrow[d]\times \FCoh_{Q, 0,\dd}(\DD) \arrow[r] & \FCoh_{Q, 0,\dd}(\DD)
    \\
    \Hecke^{k}_{\dd}(\DD, N/G)\arrow[r, hookrightarrow, "a^k"]  &  D^k_{\dd} \arrow[r, "b^k"] & \Map(\DD_{k}, N/G) \times \FCoh_{Q, 0,\dd}(\DD) \arrow[ur] 
\end{tikzcd}
\]
where all the squares are Cartesian.
All the spaces in the bottom row are algebraic stacks of finite type. Although in the course of the proof we will work formally over $\DD$, all our constructions are relative to $\Map(\DD, N/G)$, are well defined on finite approximations $\Map(\DD_{k}, N/G)$ for $k\ggg 0$, and are always to be understood in terms of these finite approximations.

We conclude this preparatory section by pointing out that, with slight abuse of notation, we will systematically drop the information about the sheaves when describing closed points in all the moduli spaces in \eqref{eq: alternative factorization map hall Coulomb after quotienting} as well as in their variants defined in the course of the next section. For instance, we will denote a point in $\Hecke_{\dd}(\DD, N/G)$ as $(g, n', n)$ as opposed to $(\Sh{V}, \Sh{V'}, g, n', n)$, and a point in $D_{\dd}$ as $(h, n', \beta)$ as opposed to $(\Sh{V}', Q, h, n', \beta)$. In other words, it will be implicitly understood that all the sheaves will be allowed to vary in moduli.

\subsection{Strategy}
Both multiplication maps $\vmult^{\Xi_Q^{\onil}}$ and $\eta^+$ in the statement of Theorem \ref{thm: main thm positive Coulomb hall app} are constructed via a pull-push in Borel-Moore homology. In the proof, we prove that both the pull operation and the push operation are compatible with the map $\theta$. Fix $\dd\in \BoN^{Q_0}$ and a decomposition $\dd=\dd'+\dd''$. Let $\Hecke^{(2)}(\DD, N/G)$ be the moduli stack parametrizing tuples $(g,h, n,n', n'')$ where $g=\{g_i\}_{i\in Q_0}$ and $h=\{h_i\}_{i\in Q_0}$ are collections of injective morphism of locally free sheaves on $\DD$
\[
\Sh{V}_i\xhookrightarrow{g_i}\Sh{V}_i' \xhookrightarrow{h_i}\Sh{V}''_i
\]
with rank zero torsion co-kernels and $ n=\{n_a\}_{a\in Q_1}$, $ n'=\{n'_a\}_{a\in Q_1}$, and $ n''=\{n''_a\}_{a\in Q_1}$ are three collections of morphisms of locally free sheaves 
\[
n_a: \Sh{V}_{s(a)}\to \Sh{V}_{t(a)}
\qquad 
n_a': \Sh{V}'_{s(a)}\to \Sh{V}'_{t(a)}
\qquad 
n''_a: \Sh{V}''_{s(a)}\to \Sh{V}''_{t(a)}
\]
such that the following diagrams commute: 
\begin{equation}
\label{eq: data parametrized by Filt^2}
    \begin{tikzcd}
    \Sh{V}_{s(a)} \arrow[r, "g_{s(a)}"]\arrow[d, "n_a"] & \Sh{V}'_{s(a)} \arrow[r, "h_{s(a)}"]\arrow[d, "n'_a"]& \Sh{V}''_{s(a)}\arrow[d, "n''_a"] 
    \\
    \Sh{V}_{t(a)} \arrow[r, "g_{t(a)}"] & \Sh{V}'_{t(a)}\arrow[r, "h_{t(a)}"] & \Sh{V}''_{t(a)}
\end{tikzcd}
\end{equation}
We denote by $\Hecke^{(2)}_{\dd',\dd''}\subset \Hecke^{(2)}$ the component such that $\Sh{V}'/\Sh{V}$ and $\Sh{V}''/\Sh{V}'$ have, respectively, length vector $\dd'$ and $\dd''$. The stack $\Hecke^{(2)}_{\dd',\dd''}$ admits an explicit presentation as a global quotient, and it is easy to check that 
\begin{equation}
    \label{eq: var triples and filt 2}
    \Hecke^{(2)}_{\dd',\dd''}(\DD, N/G)=G\fps\backslash \left(q^+((p^+)^{-1}(\Rt^+_{\dd'}\times \Rt^+_{\dd''}))\right)
\end{equation}
where $q^+((p^+)^{-1}(\Rt^+_{\dd'}\times \Rt^+_{\dd''})$ was the auxiliary space introduced in \S\ref{subsec: Positive part of the Coulomb branch} for the definition of the convolution product on $\HO^{\BM}_{G\fps}( \Rt^+, \DD\BoQ_{\Rt^+}, \BoQ)$. Let also $\FCoh_{Q,0, \dd',\dd''}(\DD)$ be the stack of short exact sequences of type $\dd',\dd''$ in $\Coh_Q(\DD)$, namely the stack parametrizing tuples tuples of short exact sequences $0\to Q'_i\to Q_i\to Q''_i\to 0$ of torsion sheaves such that $\text{length}(Q'_i)=\dd_i'$ and $\text{length}(Q''_i)=\dd_i''$ together with tuples of morphisms $\beta=\{\beta_a\}_{a\in Q_1}$, $\beta'=\{\beta'_a\}_{a\in Q_1}$ and $\beta''=\{\beta''_a\}_{a\in Q_1}$ fitting in commutative diagrams
\begin{equation*}
    \begin{tikzcd}
        0\arrow[r] & Q'_{s(a)}\arrow[r]\arrow[d, "\beta'_{a}"] & Q_{s(a)}\arrow[r]\arrow[d, "\beta''_{a}"] & Q''_{s(a)}\arrow[r]\arrow[d, "\beta''_{a}"] & 0
        \\
        0\arrow[r] & Q'_{t(a)}\arrow[r] & Q_{t(a)}\arrow[r] & Q''_{t(a)}\arrow[r] & 0
    \end{tikzcd}
\end{equation*}
Notice that we have a canonical morphism 
\[
d: \Hecke^{(2)}_{\dd',\dd''}(\DD, N/G)\to \FCoh_{Q,0,\dd',\dd''}(\DD)
\]
sending a chain of inclusions $\Sh{V}_i\xhookrightarrow{g_i}\Sh{V}_i' \xhookrightarrow{h_i}\Sh{V}''_i$ to the short exact sequence of rank zero sheaves
    \[
    0\to \Sh{V}'_i/\Sh{V}_i \to \Sh{V}''_i/\Sh{V}_i\to \Sh{V}''_i/\Sh{V}_i'\to 0
    \]
    and the sections $(n,n',n'')$ to induced maps on the quotient factors
    \begin{equation}
    \label{eq: induced maps quotients Filt2}
        \begin{tikzcd}
            0\arrow[r] &\Sh{V}_{s(a)}'/\Sh{V}_{s(a)}\arrow[r]\arrow[d, "\beta'_a"] & \Sh{V}_{s(a)}''/\Sh{V}_{s(a)}\arrow[r]\arrow[d, "\beta_a"] & \Sh{V}_{s(a)}''/\Sh{V}_{s(a)}'\arrow[d, "\beta''_a"] \arrow[r]& 0
            \\
            0\arrow[r] &\Sh{V}_{s(a)}'/\Sh{V}_{t(a)}\arrow[r] & \Sh{V}_{t(a)}''/\Sh{V}_{t(a)}\arrow[r]& \Sh{V}_{t(a)}''/\Sh{V}_{t(a)}'\arrow[r]& 0
        \end{tikzcd}
    \end{equation}
Now consider the diagrams
\begin{equation}
    \label{eq: diagram Coulomb hall is mor algebra 1}
    \begin{tikzcd}    \Hecke_{\dd'}(\DD, N/G)\times\Hecke_{\dd''}(\DD, N/G)\arrow[d, swap, "\tau^{(23)}\circ(c\times c)"]  & \Hecke^{(2)}_{\dd',\dd''}(\DD, N/G)\arrow[l, swap,  "q_1\times q_2"]\arrow[d, "l'\times l''\times d"]
   \\
   \Map(\DD, N/G)^2\times \FCoh_{Q,0,\dd'}(\DD)\times \FCoh_{Q,0, \dd''}(\DD)\arrow[d, "\pr_2\times \id\times \id"] &  \Map(\DD, N/G)^{2}\times\FCoh_{Q,0,\dd',\dd''}(\DD)\arrow[l, swap, "\id^2\times \bar q"]\arrow[d, "\pr_2\times \id "] 
   \\
   \Map(\DD, N/G)\times \FCoh_{Q,0,\dd'}(\DD)\times \FCoh_{Q,0,\dd''}(\DD) & \Map(\DD, N/G)\times\FCoh_{Q,0,\dd',\dd''}(\DD)\arrow[l, swap, "\id\times \bar q"]
\end{tikzcd}
\end{equation}
\begin{equation}
    \label{eq: diagram Coulomb hall is mor algebra 2}
\begin{tikzcd}
    \Hecke^{(2)}_{\dd',\dd''}(\DD, N/G)\arrow[d, swap, "\gr^{(2)}"] \arrow[r, "\tilde m"] & \Hecke^{(2)}_{\dd', \dd''}(\DD, N/G)\arrow[d, "l''\times c"]
    \\
    \Map(\DD, N/G)\times\FCoh_{Q,0,\dd',\dd''}(\DD)\arrow[r, "\bar p\times \id"] & \Map(\DD, N/G)\times \FCoh_{Q,0, \dd}(\DD)
\end{tikzcd}
\end{equation}
where 
\begin{itemize}
    \item $\tau^{(23)}$ swaps the second and third factor of the fiber product.
    \item the map $q_1$ sends the tuple $(g, h, n, n', n'')$ to $(g, n, n')$ and, similarly,  $q_2$ sends a tuple $(g, h, n, n', n'')$ to $(h, n', n'')$.
    \item The map $l'$ (resp. $l''$) sends a tuple $(g, h, n'', n', n)$ to $n'$ (resp. to $n''$).
    \item The maps $\bar q$ and $\bar p$ send, respectively a short exact sequence of $Q$-valued coherent sheaves to its associated graded and to its middle term. They coincide with the maps $\bar q$ and $\bar p$ in \eqref{eq: diagram coha Xi_Q} under the identification $\FM^{T, \onil}_{\dd}(\Xi_Q)\cong \FCoh_{Q,0,\dd}(\DD)$.
\end{itemize}
The diagrams are commutative. Notice that the top rows of the diagrams above are the maps defining the convolution product on the positive Coulomb branch $ \HO^{\BM}(\Hecke_{\dd}(\DD, N/G), \BoQ^{\vir})$. On the other hand, the bottom rows are, up to the cohomologically trivial factor $\Map(\DD, N/G)$, the maps entering in the definition of the Hall product on $\HO^{\BM}(\FCoh_{Q,0,\dd}(\DD), \BoQ)$. In essence, we generalize the construction of \S\ref{subsec: From hall to Coulomb}. 

Set $c^{(2)}=\pr_1\circ (l'\times l'' \times d)=l''\times d$. As a first step of the proof, we construct a virtual pullback 
\begin{equation}
    \label{eq: virtual pull filt version}
    (c^{(2)})^!: \HO^{\BM}(\Map(\DD, N/G) \times \Filt_{\dd',\dd''}(\FCoh_{Q, 0,\dd}(\DD)), \BoQ)\to \HO^{\BM}(\Hecke^{(2)}(\DD, N/G), \BoQ)
\end{equation}
generalizing the map $c^!=a^!\circ b^*$ from diagram \eqref{eq: diagram virtaul pull hall Coulomb app}.
Let $D^{(2)}$ be the moduli space parametrizing the data $(g,h,n,\beta', \beta, \beta'')$, where $g$ and $h$ are inclusions of vector bundles $\Sh{V}_i\xhookrightarrow{g_i} \Sh{V}'_i \xhookrightarrow{h_i} \Sh{V}''_i$ of the same rank, $n=\{n_a\}_{a\in Q_1}$ are morphisms of coherent sheaves $n_a:\Sh{V}_{s(a)}\to\Sh{V}_{t(a)}$ and $\beta, \beta',\beta''$ are collections of morphism fitting in the diagram 
\begin{equation}
    \label{eq: diagram maps associated greded app}
\begin{tikzcd}
    \Sh{V}'_{s(a)}/\Sh{V}_{s(a)}\arrow[r] \arrow[d, "\beta'_a"]& \Sh{V}''_{s(a)}/\Sh{V}_{s(a)}  \arrow[d, "\beta_a"]\arrow[r] & \Sh{V}''_{s(a)}/\Sh{V}'_{s(a)}\arrow[d, "\beta''_a"]
    \\
    \Sh{V}'_{t(a)}/\Sh{V}_{t(a)}\arrow[r] & \Sh{V}''_{t(a)}/\Sh{V}_{t(a)} \arrow[r] & \Sh{V}''_{t(a)}/\Sh{V}'_{t(a)}
\end{tikzcd}
\end{equation}
No compatibility between $n''$ and the maps $\beta, \beta',\beta''$ is assumed here. Let $D^{(2)}_{\dd',\dd''}$ denote the component such that the length vectors of $\Sh{V}''/\Sh{V}'$ and $\Sh{V'}/\Sh{V}$ are, respectively $\dd'$ and $\dd''$. 
Consider the following factorization of the map $c^{(2)}$:
\begin{equation}
    \label{eq: factorization c^{(2)}}
    \begin{tikzcd}
    \Hecke^{(2)}_{\dd',\dd''}(\DD, N/G)\arrow[r, "a^{(2)}"]\arrow[dr, "c^{(2)}"]  &  D^{(2)}_{\dd',\dd''} \arrow[d, "b^{(2)}"] 
    \\
    & \Map(\DD, N/G) \times \FCoh_{Q, 0, \dd',\dd''}(\DD) 
\end{tikzcd}
\end{equation}
Here,
\begin{enumerate}
    \item $a^{(2)}$ sends a closed point $(g,h, n,n', n'')\in  \Hecke^{(2)}_{\dd',\dd''}(\DD, N/G)$ to $(g,h,n,\beta', \beta, \beta'')$, where the maps $(\beta', \beta, \beta'')$ are induced by $(n,n', n'')$ by taking iterated quotients.
    \item the projection on $b^{(2)}$ on the first factor sends a closed point $(g,h,n,\beta', \beta, \beta'')\in D^{(2)}_{\dd',\dd''}$ to $n''\in \Map(\DD, N/G) $ and the projection onto the second factor sends $(g,h,n,\beta', \beta, \beta'')$ to the data \eqref{eq: diagram maps associated greded app}, which define a closed point in $\FCoh_{Q, 0,\dd',\dd''}(\DD)$.
\end{enumerate}

We will show that, similarly to \eqref{eq: alternative factorization map hall Coulomb after quotienting}, the map $b^{(2)}$ is smooth and $a^{(2)}$ is a quasi-smooth embedding. The general machinery will then allow us to define \eqref{eq: virtual pull filt version} as the composition $(c^{(2)})^!=(a^{(2)})^!\circ (b^{(2)})^*$

\begin{lemma}
\label{lemma: first lemma factorization c^{(2)}}
    The morphism $b^{(2)}$ is smooth of relative dimension $\vv^T\dd$.
\end{lemma}
\begin{proof}
    It suffices to prove the case $Q=A^1$ (and hence $N=0$). Let $\Gr^+$ be the moduli space parametrizing pairs of embeddings $ \Sh{V}\xhookrightarrow{g}  \Sh{V}'_i \xhookrightarrow{h} \Sh{O}^{\vv}_{\DD}$ of rank $\vv$ locally free sheaves on the disk $\DD$. Here $\Sh{O}^{\vv}_{\DD}$ is fixed and does not vary in moduli. Let $\Gr^+_{\dd', \dd''}$ be the component such that $\Sh{O}_{\DD}^{\vv}/\Sh{V}'$ and $\Sh{V}'/\Sh{V}$ have length $\dd'$ and $\dd''$, respectively\footnote{The moduli space $\Gr^+_{\dd', \dd''}$ is isomorphic to the nested quot scheme parametrizing quotients $\Sh{O}^{\vv}_{\DD}\twoheadrightarrow Q\twoheadrightarrow Q'' $, where $Q$ has length $\dd=\dd'+\dd''$ and $Q''$ has length $\dd'$ Informally speaking, the isomorphism is obtained by trading quotients with their kernels.}. 
    We have a commutative diagram
    \begin{equation}
        \label{eq: diagram proof b^{(2)} is smooth}
        \begin{tikzcd}
        \Gr^+_{\dd', \dd''}\arrow[d] \arrow[r]&  \Gr^+_{\dd} \arrow[d]
        \\ 
        \FCoh_{0,\dd',\dd''}(\DD)\arrow[r] & \FCoh_{0,\dd}
    \end{tikzcd}
    \end{equation}
    where the top horizontal and left vertical maps send a chain of embeddings $\Sh{V}\xhookrightarrow{g}  \Sh{V}'_i \xhookrightarrow{h} \Sh{O}^{\vv}_{\DD}$ to the composition $\Sh{V}\xhookrightarrow{h\circ g} \Sh{O}^{\vv}_{\DD}$ and to the iterated quotient $0\to \Sh{V}'/\Sh{V}\to  \Sh{O}^{\vv}/\Sh{V}\to  \Sh{O}^{\vv}/\Sh{V}\to 0$, respectively. The diagram is Cartesian as for given any short exact sequence of rank zero sheaves  $0\to Q\to Q'\to Q''\to 0$ together with a quotient map $\epsilon: \Sh{O}^{\vv}_{\DD}\twoheadrightarrow Q$, there are unique locally free sheaves $\Sh{V}'$, $\Sh{V}''$ and injective morphisms of locally free sheaves $i: \Sh{V}\hookrightarrow \Sh{V}'$ and $j: \Sh{V}'\hookrightarrow \Sh{O}_{\DD}^{\vv}$ fitting in the commutative diagram 
    \[
    \begin{tikzcd} 
        &  0\arrow[d]& 0\arrow[d] & 0\arrow[d]
        \\
        &\Sh{V} \arrow[r, equal]\arrow[d, "j"] & \Sh{V}\arrow[d]\arrow[r, hookrightarrow, "j"] & \Sh{V}'\arrow[d, "i"] 
        \\
        & \Sh{V}' \arrow[r, hookrightarrow, "i"]\arrow[d] & \Sh{V}''\arrow[d]\arrow[r, equal] & \Sh{V}''\arrow[d] 
        \\
        0\arrow[r]& Q'\arrow[d] \arrow[r] & Q \arrow[r]\arrow[d] & Q'' \arrow[r]\arrow[d] & 0
        \\
        &0 & 0 & 0
    \end{tikzcd}
    \]
    By Lemma \ref{lemma: gr to coh is smooth}, the right vertical map in \eqref{eq: diagram proof b^{(2)} is smooth} is smooth so the same is true for the left vertical map. The quotient of the latter by the $G\fps$-action acting trivially on the target and reparametrizing $\Sh{O}_{\DD}^{\vv}$ gives the map $b^{(2)}$ in the special case $Q=A^1$. This completes the proof.
    \end{proof}

    We now proceed with the analysis of the map $a^{(2)}$. Consider the following commutative diagram, 
\begin{equation}
    \label{eq: diagram appendix hall Coulomb 1}
    \begin{tikzcd}
    &  D^{(2)}_{\dd',\dd''}\arrow[ddl, swap,  "p'"] \arrow[d, swap, "b^{(2)}"]& \arrow[l, swap,  "a^{(2)}"]\Hecke^{(2)}_{\dd',\dd''}(\DD, N/G)\arrow[dl, swap, "c^{(2)}"]\arrow[ddl, "\tilde m^+"]
    \\
    & \FCoh_{Q,0,\dd',\dd''}(\DD))\times \Map(\DD, N/G)\arrow[ddl, near start, "\bar p\times \id"]
    \\
    D_{\dd} \arrow[d, swap, "b"]& \arrow[l, swap,  "a"]\Hecke_{\dd}(\DD, N/G)\arrow[dl, "c"]
    \\
    \Map(\DD, N/G)\times \FCoh_{Q,0,\dd}(\DD)
\end{tikzcd}
\end{equation}
where $\bar p$ and $\tilde m$ are as in diagram \eqref{eq: diagram Coulomb hall is mor algebra 2} and $p'$ sends $(g, h, n, \beta', \beta, \beta'')\in D^{(2)}_{\dd',\dd''}$ to $(h\circ g, n, \beta', \beta, \beta'')\in D_{\dd}$. The proof of the following lemma follows from the modular definition of $D^{(2)}_{\dd',\dd''}$ and the analog description of $D_{\dd}$, cf. \S\ref{subsec: From hall to Coulomb}.

\begin{lemma}
\label{lemma: lemma comp hecke Coulomb 1}
    The leftmost square of diagram \eqref{eq: diagram appendix hall Coulomb 1}, which consists of the maps $b,b^{(2)},\bar p\times \id$, and $p'$, is Cartesian.
\end{lemma}

Recall from \S\ref{subsec: From hall to Coulomb} that we have a vector bundle $E_{\dd}\to D_{\dd}$ together with a section $s: D_{\dd}\to E_{\dd}$ such that $s^{-1}(0)=\Hecke_{\dd}(\DD, N/G)$.  Let $E^{(2)}_{\dd',\dd''}$ be the pullback of $E_{\dd}$ along $p'$ and let $s^{(2)}: D^{(2)}_{\dd',\dd''}\to E^{(2)}_{\dd',\dd''}$ be the pulled-back section. Let also $o^{(2)}: D^{(2)}_{\dd',\dd''}\to E^{(2)}_{\dd',\dd''}$ be the zero section.

\begin{lemma}
\label{lemma: second lemma factorization c^{(2)}}
    The morphism $a^{(2)}: \Hecke^{(2)}_{\dd',\dd''}(\DD, N/G)\to D^{(2)}_{\dd',\dd''}$ coincides with the inclusion $(s^{(2)})^{-1}(0)\hookrightarrow D^{(2)}_{\dd',\dd''}$. Moreover, we have a commutative diagram 
    \[
    \begin{tikzcd}
        & E^{(2)}_{\dd',\dd''} \arrow[dl]&   D^{(2)}_{\dd',\dd''}\arrow[l, swap, "o"]\arrow[dl, "p'"]
        \\
        E_{\dd} & D_{\dd}\arrow[l, swap, "o^{(2)}"]
        \\
        & D^{(2)}_{\dd}\arrow[uu, shift left=2, near start, "s^{(2)}"]\arrow[dl] & \Hecke^{(2)}_{\dd',\dd''}(\DD, N/G)\arrow[uu, swap, "a^{(2)}"]\arrow[dl]\arrow[l, swap, "a^{(2)}"]
        \\
        D_{\dd}\arrow[uu, "s"] & \Hecke_{\dd}(\DD, N/G)\arrow[l, swap, "a"]\arrow[uu, shift right=2, near start, "a"]
    \end{tikzcd}
    \]
    with all the squares Cartesian. 
\end{lemma}
\begin{proof}
By definition of $\Hecke^{(2)}(\DD, N/G)$ it follows that $s^{(2)}\circ a^{(2)}=o^{(2)}\circ a^{(2)}$. Hence, by the universal property of the fiber product, there exists a canonical morphism $\Hecke^{(2)}(\DD, N/G)\to (s^{(2)})^{-1}(0)$. We will show that this is an isomoprhism.
The zero locus $(s^{(2)})^{-1}(0)$ consists of those elements $(g,h,n, \beta', \beta, \beta'')\in B^{(2)}$ such that there exists morphisms $n''_a: \Sh{V}''_{s(a)}\to \Sh{V}_{t(a)}''$ making the following diagram commute 
\[
\begin{tikzcd}
     0\arrow[r] &\Sh{V}_{s(a)}\arrow[r, "h_{s(a)}g_{s(a)}"] \arrow[d, "n_a"]& \Sh{V}''_{s(a)} \arrow[d, "n''_a"]\arrow[r] & \Sh{V}''_{s(a)}/\Sh{V}_{s(a)}\arrow[d, "\beta_a"]\arrow[r] & 0 
    \\
     0\arrow[r] & \Sh{V}_{t(a)}\arrow[r, "h_{t(a)}g_{t(a)}"] & \Sh{V}''_{t(a)} \arrow[r] & \Sh{V}''_{t(a)}/\Sh{V}_{t(a)}\arrow[r] & 0 
\end{tikzcd}
\]
Recall that by definition the maps $\beta', \beta$ and $\beta''$ fit in the diagram of short exact sequences
\begin{equation}
\label{eq: diagram ses quotients}
\begin{tikzcd}
    0\arrow[r] & \Sh{V}'_{s(a)}/\Sh{V}_{s(a)}\arrow[r] \arrow[d, "\beta'_a"]& \Sh{V}''_{s(a)}/\Sh{V}_{s(a)}  \arrow[d, "\beta_a"]\arrow[r] & \Sh{V}''_{s(a)}/\Sh{V}'_{s(a)}\arrow[d, "\beta''_a"]\arrow[r] & 0 
    \\
     0\arrow[r] & \Sh{V}'_{t(a)}/\Sh{V}_{t(a)}\arrow[r] & \Sh{V}''_{t(a)}/\Sh{V}_{t(a)} \arrow[r] & \Sh{V}''_{t(a)}/\Sh{V}'_{t(a)} \arrow[r] & 0 
\end{tikzcd}
\end{equation}
where the rows are induced by taking iterated quotients of the inclusions $\Sh{V}_i\hookrightarrow \Sh{V''}_i\hookrightarrow \Sh{V}''_i$.
Therefore, the composition of $h_{s(a)}: \Sh{V}'_{s(a)}\to \Sh{V}''_{s(a)}$ with $n''_a: \Sh{V}''_{s(a)}\to \Sh{V}''_{t(a)}$ and then further with the quotient $\Sh{V}''_{t(a)}\to \Sh{V}''_{s(a)}/\Sh{V}'_{t(a)}$ is zero. Hence, there exists a map $n'_{s(a)}: \Sh{V}'_{s(a)}\to \Sh{V}'_{t(a)}$ fitting in the commutative diagram 
\[
\begin{tikzcd}
    0\arrow[r] & \Sh{V}'_{s(a)}\arrow[r, "h_{s(a)}"] \arrow[d, "n'_a"]& \Sh{V}''_{s(a)}  \arrow[d, "n''_{a}"]\arrow[r] & \Sh{V}''_{s(a)}/\Sh{V}'_{s(a)}\arrow[d, "\beta''_a"]\arrow[r] & 0 
    \\
     0\arrow[r] & \Sh{V}'_{t(a)}\arrow[r] & \Sh{V}''_{t(a)} \arrow[r] & \Sh{V}''_{t(a)}/\Sh{V}'_{t(a)} \arrow[r] & 0 
\end{tikzcd}
\]
Moreover, it is easy to check that the following diagram commutes 
\[
\begin{tikzcd}
     \Sh{V}_{s(a)}\arrow[r, "g_{s(a)}"] \arrow[d, "n_a"]& \Sh{V}'_{s(a)} \arrow[d, "n'_a"]\arrow[r, "h_{s(a)}"] & \Sh{V}''_{s(a)} \arrow[d, "n''_a"]
    \\
     \Sh{V}_{t(a)}\arrow[r, "g_{t(a)}"] & \Sh{V}''_{t(a)} \arrow[r, "h_{t(a)}"] & \Sh{V}''_{t(a)}/\Sh{V}_{t(a)}
\end{tikzcd}
\]
and that it recovers \eqref{eq: diagram ses quotients} by taking iterated quotients. 
Therefore this construction provides a map $s:(s^{(2)})^{-1}(0)\to \Hecke^{(2)}(\DD, N/G)$ which is manifestly an inverse of the morphism $\Hecke^{(2)}(\DD, N/G)\to (s^{(2)})^{-1}(0)$ defined above. This proves the first statement of the lemma. For the second statement, it suffices to observe that the front square is Cartesian by Proposition \ref{prop: virtual pullback data hall and Coulomb}, the back square is Cartesian by the argument above, and the left and top squares are Cartesian by definition. Therefore, the right square is also Cartesian.
\end{proof}

Combining Lemma \ref{lemma: first lemma factorization c^{(2)}} and Lemma \ref{lemma: second lemma factorization c^{(2)}}, we deduce that \eqref{eq: factorization c^{(2)}} provides a factorization of $c^{(2)}$ into the smooth morphism $b^{(2)}$ and the quasi-smooth embedding $a^{(2)}$. Therefore, we can define a canonical pullback morphism $\theta^{(2)}_{\dd',\dd''}$ fitting in the diagram
\begin{equation*}
    \begin{tikzcd}
    \HO^{\BM}(\FCoh_{Q,0,\dd',\dd''}(\DD), \BoQ)\arrow[r, equal]\arrow[d, "\theta^{(2)}_{\dd',\dd''}"] & \HO^{\BM}(\Map(\DD, N/G)\times \FCoh_{Q,0,\dd',\dd''}(\DD), \BoQ)\arrow[d, "(b^{(2)})^*"]
    \\
     \HO^{\BM}(\Hecke^{(2)}_{\dd', \dd''}(\DD, N/G), \BoQ^{\vir}) & \arrow[l, swap, "(a^{(2)})^!"] \HO^{\BM}(D_{\dd',\dd''}, \BoQ[-\vv^T\dd])
\end{tikzcd}
\end{equation*}
It generalizes \eqref{eq: diagram virtaul pull hall Coulomb app}. Here $\BoQ^{\vir}$ stands for the shifted constant sheaf $\BoQ[2s_{\dd}]$. We remark that, as usual, all these morphisms are understood in terms of sufficiency large approximations of the mapping stack $\Map(\DD, N/G)$, cf. \S\ref{subsec: cohomology of pro and ind stacks}. Once a finite type approximation $\Map(\DD_k, N/G)=N^k\fps/G^k\fps$ is fixed, it determines finite type approximations of the spaces and maps involved in the definitions of the constituents of the diagram above.  

The next two lemmas are the key technical ingredients for the proof of Theorem \ref{thm: main thm positive Coulomb hall app}. We refer to diagrams \eqref{eq: diagram Coulomb hall is mor algebra 1} and  \eqref{eq: diagram Coulomb hall is mor algebra 2} for the notation.

\begin{lemma}
\label{lemma: comp hecke Coulomb 3}
    The following diagram commutes
    \begin{equation*}
        \begin{tikzcd}
            \HO^{\BM}(\FCoh_{Q,0,\dd',\dd''}(\DD), \BoQ) \arrow[d, "\theta^{(2)}_{\dd',\dd''}"]\arrow[r, "\bar q_*\otimes \id"] & \HO^{\BM}(\FCoh_{Q,0,\dd}(\DD), \BoQ) \arrow[d, "\theta_{\dd}"] 
            \\
           \HO^{\BM}(\Hecke^{(2)}_{\dd',\dd''}(\DD, N/G),\BoQ^{\vir})\arrow[r, "(\tilde m^+)_*"] & \HO^{\BM}(\Hecke_{\dd}(\DD, N/G), \BoQ^{\vir}).
        \end{tikzcd}
    \end{equation*}
\end{lemma}
\begin{proof}
    Recall that $\theta_{\dd}=a^!\circ b^*$ and $\theta^{(2)}_{\dd',\dd''}=(a^{(2)})^!\circ (b^{(2)})^*$. By Lemma \ref{lemma: lemma comp hecke Coulomb 1} we get $b^* \circ \bar p_*= p'_*\circ (b^{(2)})^*$ Moreover, Lemma \ref{lemma: second lemma factorization c^{(2)}} implies that the virtual pullback data $(E_{\dd}, s)$ and $(E^{(2)}_{\dd',\dd''}, s^{(2)})$ are compatible, so Proposition \ref{prop: compat gysin and pull push Cartesian squares} gives $s^!\circ p'_*=(\tilde m_+)_*\circ (s^{(2)})^!$. The proof follows.
\end{proof}
In order to state the second lemma, let $(q_1\times q_2)^!$ be the composition $(q_1\times q_2)^!\coloneqq (\tilde q^{+,*})^{-1}\circ (\tilde p^+)^!$, where $(\tilde q^+)^*$ and $(\tilde p^+)^!$ are the morphisms introduced in \eqref{eq: gysin pullback Coulomb +} and \eqref{eq: gysin pullback Coulomb + 2}. To match the notation of this appendix with that in the references, one uses \eqref{eq: var triples and filt} and \eqref{eq: var triples and filt 2}.
\begin{lemma}
\label{lemma: comp hecke Coulomb 4}
    Then the following diagram commutes
    \begin{equation*}
        \begin{tikzcd}
            \HO^{\BM}(\FCoh_{Q,0,\dd'}(\DD),\BoQ)\otimes  \HO^{\BM}(\FCoh_{Q,0,\dd''}(\DD), \BoQ)\arrow[d, "\theta_{\dd'}\otimes \theta_{\dd''}"]\arrow[r, "\bar q^! \otimes \id"]& \HO^{\BM}(\FCoh_{Q,0,\dd',\dd''}(\DD), \BoQ)\arrow[d, "\theta^{(2)}_{\dd',\dd''}"]
            \\
            \HO^{\BM}(\Hecke_{\dd'}(\DD, N/G), \BoQ^{\vir}) \otimes \HO^{\BM}(\Hecke_{\dd''}(\DD, N/G), \BoQ^{\vir})  \arrow[r, "(q_1\times q_2)^!"]&  \HO^{\BM}(\Hecke^{(2)}_{\dd',\dd''}(\DD, N/G), \BoQ^{\vir})
        \end{tikzcd}
    \end{equation*}    
\end{lemma}
Before proving the lemma, we can easily observe how it implies Theorem \ref{thm: main thm positive Coulomb hall}. 

\begin{proof}[Proof of Theorem \ref{thm: main thm positive Coulomb hall}]
Combining Lemma \ref{lemma: comp hecke Coulomb 3} and Lemma \ref{lemma: comp hecke Coulomb 4} we obtain a commutative square 
\[
\begin{tikzcd}
			 \HO^{\BM}(\FCoh_{Q,0,\dd'}(\DD), \BoQ)\otimes \HO^{\BM}(\FCoh_{Q,0,\dd''}(\DD), \BoQ)\arrow[d, "\theta_{\dd'}\otimes \theta_{\dd''}"]\arrow[r, "\bar p_* \circ \bar q^!"] & \HO^{\BM}(\FCoh^0_{Q,\dd}(\DD), \BoQ)
        \arrow[d, "\theta_{\dd'+\dd''}"]
        \\
                \HO^{\BM}(\Hecke_{\dd'}(\DD, N/G), \BoQ^{\vir}) \otimes \HO^{\BM}(\Hecke_{\dd''}(\DD, N/G), \BoQ^{\vir})\arrow[r] & \HO^{\BM}(\Hecke_{\dd}(\DD, N/G), \BoQ^{\vir})
			\end{tikzcd}
\]
By definition, we have $\bar p_* \circ \bar q^!=\vmult^{\Xi_Q^{\onil}}$. On the other hand  the bottom horizontal map is the composition $(\tilde m^+)_*\circ (q_1\times q_2)^!$, which is, by definition, the Coulomb branch multiplication $\eta^+$.
\end{proof}

\subsection{Proof of Lemma \ref{lemma: comp hecke Coulomb 4}}

As a warm up, consider the proof in the case where $Q=A_1$, so that $\FCoh_{Q,0}(\DD)=\FCoh_{0}(\DD)$ and $\Bun_Q(\DD)=\Bun(\DD)$. In this case, all the morphisms in the outer frame of diagram \eqref{eq: diagram Coulomb hall is mor algebra 1} are smooth. Therefore, the Gysin pullbacks in the statement of Lemma \ref{lemma: comp hecke Coulomb 4} are ordinary pullbacks along smooth morphisms. Hence, the statement of the lemma follows from functoriality of smooth pullback in Borel-Moore cohomology. 

We now consider the general case. To remove clutter, we set $\Map=\Map(\DD, N/G)$, $\FCoh_{Q,0}(\DD)=\FCoh_{Q,0}$, $\Hecke(\DD, N/G)=\Hecke$, and similarly for the other spaces carrying $N/G$ as part of the notation. Firstly, we refine diagram \eqref{eq: diagram Coulomb hall is mor algebra 1} as follows
\begin{equation}
\label{eq: split diagram in two}
    \begin{tikzcd}             \Hecke^+_{\dd'}\times\Hecke_{\dd''}\arrow[d, swap, "c\times c"]  & & \Hecke^{(2)}_{\dd',\dd''}\arrow[ll, swap,  "q_1\times q_2"]\arrow[dd, "c^{(2)}"]\arrow[dl, "\bar c"]
        \\
        \Map\times \FCoh_{Q,0, \dd'}\times \Map \times \FCoh_{Q,0, \dd''} \arrow[d, "\pr_{134}"] &\Hecke_{\dd'}\times \FCoh_{Q,0,\dd''} \arrow[dl, "c\times \id"]& 
        \\
        \Map \times \FCoh_{Q,0, \dd'}\times \FCoh_{Q,0, \dd''}  & & \Map\times (\FCoh_{Q,0,\dd',\dd''})\arrow[ll, swap, "\id\times \bar q"]
   \end{tikzcd}
\end{equation}
The projection of $\bar c$ to $\Hecke_{\dd'}$ sends the a tuple $(g, h, n, n', n'')\in\Hecke^{(2)}_{\dd',\dd''}$ as in \eqref{eq: data parametrized by Filt^2} to $(g, n, n')$.
The projection of $\bar c$ to $\FCoh_{Q,0, \dd''}$ sends $(g, h, n, n', n'')$ to the induced quotient maps $\beta''=\{\beta''_a\}_{a\in Q_1}$, where $\beta''_a: \Sh{V}''_{s(a)}/\Sh{V}'_{s(a)}\to\Sh{V}''_{t(a)}/\Sh{V}'_{t(a)}$. The diagram is commutative.  By the usual strategy, we now define a virtual pullback 
\[
\bar c^!: \HO^{\BM}(\Hecke_{\dd'}\times \FCoh_{Q,0, \dd''}, \BoQ^{\vir})\to \HO^{\BM}(\Hecke^{(2)}_{\dd',\dd''}, \BoQ^{\vir}).
\]
In brief, we factor $\bar c$ as 
\begin{equation}
    \label{eq: factorization bar gr}
    \begin{tikzcd}
\bar E &  \bar D \arrow[l, swap, "\bar o"] 
\\
\bar D \arrow[u, "\bar s"]\arrow[d, swap, "\bar b"]& \arrow[l, swap,  "\bar a"]\Hecke^{(2)}_{\dd',\dd''}\arrow[dl, "\bar c"]\arrow[u, "a"]\\
   \Hecke_{\dd'}\times \FCoh_{Q,0}
\end{tikzcd}
\end{equation}
where $\bar D$ parametrizes maps $(g,h, n, n',n'', \beta'')$ fitting in the diagrams 
\begin{equation}
\label{eq: diagrams bar D}
  \begin{tikzcd}
     \Sh{V}_{s(a)}\arrow[r, "g_{s(a)}"] \arrow[d, "n_a"]& \Sh{V}'_{s(a)} \arrow[d, "n'_a"]
    \\
     \Sh{V}_{t(a)}\arrow[r, "g_{t(a)}"] & \Sh{V}''_{t(a)}
\end{tikzcd}
\qquad 
\begin{tikzcd}
   0\arrow[r]& \Sh{V}'_{s(a)}\arrow[r,  "h_{s(a)}"] & \Sh{V}''_{s(a)}  \arrow[r] \arrow[d, "n_a''"]& \Sh{V}''_{s(a)}/\Sh{V}'_{s(a)}\arrow[d, "\beta''_a"]\arrow[r] & 0
    \\
       0\arrow[r]& \Sh{V}'_{t(a)}\arrow[r, "h_{t(a)}"] & \Sh{V}''_{t(a)}\arrow[r] & \Sh{V}''_{t(a)}/\Sh{V}'_{t(a)}\arrow[r] & 0
\end{tikzcd}  
\end{equation}
and $\bar E\to \bar D$ is a vector bundle with fiber $\bigoplus_{a\in Q_1}\Hom(\Sh{V}_{s(a)}'', \Sh{V}_{t(a)}''/\Sh{V}_{t(a)}')$. The right diagram is not required to commute. For any $a\in Q_1$, the projection of the section $\bar s$ onto the summand $\Hom(\Sh{V}_{s(a)}'', \Sh{V}_{t(a)}''/\Sh{V}_{t(a)}')$ sends a closed point $(g,h, n, n', \beta'')\in \bar D$ to $\beta_{a}'' q_{s(a)}-q_{t(a)}n_{a}''$, where $q_i: \Sh{V}_i''\to  \Sh{V}_i''/ \Sh{V}_i'$ are the quotient maps induced by the inclusions $h_i$.
More formally, one can construct $\bar E\to \bar D$ arguing along the lines of Proposition \ref{prop: virtual pullback data hall and Coulomb}. We leave the details argument to the interested reader.
The maps $\bar a$ and $\bar b$ are the obvious ones and $\bar o$ is the zero section. The top square is Cartesian and the map $\bar b$ is smooth. Therefore, the pullback $\bar c^!$ is obtained as the composition $\bar c^!=\bar a^!\circ \bar b^*$ where $\bar b^*$ is smooth pullback and $\bar a^!$ is the virtual pullback from \S\ref{app: Refined Gysin pullback}.

The proof of the lemma reduces to proving the following two equalities
\begin{align}
    &(q_1\times q_2)^!\circ (c^!\boxtimes c^!) \circ (\pr_{134})^*=\bar c^!\circ (c^!\boxtimes \id) \label{eq: first of two equalities}
    \\
    & \bar c^!\circ (c^!\boxtimes \id)= (c^{(2)})^!\circ (\id\boxtimes \bar q^!), \label{eq: second of two equalities}
\end{align}
which are dealt with in Lemma \ref{lemma: first of two equalities} and Lemma \ref{lemma: second of two equalities} below.

\begin{lemma}
\label{lemma: first of two equalities}
=    Equality \eqref{eq: first of two equalities} holds.
\end{lemma}
\begin{proof}
Consider the quotient $G\fps\backslash \left(G\lfps\times_{G\fps} \Rt^+ \right)$. It is the moduli space parametrizing data $(g,h,n'', \beta'')$ as in the right diagram of \eqref{eq: diagrams bar D}, but now with the constrain that the rightmost square commutes. We denote by $G\fps\backslash \left(G\lfps^+_{\dd'}\times_{G\fps} \Rt^+_{\dd''} \right)$ the component such that the length vector of the cokernels of $g$ and $h$ are $\dd'$ and $\dd''$, respectively. The top diagonal triangle of \eqref{eq: split diagram in two} can be further expanded as 
\[
\begin{tikzcd}
& G\fps\backslash \Tt^+_{\dd'} \times\Hecke_{\dd''}\arrow[dd, near start, shift right=2, "\id\times c"]& G\fps\backslash \left(G\lfps^+_{\dd'}\times_{G\fps} \Rt^+_{\dd''} \right) \arrow[l, swap, "\tilde q_1\times \tilde q_2"]\arrow[dd, "\bar c'"]
\\
\Hecke_{\dd'}\times\Hecke^{+}_{\dd''} \arrow[dd, "\id\times c"]\arrow[ur, hookrightarrow] & \Hecke^{(2)}_{\dd',\dd''}\arrow[l, near end, swap, "q_1\times q_2"]\arrow[ur, hookrightarrow]\arrow[dd, near start, shift left=2, "\bar c"]
\\
& G\fps\backslash \Tt^+_{\dd'}\times \Map\times \FCoh_{Q,0, \dd''} & G\fps \backslash \Tt^+_{\dd'}\times \FCoh_{Q,0, \dd''}\arrow[l, near start, swap, "\delta'\times \id"]
\\
\Hecke_{\dd'}\times \Map\times \FCoh_{Q,0, \dd''} \arrow[d, swap, "c\times p_2"]\arrow[ur, hookrightarrow] &\Hecke_{\dd'}\times \FCoh_{Q,0, \dd''}\arrow[l, swap, "\delta\times \id"]\arrow[dl, "c\times \id"]\arrow[ur, hookrightarrow]
\\
\Map\times \FCoh_{Q,0, \dd'}\times \FCoh_{Q,0, \dd''} 
\end{tikzcd}
\]
Here:
\begin{enumerate}
    \item all the upward pointing maps are the canonical inclusions induced by 
    \[
    \Hecke_{\dd}=G\fps\backslash \Rt^+_{\dd}\hookrightarrow G\fps\backslash \Tt^+_{\dd}.
    \] 
    \item The map $\delta$ sends a tuple $(g, n, n')\in\Hecke_{\dd'}$ to $((g, n, n'), n')\in\Hecke_{\dd'}\times \Map$. 
    \item  The map $\delta'$ sends a tuple $(g, n')\in G\fps\backslash \Tt^+_{\dd'}$ to $((g, n'), n')$.
    \item The map $\bar c'$ sends a tuple $(g,h,n'', \beta'')\in G\fps\backslash \left(G\lfps^+_{\dd'}\times_{G\fps} \Rt^+_{\dd''} \right) $ to $((g, n''), \beta'')$.
    \item All the other maps are defined after diagrams \eqref{eq: diagram Coulomb hall is mor algebra 1} and \eqref{eq: diagram Coulomb hall is mor algebra 2}.
\end{enumerate}
 It is easy to check that the diagram is commutative. Moreover, it is clear that 
 \[
 (q_1\times q_2)^!\circ (c^!\boxtimes c^!) \circ (\pr_{134})^*=(q_1\times q_2)^!\circ (\id\boxtimes c^!) \circ (c^!\boxtimes p_2^*),
 \]
so \eqref{eq: first of two equalities} is equivalent to 
\begin{equation}
\label{eq: alt first of two equalities}
    (q_1\times q_2)^!\circ (\id\boxtimes c^!) \circ (c^!\boxtimes p_2^*)=\bar c^!\circ (c^!\boxtimes \id),
\end{equation}
which we now proceed to prove. Firstly, observe that since $\Map=N\fps/G\fps$ is topologically trivial, arguing as in Lemma \ref{lemma: Lemma pullback BFN} we deduce that 
\begin{equation*}
    \delta^* \DD\BoQ_{\Hecke_{\dd'}\times \Map\times \FCoh_{Q,0,\dd''}}=\DD\BoQ_{\Hecke_{\dd'}\times \FCoh_{Q,0,\dd''}}
\end{equation*}
Therefore, we have a well-defined pullback morphism 
\begin{equation}
    \label{eq: can pullback rgardless delta}
    \delta^*: \HO^{\BM}(\Hecke_{\dd'}\times \Map\times \FCoh_{Q,0,\dd''}, \BoQ) \xrightarrow{} \HO^{\BM}(\Hecke_{\dd'}\times \FCoh_{Q,0,\dd''}, \BoQ[-2\infty_{G,N}]).
\end{equation}
Moreover, the composition 
\[
\Hecke_{\dd'}\xrightarrow{\delta}\Hecke_{\dd'}\times \Map\xrightarrow{}\Hecke_{\dd'}
\]
is the identity, so we deduce that $(c^!\boxtimes \id)=(\delta^*\boxtimes \id )\circ(c^!\boxtimes p_2^*)$. Therefore, to prove \eqref{eq: alt first of two equalities} it suffices to show that
\begin{equation}
    \label{eq: alt 2 first of two equalities}
    (q_1\times q_2)^!\circ (\id\boxtimes c^!) =\bar c^!\circ (\delta^*\boxtimes \id ).
\end{equation}
The four pullbacks in this equation correspond to the four sides of the front square of the diagram above. To deduce the statement, we first prove a sheafified version of the statement for the back square and obtain the result by restricting to the front square and applying derived global sections. The precise argument goes as follows. It is easy to check that $\bar c'$ admits a two step factorization $\bar c'=\bar b'\circ \bar a'$ fitting in the diagram
\[
\begin{tikzcd}
    &  \bar D' \arrow[ddl, swap,  "\delta''"] \arrow[d, swap, "\bar b'"]& \arrow[l, swap,  "\bar a'"]G\fps\backslash \left(G\lfps^{+}_{\dd'}\times_{G\fps} \Rt^+_{\dd''} \right)\arrow[dl, swap, "\bar c' "]\arrow[ddl,  "\tilde q_1\times \tilde q_2"]
    \\
    & G\fps \backslash \Tt^+_{\dd'}\times \FCoh_{Q,0, \dd''}\arrow[ddl, near start, swap, "\delta' \times \id"]
    \\
    G\fps\backslash \Tt^+_{\dd'}\times D_{\dd''} \arrow[d, swap, "\id \times b"]& \arrow[l, swap,  "\id \times a"]G\fps\backslash\Tt^+_{\dd'}\times\Hecke_{\dd''}\arrow[dl, "\id \times c"]
    \\
    G\fps\backslash \Tt^+_{\dd'}\times \Map\times \FCoh_{Q,0, \dd''}
\end{tikzcd}
\]
All squares are Cartesian and, additionally, the map $\bar a'$ has virtual pullback data given by $(\delta'')^*(G\fps\backslash \Tt^+_{\dd'}\times E_{\dd''} , \id \times s)$, where $(E_{\dd},s)$ are the virtual pullback data for $c:\Hecke_{\dd''}\to \Map\times \FCoh_{Q,0, \dd''}$, see \S\ref{subsec: From hall to Coulomb} and the discussion at the beginning of \S\ref{subsec: APP Statement and conventions}. Therefore, $(\delta'')^*(G\fps\backslash \Tt^+_{\dd'}\times E_{\dd''} , \id \times s)$ and $G\fps\backslash \Tt^+\times E_{\dd''} , \id \times s)$ are compatible along $\delta''$, cf. Appendix \ref{app: Refined Gysin pullback}. Arguing as in the proof of Proposition \ref{prop: compat gysin and pull push Cartesian squares}, one obtains
\begin{equation}
    \label{eq: square 1}
\begin{tikzcd}
     \DD\BoQ_{ \Tt^+_{\dd'}\times D_{\dd''}}\arrow[d]\arrow[r] &(\delta'')_*\DD\BoQ_{\bar D'}[2\infty_{G,N}]\arrow[d]
     \\
    (\id\times a)_*\DD\BoQ_{\Tt^+_{\dd'}\times\Hecke_{\dd''}}\arrow[r] & (\delta'')_*(\bar a')_*\DD\BoQ_{ \left(G\lfps^{+}_{\dd'}\times_{G\fps} \Rt^+_{\dd''} \right)}[2\infty_{G,N}]
\end{tikzcd}
\end{equation}
Here, the horizontal maps are induced by the units $\id\to (\delta'')_* (\delta'')^*$ and $\id\to (\tilde q_1\times \tilde q_2)_*(\tilde q_1\times \tilde q_2)^*$ and the vertical ones are the sheaf version of the refined pullbacks \eqref{eq: sheafy refined Gysin} for $(G\fps\backslash \Tt^+_{\dd'}\times E_{\dd''} , \id \times s))$ and, $(\delta'')^*(G\fps\backslash \Tt^+_{\dd'}\times E_{\dd''} , \id \times s))$ respectively. Notice that by \eqref{eq: can pullback rgardless delta} and the Cartesian diagrams above (or by a direct argument as in Lemma \ref{lemma: Lemma pullback BFN}) it follows that 
\[
(\delta'')^* \DD\BoQ_{ \Tt^+_{\dd'}\times D_{\dd''}}=\DD\BoQ_{\bar B'}[2\infty_{G,N}],
\]
so the proof of Proposition \ref{prop: compat gysin and pull push Cartesian squares} generalizes to this case even though $\delta''$ is not smooth—as it would be required in its assumptions. Now, by functoriality, we have an obvious commutative diagram 
\begin{equation}
    \label{eq: square 2}
    \begin{tikzcd}
    \DD\BoQ_{ \Tt^+_{\dd'}\times \Map\times \FCoh_{Q,0, \dd''}}\arrow[r]\arrow[d] & (\delta'\times \id)_*\DD\BoQ_{ \Tt^+_{\dd''}\times \FCoh_{Q,0, \dd''}}[2\infty_{G,N}]\arrow[d]
    \\
    (\id\times b)_*\DD\BoQ_{ \Tt^+_{\dd'}\times D_{\dd''}}\arrow[r] &(\delta'\times \id)_*(\bar b')_*\DD\BoQ_{\bar D'}[2\infty_{G,N}]
\end{tikzcd}
\end{equation}
All the maps are induced by the units $\id \to x_*x^*$ for $x=b, \delta'', \delta'\times \id, \bar b'$. Applying $(\id\times b)_*$ to \eqref{eq: square 1} and composing with \eqref{eq: square 2} we get the diagram 
\begin{equation*}
    \begin{tikzcd}
    \DD\BoQ_{ \Tt^+_{\dd'}\times \Map\times \FCoh_{Q,0, \dd''}}\arrow[r]\arrow[d] & (\delta'\times \id)_*\DD\BoQ_{\backslash \Tt^+_{\dd'}\times \FCoh_{Q,0, \dd''}}[2\infty_{G,N}]\arrow[d]
    \\
    (\id\times b)_* (\id\times a)_*\DD\BoQ_{G\fps\backslash\Tt^+_{\dd'}\times\Hecke_{\dd''}}\arrow[r] & (\delta'')_*(\bar a')_*\DD\BoQ_{\fps\backslash \left(G\lfps^{+}_{\dd'}\times_{G\fps} \Rt^+_{\dd''} \right)}[2\infty_{G,N}]
\end{tikzcd}
\end{equation*}
Applying the functor 
\[
\left(\Hecke_{\dd'} \times \Map\times \FCoh_{Q,0, \dd''}\hookrightarrow G\fps\backslash \Tt^+_{\dd'}\times \Map\times \FCoh_{Q,0, \dd''}\right)^!
\]
and base changing, we get
\begin{equation*}
    \begin{tikzcd}
    \DD\BoQ_{\Hecke_{\dd'} \times \Map\times \FCoh_{Q,0, \dd''}}\arrow[r]\arrow[d] & (\delta'\times \id)_*\DD\BoQ_{\Hecke_{\dd'} \times \FCoh_{Q,0, \dd''}}[2\infty_{G,N}]\arrow[d]
    \\
    (\id\times b)_* (\id\times a)_*\DD\BoQ_{\Hecke_{\dd'}\times\Hecke^+_{\dd''}}\arrow[r] & (\delta'')_*(\bar a)_*\DD\BoQ_{\Hecke^{(2)}_{\dd',\dd''}}[2\infty_{G,N}]
\end{tikzcd}
\end{equation*}
We now take derived global sections. The bottom row gives, by definition, $(q_1\times q_2)^!$. The top row gives $\delta^*\boxtimes \id$, where $\delta$ is as in \eqref{eq: can pullback rgardless delta}. With a standard base change argument, one can check that the vertical maps give $c^!$ and $\bar c^!$. Therefore, commutativity of the square above implies \eqref{eq: alt 2 first of two equalities}. This completes to proof of the lemma.

\end{proof}

\begin{lemma}
\label{lemma: second of two equalities}
    Equality \eqref{eq: second of two equalities} holds.
\end{lemma}
\begin{proof}
    Consider the commutative diagram
    \begin{equation}
\label{eq: split diagram in two}
    \begin{tikzcd}              
   \Hecke_{\dd'}\times \FCoh_{Q,0, \dd''} \arrow[d, "c\times \id"]& \Hecke^{(2)}_{\dd',\dd''}\arrow[l, swap,  "\bar c"]\arrow[d, "c^{(2)}"]
        \\
        \Map \times \FCoh_{Q,0, \dd'}\times \FCoh_{Q,0, \dd''}  & \Map\times (\FCoh_{Q,0,\dd',\dd''})\arrow[l, swap, "\id\times \bar q"]
   \end{tikzcd}
\end{equation}
Each morphism is equipped with virtual pullback data inducing virtual pullbacks and we need to prove that $\bar c^!\circ (c^!\boxtimes \id)= (c^{(2)})^!\circ (\id\boxtimes \bar q^!)$. All the moduli spaces are of finite type and admit compatible derived enhancements inducing the virtual pullbacks in the equation above. Therefore, the equality follows by \ref{prop: quasismooth gysin functorial}.

\begin{remark}
	The compatibility of  the derived structures discussed in the proof of Lemma \ref{lemma: second of two equalities} can be seen very explicitly over a point in $ \Hecke^{(2)}_{\dd',\dd''}$. Indeed all the derived enhancements of the morphisms involved in the statement of the Lemma are quasi-smooth and the associated relative cotangent complexes, which are concentrated in degree $[-1,0]$ have fibers over $(g,h,n, n', n'')\in  \Hecke^{(2)}_{\dd',\dd''}$ given by the deformation/obstruction theories
	\begin{align*}
		\LL_{\id\times \bar q}&=[\Ext^1(\Sh{V}_{s(e)}''/\Sh{V}_{s(e)}', \Sh{V}'_{t(e)}/\Sh{V}_{t(e)})\xrightarrow{0} \Hom(\Sh{V}_{s(e)}''/\Sh{V}_{s(e)}', \Sh{V}'_{t(e)}/\Sh{V}_{t(e)})]
		\\
		\LL_{\bar c\times \id}&=[\Hom(\Sh{V}'_{s(e)}, \Sh{V}'_{t(e)}/\Sh{V}_{t(e)})\to 0]
		\\
		\LL_{\bar c^{(2)}}&=[\Hom(\Sh{V}''_{s(e)}, \Sh{V}''_{t(e)}/\Sh{V}_{t(e)})\to 0]
		\\
		\LL_{\bar c}&=[\Hom(\Sh{V}''_{s(e)}, \Sh{V}''_{t(e)}/\Sh{V}'_{t(e)})\to 0]
	\end{align*}
	In terms of this explicit presentation, the compatibility of the derived structures reduces to the classical exact sequences
	\begin{equation*}
		0\to \Hom(\Sh{V}''_{s(e)}, \Sh{V}'_{t(e)}/\Sh{V}_{t(e)})\to \Hom(\Sh{V}''_{s(e)}, \Sh{V}''_{t(e)}/\Sh{V}_{t(e)})\to \Hom(\Sh{V}''_{s(e)}, \Sh{V}''_{t(e)}/\Sh{V}'_{t(e)})\to 0
	\end{equation*}
	and 
	\begin{align*}
	0&\to \Hom(\Sh{V}_{s(e)}''/\Sh{V}_{s(e)}', \Sh{V}'_{t(e)}/\Sh{V}_{t(e)})\to\Hom(\Sh{V}_{s(e)}'', \Sh{V}'_{t(e)}/\Sh{V}_{t(e)}) \to \Hom(\Sh{V}_{s(e)}', \Sh{V}'_{t(e)}/\Sh{V}_{t(e)}) 
	\\
	&\to \Ext^1(\Sh{V}_{s(e)}''/\Sh{V}_{s(e)}', \Sh{V}'_{t(e)}/\Sh{V}_{t(e)})\to 0.
    \end{align*}
\end{remark}

\end{proof}   

\newpage

\section{The derived structure of the moduli space of quasimaps}
\label{app: derived quasimaps}

\subsection{Quasimaps as a derived stack}
Fix a smooth connected curve $C$ and a derived Artin stack $\FX$. Let $\DMap(C, \FX)$ be the derived stack of maps from $C$ to $\FX$. $\Map(C, \FX)$ is recovered as the classical truncation
\[
\Map(\bbP^1, \FX)=t_0(\DMap(C, \FX))\to \DMap(C, \FX).
\]
The functor $\DMap(C, \FX)$ assigns to any derived scheme $B$ the $\infty$-groupoid $\Map(C\times B, \FX)$. The derived stack $\DMap(C, \FX)$ is algebraic, see \cite[Sect. 2.2.6.3]{ToenVezzosi}.

Assume that $\FX$ is the derived enhancement of the GIT quotient of a possibly singular affine affine variety $Y$ by a reductive group $G$. We denote the GIT stability by $\zeta$ so that $X=Y^{\zeta\sst}/G \subset \FX$. Let also $\eta\in C$ be the generic point. We define the derived moduli stack  $\DQM(C, \FX)$ as the substack of $\DMap(C, \FX)$ whose functor of points $\DQM(C, \FX)(B)$ consists of those maps $f\in \Map(C, \FX)(B)$ such that for any geometric point $b:\Spec(\BoC)\to B$, the evaluation of $f$ factors through the semistable locus, i.e.
\begin{equation}
    \label{eq: quasimap stability}
    \begin{tikzcd}
        \Spec(\BoC)\times \eta \arrow[r, "b"]\arrow[drr, dashed,"\exists"]& B \times C  \arrow[r, "f"] & \FX
    \\
    && X\arrow[u]
    \end{tikzcd}
\end{equation}
 Notice that since any geometric point $b$ of $B$ factors as $\Spec(\BoC)\to t_0(B)\to B$, the quasimap condition is effectively a condition on the classical truncation of $\DQM(C, \FX)$. By definition, we have 
\[
\QM(C, X)=t_0(\DQM(C, X))
\]
and a fiber square 
\[
\begin{tikzcd}
    \QM(C, X) \arrow[d]\arrow[r] & \DQM(C, X)\arrow[d]
    \\
    \Map(C, \FX)\arrow[r] &\DMap(C, \FX).
\end{tikzcd}
\]
Fix closed points $\Spec(\BoC)\xrightarrow{p} C$ and $\Spec(\BoC)\xrightarrow{x} X$. We will be also interested in the derived  stack of based maps $\DMap_{p\to x}(C, \FX)$. Formally, the former is defined as the (derived) fiber product
\[
\DMap_{p \to x} \coloneq \DMap(C, \FX)\times_{\FX} x
\]
associated to the composition $\Spec(\BoC)\xrightarrow{x} X\to \FX$ and the evaluation $\ev_p: \DMap_{p\to x}(C, \FX)\to \FX$. The substack of based quasimaps $\DQM_{p\to x}(C, X)$ is defined in the same way.

Below we will focus on the case where $C=\bbP^1$ and the target $\FX$ is a quiver variety, namely when
 \[
 \FX=\mu^{-1}(0)/\GL_{\vv}\subset (T^*\Rep_Q(\vv,\ww))/\GL_{\vv}\qquad X=(\mu^{-1}(0))^{\zeta\sst}/\GL_{\vv}.
 \]
Explicitly, a object $f\in \DMap(\bbP^1\times B, \FX)$ consists of a tuple $(\Sh{V}_i, \Sh{I}_i, \Sh{J}_i, \Sh{B}_{a}, \Sh{B}_{a^*})$ where $\Sh{V}_i$ are vector bundles on $B \times \bbP^1$ and the remaining data are morphisms of sheaves
\begin{align*}
    &\Sh{B}_{a}: \Sh{V}_{s(a)}\to \Sh{V}_{t(a)} \qquad \Sh{J}_i: \Sh{V}_i \to W_i \otimes \Sh{O}_{B\times \bbP^1}\\
    & \Sh{B}_{a^*}: \Sh{V}_{t(a)}\to \Sh{V}_{s(a)}  \qquad \Sh{I}_i:   \Sh{O}_{B\times \bbP^1}^{\ww_i}\to \Sh{V}_i
\end{align*}
together with a choice of nullhomotopies 
\[
\sum_{a} \comm{\Sh{B}_{a}}{\Sh{B}_{a^*}}_i+\Sh{I}_i \Sh{J}_i\simeq 0
\]
in the space of maps $\Map(\Sh{V}_i,\Sh{V}_i)$, for all $i\in Q_0$. In the case of a based quasimap $\DMap_{\infty\to x}(\bbP^1, \FX)$ to a point $\xi =(\xi_2, \xi_3, \iota, \jmath ) \in X$, an object is further specified by choices of trivializations 
\begin{equation}
\label{eq: based quasimaps trivializes byndles at infinity}
    \Sh{V_i}|_{\{\infty\}\times \bbP^1}\simeq V_i\otimes \Sh{O}_B
\end{equation}
identifying the restrictions of the tuple $(\Sh{I}_i, \Sh{J}_i, \Sh{B}_{a}, \Sh{B}_{a^*})|_{\{\infty\}\times \bbP^1}$ with $(\xi_2, \xi_3, \iota, \jmath) \otimes\id$.

Assume now that we fix the stability condition $\zeta$ on the quiver variety and consider the derived stack $\DQM(\bbP^1, X)$. Then an object $f\in \DQM(\bbP^1, X)(B)$ is just a tuple $(\Sh{V}_i,\Sh{I}_i, \Sh{J}_i, \Sh{B}_{a}, \Sh{B}_{a^*})$ satisfying the quasimap stability \eqref{eq: quasimap stability}. In particular, for $\zeta=\zeta^-$ as in \S\ref{subsec: Nakajima quiver varieties}, for any geometric point $\Spec(\BoC)\xrightarrow{p} B$ we have
\begin{equation}
    \label{eq: positive homology vaishing quasimaps}
    H^0(\bbP^1, \Sh{V}_i|_{b\times \bbP^1}(-k))=0
\end{equation}
for all $i\in Q_0$ and $k > 0$.
\subsection{Quiver quasimaps as a derived stack}

Consider the smooth Artin stack,
\[
\FM_{Q^{\QM}_{\ww,\vv}}(\mathbf{d}) = \rep_{Q^{\QM}_{\ww,\vv}}(\dd,1)/\GL(\dd) \ltimes \Hom(\BoC^{\dd}, \BoC^{\vv}),
\]
which in this section we regard as a derived stack. Let $ \FM_{Q^{\QM}_{\ww,\vv}}^{\zeta^+\sst}(\mathbf{d})$ be the open substack of $\zeta^+$-semistable representations as is \S\ref{subsec: quasimaps to quiver varieties}.
Let $\Dcrit_{\dd}(W)$ be derived critical locus of the quasimap potential $W: \FM_{Q^{\QM}_{\ww,\vv}}(\mathbf{d})\to \BoC$. Formally, $\Dcrit_{\dd}(W)$ is the derived fiber product 
\[
\Dcrit_{\dd}(W)=\FM_{Q^{\QM}_{\ww,\vv}}(\mathbf{d}) \times_{T^* \FM_{Q^{\QM}_{\ww,\vv}}(\dd)} \FM_{Q^{\QM}_{\ww,\vv}}(\dd)
\] 
defined with respect to the sections $dW, 0: \FM_{Q^{\QM}_{\ww,\vv}}(\mathbf{d})\to T^* \FM_{Q^{\QM}_{\ww,\vv}}(\mathbf{d})$.
Let also $\Dcrit_{\dd}(W)^{\zeta^+\sst}$ be the substack of semistable points, fitting in the fiber diagram
\[
\begin{tikzcd}
    \Dcrit_{\dd}(W)^{\zeta^+\sst} \arrow[r]\arrow[d] & \Dcrit_{\dd}(W)\arrow[d]
    \\
    \FM_{Q^{\QM}_{\ww,\vv}}(\mathbf{d})^{\zeta\sst} \arrow[r] & \FM_{Q^{\QM}_{\ww,\vv}}(\mathbf{d})
\end{tikzcd}
\]
As a functor of points, $\Dcrit_{\dd}(W)$ assigns to a derived scheme $B$ the groupoid $\FM_{Q^{\QM}_{\ww,\vv}}(\mathbf{d})(S)$ whose objects consists of the following data 
\begin{enumerate}
    \item \label{quiver quasimap data 1} Vector bundles $\Sh{F}_i$ of rank $\dd_i + \vv_i$ on $B$ together with a morphism $\Sh{O}_B^{\vv_i}\to \Sh{F}_i$ such that $\Sh{U}_i\coloneqq\cofib(\Sh{O}_B^{\vv_i}\to \Sh{F}_i)$ is a rank $\dd_i$ vector bundle.
    \item \label{quiver quasimap data 2} Morphisms 
    \[
    \mathbb{B}_{a}: \Sh{F}_{s(a)}\to \Sh{F}_{t(a)}\qquad \mathbb{B}_{a^*}: \Sh{F}_{s(a^*)}\to \Sh{F}_{t(a^*)}\qquad \zeta_i: \Sh{F}_i\to \Sh{U}_i\qquad  \mathbb{J}_i: \Sh{F}_i\to \Sh{O}_B^{\ww_i}
    \]
    satisfying a number of relations. These relations are the same as those in \cite[Lemma 3.4]{quasimapcrit}, except they are now understood in the homotopical sense.
\end{enumerate}

\subsection{An equivalence of derived stacks}
Fix $\infty \in \bbP^1$ and a chosen point $\xi = (\xi_2, \xi_3, \iota, \jmath) \in X$. The main result of \cite{quasimapcrit} constructs an equivalence of classical stacks
\begin{equation}
    \label{eq: spencer's result classical}
    \QM_{\infty\to \xi}(\bbP^1, X)\simeq  \crit_{\dd}(W)^{\zeta^+\sst}
\end{equation}
In this section, we prove that this equivalence lifts to the derived enhancements of these moduli spaces described above. More precisely, we prove the following.
\begin{theorem}
\label{thm: spencer's result derived}
    There exists an equivalence of derived stacks 
    \[
     \DPsi: \DQM_{\infty\to \xi}(\bbP^1, X)\simeq  \Dcrit_{\dd}(W)^{\zeta^+\sst}
    \]
    recovering \eqref{eq: spencer's result classical} upon taking classical truncations.
\end{theorem}
\begin{proof}
    Firstly, we construct a morphism $\DPsi: \DQM_{\infty\to \xi}(\bbP^1, X)\to  \Dcrit_{\dd}(W)^{\zeta^+\sst}$ recovering \eqref{eq: spencer's result classical} when passing to the classical truncations. This is an adaptation of the argument in \cite{quasimapcrit}. Fix a derived scheme $B$ and consider the projections $\pi_{12}, \pi_{13}: B\times \bbP^1\times \bbP^1\to B\times \bbP^1$. We further introduce projections $p: B\times \bbP^1\to B$ and $q: B\times \bbP^1\to \bbP^1$. 
    Consider the resolution of the diagonal 
    \begin{equation}
        \label{eq: res diag}
        \Sh{O}_{\bbP^1}(-1)\boxtimes \Sh{O}_{\bbP^1}(-1)\to \Sh{O}_{\bbP^1\times \bbP^1}\to \Sh{O}_{\Delta}
    \end{equation}
    seen as a triangle in the stable $\infty$-category $\Coh(\bbP^1)$. Fix a vector bundle $\Sh{V}$ on $B\times \bbP^1$ and denote by $\Sh{V}(-k)$ the twist by $q^*\Sh{O}_{\mathbb{P}^1}(-k)$. Pulling back \eqref{eq: res diag} by $B\times \bbP^1\times \bbP^1\to \bbP^1\times \bbP^1$ and applying the functor 
    \[
    (\pi_{13})_*((\pi_{12})^*\Sh{V}(-1)\otimes -)
    \]
    we obtain, via base-change, a distinguished triangle 
    \[
    p_*(\Sh{V}(-2))\boxtimes \Sh{O}_{\bbP^1}(-1) \to p_*(\Sh{V}(-1))\boxtimes \Sh{O}_{\bbP^1} \to \Sh{V}(-1)
    \]
    which we further twist to obtain the triangle
    \[
    p_*(\Sh{V}(-2))\boxtimes \Sh{O}_{\bbP^1} \to p_*(\Sh{V}(-1))\boxtimes \Sh{O}_{\bbP^1}(1) \to \Sh{V}
    \]
    Notice that since the morphism  $p$ is perfect of relative dimension $1$, the pushforwards $p_*(\Sh{V}(-2))$ and $p_*(\Sh{V}(-1))$ are perfect complexes of amplitude $[0,1]$. We can now define the functor $\DPsi$. Fix an object $(\Sh{V}_i,\Sh{I}_i, \Sh{J}_i, \Sh{B}_{a}, \Sh{B}_{a^*})$ in the groupoid $\DQM_{\infty\to \xi}(\bbP^1, X)(B)$. By \eqref{eq: positive homology vaishing quasimaps} we deduce that each vector bundle $\Sh{V}_i$ fits in an exact triangle 
    \[
    \Sh{V}_i\to p_*(\Sh{V}_i(-2))[-1]\boxtimes \Sh{O}_{\bbP^1} \to p_*(\Sh{V}_i(-1))[-1]\boxtimes \Sh{O}_{\bbP^1}(1)
    \]
    where $p_*(\Sh{V_i}(-k))[-1]$ is a vector bundle (concentrated in degree zero) for $k=1,2$. Pulling back this triangle to $B\times \{\infty\}$ and using \eqref{eq: based quasimaps trivializes byndles at infinity}, we get a triangle 
    \[
    V_i\otimes\Sh{O}_B\to \Sh{F}_i\to\Sh{U}_i
    \] 
    where each term is a vector bundle. Hence, to any $(\Sh{V}_i,\Sh{I}_i, \Sh{J}_i, \Sh{B}_{a}, \Sh{B}_{a^*})$ in $\DQM_{\infty\to \xi}(\bbP^1, X)(B)$ we can associate a morphism $V_i\otimes\Sh{O}_B\to \Sh{F}_i$ satisfying the conditions of \eqref{quiver quasimap data 1}. Moreover, using the fact the the definition of the triangle above is functorial, the data $\Sh{I}_i, \Sh{J}_i, \Sh{B}_{a}, \Sh{B}_{a^*}$ induce morphisms $\mathbb{B}_{a},\mathbb{B}_{a^*},\zeta_i,\mathbb{J}_i$ satisfying the condition of \eqref{quiver quasimap data 2}. The precise construction is a straightforward adaptation of the argument of \cite[Prop 3.5]{quasimapcrit}, so we skip it. 
   This analysis produces a functor $\DPsi: \DQM_{\infty\to \xi}(\bbP^1, X)\to \Dcrit(W)$ which restricts to $\Psi: \QM_{\infty\to \xi}(\bbP^1, X)\to \crit(W)$ on the classical truncations. Therefore, since semistability is a condition on classical points, we conclude that $\DPsi$ factors through $\Dcrit(W)^{\zeta^+\sst}$, and hence defines the sought after functor 
   \[
   \DPsi: \DQM_{\infty\to \xi}(\bbP^1, X)\to \Dcrit(W)^{\zeta^+\sst}
   \]
   Additionally, since $\Psi$ is an equivalence, to deduce that $\DPsi$ is also an equivalence, it suffices to check that $\DPsi$ is smooth, i.e. that the relative cotangent complex $\LL_{\DPsi}$ is perfect of amplitude $[0,0]$. In Lemma \ref{lemma: derived etale} below, which completes the proof, we actually prove the stronger statement that $\DPsi$ is \'etale.
\end{proof}

\begin{proposition}
\label{prop: derived etale}
    The morphism $\Psi$  is \'etale, i.e. $\LL_{\DPsi}=0$.
\end{proposition}
\begin{proof}
It suffices to show that that the leftmost morphism fitting in the triangle 
\[
\DPsi^*\LL_{\Dcrit(W)^{\zeta^+\sst}}\to \LL_{ \DQM_{\infty\to \xi}(\bbP^1, X)}\to  \LL_{\DPsi}
\]
is a quasi-isomorphism. Equivalently, we will show that the dual 
\begin{equation} \label{eq: comparison tangent complexes at geometric points}
d\DPsi: \TT_{\DQM_{\infty \to \xi}(\bbP^1, X)} \to \DPsi^* \TT_{\Dcrit(W)^{\zeta^+\sst}}
\end{equation}
is an isomorphism. It suffices to verify that $d\DPsi$ induces a quasi-isomorphism between tangent complexes at all geometric points $\Spec(\BoC)\to \DQM_{\infty\to \xi}(\bbP^1, X)$. This amounts to a straightforward analysis of the deformation theory of quasimaps that will occupy the rest of this appendix.
\end{proof}

Everything that follows is written at an arbitrary geometric point of $\DQM$, at which we have fixed quasimap and quiver data $\Psi(\Sh{V}, \Sh{B}_2, \Sh{B}_3, \Sh{I}, \Sh{J}) = (B_1, I, B_2, J_3, B_3, J_2, \alpha)$.

\subsection{Preliminaries}
Recall from \cite{quasimapcrit} that, for the bundle $\Sh{V}$ entering the quasimap data, we have a resolution 
\begin{equation} \label{eq: resolve V}
\Sh{V}(-1) \simeq \Bigg[\begin{tikzcd} F \otimes \Sh{O}_{\mathbb{P}^1}(-1) \arrow[r] &  U \otimes \Sh{O}_{\mathbb{P}^1} \end{tikzcd} \Bigg].
\end{equation}
It will be useful to note that this resolution is acyclic for the functor of global sections $\Gamma = \Gamma(\mathbb{P}^1, -)$, because $H^1(\mathbb{P}^1, \Sh{O}_{\mathbb{P}^1}) = H^1(\mathbb{P}^1, \Sh{O}_{\mathbb{P}^1}(-1)) = 0$. 

Applying the exact functor $\mathscr{H}om(\Sh{V}, -)$ to the above, we find another resolution 
\begin{equation} \label{eq: resolve Hom}
\mathscr{H}om(\Sh{V}, \Sh{V}(-1)) \simeq \Bigg[\begin{tikzcd} F \otimes \Sh{V}^*(-1) \arrow[r] &  U \otimes \Sh{V}^* \end{tikzcd} \Bigg].
\end{equation}
Note this is also acyclic for $\Gamma(\mathbb{P}^1, -)$, because $H^1(\mathbb{P}^1, \Sh{V}^*(-1)) \simeq H^0(\mathbb{P}^1, \Sh{V}(-1))^* = 0$ and $H^1(\mathbb{P}^1, \Sh{V}^*) \simeq H^0(\mathbb{P}^1, \Sh{V}(-2))^* = 0$ by Serre duality and Lemma \ref{lemma negative vector bundle}.

\subsection{Proof of Proposition \ref{prop: derived etale}}
This will amount to a straightforward analysis of the quasimap tangent complex $\TT_{\DQM_{\infty \mapsto \xi}}$ using the resolutions \eqref{eq: resolve V} and \eqref{eq: resolve Hom}. 

\subsubsection{Quasimap tangent complex}
Recall $\DQM$ fits into a diagram 
\[
\begin{tikzcd}
    \DQM(\mathbb{P}^1, X) \times \mathbb{P}^1 \arrow[d, "p"] \arrow[r, "u"] & \FX \\
    \DQM(\mathbb{P}^1, X)
\end{tikzcd}
\]
where $p$ is the canonical projection and $u$ is the universal evaluation. Since $\DQM$ is open in the mapping stack $\DMap(\mathbb{P}^1, \FX)$, its tangent complex may be presented as 
\[
\TT_{\DQM} = Rp_* Lu^* \TT_{\FX}. 
\]
As $\FX$ is a symplectic reduction of $T^* \text{Rep}_{Q}(\vv, \ww)$, we have the following model for its tangent complex: 
\[
\TT_{\FX} \simeq \Bigg[ \begin{tikzcd}
    \mathfrak{g}_{\FX} \arrow[r, "\text{act}"] & T_{T^*\text{Rep}_Q(\vv, \ww)} \arrow[r, "d\mu"] & \mathfrak{g}^*_{\FX}
\end{tikzcd} \Bigg]
\]
where $\mathfrak{g}_{\FX}$ is the bundle on $\FX$ associated to the adjoint representation of $\mathfrak{g} = \text{Lie} (G_{\vv})$ (likewise, $\mathfrak{g}_{\FX}^*$ with the coadjoint), the terms are in degrees $[-1, 0, 1]$, the first map is the differential of the action map, and the second map is the differential of the moment map. Passing to based quasimaps $\DQM_{\infty \mapsto \xi}$ and base changing to our given geometric point in $\DQM_{\infty \mapsto \xi}$, using the above model for $\TT_{\FX}$ to compute the derived pullback we find for the pointwise tangent complex: 
\begin{equation} \label{eq: quasimap tangent complex}
\TT_{\DQM_{\infty \mapsto \xi}} \simeq R \Gamma \Bigg[ \begin{tikzcd} \mathscr{H}om(\Sh{V}, \Sh{V}(-1)) \arrow[r] &  T^* \mathscr{R}ep \otimes \Sh{O}_{\mathbb{P}^1}(-1) \arrow[r] & \mathscr{H}om(\Sh{V}, \Sh{V}(-1)) \end{tikzcd} \Bigg]
\end{equation}
where $\Gamma = \Gamma(\mathbb{P}^1, -)$ is the functor of global sections on $\mathbb{P}^1$. We have introduced the abbreviation ($f$ denotes the quasimap regarded as a map to $\FX$, the pointwise restriction of $u$):
\[
T^* \mathscr{R}ep \coloneqq f^*T_{T^* \text{Rep}_Q(\vv, \ww)} \simeq \Sh{V} \otimes W^*  \oplus \Sh{V}^* \otimes W \oplus \mathscr{H}om(\Sh{V}, \Sh{V}) \otimes \Sh{O}_{\mathbb{P}^1}^{\oplus 2}. 
\]

\subsubsection{Quasimap tangent complex in quiver language}
Now we use the acyclic resolutions \eqref{eq: resolve V}, \eqref{eq: resolve Hom} to compute $R \Gamma$. We can take direct sums of \eqref{eq: resolve V}, \eqref{eq: resolve Hom} to get a resolution of $T^* \mathscr{R}ep(-1)$ (note we do not need to resolve the summand $W \otimes \Sh{V}^*(-1)$ in $T^* \mathscr{R}ep(-1)$ since it is already $\Gamma$-acyclic). Using, from Serre duality and the defintion of $\Psi$ in \cite{quasimapcrit}, 
\[
\begin{split}
    H^0(\mathbb{P}^1, \Sh{V}^*(-1)) & \simeq H^1(\mathbb{P}^1, \Sh{V}(-1))^* = U^* \\ 
    H^0(\mathbb{P}^1, \Sh{V}^*) & \simeq H^1(\mathbb{P}^1, \Sh{V}(-2))^* = F^*
\end{split}
\]
we find the following double complex model for $\TT_{\DQM_{\infty \mapsto \xi}}$:
\begin{equation} \label{eq: quasimap tangent double complex}
\TT_{\DQM_{\infty \mapsto \xi}} \simeq \Bigg[
\begin{tikzcd}
    \text{Hom}(F, U) \arrow[r] & \text{Hom}(W, U) \oplus \text{Hom}(F, U) \otimes \mathbb{C}^2 \arrow[r] & \text{Hom}(F, U) \\
    \text{Hom}(U, F) \arrow[u] \arrow[r] & \text{Hom}(U, W) \oplus \text{Hom}(U, F) \otimes \mathbb{C}^2 \arrow[u] \arrow[r] & \text{Hom}(U, F) \arrow[u]
\end{tikzcd} \Bigg].
\end{equation}
The vertical differential is read off from \eqref{eq: resolve V}, \eqref{eq: resolve Hom}, and the horizontal differential may be computed following the logic of \cite[Lemma 2.2]{quasimapcrit}. Computation of this differential will essentially imply Proposition \ref{prop: derived etale}. As a preliminary observation, note that if we regard this double complex as an ordinary complex by summing along the antidiagonals, the degree 0 and degree 1 terms may be identified canonically with 
\[
\begin{split}
T \text{Rep}_{Q^{\QM}} & \simeq \text{Hom}(F, U) \oplus \text{Hom}(U, F) \otimes \mathbb{C}^2 \oplus \text{Hom}(U, W) \\ 
T^* \text{Rep}_{Q^{\QM}} & \simeq \text{Hom}(U, F) \oplus \text{Hom}(F, U) \otimes \mathbb{C}^2 \oplus \text{Hom}(W, U).
\end{split}
\]

\subsubsection{Quiver tangent complex}
Because $\Dcrit(W)$ is the derived critical locus of a function on a quotient stack, we have a canonical chain complex model for $\TT_{\Dcrit(W)}$. It is given by the following:
\begin{equation} \label{eq: quiver tangent complex}
\TT_{\Dcrit(W)} \simeq \Bigg[ \begin{tikzcd} \text{Hom}(U, F) \arrow[r] &  T \text{Rep}_{Q^{\QM}} \arrow[r, "\text{Hess}(W)"] &  T^* \text{Rep}_{Q^{\QM}} \arrow[r] & \text{Hom}(F, U) \end{tikzcd} \Bigg].
\end{equation}
Here we view $\text{Hom}(U, F)$ as the Lie algebra of the group $GL(U) \ltimes \text{Hom}(U, V)$, the first map in the complex is the infinitesimal action map, and the last map is the dual of the infinitesimal action map. This complex sits in degrees $[-1, \dots, 2]$. 

\subsubsection{Comparing the complexes}
The following implies Proposition \ref{prop: derived etale}. 

\begin{proposition} \label{prop: quiver tangent complex=quasimap tangent complex}
The total complex of the double complex \eqref{eq: quasimap tangent double complex} is isomorphic to \eqref{eq: quiver tangent complex}. 
\end{proposition}
Explicitly, the isomorphism just makes some sign adjustments. The proof will be a little lengthy, but amounts to no more than unwinding definitions using the quiver construction set up in \cite{quasimapcrit}, to explicitly compare the differentials in \eqref{eq: quasimap tangent double complex}, \eqref{eq: quiver tangent complex}. 

\begin{proof}
Let's first deal with the vertical differentials in \eqref{eq: quasimap tangent double complex}. Recall from \cite{quasimapcrit} that the differentials in \eqref{eq: resolve V}, \eqref{eq: resolve Hom} are described in the quiver language as the natural maps induced by
\[
\begin{pmatrix} B_1 && I\end{pmatrix} - z \pi
\]
where $\pi: F \to U$ is the natural projection, and $z$ is the affine coordinate on $\mathbb{P}^1$ regarded as a global section of $\Sh{O}_{\mathbb{P}^1}(1)$. The vertical differential in \eqref{eq: quasimap tangent double complex} is obtained by applying the functor $H^0(\mathbb{P}^1, -)$ to these morphisms of sheaves and taking direct sums. Using the results of \cite{quasimapcrit}, these can be read off as follows. The leftmost vertical map is 
\[
\begin{split}
\text{Hom}(U, F) \ni \begin{pmatrix} X \\ Y \end{pmatrix} & \mapsto   \begin{pmatrix} B_1 && I \end{pmatrix} \begin{pmatrix} X \\Y  \end{pmatrix} \begin{pmatrix} 1 && 0 \end{pmatrix} - \begin{pmatrix} 1 && 0 \end{pmatrix} \begin{pmatrix} X \\ Y\end{pmatrix} \begin{pmatrix} B_1 && I\end{pmatrix} \\
& = \begin{pmatrix} \comm{B_1}{X} + IY && - XI \end{pmatrix} \in \text{Hom}(F, U)
\end{split}
\]
where we are recalling from \cite{quasimapcrit} that the projection $\pi = \begin{pmatrix} 1 & 0 \end{pmatrix}$ in coordinates, and is the map induced in $H^1$ by the tautological inclusion $\Sh{V}(-1) \hookrightarrow \Sh{V}$, while multiplication by $z: \Sh{V}(-2) \to \Sh{V}(-1)$ induces $\begin{pmatrix} B_1 & I \end{pmatrix}$ in $H^1$; we get the dual version of these maps since we are in $H^0$ of $\Sh{V}^*$ rather than $H^1$ of $\Sh{V}$, by Serre duality as above. 

Let us introduce the notation 
\[
(\delta \alpha, \delta B_2, \delta J_3, \delta B_3, \delta J_2) \in \text{Hom}(U, W) \oplus \text{Hom}(U, F) \oplus \text{Hom}(U, F)
\]
for an element of the middle term in the bottom row of \eqref{eq: quasimap tangent double complex}. By the same rationale, the vertical map at the middle entry is the direct sum of maps sending
\[
\begin{split}
    \text{Hom}(U, F) \ni \begin{pmatrix} \delta B_2 \\ \delta J_3 \end{pmatrix} & \mapsto \begin{pmatrix} \comm{B_1}{\delta B_2} + I \delta J_3 && - \delta B_2 I \end{pmatrix} \in \text{Hom}(F, U) \\
    \text{Hom}(U, F) \ni \begin{pmatrix} \delta B_3 \\ -\delta J_2 \end{pmatrix} & \mapsto -\begin{pmatrix} \comm{\delta B_3}{B_1} + I \delta J_2 &&  \delta B_3 I \end{pmatrix} \in \text{Hom}(F, U) 
\end{split}
\]
and acts by $0$ on $\delta \alpha$. The sign on $\delta J_2$ is to match the conventions of \cite{quasimapcrit}. Finally, the rightmost vertical map is literally the same as the leftmost vertical map. 

Now we deal with the horizontal maps; first we unwind the definition. The action map 
\[
\mathfrak{g}_{\FX} \to T_{T^* \Rep_{Q}(\vv, \ww)}
\]
induces, via $f^*$, a morphism of sheaves 
\[
\mathscr{H}om(\Sh{V}, \Sh{V}(-1)) \to T^* \mathscr{R}ep(-1)
\]
fitting into a commutative diagram 
\[
\begin{tikzcd}
    \mathscr{H}om(\Sh{V}, \Sh{V}(-1)) \arrow[r] \arrow[d] & \mathscr{H}om(\Sh{V}, F \otimes \Sh{O}_{\mathbb{P}^1}(-1)) \arrow[r] \arrow[d] & \mathscr{H}om(\Sh{V}, U \otimes \Sh{O}_{\mathbb{P}^1}) \arrow[d] \\ 
    T^* \mathscr{R}ep \otimes \Sh{O}_{\mathbb{P}^1}(-1) \arrow[r] & \langle \Sh{V}, W \rangle \oplus \langle W, F \rangle \oplus \langle \Sh{V}, F \rangle \otimes \mathbb{C}^2 \arrow[r] & \langle W, U \rangle \oplus \langle \Sh{V}, U \rangle \otimes \mathbb{C}^2
\end{tikzcd}
\]
where we have abbreviated e.g. $\langle \Sh{V}, F \rangle \coloneqq \mathscr{H}om(\Sh{V}, F \otimes \Sh{O}_{\mathbb{P}^1}(-1))$. The horizontal arrows on the left hand side of \eqref{eq: quasimap tangent double complex} are obtained by applying $H^0(\mathbb{P}^1, -)$ to the second and third downward arrows in the above diagram. On sheaf level (i.e. for local sections), the action map is $\epsilon \mapsto (\comm{\Sh{B}_2}{\epsilon}, \comm{\Sh{B}_3}{\epsilon}, - \epsilon \Sh{I}, \Sh{J} \epsilon)$ in terms of the quasimap sections $(\Sh{B}_2, \Sh{B}_3, \Sh{I}, \Sh{J})$. Unravelling the definitions and using the matrix descriptions of the quasimap sections from \cite{quasimapcrit}, we arrive at the following description of the bottom left horizontal differential in \eqref{eq: quasimap tangent double complex}. It has $\delta \alpha$ component 
\[
\text{Hom}(U, F) \ni \begin{pmatrix} X \\ Y \end{pmatrix} \mapsto \begin{pmatrix} \alpha && \jmath \end{pmatrix} \begin{pmatrix} X \\ Y\end{pmatrix} = \alpha X + \jmath Y \in \text{Hom}(U, W)
\]
and $(\delta B_2, \delta J_3, \delta B_3, \delta J_2)$ components 
\[
\begin{split}
\text{Hom}(U, F) \ni \begin{pmatrix} X \\ Y \end{pmatrix} & \mapsto \begin{pmatrix} B_2 && 0 \\ J_3 && \xi_2 \end{pmatrix} \begin{pmatrix} X \\ Y \end{pmatrix} - \begin{pmatrix} X \\ Y \end{pmatrix} B_2 \\ 
& = \begin{pmatrix} \comm{B_2}{X} \\ J_3 X + \xi_2 Y - Y B_2 \end{pmatrix} \in \text{Hom}(U, F) \\ 
\text{Hom}(U, F) \ni \begin{pmatrix} X \\ Y \end{pmatrix} & \mapsto \begin{pmatrix} B_3 && 0 \\ -J_2 && \xi_3 \end{pmatrix} \begin{pmatrix} X \\ Y \end{pmatrix} - \begin{pmatrix} X \\ Y\end{pmatrix}B_3 \\
& = \begin{pmatrix} \comm{B_3}{X} \\ -(J_2 X - \xi_3 Y + Y B_3) \end{pmatrix} \in \text{Hom}(U, F). 
\end{split}
\]
Likewise, the top left horizontal differential in \eqref{eq: quasimap tangent double complex} has $\text{Hom}(W, U)$-component 
\[
\text{Hom}(F, U) \ni \begin{pmatrix} \delta B_1 && \delta I\end{pmatrix} \mapsto -\begin{pmatrix} \delta B_1 && \delta I \end{pmatrix} \begin{pmatrix} 0 \\ \iota \end{pmatrix} = -\delta I \iota \in \text{Hom}(W, U)
\]
and $\text{Hom}(F, U) \otimes \mathbb{C}^2$-components 
\[
\begin{split}
\text{Hom}(F, U) \ni \begin{pmatrix} \delta B_1 && \delta I \end{pmatrix} & \mapsto B_2 \begin{pmatrix} \delta B_1 && \delta I \end{pmatrix} - \begin{pmatrix} \delta B_1 && \delta I \end{pmatrix} \begin{pmatrix} B_2 && 0 \\ J_3 && \xi_2 \end{pmatrix} \\
& = -\begin{pmatrix} \comm{\delta B_1}{B_2} + \delta I J_3 && -B_2\delta I + \delta I \xi_2 \end{pmatrix} \in \text{Hom}(F, U) \\
\text{Hom}(F, U) \ni \begin{pmatrix} \delta B_1 && \delta I \end{pmatrix} & \mapsto B_3 \begin{pmatrix} \delta B_1 && \delta I\end{pmatrix} - \begin{pmatrix} \delta B_1 && \delta I \end{pmatrix} \begin{pmatrix} B_3 && 0 \\ -J_2 && \xi_3\end{pmatrix} \\
& = \begin{pmatrix} \comm{B_3}{\delta B_1} + \delta I J_2 && B_3 \delta I - \delta I \xi_3 \end{pmatrix} \in \text{Hom}(F, U). 
\end{split}
\]
In a similar fashion, the differential of the moment map 
\[
T_{T^* \text{Rep}_Q(\vv, \ww)} \to \mathfrak{g}^*_{\FX}
\]
induces a morphism of sheaves 
\[
T^* \mathscr{R}ep(-1) \to \mathscr{H}om(\Sh{V}, \Sh{V}(-1))
\]
fitting into a similar commutative diagram as above. On local sections, the differential of the moment map acts of course by $(\delta \Sh{B}_2, \delta \Sh{B}_3, \delta \Sh{I}, \delta \Sh{J}) \mapsto \comm{\delta \Sh{B}_2}{\Sh{B}_3} + \comm{\Sh{B}_2}{\delta \Sh{B}_3} + \delta \Sh{I} \Sh{J} + \Sh{I} \delta \Sh{J}$. Following a similar line of reasoning as above, we see that the bottom right horizontal differential in \eqref{eq: quasimap tangent double complex} acts by 
\[
\begin{split}
(\delta \alpha, \delta B_2, \delta J_3, \delta B_3, \delta J_2) & \mapsto \begin{pmatrix} \delta B_2 \\ \delta J_3 \end{pmatrix} B_3 - \begin{pmatrix} B_3 && 0 \\ -J_2 && \xi_3 \end{pmatrix} \begin{pmatrix} \delta B_2 \\ \delta J_3\end{pmatrix}  \\
& + \begin{pmatrix} B_2 && 0 \\ J_3 && \xi_2 \end{pmatrix} \begin{pmatrix} \delta B_3 \\ - \delta J_2 \end{pmatrix} - \begin{pmatrix} \delta B_3 \\ - \delta J_2 \end{pmatrix} B_2 + \begin{pmatrix} 0 \\ \iota \end{pmatrix} \delta \alpha \\
& = \begin{pmatrix} \comm{\delta B_2}{B_3} + \comm{B_2}{\delta B_3} \\
\sum_{m = 1, 2} (\delta J_m B_m + J_m \delta B_m - \xi_m \delta J_m) + \iota \delta \alpha
\end{pmatrix} \in \text{Hom}(U, F). 
\end{split}
\]
Finally, let's introduce the notation 
\[
(\eta_\alpha, \eta_{B_2}, \eta_{J_3}, \eta_{B_3}, \eta_{J_2}) \in \text{Hom}(W, U) \oplus \text{Hom}(F, U) \oplus \text{Hom}(F, U)
\]
for an element of the top middle entry of \eqref{eq: quasimap tangent double complex}. Then, reasoning as above, the horizontal differential at this term acts by 
\[
\begin{split}
(\eta_\alpha, \eta_{B_2}, \eta_{J_3}, \eta_{B_3}, \eta_{J_2}) & \mapsto \begin{pmatrix} \eta_{B_3} && -\eta_{J_2} \end{pmatrix} \begin{pmatrix} B_3 && 0 \\ -J_2 && \xi_3 \end{pmatrix}- B_3 \begin{pmatrix} \eta_{B_3} && - \eta_{J_2} \end{pmatrix} \\ 
& + B_2 \begin{pmatrix} \eta_{B_2} && \eta_{J_3} \end{pmatrix} - \begin{pmatrix} \eta_{B_2} && \eta_{J_3} \end{pmatrix} \begin{pmatrix} B_2 && 0 \\ J_3 && \xi_2 \end{pmatrix} + \eta_{\alpha} \begin{pmatrix} \alpha && \jmath \end{pmatrix} \\
& = \begin{pmatrix} \comm{B_2}{\eta_{B_2}} + \comm{\eta_{B_3}}{B_3} + \eta_{J_2} J_2 - \eta_{J_3} J_3 + \eta_\alpha \alpha \\
B_3 \eta_{J_2} - \eta_{J_2} \xi_3 + B_2 \eta_{J_3} - \eta_{J_3} \xi_2  + \eta_\alpha \jmath \end{pmatrix}^T \in \text{Hom}(F, U). 
\end{split}
\]
Now that we have a complete understanding of \eqref{eq: quasimap tangent double complex}, we may compare with \eqref{eq: quiver tangent complex}. First we replace $\eqref{eq: quasimap tangent double complex}$ with a quasi-isomorphic complex, by making the following sign replacements. Reverse the sign in the $\text{Hom}(W, U)$ factor and one of the $\text{Hom}(F, U)$ factors of the top middle entry (in the above notation, $\eta_{\alpha} \mapsto - \eta_{\alpha}$, $\eta_{B_3} \mapsto -\eta_{B_3}, \eta_{J_2} \mapsto - \eta_{J_2}$), and reverse the sign of the middle vertical differential so that both squares anticommute rather than commute. Finally, reverse the overall sign in the top right $\text{Hom}(F, U)$ entry. Then, combining the results above, the differential $D$ in the total complex acts as follows. In degree $-1$, 
\[
\begin{split}
D \begin{pmatrix} X \\ Y \end{pmatrix} & = \Bigg( \alpha X + jY, \comm{B_2}{X}, J_3 X + \xi_2 Y - Y B_2, \comm{B_3}{X}, J_2 X - \xi_3 Y + Y B_3, \\
& \comm{B_1}{X} + IY, -XI \Bigg) \in \text{Hom}(U, W) \oplus \text{Hom}(U, F) \oplus \text{Hom}(U, F) \oplus \text{Hom}(F, U)
\end{split}
\]
where the first line is from the horizontal differential and the second line is from the vertical differential. In degree $0$, using an abbreviated notation, 
\[
\begin{split}
D ( \delta \alpha, \delta B, \delta J, \delta B_1, \delta I) & = \Bigg( \delta I \iota, \comm{B_3}{\delta B_1} + \delta I J_2 + \comm{\delta B_3}{B_1} + I \delta J_2, B_3 \delta I - \delta I \xi_3 + \delta B_3 I  \\
& \comm{\delta B_1}{B_2} + \delta I J_3 + \comm{B_1}{\delta B_2} + I \delta J_3, B_2 \delta I - \delta I \xi_2 + \delta B_2 I, \\
& \comm{\delta B_2}{B_3} + \comm{B_2}{\delta B_3}, \sum_m (\delta J_m B_m + J_m \delta B_m - \xi_m \delta J_m) + \iota \delta \alpha \Bigg) \\
& \in \text{Hom}(W, U) \oplus \text{Hom}(F, U) \oplus \text{Hom}(F, U) \oplus \text{Hom}(U, F)
\end{split}
\]
where we have implemented the sign replacements above (and recalled that the definition of $\eta_{J_2}$ comes with an overall minus sign). Finally, in degree $1$, 
\[
\begin{split}
D( \eta_\alpha, \eta_B, \eta_J, \eta_{B_1}, \eta_I) & = \begin{pmatrix} \comm{\eta_{B_2}}{B_2} + \comm{\eta_{B_3}}{B_3} + \eta_{J_2} J_2 + \eta_{J_3} J_3 + \eta_\alpha \alpha + \comm{\eta_{B_1}}{B_1} - I \eta_I  \\ B_3 \eta_{J_2} - \eta_{J_2} \xi_3 + \eta_{J_3} \xi_2 - B_2 \eta_{J_3} + \eta_\alpha \jmath + \eta_{B_1}I\end{pmatrix}^T \\
& \in \text{Hom}(F, U)
\end{split}
\]
where we introduced the notation $\begin{pmatrix} \eta_{B_1} & \eta_I \end{pmatrix} \in \text{Hom}(U, F)$ for an element of the bottom right entry in \eqref{eq: quasimap tangent double complex}. 

Now from the above formulas, it is obvious that $D$ is exactly the differential in \eqref{eq: quiver tangent complex}, once we recall the formula \cite{quasimapcrit} for the $GL(U) \ltimes \text{Hom}(U, V)$ action on $\text{Rep}_{Q^{\QM}}$ and explicitly compute the second derivatives of the function 
\[
W = \tr B_1 \comm{B_2}{B_3} + \tr J_2(B_2I - I \xi_2) + \tr J_3 (B_3 I - I \xi_3) + \tr I \iota \alpha. 
\]
\end{proof}
Notice that the basic principle entering the quiver description of quasimaps is a certain functoriality of Beilinson resolutions, and it is this same principle that allows us to analyze the deformation theory as above. 

\subsection{The derived structure in the Čech model}

We would now like to show that Theorem \ref{thm: Čech critical quasimap} also holds at the derived level. Arguing as above in the quiver case, this consideration is reduced to a comparison of the tangent complex of $\QM^\xi(X)$ to the natural tangent complex on the derived critical locus $\Dcrit(\Sh{W}^\xi)$. 

We follow the notations and conventions of Section \ref{subsubsect: Čech review}. Introduce in addition the notation, for a vector bundle $\mathscr{E}$ on $\mathbb{P}^1$, 
\begin{align*}
    \Gamma_+(\mathscr{E}) & \coloneqq \Gamma(\mathbb{D}_0^\times, \mathscr{E})/\Gamma(U_\infty, \mathscr{E}) \\
    \Gamma_-(\mathscr{E}) & \coloneqq \Gamma(\mathbb{D}_0^\times, \mathscr{E})/\Gamma(\mathbb{D}_0, \mathscr{E}). 
\end{align*}
Notice that if we use local trivializations and transition functions as in Section \ref{subsubsect: Čech review}, then $\Gamma_-$ is explicitly vector-valued Laurent series $s(t)$ modulo those of the form $g^{-1}(t) \cdot s_0(t)$, where $s_0(t)$ is a vector-valued Taylor series. On the other hand, using the trivializations from Section \ref{subsubsect: Čech review}, the quotient in $\Gamma_+$ involves no transition function, and consists of vector-valued Laurent series $s(t)$ modulo vector-valued polynomials $s_\infty(t)$ in $t^{-1}$.  

The following provides equivalent descriptions of the Čech complex, and is an elementary consequence of the definitions.

\begin{lemma} \label{lemma: equivalent Čech complexes}
Let $\mathscr{E}$ be a vector bundle on $\mathbb{P}^1$. Consider the following diagram of two-term complexes 
\begin{equation}
\begin{tikzcd}
    \Gamma(U_\infty, \mathscr{E}) \oplus \Gamma(\mathbb{D}_0, \mathscr{E}) \arrow[r] \arrow[d, two heads] & \Gamma(\mathbb{D}^\times_0, \mathscr{E}) \arrow[d, two heads] \\ 
    \Gamma(U_\infty, \mathscr{E}) \arrow[r] \arrow[d, hook] & \Gamma_-(\mathscr{E}) \arrow[d, hook] \\
    \Gamma(\mathbb{D}_0^\times, \mathscr{E}) \arrow[r] & \Gamma_-(\mathscr{E}) \oplus \Gamma_+(\mathscr{E}).
\end{tikzcd}
\end{equation}
All vertical maps, taken as maps of complexes, are quasi-isomorphisms. 
\end{lemma}

\subsubsection{Quasimap tangent complex in Čech language}
We now proceed as above; the strategy will again be to use \eqref{eq: quasimap tangent complex}, though this time using Čech resolutions to compute $R\Gamma$. In the notation used in Section \ref{subsubsect: Čech review},  
\[
T^* \mathscr{R}ep = \Sh{N} \oplus \Sh{N}^*. 
\]
Recall $\Sh{N}$ is the bundle on $\mathbb{P}^1$ associated to $\text{Rep}_Q(\vv, \ww)$. Now resolve the first $\mathscr{H}om(\Sh{V}, \Sh{V}(-1))$ and $\Sh{N}(-1)$ in \eqref{eq: quasimap tangent complex} by the usual Čech complex, and $\Sh{N}^*(-1)$ and the second $\mathscr{H}om(\Sh{V}, \Sh{V}(-1))$ by the bottom line of the diagram in Lemma \ref{lemma: equivalent Čech complexes}. We arrive at the following double complex for $\TT_{\DQM^\xi(X)}$: 
\begin{equation} \label{eq: quasimap tangent double complex Čech}
\begin{tikzcd}
    \Gamma(\mathbb{D}_0^\times, \Sh{E}) \arrow[r] & \Gamma(\mathbb{D}_0^\times, \Sh{N}(-1)) \oplus \Gamma_-(\Sh{N}^*(-1)) \oplus \Gamma_+(\Sh{N}^*(-1)) \arrow[r] & \Gamma_-(\Sh{E}) \oplus \Gamma_+(\Sh{E}) \\
    \Gamma(U_\infty, \Sh{E}) \oplus \Gamma(\mathbb{D}_0, \Sh{E}) \arrow[u] \arrow[r] & \Gamma(U_\infty, \Sh{N}(-1)) \oplus \Gamma(\mathbb{D}_0, \Sh{N}(-1)) \oplus \Gamma(\mathbb{D}_0^\times, \Sh{N}^*(-1)) \arrow[u] \arrow[r] & \Gamma(\mathbb{D}_0^\times, \Sh{E}). \arrow[u] 
\end{tikzcd} 
\end{equation}
We have abbreviated $\Sh{E} \coloneqq \mathscr{H}om(\Sh{V}, \Sh{V}(-1))$. 

\subsubsection{Tangent complex of $\Dcrit(\Sh{W}^\xi)$} The objective is to compare \eqref{eq: quasimap tangent double complex Čech} with the natural tangent complex on the derived critical locus $\Dcrit(\Sh{W}^\xi)$. Recall that $\Sh{M}^\xi_{G, N}$ defined in \eqref{eq: define MxiGN} is explicitly a quotient by $G(\Sh{O}) \times G_1[t^{-1}]$. We have a natural isomorphism 
\[
\text{Lie}(G(\Sh{O}) \times G_1[t^{-1}]) \simeq \Fg(\Sh{O}) \oplus t^{-1} \Fg[t^{-1}]
\]
and, via the \textit{left} trivialization, an isomorphim 
\[
T G(\Sh{K}) \simeq G(\Sh{K}) \times \Fg(\Sh{K}). 
\]
Via these isomorphisms, we have the following model for the pointwise tangent complex of $\Sh{M}^\xi_{G, N}$, concentrated in degrees $[-1, 0]$:
\[
\TT_{\Sh{M}^\xi_{G, N}} \simeq \Bigg[ 
\begin{tikzcd}
    \Fg(\Sh{O}) \oplus t^{-1} \Fg[t^{-1}] \arrow[r] & \Fg(\Sh{K}) \oplus t^{-1}N[t^{-1}] \oplus N(\Sh{O}) \oplus N^*(\Sh{K})
\end{tikzcd} \Bigg].
\]
The standard model for the tangent complex to a derived critical locus reads, as usual:
\begin{equation} \label{eq: Čech tangent complex}
\TT_{\Dcrit(\Sh{W}^\xi)} \simeq \Bigg[ 
\begin{tikzcd}
    \TT^{-1}_{\Sh{M}^\xi_{G, N}} \arrow[r] & \TT^0_{\Sh{M}^\xi_{G, N}} \arrow[r, "\text{Hess}(\Sh{W})"] & (\TT^0_{\Sh{M}^\xi_{G, N}})^\vee \arrow[r] & (\TT^{-1}_{\Sh{M}^\xi_{G, N}})^\vee
\end{tikzcd} \Bigg].
\end{equation}

\subsubsection{Comparing the complexes}
Now we explicitly compare \eqref{eq: quasimap tangent double complex Čech} and \eqref{eq: Čech tangent complex}. First we compare the individual terms, and then the differentials.

Because $\mathscr{H}om(\Sh{V}, \Sh{V}(-1))$ is viewed as the sheaf of sections of $\mathscr{E}nd(\Sh{V})$ vanishing to first order at $\infty \in \mathbb{P}^1$, trivializing $\Sh{V}$ over $U_\infty$ and $\mathbb{D}_0$ provides isomorphisms 
\begin{align*}
t^{-1} \Fg[t^{-1}] & \simeq \Gamma(U_\infty, \mathscr{H}om(\Sh{V}, \Sh{V}(-1))) \\ 
\Fg(\Sh{O}) & \simeq \Gamma(\mathbb{D}_0, \mathscr{H}om(\Sh{V}, \Sh{V}(-1)))
\end{align*}
identifying the degree $-1$ terms in the total complex of \eqref{eq: quasimap tangent double complex Čech} and \eqref{eq: Čech tangent complex}. Trivializing over the patches likewise gives isomorphisms 
\begin{align*}
\Fg(\Sh{K}) & \simeq \Gamma(\mathbb{D}_0^\times, \mathscr{H}om(\Sh{V}, \Sh{V}(-1))) \\ 
t^{-1} N[t^{-1}] & \simeq \Gamma(U_\infty, \Sh{N}(-1)) \\ 
N(\Sh{O}) & \simeq \Gamma(\mathbb{D}_0, \Sh{N}(-1)) \\ 
N^*(\Sh{K}) & \simeq \Gamma(\mathbb{D}_0^\times, \Sh{N}^*(-1)) 
\end{align*}
identifying the degree 0 part of the total complex of \eqref{eq: quasimap tangent double complex Čech} with the degree 0 part of \eqref{eq: Čech tangent complex}. 

In a similar fashion, one finds 
\begin{align*}
    (N^*(\Sh{K}))^\vee & \simeq N(\Sh{K}) \simeq \Gamma(\mathbb{D}_0^\times, \Sh{N}(-1)) \\ 
    (N(\Sh{O}))^\vee & \simeq N^*(\Sh{K})/g^{-1}(t) \cdot N^*(\Sh{O}) \simeq \Gamma_-(\Sh{N}^*(-1)) \\
    (t^{-1}N[t^{-1}])^\vee & \simeq N^*(\Sh{K})/t^{-1} N^*[t^{-1}] \simeq \Gamma_+(\Sh{N}^*(-1)) \\
    (\Fg(\Sh{K}))^\vee & \simeq \Fg(\Sh{K}) \simeq \Gamma(\mathbb{D}_0^\times, \mathscr{H}om(\Sh{V}, \Sh{V}(-1))). 
\end{align*}
The first isomorphism in each line comes from the residue pairing and trace, while the second isomorphism follows from trivializing over patches. Notice that the transition function appears in the expression for $N(\Sh{O})^\vee$ because the duality must be $G(\Sh{O}) \times G_1[t^{-1}]$-equivariant, and our convention is that $G_1[t^{-1}]$ acts on $N^*(\Sh{K})$. For a similar reason, the left trivialization of $TG(\Sh{K})$ is preferred above. This identifies the degree $1$ part of the total complex of \eqref{eq: quasimap tangent double complex Čech} with the degree $1$ part of \eqref{eq: Čech tangent complex}. The degree 2 parts are also identified following this logic. 

In fact, the differentials agree as well, as expressed via the following

\begin{proposition} \label{prop: Čech tangent complex=quasimap tangent complex}
The total complex of the double complex \eqref{eq: quasimap tangent double complex Čech} is isomorphic to \eqref{eq: Čech tangent complex}. 
\end{proposition}
Once again, the proof is by a direct computation of the differentials in both complexes. 

\begin{proof}
    We begin at the bottom left corner of \eqref{eq: quasimap tangent double complex Čech}. That the horizontal differential at this entry agrees with \eqref{eq: Čech tangent complex} is clear by construction (namely, the horizontal differential in \eqref{eq: quasimap tangent double complex Čech} at the leftmost column is defined as the action map, and all isomorphims above are $ G(\Sh{O}) \times G_1[t^{-1}]$-equivariant). It remains to check that the vertical Čech differential in \eqref{eq: quasimap tangent double complex Čech} agrees with the infinitesimal action of $G(\Sh{O}) \times G_1[t^{-1}]$ on $G(\Sh{K})$ in the left trivialization from \eqref{eq: Čech tangent complex}. Recalling that $G(\Sh{O}) \times G_1[t^{-1}] \ni (g_0(t), g_\infty(t))$ acts on $G(\Sh{K})$ by $g(t) \mapsto g_0^{-1}(t) g(t) g_\infty(t)$, the infinitesimal action in the left trivialization is 
    \[
    \begin{split}
    \Fg(\Sh{O}) \oplus t^{-1} \Fg[t^{-1}] & \ni (\epsilon_\infty(t), \epsilon_0(t))  \\
    & \mapsto g^{-1}(t) ( - \epsilon_0(t) g(t) + g(t) \epsilon_\infty(t))
    = \epsilon_\infty(t) - g^{-1}(t) \epsilon_0(t) g(t)
    \end{split}
    \]
    which clearly agrees with the Čech differential on $\mathscr{H}om(\Sh{V}, \Sh{V}(-1))$.

    Moving to degree $0$, the differential is a sum of two pieces. Denote an element of the middle entry of the bottom row of \eqref{eq: quasimap tangent double complex Čech} by $(\delta \Sh{B}_{3, \infty}, \delta \Sh{J}_\infty, \delta \Sh{B}_{3, 0}, \delta \Sh{J}_0, \delta \Sh{B}_2, \delta \Sh{I})$, the vertical differential is a direct sum of maps sending 
    \begin{equation} \label{eq: degree 0 Čech vertical differential}
    \begin{split}
    & \Gamma(U_\infty, \Sh{N}(-1)) \oplus \Gamma(\mathbb{D}_0, \Sh{N}(-1)) \ni (\delta \Sh{B}_{3, \infty}(t), \delta \Sh{J}_\infty(t), \delta \Sh{B}_{3, 0}, \delta \Sh{J}_0(t)) \\
    & \mapsto (\delta \Sh{B}_{3, \infty}(t) - g^{-1}(t) \delta \Sh{B}_{3, 0} g(t), \delta \Sh{J}_\infty(t) - \delta \Sh{J}_0(t) g(t)) \in \Gamma(\mathbb{D}_0^\times, \Sh{N}(-1)) \\
    & \Gamma(\mathbb{D}_0^\times, \Sh{N}^*(-1)) \ni (\delta \Sh{B}_2(t), \delta \Sh{I}(t)) \\
    & \mapsto ([\delta \Sh{B}_2(t), \delta \Sh{I}(t)]_-, [\delta \Sh{B}_2(t), \delta \Sh{I}(t)]_+) \in \Gamma_-(\Sh{N}^*(-1)) \oplus \Gamma_+(\Sh{N}^*(-1)). 
    \end{split}
    \end{equation}
    The horizontal differential in degree 0 sends 
    \begin{equation} \label{eq: degree 0 Čech horizontal differential}
    \begin{split}
        & \Gamma(\mathbb{D}_0^\times, \mathscr{H}om(\Sh{V}, \Sh{V}(-1))) \ni \epsilon(t) \\
        & \mapsto (\comm{\Sh{B}_{3, \infty}(t)}{\epsilon(t)}, \Sh{J}_\infty(t) \epsilon(t), [ \comm{\Sh{B}_2(t)}{\epsilon(t)}, - \epsilon(t) \Sh{I}(t)]_-, 0) \\
        & \in \Gamma(\mathbb{D}_0^\times, \Sh{N}(-1)) \oplus \Gamma_-(\Sh{N}^*(-1)) \oplus \Gamma_+(\Sh{N}^*(-1)) \\
        & \Gamma(U_\infty, \Sh{N}(-1)) \oplus \Gamma(\mathbb{D}_0, \Sh{N}(-1)) \oplus \Gamma(\mathbb{D}_0^\times, \Sh{N}^*(-1))) \ni (\delta \Sh{B}_{3, \infty}, \delta \Sh{J}_\infty, \delta \Sh{B}_{3, 0}, \delta \Sh{J}_0, \delta \Sh{B}_2, \delta \Sh{I}) \\
        & \mapsto \comm{\delta \Sh{B}_2}{\Sh{B}_{3, \infty}} + \comm{\Sh{B}_2}{\delta \Sh{B}_{3, \infty}} + \delta \Sh{I} \Sh{J}_\infty + \Sh{I} \delta \Sh{J}_\infty \\ 
        & \in \Gamma(\mathbb{D}_0^\times, \mathscr{H}om(\Sh{V}, \Sh{V}(-1))). 
    \end{split}
    \end{equation}
    The first three lines describe the differential at the top left of \eqref{eq: quasimap tangent double complex Čech}, while the next three lines describe the differential at the middle entry of the bottom row. It may be somewhat surprising that the $\Gamma_+$ component of the top map vanishes, but this is required e.g. to make the left square of \eqref{eq: quasimap tangent double complex Čech} commute. 
    
    Onward to degree $1$, the differential is a sum of two contributions. The rightmost vertical map in \eqref{eq: quasimap tangent double complex Čech} is simply the natural projection
    \begin{equation} \label{eq: degree 1 Čech vertical differential}
    \Gamma(\mathbb{D}_0^\times, \mathscr{H}om(\Sh{V}, \Sh{V}(-1)) \ni \eta_g(t) \mapsto ( [\eta_g(t)]_-, [\eta_g(t)]_+). 
    \end{equation}
    The horizontal differential at the top middle entry in \eqref{eq: quasimap tangent double complex Čech} acts as follows. Introduce notation 
    \[
    \begin{split}
        (\eta_2(t), \eta_{\Sh{I}}(t)) & \in \Gamma(\mathbb{D}_0^\times, \Sh{N}(-1)) \\
        (\tilde \eta_{3, 0}(t),  \tilde \eta_{\Sh{J}_0}(t)) & \in \Gamma_-(\Sh{N}^*(-1)) \\
        (\eta_{3, \infty}(t), \eta_{\Sh{J}_\infty}(t)) & \in \Gamma_+(\Sh{N}^*(-1)). 
    \end{split}
    \]
    Then the differential sends such a tuple to 
    \begin{equation} \label{eq: degree 1 Čech horizontal differential}
    \begin{split}
        \Big( [ \comm{\Sh{B}_2}{\eta_2} + \comm{\tilde \eta_{3, 0}}{\Sh{B}_{3, \infty}} + \Sh{I} \eta_{\Sh{I}} + \tilde \eta_{\Sh{J}_0} \Sh{J}_\infty ]_-, [ \comm{\eta_{3, \infty}}{\Sh{B}_{3, \infty}} + \eta_{\Sh{J}_\infty} \Sh{J}_\infty]_+ \Big). 
    \end{split}
    \end{equation}
    Notice that this is well-defined by the equations $\Sh{B}_{3, \infty} = g^{-1} \Sh{B}_{3, 0} g$, $\Sh{J}_\infty = \Sh{J}_0 g$, and that the right square of \eqref{eq: quasimap tangent double complex Čech} commutes using the other critical locus conditions. 

    Now we make sign replacements in \eqref{eq: quasimap tangent double complex Čech}. Reverse the sign of the middle vertical differential, and reverse the sign of every factor in the top middle entry (in the above notations, every $\eta$ variable is replaced with $-\eta$). Then reverse the sign in the $\Gamma_+$ factor at the top right corner of \eqref{eq: quasimap tangent double complex Čech}. The total differential $D$ in \eqref{eq: quasimap tangent double complex Čech} now takes the form 
    \[
    \begin{split}
        &D( \epsilon_0(t), \epsilon_\infty(t))  = (\epsilon_\infty(t) - g^{-1}(t) \epsilon_0(t) g(t); \comm{\Sh{B}_{3, \infty}(t)}{\epsilon_\infty(t)}, \Sh{J}_\infty(t) \epsilon_\infty(t); \\
        & \comm{\Sh{B}_{3, 0}(t)}{\epsilon_0(t)}, \Sh{J}_0(t) \epsilon_0(t); \comm{\Sh{B}_2(t)}{\epsilon_\infty(t)}, - \epsilon_\infty(t) \Sh{I}(t)) \\
        & \in \Gamma(\mathbb{D}_0^\times, \Sh{E}) \oplus \Gamma(U_\infty, \Sh{N}(-1)) \oplus \Gamma(\mathbb{D}_0, \Sh{N}(-1)) \oplus \Gamma(\mathbb{D}_0^\times, \Sh{N}^*(-1)). \\
        & D(\epsilon; \delta \Sh{B}_{3, \infty}, \delta \Sh{J}_\infty; \delta \Sh{B}_0, \delta \Sh{J}_0; \delta \Sh{B}_2, \delta \Sh{I})  = ( \delta \Sh{B}_{3, \infty} - g^{-1} \delta \Sh{B}_{3, 0} g - \comm{\Sh{B}_{3, \infty}}{\epsilon}, \delta \Sh{J}_\infty - \delta \Sh{J}_0 g - \delta \Sh{J}_\infty \epsilon; \\
        & [\delta \Sh{B}_2 - \comm{\Sh{B}_2}{\epsilon} , \delta \Sh{I} + \epsilon \Sh{I}]_-; [\delta \Sh{B}_2, \delta \Sh{I}]_+; \comm{\delta \Sh{B}_2}{\Sh{B}_{3, \infty}} + \comm{\Sh{B}_2}{\delta \Sh{B}_{3, \infty}} + \delta \Sh{I} \Sh{J}_\infty + \Sh{I} \delta \Sh{J}_\infty) \\
        & \in \Gamma(\mathbb{D}_0^\times, \Sh{N}(-1)) \oplus \Gamma_-(\Sh{N}^*(-1)) \oplus \Gamma_+(\Sh{N}^*(-1)) \oplus \Gamma(\mathbb{D}_0^\times, \Sh{E}). \\
        & D( \eta_2, \eta_{\Sh{I}}; \tilde \eta_{3, 0}, \tilde \eta_{\Sh{J}_0}; \eta_{3, \infty}, \eta_{\Sh{J}_\infty}; \eta_g) = ( [\eta_g + \comm{\eta_2}{\Sh{B}_2} - \Sh{I} \eta_{\Sh{I}} - \comm{\tilde \eta_{3, 0}}{\Sh{B}_{3, \infty}} - \tilde \eta_{\Sh{J}_0} \Sh{J}_\infty]_-; \\
        & [ \comm{\eta_{3, \infty}}{\Sh{B}_{3, \infty}} + \eta_{\Sh{J}_\infty} \Sh{J}_\infty - \eta_g]_+) \in \Gamma_-(\Sh{E}) \oplus \Gamma_+(\Sh{E})
    \end{split}
    \]
    where we recall the shorthand $\Sh{E} = \mathscr{H}om(\Sh{V}, \Sh{V}(-1))$. This is essentially manifestly the same as \eqref{eq: Čech tangent complex} upon using the critical locus condition $\delta \Sh{W} = 0$ to rewrite $\text{Hess}(\Sh{W})$ as follows:
    \[
    \begin{split}
        \delta(\Sh{B}_{3, \infty} - g^{-1} \Sh{B}_{3, 0} g) & = \delta\Sh{B}_{3, \infty} - g^{-1} \delta \Sh{B}_{3, 0} g - \comm{g^{-1} \Sh{B}_{3, 0}g}{g^{-1} \delta g} \\ 
        & = \delta \Sh{B}_{3, \infty} - g^{-1} \delta \Sh{B}_{3, 0} g - \comm{\Sh{B}_{3, \infty}}{g^{-1} \delta g}
    \end{split}
    \]
    and
    \[
    \begin{split}
        \delta (g \Sh{B}_{2} g^{-1}, g \Sh{I}) &  \mod N^*(\Sh{O}) \\
        = (g \delta \Sh{B}_2 g^{-1} + g( \comm{g^{-1} \delta g}{\Sh{B}_2})g^{-1}, g (\delta \Sh{I} + g^{-1} \delta g \Sh{I}) ) & \mod N^*(\Sh{O}) \\
        \mapsto  (\delta \Sh{B}_2 + \comm{g^{-1} \delta g}{\Sh{B}_2}, \delta \Sh{I} + g^{-1} \delta g \Sh{I}) &  \in \Gamma_-(\Sh{N}^*(-1)).
    \end{split}
    \]
    Notice moreover that the naive obstruction coordinate $\eta'_g$ is related to our $\eta_g$ by $\eta'_g  = \eta_g + \comm{\eta_2}{\Sh{B}_2} - \Sh{I} \eta_{\Sh{I}}$, accounting for the expressions we use for the differential at the degree $1$ term. This completes the proof.
\end{proof}

\begin{remark}
Notice that the tangent complex \eqref{eq: quasimap tangent double complex Čech} is in turn quasi-isomorphic to the tangent complex to any sufficiently large finite approximation $\Dcrit(\Sh{W}^{l, m})$, essentially by a linearized version of the proof of Proposition \ref{prop: approx crit Čech}. 
\end{remark}

\subsection{Orientation data} \label{ap: orientations}
The intrinsic DT sheaf $\varphi_{\QM^\xi(X)}$ on the moduli space of quasimaps depends not only on its $(-1)$-shifted symplectic structure, but also on a choice of \textit{orientation data}. An orientation data is a choice of square root of the virtual canonical bundle 
\[
\Sh{K}_{\text{vir}}  = \det(\mathbb{T}_{\QM^\xi(X)})^{-1}
\]
in $\text{Pic}(\QM^\xi(X))$. A canonical choice is provided by Serre duality and the polarization $N = \text{Rep}_Q(\vv, \ww) \subset T^* \text{Rep}_Q(\vv, \ww)$ afforded by the zero section. From \eqref{eq: quasimap tangent complex}, we have 
\[
\det(\mathbb{T}_{\QM^\xi(X)})^{-1} \simeq \frac{(\det \text{Ext}^\bullet(\Sh{V}, \Sh{V}(-1))^2}{\det H^\bullet(\mathbb{P}^1, \Sh{N}(-1)) \otimes \det H^\bullet(\mathbb{P}^1, \Sh{N}^*(-1))}. 
\]
By Serre duality, 
\[
\det H^\bullet(\mathbb{P}^1, \Sh{N}^*(-1)) = \det H^\bullet(\mathbb{P}^1, \Sh{N}(-1)),
\]
therefore we have a canonical choice of square root\footnote{Notice this also holds equivariantly since we work in the Calabi-Yau specialization throughout this paper.}
\[
\Sh{K}_{\text{vir}}^{1/2} \coloneqq \frac{\det \text{Ext}^\bullet(\Sh{V}, \Sh{V}(-1))}{\det H^\bullet(\mathbb{P}^1, \Sh{N}(-1))}. 
\]
On the other hand, consider the derived critical locus $\Dcrit(W)$ of a function $W: M \to \mathbb{A}^1$ on a smooth stack $M$. In this case a canonical orientation data is given by 
\[
\Sh{K}_{\text{vir}}^{1/2} \coloneqq \det( \mathbb{T}_M |_{\Dcrit(W)})^{-1}.
\]
We would like to show that the canonical orientation data in the quiver or Čech model agrees with the intrinsic orienation fixed by $\Sh{N}$. 

\begin{lemma}
In the quiver model, both orientations agree. 
\end{lemma}

In this proof, we explicitly indicate the equivariant parameters $t_1, t_a, t_{a^*}$. 

\begin{proof}
    Notice that, using Beilinson resolution to compute derived global sections as in \eqref{eq: quasimap tangent double complex} above gives us relations 
    \[
    \begin{split}
        \text{Ext}^\bullet(\Sh{V}, \Sh{V}(-1)) & = \sum_{i \in Q_0} \text{Hom}(U_i, F_i) - t_1 \text{Hom}(F_i, U_i) \\ 
        H^\bullet(\mathbb{P}^1, \Sh{N}(-1)) & = \sum_{i \in Q_0} \text{Hom}(U_i, W_i) + \sum_{a \in Q_1} t_a(\text{Hom}(U_{s(a)}, F_{t(a)}) - t_1 \text{Hom}(F_{s(a)}, U_{t(a)}))
    \end{split}
    \]
    in $K_T(\QM^\xi_{\dd}(X))$. On the other hand for ambient space $M = \text{Rep}_{Q^{\QM}_{\vv, \ww}}(\dd)/GL(U) \ltimes \text{Hom}(U, V)$ we have 
    \[
    \begin{split}
    \mathbb{T}_M & = \sum_{i \in Q_0} t_1 \text{Hom}(F_i, U_i) - \text{Hom}(U_i, F_i) \\
    & + \sum_{i \in Q_0} \text{Hom}(U_i, W_i) + \sum_{a \in Q_1} (t_a \text{Hom}(U_{s(a)}, F_{t(a)}) + t_{a^*} \text{Hom}(U_{t(a)}, F_{s(a)}). 
    \end{split}
    \]
    Using that we are at the Calabi-Yau specialization so that $t_1 t_a t_{a^*} = 1$ for each edge $a$ with opposite $a^*$, we see that the ratio of orientations is 
    \[
    \prod_{a \in Q_1} \det t_{a^*} \text{Hom}(U_{t(a)}, F_{s(a)}) \otimes \det t_1t_a \text{Hom}(F_{s(a)}, U_{t(a)}) = \Sh{O}_{\QM^\xi(X)}.
    \]
\end{proof}

\begin{lemma}
    In the Čech model, both orientations agree. 
\end{lemma}

During the course of the proof, we will manipulate determinant lines of formally infinite rank vector bundles. The may be understood via finite type approximations as in Section \ref{subsec: Finite type approximations of the Čech model}. 

\begin{proof}
    For the ambient space in the Čech model we have 
    \[
    \begin{split}
    \mathbb{T}_{\Sh{M}^\xi_{G, N}} & = \Gamma(U_\infty, \Sh{N}(-1)) + \Gamma(\mathbb{D}_0, \Sh{N}(-1)) + \Gamma(\mathbb{D}_0^\times, \Sh{N}^*(-1)) \\
    & + \Gamma(\mathbb{D}_0^\times, \mathscr{E}nd(\Sh{V})(-1)) - \Gamma(U_\infty, \mathscr{E}nd(\Sh{V})(-1)) - \Gamma(\mathbb{D}_0, \mathscr{E}nd(\Sh{V})(-1)).
    \end{split}
    \]
    Via residue we have 
    \[
    \Gamma(\mathbb{D}_0^\times, \Sh{N}^*(-1)) \simeq \Gamma(\mathbb{D}_0^\times, \Sh{N}(-1))^\vee
    \]
    whence 
    \[
    \det \mathbb{T}_{\Sh{M}^\xi_{G, N}} \eval_{\Dcrit(\Sh{W}^\xi)} = \frac{\det H^\bullet(\mathbb{P}^1, \Sh{N}(-1))}{\det \text{Ext}^\bullet(\Sh{V}, \Sh{V}(-1))}
    \]
    as desired. 
    
\end{proof}

\end{appendix}

\newpage

\printbibliography

\end{document}